\documentclass[leqno,12pt]{book}

\usepackage{makeidx}

\usepackage[noautomatic]{imakeidx}\makeindex[title=Index of Terminology]
\makeindex[name=Oct2026-sorted-notation, title=Index of Notations]
\usepackage[bbgreekl]{mathbbol}
\usepackage{amsfonts}
\DeclareSymbolFontAlphabet{\mathbb}{AMSb} 
\DeclareSymbolFontAlphabet{\mathbbl}{bbold}
\usepackage{lmodern}
\usepackage{amsmath}
\usepackage{accents}
\usepackage{bm}
\RequirePackage{etoolbox}
\usepackage{amscd}

\usepackage{stmaryrd}
\usepackage[all]{xy}
\xyoption{rotate}
\xyoption{pdf}
\usepackage{tikz-cd}
\tikzset{mystyle/.style={sloped, anchor=south}}
\usepackage{amssymb}
\usepackage{graphicx} 
\usepackage{mathrsfs}
\usepackage{amsthm}
\usepackage{enumerate}
\usepackage{leftidx}
\usepackage{extarrows}
\usepackage{wasysym}
\usepackage{indentfirst}
\usepackage[OT2,T1]{fontenc}
\usepackage{epstopdf}
\usepackage {hyperref}
\hypersetup{
	colorlinks=false,       
	linkcolor=red,          
	citecolor=brown,        
}

\numberwithin{equation}{subsection}
\theoremstyle{theorem}

\newtheorem{thm}[subsection]{Theorem}
\newtheorem{introthm}{Theorem}
\newtheorem{introprop}{Proposition}
\newtheorem*{thm*}{Theorem}
\newtheorem{cor}[subsection]{Corollary}
\newtheorem{lem}[subsection]{Lemma}

\newtheorem{conj}[subsection]{Conjecture}
\newtheorem{introconj}{Conjecture}
\newtheorem{prop}[subsection]{Proposition}
\theoremstyle{definition}
\newtheorem{defn}[subsection]{Definition}

\theoremstyle{remark}
\newtheorem{rem}[subsection]{Remark}

\DeclareMathOperator{\adm}{a}
\DeclareMathOperator{\ab}{ab}
\DeclareMathOperator{\ad}{ad}
\DeclareMathOperator{\alg}{alg}
\DeclareMathOperator{\Aut}{Aut}
\DeclareMathOperator{\B}{B}

\DeclareMathOperator{\bij}{{\mathscr B} }
\DeclareMathOperator{\clsc}{\mathrm{cla}}
\DeclareMathOperator{\cls}{cls}
\DeclareMathOperator{\can}{can}
\DeclareMathOperator{\ch}{ch}
\DeclareMathOperator{\cok}{coker}

\DeclareMathOperator{\cris}{cris}
\DeclareMathOperator{\Crys}{Crys}

\DeclareMathOperator{\mot}{mot}
\DeclareMathOperator{\neu}{n}
\DeclareMathOperator{\ntr}{neu}

\DeclareMathOperator{\diag}{diag}
\DeclareMathOperator{\disc}{disc}

\DeclareMathOperator{\dR}{dR}

\DeclareMathOperator{\et}{\acute{e}t }
\DeclareMathOperator{\Gal}{Gal}
\DeclareMathOperator{\coh}{\mathbf {H}}
\DeclareMathOperator{\tatecoh}{\widehat{\mathbf {H}}}
\DeclareMathOperator{\Hom}{Hom}

\DeclareMathOperator{\id}{id}
\DeclareMathOperator{\idfunc}{\mathbbl 1}
\DeclareMathOperator{\Int}{Int}

\DeclareMathOperator{\SIsom}{\underline{Isom}}
\DeclareMathOperator{\im}{im}

\DeclareMathOperator{\Mod}{Mod}
\DeclareMathOperator{\Mot}{Mot}
\DeclareMathOperator{\RZ}{RZ}
\DeclareMathOperator{\N}{N}

\DeclareMathOperator{\topo}{top}
\DeclareMathOperator{\Pic}{Pic}

\DeclareMathOperator{\Res}{Res}

\DeclareMathOperator{\sa}{sa}

\DeclareMathOperator{\Sh}{Sh}

\DeclareMathOperator{\tors}{tors}
\DeclareMathOperator{\Spec}{Spec}

\DeclareMathOperator{\SL}{SL}
\DeclareMathOperator{\cpt}{cpt}
\DeclareMathOperator{\inv}{inv}
\DeclareMathOperator{\reg}{reg}

\DeclareMathOperator{\ur}{ur}

\DeclareMathOperator{\sconn}{sc}
\DeclareMathOperator{\semi}{ss}
\DeclareMathOperator{\der}{der}

\DeclareMathOperator{\vol}{vol}

\DeclareMathOperator{\Fil}{Fil}
\DeclareMathOperator{\spec}{spec}

\DeclareMathOperator{\Out}{Out}
\DeclareMathOperator{\geom}{geom}

\DeclareMathOperator{\GL}{GL}

\DeclareMathOperator{\Mult}{Mult}
\DeclareMathOperator{\Tori}{Tori}
\DeclareMathOperator{\Rep}{Rep}

\newcommand{\asim}[1]{\accentset{\sim}{#1}} 
\newcommand{\AM}{\mathcal{AM}}

\newcommand{\bD}{\mathbf{D}}

\newcommand{\ica}{{\mathbb I}}
\newcommand{\cca}{{\mathbb J}}
\newcommand{\lprod}[1]{\langle #1 \rangle}
\newcommand{\Abshv}{\mathcal{A}b\mathcal{S}hv}
\newcommand{\Ch}{\mathcal{C}h^{\flat}}
\newcommand{\dercat}[1]{\mathcal{D}^b(#1)}
\newcommand{\dmu}{\bbmu}

\newcommand{\idrho}{\sideset{_\rho}{}{\idfunc}}
\newcommand{\idlambdarho}{\sideset{_{\rho^{\lambda}}}{}{\idfunc}}

\newcommand{\idrholambda}{\sideset{_\rho}{^\lambda}{\idfunc}}

\newcommand{\tauhome}[1]{I_{#1} (\A_f) \backslash I_{#1} ^{\ad} (\A_f)/ I_{#1} ^{\ad} (\QQ)}

\newcommand{\adele}{\mathbb A}
\newcommand{\Qbar}{\overline \QQ}
\newcommand{\abs}[1]{\left\vert#1\right\vert}
\newcommand{\set}[1]{\left\{#1\right\}}

\newcommand{\tmar}{{\mathcal {M}ar}}

\newcommand{\Modfp}{{\mathrm {Mod}}^{\mathrm{fp}}}
\newcommand{\motinv}{\mathscr{M}}
\newcommand{\zzeta}{\zeta^{p,\infty}}

\newcommand{\site}{\mathsf{S}}
\newcommand{\spK}{\mathbb{K}}
\newcommand{\gbBC}{{\mathrm{PB}}}

\newcommand{\To}{\longrightarrow}
\newcommand{\isom}{%
	\ifbool{@display}{\overset{\sim}{\longrightarrow}}{\xrightarrow\sim}%
}
\newcommand{\lix}{\leftidx}
\newcommand{\A}{\mathbb{A}}

\newcommand{\NN}{\ZZ_{\geq 1} }

\newcommand{\ZZ}{\mathbb{Z}}
\newcommand{\QQ}{\mathbb{Q}}
\newcommand{\RR}{\mathbb{R}}
\newcommand{\II}{\mathscr {I}}
\newcommand{\JJ}{\mathscr {J}}
\newcommand{\FF}{\mathbb{F}}

\newcommand{\CC}{\mathbb{C}}
\newcommand{\GG}{\mathbb{G}}
\newcommand{\DD}{\mathbb{D}}
\newcommand{\oo}{\mathcal{O}}
\renewcommand{\O}{\mathcal{O}}
\newcommand{\D}{\mathfrak{D}}
\newcommand{\G}{\mathfrak{G}}
\newcommand{\E}{\mathfrak{E}}
\newcommand{\F}{\mathfrak{scc}}

\newcommand{\Qf}{\mathfrak{Q}}
\newcommand{\gP}{\mathfrak{P}}
\newcommand{\gQ}{\mathfrak{Q}}
\newcommand{\gG}{\mathfrak{G}}
\newcommand{\gS}{\mathfrak{S}}
\newcommand{\gH}{\mathfrak{H}}
\newcommand{\gM}{\mathfrak{M}}
\newcommand{\Shh}{\mathscr{S}}
\newcommand{\SigmaREllip}{\Sigma_{\RR \text{-} \el}}

\newcommand{\Newton}{\mathcal N}
\newcommand{\Gm}{\mathbb G_m}

\newcommand{\LL}{\breve{\QQ}_p}

\newcommand{\wa}[1]{\mathrm{MF}^{\varphi, \mathrm{a}} _{#1}}
\newcommand{\ignore}[1]{}
\newcommand{\ggp}{\mathcal{LRP}^{\mathrm{gg}}}
\newcommand{\wgp}{\mathcal{LRP}}
\newcommand{\spd}{\mathcal{SPD}}
\newcommand {\lps}{[\![}
\newcommand {\rps}{]\!]}

\newcommand{\Tstr}{{\mathfrak T}} 
\newcommand{\KTstr}{{\mathfrak{KT}}^{\mathrm{str}}}

\newcommand{\KT}{{\mathfrak{KT}}}

\newcommand{\KP}{{\mathfrak{KP}}}

\newcommand{\redgb}{\mathcal{G}rb}
\newcommand{\proredgb}{{\mathrm{pro}\text{-}\mathcal{G}rb}}
\newcommand{\proR}{{\mathrm{pro}\text{-}\mathcal{R}}}
\DeclareSymbolFont{cyrletters}{OT2}{wncyr}{m}{n}
\DeclareMathSymbol{\Sha}{\mathalpha}{cyrletters}{"58}

\newcommand{\cB}{{\mathcal{B}}}

\newcommand{\cD}{{\mathcal{D}}}
\newcommand{\cE}{{\mathcal{E}}}
\newcommand{\cF}{{\mathcal{F}}}

\newcommand{\cR}{{\mathcal{R}}}
\newcommand{\cS}{{\mathcal{S}}}
\newcommand{\cZ}{{\mathcal{Z}}}
\newcommand{\cA}{{\mathcal{A}}}

\newcommand{\fkc}{{\mathfrak c}}
\newcommand{\fkH}{{\mathfrak H}}
\newcommand{\fkD}{{\mathfrak D}}
\newcommand{\fke}{{\mathfrak e}}

\newcommand{\fkg}{{\mathfrak g}}
\newcommand{\fkk}{{\mathfrak k}}
\newcommand{\fkm}{{\mathfrak m}}

\newcommand{\fkp}{{\mathfrak p}}

\newcommand{\fkE}{{\mathfrak E}}
\newcommand{\fkK}{{\mathfrak K}}
\newcommand{\fks}{{\mathfrak s}}

\newcommand{\fkX}{{\mathfrak X}}
\newcommand{\fkZ}{{\mathfrak Z}}

\newcommand{\mL}{{\mathscr L}}

\newcommand{\mR}{{\mathscr R}}

\newcommand{\Z}{\mathbb{Z}}
\newcommand{\Q}{\mathbb{Q}}

\newcommand{\R}{\mathbb{R}}

\newcommand{\el}{\mathrm{ell}}

\newcommand{\dertau}{\underline {\tau}}
\newcommand{\dersigma}{\underline {\sigma}}
\newcommand{\Cent}{\mathrm{Cent}}
\newcommand{\tomega}{\Upsilon}

\newcommand{\circsim}{\overset{\circ}{\sim}}

\newcommand{\bs}{\backslash}

\newcommand{\ltensor}{\overset{\mathrm{L}}{\otimes}}
\newcommand{\RX}{\mathrm{R}X}

\newcommand{\tr}{{\mathrm {tr}}\,}

\newcommand{\Isoc}{{\mathrm{Isoc}}}
\newcommand{\GIsoc}{{\mathcal G\text{-}\mathrm{Isoc}}}

\newcommand{\GMod}{{\Gamma\text{-} \mathrm{Mod}}}

\newcommand{\GFMod}{{\Gamma_F\text{-}\mathrm{Mod}}}

\newcommand{\cM}{{\mathcal{M}}}
\newcommand{\cP}{{\mathcal{P}}}
\newcommand{\cT}{{\mathcal{T}}}

\newcommand{\cG}{{\mathcal{G}}}
\newcommand{\cH}{{\mathcal{H}}}
\newcommand{\cI}{\mathcal{I}}

\newcommand{\cO}{\mathcal{O}}
\newcommand{\cV}{\mathcal{V}}

\newcommand{\Frac}{\mathrm {Frac\,}}
\newcommand{\Fpbar}{\overline{\mathbb F}_p}
\newcommand{\Qpbar}{\overline{\mathbb Q}_p}

\newcommand{\Qpur}{\Q_p^{\ur}}
\newcommand{\Zpur}{\Z_p^{\ur}}
\newcommand{\Qlbar}{\overline{\mathbb Q}_l}

\newcommand{\Fix}{\mathcal {FIX}}

\renewcommand{\hat}{\widehat}

\renewcommand{\lg}{\langle}
\newcommand{\rg}{\rangle}
\newcommand{\hra}{\hookrightarrow}
\newcommand{\ra}{\rightarrow}

\newcommand{\ol}{\overline}
\renewcommand{\tilde}{\widetilde}

\newcommand{\benu}{\begin{enumerate}}
	\newcommand{\eenu}{\end{enumerate}}
\newcommand{\beq}{\begin{equation}}
	\newcommand{\eeq}{\end{equation}}
\newcommand{\bit}{\begin{itemize}}
	\newcommand{\eit}{\end{itemize}}

\title{The stable trace formula for Shimura varieties of abelian type} 
\author{Mark Kisin, Sug Woo Shin, Yihang Zhu}
\date{November 2026}
\begin{document}
	
\maketitle
 
  \frontmatter 
  
  \chapter*{Abstract}
  \addcontentsline{toc}{chapter}{Abstract}

We express the Frobenius--Hecke traces on the compactly supported cohomology of a Shimura variety of abelian type in terms of elliptic parts of stable Arthur--Selberg trace formulas for the endoscopic groups. This confirms predictions of Langlands and Kottwitz at primes where the level is hyperspecial.

  \cleardoublepage          
  \thispagestyle{empty}     
  \vspace*{\fill}
  \begin{flushright}
  	\begin{minipage}{0.65\textwidth}
  		\itshape
  	 	{For various reasons it is to be expected that the L-functions
  		associated to these Shimura varieties can be expressed in terms of the L-functions associated
  		to automorphic forms on the group defining the variety and on certain related groups.}
  		\par\medskip
  		\normalfont\hfill --- Robert P. Langlands, \textit{Some Contemporary Problems with Origins in the Jugendtraum}
  	\end{minipage}
  \end{flushright}
  \vspace*{\fill}
  \cleardoublepage
   
  \chapter*{Preface}

  This book grew out of a collaboration that began in the spring of 2014 in Cambridge, Massachusetts, where two of us were faculty members and the third a graduate student. The initial manuscript was submitted in October 2021, and revisions were finalized in 2026, with the three of us spread across 
  two sides of the North American continent and of the Pacific. 
  
  During the revision period, we uncovered and corrected numerous issues beyond those identified by the referees, often finding that resolving them led to clearer and more conceptual proofs. Beyond correcting flaws, we also discovered several major simplifications of existing arguments, enabling us to dispense with cumbersome constructions. We hope that these improvements go some way toward justifying the extended revision process.
  
  Spanning over a decade and across widely scattered locations, this collaborative endeavor accompanied major transitions in our personal and professional lives, and continued to be a constant factor in whatever measure in a changing world.  Arriving at its completion in 2026, one also senses that our subject itself may be entering an entirely new era. Cherished as it has been, this unforgettable chapter has reached a long-awaited close, leaving us to follow the road that surely continues, wherever it may lead.

  \begin{flushright}
  	The Authors\\
  	November 2026
  \end{flushright}

  \chapter*{Acknowledgements}
  \addcontentsline{toc}{chapter}{Acknowledgements}
  We thank R.~Zhou for many inspiring discussions and T.~Haines for answering our questions on parahoric group schemes. We are grateful to M.~Rapoport for intriguing conversations on the Langlands--Rapoport Conjecture and for useful suggestions. In addition, we thank M.~Harris for his sustained interest in this project and his encouragement. Y.~Z.~would like to thank W.-H.~Yang for motivating the completion of the final revisions. We thank the referees for carefully reading the manuscript and providing valuable corrections and suggestions, and the editors for their patience during the exceptionally long revision period.
  
  M.~K.~is partially supported by NSF grant DMS-1902158. S.~W.~S.~is partially supported by~NSF grant DMS-1802039/2101688, NSF RTG grant DMS-1646385, and a Miller Professorship.
  Y.~Z.~is partially supported by NSF grant DMS-1802292, a startup grant at the University of Maryland, and a startup grant at Tsinghua University.

  \mainmatter

  \setcounter{tocdepth}{1}
  \tableofcontents

\chapter*{Introduction}
\addcontentsline{toc}{chapter}{Introduction}
\numberwithin{equation}{section}
\section{The main results}
Shimura varieties have provided a testing ground for many conjectures in the Langlands Program, and have been indispensable in the (partial) solutions of some of these conjectures. Motivated by the work of Eichler, Shimura, Kuga, Sato, and Ihara, Langlands formulated the problem of expressing the Hasse--Weil zeta function of a Shimura variety in terms of automorphic $L$-functions. 
This question is itself a special case of Langlands' conjecture that all motivic $L$-functions are automorphic.

In a series of papers \cite{LanglandsAntwerp,Lan76,langlands1977shimura,langlandsmarchen,langlands1979zeta}, Langlands developed the idea of systematically using trace formulas to attack this problem. In his initial investigations, he encountered the phenomenon of $L$-indistinguishability, which motivated the theory of endoscopy. Based on the latter, Langlands predicted that one should be able to compare a  Lefschetz-type trace formula for the Shimura variety with trace formulas arising in the theory of automorphic representations \emph{after} stabilizing both types of the formulas. This prediction was formulated as a precise conjecture in Kottwitz's paper \cite{Kot90}. 

The main result of the present work is the confirmation of the elliptic part of this conjecture for Shimura varieties of abelian type: We prove an identity between a Grothendieck--Lefschetz--Verdier trace formula on the Shimura variety and elliptic parts of stable Arthur--Selberg trace formulas for the endoscopic groups. 

Let $(G,X)$ be a Shimura datum with reflex field $E$. Fix a prime $\ell$, and let $\xi$ be an algebraic representation of $G$ over $\ol \QQ_{\ell}$. We set  
$$\coh^i_c (\Sh, \xi): = \varinjlim_K \coh^i_c(\Sh_K(G,X)_{\overline E}, \mL_{\xi}),$$
where $K$ runs through all sufficiently small compact open subgroups of $G(\A_f)$, $\Sh_K(G,X)$ denotes the Shimura variety at level $K$, and  $\mL_{\xi}$ is the automorphic $\ol \QQ_{\ell}$-local system associated with $\xi$.\footnote{A 
	technical condition on $\xi$ is required for $\mL_\xi$ to be well defined; for the 
	purposes of this Introduction, the reader may simply take $\xi$ to be trivial, 
	in which case $\mL_\xi = \ol\QQ_\ell$.}  Then $\coh^i_c(\Sh, \xi)$ carries commuting actions by $\Gal(\overline E/E)$ and $G(\A_f)$.  

Next, let $p \neq \ell$ be a prime, and let $\Phi_\fkp \in \Gal(\ol E/E)$ be a geometric 
Frobenius element at a place $\fkp$ of $E$ above $p$. Let $f$ be an element of the 
Hecke algebra $\cH(G(\A_f))$. We assume throughout that $f$ factors as 
$f = 1_{K_p} f^p$, where $f^p \in \cH(G(\A_f^p))$ and $1_{K_p}$ is 
the characteristic function of a hyperspecial subgroup $K_p \subset G(\Q_p)$. (Any fixed $f$ admits such a factorization for almost all primes $p$.) For a positive integer $m$, 
we define
\begin{equation}\label{eq:trace-def}
	T(m,f) := \sum_i (-1)^i \tr\bigl( f \times \Phi_\fkp^m \bigm| \coh^i_c(\Sh, \xi)\bigr).
\end{equation}
As $(\fkp, m, f)$ vary, these quantities completely characterize\footnote{Because 
	local Galois data at ramified primes are omitted, this characterization is understood 
	in the sense that a semi-simple global Galois representation is uniquely determined 
	by the characteristic polynomials of Frobenius elements at almost all  places by the Chebotarev density theorem and the 
	Brauer--Nesbitt theorem.} the virtual representation $\sum_i (-1)^i [\coh^i_c(\Sh, \xi)]$ 
in the Grothendieck group of $\Gal(\ol E/E) \times G(\A_f)$, and in particular 
determine the good Euler factors of the Hasse--Weil zeta function of the Shimura 
varieties $\Sh_K$.

\begin{introthm}[see Theorem \ref{thm:announcement 2}]\label{thm:intro main}
	Assume that $(G,X)$ is of abelian type. For all sufficiently large $m$ we have 
	\begin{align}\label{introeq:STF}
		T(m,f) =  \sum_{\fke} \iota(\fke) ST^{H_1}_{\el, \chi_{H_1}}(f^{H_1}) ,
	\end{align}
	where $\fke$ runs through the elliptic endoscopic data for $G$ up to isomorphism, $\iota(\fke) \in \QQ$ is a constant depending only on $\fke$, and $ST^{H_1}_{\el, \chi_{H_1}}$ is an elliptic stable distribution associated with $\fke$ (defined on a central extension $H_1$ of the endoscopic group in $\fke$). \end{introthm}

If the derived subgroup $G_{\der}$ of $G$ is simply connected and the center $Z_G$ of $G$ is cuspidal (i.e., having equal $\QQ$-rank and $\RR$-rank), then $ST^{H_1}_{\el, \chi_{H_1}}$ is the elliptic part of the stable trace formula for the endoscopic group in  $\fke$. Without these assumptions, the definition involves $z$-extensions and fixed central characters. The functions $f^{H_1}$ will be explained after we state Theorem \ref{thm:intro stabilization} below. The requirement that $m$ is sufficiently large is needed to ensure that the local terms in the Grothendieck--Lefschetz--Verdier trace formula can be calculated naively; this is a special case of Deligne's conjecture,\index{Deligne's conjecture} which has been proved in general by Fujiwara \cite{fujiwara}  and Varshavsky \cite{varshavsky}. For applications this restriction turns out to be harmless. Note that by contrast, knowing (\ref{introeq:STF}) only for all sufficiently \emph{divisible} $m$ would be insufficient for most applications.  

Kottwitz \cite[\S 3, \S 7]{Kot90} conjectured the equality in Theorem \ref{thm:intro main}   for general Shimura varieties (assuming $G_{\der}$ is simply connected and $Z_G$ is cuspidal), and proved it in the case of PEL type A or C in \cite[\S 19]{kottwitz1992points} and \cite[Thm.~7.2]{Kot90}. 
By results of Matsushima \cite{Matsushima} and Franke \cite{Fra98}, the $G(\A_f)$-action on $\coh^i _c(\Sh, \xi )$ can be understood in terms of automorphic representations of $G.$ It is expected that the equality in Theorem \ref{thm:intro main} should lead to a description of $\coh^i_c(\Sh, \xi)$ as a $\Gal(\ol E/E)\times G(\A_f)$-module, at least when the Shimura varieties are compact. This description should involve the global Langlands correspondence between (algebraic) automorphic representations and Galois representations, as well as Arthur's conjectures on automorphic multiplicities. Such a description would further lead to an expression of the Hasse--Weil zeta function in terms of automorphic $L$-functions. See \cite[Part II]{Kot90} for an explanation of this circle of ideas. In the non-compact case, one expects that it is the intersection cohomology of the Baily--Borel compactification, rather than $\coh^i _c(\Sh, \xi )$, that admits a description similar to the compact case. For this, one needs an intersection cohomology variant of Theorem \ref{thm:intro main}, equating traces on the intersection cohomology with \emph{full} stable trace formulas (as opposed to elliptic parts). In the present work we do not consider this variant. However,  Theorem \ref{thm:intro main} and the point counting formula in Theorem \ref{thm:intro point counting} below are expected to play a crucial role in the proof of such a variant; see for instance \cite{Mor10,Mor11,Zhu-orthogonal}.

The proof of Theorem \ref{thm:intro main} consists of two steps. The first step is to prove a ``point counting formula'', expressing $T(m,f)$ in terms of orbital integrals and twisted orbital integrals on $G$ in a way resembling the geometric side of the Arthur--Selberg trace formula. The second step is stabilization, which relates the (twisted) orbital integrals on $G$ with the terms constituting $  ST^{H_1}_{\el, \chi_{H_1}}(f^{H_1})$.  

When $G_{\der}$ is simply connected and $Z_G$ is cuspidal, the point counting formula was already conjectured by Kottwitz 
\cite[\S 3]{Kot90}. Let $q = p^r$ be the cardinality of the residue field of $\fkp$. For $m$ sufficiently large, the conjecture states that 
\begin{multline}\label{eq:intro point counting}
	T(m,f) =  \\  \sum_{\substack{(\gamma_0,\gamma,\delta) \in \KP_{\clsc}(q^{m})/{\sim},\\ \alpha(\gamma_0,\gamma,\delta) =0}} c_1(\gamma_0,\gamma,\delta) c_2(\gamma_0) O_{\gamma}(f^p) TO_{\delta}(\phi_{mr}) \tr\xi(\gamma_0).
\end{multline} Here $\KP_{\clsc}(q^m)$ consists of  triples $(\gamma_0,\gamma, \delta)\in G (\Q) \times G(\A_f^p) \times G(\Q_{q^m})$  such that $\gamma_0$ is $\RR$-elliptic, stably conjugate to $\gamma$, and stably conjugate to the degree $mr$ norm of $\delta$. There is also a technical assumption on $\delta$ which we omit here. (The notation $\KP_{\clsc}$ stands for ``classical Kottwitz parameters''.) The equivalence relation $\sim$ on $\KP_{\clsc}(q^m)$ is given by stable conjugacy on the first factor, conjugacy on the second factor, and $\sigma$-conjugacy on the third factor.  Kottwitz defines a Galois cohomological invariant $\alpha(\gamma_0,\gamma,\delta)$ for each $(\gamma_0,\gamma,\delta) \in \KP_{\clsc}(q^m)$, and in (\ref{eq:intro point counting}) the summation is subject to the condition  $\alpha (\gamma_0,\gamma,\delta)=0$. In each summand, we have an orbital integral $O_{\gamma}(f^p)$ on $G(\A_f^p)$, a twisted orbital integral $TO_{\delta}(\phi_{mr})$ on $G(\Q_{q^m})$ (where $\phi_{mr}$ is an explicit function on $G(\Q_{q^m})$), the character $\tr \xi $ of $\xi$ evaluated at $\gamma_0$, a volume term $c_1(\gamma_0,\gamma,\delta)$, and a term $c_2(\gamma_0)$ defined via Galois cohomology.  

In the conjectural formula (\ref{eq:intro point counting}), the assumption that $G_{\der}$ is simply connected is quite serious. Without it, Kottwitz's construction of the invariant $\alpha(\gamma_0,\gamma,\delta)$ for $(\gamma_0,\gamma,\delta) \in \KP_{\clsc}(q^m)$ no longer works, and  the volume term $c_1(\gamma_0,\gamma,\delta)$ is not well defined. These problems are caused by the possible disconnectedness of $G_{\gamma_0}$ for a semi-simple $\gamma_0 \in G(\QQ)$. In the following theorem, the point counting formula we prove is a generalization of (\ref{eq:intro point counting}) without any assumptions on $G_{\der}$ and $Z_G$.
\begin{introthm}[see Theorem \ref{thm:announcement 1}]\label{thm:intro point counting}
	If $(G,X)$ is of abelian type, then for all sufficiently large $m$ we have \begin{multline}
		\label{eq:intro refined point counting}	 T(m,f) \\ = 	\sum_{\gamma_0 \in \Sigma} \sum_{\substack{\fkc \in \KP(\gamma_0, q^m) \\ \alpha(\fkc) = 0 }} \abs{(G_{\gamma_0}/G_{\gamma_0}^0)(\Q)}^{-1} c_1(\fkc) c_2(\gamma_0) O_{\fkc}(f^p) TO_{\fkc}(\phi_{mr}) \tr\xi (\gamma_0), 
	\end{multline}
	where the terms $c_1(\fkc), c_2(\gamma_0), O_{\fkc}(f^p), TO_{\fkc}(\phi_{mr})$ are defined analogously as the terms in (\ref{eq:intro point counting}). 
\end{introthm}

The most significant new feature of the form of  (\ref{eq:intro refined point counting}) is that the summation index set $\KP_{\clsc} (q^m)/{\sim}$ in (\ref{eq:intro point counting}) has been replaced by a more refined set $$\coprod_{\gamma_0 \in \Sigma} \KP(\gamma_0, q^m)$$ which admits a map to the former. Here $\Sigma$ is a certain subset of the set of $\RR$-elliptic elements of $G(\QQ)$, and for each $\gamma_0 \in \Sigma$ the definition of $\KP(\gamma_0, q^m)$ is Galois cohomological in nature. (We also allow $Z_G$ to be non-cuspidal, in which case $\Sigma$ depends on the choice of a compact open subgroup $K^p \subset G(\A_f^p)$ such that $f^p$ is $K^p$-bi-invariant.) For each $\fkc \in \KP(\gamma_0, q^m)$, we define an invariant $\alpha(\fkc)$ lying in an abelian group that depends only on $G_{\gamma_0}^0$ and $G$. This definition specializes to Kottwitz's invariant $\alpha(\gamma_0,\gamma,\delta)$ when $G_{\der}$ is simply connected. In (\ref{eq:intro refined point counting}) the condition $\alpha (\fkc)=0$ is imposed, similarly 
to (\ref{eq:intro point counting}).

Once Theorem \ref{thm:intro point counting} is proved, in order to prove Theorem \ref{thm:intro main} we need to stabilize the right hand side of (\ref{eq:intro refined point counting}). We prove this stabilization in general as in the next theorem, without assuming that $(G,X)$ is of abelian type. 
\begin{introthm}[see Theorem \ref{thm:end-of-stabilization}]\label{thm:intro stabilization}
	The right hand side of (\ref{eq:intro refined point counting}) is equal to $$ \sum_{\fke} \iota(\fke) ST^{H_1}_{\el, \chi_{H_1}}(f^{H_1}). $$
\end{introthm}

Theorem \ref{thm:intro main} follows immediately from Theorems \ref{thm:intro point counting} and \ref{thm:intro stabilization}. 
The proof of Theorem \ref{thm:intro stabilization} follows the general outline of \cite[\S7]{Kot90}. Namely, we  turn the right hand side of (\ref{eq:intro refined point counting})  into a sum of $\kappa$-orbital integrals (twisted at $p$) over adelic conjugacy classes, and then  rewrite the sum in terms of stable distributions on endoscopic groups via the Kottwitz--Langlands--Shelstad transfer and the fundamental lemma. More precisely, $f^{H_1}$ away from $\{p,\infty\}$ is obtained from $f^p$ via the usual untwisted transfer, whereas $f^{H_1}$ at $p$ is a twisted transfer of $\phi_{mr}$, and $f^{H_1}$ at $\infty$ is constructed explicitly as a finite linear combination of certain stably cuspidal functions.

Although we largely follow the strategy outlined in \textit{loc.~cit.}, our approach incorporates two technical improvements to the stabilization process. First, we dispense with the simplifying hypotheses that $G_{\der}$ is simply connected and $Z_G$ is cuspidal, achieved by systematically employing $z$-extensions and fixed central characters. Second, in contrast to \cite[Thm.~7.2]{Kot90}, Theorem \ref{thm:intro main} imposes no $(G,H)$-regularity conditions on the stable distributions. The reason is that the non-$(G,H)$-regular semisimple terms vanish, as shown by Morel (see \S\ref{h-infty} and Lemma \ref{lem:adelic-orb-int-identity} below). 

We refer the reader to \S\ref{subsec:intro stabilization} for a more detailed summary of the proof of Theorem \ref{thm:intro stabilization}.

\section{Applications}
\label{sec:intro applications}

The following results have been obtained using 
Theorem \ref{thm:intro main} or its proof (Theorem \ref{introthm:LRtau} below). 

\begin{itemize}
	\item One of us (Shin) has, in joint work with Kret, constructed the automorphic to Galois direction of the Langlands correspondence for 
	$\mathrm{GSp}_{2n}$ and (a form of) $\textup{GSO}_{2n}$ over totally real fields, under a technical local hypothesis 
	\cite{KS-GSp,KS-GSO}. This involves constructing Galois representations into $\mathrm{GSpin}$ groups, 
	cf. \cite[p.~268]{faltingschai}. 
	\item One of us (Zhu) has given a description of the Hasse--Weil zeta function and the Hecke--Galois action on the intersection cohomology of the Baily--Borel compactification of Shimura varieties for some global forms of $\textup{SO}(N,2)$ in terms of automorphic representations \cite{Zhu-orthogonal}. This completes the Langlands--Kottwitz program for these varieties at almost all primes.  
	\item Youcis \cite{Youcis} has extended Scholze's version \cite{Scholze:LK} of the Langlands--Kottwitz method for Shimura varieties 
	with bad reduction from the case of PEL type to the case of abelian type.
	\item Mack-Crane and one of us \cite{Shin-MC} have obtained a trace formula for Igusa varieties of Hodge type which is 
	analogous to Theorem \ref{thm:intro point counting} in the case of Hodge type. This generalizes \cite{Shi09}. 
	(A generalization to the case of abelian type, as well as a stabilization analogous to Theorem \ref{thm:intro main} is expected, cf.~\cite{Shin10}.) 
\end{itemize}

We stress that the Shimura data appearing in concrete applications, as in the first two items, are typically of abelian type but \emph{not} of Hodge type. The same is true of the three applications below. As we will explain in \S \ref{subsec:variants} below, the proof of Theorem \ref{thm:intro main} in the case of Hodge type is substantially easier, but this does not suffice for many applications. 

In general, 
Theorem~\ref{thm:intro main} is the key to determining the fundamental virtual $G(\A_f)\times \Gal(\ol E/E)$-module $[\coh_c(\Sh,\xi)]:=\sum_{i} (-1)^i \coh^i_c(\Sh, \xi)$ in terms of automorphic representations. A notable consequence would be an expression of the Hasse--Weil zeta function of $\Sh_{K}$ as an alternating product of automorphic $L$-functions for sufficiently small compact open subgroups $K\subset G(\A_f)$, possibly up to Euler factors at finitely many primes.
We intend to work out the details in a forthcoming follow-up paper, and obtain unconditional results in various special cases, which will provide important ingredients for some remarkable arithmetic results:
\begin{itemize}
	\item The Beilinson--Bloch--Kato conjecture for Rankin--Selberg motives using  
	Shimura varieties of unitary groups \cite{LTXZZ}.
	\item Higher-dimensional Gross--Zagier formula (i.e., arithmetic inner product formula) using  
	Shimura varieties of unitary groups \cite{LiLiu}.
	\item Euler systems arising from Shimura varieties of $\textup{SO}(2n-1,2)$ \cite{Cornut}.
\end{itemize} Related to the applications in \cite{LTXZZ, LiLiu}, 
it is worth pointing out that when one passes between a Shimura datum of abelian type and an ``isogenous'' Shimura datum of Hodge type, the reflex field is often not preserved. Thus even if one is just interested in constructing representations of $\Gal(\ol E/E)$ using the cohomology of a Shimura variety of a unitary group with reflex field $E$ (which is of abelian type but not of Hodge type), one cannot pass to a Shimura variety of Hodge type without having to enlarge $E$ in general.

To understand the structure of $[\coh_c(\Sh,\xi)]$ in the general case of abelian type, there are two main obstacles to proving an unconditional theorem. Let us briefly address them, leaving the details to \S\ref{sub:speculation} below.

When $G$ is anisotropic modulo center over $\Q$, or equivalently when the finite-level Shimura varieties $\Sh_K(G,X)$ are projective, 
the first problem is to show that 
$$ST^{H_1}_{\el,\chi_{H_1}}(f^{H_1})=ST^{H_1}_{\chi_{H_1}}(f^{H_1})$$ 
in the summand of Theorem \ref{thm:intro stabilization}, namely that the non-elliptic terms in the stable distribution cancel each other out. The resolution of this problem is within reach at least in various special cases that are of interest for applications, and this will be treated in our forthcoming follow-up paper. The second problem is that the endoscopic classification for automorphic representations is not available for general reductive groups.

If $G$ is isotropic modulo center, the second problem remains the same. In place of the first problem, however, it is desirable to promote Theorem~\ref{thm:intro main} by proving an equality where the compactly supported cohomology and $ST^{H_1}_{\el,\chi_{H_1}}$ are replaced with the intersection cohomology of the Baily--Borel compactification and $ST^{H_1}_{\chi_{H_1}}$, respectively. As mentioned above, such an upgrade is obtained for $\mathrm{SO}(N,2)$ in \cite{Zhu-orthogonal}. Some results were previously known for Shimura varieties of PEL type A and C \cite{Montreal92,Mor08,Mor10,Mor11}.

Another application of our work would be the analogues of Theorems \ref{thm:intro main} and \ref{thm:intro point counting} for Shimura varieties of parahoric level at $p$, in light of recent advances on the Haines--Kottwitz test function conjecture \cite{HainesRicharz1,HainesRicharz2} and the Langlands--Rapoport conjecture in the parahoric case \cite{ZhouParahoric,vanHoftenLR}. The latter takes as an input the hyperspecial case through the earlier work \cite{kisin2012modp}; a strengthening should be possible by appealing to our improvement (Theorem \ref{introthm:LRtau}) instead.

\section{Variants of the Langlands--Rapoport Conjecture}\label{subsec:variants}
We now discuss the proof of Theorem \ref{thm:intro point counting}. 
To simplify the exposition we assume that $\xi$ is trivial (so  $\mL_{\xi} = \overline \QQ_{\ell}$), and that $f^p = 1_{K^p}$ for a sufficiently small compact open subgroup $K^p \subset G(\A_f^p)$, while continuing to assume that $K_p$ is hyperspecial. 
If $\Sh_{K_pK^p}$ is proper over $E,$ then one expects that there exists a proper smooth integral model $\Shh_{K_pK^p}$ 
of $\Sh_{K_pK^p}$ over $\oo_{E,(\fkp)}$. In this case we have (with the Haar measure normalized such that $K_p K^p$ has volume $1$) 
\begin{align}\label{eq:intro Lefschetz}
	T(m,f)  = \# \Shh_{K_pK^p}(\FF_{q^m}).
\end{align} Even outside the proper case, one still conjectures that there exists a canonical smooth integral model $\Shh_{K_pK^p}$ satisfying (\ref{eq:intro Lefschetz}) (among other conditions). Hence in all cases we seek a formula for $\# \Shh_{K_pK^p}(\FF_{q^m})$, thus the name ``point counting formula''. 

For Shimura varieties of Hodge type, it is possible to establish a point counting formula by generalizing the considerations of Kottwitz \cite{kottwitz1992points} in the PEL case, with the aid of the results from \cite{kisin2012modp}. In this approach one attaches group-theoretic invariants to isogeny classes over a fixed finite field $\FF_{q^m}$, and then counts the number of $\FF_{q^m}$-points having the same invariants; see \cite{Lee}. However, it 
does not seem possible to derive Theorem \ref{thm:intro point counting} for general Shimura varieties of abelian type from such an approach.

In the present work, we adopt the perspective of Langlands--Rapoport \cite{langlands1987gerben}, aiming to relate $\ol \FF_q$-isogeny classes to certain group-theoretic data defined via Galois gerbes. While establishing these results in the Hodge-type setting requires greater technical effort, this approach has the advantage that one can then bootstrap these statements to the case of abelian type, from which the point counting formula for abelian type would ultimately follow.

Write $\Shh_{K_p}$ for $\varprojlim_{K^p} \Shh_{K_p K^p}$, and write $\Phi$ for the geometric $q$-Frobenius. The Langlands--Rapoport Conjecture states that there is a $G(\A_f^p)\times \Phi^{\ZZ}$-equivariant bijection
$$ \Shh_{K_p } (\overline \FF_q) \isom \coprod_{\phi} \varprojlim_{K^p} I_{\phi} (\QQ) \backslash X(\phi)/K^p. $$ Here $\phi$ runs through conjugacy classes of \emph{admissible morphisms} $\Qf \to \G_G$, which are objects of purely group-theoretic nature. Associated with each $\phi$ there is a reductive group $I_\phi$ over $\Q$ and a set $X(\phi)$ equipped with commuting actions by $I_\phi(\A_f)$, $\Phi^{\ZZ}$, and $G(\A_f^p)$.

Currently, the Langlands--Rapoport Conjecture is open even for the Siegel modular varieties. (For some quaternionic Shimura varieties the conjecture has been proved by Reimann \cite{reimann1997zeta}.) In \cite{kisin2012modp}, a weaker version of the conjecture is proved for the canonical integral models of Shimura varieties of abelian type for $p>2$; these integral models were previously constructed in \cite{kisin2010integral} for $p>2 $ and in \cite{KMP16} for $p=2$, and shown to satisfy  (\ref{eq:intro Lefschetz}) in  \cite{lanstrohII}.   In this weaker version, the set 
$$\varprojlim_{K^p} I_{\phi}(\QQ) \backslash X(\phi)/K^p $$ is replaced by 
$$\varprojlim_{K^p} I_{\phi}(\QQ)_{\tau(\phi)} \backslash X(\phi)/K^p, $$ where $\tau(\phi)$ is an unspecified element of $I_{\phi}^{\ad}(\A_f)$, and $I_{\phi}(\QQ)_{\tau(\phi)} $ is the image of $I_{\phi}(\QQ)$ under $$ I_{\phi}(\QQ) \hookrightarrow I_{\phi}(\A_f) \xrightarrow{\Int(\tau(\phi))} I_\phi (\A_f). $$ In order to deduce (\ref{eq:intro refined point counting}) from such a weaker statement, one must have some control of the elements $\tau(\phi)$. 

We introduce some definitions in order to state the desired control of $\tau(\phi)$. For each admissible morphism $\phi$, we have the \emph{algebraic part} $\phi^{\Delta}$ of $\phi$, which is a $\Qbar$-homomorphism from a pro-torus $\Qf^{\Delta}$ to $G_{\Qbar}$. Different
admissible morphisms whose algebraic parts are conjugate under $G(\Qbar)$ can be obtained from each other by a Galois cohomological twisting. The double quotient set $$\cH(\phi):=\tauhome{\phi}$$ is an abelian group, and up to canonical isomorphism it depends only on the $G(\Qbar)$-conjugacy class of $\phi^{\Delta}$. For each maximal torus $T$ in $I_{\phi}$, write $\cH(\phi,T)$ for the cokernel\footnote{The notation $\cH(\phi,T)$ is used only in the Introduction.} of the localization map $$\ker\big(\coh^1(\Q,T) \to \coh^1(\RR,T) \oplus \coh^1_{\ab}(\Q,G)\big) \To  \coh^1(\A_f,T). $$ There is a natural homomorphism $ \cH(\phi) \to \cH(\phi,T)$ as explained in  Definition \ref{defn:rectification}. 
\begin{introconj}[``Langlands--Rapoport--$\tau$''; see Conjecture \ref{hypo about LR}]\label{introconj:LRtau} There is a  bijection
	\begin{align}\label{eq:intro LRbij}
		\Shh_{K_p} (\overline \FF_q) \isom \coprod_{\phi} \varprojlim_{K^p} I_{\phi} (\QQ)_{\tau(\phi)} \backslash X(\phi)/K^p,
	\end{align}  which is $G(\A_f^p)\times \Gal(\overline \FF_{q}/\FF_q)$-equivariant, with   $\tau(\phi) \in I_{\phi}^{\ad}(\A_f)$ satisfying the following conditions. 
	\begin{enumerate}
		\item The image of $\tau(\phi)$ in $\cH(\phi)$ depends only on the $G(\Qbar)$-conjugacy class of $\phi^{\Delta}$. 
		\item For each maximal torus $T$ in $I_{\phi}$, the image of $\tau(\phi)$ in $\cH(\phi)$ lies in the kernel of $\cH(\phi) \to \cH(\phi,T)$. 
	\end{enumerate}
\end{introconj} While Conjecture \ref{introconj:LRtau} allows for the elements $\tau(\phi)$ and is therefore a weakening of the original Langlands--Rapoport Conjecture, it nonetheless imposes precise conditions controlling them. No such control is proved in \cite{kisin2012modp}.

\begin{introthm}[see Theorem \ref{thm:main thm in point counting}] \label{introthm:one conj implies the other} 
	Conjecture \ref{introconj:LRtau} implies the point counting formula (\ref{eq:intro refined point counting}). 
\end{introthm}

The proof of Theorem \ref{introthm:one conj implies the other} is group-theoretic 
in nature, with key ingredients 
drawn from \cite{kottwitztwisted} and \cite[\S 5]{langlands1987gerben}. In \textit{loc.~cit.}, Langlands and Rapoport laid the groundwork for the transition from their conjecture 
to a formula of the form \eqref{eq:intro point counting} under the hypothesis that $G_{\der}$ 
is simply connected; see also \cite{milne92}. To establish Theorem~\ref{introthm:one conj implies the other}, 
we must incorporate the possibly non-trivial elements $\tau(\phi)$ and extend the 
method to arbitrary Shimura data without hypotheses on $G_{\der}$ or $Z_G$, 
which eventually leads to the group-theoretically more sophisticated formula 
\eqref{eq:intro refined point counting}. \footnote{A further technical complication 
	is that the original construction of the quasi-motivic gerb in \cite{langlands1987gerben} 
	contained an error, which was corrected by Reimann \cite{reimann1997zeta}. Consequently, 
	several key definitions and arguments from \cite[\S 5]{langlands1987gerben} require 
	substantial modification in our setting.}

The logical relations between the various conjectures are depicted in the following diagram.  
$$
\xymatrixrowsep{3pc} \xymatrixcolsep{1.4pc}
\xymatrix{
	\boxed{\text{LR Conjecture}} \ar@{=>}[r]
	& \boxed{\begin{array}{c} \text{LR--$\tau$ Conjecture} \\ \text{(Conjecture~\ref{introconj:LRtau})} \end{array}}
	\ar@{=>}[r]
	\ar@{=>}[d]^{\text{Theorem \ref{introthm:one conj implies the other}}}
	& \boxed{\begin{array}{c} \text{weak LR Conjecture} \\ \text{as in \cite{kisin2012modp}} \end{array}} \\
	& \boxed{\begin{array}{c} \text{point counting formula} \\ \text{(\ref{eq:intro refined point counting})} \end{array}}
	\ar@{=>}[d]^{\text{Theorem \ref{thm:intro stabilization}}} \\
	& \boxed{\begin{array}{c} \text{stable trace formula} \\ \text{(\ref{introeq:STF})} \end{array}}
}
$$ 
\newpage
\begin{introthm}[see Theorem \ref{thm:LRtau for abelian type}]\label{introthm:LRtau} If $(G,X)$ is of abelian type, then Conjecture \ref{introconj:LRtau} holds with respect to the canonical integral models. 
\end{introthm}
Theorem \ref{thm:intro point counting} follows from Theorems \ref{introthm:one conj implies the other} and \ref{introthm:LRtau}. We now  discuss the proof of Theorem \ref{introthm:LRtau}.

\section{Proof of the Langlands--Rapoport--\texorpdfstring{$\tau$}{tau} Conjecture} \label{subsec:intro Hodge type}

In Conjecture \ref{introconj:LRtau}, we have imposed the minimal conditions on the elements $\tau(\phi)$ required to deduce the point-counting formula (\ref{eq:intro refined point counting}), regardless of the type of the Shimura datum (see Theorem \ref{introthm:one conj implies the other}). However, one may strengthen the conjecture by requiring compatibility with 
additional structures on the Shimura variety---namely, the set of connected 
components and a certain symmetry arising from the Galois cohomology of the center. 
For ease of reference in this Introduction, we shall refer to the latter as 
\emph{$\mathscr{A}$-symmetry}.\footnote{The precise formulation of this symmetry is somewhat involved; see \S \ref{subsec:LR str}.} Keeping track of these extra structures on  Shimura varieties of Hodge type is precisely what allows one to reconstruct Shimura varieties of abelian type (both for the generic fiber and the integral model). Our strategy therefore proceeds in two steps. First, we establish a strengthened version of Conjecture \ref{introconj:LRtau} incorporating connected components and $\mathscr A$-symmetry for Shimura varieties of Hodge type. Second, we bootstrap this result to derive the original Conjecture \ref{introconj:LRtau} for general Shimura varieties of abelian type.

\bigskip
\textbf{(I) Proving a strengthening of Conjecture \ref{introconj:LRtau} for Hodge type.}
\bigskip 

By the theory of integral models \cite{kisin2010integral}, after fixing a  Siegel embedding of the Shimura datum of Hodge type, for each $x\in \Shh_{K_p}(\overline \FF_q)$ we obtain an abelian variety $\cA_x$ (up to prime-to-$p$ isogeny) over the residue field of $x$, together with tensors on $\coh^1_{\et}(\cA_x, \A_f^p)$ and on $\coh^1_{\cris}(\cA_x/\Zpur) : = \coh^1_{\cris}(\cA_x/W(k))\otimes_{W(k)} \Zpur$ (where $k$ is the finite residue field of the image of $x$ to some finite level), whose stabilizer group schemes are canonically identified, up to conjugation by $\A_f^p$-points and $\Zpur$-points, with $G_{\A_f^p}$ and $\cG_{\Zpur}$ respectively. Here  $\cG$ denotes the reductive group scheme over $\Z_p$ corresponding to the hyperspecial subgroup $K_p$. The set $\Shh_{K_p}(\overline \FF_q)$ is partitioned into  \emph{isogeny classes}, where two points $x$ and $x'$ are defined as  isogenous if there exists a quasi-isogeny $\cA_x \to \cA_{x'} $ preserving the tensors.  

For $x\in\Shh_{K_p}(\overline \FF_q)$, the relative Frobenius  on $\coh^1_{\et}(\cA_x, \A_f^p)$ and the absolute Frobenius on  $\coh^1_{\cris}(\cA_x/\Zpur)$ give rise to 
an element of\footnote{The notation $G^\sim$ is used only in the Introduction.}
$$  G^\sim : = \bigg(\varinjlim_{n\in \mathbb Z_{\geq 1}} G(\A_f^p) \bigg) \times \set{\cG(\Z_p^{\ur}) \text{-}\sigma\text{-conjugacy classes in } G(\Qpur)}, $$
where the direct limit is with respect to the transition maps $\gamma \mapsto \gamma^{n'/n}$
for $n| n'$.   Write $\pi(G,X)$ for the set of connected components of the scheme $ \Shh_{K_p, \O_{\overline {E_\fkp}}}$, which has a group-theoretic description. We then have a natural map $$\Shh_{K_p}(\overline \FF_q) \To G^\sim \times \pi(G,X). $$ 
Analogously, for each admissible morphism $\phi$, we have a natural map $$X(\phi) \To G^{\sim} \times \pi(G,X).$$

We denote by $\ica$ the set of isogeny classes in $\Shh_{K_p}(\overline \FF_q)$, and by $\cca$ the set of admissible morphisms $\Qf \to \G_G$ up to conjugacy. A notion of central importance in our argument is that of \emph{amicable pairs}:\index{amicable pair} these are pairs $(\II, \phi) \in \ica \times \cca$ that are eligible to be matched under the bijection \eqref{eq:intro LRbij}. The defining feature of an amicable pair is the existence of a \emph{marking}. In the simplest setting, a marking consists of a pair of points $(y,y') \in \II \times X(\phi)$ sharing the same image in $G^\sim \times \pi(G,X)$. In general, however, the marking points are allowed to lie in certain canonical enlargements of $\II$ and $X(\phi)$, with the requirement that their images in a canonical enlargement of $G^\sim \times \pi(G,X)$ agree. We shall return to these enlargements shortly. The following result is a prerequisite for constructing \eqref{eq:intro LRbij}.

\begin{introthm}[see Theorem \ref{thm:matchability of bij}]\label{introprop:Hodge}
	There exists a bijection $\bij: \cca \isom \ica$ whose graph is contained in the set of amicable pairs. Moreover  
	$\bij$ is equivariant with respect to Galois cohomological twistings on both sides.
\end{introthm}
We briefly explain the last requirement on $\bij$. The Galois cohomological twisting on $\ica$ permutes transitively the $\QQ$-isogeny classes in a $\Qbar$-isogeny class, while the twisting on $\cca$ permutes  transitively the conjugacy classes of admissible morphisms having conjugate algebraic part. For an arbitrary amicable pair $(\II,\phi)$, the set of Galois cohomology classes that can be used to twist $\phi$ is canonically identified with the analogous set for $\II$. The requirement is that for every such class $\beta$, we have $\II^{\beta} = \bij (\phi^{\beta})$ provided that $\II = \bij (\phi)$. 

We now explain the proof of Theorem \ref{introprop:Hodge}. We make use of the following diagram:
\begin{align*}
	\xymatrixcolsep{5pc} \xymatrix{& \{ \text{special point data}\}  \ar@/_1pc/[ldd]  \ar@/^1pc/[rdd]  \ar@{-->}[d] \\ &  \{ \text{amicable pairs} \}   \ar[ld] \ar[rd] \\ \cca && \ica }
\end{align*}
Here a\index{special point datum} \emph{special point datum} refers to a triple $(T,i,h)$, where $(T,h)$ is a Shimura datum on a torus, and $i$ is an embedding of Shimura data $(T,h) \to (G,X)$ such that $i(T)$ is a maximal torus in $G$. Given $(T,i,h)$, we obtain a special point $\tilde x_{(T,i,h)}$ in the generic fiber of $\Shh_{K_p}$, which specializes to a point  $x _{(T,i,h)}\in \Shh_{K_p}(\overline \FF_q)$. Define $ \II(T,i,h)$ to be the isogeny class of $x_{(T,i,h)}$. Analogously, $(T,i,h)$ also gives rise to an admissible morphism $\phi(T,i,h)$. We thus have maps from the set of special point data to $\ica$ and $\cca$. The following proposition asserts that we can fill in the dashed arrow in the above diagram. 

\begin{introprop}[see Corollary \ref{cor:standard amicable}]\label{introprop:amicable}
	Let $(T,i,h)$ be a special point datum. Then the pair $(\II(T,i,h), \phi(T,i,h))$ is amicable. 
\end{introprop}

The key arithmetic ingredient in the proof of Proposition~\ref{introprop:amicable} (and in several subsequent arguments) is the construction of certain \emph{integral special points},\index{integral special point}  
together with the computation of their images in a canonical enlargement $G^\sim \times \pi^*(G,X)$ 
of $G^\sim \times \pi(G,X)$. The main complication arises from an incompatibility 
between integral models: if $\cT^{\circ}$ denotes the connected N\'eron model of 
$T_{\Q_p}$ over $\Z_p$, the morphism $i \colon T_{\Q_p} \to G_{\Q_p}$ does not in general 
extend to a morphism $\cT^\circ \to \cG$. Consequently, the image of the ``standard'' 
special point $x_{(T,i,h)} \in \II(T,i,h)$ in $G^\sim \times \pi(G,X)$ cannot easily 
be computed directly.\footnote{If the second factor of $G^\sim$ (the component at $p$) 
	is considered only up to $G(\Qpur)$-$\sigma$-conjugacy, this image is well understood. 
	For our purposes, however, we require the finer precision modulo $\cG(\Zpur)$-$\sigma$-conjugacy.} The integral special points are designed to respect the integral structure provided 
by $\cT^\circ$ in the same way that $x_{(T,i,h)}$ respects that  provided 
by $\cG$. Due to the 
same incompatibility, such points cannot in general be constructed inside $\II(T,i,h)$ 
itself; instead, they are constructed within a canonical enlargement 
$\II(T,i,h)^* \supset \II(T,i,h)$ introduced in a purely formal way. Crucially, the 
map $\II(T,i,h) \to G^\sim \times \pi(G,X)$ extends to a map
\begin{equation}\label{eq:extended map}
	\II(T,i,h)^* \longrightarrow G^\sim \times \pi^*(G,X),
\end{equation}
under which the images of these integral special points can be explicitly computed. This extended map and its counterpart for admissible morphisms are precisely the maps used in the general definition of an amicable pair alluded to above.

To carry this out, we first state a general construction from $p$-adic Hodge theory. Let $A$ be an abelian variety over a finite extension $F$ of $\QQ_p$ with good reduction and special fiber $\bar A$. Write $$\Lambda_A := \coh^1_{\et}(A_{\ol F} , \Z_p), \quad M_A:= \coh^1_{\cris}(\bar A/\Zpur) = \coh^1_{\cris}(\bar A/W(k_F)) \otimes_{W(k_F)} \Zpur. $$ Let $P$ be a reductive $\Q_p$-subgroup of $\GL(\Lambda_A[1/p])$ such that $P(\Q_p)$ contains the image of $\Gal(\overline F/F)$, and let $\cP(\Z_p)$ be a parahoric subgroup of $P(\Q_p)$ stabilizing the lattice $\Lambda_A \subset \Lambda_A[1/p]$. Under the hypotheses that $\cP \to \GL(\Lambda_A)$ is a closed immersion of $\Z_p$-group schemes and that $\cP$ is a connected reductive group scheme, a fundamental  construction of  \cite{kisin2010integral,kisin2012modp} provides a $\Zpur$-linear isomorphism 
\begin{align}\label{eq:eta}
	\eta_{A, \cP} :  \Lambda_A \otimes_{\Z_p} \Zpur \isom M_A,
\end{align} 
canonical up to the $\cP(\Zpur)$-action on the source; this is the foundation for the construction of the crystalline tensors on $\cA_x$ for $x\in \Shh_{K_p}(\overline \FF_q)$. In the present work, we 
generalize this construction by removing the two hypotheses, which will be outlined in  \S\ref{subsec:intro p-adic Hodge} below.

We now apply this machinery to sketch the construction of the integral special point. Associated with the special point $\tilde x_{(T,i,h)}$ in the generic fiber is a CM abelian variety $A$ defined over a number field, equipped with a canonical $\cG$-action on $\Lambda :  = 
\Lambda_A$. 
In fact, $\Lambda$ can be canonically identified with the integral symplectic representation of $\cG$ fixed in the ambient Hodge-type setting. Via the morphism $i :  T_{\Q_p} \to G_{\Q_p}$, the torus $T_{\Q_p}$ acts on $\Lambda[1/p]$. Fix a $\cT^\circ$-stable $\Z_p$-lattice $\Lambda'$ in $\Lambda[1/p]$, which then corresponds to a $\Q_p$-isogeny $\iota : A
' \to A$ between CM abelian varieties. Because of the aforementioned incompatibility, we cannot in general choose $\Lambda' = \Lambda$, and hence $\iota$ need not be an isomorphism. Choose a finite extension $F/\Q_p$ over which $A$, $A'$, and $\iota$ are defined, and such that $A$ and $A'$ admit good reduction. Reduction of $\iota$ then induces an isomorphism $\iota^* \colon M_{ A}[1/p] \isom M_{ A'}[1/p]$.

Because $A$ arises from the special point $\tilde x_{(T,i,h)}$, the $\Gal(\overline F/F)$-action on $\Lambda[1/p]$ factors through $T(\Q_p)$. Since $\cG$ and $\cT^{\circ}$ are both parahoric group schemes, specializing~\eqref{eq:eta} yields isomorphisms\footnote{Note that in general, $\cT^{\circ}$ is not reductive due to ramification of $T_{\Q_p}$, and the morphism $\cT^\circ \to \GL(\Lambda')$ is not a closed immersion (and indeed, there may not exist a $\cT^\circ$-stable lattice $\Lambda'$ in the $T_{\Q_p}$-representation $\Lambda[1/p]$ having this property). Therefore we genuinely require the full generality of~\eqref{eq:eta}.}
\begin{align*}
	\eta= 	\eta_{A, \cG} & \colon \Lambda \otimes_{\Z_p} \Zpur \isom M_A, \\ 
	\eta'= 	\eta_{A',\cT^\circ} & \colon \Lambda' \otimes_{\Z_p} \Zpur \isom M_{A'}.
\end{align*}  Now we define $g\in \GL(\Lambda\otimes_{\Z_p} \Qpur)$ to be the composite isomorphism 
\[
\xymatrixcolsep{3pc}
\xymatrix{
	\Lambda \otimes \Qpur \cong \Lambda' \otimes \Qpur \ar[r]^-{\eta'[1/p]}_-{\sim} & 
	M_{A'}[1/p] \ar[r]^-{(\iota^*)^{-1}}_-{\sim} & 
	M_A[1/p] \ar[r]^-{\eta^{-1}[1/p]}_-{\sim} & 
	\Lambda \otimes \Qpur.
}
\]
It turns out that $g$ lies in the subgroup $G(\Qpur)$. Since $\eta$ and $\eta'$ are canonical up to the actions of $\cG(\Zpur)$ and $\cT^\circ(\Zpur)$ on the source respectively, the double coset $\cG(\Zpur) g\cT^{\circ}(\Zpur)$ is independent of the choices of $\eta$ and $\eta'$. In fact it is also independent of the choice of $\Lambda'$. 

If the coset $g \cG(\Zpur) \in G(\Qpur)/\cG(\Zpur)$ lies in the relevant affine 
Deligne--Lusztig variety, it defines a point in the isogeny class $\II(T,i,h)$ 
by virtue of the uniformization of the latter by this variety (a key construction 
in \cite{kisin2012modp}, ignoring structures away from $p$). In general, however, 
this condition fails. To circumvent this, the formal enlargement 
$\II(T,i,h)^* \supset \II(T,i,h)$ is defined so as to be uniformized by the full 
quotient $G(\Qpur)/\cG(\Zpur)$, compatibly extending the uniformization of the subset 
$\II(T,i,h)$ by the affine Deligne--Lusztig variety. In this way, $g \cG(\Zpur)$ 
specifies the desired integral special point $x_g \in \II(T,i,h)^*$.

The crux of the proof of Proposition~\ref{introprop:amicable} then reduces to computing the image of $x_g$ under the extended map (\ref{eq:extended map}). This computation rests on the following result, whose proof is deferred to \S\ref{subsec:intro p-adic Hodge}.

\begin{introprop}\label{introlem:1}  The following statements hold.
	\begin{enumerate}
		\item  The crystalline Frobenius on $M_{A'}[1/p]$ corresponds via $\eta'$ to $\delta \sigma,$ with $\delta \in T(\Qpur) \subset \GL(\Lambda' \otimes_{\ZZ_p} \Qpur)$. Moreover,  the image of $\delta$ in $X_*(T)_{\Gamma_{p,0}}$ (where $\Gamma_{p,0}$ denotes the inertia subgroup of $\Gal(\Qpbar/\Q_p)$) under the Kottwitz homomorphism is equal to the negative of the natural image of $\mu_h$, the Hodge cocharacter of $h: \mathbb S \to T_{\RR}$. 
		\item  The coset $g \cG(\Zpur)$ contains an element in the image of $G_{\sconn} (\Qpur)$. 
	\end{enumerate} 
\end{introprop} Using part (i) of the proposition and the Shimura--Taniyama reciprocity law, we obtain an explicit description of the image of $x_g$ in $G^\sim$. Using part (ii), we prove that the image of $x_g$ in $\pi^*(G,X)$ coincides with the natural image of the point $i\circ h \in X$. 

For Galois gerbs there is a parallel construction of integral special points. By comparing the images of integral special points in $G^\sim \times \pi^*(G,X)$ in the two contexts, we obtain Proposition \ref{introprop:amicable}. 

Given Proposition \ref{introprop:amicable}, we prove Theorem \ref{introprop:Hodge}  in two stages. In the first stage, we construct subsets $\cca' \subset \cca$ and $\ica' \subset \ica$ and a bijection $\bij: \cca' \isom \ica'$ such that for each $\phi \in \cca'$ there exists a special point datum which induces both $\phi$ and $\bij(\phi)$. By Proposition \ref{introprop:amicable} we know that $( \bij(\phi) ,\phi)$ is amicable for such $\phi$. In the second stage, we use Galois cohomological twistings on both sides to extend $\bij$ to a bijection $\cca \isom \ica$. For this, we need the following result.

\begin{introprop}[see Theorem \ref{thm:matchable pair 2}]\label{introprop:equivariance}
	If $( \II, \phi)$ is an amicable pair, then the twists of $\II$ and $\phi$ under a common Galois cohomology class still constitute an amicable pair. 
\end{introprop} To show this we again utilize integral special points. More precisely, we use the fact that for each amicable pair $(\II, \phi)$ and each maximal torus $T \subset I_\phi$, there exist two special point data of the form $(T,i,h)$ and $(T,i',h)$, with $i$ conjugate to $i'$ under $G(\Qbar)$, such that 
\begin{align}\label{eq:abundance}
	\II= \II(T,i,h) \quad  \text{ and } \quad  \phi = \phi(T,i',h) \text{ (up to conjugacy)}.
\end{align}  This essentially follows from the CM lifting theorem in \cite{kisin2012modp} and an analogous result for Galois gerbs. Using this fact, we are able to reduce Proposition \ref{introprop:equivariance} to studying the twistings of integral special points, which can be done explicitly (see Proposition \ref{prop:LD2}). 

This concludes our discussion of the proof of Theorem~\ref{introprop:Hodge}. At the time of writing, it remains an open question whether the bijection $\cca \isom \ica$ can be chosen compatibly with the natural maps from the set of special point data to $\cca$ and $\ica$. Since the latter two maps are surjective, such a compatible bijection is necessarily unique, and establishing its existence is expected to provide deeper insight into the Langlands--Rapoport conjecture.

Having shown Theorem \ref{introprop:Hodge}, we proceed to prove the strengthened version of Conjecture \ref{introconj:LRtau} in the case of Hodge type as follows. Fix a bijection $\bij$ as in Theorem \ref{introprop:Hodge}. For each $\phi \in \cca$, since $(\bij(\phi),\phi)$ is amicable, after choosing some ``marking points'' we can find an element $\tau(\phi) \in I_{\phi}^{\ad}(\A_f)$ and a $G(\A_f^p)\times \Gal(\overline \FF_{q}/\FF_q)$-equivariant bijection $$ \mathscr B_{\phi}: \bij(\phi) \isom \varprojlim_{K^p} I_{\phi}(\Q)_{\tau(\phi)} \backslash X(\phi)/K^p$$ which commutes with the natural maps from both sides to $\pi(G,X)$; see Proposition \ref{prop:marked bij}. Here both $\tau(\phi)$ and $\mathscr B_{\phi}$ depend on choices, but the image of $\tau(\phi)$ in the set $$ \fkH(\phi) : =  (\text{kernel of $I_\phi(\A_f)$ acting on $\pi(G,X)$}) \backslash I_{\phi}^{\ad}(\A_f)/ I_{\phi}^{\ad} (\Q)$$  depends only on the pair $(\bij(\phi),\phi)$; see Proposition \ref{prop:well def of tau}. The set  $\fkH(\phi)$ turns out to be an abelian quotient group of $I_{\phi}^{\ad}(\A_f)$, and the map  $I_{\phi}^{\ad}(\A_f) \to \cH(\phi)$ factors through $ \fkH(\phi)$. Moreover, just like $\cH(\phi)$, the group $\fkH(\phi)$ depends only on the conjugacy class of $\phi^{\Delta}$ up to canonical isomorphism. 

By grouping together the bijections $\mathscr B_\phi$, we obtain a $G(\A_f^p)\times \Gal(\overline \FF_{q}/\FF_q)$-equivariant bijection as in (\ref{eq:intro LRbij}), satisfying the extra condition that it commutes with the maps from both sides to $\pi(G,X)$. With more care, we can also arrange the bijection to respect $\mathscr A$-symmetry. In addition, we obtain the following strengthened versions of conditions (i) and (ii) in Conjecture \ref{introconj:LRtau}.

Strengthening condition (i), we prove that the image of $\tau(\phi)$ in the more refined group $\fkH(\phi)$ depends only on the conjugacy class of $\phi^{\Delta}$. This is a byproduct of the proof of Proposition \ref{introprop:equivariance}.

Strengthening condition (ii), we prove that for each maximal torus $T$ in $I_\phi$ the image of $\tau(\phi)$ in $\fkH(\phi)$ lies in the kernel of the top row of a natural commutative diagram 
$$ \xymatrix{\fkH(\phi) \ar[d] \ar[r] & \fkH(\phi, T) \ar[d]\\ \cH(\phi) \ar[r] & \cH(\phi,T).}$$  
Roughly speaking, $\fkH(\phi, T)$ is defined such that its dependence on $G$ is only via $G_{\der}$ modulo an error from $\cG^{\ab}(\Z_p)$, the hyperspecial subgroup of $G^{\ab} = G/G_{\der}$. For its precise definition, see \S \ref{subsec:refined tori-rational}. In order to prove the vanishing of $\tau(\phi)$ in $\fkH(\phi,T)$, we can again arrange (\ref{eq:abundance}) for the  amicable pair $(\II = \bij(\phi), \phi)$ and the given $T$. The proof then follows from an explicit computation using integral special points; see Proposition \ref{prop:LD1}.   

\bigskip
\textbf{(II) Proving the original  Conjecture \ref{introconj:LRtau} for abelian type.}
\bigskip

We now explain how we bootstrap the strengthened version of Conjecture \ref{introconj:LRtau} in the case of  Hodge type to deduce the conjecture in the case of abelian type.  As the bijection (\ref{eq:intro LRbij}) preserves the extra structures in the Hodge-type setting sufficient for reconstructing abelian-type Shimura varieties, it can be bootstrapped to yield a bijection in the abelian-type case, while simultaneously ensuring that condition (i) remains satisfied.

For condition (ii) in the abelian-type setting, however, reducing it to even the strengthened version of condition (ii) in the Hodge-type setting (namely the vanishing of $\tau(\phi)$ in $\fkH(\phi,T)$) meets an obstruction. This obstruction arises because in the Hodge-type setting, the dependence of the group $\fkH(\phi,T)$ on $G$ is not solely via $G_{\der}$. Nonetheless, this obstruction vanishes if the Shimura datum $(G,X)$ in the Hodge-type setting satisfies the following hypothesis:
\begin{itemize}
	\item The $\Gal(\Qbar/\QQ)$-module $X_*(G^{\ab})$ is generated by the element induced by the Hodge cocharacter of $(G,X)$. 
\end{itemize} 
Finally, using a construction due to Deligne \cite{deligne1982motifs}, we show that every Shimura variety of abelian type can be reconstructed from one of Hodge type satisfying this  hypothesis. This completes the proof of Theorem \ref{introthm:LRtau}.

\section{Integral \texorpdfstring{$p$}{p}-adic Hodge theory}\label{subsec:intro p-adic Hodge} 
We now explain the ingredients from $p$-adic Hodge theory that go into the construction and study of the integral special points in 
\S \ref{subsec:intro Hodge type}.

Let $P$ be a reductive group over $\QQ_p$, and $\cP$ a parahoric group scheme over $\ZZ_p$ with generic fiber $P$. Write $\Rep_{\Z_p} \cP$ for the category of $\cP$-representations on finite free $\ZZ_p$-modules.   Let $F/\QQ_p$ be a finite extension with residue field $k$. Write $W$ for $W(k)$. Fix a uniformizer in $F$ and let $E \in W[u]$ be its minimal polynomial over $W$. Let $\gS = W\lps u \rps $.  
Fix a crystalline representation $$\rho : \Gal(\ol F/F) \To \cP(\Z_p).$$ Then for each $L \in \Rep_{\Z_p}  \cP$, we can view $L$ as a $\Gal(\ol F/F)$-stable lattice in the crystalline representation $(L[1/p], \rho)$. Using the functor $\gM$ in \cite{kisin2006crystalline} we attach to $L$ a BK module $\gM(L)$ over $\gS$.  
Using a purity result\footnote{The special cases for connected reductive group schemes and  parahoric group schemes with tamely ramified generic fibers were previously proved in \cite{kisin2010integral} and \cite{KisinPappas} respectively.} proved recently by Ansch\"utz \cite{Anschutz}, we prove the following result.
\begin{introprop}[see \S \ref{para:KL} and Lemma \ref{keylemma}]\label{introprop:KL}
	The $\otimes$-functor $\omega_{\rho}$ sending $L \in \Rep_{\Z_p} \cP$ to the underlying $\gS$-module of $\gM(L)$ is isomorphic to the standard fiber functor $\idfunc_{\gS}: L \mapsto L \otimes_{\ZZ_p} \gS$.
\end{introprop} We choose such an isomorphism $\idfunc_{\gS} \isom \omega_{\rho}$. The Frobenius on the BK module $\gM(L)$ is then given by an element $\delta_{\gS} \in \cP(\gS[1/E])$ (which is independent of $L$). By considering the de Rham realization of $\gM$ and choosing a splitting of the filtration, we obtain a cocharacter $\mu: \mathbb G_{m,F} \to P_F$.

\begin{introprop}[see Corollary \ref{cor:comparing delta and mu}] \label{introprop:p-adic Hodge}
	Let  $\delta$ be the image of $\delta_{\gS}$ under the map $\cP(\gS[1/E])  \to P(W[1/p])$ induced by $\gS[1/E] \to W[1/p], u \mapsto 0 $. Then the image of $\delta$ under the Kottwitz homomorphism $P(\Qpur) \to \pi_1(P)_{\Gamma_{p,0}}$ is equal to the image of $\mu$. 
\end{introprop}  

The proof of Proposition~\ref{introprop:p-adic Hodge} proceeds in two steps. First, using fundamental properties of the functor $\gM$, we show that the image of $\delta_{\gS}$ in $ P(\widehat \gS_{(E)}[1/E])$ lies in the double coset\footnote{Here $\widehat \gS_{(E)}$ is a complete DVR with residue field $F$ and uniformizer $E$. There is a unique $W[1/p]$-linear multiplicative lifting $F \to \widehat \gS_{(E)}$, making $\widehat \gS_{(E)}$ an $F$-algebra. Hence $\mu(E) \in P(\widehat \gS_{(E)}[1/E])$ is defined.} $P(\widehat \gS_{(E)}) \mu(E) P(\widehat \gS_{(E)}) $. Second, we develop a family version of the Kottwitz map, whose key properties 
are summarized in the proposition below and may be of independent interest. 
Applying both parts of this result, we deduce Proposition~\ref{introprop:p-adic Hodge} 
directly from the double coset containment.

\begin{introprop}[see \S \ref{subsec:Kottwitz hom}]\label{introprop:Kottwitz hom} Let $F$ be a field of characteristic $0,$ $R$ an 
	$F$-algebra, and $v$ a discrete valuation on $R.$ Then for every reductive group $P$ over $F$ there is a natural map 
	$$\kappa_P^v: P(R) \To \pi_1(P)_{\Gal(\overline F/F)},$$ 
	which is functorial in $P$ and  in $(R,v)$. For $P = \GG_m$, it is given by $v : R^{\times} \to \ZZ$. 
	Moreover the following statements hold for all $P$. 
	\begin{enumerate}
		\item Assume that $\Spec R$ has trivial Picard group. For discrete valuations $v_1,\ldots, v_n$ on  $R$ and integers $a_1,\ldots, a_n$ such that $ \sum_i a_i v_i$ vanishes on $R^\times$, we have $$\sum_i a_i \kappa_P^{v_i}(g)= 0 $$ for all $g\in P(R)$. 
		\item Let $v$ be a discrete valuation on $F$ with valuation ring $\oo_F$. For every smooth affine group scheme $\mathcal P$ over $\oo_F$ with connected fibers extending $P$, we have $$\kappa_P^{v}(\mathcal P(\oo_F)) = 0. $$   
	\end{enumerate}
\end{introprop}

We end this section by explaining the following applications of the above local theory in the global situation of \S \ref{subsec:intro Hodge type}, following the notation there.

\textbf{(A) The construction of (\ref{eq:eta}).} Up to enlarging $F$, the $\Gal(\overline F / F)$-action on $\Lambda_A$ factors through a crystalline representation 
$$
\rho : \Gal(\overline F/ F) \to \cP(\Z_p).
$$ By Proposition \ref{introprop:KL} we can choose an isomorphism $\idfunc_{\gS} \isom \omega_{\rho}$. Evaluating at $\Lambda_A$ and reducing modulo $u$, we obtain an isomorphism $$\Lambda_ A \otimes_{\Z_p} W \isom \gM(\Lambda_A)/u\gM(\Lambda_A). 
$$ 
On the other hand, the integral comparison theorem provides an isomorphism  
$$ 
\gM(\Lambda_A)/u\gM(\Lambda_A) \isom \coh^1_{\cris}(\bar A/W). 
$$ 
Composing the above two isomorphisms and base changing from $W$ to $\Zpur$ we obtain the desired isomorphism~\eqref{eq:eta}.\footnote{By a result of T.~Liu \cite{Liucompatibility} this construction is insensitive to enlarging $F$.}

\textbf{(B) The proof of Proposition \ref{introlem:1}.}  Let $\rho: \Gal(\overline F/F) \to \cG(\Z_p)$ and $\rho_T: \Gal(\overline F/F) \to \cT^{\circ}(\Z_p)$ be the crystalline representations giving rise to the $\Gal(\overline F/F)$-action on $\Lambda[1/p]$ (up to enlarging $F$). Part (i) of Proposition \ref{introlem:1} follows from Proposition \ref{introprop:p-adic Hodge} applied to $(\rho_T, \cT^{\circ})$ and the fact that $\mu = -\mu_h$. For part (ii), we compare $\cT^\circ$ and $\cG$ indirectly via the canonical $\Z_p$-homomorphisms $\cG \to \cG^{\ab}$ and $\cT^\circ \to \cG^{\ab}$, where $\cG^{\ab}$ denotes the connected N\'eron model over $\Z_p$ of $G^{\ab}$. Using the fact that $\rho$ and $\rho_T$ induce the same crystalline representation $\Gal(\overline F /F) \to \cG^{\ab}(\Z_p)$ and the functoriality of the local theory, one immediately obtains the weaker version of the desired statement with $G_{\sconn}(\Qpur)$ replaced by $G_{\der}(\Qpur)$. To obtain the full statement, we proceed by lifting $\rho_T$ to a $z$-extension of $G$ (again up to enlarging $F$); see Proposition~\ref{prop:intstrcomp}. This strategy is made possible by the fact that the local theory involves only crystalline representations, not more rigid geometric objects such as abelian varieties or $p$-divisible groups.

\section{Stabilization} \label{subsec:intro stabilization} 

We close the introduction with a summary of the stabilization, namely the proof of Theorem \ref{thm:intro stabilization}.  

A key idea in the initial steps (\S\ref{sub:initial}) is to perform an inverse Fourier transform on the finite abelian group $\mathfrak K(I_0/\Q)$ (defined in terms of $\gamma_0$), of which the Kottwitz invariant $\alpha(\fkc)$ is a character. This takes away the delicate global condition $\alpha(\fkc)=0$ from \eqref{eq:intro refined point counting}, and the global coefficients in the formula turn out to depend only on $\gamma_0$. 
The outcome of \S\ref{sub:initial} is a sum over the set of $\gamma_0$ and $\kappa\in \mathfrak K(I_0/\Q)$ of a summand which is a product of purely local terms, setting the stage for local transfer of orbital integrals in \S\ref{sub:local-transfer}.

Up to a Kottwitz sign and a volume term at $\infty$, the summand for each $(\gamma_0,\kappa)$ is the sum of 
$$\langle \alpha(\gamma_0,a,[b]),\kappa \rangle 
O_{\gamma}(f^p) TO_{\delta}(\phi_{mr}) \tr\xi(\gamma_0)$$
over $(a,[b])$ making $(\gamma_0,a,[b])$ a Kottwitz parameter. Naturally this is broken into three parts: (i) away from $\{p,\infty\}$: $\kappa$-orbital integrals, namely a linear combination of orbital integrals over the (local) stable conjugacy class of $\gamma_0$ weighted by certain character values, (ii) at $p$: a twisted analogue of a $\kappa$-orbital integral, and (iii) at $\infty$: $\tr\xi(\gamma_0)$ weighted by a character value (up to a volume factor). 
Rewriting (i) and (ii) as stable orbital integrals of $f^{H,p,\infty}$ and $f^H_p$ on endoscopic groups $H$ relies on the fundamental lemma and its twisted version, together with relevant results on endoscopic transfer. There is an extra subtlety  at $p$ coming from the $p$-part of $\langle \alpha(\gamma_0,a,[b]),\kappa \rangle$, which was handled by Kottwitz in \cite[App.~A]{Mor10}. The term (iii) at $\infty$ is identified with a stable orbital integral on an endoscopic group by constructing $f^H_\infty$ explicitly in terms of pseudo-coefficients averaged over discrete series $L$-packets.

Throughout \S\ref{sub:local-transfer} there is an extra layer of complexity caused by $z$-extensions and fixed central characters, so we really end up with stable orbital integrals of $f^{H_1}=f^{H_1,p,\infty}f^{H_1}_p f^{H_1}_\infty$ on central extensions $H_1$ of endoscopic groups $H$. 
Here a useful fact is that once a $z$-extension $G_1$ of $G$ is fixed, it induces central extensions $H_1$ of endoscopic groups $H$ for $G$. To transfer functions with fixed central characters (thus the functions are not compactly supported in general), the main point is that the transfer factors enjoy an equivariance property with respect to the translation by central elements.

The main idea in the final steps (\S\ref{sub:final}) is to reparametrize the pairs $(\gamma_0,\kappa)$ in terms of elliptic endoscopic groups and stable conjugacy classes on them. Doing so requires a careful analysis of global invariants, including the obstruction to an adelic conjugacy class containing a rational element; cf.~proof of Lemma \ref{lem:Q-conj-in-A-conj}. The arguments in \S\ref{sub:final} are mostly summarized by the chain of equalities in the proof of Theorem \ref{thm:end-of-stabilization}.

\section*{Postscript: AI and \texorpdfstring{$\tau$}{tau}}
\phantomsection
\addcontentsline{toc}{section}{Postscript: AI and \texorpdfstring{$\tau$}{tau}}
This project was initiated in 2014. The initial manuscript was submitted in 2021, and the main revisions were largely completed by January 2026.

In the fall of 2026, we began a final revision of the manuscript, just as models such as GPT Astra became available. It was a delight to one of us (M.K.) that Astra quickly reproduced calculations we had done thirteen years earlier, but for which we had kept no notes. These calculations showed that, for many groups $I$, the vanishing of the class of $\tau$ when restricted to maximal tori (cf.~Theorem \ref{thm:matchable pair 1}) forces the class itself to vanish. They also gave counterexamples showing that the condition on tori was insufficient for certain unitary inner forms of even degree. Unfortunately, these calculations are only enough to prove the Langlands--Rapoport conjecture for certain special $G$ and $p.$ Even for $\mathrm{GSp}_{2g}$ this obstruction arises when $g> 1.$

We then used Astra to proofread the manuscript. It identified a large number of typos and mathematical inaccuracies, all minor and local in nature, which the authors corrected.  

On October 6, 2026, OpenAI announced a proof of the Hodge conjecture for CM abelian varieties. This case of the Hodge conjecture was already known to imply the Tate conjecture for abelian varieties over finite fields. With this result in hand, it becomes possible to define the motivic Galois gerb using a Tannakian category of motives generated by abelian varieties. This should prove that $\tau$ is trivial and allow a much more natural construction of the correspondence in \S\ref{sec:Shimura Hodge type}. Given the late stage of this project at the time of the announcement, we have chosen to preserve our original treatment in \S\ref{sec:Shimura Hodge type}.

\section*{Organization of the monograph}
\phantomsection
\addcontentsline{toc}{section}{Organization of the monograph}
In \S \ref{sec:pcf}, after some group-theoretic preparations, we state the conjectural point counting formula as in (\ref{eq:intro refined point counting}) for a general Shimura datum (with hyperspecial level at $p$; see Conjecture \ref{conj:point counting formula}). 
Of note is \S \ref{subsec:Kottwitz hom}, in which we generalize the Kottwitz homomorphism to families, as discussed in Proposition \ref{introprop:Kottwitz hom} above. 

In \S \ref{sec:LRtau}, we first recall the formalism of Galois gerbs and the Langlands--Rapoport Conjecture, and then state the Langlands--Rapoport--$\tau$ Conjecture (see Conjecture \ref{hypo about LR}). For the formulation of the conjecture we need a general twisting construction for admissible morphisms. We study this in 
\S \ref{subsec:preparation LR}. 

In \S \ref{sec:study of pairs}, we show that the Langlands--Rapoport--$\tau$ Conjecture implies the desired point counting formula. The key step is the assignment of a Kottwitz parameter (i.e., a summation index in the point counting formula) to a gg LR pair. We study this construction in \S \ref{subsec:attach} 
in the presence of general twisting elements $\tau(\phi)$. In \S \ref{subsec:axiomatic setting}, we study the special case when the elements $\tau(\phi)$ are controlled as in the Langlands--Rapoport--$\tau$ Conjecture. Roughly speaking, we show that the set of Kottwitz parameters that can arise from this construction is unaffected by the presence of these twisting elements. 

In \S \ref{sec:crystalline}, we develop the input from integral $p$-adic Hodge theory. We define\index{F-isocrystal@$F$-isocrystal!integral with $\cG$-structure} \emph{integral $F$-isocrystals with $\cG$-structure} in \S \ref{subsec:intFisoc}, and attach them to $\cG(\Z_p)$-valued crystalline representations in \S \ref{subsec:crystalline lattices} (for $\cG$ a $\ZZ_p$-group scheme satisfying ``$\gS$-purity''). Most of the content in this chapter is local in nature, with the exception of the latter parts of \S \ref{subsec:crystalline-representations-with-g-structure}, which treat the global abelian Galois representations related to the Shimura--Taniyama reciprocity law. 

In \S \ref{sec:Shimura Hodge type} and \S \ref{sec:abelian type}, we prove the Langlands--Rapoport--$\tau$ Conjecture for Shimura varieties of abelian type. We prove in \S \ref{subsec:LR str} the strengthened version of the conjecture in the case of Hodge type, and then bootstrap this result and deduce the Langlands--Rapoport--$\tau$ Conjecture in the case of abelian type in \S \ref{sec:abelian type}.   Chapter  \ref{sec:Shimura Hodge type} is the core of the arithmetic part of this work, and is arguably the most technical. We refer the reader to its own introduction for its structure.   

We devote 
\S\S \ref{sec:prelim-for-stab}--\ref{sec:stabilization} to stabilization. 
In \S \ref{sec:prelim-for-stab}, we have preparatory discussions on central character data in \S\ref{sec:TF-fixed-central}, endoscopic data and $z$-extensions in \S \ref{sub:endoscopic-data}, Galois cohomology invariants in \S \ref{sub:cohomological}, and the Langlands--Shelstad--Kottwitz transfer in \S \ref{sub:LSK-transfer}. By implementing central characters and $z$-extensions, we make it unnecessary to assume any undesirable technical hypothesis such as cuspidality of $Z_G$ or simple connectedness of $G_{\der}$.

In \S \ref{sec:stabilization}, we present the stabilization in three steps following Kottwitz. We rewrite the point counting formula in terms of adelic $\kappa$-orbital integrals (\S\ref{sub:initial}), transfer $\kappa$-orbital integrals to stable orbital integrals on central extensions of endoscopic groups (\S\ref{sub:local-transfer}), and then finish by reorganizing the terms into the sum of stable distributions intrinsic to the endoscopic groups and their central extensions (\S\ref{sub:final}).

Finally, in \S \ref{sec:spectral interpretation}, after recalling the general stable trace formula, we indicate what extra information and steps are needed, in addition to our main results, for understanding the cohomology of Shimura varieties unconditionally.

\newpage
\section*{Conventions}
\phantomsection
\addcontentsline{toc}{section}{Conventions}
\begin{itemize}
	\item When $G$ is a group acting on a set $X$ on the left, we denote the action map $G \times X \to X$ by  $(\rho, x) \mapsto \leftidx^{\rho} x$, or $(\rho, x) \mapsto \rho (x)$. We shall \emph{not} use the notation $x^{\rho}$.
	\item If $g, h$ are elements of a group, we define $\Int(g) h $ to be $ghg^{-1}. $
	\item Given a field $F$, we denote by $\ol F$ a fixed algebraic closure, and write $\Gamma_F$\index[Oct2026-sorted-notation]{zz03 gamma 0 00 00 gamma subscript f@$\Gamma_F$} for $\Gal(\ol F/F)$. 
	\item For every place $v $ of $\QQ$ we fix $\ol {\QQ}_v$ (with $\ol {\QQ}_\infty = \CC$) and fix an embedding ${\overline {\QQ}} \to \overline{\QQ}_v$. The Galois groups will sometimes be abbreviated as follows:\index[Oct2026-sorted-notation]{zz03 gamma 0 00 00 gamma@$\Gamma$}\index[Oct2026-sorted-notation]{zz03 gamma 0 00 00 gamma subscript v@$\Gamma_v$}\index[Oct2026-sorted-notation]{zz03 gamma 0 00 00 gamma subscript v 0@$\Gamma_{v,0}$} $$\Gamma = \Gamma_\QQ  = \Gal({\overline {\QQ}} /\QQ),\quad \Gamma_v = \Gamma_{\QQ_v} = \Gal (\overline{\QQ}_v / \QQ_v).$$ Via the fixed embedding $\ol \QQ \to \ol \QQ_v$, we view $\Gamma_v$ as a subgroup of $\Gamma$. For $v$ non-archimedean, we write $\Gamma_{v,0}$ for the inertia subgroup of $\Gamma_{v}$.
	
	\item If $L$ is a subfield of $\overline \QQ$ and $v$ is a place of $\QQ$, we denote by $L_v$ the completion of $L$ inside $\overline \QQ_v$, with respect to the fixed embedding $\overline \QQ \to \overline \QQ_v$.
	\item By a \emph{CM field} we mean a totally imaginary quadratic extension of a totally real number field, contained in the fixed $\Qbar$.
	\item For a prime $p$, we denote by $\QQ_{p}^{\ur}$\index[Oct2026-sorted-notation]{q 0 06 00 q subscript p superscript ur@$\mathbb Q_p^{\ur}$} the maximal unramified extension of $\QQ_p$ in $\ol\QQ_p$. We write $\QQ_{p^n}$\index[Oct2026-sorted-notation]{q 0 06 00 q subscript p superscript n@$\mathbb Q_{p^n}$} for the degree $n$ unramified extension of $\QQ_p$ in $\QQ_{p}^{\ur}$, and write $\ZZ_{p^n}$\index[Oct2026-sorted-notation]{z 0 06 00 z subscript p superscript n@$\mathbb Z_{p^n}$} for its ring of integers. We write $\ZZ_p ^{\ur}$\index[Oct2026-sorted-notation]{z 0 06 00 z subscript p superscript ur@$\mathbb Z_{p}^{\ur}$} for the strict henselization of $\ZZ_p$, i.e., $\ZZ_p^{\ur}= \bigcup _{n} \ZZ_{p^n}$.
	
	\item For a prime $p$, we denote by $\LL$\index[Oct2026-sorted-notation]{q 0 06 21 breve q subscript p@$\breve {\mathbb Q}_p$} the completion of $\QQ_{p}^{\ur}$, and denote by $\breve{\ZZ}_p$\index[Oct2026-sorted-notation]{z 0 06 21 breve z subscript p@$\breve{\mathbb Z}_p$} its ring of integers. We fix an embedding $\Qpbar \hookrightarrow \overline {\breve \QQ}_p$, and thereby identify $\Gamma_{\LL} = \Gal(\overline {\breve \QQ}_p / \LL) $ with the inertia subgroup $\Gamma_{p,0}$ of $\Gamma_p$. We denote by $\sigma$\index[Oct2026-sorted-notation]{zz18 sigma 1 00 00 sigma@$\sigma$} the arithmetic $p$-Frobenius in $\Aut(\LL)$.
	\item We denote by $\A, \A_f, \A_f^p$\index[Oct2026-sorted-notation]{a 0 06 00 a@$\mathbb A$}\index[Oct2026-sorted-notation]{a 0 06 00 a subscript f@$\mathbb A_f$}\index[Oct2026-sorted-notation]{a 0 06 00 a subscript f superscript p@$\mathbb A_f^p$}\index[Oct2026-sorted-notation]{a 0 06 00 a subscript f superscript star@$\mathbb A_f^*$} respectively the adeles over $\QQ$, the finite adeles over $\QQ$, and the finite adeles away from $p$. We also denote by $\A^*_f$ the product $\A^p_f \times \QQ_p^{\ur}$, when $p$ is clear in the context. 
	\item For a perfect field $k$, we write $W(k) $ for the ring of Witt vectors over $k$. When $k = \FF_{p^n}$ we identify $W(k)$ with $\ZZ_{p^n}$.
	
	\item Reductive groups over a field are always assumed connected. 
	\item For a reductive  group $I$ over a field, we write $I_{\der}, I_{\sconn}, I^{\ad}$ for the derived subgroup, the simply connected cover of the derived subgroup, and the adjoint group respectively. We define $I^{\ab} $ to be $I/I_{\der}$, the maximal torus quotient of $I$. We often use the symbol $\varrho$ to denote the natural map $I_{\sconn} \to I$. 
	\item For a reductive group $I$ over $\RR$, we write\index[Oct2026-sorted-notation]{0000 plus 1 00 00 superscript and subscript@$+$ (superscript and subscript)} $I(\R)^+$ or $I(\RR)^0$ for the identity component of the real Lie group $I(\R)$, and write $I(\R)_+$ for the preimage of $I^{\ad}(\R)^+$ under $I(\R) \to I^{\ad}(\R)$. If $I$ is defined over $\QQ$, we write $I(\Q)^+$ for $I(\Q) \cap I(\R)^+$, and write $I(\Q)_+$ for $I(\Q) \cap I(\R)_+$.
	\item All group cohomology classes or cocycles for profinite groups (e.g.~Galois groups) are understood in the continuous sense. We shall denote a $1$-cochain by $(g_{\rho})_\rho$, or $\rho \mapsto g_{\rho}$, or simply $g_{\rho}$. The cocycle relation is  $g_{\rho \sigma} = g_{\rho} \leftidx^{\rho} g_{\sigma}$.
	\item
	Let $I$ be a linear algebraic group over a field $F$ of characteristic zero. Let $F'/F$ be a Galois extension. We denote by $Z^1(F'/F,I(F'))$\index[Oct2026-sorted-notation]{z 0 00 00 z superscript 1 1 minus cocycles@$Z^1$ (1-cocycles)} the set of continuous $1$-cocycles $\Gal(F'/F) \to I(F')$. Denote by $\coh^1(F'/F, I(F'))$ the corresponding cohomology set. When $F' = \overline F$ is an algebraic closure of $F$, we write $Z^1(F, I)$ and $\coh^1(F,I)$. When $I$ is reductive, for $\tau \in I^{\ad}(F)$, its \emph{image} in $\coh^1(F,Z_I)$ is understood to be the  class of the cocycle $\rho \mapsto \tilde \tau^{-1}\leftidx^{\rho} \tilde \tau$, where $\tilde \tau \in I(\overline F)$ is an (arbitrary) lift of $\tau$.
	\item When $I$ is a reductive  group over $F$, we denote by\index[Oct2026-sorted-notation]{ss 1 01 00 ss semisimple@$\mathrm{ss}$ (semisimple)} $I(F)_{\semi}$ the set of semi-simple elements of $I(F)$.
	\item We denote by $\mathbb S$\index[Oct2026-sorted-notation]{s 0 06 00 s@$\mathbb S$} the Deligne torus $\Res_{\CC/\R} \Gm.$
\end{itemize}

\cleardoublepage

 \numberwithin{equation}{subsection}
  \part{Axiomatic point counting}
 
 \cleardoublepage
 \chapter{The point counting formula}\label{sec:pcf} 
 The goal of this chapter is to state the point counting formula (\ref{pcf}). This generalizes the formula conjectured by Kottwitz (and proved in some PEL cases) by removing the group-theoretic constraints originally imposed. 
 	
 	In \S \ref{subsec:RX_*}--\S \ref{subsec:decent}, we make group-theoretic preparations that will be used throughout the book. Following a review of the setting of Shimura varieties in \S \ref{sec:SV}, we introduce the primary components of the formula in \S \ref{coh Kott} and \S \ref{Kottinv}: \emph{Kottwitz parameters} and their \emph{invariants}. Finally, we state the point counting formula in \S \ref{subsec:stating pcf}.
 	
 	Throughout this chapter, $F$ denotes a field of characteristic zero, $\overline F$ denotes a fixed algebraic closure, and $\Gamma_F$\index[Oct2026-sorted-notation]{zz03 gamma 0 00 00 gamma subscript f@$\Gamma_F$} denotes $\Gal(\overline F/F)$. We will often specialize to the cases where $F$ is  a local or global field. When $F$ is clear from the context, we simply write $\Gamma$ for $\Gamma_F$.

 \section{Abelianized Galois cohomology}\label{subsec:RX_*}
 We review the theory of\index{abelianized Galois cohomology} \emph{abelianized Galois cohomology} of reductive groups and fix notation. 
 \subsection{}\label{para:setting for Gamma modules}
 Let 
 $\GMod$\index[Oct2026-sorted-notation]{zz03 gamma 0 00 00 gamma minus mod@$\Gamma\text{-}\mathrm{Mod}$} (resp.~$\GMod_f$\index[Oct2026-sorted-notation]{zz03 gamma 0 00 00 gamma minus mod subscript f@$\Gamma\text{-}\mathrm{Mod}_f$}) be the  category of discrete $\ZZ[\Gamma]$-modules that are finitely generated (resp.~finite free) over $\ZZ.$ 
 Let $\Mult(F)$\index[Oct2026-sorted-notation]{mult 0 01 00 mult f@$\mathrm{Mult}(F)$} be the category of algebraic groups of multiplicative type over $F,$ and $\Tori(F) \subset \Mult(F)$\index[Oct2026-sorted-notation]{tori 0 01 00 tori f@$\mathrm{Tori}(F)$} the full subcategory of tori. 
 Each of the above four categories is naturally an exact category, and we have the corresponding bounded and unbounded 
 derived categories; see for instance \cite{keller}.	
 
 We have an exact anti-equivalence $$X^* : \Mult(F) \isom  \GMod  $$ sending $G\in \Mult(F)$ to its group of characters $X^*(G) = \Hom_{\ol F}(G_{\ol F}, \mathbb G_m)$. This restricts to an anti-equivalence $\Tori(F) \isom \GMod_f$. We also have the functor $$X_*: \Mult(F) \To \GMod_f $$ sending $G$ to its group of cocharacters $X_*(G) = \Hom_{\ol F} (\mathbb G_m, G_{\ol F})$. This functor is canonically identified with $\Hom_{\ZZ} (\cdot , \ZZ) \circ X^* $, and it restricts to an exact equivalence $\Tori(F) \isom \GMod_f$.
 
 \begin{prop}\label{prop:commdiagder} We have a commutative diagram of equivalences of triangulated categories 
 	$$\xymatrix{
 		\cD^b(\Tori(F)) \ar[r]^\sim \ar[d]^{\sim} &  \cD^b(\GMod_f) \ar[d]^{\sim} \\
 		\cD^b(\Mult(F)) \ar[r]^\sim & \cD^b(\GMod)}	$$
 	where the top functor is induced by $X_*$, the vertical functors are induced by the inclusions $\Tori(F) \hookrightarrow \Mult(F)$ and  $\GMod_f \hookrightarrow \GMod$, and the bottom functor realizes the derived functor $\RX_*$ of $X_*: \Mult(F) \to \GMod$. 
 \end{prop} 
 
 \begin{proof} We first show that inclusion induces an equivalence $\cD^+(\Tori(F)) \isom  \cD^+(\Mult(F))$. By the criterion in \cite[\S 12]{keller}, it suffices to check the following properties:
 	\begin{enumerate}
 		\item Every $A \in \Mult(F)$ admits an injection to some $B \in \Tori(F)$.
 		\item For every surjection $B\to A$ in $\Mult(F)$ with $B \in \Tori(F)$, we have $A\in \Tori(F)$. 
 		\item The subcategory $\Tori(F) \subset \Mult(F)$ is closed under extensions. 
 	\end{enumerate} Via the anti-equivalence $X^*$, (i) follows from the observation that 
 	every object in $\GMod$ is a $\ZZ[\Gamma/U]$-module for some open normal subgroup $U \subset \Gamma$, and hence admits a surjection from $\ZZ[\Gamma/U]^{\oplus r} \in \GMod_f$ for some $r$. (ii) and (iii) correspond to the fact that $\GMod_f$ is closed under subobjects and extensions in $\GMod$. 
 	
 	Now in order to get the equivalence $\cD^b(\Tori(F)) \isom \cD^b(\Mult(F))$, it suffices to show that the acyclicity of a complex in $\Tori(F)$ at a given degree (in the sense of \cite[\S 11]{keller}) is equivalent to acyclicity in $\Mult(F)$. This follows from (ii) above. 
 	
 	Applying the exact anti-equivalence $X^*$ and the commutation of $D^b(\cdot)$ with $(\cdot)^{\mathrm{op}}$, we see that $D^b(\GMod_f) \to D^b(\GMod)$ induced by inclusion is an equivalence.  
 	
 	Since $X_*: \Tori(F) \rightarrow \GMod_f$ is an exact equivalence, the functor on the top is an equivalence. 
 	Now we can fill in the equivalence in the bottom row. The exactness also implies that $X_*: \Tori(F) \rightarrow \GMod_f$ preserves 
 	acyclicity of bounded complexes, which implies that the bottom row realizes the derived functor of $X_*,$ by 
 	\cite[\S I, Thm.~5.1]{HartsRD}.
 \end{proof}

 \subsection{} \label{para:fundgp}
 Let $G$ be a reductive group over $F$. Let $Z_G$ and $Z_{G_{\sconn}}$ denote the centers of $G$ and $ G_{\sconn}$  respectively.\index[Oct2026-sorted-notation]{z 0 00 00 z subscript g@$Z_G$} We write $\mathscr Z_G$\index[Oct2026-sorted-notation]{z 0 05 00 z subscript g@$\mathscr Z_G$} for the complex ($Z_{G_{\sconn}} \to Z_G)$ in $\cD^b(\Mult(F))$, at degrees $-1, 0$. Recall that the algebraic fundamental group\index{algebraic fundamental group}   $\pi_1(G)$\index[Oct2026-sorted-notation]{zz16 pi 1 00 00 pi subscript 1 g@$\pi_1(G)$} is defined to be the $\Gamma$-module $X_*(T)/X_*(\tilde T)$, where $T$ is a maximal torus in $G$ and $\tilde T$ denotes its inverse image in $G_{\sconn}.$ It is independent of the choice of $T$ up to the canonical isomorphism induced by any inner automorphism over $\ol F$ relating the different choices of $T$. The natural map from $\mathscr Z_G$ to the complex $(\tilde T \to T)$ at degrees $-1,0$ is a quasi-isomorphism. Hence we obtain an isomorphism 
 $ \RX_*(\mathscr Z_G) \isom   \RX_*(\tilde T \to T) = X_*(T)/X_*(\tilde T)$ in $\dercat{\GMod}$. The resulting isomorphism $$\RX_*(\mathscr Z_G) \isom  \pi_1(G)$$ is independent of $T$.

 \subsection{}\label{para:coh of mult} We now review the theory of abelianized Galois cohomology, developed by Borovoi \cite{borovoi} and Labesse \cite{Lab99}.
 
 For the rest of this section $F$ will be a local or global field of characteristic zero. We introduce a symbol $?$ as follows. When $F$ is local, $?$ denotes $F$.  When $F$ is global, $?$ denotes one of $F$, $\adele_F/F$, or $\adele_F^S,$ where $S$ is a finite  set of places of $F$, and $\adele_F^S$\index[Oct2026-sorted-notation]{a 0 06 00 a subscript f superscript s@$\mathbb A_F^S$} is the ring of adeles away from $S$. Let $\Gamma = \Gal(\overline F/F)$, and define the discrete $\Gamma$-module\index[Oct2026-sorted-notation]{d 0 02 00 d subscript@$\mathbf D_?$}
 $$ \bD_? = \begin{cases}
 	\overline F ^{\times}, & \text{if } ? = F ,\\
 	(\bar \adele_F^S)^{\times}, & \text{if } ? = \adele_F^S ,\\
 	(\bar \adele_F)^{\times}/ \overline F^{\times}, & \text{if } ? = \adele_F/F.
 \end{cases}$$
 Here $\bar \A_F^S$ denotes $\overline F\otimes _F\A_F^S$.
 
 For any bounded complex $C^{\bullet}$ in $\Mult(F)$, we define the abelian groups 
 $$\coh^i (?, C^{\bullet}) :  = \coh^i (\Gamma, \RX_*(C^{\bullet})\ltensor_{\ZZ}\bD_?), \quad i\in \ZZ ,$$ 
 cf.~\cite[p.~22, p.~26]{Lab99}.  Here the term on the right denotes the continuous group cohomology of the profinite group $\Gamma$.
 
 \subsection{}\label{para:ab Gal coh}
 Let $G$ be a  reductive group over $F$. We define (cf.~\cite[\S 1.6]{Lab99}) \index[Oct2026-sorted-notation]{h 0 02 00 h superscript i subscript ab g@$\mathbf H^i_{\ab}(?,G)$}
 $$\coh^i_{\ab}(?, G) : =  \coh^i(?, \mathscr Z_G) \cong  \coh^i (\Gamma, \pi_1(G) \ltensor_{\ZZ}\bD_? ) .$$   This is a functor in $G$.
 
 If $I \to G$ is a 	homomorphism  of reductive  groups over $F$, we denote by $\pi_1(I\to G)$\index[Oct2026-sorted-notation]{zz16 pi 1 00 00 pi subscript 1 i to g@$\pi_1(I\to G)$} the complex $(\pi_1(I) \to \pi_1(G))$ at degrees $-1,0$, in $\cD^b(\GMod)$. We set (cf.~\cite[\S 1.8]{Lab99})  \index[Oct2026-sorted-notation]{h 0 02 00 h superscript i subscript ab i to g@$\mathbf H^i_{\ab}(?, I \to G) $}
 $$ \coh^i_{\ab}(? , I \to G) : = \coh^i (\Gamma, \pi_1(I \to G) \ltensor_{\ZZ}\bD_? ) .$$ We have a long exact sequence 
 $$ \cdots \to \coh^i_{\ab}(?, I) \to \coh^i_{\ab}(?, G) \to \coh^i_{\ab} (?, I \to G) \to \coh^{i+1}_{\ab}(?, I) \to \cdots  $$

 When $?$ is $F$ or $\A_F^S$, for an arbitrary linear algebraic group $H$ over $F$ we have the usual Galois or adelic cohomology
 $$\coh^i (?, H), \quad i = 0,1,$$ defined to be the continuous cohomology of $\Gamma$ acting on $H(\overline F)$ or $H(\bar \A_F^S)$ respectively. This is a group for $i = 0$ and a pointed set for $i =1$. For reductive $G$, we have natural ``abelianization'' maps\index[Oct2026-sorted-notation]{ab 1 01 00 ab superscript i@$\mathrm{ab}^i$}
 \begin{align*}\
 	\ab^i_? : \coh^i (?, G) \To  \coh^i_{\ab} (?,G), \quad i=0,1,
 \end{align*}
 which is a group homomorphism for $i=0$ and a map of pointed sets for $i=1$.
 By \cite[Prop.~1.6.7]{Lab99}, the map $\ab^1_F: \coh^1(F,G) \to \coh^1_{\ab}(F,G)$ is surjective, and is bijective when $F$ is local non-archimedean. In particular, when $F$ is local non-archimedean, $\coh^1(F,G)$ has a canonical structure of an abelian group.
 
 When $F$ is global, we have
 \begin{align} \label{eq:adelic coh}
 	\coh^1(\adele_F^S, G ) \cong \prod'_{\substack{\text{places $v$ of $F$,} \\ v \notin S}}  \coh^1(F_v, G),\end{align}
 where the restricted product is with respect to the trivial elements; see for instance \cite[p.~298, Cor.~1]{plantonov-rapinchuk}. Analogously we have
 \begin{align}\label{eq:adelic ab coh}
 	\coh^1_{\ab}(\adele_F^S, G ) \cong \bigoplus_{\substack{\text{places $v$ of $F$,} \\ v \notin S}}  \coh^1_{\ab}(F_v, G).
 \end{align} The decompositions (\ref{eq:adelic coh}) and (\ref{eq:adelic ab coh}) are compatible with $\ab^1_{\A_F^S}$ and $\ab^1_{F_v}$. If $S$ contains all archimedean places of $F$, then $\ab^1_{\adele_F^S} : \coh^1(\adele_F^S, G) \to \coh^1_{\ab}(\adele_F^S, G)$ is bijective, giving $\coh^1(\adele_F^S, G)$ the canonical structure of an abelian group.
 
 \subsection{}\label{para:two groups}  
 Assume $?$ is $F$ or $\A_F^S$. We adopt the following notation from \cite{Lab99}: Let $G$ be a linear algebraic group over $F$ and $I$ a subgroup. We write  \index[Oct2026-sorted-notation]{d 0 04 00 d i g@$\mathfrak D(I,G,?)$} $$ \D(I,G ;?     ) := \ker \big(\coh^1(?, I) \to \coh^1(?, G )\big),  $$ which is a pointed set. When $G$ and $I$ are reductive, we write 	\index[Oct2026-sorted-notation]{e 0 04 00 e i g@$\mathfrak E(I,G,?)$}  $$ \E(I,G ; ?)  : =  \ker \big(\coh^1_{\ab}(?, I) \to \coh^1_{\ab}(?, G )\big),$$ which is an abelian group.
 In this case we have a map of pointed sets
 $$\D(I,G ; ?  ) \To \E(I,G ; ?)$$ induced by $\ab^1_?$. This map is bijective in the following two cases:
 \begin{itemize}
 	\item $F$ is local non-archimedean and $?=F$.
 	\item $F$ is global, $? = \adele_F^S$, and $S$ contains all the archimedean places of $F$.
 \end{itemize} This is because in these two cases, $\coh^1(?, U) \cong \coh^1_{\ab}(?,U)$ for $U = I$ or $G$,  as discussed in \S \ref{para:ab Gal coh}.

 \subsection{}\label{para:setting for Borovoi}
 Let $M \in \Gamma_F\text{-}\mathrm{Mod}.$
 When $F$ is global or local non-archimedean,  we set $\mathcal A_F(M) = M_{\Gamma_F, \tors},$\index[Oct2026-sorted-notation]{a 0 03 00 a subscript f cdot@$\mathcal A_F(\cdot)$} the torsion subgroup of the group of coinvariants $M_{\Gamma}.$ 
 When $F$ is local archimedean, we set $\mathcal A_F(M) = \tatecoh^{-1}(\Gamma_F, M).$ (As usual, $\tatecoh^i$ denotes Tate cohomology.) Note that in both cases, there is a canonical embedding $\mathcal A_F(M) \hookrightarrow M_{\Gamma_F, \tors}$.
 
 Now suppose $F$ is global. For each place $v$ of $F$ we fix an embedding $\overline F \to \overline F_v$, and thereby view $\Gamma_{F_v}$ as a subgroup of $\Gamma_F$. For  $M \in \Gamma_F \text{-}\mathrm{Mod},$ set \index[Oct2026-sorted-notation]{b 0 03 00 b subscript f cdot@$\mathcal B_F(\cdot)$}
 \begin{align*}
 	\mathcal B_F(M) :=  \bigoplus_{\substack{\text{places $v$ of $F$} \\ }}    \mathcal A_{F_v} (M),
 \end{align*}  
 and define the map\index[Oct2026-sorted-notation]{p 0 05 00 p cdot@$\mathscr P(\cdot)$} $$\mathscr P(M):\cB_F (M)  \To \mathcal A_F(M) $$ to be the direct sum over $v$ of the maps $$\mathcal A_{F_v} (M) \hookrightarrow M_{\Gamma_v, \tors} \overset{\dagger}{\To}  M_{\Gamma_F, \tors} = \mathcal A_F(M), $$ where $\dagger$ is induced by the identity on $M$. Then $\mathscr P: \mathcal B_F \to \mathcal A_F$ is a natural transformation between functors from $\Gamma_F \text{-}\mathrm{Mod}$ to the category of  abelian groups.
 
 More generally, if $S$ is a finite set of places of $F$, we set \index[Oct2026-sorted-notation]{b 0 03 00 b subscript f superscript s cdot@$\mathcal B_F^S(\cdot)$}	\begin{align*}
 	\mathcal B_F^S(M) =  \bigoplus_{\substack{\text{places $v$ of $F$,} \\ v\notin S }}    \mathcal A_{F_v} (M). \end{align*}

 \begin{prop}\label{prop:apply TN} Let $F$ be local or global, and let $I \to G$ be a homomorphism of reductive  groups over $F$. Assume that the induced map $\pi_1(I) \to \pi_1(G)$ is surjective, and denote its kernel by $K.$ 
 	\begin{enumerate}
 		\item Let $? = F$ when $F$ is local, and let $? = \adele_F/F$ when $F$ is global. Then the exact sequence
 		$$ \coh^0_{\ab}(?, I\to G)  \To  \coh^1_{\ab}(? , I) \To \coh^1_{\ab}(?, G) $$ is canonically isomorphic to the natural sequence
 		$$ \mathcal A_F(K) \To \mathcal A_F(\pi_1(I)) \To \mathcal A_F(\pi_1(G)). $$ If $F$ is global and $? = \A_F^S$ for a finite set $S$ of places of $F$, then the above statement holds with $\cA_F$ replaced by $\cB_F^S$. 
 		\item  When $F$ is global, the commutative diagram with exact rows
 		$$ \xymatrix{\coh^0_{\ab} (\adele_F, I \to G) \ar[r] \ar[d] & \coh^1_{\ab}(\adele_F, I) \ar[d] \ar[r] & \coh^1_{\ab}(\adele_F, G) \ar[d] \\ \coh^0_{\ab} (\adele_F/F, I \to G) \ar[r]  & \coh^1_{\ab}(\adele_F/F, I) \ar[r] & \coh^1_{\ab}(\adele_F/F, G) }$$ is canonically identified with the natural commutative diagram
 		$$\xymatrix{\mathcal B_F(K) \ar[r] \ar[d] ^{\mathscr P(K)} & \mathcal B_F(\pi_1(I)) \ar[d] ^{\mathscr P(\pi_1(I))} \ar[r] & \mathcal B_F(\pi_1(G)) \ar[d]^{\mathscr P(\pi_1(G))} \\\mathcal A_F(K) \ar[r]  &\mathcal A_F(\pi_1(I)) \ar[r] &  \mathcal A_F(\pi_1(G)) .} $$ 
 	\end{enumerate}
 \end{prop}
 \begin{proof} This follows from the results of \cite[\S4]{borovoi} and the fact that $K[1]$ is quasi-isomorphic to $\pi_1(I\to G).$ 
 \end{proof}
 The following lemma will be used in the proof of Proposition \ref{computing E} below.

 \begin{lem}\label{lem:eq of kernels} Let $F$ be local or global. Let $M\to M'$ be an injection in $\GMod$. The natural inclusions
 	\begin{align*}
 		\ker \!\bigg( \mathcal A_F(M) \to \mathcal A_F(M')\bigg)  \hookrightarrow \ker \!\bigg(M_{\Gamma,\tors} \to M'_{\Gamma,\tors}\bigg)   \hookrightarrow \ker \!\bigg(M_{\Gamma} \to M'_{\Gamma}\bigg) 
 	\end{align*}
 	are equalities.
 \end{lem}
 \begin{proof} As in the proof of \cite[Prop.~B.4]{milne92}, when $F$ is global or local non-archimedean, there exists a finite Galois extension $E/F$ such that the actions of $\Gamma$ on $M$ and $M'$ factor through $\Gal(E/F)$ and such that
 	$$	\mathcal A_F(M)  = \tatecoh^{-1}(\Gal(E/F),M), \quad \mathcal A_F(M')  = \tatecoh^{-1}(\Gal(E/F),M'). $$ When $F$ is local archimedean, the above is still true by the definition of $\mathcal A_F$.   Thus the first and third groups are equal and both canonically identified with the cokernel of $$\tatecoh^{-2}(\Gal(E/F), M')  \To  \tatecoh^{-2}(\Gal(E/F), M'/M).$$  
 \end{proof}

 \section{Cohomology with local triviality conditions}\label{sec:inner twisting}
 We review the notion of an \emph{inner twisting} and introduce generalizations that allow  for larger target groups. We define kernels of certain maps between Galois cohomology groups with local triviality conditions and establish several results to be utilized in the sequel.
 
 We begin by letting $F$ be an arbitrary field of characteristic zero. 
 
 \begin{defn}\label{defn:inner twistings}
 	Let $H, H_1$ be algebraic groups over $F$. By an\index{inner twisting} \emph{inner twisting} from $H$ to $H_1$, we mean an $\ol F$-group isomorphism $\psi: H_{\overline F} \isom H_{1, \overline F}$ such that for each $\rho \in \Gamma$ the automorphism $(\leftidx ^{\rho} \psi  )^{-1} \psi $ of $H_{\overline F}$ is inner. Two inner twistings from $H$ to $H_1$ are called \emph{equivalent}, if they differ by an inner automorphism of $H_{\overline F}$. 
 \end{defn}
 \begin{defn}\label{defn:inner form}
 	Let\index{inner form} $H$ be an algebraic group over $F$. By an \emph{inner form of $H$}, we mean a pair $(H_1, [\psi])$, where $H_1$ is an algebraic group over $F$ and $[\psi]$ is an equivalence class of inner twistings from $H$ to $H_1$. When $[\psi]$ is understood from the context we simply say that $H_1$ is an inner form of $H$. By an \emph{isomorphism} between two inner forms $(H_1, [\psi])$ and $(H_1', [\psi'])$ of $H$, we mean an $F$-group isomorphism $H_1 \isom H_1'$ under which $[\psi]$ is identified with $[\psi']$.  
 \end{defn}
 \begin{rem}\label{rem:inner form}There is a bijection from the set of isomorphism classes of inner forms of $H$ to $\coh^1(F, H^{\ad})$, sending $(H_1, [\psi])$ to the class of the cocycle $\rho \mapsto \psi^{-1} \circ \lix^{\rho} \psi$. 
 \end{rem}
 
 The following notions generalize Definition \ref{defn:inner twistings}. 
 \begin{defn}\label{defn:inner transfer} Let $I$ and $G$ be linear algebraic  groups over  $F$. Let $f:I_{\overline F} \to G_{\overline F}$ be a homomorphism over $\ol F$. 
 	\begin{enumerate}
 		\item We call $f$ a\index{pre-inner twisting} \emph{pre-inner twisting}, if $f$ is injective, and for every $\tau \in \Gamma_F$ the homomorphisms $f$ and $\lix^{\tau} f$ differ by post-composing with an inner automorphism of $G_{\ol F}$.
 		\item When $I$ and $G$ are reductive, we say that $f$ is \emph{$F$-rational at $\pi_1$-level}, if the map $\pi_1(I) \to \pi_1(G)$ induced by $f$ is $\Gamma_F$-equivariant. 
 	\end{enumerate}
 \end{defn}
 
 \begin{lem}\label{lem:rational pi_1}
 	Every pre-inner twisting between reductive groups is $F$-rational at $\pi_1$-level.
 \end{lem}
 \begin{proof}
 	This follows from the fact that any inner automorphism over $\ol F$ induces the identity on $\pi_1$. 
 \end{proof}
 \begin{lem}\label{lem:tilde rat}
 	Let $I$ be a reductive group over $F$. Let $h: I' \to I_{\ol F}$ be an $\ol F$-homomorphism from a reductive group over $\ol F$ with finite central kernel inducing an isomorphism between adjoint groups. Then there exists at most one $F$-rational structure on $I'$ such that $h$ is defined over $F$. This exists if and only if the image of the induced map $\pi_1(I') \to \pi_1(I)$ is $\Gamma_F$-stable. 
 \end{lem}
 \begin{proof}
 	First note that if $\phi: I' \to I'$ is an $\ol F$-automorphism over $I_{\ol F}$, then $I' \to I', g \mapsto \phi(g)g^{-1}$ factors through $\ker(h)$ and hence must be constantly $1$. It follows that $\phi= \id $. This implies the uniqueness of $F$-descent datum on $I'$ that is compatible with the $F$-rational structure of $I$. 
 	
 	We claim that if $h_1:I_1' \to I_{\ol F}$ is another homomorphism satisfying the same assumptions as $h$ such that $\pi_1(I_1') \to \pi_1(I)$ and $\pi_1(I') \to \pi_1(I)$ have the same image, then there exists an isomorphism $I_1' \isom I'$ over $I_{\ol F}$. Indeed, by a simple devissage one reduces to the cases where either $I$ is a torus or $\cok(\pi_1(I') \to \pi_1(I))$ is finite. In the first case $h$ and $h_1$ are morphisms between tori with finite kernels, and the claim is clear. In the second case $h$ and $h_1$ are central isogenies between reductive groups, and the claim follows from the classification of these in terms of root data. 
 	
 	Now suppose that the image of $\pi_1(I') \to \pi_1(I)$ is $\Gamma_F$-stable. By the above claim, for each $\tau \in \Gamma_F$ there exists an isomorphism $a_\tau: \tau^*I' \isom I'$ commuting with $\tau^* h: \tau^* I' \to \tau^* I \cong I$ and $h$. By our first observation about automorphisms, $a_\tau$ is unique, and $(a_\tau)_{\tau}$ automatically satisfies the cocycle condition and provides the descent datum. 
 \end{proof}
 \begin{lem}\label{lem:tilde conn}
 	Let $I \hookrightarrow G$ be an injective homomorphism between reductive groups over $F$. Then the inverse image $\tilde I$ of $I$ in $G_{\sconn}$ is connected if and only if $\cok(\pi_1(I) \to \pi_1(G))$ is torsion-free. When this is the case, the group $I\cap G_{\der}$ is also connected, and we have exact sequences \begin{align*}
 		&	0 \To \pi_1(\tilde I) \To \pi_1(I) \To \pi_1(G)   ,\\ &  0 \To \pi_1(I\cap G_{\der}) \To \pi_1(I) \To \pi_1(G^{\ab})  .
 	\end{align*}
 \end{lem}
 \begin{proof} Only the first statement is non-trivial. To prove it we first reduce to the case where $G$ is semi-simple. 
 	If $\tilde I$ is connected, then so is $I \cap G_{\der}$. Now  $I \cap G_{\der}$ is the inverse image of $\ker(I^{\ab} \to G^{\ab})$ along $I \to I^{\ab}$, so it is connected if and only if $\ker(I^{\ab} \to G^{\ab})$ is connected, if and only if $\pi_1(I) \to \pi_1(G^{\ab})$ has torsion-free cokernel. When these conditions hold, we have $\pi_1(I\cap G_{\der}) = \ker (\pi_1(I) \to \pi_1(G^{\ab}))$. On the other hand, $\cok(\pi_1(I) \to \pi_1(G))$ is torsion-free if and only if $\ker (\pi_1(I) \to \pi_1(G^{\ab}))$ maps onto $\pi_1(G_{\der})$ and $\cok(\pi_1(I) \to \pi_1(G^{\ab}))$ is torsion-free. Thus we have reduced to showing that for a (connected) reductive subgroup $J\subset G_{\der}$, its inverse image 
 	$\tilde J$ in $G_{\sconn}$ is connected if and only if $\pi_1(J) \to \pi_1(G_{\der})$ is surjective.

 	Let $S$ be a maximal torus in $J$ and write $\tilde S$ for its inverse image in $\tilde J$.    Write $Z$ for $\ker (G_{\sconn} \to G_{\der})$. If $\tilde J$ is connected, then $\tilde J\to J$ is a central isogeny between reductive groups, from which it follows that $\tilde S$ is a maximal torus in $\tilde J$ and we have $\tilde S = \tilde S^\circ \supset Z_{\tilde J} \supset Z$. Conversely, suppose $\tilde S^{\circ} \supset Z$. Then $Z \subset \tilde J^{\circ}$. Since $\tilde J^{\circ} \to J$ is surjective, we have $\tilde J = \tilde J^\circ$. Thus we are left to show that $\tilde S^{\circ} \supset Z$ if and only if $\pi_1(J) \to \pi_1(G_{\der})$ is surjective. For this, let $T$ be a maximal torus in $G_{\der}$ containing $S$ with inverse image $\tilde T$ in $G_{\sconn}$. Then $X_*(\tilde S^\circ), X_*(S), X_*(\tilde T)$ are all identified with subgroups of $X_*(T)$, and we have $X_*(\tilde S^\circ) = X_*(S) \cap  X_*(\tilde T)$. Let $d$ be the order of $Z$. Now we have $\tilde S^{\circ} \supset Z$ if and only if the kernel of $\tilde S^\circ \to S$ has order $d$, if and only if the index of $X_*(\tilde S^\circ)$ in $X_*(S)$ is $d$. But $d$ is the index of $X_*(\tilde T)$ in $X_*(T)$, so the above  condition is equivalent to the surjectivity of $X_*(S) \to X_*(T)/X_*(\tilde T)$, or equivalently the surjectivity of $\pi_1(J) \to \pi_1(G_{\der})$.  
 \end{proof}
 \subsection{}\label{para:construction of tilde} Using the above two lemmas, we make a general construction that will be used later in several different contexts. Let $I$ and $G$ be reductive groups over $F$, and let $f: I_{\ol F} \to G_{\ol F}$ be an injective homomorphism that is $F$-rational at $\pi_1$-level (e.g., a pre-inner twisting by Lemma \ref{lem:rational pi_1}). Assume that $\pi_1(I) \to \pi_1(G)$ has torsion-free cokernel. Let $f^{\dagger}: I^{\dagger}_{\ol F} \to G_{\der,\ol F}$ and $\tilde f: \tilde I_{\ol F} \to G_{\sconn,\ol F}$ be the pull-back of $f$ along $G_{\der} \to G$ and $G_{\sconn} \to G$ respectively. By Lemma \ref{lem:tilde conn} (applied to $\ol F$), $f^\dagger$ and $\tilde f$ are injections between (connected) reductive groups, and the images of $\pi_1(I^{\dagger})$ and $\pi_1(\tilde I)$ in $\pi_1(I)$ are $\Gamma_F$-stable subgroups since they are kernels of $\Gamma_F$-equivariant maps out of $\pi_1(I)$. By Lemma \ref{lem:tilde rat}, there exist unique $F$-rational structures on $I^{\dagger}_{\ol F}$ and $\tilde I_{\ol F}$ such that the maps $ I_{\ol F}^{\dagger} \to I_{\ol F}$ and $\tilde I_{\ol F} \to I_{\ol F}$ are $F$-rational. We denote the resulting $F$-groups by $I^{\dagger}$ and $\tilde I$. Clearly $f^{\dagger}$ and $\tilde f$ are $F$-rational at $\pi_1$-level. 
 \subsection{}
 \label{subsubsec:notation of Sha}
 Assume that $F$ is  local or global, and let the symbol $?$ be as in \S \ref{para:coh of mult}. Let $I$ and $G$ be reductive  groups over $F$, and let $f:I_{\ol F} \to G_{\ol F}$ be an $\ol F$-homomorphism that is $F$-rational at $\pi_1$-level. Then the $\Gamma_F$-equivariant map $\pi_1(I) \to \pi_1(G)$ induces a homomorphism
 $$ \coh^i_{\ab}(?, I)  \To \coh^i_{\ab}(?, G), $$ which depends on $f$ only via the $G(\ol F)$-conjugacy class of $f$.
 
 Now assume that $F$ is global, and let $S$ be a (finite or infinite) set of places of $F$ containing all the archimedean places. For any linear algebraic group $H$ over $F$, we define \index[Oct2026-sorted-notation]{zzz sha 0 00 00 sha superscript s@$\Sha^S$} $$ 
 \Sha^S (F, H) : = \ker \big(\coh^1(F, H) \to \prod_{v\in S} \coh^1(F_v, H) \big) .$$ If $H$ is reductive , we also define \index[Oct2026-sorted-notation]{zzz sha 0 00 00 sha superscript s subscript ab@$\Sha^S_{\ab}$} $$ \Sha^S_{\ab} (F, H) = \ker \big(\coh^1_{\ab}(F, H) \to \prod_{v\in S} \coh^1_{\ab}( F_v, H) \big ).$$ In this case, by \cite[Thm.~5.12 (i)]{borovoi}, we have a canonical isomorphism
 \begin{align}\label{eq:Sha ab}
 	\Sha^S (F, H) \cong \Sha ^S_{\ab} (F,H )\end{align}
 induced by $\ab^1_F$. In the sequel we shall often make this identification implicitly.  
 
 For $f: I_{\ol F} \to G_{\ol F}$ as above, the induced  homomorphism $\coh^1_{\ab}(F, I) \to \coh^1_{\ab}(F, G)$ restricts to a homomorphism 
 $$  \Sha ^S_{\ab} (F, I) \To  \Sha^S_{\ab} (F, G).$$ We  denote its kernel by \index[Oct2026-sorted-notation]{zzz sha 0 00 00 sha superscript s subscript g@$\Sha^S_G$} $$\Sha^S_G (F, I).$$ Via (\ref{eq:Sha ab}) we also view $\Sha^S_G(F, I)$ as a subset of $\coh^1(F, I)$. 
 If $S$ is the set of all places then we omit it from the notation, and write $\Sha(F,\cdot)$, $\Sha _G(F,\cdot)$, etc.

 \begin{prop}\label{prop:lift to sc} Keep the setting and notation of \S \ref{para:construction of tilde} with $F= \QQ$. Assume that $\pi_1(I) \to \pi_1(G)$ is surjective. Then we have $\Sha _{G_{\sconn}} ^{\infty} (\QQ, \tilde I) \cong \Sha^{\infty}_{\ab} (\QQ, \tilde I)$, and the natural map $\Sha _{G_{\sconn}} ^{\infty} (\QQ, \tilde I) \to \Sha_G^{\infty}(\QQ, I)$ is surjective.   \end{prop}
 \begin{proof} The first statement follows from the triviality of $\coh ^1_{\ab} (\QQ, G_{\sconn})$. In view of the first statement, the second statement amounts to the exactness of 
 	$$ \Sha^{\infty}_{\ab}(\QQ, \tilde I) \To \Sha^{\infty}_{\ab}(\QQ,  I) \To \Sha^{\infty}_{\ab}(\QQ,G).$$ 
 	
 	By \cite[Lem.~5.5.5]{BKH}, for every reductive group $H$ over $\QQ$ we have a functorial isomorphism 
 	$$\Sha^{\infty}_{\ab}(\QQ,H) \cong \big( \pi_1(H)[V_f(E)]_0 \big)_{\Gal(E/\Q), \tors}, $$ where $E/\QQ$ is any sufficiently large finite Galois totally complex extension such that the $\Gamma_{\QQ}$-action on $\pi_1(H)$ factors through $\Gal(E/\QQ)$, $V_f(E)$ denotes the set of finite places of $E$, and $\pi_1(H)[V_f(E)]_0$ denotes the group of finite formal linear combinations $\sum_{v\in V_f(E)} m_v [v]$ with $m_v \in \pi_1(H)$ such that $\sum_v m_v = 0$. 
 	
 	By our assumption and Lemma \ref{lem:tilde conn}, we have a short exact sequence $0 \to \pi_1(\tilde I) \to \pi_1(I) \to \pi_1(G) \to 0$. Choose a common extension $E/\QQ$ as above for the three terms. Clearly 
 	$$ 0 \To \pi_1(\tilde I)[V_f(E)]_0 \To \pi_1(I)[V_f(E)]_0 \To \pi_1(G)[V_f(E)]_0 \To 0$$ remains exact. The associated long exact sequence of $\Gal(E/\QQ)$-homology and the fact that $\coh_1(\Gal(E/\QQ), \cdot)$ is killed by $[E:\QQ]$ imply the 
 	desired exact sequence. 
 \end{proof}
 \begin{rem}
 	In the present work, we only require certain special cases of the above proposition, 
 	for which the argument reduces to the setting where $f: I \hookrightarrow G$ is the inclusion of a full-rank 
 	$\QQ$-subgroup. In that case, the surjectivity assertion admits an elementary proof by working directly with cocycles and applying real approximation.
 \end{rem}
 
 \begin{prop}\label{prop:G(Q)_+}
 	Let $G$ be a reductive  group over $\QQ$, and $T $ a maximal torus in $G$. Let $$\partial : G^{\ab} (\QQ) \To \coh^1(\QQ, Z_{G_{\der}})$$ be the coboundary map associated with the short exact sequence $1 \to Z_{G_{\der}} \to Z_G \to G^{\ab} \to 1$. Then 	  the image of $G(\Q)_+$ (i.e., the inverse image of $G^{\ad}(\RR)^+$ in $G(\QQ)$) under the composite map $$G(\QQ) \To G^{\ab} (\QQ) \overset{\partial}{\To} \coh^1(\QQ, Z_{G_{\der}})$$ is contained in the kernel of $\coh^1(\QQ, Z_{G_{\der}}) \to \coh^1(\QQ, G_{\der}) \times \coh^1(\RR, Z_{G_{\der}})$. 
 \end{prop}
 \begin{proof}
 	For $g\in G(\QQ)$, write $g = g_0 z$ with $g_0 \in G_{\der}(\Qbar)$ and $z \in Z_G(\Qbar)$. Then $\partial (g)$ is represented by the cocycle $\Gamma \ni \rho \mapsto z^{-1}  \lix ^{\rho} z$ in $Z^1(\QQ, Z_{G_{\der}})$. Since $\lix^{\rho} g = g$, this cocycle is equal to $\rho \mapsto  \lix ^{\rho} g_0^{-1} g_0$, and this is equal to $g_0 \lix ^{\rho} g_0^{-1}$ since it commutes with $g_0$.  It follows that $\partial (g)$ has trivial image in $\coh^1(\QQ, G_{\der})$. To show that it has trivial image in $\coh^1(\RR, Z_{G_{\der}})$, note that the image of $g_0$ in $G^{\ad}(\CC)$ is equal to that of $g$, and hence lies in $G^{\ad}(\RR)^{+} \subset \im (G_{\der}(\RR) \to G^{\ad}(\RR))$. Hence $g_0 \in G_{\der}(\RR) Z_{G_{\der}}(\CC)$, and so the cocycle $\rho \mapsto  g_0 \lix ^{\rho} g_0^{-1} $ is trivial in $\coh^1(\RR, Z_{G_{\der}})$. 
 \end{proof}

 \section{The Kottwitz homomorphism}
 \label{subsec:Kottwitz hom}  Let $G$ be a reductive group over $F.$ Recall that, when $F$ is a complete discretely valued field with algebraically closed residue field, the \emph{Kottwitz homomorphism} is a homomorphism $ G(F) \rightarrow \pi_1(G)_{\Gamma_F}$
 which is functorial in $G,$ and for $G = \mathbb G_m$ it is the valuation map $F^{\times} \to \ZZ$. The original construction in \cite[\S 7]{kottwitzisocrystal2} relies on Steinberg's theorem for $F$.
 
 In this section we generalize the construction of the Kottwitz homomorphism to an arbitrary field $F$ of characteristic zero. We shall obtain a homomorphism \index[Oct2026-sorted-notation]{zz10 kappa 1 00 00 kappa subscript g superscript r v@$\kappa_G^{R,v}$} $\kappa_G^{R,v}: G(R) \to \pi_1(G)_{\Gamma_F}$, where $R$ is any $F$-algebra equipped with a discrete valuation $v$. This will allow us to show that the Kottwitz homomorphism
 is constant in certain families.

 \subsection{}Let $\site$ be the big fpqc site of $\Spec F$. Let $\Abshv (\site)$\index[Oct2026-sorted-notation]{a 0 03 00 a b s hv s@$\mathcal A b \mathcal S hv (\site)$} be the category of abelian sheaves on $\site$. We shall view $\Mult(F)$ and $\Tori(F)$ (see \S \ref{para:setting for Gamma modules}) as full subcategories of $\Abshv(\site)$. Let $\cD(\site)$\index[Oct2026-sorted-notation]{d 0 03 00 d s@$\mathcal D(\site)$}\index[Oct2026-sorted-notation]{d 0 03 00 d superscript minus 1 0 s leftbracket rightbracket@$\mathcal D^{[-1,0]}(\site)$} be the derived category of $\Abshv(\site)$, and let $\cD^{[-1,0]}(\site)$ be the full subcategory consisting of $L \in \cD(\site)$ such that $\coh^i(L) = 0$ unless $i\in \set{-1,0}$.
 
 We shall need the formalism of\index{Picard stack} Picard stacks\footnote{We omit the adjective ``strictly commutative'', as that will always be understood.} on $\site$, as in \cite[Expos\'e XVIII, \S 1.4]{SGA4-3}. Following \textit{loc.~cit.}, let $\Ch(\site)$\index[Oct2026-sorted-notation]{c 0 03 00 c h superscript flat s@$\mathcal C h^{\flat} (\site)$} be the category whose objects are the small Picard stacks on $\site$, and whose morphisms are the isomorphism classes of additive functors between Picard stacks. By \cite[Expos\'e XVIII, Prop.~1.4.15]{SGA4-3}, we have an equivalence of categories\index[Oct2026-sorted-notation]{ch 1 01 00 ch@$\mathrm {ch}$} $$\ch: \cD^{[-1,0]}(\site) \isom \Ch(\site).$$ For a complex $C^{\bullet} = (C^{-1} \to C^0)$ in $\Mult(F)$ (at degrees $-1, 0$), $\ch(C^{\bullet})$ is given by the quotient stack $[C^{-1} \backslash C^0]$.

 \subsection{}
 As in Proposition \ref{prop:commdiagder}, we have an equivalence of categories
 $$\RX_*: \dercat{\Mult(F)} \To \dercat{\GFMod}.$$
 Choose a quasi-inverse $\mathcal Y$, with unit $\eta: \id \to \mathcal Y \circ \RX_* $ and counit  $\epsilon: \RX_* \circ \mathcal Y \to \id$ satisfying the triangle identities. 
 
 Let $G$ be a reductive group over $F$. By \S \ref{para:fundgp}, we have a canonical isomorphism $\RX_*(\mathscr Z_G) \cong \pi_1(G)$ in $\dercat{\GFMod}$. Let $\mathrm {pr}: \pi_1(G) \to \pi_1(G)_{\Gamma_F}$ be the canonical map, viewed as a morphism in $\dercat{\GFMod}$. Let \index[Oct2026-sorted-notation]{z 0 05 00 z subscript g superscript f@$\mathscr Z_G^F$} $$\mathscr Z_G^F : = \mathcal Y(\pi_1(G)_{\Gamma_F}) \in \dercat{\Mult(F)}, $$ and let $$\kappa_0: \mathscr Z_G \To \mathscr Z_G^F$$ be given by the composite
 $$  \mathscr Z_G \overset{\eta}{\To} \mathcal Y (\RX_*(\mathscr Z_G)) \xrightarrow{\mathcal Y(\mathrm {pr})}  \mathcal Y (\pi_1(G)_{\Gamma_F}) = \mathscr Z_G^F.$$ Thus we have a canonical isomorphism $$\epsilon : \RX_*(\mathscr Z_G^F) \isom \pi_1(G)_{\Gamma_F}.$$
 
 \begin{lem}\label{lem:rep by tori}
 	In $\dercat{\Mult(F)}$, $\mathscr Z_G^F$ is isomorphic to a complex of the form $$T^{-1} \To T^0, $$ where $T^{-1}$ and $T^0$ are split tori over $F$ and are located in degrees $-1 $ and $0$. In particular, the image of $\mathscr Z_G^F$ in $\cD(\site)$ lies in $\cD^{[-1,0]} (\site)$.
 \end{lem}
 \begin{proof}
 	This follows from the fact that $\pi_1(G)_{\Gamma_F}$ is isomorphic to a complex $L^{-1} \to L^0$ in $\dercat{\GFMod}$, where $L^{-1}$ and $L^0$ are finite free $\ZZ$-modules with the trivial $\Gamma$-action.
 \end{proof}
 \subsection{}\label{para:Kott stack}
 We write $\mathcal K^F_G$\index[Oct2026-sorted-notation]{k 0 03 00 k subscript g superscript f@$\mathcal K_G^F$} for $\ch(\mathscr Z_G^F) \in \Ch(\site)$, and call it the\index{Kottwitz stack} \emph{Kottwitz stack} of $G$ over $F$. The inclusions induce a canonical equivalence between quotient stacks $[Z_{G_{\sconn}} \backslash Z_G] \to [G_{\sconn} \backslash G]$. Thus we obtain a functor between stacks on $\site$:\index[Oct2026-sorted-notation]{zz10 kappa 1 00 00 kappa subscript g superscript can@$\kappa_G^{\can}$}
 $$\kappa_G^{\can}: G \To [G_{\sconn} \backslash G] \isom [Z_{G_{\sconn}} \backslash Z_G] \cong \ch(\mathscr Z_G) \xrightarrow{\ch(\kappa_0)} \ch(\mathscr Z_G^F)  = \mathcal K^F_G,$$ which is canonical up to isomorphism.
 
 For any small Picard category $P$ (strictly commutative, as always), we denote by $\pi_0(P)$ the set of isomorphism classes of $P$, which is naturally an abelian group. Then $\kappa_G^{\can}$ induces a morphism 
 $$ \pi_0(\kappa_G^{\can}) : G \To \pi_0( \mathcal K^F_G (\cdot))$$ of presheaves in groups\footnote{We caution the reader that the right hand side is not a sheaf.} on $\site$. 
 
 \subsection{}\label{para:setting for gamma} Now choose $T^{\bullet} = (T^{-1}\to T^0)$ as in Lemma \ref{lem:rep by tori}, and choose an isomorphism $f : T^{\bullet} \isom \mathscr Z_G^F$ in $\dercat{\Mult(F)}$.  If $R$ is any $F$-algebra with $\Pic(R) = 0 $, then we have canonical isomorphisms of abelian groups 
 \begin{multline}\label{eq:ch to coh0}
 	\pi_0 \big( \ch(T^\bullet) (R) \big ) \cong  T^0(R)/ T^{-1}(R) \\ \cong  (X_*(T^0)/X_*(T^{-1})) \otimes_{\ZZ} R^\times  \cong \coh^0(\RX_*(T^{\bullet})) \otimes_{\ZZ} R^\times , \end{multline}  since $T^{-1}$ is a split torus.  In this case, consider the composite isomorphism:\index[Oct2026-sorted-notation]{zz03 gamma 1 00 00 gamma subscript r@$\gamma_R$}
 \begin{multline*}
 	\gamma_R:  \pi_0 (\mathcal K^F_G (R) )  \xrightarrow{f^{-1}}  \pi_0 \big( \ch(T^{\bullet}) (R) \big )  \cong \coh^0(\RX_*(T^{\bullet})) \otimes R^{\times } \\  \overset{f}{\To} \coh^0(\RX_* (\mathscr Z^F_G)) \otimes  R^{\times} \xrightarrow{\epsilon \otimes \id } \pi_1(G)_{\Gamma_F} \otimes R^\times.
 \end{multline*}
 Then $\gamma_R$ is independent of the choice of $(T^{\bullet}, f)$, by the functoriality of (\ref{eq:ch to coh0}) in $T^{\bullet}$.

 \subsection{}\label{para:kottmapconst}
 If $R$ is a commutative ring, by a \emph{discrete valuation on $R$}, we mean a function $v: R \to \ZZ \cup \set{\infty}$ satisfying $v(0) = \infty$, $v(1) = 0$, $v(ab) = v(a) +v(b)$, and $v(a +b) \geq \min(v(a), v(b))$, for all $a,b\in R$. (Here $\infty >n, \forall n\in \ZZ$, and we do not require $v(a) =\infty \Rightarrow a=0$.)
 
 Now consider an $F$-algebra $R$ satisfying $\Pic(R) = 0$, and a discrete valuation $v$ on $R$. Composing the canonical map $\pi_0(\kappa_G^{\can})$ in \S \ref{para:Kott stack} with the canonical map $\gamma_R$ in \S \ref{para:setting for gamma}, we obtain the canonical map \index[Oct2026-sorted-notation]{zz10 kappa 1 00 00 kappa subscript g superscript r@$\kappa_G^R$}
 $$ \kappa_G^R:   G(R) \xrightarrow{\pi_0(\kappa_G^{\can})}   \pi_0 (\mathcal K_G^F (R) )  \overset{\gamma_R} {\To} \pi_1(G)_{\Gamma_F} \otimes R^{\times}  ,$$ which is a group homomorphism.
 On composing the above with $v: R^{\times} \to \ZZ$, we obtain the group homomorphism 
 \begin{align}\label{eq:kappa^R,v}
 	\kappa_G^{R,v} : G(R) \xrightarrow{\kappa_G^R}  \pi_1(G)_{\Gamma_F} \otimes R^{\times}  \overset{v}{\To}  \pi_1(G)_{\Gamma_F}  .
 \end{align}
 We now extend the definition of $\kappa_G^{R,v}$, dropping the hypothesis $\Pic(R) =  0 $.

 \begin{defn}\label{defn:Kott map} Let $R$ be an arbitrary $F$-algebra, and let $v$ be a discrete valuation on $R$. The elements $r \in R$ with $v(r) = \infty$
 	form a prime ideal $\mathfrak p$, and $v$ factors as $R \to \Frac(R/\mathfrak p) \xrightarrow{\bar v} \ZZ \cup \set{\infty}$. We define $\kappa_G^{R,v}$\index[Oct2026-sorted-notation]{zz10 kappa 1 00 00 kappa subscript g superscript r v@$\kappa_G^{R,v}$} to be the composition $$  G(R) \To G(\Frac (R/\mathfrak p)) \xrightarrow{\kappa_G^{\Frac(R/\mathfrak p), \bar v} }  \pi_1(G)_{\Gamma_F}  .$$ We call $\kappa_G^{R,v}$ the\index{Kottwitz homomorphism} \emph{Kottwitz homomorphism}. We often simply write $\kappa_G^v$ for $\kappa_G^{R,v}$.
 \end{defn}
 \subsection{}\label{para:functoriality of Kott map}
 In Definition \ref{defn:Kott map}, if $\Pic(R) = 0 $, then  one checks that the definition of $\kappa_G^{R,v}$ agrees with the previous definition (\ref{eq:kappa^R,v}). Moreover, the generally defined $\kappa_G^{R,v}$   is functorial in the pair $(R,v)$ in the following sense: Let $R'$ be another $F$-algebra equipped with a discrete valuation $v'$. Suppose there is an $F$-algebra map $h: R\to R'$ such that $v$ is the pull-back of $v'$ along $h$. Then $\kappa_G^{R,v}$ equals the composition $$ G(R) \overset{h}{\To} G(R') \xrightarrow{\kappa _G^{R',v'}}   \pi_1(G)_{\Gamma_F}  . $$
 
 For fixed $(R,v)$, the homomorphism $\kappa_G^{R,v}$ is also functorial in the reductive group $G$ over $F$. When $(F,v)$ is a complete discretely valued field with algebraically closed residue field and $v$ is normalized, one easily checks that $\kappa_G^{F,v}$ agrees with Kottwitz's original construction in \cite[\S 7]{kottwitzisocrystal2} by checking that the former satisfies Kottwitz's unique characterization in terms of functoriality and behavior on tori. 
 
 \begin{prop}\label{prop:sumkottmaps} Let $R$ be an $F$-algebra, and let $v_1, \dots, v_n$ be a collection of discrete valuations on $R$. Let $a_1,\dots, a_n \in \mathbb Z$. Assume the following conditions.
 	\begin{enumerate}
 		\item $\sum_{i=1}^n a_i v_i(u) = 0$ for all $u\in R^{\times}$. 
 		\item $\Pic(R) = 0 $.
 	\end{enumerate}
 	Then $\sum_{i=1}^n a_i\kappa_G^{v_i}(g) = 0$ for all  $g \in G(R).$
 \end{prop}
 \begin{proof} By condition (ii), each $\kappa_G^{v_i}$ factors as in (\ref{eq:kappa^R,v}). Thus the map $\sum_i  a_i \kappa_G^{v_i}$ factors through the map $\id \otimes \sum a_i v_i : \pi_1(G)_{\Gamma_F} \otimes R^{\times} \to  \pi_1(G)_{\Gamma_F}$, which is zero by (i).
 \end{proof}
 
 \begin{prop}\label{prop:intkottmaps} Let $R \supset F$ be a domain, and $v_1, \ldots , v_n$
 	a collection of discrete valuations on $R$. Let $a_1,\dots, a_n \in \mathbb Z$. Suppose that $R = R^\circ[1/f_j]_{j=1}^m,$ where $R^\circ \subset R$ is a subring and $f_j \in R^{\circ}$ are non-zero prime elements, satisfying the following conditions.
 	\begin{enumerate}
 		\item The ring $R^{\circ}$ is a noetherian locally factorial domain.
 		\item For $i=1, \dots, n$, we have $v_i(R^\circ) \subset \ZZ_{\geq 0} \cup \set{\infty}$. 
 		\item  For $j=1,\ldots, m$, we have $\sum_{i=1}^n a_iv_i(f_j) = 0$. 	\end{enumerate}
 	Then $\sum_{i=1}^n a_i\kappa^{v_i}_G(h) = 0$ for all $h \in G(R)$.
 \end{prop}
 \begin{proof} Let $R^{\circ\prime}$ be the ring obtained from $R^{\circ}$ by inverting all elements $f$
 	such that $v_i(f) = 0$ for all $i.$ The conditions of the proposition continue to hold if we replace $R^{\circ}$ (resp.~$R$) by $R^{\circ\prime}$ (resp.~$R^{\circ\prime}[1/f_j]_{j=1}^m$), and omit from the list of $f_j$ those elements such that $v_i(f_j) = 0$ for all $i$ (as they become units in $R^{\circ\prime}$). By the functoriality of the Kottwitz map as discussed in \S \ref{para:functoriality of Kott map}, we reduce to the case where $R^{\circ}= R^{\circ\prime}.$ Then $R^{\circ}$ is semi-local, as each proper ideal is contained in one of the prime ideals $\mathfrak p_i = \set{x\in R^{\circ}\mid v_i(x) >0}, i = 1,\dots,n$. 
 	
 	Since $R^{\circ}$ is noetherian and locally factorial, the natural map $\Pic( R^{\circ}) \to \Pic( R)$ is surjective (see \cite[Cor.~21.6.11]{EGA4-4}). Since $R^{\circ}$ is semi-local, we have $\Pic (R^{\circ}) = 0 $, and so $ \Pic(R) = 0 $.
 	
 	Now since the $f_j$ are prime in $R^{\circ}$, any unit in $R$ has the form $u = w f_1^{e_1}\dots f_m^{e_m}$ where
 	$w \in R^{\circ\times}, e_j \in \ZZ.$ Since $v_i(w) = 0$ for all $i$, we have
 	$$ \sum_{i=1}^n a_iv_i(u) =  \sum_{j=1}^m e_j \sum_{i=1}^n a_iv_i(f_j) = 0.$$
 	The proposition now follows from Proposition \ref{prop:sumkottmaps}.
 \end{proof}
 
 \begin{cor}\label{cor:kottparahoric} Suppose that $F$ is equipped with a surjective discrete valuation $v_F: F \rightarrow \mathbb Z\cup\{\infty\}$ with ring of integers $\O_F,$ and that $\cG$ is a smooth affine group scheme over $\O_F$ extending $G$ with connected special 
 	fiber. Then $\kappa_G^{v_F}: G(F) \to \pi_1(G)_{\Gamma_F}$ is zero on $\cG(\O_F)  \subset G(F)$. 
 \end{cor}
 \begin{proof}
 	Let $R^{\circ}$ be the affine ring of $\cG$. Let $\pi_F \in \oo_F$ be a uniformizer.
 	Our conditions imply that $\pi_F$ is a prime element of the noetherian domain $R^{\circ},$ and hence $R^{\circ}_{(\pi_F)}$ is a DVR. Let $v_0$ be the pull-back to $R^{\circ}$ of the canonical normalized discrete valuation on $R^{\circ}_{(\pi_F)}$. For each $g \in \cG(\O_F)$, let $v_{g}$ be the valuation 
 	$$R^{\circ} \overset{g}{\To} \O_F \xrightarrow{v_F} \ZZ \cup \set{\infty}. $$ 	Note that $v_0(\pi_F) = v_{g}(\pi_F) = 1$. Hence $v_0$ and $v_g$ extend to $R: = R^{\circ}[1/\pi_F]$.  By Proposition \ref{prop:intkottmaps} applied to $R^{\circ}$ and the relation   $v_0(\pi_F) - v_g(\pi_F) = 0$, we have 
 	$\kappa^{v_0}_G(h) - \kappa_G^{v_{g}}(h) =0 $ for all $h \in G(R)$.
 	Applying this to the universal point $g_u \in \cG(R^{\circ}) \subset G(R)$, we get
 	$$ \kappa^{v_0}_G(g_u) = \kappa_G^{v_{g}}(g_u) = \kappa_G^{v_F}(g),$$ where the second equality follows from functoriality (\S \ref{para:functoriality of Kott map}).
 	It follows that $\kappa_G^{v_F}$ is constant, and therefore zero, on $\cG(\O_F)$.  
 \end{proof}
 \begin{rem}
 	In Corollary \ref{cor:kottparahoric}, if the discretely valued field $(F,v_F)$ is strictly henselian and $\cG$ is a parahoric group scheme, then the conclusion follows from \cite[Prop.~3]{HainesRapoport}.
 \end{rem}
 \subsection{}\label{para:setting for Cartan decomp}
 Keep the setting and notation of Corollary \ref{cor:kottparahoric}. Assume that $\cG$ is reductive,
 and $F$ is complete. Let $\pi _F \in F$ be a uniformizer. Let $\cS \subset \cG$ be a maximal split torus. Then we have the Cartan--Iwahori--Matsumoto decomposition
 $$ G(F) = \bigcup_{\mu \in X_*(\cS)} \cG(\O_F)\mu(\pi_F)\cG(\O_F).$$ (The union is not disjoint.) The decomposition in this generality is proved in \cite[Thm.~1.3, Rmk.~3.5]{CIMdecomp}.

 \begin{cor}\label{cor:kottCartan} 	If $g \in G(F)$ belongs to the double coset indexed by $\mu \in X_*(\cS),$ then $\kappa_G^{v_F}(g) = [\mu]$, where $[\mu]$ is the image of $\mu$ under the natural map $X_*(\cS) = \pi_1(\cS_{F}) \to \pi_1(G) \to \pi_1(G)_{\Gamma_F}$.
 \end{cor}
 \begin{proof} By Corollary \ref{cor:kottparahoric}, it suffices to show that $\kappa_G^{v_F}(\mu(\pi_F)) = [\mu].$ But this follows from the functoriality of the Kottwitz map in the group $G$. \end{proof}

 \section{Decent elements and twisting}\label{subsec:decent}
 Throughout this section, we fix a prime $p$, and denote by $\sigma$ the arithmetic $p$-Frobenius in $\Aut(\LL)$. 	Let $G$ be a reductive  group over $\QQ_p$.  We recall various definitions and constructions (especially a \emph{twisting construction}) concerning $\sigma$-conjugacy classes in $G(\breve \QQ_p)$. 
 
 \subsection{}\label{isocdefns}
 For $b\in G(\breve \Q_p)$, we write $\nu_b$\index[Oct2026-sorted-notation]{zz13 nu 1 00 00 nu subscript b@$\nu_b$}  for the\index{Newton cocharacter} Newton cocharacter of $b$, which is a fractional cocharacter of $G_{\LL}$; see \cite[\S 4]{kottwitzisocrystal} and \cite[\S 1.7]{RZ96}. Following \cite[Def.~1.8]{RZ96}, we say that $b$ is\index{decent} \emph{decent}, if there exists $n \in \ZZ_{\geq 1}$ such that $n\nu_b$ is a cocharacter of $G_{\LL}$ and  
 \begin{align}\label{eq:n-decency}
 	b \sigma(b ) \cdots \sigma^{n-1} (b) = (n \nu_b) (p).
 \end{align}
 In this case, we also say that $b$ is \emph{$n$-decent}.
 
 For any $n$-decent $b \in G(\LL)$, it is shown in \cite[Cor.~1.9]{RZ96} that $b \in G(\QQ_{p^n})$, and that $\nu_b$ is defined over $\QQ_{p^n}$. In particular, if $b$ is $n$-decent, then it is also $n'$-decent for $n|n'$. Clearly the condition that an element of $G(\Q_{p^n})$ is $n$-decent is invariant under $\sigma$-conjugation by $G(\Q_{p^n})$. Conversely, if $b,b'$ are $n$-decent and if $b' = g b \sigma(g)^{-1}$ for some $g\in G(\breve \Q_p)$, then necessarily $g\in G(\Q_{p^n})$; see \cite[Cor.~1.10]{RZ96}.

 \subsection{}\label{para:B(G)}
 We denote by $\B(G)$\index[Oct2026-sorted-notation]{b 0 01 00 b g@$\mathrm B(G)$} the set of $\sigma$-conjugacy classes in $G(\LL)$. For $b \in G(\LL)$, we denote its class in $\B(G)$ by $[b]$. We recall Kottwitz's classification of elements of $\B(G)$. Let $\Newton(G)$\index[Oct2026-sorted-notation]{n 0 03 00 n g@$\mathcal N(G)$} denote the set of $\sigma$-stable $G(\LL)$-conjugacy classes of fractional cocharacters of $G_{\LL}$. The association $b\mapsto \nu_b$ descends to the\index{Newton map} \emph{Newton map} \index[Oct2026-sorted-notation]{zz13 nu 1 00 20 bar nu@$\bar \nu$}
 $\bar \nu : \B(G) \to \Newton (G)$. As in \S\ref{subsec:Kottwitz hom}, we have the Kottwitz homomorphism \index[Oct2026-sorted-notation]{w 1 00 00 w subscript g@$w_G$} $w_G: G(\breve \Q_p) \to \pi_1(G)_{\Gamma_{p,0}}$ associated with the $p$-adic valuation on $\breve \Q_p$. By \cite[\S 7]{kottwitzisocrystal2}, $w_G$ is surjective, and descends to a map\index[Oct2026-sorted-notation]{zz10 kappa 1 00 00 kappa subscript g@$\kappa_G$} $\kappa_G: \B(G) \to \pi_1(G)_{\Gamma_p} ,$ called the\index{Kottwitz map} \emph{Kottwitz map}. By \cite[\S 4.13]{kottwitzisocrystal2}, the map $$(\bar \nu, \kappa_G): \B(G) \To \Newton(G) \times \pi_1(G)_{\Gamma_p}$$ is injective.
 
 It is proved by Kottwitz \cite[\S 4.3]{kottwitzisocrystal} (cf.~\cite[\S 1.11]{RZ96}) that every $\sigma$-conjugacy class in $G(\LL)$ is represented by a decent element. Thus $\B(G)$ is in natural bijection with the set of $G(\Q_p^{\ur})$-orbits in the set of decent elements of $G(\Q_p^{\ur})$, where $G(\Q_p^{\ur})$ acts by $\sigma$-conjugation.
 
 \subsection{}\label{para:J_b}
 Let $b \in G(\LL)$. The functor sending any $\QQ_p$-algebra $R$ to the group \index[Oct2026-sorted-notation]{j 0 00 00 j subscript b@$J_b$}
 $$ J_b (R) : = \set{g\in G(R\otimes _{\QQ_p} \LL) \mid gb = b \sigma (g)}$$ is represented by a reductive group $J_b$ over  $\QQ_p$. We shall also write $J_b^G$ for $J_b$, to make the presence of $G$ explicit. If $b$ is decent (so $b \in G(\Q_p^{\ur})$), then by \cite[Cor.~1.14]{RZ96}, there is a canonical $\Q_p^{\ur}$-isomorphism from $J_{b, \Q_p^{\ur}}$ to the centralizer $G_{\Q_p^{\ur}, \nu_b}$ of $\nu_b$ in $G_{\Q_p^{\ur}}$. (In this case, for any $\Qpur$-algebra $R$ we have $J_b(R)\subset G(R \otimes_{\Q_p} \Q_p^{\ur})$, and the embedding $J_{b,\Qpur} \to G_{\Qpur}$ is induced by the natural map $R \otimes_{\Q_p} \Q_p^{\ur} \to R$.) 
 In this case the action of $\sigma$ on $G_{\Q_p^{\ur}, \nu_b}(\Q_p^{\ur})$ with respect to the $\Q_p$-form $J_b$ is given by $g \mapsto b \sigma (g) b^{-1}$, where $\sigma(\cdot)$ is defined with respect to the $\Q_p$-form $G$.
 
 If $b$ is decent and if 
 \begin{align}\label{eq:W non-empty}
 	\set{c \in  G(\Q_p^{\ur}) \mid c \nu_b c^{-1}\text{ is defined over } \Q_p} \neq \emptyset, 
 \end{align} then the canonical embedding $J_{b , \Q_p^{\ur}} \to G_{\Q_p^{\ur}}$ is a pre-inner twisting between $\QQ_p$-groups  (Definition \ref{defn:inner transfer}), as for any $c$ in the above set $\Int(c): J_{b, \Qpur} \to \Int(c)(J_{b,\Qpur})$ is an inner twisting between $\Q_p$-groups. In particular, we  obtain a canonical map
 \begin{align}\label{eq:from coh J_b}
 	\coh^1(\Q_p, J_b) \cong \coh^1_{\ab}(\Q_p, J_b) \to \coh^1_{\ab}(\Q_p, G) \cong \coh^1(\Q_p, G) 
 \end{align}
 as discussed in \S \ref{sec:inner twisting}. 
 Note that for $b$ decent, (\ref{eq:W non-empty}) holds in the following two cases: 
 \begin{enumerate}
 	\item $G$ is quasi-split over $\Q_p$.
 	\item $b$ is basic in $G$, i.e., $\nu_b$ is central.
 \end{enumerate}
 Indeed, the $G(\Q_p^{\ur})$-conjugacy class of $\nu_b$ is always stable under $\Gal(\Q_p^{\ur}/\Q_p)$. In case (i), this conjugacy class must contain a member defined over $\Q_p$. In case (ii), $\nu_b$ itself is already defined over $\QQ_p$.
 
 In case (ii), the pre-inner twisting is an inner twisting, so $J_b$ has the structure of an inner form of $G$. The isomorphism class of this inner form of $G$ (see Definition \ref{defn:inner form}) depends only on $[b] \in \B(G)$, not on the decent representative $b$.

 \subsection{}\label{para:setting for delta(b,beta)} Let $G$ be a reductive  group over $\Q_p$. Let $b \in G(\LL)$ be a decent element, and fix an element $\beta \in \coh^1(\Q_p, J_b)$. By Steinberg's theorem, $\beta$ is represented by a cocycle $ (a_{\rho})_{\rho} \in Z^1(\Q_p^{\ur}/\Q_p, J_b(\Q_p^{\ur}))$. Under the canonical isomorphism $J_{b, \Q_p^{\ur}} \isom G_{\Q_p^{\ur}, \nu_b}$, we view $a_{\sigma} \in J_b(\Q_p^{\ur})$ as an element of $G(\Q_p^{\ur})$, and define $b' : = a_{\sigma} b \in G(\Q_p^{\ur})$. It is easy to see that the $\sigma$-conjugacy class of $b'$ in $G(\Q_{p}^{\ur})$ depends only on $b, \beta$, not on $(a_{\rho})_{\rho}$. We shall say that (the $\sigma$-conjugacy class of) $b'$ is the twist of $b$ by $\beta$.  
 \begin{prop}\label{prop:twist p}
 	In the setting of \S \ref{para:setting for delta(b,beta)}, we have $\nu_{b'} = \nu_b$, and $b'$ is decent. Moreover, if (\ref{eq:W non-empty}) holds for $b$ (e.g., if $G$ is quasi-split), then $\kappa_G([b']) - \kappa _G([b])$ is equal to the image of $\beta$ under
 	$$ \coh^1(\Q_p, J_b) \To  \coh^1_{\ab}(\Q_p, G) \isom \pi_1(G)_{\Gamma_p, \tors}, $$ where the first map is (\ref{eq:from coh J_b}), and the second isomorphism is as in Proposition \ref{prop:apply TN}.
 \end{prop}
 \begin{proof}Choose $n\in\NN$ to be divisible enough such that $$(a_{\rho})_{\rho} \in Z^1(\Q_{p^n}/\Q_p, J_b(\Q_{p^n})), $$ and such that $b$ is $n$-decent. Using that $(a_{\rho})_{\rho}$ is a cocycle, one shows by induction that
 	$$ a_{\sigma^i} = b' \sigma (b') \cdots \sigma ^{i-1} (b') \sigma ^{i-1}(b^{-1}) \cdots \sigma (b^{-1}) b^{-1} ,$$ for each $ i \in \ZZ_{\geq 1}$. Since $a_{\sigma^n} =1$, we have
 	$$ b' \sigma (b') \cdots \sigma ^{n-1} (b') = b \sigma (b) \cdots \sigma ^{n-1} (b). $$
 	Since $b$ is $n$-decent, the right hand side is equal to $p^{n\nu_b}$. By the characterization of $\nu_{b'}$ (see \cite[\S 4]{kottwitzisocrystal}), we conclude that $\nu_{b'} = \nu_b$, and that $b'$ is $n$-decent.
 	
 	We now prove the statement about $\kappa_G([b']) - \kappa_G([b])$. Since (\ref{eq:W non-empty}) holds, we can replace $b$ by a $\sigma$-conjugate in $G(\Q_{p}^{\ur})$ and assume that $\nu_b$ is defined over $\Q_p$. Then we can replace $G$ by $G_{\nu_b}$ and reduce to the case where $\nu_b$ is central. To finish the proof we only need to show that the image of $\beta$ under the composite isomorphism $\coh^1(\Q_p, J_b) \cong \coh^1_{\ab}(\Q_p, J_b) \cong \pi_1(J_b)_{\Gamma_p, \tors}$ is equal to $\kappa_{J_b} (a_{\sigma}).$ This follows from \cite[Rmk.~5.7]{kottwitzisocrystal} applied to $J_b$ (cf.~\cite[Rmk.~2.2 (iv)]{RV14}).	
 \end{proof}
 
 \section{Shimura varieties and their cohomology}\label{sec:SV}
 We first establish a general formalism that will allow us to discuss the \emph{Hecke action on the $\ell$-adic cohomology}; this will be applied both to the tower of canonical models of Shimura varieties, and to the tower of canonical integral models of Shimura varieties (of a fixed hyperspecial level at $p$). We then review the definition of the \emph{canonical model} of a Shimura variety and fix the sign convention in the reciprocity law. In order to clarify the condition for an algebraic representation to give rise to an $\ell$-adic local system, we give a general discussion of \emph{cuspidality} and the \emph{anti-cuspidal part} of a torus. Finally we define the \emph{Frobenius--Hecke trace} on the cohomology of the Shimura variety; the study of this quantity is the central objective of this book. 
 \begin{defn}\label{defn:adm action} Let $H$ be a locally profinite group admitting a countable neighborhood basis of the identity. Let $S$ be  a scheme equipped with a right $H$-action. We say that the action is\index{admissible!action} \emph{admissible} if there exists a non-empty class $\mathscr K$ of compact open subgroups of $H$ satisfying the following conditions. 
 	\begin{enumerate}
 		\item If $K \in \mathscr K$, then every open subgroup of $K$ belongs to $\mathscr K$. 
 		\item For $K \in \mathscr K$, the categorical quotient $S/K$ in the category of schemes exists and is locally noetherian.
 		\item For $K_1, K_2 \in \mathscr K$  such that $K_1 \subset K_2$, the natural map $S/K_1 \to S/K_2$ is finite \'etale.
 		\item The natural map  $S  \to \varprojlim_{K\in \mathscr K} S/K$ is an isomorphism. Here the right hand side is defined in the category of schemes since the transition maps are affine.  		
 	\end{enumerate}
 \end{defn}
 \subsection{}\label{para:sheaf and Hecke} Let $H$, $S$, and $\mathscr K$ be as in Definition \ref{defn:adm action}. For $K \in \mathscr K$, we write $S_K$ for $S/K$.
 Let $\ell$ be a prime number. The construction in \cite[\S III.2]{harristaylor} can be generalized to define $\ell$-adic sheaves on $S_K$ and the Hecke action on cohomology. We explain this in the following.

 Let $K \in \mathscr{K}$ and let $U$ be an open normal subgroup of $K$. The finite group $K/U$ acts on $S_U$ via $S_K$-automorphisms. We denote by\index[Oct2026-sorted-notation]{k 0 00 00 k u superscript natural@$(K/U)^{\natural}$} $(K/U)^{\natural}$ the maximal quotient of $K/U$ that acts faithfully on $S_U$. Because $S_U \to S_K$ is \'etale, the criterion of \cite[Exp.~V, Prop.~3.1]{SGA1} ensures that the action of $(K/U)^{\natural}$ on $S_U$ is admissible in the sense of \S1 of that expos\'e, and the quotient scheme constructed there agrees with  $S_K$, as the latter is the categorical quotient. We make the following crucial assumption on $S$:
 \begin{enumerate}
 	\item[($\dagger$)] For all $K,U$ as above and every connected component $C$ of $S_U$, the stabilizer of $C$ in $(K/U)^{\natural}$ acts faithfully on $C$.
 \end{enumerate}
 In view of \cite[Exp.~V, Cor.~2.4, Prop.~2.6]{SGA1}, condition~($\dagger$) and the \'etaleness of $S_U \to S_K$ imply that the admissible  $(K/U)^{\natural}$-action on $S_U$ has no inertia, and therefore $S_U \to S_K$ is a torsor under $(K/U)^{\natural}$.

 For each $K \in \mathscr K$, we define the profinite group\index[Oct2026-sorted-notation]{k 0 00 00 k superscript natural@$K^{\natural}$}
 \begin{align}\label{eq:profinite Galois}
 	K^{\natural} : = \varprojlim_{ U \mathrel{\unlhd} K \text{ open}} (K/U)^{\natural} .
 \end{align} When we need to record the dependence on $S$, we write $(K/U)^{\natural, S}$ and $K^{\natural, S}$  for $(K/U)^{\natural}$ and $K^{\natural}$, respectively. 
 
 By our assumption on $H$, there exists a decreasing sequence of open normal subgroups of $K$ $$\cdots \subset U_i \subset U_{i-1} \subset \cdots \subset U_1 \subset K$$ forming a neighborhood basis of $1$.   Hence we have a presentation 	\begin{align}\label{eq:Gal(S/S_K)}
 	K^{\natural} \cong \varprojlim_i (K/U_i)^{\natural}.
 \end{align}

 We are thus in a situation similar to the general setting at the beginning 
 of \cite[\S III.2]{harristaylor}, and essentially the same argument as there 
 produces a lisse $\Lambda$-sheaf on $S_K$ associated with any continuous 
 representation of $K^{\natural}$ on a finite free $\Lambda$-module (where $\Lambda = \ol{\QQ}_{\ell}$, or $\Lambda = L$ or $\mathcal{O}_L$ 
 for a finite extension $L/\QQ_{\ell}$). We take this opportunity to record an oversight 
 in that reference. In our notation, \textit{loc.~cit.}~attaches a lisse sheaf on 
 $S_K$ to a continuous representation of $\varprojlim_i \Gal(S_{U_i}/S_K)$. However, 
 owing to the possible disconnectedness of $S_{U_i}$ (an essential feature in applications 
 to Shimura varieties), the morphism $S_{U_i} \to S_K$ is generally not a 
 $\Gal(S_{U_i}/S_K)$-torsor, nor is there an intrinsic group of $S_K$-automorphisms 
 making it one. Consequently, the resulting lisse sheaf may have strictly smaller rank 
 than that of the representation, contrary to what is intended. The construction must 
 instead exploit the torsor structure: the argument of \textit{loc.~cit.}~goes through 
 verbatim upon replacing $\Gal(S_{U_i}/S_K)$ with $(K/U_i)^{\natural}$, which acts on 
 $S_{U_i} \to S_K$ as a torsor.
 
 In what follows we make the following assumption in addition to ($\dagger$), which will be satisfied in the application to Shimura varieties.
 \begin{enumerate}
 	\item[($\ddagger$)] For each $K\in \mathscr K$, there exists a subgroup $\alpha_K \subset K$ central in $H$ such that for any open normal subgroup $U \subset K$ the kernel of $K/U \to (K/U)^{\natural}$ is $\alpha_K U / U$. 
 \end{enumerate}  By this assumption and (\ref{eq:Gal(S/S_K)}), we have $ K^{\natural} \cong \varprojlim_i K / \alpha_K U_i. $ Since $K$ is profinite, we also have $K \cong  \varprojlim_i K/U_i$.  Hence the natural map $K\to K^{\natural}$ is surjective by the Mittag-Leffler condition. Moreover, the kernel is the closure of $\alpha_K$ in $K$. 
 
 Now let $\xi$ be a continuous representation of $H$ on a finite-dimensional $\ol \QQ_{\ell}$-vector space $W$ satisfying the following assumption:
 \begin{itemize}
 	\item For any $K\in \mathscr K$, the restriction $\xi|_K$ is trivial on $\ker(K\to K^\natural)$. 
 \end{itemize}  Then by the surjectivity of $K \to K^{\natural}$, we have an induced continuous representation of $K^{\natural}$ on $W$. As discussed above, this gives rise to a lisse $\ol \QQ_{\ell}$-sheaf on $S_K$, which we denote by $\mL_{\xi, K}$.
 
 Consider $K_1,K_2 \in \mathscr K$ and $g \in H$ such that $g^{-1} K_1 g \subset K_2$. The action of $g$ on $S$ induces a map $g: S_{K_1} \to S_{K_2}$. As on p.~96 of \cite{harristaylor}, associated with $g$ there is a canonical morphism
 $$ \overrightarrow{g} ^* : g^* \mL_{\xi, K_2} \To \mL_{\xi, K_1}$$
 between $\ol \QQ_{\ell}$-sheaves on $S_{K_1}$. 
 
 Now we assume in addition that $S$ is a scheme over a base scheme $B$, that $H$ acts on $S$ via $B$-automorphisms, and that $S_K$ is separated of finite type over $B$ for each $K \in \mathscr K$. Fix a closed point of $B$ with residue field $k$ such that $\ell$ is invertible on $k$, and fix an algebraic closure $\bar k$ of $k$. In this case, $\overrightarrow{g} ^*$ as above induces a map\index[Oct2026-sorted-notation]{g 1 00 10 g subscript k subscript 2 k subscript 1 leftbracket rightbracket@$[g]_{K_2, K_1}$}
 $$ [g]_{K_2, K_1}: \coh^i_{c} (S_{K_2, \bar k} , \mL_{\xi, K_2}) \To \coh^i_{c} (S_{K_1, \bar k} , \mL_{\xi, K_1}). $$
 Define\index[Oct2026-sorted-notation]{h 0 02 00 h superscript i subscript c s subscript bar k xi@$\mathbf H^i_c (S_{\bar k}, \xi)$} $$\coh^i_c (S_{\bar k}, \xi) : = \varinjlim_{K \in \mathscr K} \coh^i_c (S_{K,\bar k} , \mL_{\xi,K}) , $$ where the transition maps are given by $[1]_{K_2, K_1}$ for $K_1 \subset K_2$. The maps $[g]_{K_2, K_1}$ for varying $g, K_1, K_2$ give rise to a left $H$-action on $\coh^i_c(S_{\bar k}, \xi)$, called the \emph{Hecke action}. As on p.~97 of \cite{harristaylor}, $\coh^i_c(S_{\bar k}, \xi)$ is an admissible $H$-module over $\ol \QQ_{\ell}$. Indeed, using the Hochschild--Serre spectral sequence and the fact that the group cohomology of a finite group acting on a $\ol \QQ_{\ell}$-vector space vanishes in positive degrees, we see that for each $K \in \mathscr K$ the natural map $\coh^i_c(S_{K,\bar k}, \mL_{\xi,K}) \to \coh^i_c(S_{\bar k}, \xi)$ is injective and has image equal to the subspace of $K$-invariants. The admissibility then follows from the fact that $\coh^i_c(S_{K,\bar k}, \mL_{\xi,K})$ is finite-dimensional. We also have an action of $\Gal(\bar k / k)$ on $\coh^i_c(S_{\bar k}, \xi)$ commuting with the $H$-action.

 \subsection{}\label{para:Shimura variety}
 Let $(G,X)$ be a Shimura datum with\index{reflex field} reflex field $E= E(G,X) \subset \CC$,\index[Oct2026-sorted-notation]{e 0 00 00 e g x@$E(G,X)$} as in \cite[\S\S 2.1, 2.2]{deligne1979varietes}. By the theory of\index{canonical model} canonical models established by Shimura \cite{Shimura63,Shimura64,Shimura65,Shimura66,Shimura67b,Shimura67, Shimura70, Shimura70b}, Deligne \cite{deligne1971traveaux,deligne1979varietes}, Milne \cite{Mil83} (cf.~\cite{Mil88}), and Borovoi \cite{borovoi83}, we have a canonical $E$-scheme $\Sh = \Sh(G,X)$ equipped with a right $G(\A_f)$-action that is admissible in the sense of Definition \ref{defn:adm action}; the class $\mathscr K$ can be taken to be the class of neat compact open subgroups of $G(\A_f)$ (as defined in \cite[\S 0.6]{pink1989compactification}). For each $K \in \mathscr K$, we denote $\Sh(G,X)/K$ by $\Sh_K = \Sh_K(G,X)$.\index[Oct2026-sorted-notation]{sh 0 01 00 sh subscript k@$\mathrm{Sh}_K$} This is a smooth, quasi-projective $E$-scheme, whose analytification over $\CC$ is identified with the hermitian locally symmetric variety $G(\QQ) \backslash X \times G(\A_f) /K$.

 If $G=T$ is a torus, then $X$ consists of a single $\RR$-homomorphism $h: \mathbb S \to T_{\RR}$, and $\Sh_K(\overline E)$ is identified with the finite set $\Sh_K(\CC) = T(\QQ) \backslash T(\A_f)/ K$. The action of $\Gal(\overline E/E)$ on this finite set is given by the\index{reciprocity law} \emph{reciprocity law}, which we now review in order to fix the sign convention (cf.~\cite[\S 2.2.3]{deligne1979varietes}).
 Let\index[Oct2026-sorted-notation]{zz12 mu 1 00 00 mu subscript h@$\mu_h$} $\mu = \mu_h$, which is a cocharacter of $T$ defined over $E$.\footnote{The convention used here is the same as in \cite{deligne1979varietes}. If $T_{\RR}= \GG_{m, \RR}$ and $h$ is given by $h(z) = z^p \bar z^q$, then $\mu_h$ is given by $\mu_h(z) = z^{-p}$.} Consider the composite  homomorphism of $\Q$-algebraic groups\index[Oct2026-sorted-notation]{r 1 00 00 r mu superscript alg@$r(\mu)^{\alg}$}
 \begin{align*}
 	r(\mu)^{\alg}:  \Res_{E/\Q} \GG_m \xrightarrow {\Res_{E/\Q} \mu}   \Res_{E/\Q} T \xrightarrow{\N_{E/\Q}} & T  .
 \end{align*} This induces a group homomorphism
 $$\pi_0(E^{\times} \backslash \adele_E^{\times} ) \To \pi_0(T(\Q) \backslash T(\A)).$$ Now the left hand side is identified with $\Gal(E^{\ab}/E)$ under the global Artin map (normalized geometrically, i.e., uniformizers correspond to geometric Frobenius elements at the finite places), while the right hand side admits a natural map to $T(\Q) \backslash T(\A_f)/K$ (which is finite and discrete) induced by the projection $T(\adele) \to T(\adele_f)$.  We thus obtain a group homomorphism $$ r: \Gal(E^{\ab}/E)  \To T(\Q) \backslash T(\A_f)/K. $$ For $\sigma \in \Gal(\overline E/E)$ and $x \in \Sh_K(\overline E)\cong T(\QQ)\backslash T(\A_f)/K$, we have the\index{reciprocity law} \emph{reciprocity law} $$\sigma (x) = r(\sigma) \cdot x. $$ This uniquely determines the $E$-scheme structure of $\Sh_K$. Note that the above reciprocity law differs from \cite{deligne1979varietes} in the sign of $\mu$. Thus the $E$-scheme which we call the canonical model for $(T,\set{h})$ would be called the canonical model for $(T, \set{h^{-1}})$ according to \textit{loc.~cit.} Our sign convention is the same as that used by Pink \cite{pink1989compactification, pink1992ladic} and Morel \cite{Mor10}. 
 
 For general Shimura data, the canonical models are uniquely characterized by functoriality and the case of tori. According to our sign convention, the Siegel modular varieties which are moduli of polarized abelian varieties are canonical models for the Siegel Shimura data as specified in \cite[\S 2.1.5]{kisin2010integral}; see \cite[Thm.~4.21]{deligne1971traveaux}.\footnote{There is a discrepancy between this fact and Deligne's sign conventions in \cite{deligne1979varietes}, which was observed in \cite{milnelettertodeligne}.}
 
 We next discuss the standard construction of $\ell$-adic sheaves on Shimura varieties. For this we need some preparations. 
 
 \begin{defn}\label{defn:anti-cusp}
 	Let $T$ be a $\QQ$-torus. We denote by $T_a$\index[Oct2026-sorted-notation]{t 0 00 00 t subscript a@$T_a$} the maximal $\QQ$-anisotropic subtorus of $T$ (see 
 	\cite[Prop.~13.2.4]{Springerbook}). We denote by $T_{ac}$\index[Oct2026-sorted-notation]{t 0 00 00 t subscript ac@$T_{ac}$} the smallest $\QQ$-subgroup of $T_a$ whose base change to $\RR$ contains the maximal $\RR$-split subtorus of $T_{a,\RR}$. (Clearly $T_{ac}$ exists and is a torus.) We call $T_{ac}$ the\index{anti-cuspidal part} \emph{anti-cuspidal part} of $T$.  We say that $T$ is\index{cuspidal torus} \emph{cuspidal}, if $T$ has equal $\QQ$-split rank and $\RR$-split rank (cf.~\cite[Def.~3.1.1]{Mor10}). 
 \end{defn} 
 \begin{lem}\label{lem:cusp TFAE}
 	Let $T$ be a $\QQ$-torus. The following statements are equivalent.\begin{enumerate}
 		\item $T$ is cuspidal.
 		
 		\item $T$ is isogenous over $\QQ$ to the product of a split $\Q$-torus and a $\Q$-torus that is anisotropic over $\R$.
 		
 		\item $T_a$ is $\RR$-anisotropic. 
 		\item $T_{ac}$ is trivial. 
 		\item $T(\QQ)$ is discrete in $T(\A_f)$.
 		\item All arithmetic subgroups of $T(\QQ)$ are finite. 
 		\item $T$ satisfies the following\index{Serre condition} \emph{Serre condition} (cf.~\cite[\S 3.5.6]{kisin2012modp}). Fix a complex conjugation $\iota \in \Gal(\Qbar/\Q)$. For all $\tau \in \Gal(\Qbar/\Q ) $ and all $\mu \in X_*(T)$, we have 
 		$$ (\tau - 1) (\iota+1) \mu = (\iota+1)(\tau-1)\mu = 0.$$
 	\end{enumerate}
 	In general, $T_{ac}$ is the smallest $\QQ$-subgroup $S$ of $T$ such that $T/S$ is cuspidal. 
 \end{lem}
 \begin{proof}
 	The equivalence of (i), (ii), (iii), and (iv) follows  from \cite[Prop.~13.2.4]{Springerbook}. The equivalence of (i), (v), and (vi) is shown in  \cite[Prop.~1.4]{Gro99}. We now show that (vii) is equivalent to the other conditions. Note that (vii) is invariant under isogeny over $\QQ$, and is satisfied when $T$ is either split over $\QQ$ or anisotropic over $\RR$. Hence (ii) implies (vii). Conversely, if (vii) holds, then every $\mu \in X_*(T)^{\iota = 1}$ is fixed by $\Gal(\Qbar/\Q)$, since $(\iota+1)\mu = 2 \mu$ is fixed by $\Gal(\Qbar/\Q)$. This implies (i).
 	
 	The last assertion in the lemma is clear since (i) and (iii) are equivalent.  
 \end{proof} 
 \begin{lem}\label{lem:Z_s}
 	Let $G$ be a reductive group over $\QQ$, and $K\subset G(\A_f)$ a neat compact open subgroup. Then the intersection $Z_G(\Q) \cap K$ in $G(\adele_f)$ is contained in $(Z_G^0)_{ac}(\Q)$.
 \end{lem}
 \begin{proof}
 	By Lemma \ref{lem:cusp TFAE},  all congruence subgroups of $(Z_G^0/(Z_G^0)_{ac})(\Q)$ are finite. It follows that the image of $Z_G(\QQ)\cap K$ in $(Z_G/(Z_G^0)_{ac})(\Q)$ is finite. But this image is also neat inside $(Z_G/(Z_G^0) _{ac})(\A_f)$, so it is trivial.
 \end{proof}
 
 The following result will not be used until \S \ref{sec:Shimura Hodge type}. Recall that a Shimura datum is of \emph{Hodge type}\index{Shimura datum!Hodge type} if it admits an embedding into a Siegel Shimura datum. 
 \begin{lem}\label{lem:G=G^c}
 	For any Shimura datum $(G,X)$ of Hodge type, the $\QQ$-tori $Z_G^0$ and $G^{\ab}$ are cuspidal, and in particular $(Z_G^0)_{ac}$ is trivial. If $T$ is a maximal torus in $G$ such that some $h\in X$ factors through $T_{\RR}$, then $T$ is cuspidal and $T_{\RR}$ is an elliptic maximal torus in $G_{\RR}$. 
 \end{lem}
 \begin{proof}
 	Since $(G,X)$ is of Hodge type, its weight cocharacter $w$ is defined over $\QQ$, and $\Int(h(\sqrt{-1}))$ induces a Cartan involution on $G_{\RR}/w(\GG_{m,\RR})$ for each $h \in X$ (by directly checking the same properties for the Siegel Shimura datum). Recall from \cite[\S 2.1.1]{deligne1979varietes} that $w$ is central. Thus $(Z_G^0)_{\RR}/ w(\GG_{m,\RR})$ is anisotropic. Since $w(\GG_{m,\RR})$ is defined and split over $\QQ$, we see that $Z_G^0$ is cuspidal. Equivalently (Lemma \ref{lem:cusp TFAE}), $(Z_G^0)_{ac}$ is trivial. Since $G^{\ab}$ is isogenous to $Z_G^0$ over $\QQ$, it is also cuspidal. If $h\in X$ factors through $T_{\RR}$, then the Cartan involution $\Int(h(\sqrt{-1}))$ on $G_{\RR}/w(\GG_{m,\RR})$ restricts to the identity on $T_{\RR}/w(\GG_{m,\RR})$, forcing the latter to be anisotropic. This implies that $T$ is cuspidal and that $T_{\RR}$ is elliptic in $G_{\RR}$. 
 \end{proof}
 
 \subsection{}\label{para:aut sheaf} 
 Let $(G,X)$ be a Shimura datum with reflex field $E$. Fix a prime number $\ell$, and let $\xi$ be an irreducible algebraic representation of $G$ over $\ol \Q_{\ell}$. 
 We set $G^c : = G/ (Z_G^0)_{ac}$,\index[Oct2026-sorted-notation]{g 0 00 00 g superscript c@$G^c$} and assume that $\xi$ factors through $G^c$.\footnote{  
 	In \cite[\S III]{Mil88}, $G^c$ is defined to be $G/T$ for $T$ the maximal $\QQ$-subtorus of $Z_G^0$ that is $\QQ$-anisotropic and $\RR$-split. Note that it is assumed in \cite[\S II, (2.1.4)]{Mil88} that $Z_G^0$ splits over a CM field. Under this assumption, we have $T= (Z_G^0)_{ac}$. In general, the two can be different. }   In the following we construct $\ol \Q_{\ell}$-sheaves on $\Sh_K$ by applying the general formalism in \S \ref{para:sheaf and Hecke} to the $G(\adele_f)$-action on $\Sh$. This construction is well known, cf.~for instance \cite[\S 6]{kottwitz1992points}, \cite[\S 5]{pink1992ladic}, \cite[\S III.2]{harristaylor}, and \cite[\S 3]{lanstrohII}, which work at different levels of generality. Note that in all but the last reference, the Shimura varieties being considered satisfy $G= G^c$. 
 
 Let $K$ be a neat compact open subgroup of $G(\A_f)$ and $U$ an open normal subgroup of $K$. If $\gamma \in G(\QQ)$ has a fixed point in $X$ and lies in a neat compact open subgroup of $G(\A_f)$, then the image of $\gamma$ in $G^{\ad}(\QQ)$ is neat and lies in a compact subgroup of $G^{\ad}(\mathbb A)$, and is therefore trivial. Using this it is easy to see the following:
 \begin{align}\label{eq:(K/U)'}
 	\{k\in K/U \mid k \text{ has a fixed point in }\Sh_U(\mathbb C)\} & = \ker (K/U \to \Aut (\Sh_U)) \\ \nonumber & = (Z_G(\QQ) \cap K) U/U.
 \end{align}  
 In particular, the assumptions   ($\dagger$) and ($\ddagger$) in \S \ref{para:sheaf and Hecke} are both satisfied, and $\alpha_K = Z_G(\QQ)\cap K$.\index[Oct2026-sorted-notation]{zz01 alpha 1 00 00 alpha subscript k@$\alpha_K$}  
 Write $Z$ for $Z_G$, write $Z(\Q)^-$ for the closure of $Z(\Q)$ in $Z(\A_f)$, and write\index[Oct2026-sorted-notation]{z 0 00 00 z q superscript minus subscript k@$Z(\mathbb Q)^-_K$} $Z(\Q)^-_K$ for the intersection $Z(\Q)^- \cap K$ inside $G(\A_f)$. Note that $Z(\Q)^-_K$ is also the closure of $Z(\Q)\cap K$ inside $K$, since $K$ is open and closed in $G(\A_f)$. Hence by the discussion in \S \ref{para:sheaf and Hecke} below ($\ddagger$), we have  (cf.~\cite[\S 2.1.9]{deligne1979varietes})
 \begin{align}\label{eq:computation of Galois group}
 	K^{\natural,\Sh} \cong K/Z(\Q)_K^-.
 \end{align}  By Lemma \ref{lem:Z_s}, the natural map $K \to G^c(\A_f)$ factors through $K^{\natural,\Sh}$.\footnote{In fact, the induced map $K^{\natural,\Sh} \cong K /Z(\Q)_K^- \to G^c(\A_f)$ is never injective, if $(Z_G^0)_{ac}$ is non-trivial. Write $T$ for $(Z_G^0)_{ac}$. Then this follows from the fact that $T(\Q)(K\cap T(\A_f))$ has finite index in $T(\A_f)$ (\cite[Thm.~5.1]{Bor63}), and the fact that the closure of  $T(\Q)$ in $T(\A_f)$ has infinite index in $T(\A_f)$ (\cite[Prop.~7.13(2)]{plantonov-rapinchuk}). In \cite[\S III, Rmk.~6.1]{Mil88} it is incorrectly asserted that $K^{\natural,\Sh}$ is isomorphic to the image of $K$ in $G^c(\A_f)$, cf.~\S 3 of the updated online version of \cite{lanstrohII} and its erratum.}
 
 Via the projection $G(\A_f) \to G(\QQ_{\ell})$, $\xi$ induces a continuous representation of $G(\A_f)$ on a finite-dimensional $\ol\QQ_{\ell}$-vector space. This representation satisfies the assumption in \S \ref{para:sheaf and Hecke}, namely that its restriction to $K$ factors through $K^{\natural,\Sh}$ for all neat  $K$. By the construction in \S \ref{para:sheaf and Hecke}, we obtain a lisse $\ol \QQ_{\ell}$-sheaf $\mL_{\xi,K}$ on $\Sh_K$, and obtain the admissible $G(\A_f)$-module\index[Oct2026-sorted-notation]{h 0 02 00 h superscript i subscript c sh subscript bar e xi@$\mathbf H^i_c(\mathrm{Sh}_{\ol E}, \xi)$}
 $$ \coh^i_c(\Sh_{\ol E}, \xi) : =\varinjlim_{K} \coh^i_c(\Sh_{K, \ol E}, \mL_{\xi, K})$$ equipped with an action of  $\Gal(\ol E/E)$ commuting with $G(\A_f)$. Our main interest is to understand the virtual $G(\A_f)\times \Gal(\ol E/E)$-module 
 $$ \sum_i (-1)^i \coh^i_c (\Sh_{\overline E}, \xi). $$ More specifically, we are interested in the trace of a power of Frobenius at a ``good prime'' composed with a Hecke operator away from that prime on this virtual module. To formulate this problem we need more notations.

 \subsection{}\label{para:unramified SD}
 Let $p$ be a prime number unequal to $\ell$. Recall that $E$ is by definition a subfield of $\CC$. From the fixed embeddings $\Qbar\hookrightarrow \CC$ and $\Qbar \hookrightarrow \Qpbar$, we obtain a prime $\mathfrak p$\index[Oct2026-sorted-notation]{p 1 04 00 p@$\mathfrak p$} of $E$ over $p$. Let $q= p^r$\index[Oct2026-sorted-notation]{q 1 00 00 q p superscript r@$q= p^r$} be the cardinality of the residue field of $\mathfrak p$. 
 
 Now we assume that $G$ is unramified over $\Q_p$, and fix a connected reductive model $\cG$ of $G_{\Q_p}$ over $\Z_p$. In the sequel we shall call such a quadruple $(G,X, p , \cG)$ an\index{unramified!Shimura datum} \emph{unramified Shimura datum}. In this case $E_{\fkp}$ is unramified over $\Q_p$, so it is identified with $\Q_{p^r}$. We write\index[Oct2026-sorted-notation]{k 0 00 00 k subscript p@$K_p$} $K_p$ for the hyperspecial subgroup $\cG(\Z_p) \subset G(\Q_p)$.

 We fix notations for Hecke algebras. 	For each compact open subgroup $K^p \subset G(\A_f^p)$,
 let 
 $\cH(G(\A_f^p){\sslash} K^p)$\index[Oct2026-sorted-notation]{h 0 03 00 h g a subscript f superscript p double quotient k superscript p@$\mathcal H(G(\A_f^p){\sslash}K^p)$}\index[Oct2026-sorted-notation]{h 0 03 00 h g a subscript f superscript p double quotient k superscript p subscript q@$\mathcal H(G(\A_f^p){\sslash}K^p)_{\QQ}$} be the Hecke algebra of $\CC$-valued smooth compactly
 supported $K^p$-bi-invariant distributions on $G(\A_f^p
 )$. Let $\cH(G(\A_f^p){\sslash} K^p)_{\QQ}$ be the $\QQ$-subalgebra of $\cH(G(\A_f^p){\sslash} K^p)$ consisting of distributions that are rational on characteristic functions of compact open subsets of $G(\A_f^p)$. Elements of $\cH(G(\A_f^p){\sslash} K^p)$  can be represented as $f^p dg^p$, where $f^p$ is a $\CC$-valued smooth compactly supported $K^p$-bi-invariant function on $G(\A_f^p)$, and $dg^p$ is a Haar measure on $G(\A_f^p)$ assigning rational volumes to compact open subgroups. Elements of $\cH(G(\A_f^p){\sslash} K^p)_{\QQ}$ can similarly be represented as $f^p dg^p$, where $dg^p$ is as before and $f^p$ is $\QQ$-valued. 
 
 We have an induced action of $\cH(G(\A_f^p)\sslash K^p)_{\QQ}$ on the admissible $G(\A_f^p)$-module $\coh^i_c(\Sh_{\ol E}, \xi)^{K_p}$. Fix a decomposition subgroup $D_{\fkp} \subset \Gal(\overline E/E)$ at $\fkp$, and fix an element $\Phi_{\fkp} \in D_{\fkp}$ that lifts the geometric $q$-Frobenius. 
 For $m \in \ZZ_{\geq 1}$ and $f^p dg^p\in \mathcal H(G(\A_f^p){\sslash} K^p)_{\QQ}$, we define\index[Oct2026-sorted-notation]{t 0 00 00 t phi subscript p superscript m f superscript pdg superscript p@$T(\Phi_{\fkp}^m, f^pdg^p)$} $$T(\Phi_{\fkp}^m, f^pdg^p): = \sum_{i} (-1)^i \tr \bigg({\Phi}_{\mathfrak p}^m \times (f^p dg^p)  \mid \coh_c^i(\Sh_{\ol{E}},\xi)^{K_p} \bigg)  \in  \ol \Q_{\ell}. $$
 Our goal in the rest of this chapter is to state a conjectural formula for the above quantity. The formula will be a summation over \emph{Kottwitz parameters}, under the condition that the \emph{Kottwitz invariant} vanishes. In \S \ref{coh Kott} and \S \ref{Kottinv} below we define Kottwitz parameters and the Kottwitz invariant respectively.

 \section{Kottwitz parameters}\label{coh Kott} 
 Let $(G,X)$ be a Shimura datum and $p$ a prime. We define \emph{Kottwitz parameters} with respect to $(G,X)$ and $p$, generalizing the considerations in \cite[\S 2]{Kot90} where $G_{\der}$ is  assumed to be simply connected. While we refer to Kottwitz's original triple $(\gamma_0,\gamma,\delta)$ as a \emph{classical Kottwitz parameter}, our version is a refinement when $G_{\der}$ is not simply connected. We have chosen the terminology ``Kottwitz parameter'' in order to distinguish it from the ``Kottwitz triples'' of \cite{kisin2012modp}. Roughly speaking, the latter (to be extensively used in \S \ref{sec:Shimura Hodge type})  describe $\Fpbar$-points of a Shimura variety, while Kottwitz parameters are designed to describe points over a finite field. 
 
 \subsection{}\label{para:fkp}
 Let $E$ be the reflex field of $(G,X)$, and let $\mathfrak p$ and $q= p^r$ be as at the beginning of \S \ref{para:unramified SD}. Fix a positive multiple $n$ of $r$. Unlike \S \ref{para:unramified SD}, we do not assume here that $G$ is unramified over $\QQ_p$. 
 
 The Hodge cocharacters attached to $h\in X$ all have the same image in $\pi_1(G)$, which we denote by\index[Oct2026-sorted-notation]{zz12 mu 1 00 10 mu subscript x leftbracket rightbracket@$[\mu]_X$} $[\mu]_X \in \pi_1(G)$. 
 
 We say that an element of $G(\QQ)$ is \emph{semi-simple and $\RR$-elliptic}, if it lies in $T(\RR)$ for some elliptic maximal torus $T$ in $G_{\RR}$. Since $G$ is part of a Shimura datum, $G^{\ad}_{\RR}$ admits an inner Cartan involution. Therefore $G_{\RR}$ contains elliptic maximal tori, and our definition of $\RR$-elliptic elements agrees with the more general definition in the literature, cf.~\cite[\S 9.1]{Kot86}.

 \begin{defn}\label{defn of Kott trip}
 	A\index{Kottwitz parameter!classical} \textit{classical Kottwitz parameter of degree $n$} (with respect to $(G,X)$ and $p$) is a triple $$(\gamma_0, \gamma = (\gamma_l)_{l\neq p} , \delta) \in G(\QQ)\times G(\mathbb A_f^p)\times G(\QQ_{p^n}), $$ satisfying the following conditions.
 	\begin{description}
 		\item[CKP1] The element $\gamma_0\in G(\QQ)$ is semi-simple and $\RR$-elliptic.
 		\item[CKP2] For each prime $l\neq p$, the element $\gamma_l$ is\index{stable conjugacy} stably conjugate to $\gamma_0$ as elements of $G(\QQ_l)$. (See \cite[\S 3]{Kot82} for the standard definition of stable conjugacy.)
 		\item[CKP3] The image of $\gamma_0$ in $G(\QQ_p) $ is a degree $n$ norm of $\delta $ (see \cite[\S 5]{Kot82}).
 		\item[CKP4] The image of $\delta $ under the Kottwitz map $\kappa_{G} : \B(G_{\Q_p}) \to \pi_1 (G)_{\Gamma_p} $ is equal to the image of $ -[\mu]_X$.
 	\end{description}
 	We denote by $\KP_{\clsc}(p^n)$\index[Oct2026-sorted-notation]{kp 0 04 00 kp subscript cla p superscript n@$\mathfrak{KP}_{\mathrm{cla}}(p^n)$} the set of classical Kottwitz parameters of degree $n$.
 \end{defn}
 \subsection{}
 \label{defn of equivalence of triples}
 We define an equivalence relation $\sim$ on $ G(\QQ)\times G(\mathbb A_f^p)\times G(\QQ_{p^n}) $  by declaring $(\gamma_0, \gamma, \delta) \sim ({\gamma_0}', \gamma' , \delta')$ if the following conditions are satisfied:
 \begin{itemize}
 	\item The elements $\gamma_0$ and ${\gamma_0}'$ are stably conjugate in $G(\Q)$.
 	\item The elements $\gamma$ and $\gamma'$ are conjugate in $G(\mathbb A_f^p)$.
 	\item The elements $\delta$ and $\delta'$ are $\sigma$-conjugate in $G(\QQ_{p^n})$.
 \end{itemize}
 The subset $\KP_{\clsc}(p^n) \subset  G(\QQ)\times G(\mathbb A_f^p)\times G(\QQ_{p^n})$ is stable under $\sim$.

 \begin{defn}\label{defn of coh Kott trip}
 	A\index{Kottwitz parameter} \textit{Kottwitz parameter} (with respect to $(G,X)$ and $p$)  is a tuple $(\gamma_0, a, [b])$, consisting of:
 	\begin{itemize}
 		\item a semi-simple and $\RR$-elliptic element $\gamma_0 \in G(\Q)$,
 		\item an element $a \in \D(G_{\gamma_0}^0, G ; \A_f^p  )$ (see \S \ref{para:two groups} for $\D$),
 		\item a $\sigma$-conjugacy class $[b]$ in $G_{\gamma_0}^0 (\LL)$ (i.e., $[b] \in \B(G_{\gamma_0, \Q_p}^0)$),
 	\end{itemize}
 	satisfying the following condition.
 	\begin{description}
 		\item[KP0] Let $[b]_G$ be the image of $[b]$ in $\B(G_{\Q_p})$. Then the element $\kappa_{G}([b]_G) \in \pi_1(G)_{\Gamma_p}$ is equal to the image of  $-[\mu]_X$.\end{description}  We denote by $\KP$\index[Oct2026-sorted-notation]{kp 0 04 00 kp@$\mathfrak {KP}$} the set of all Kottwitz parameters. 
 	For $\fkc = (\gamma_0, a, [b]) \in \KP$, we write $I_0(\fkc)$\index[Oct2026-sorted-notation]{i 0 00 00 i subscript 0 c@$I_0(\mathfrak c)$} for $G_{\gamma_0}^0$. When $\fkc$ is fixed in the context we also simply write $I_0$ for $I_0(\fkc)$.
 \end{defn}
 \begin{defn}\label{pn-adm}
 	We say that $(\gamma_0, a, [b])\in \KP$ is\index{admissible!Kottwitz parameter} \textit{$p^n$-admissible}, if it satisfies the following condition.
 	\begin{description}
 		\item [KP1] For one (and hence any) representative $b \in I_0(\LL)$ of $[b]$, there exists $c\in G(\LL)$ satisfying the following. Let $\delta := c^{-1} b \sigma(c) \in G(\LL)$. We have \begin{align} \label{eqn in KP1}
 			c^{-1} \gamma_0 c = \delta \cdot \sigma(\delta) \cdots \sigma^{n-1}(\delta).
 		\end{align}
 	\end{description}  
 	We denote by $\KP_{\adm}(p^n)$\index[Oct2026-sorted-notation]{kp 0 04 00 kp subscript a p superscript n@$\mathfrak{KP}_{\mathrm{a}}(p^n)$} the set of $p^n$-admissible Kottwitz parameters.
 \end{defn}
 
 \subsection{}\label{para:setting KP1}
 Next we deduce some consequences of the condition \textbf{KP1}. We shall work in a local setting as follows. Let $\gamma_0 \in G(\Q_p)_{\semi}$ and let $I_0 : = (G_{\Q_p})_{\gamma_0}^0$. Let $[b] \in \B(I_0)$, and assume that \textbf{KP1} holds for $[b]$ and $\gamma_0$.
 
 \begin{lem}\label{b decides delta} Keep the setting of \S \ref{para:setting KP1}. Choose $b$ and $c$ as in \textbf{KP1} with respect to $[b]$ and $\gamma_0$. Let $\delta = c^{-1} b \sigma(c)$. Then we have $\delta\in G(\QQ_{p^n})$. The $\sigma$-conjugacy class of $\delta$ in $G(\QQ_{p^n})$ depends only on $[b] \in \B(I_0)$, not on the choices of $b$ and $c$. Moreover, the element $\gamma_0\in G(\Q_p)$ is a degree $n$ norm of $\delta\in G(\Q_{p^n})$.
 \end{lem}
 
 \begin{proof}
 	By (\ref{eqn in KP1}) we have $(\delta\rtimes \sigma )^n = c^{-1 } \gamma_0 c \rtimes \sigma^n $. Since $b$ and $1\rtimes \sigma$ both commute with $\gamma_0$, we know that $\delta\rtimes \sigma = c^{-1} (b\rtimes \sigma ) c$ commutes with $c^{-1} \gamma_0 c$. Hence $\delta\rtimes \sigma $ commutes with $\sigma^n$, which means that $\delta \in G(\QQ_{p^n})$.
 	
 	To prove the second statement, we first note that when $b$ is fixed, any choice of $c$ has to satisfy the equation
 	\begin{align*}
 		\gamma_0 c = b \sigma(b) \cdots \sigma^{n-1}( b) \sigma^n (c).
 	\end{align*}
 	Hence two choices of $c$ give rise to the same coset $cG(\QQ_{p^n}) $. It follows that the $\sigma$-conjugacy class of $\delta$ in $G(\QQ_{p^n})$ is independent of the choice of $c$. Now suppose we choose another representative $b' \in I_0(\LL)$ of $[b] \in \B(I_0)$. Then $b' = d b \sigma( d)^{-1}  $ for some $d \in I_0(\LL)$. Letting $c' : = dc$, we have
 	$ c'^{-1} b'  \sigma(c') = \delta ,$ and
 	$$c'^{-1} \gamma_0 c' = c^{-1}\gamma_0 c = \delta \sigma (\delta) \cdots \sigma^{n-1}(\delta). $$ Hence $b'$ still determines the same $\delta\in G(\Q_{p^n})$ up to $\sigma$-conjugacy.
 	
 	Finally, we show that $\gamma_0$ is a degree $n$ norm of $\delta$. If $G_{\der}$ is simply connected, then the statement follows from (\ref{eqn in KP1}) and the definition of norm. In general, take a $z$-extension $1 \to Z \to H \to G_{\Q_p} \to 1$ over $\Q_p$. Let $\tilde c \in H(\LL)$ be a lift of $c \in G(\LL)$, and let $\tilde \delta \in H(\Q_{p^n})$ be a lift of $\delta \in G(\Q_{p^n})$. Define $\tilde \gamma_0 \in H(\LL)$ by the equation
 	\begin{align}\label{eq:tilde norm}
 		\tilde c^{-1} \tilde \gamma_0 \tilde c = \tilde \delta \sigma (\tilde \delta) \cdots \sigma^{n-1}(\tilde \delta).
 	\end{align}
 	Then $\tilde \gamma_0 $ is a lift of $\gamma_0$. We claim that $\tilde \gamma_0 \in H(\Q_p)$. Once the claim is proved, we know that $\tilde \gamma_0$ is a degree $n$ norm of $\tilde \delta$ by (\ref{eq:tilde norm}), and it follows that $\gamma_0$ is a degree $n$ norm of $\delta$ (see \cite[\S 5]{Kot82}).
 	
 	It remains to prove the claim. Let $\tilde b = \tilde c \tilde \delta \sigma (\tilde c)^{-1} \in H (\LL)$. By (\ref{eq:tilde norm}) we have $(\tilde \delta \rtimes \sigma)^n = \tilde c^{-1} \tilde \gamma_0 \tilde c \rtimes \sigma^n$. Since $\tilde \delta \rtimes \sigma$ commutes with $\sigma^n$ (i.e., $\tilde \delta \in H(\Q_{p^n})$), it must commute with $\tilde c^{-1} \tilde \gamma_0 \tilde c$. Hence $\tilde b \rtimes \sigma = \tilde c (\tilde \delta \rtimes \sigma) \tilde c^{-1}$ commutes with $\tilde \gamma_0$. On the other hand since $\tilde b$ is a lift of $b \in I_0$, we know that $\tilde b$ commutes with $\tilde \gamma_0$ by \cite[Lem.~3.1 (1)]{Kot82}. Therefore $\sigma$ commutes with $\tilde \gamma_0$, which proves the claim.
 \end{proof}

 \begin{lem}\label{polar decomp'} Keep the setting of \S \ref{para:setting KP1}. Let $b \in I_0(\LL)$ be a decent representative of $[b] \in \B(I_0)$ (see \S \ref{para:B(G)}). Then there exists $t\in \NN$ such that $b$ is $t$-decent, and such that
 	\begin{align}\label{eq:polar'}
 		\gamma_0^t = p^{nt\nu_b}k,
 	\end{align} where $k$ lies in a bounded subgroup of $G(\LL)$ (in the sense of \cite[\S 2.2]{Tits}). 
 \end{lem}
 \begin{proof} Assume $b$ is $t_0$-decent. Let $c$ be as in \textbf{KP1}. Let $\cG$ be a parahoric model of $G_{\Q_p}$ over $\ZZ_p$, and write $\bar c$ for the image of $c$ in $G(\LL)/\cG(\breve \ZZ_p)$. By (\ref{eqn in KP1}), we have $
 	\gamma_0 \cdot \bar c = (b\rtimes \sigma )^n \cdot \bar c.$ Since $\gamma_0$ commutes with $b$ and $1 \rtimes \sigma$, for any multiple $t$ of $t_0$ we have $$\gamma_0^{t} \cdot \bar c = (b\rtimes \sigma)^{nt}  
 	\cdot \bar c = (p^{nt \nu_b} \rtimes \sigma ^{nt} ) \cdot \bar c.
 	$$ By the density of $G(\Qpur)$ in $G(\breve \QQ_p)$ for the $p$-adic topology, when $t$ is sufficiently divisible we have $\sigma^{nt} \bar c =\bar c$. Then   $k:= p^{-nt\nu_b} \gamma_0^t$ lies in $c\cG(\breve \ZZ_p)c^{-1}$, which is a bounded subgroup of $G(\LL)$.\end{proof}
 
 One can view (\ref{eq:polar'}) as a ``polar decomposition'' of $\gamma_0^t$. Such a decomposition is unique in a very strong sense, as specified in the following lemma. 
 \begin{lem}\label{uniqueness of polar decomp} Let $F$ be a complete discretely valued field. Let $H$ be a linear algebraic group over $F$, and let $\epsilon$ be a semi-simple element of $H(\overline F)$. Then there exists at most one cocharacter $\nu$ of $H_{\epsilon}$ defined over $\overline F$ with the following property: There exists a uniformizer $\pi \in F^{\times}$ such that $\epsilon ^{-1} \pi^{\nu}$ lies in an $H(\overline F  )$-conjugate of a bounded subgroup of $H(F)$.  
 \end{lem}
 \begin{proof}Let $\rho: H \to \GL_N$ be a  faithful representation of $H$ over $F$. If $k$ is an element of a bounded subgroup of $H(F)$, then $\rho(k)$ lies in a $\GL_N(F)$-conjugate of $\GL_N(\oo_F)$ (cf.~the proof of \cite[Lem.~2.3.1]{kisin2010integral}), and hence all eigenvalues of $\rho(k)$ have valuation zero. (Here and below, by eigenvalues we always mean eigenvalues in $\ol F$.) To prove the lemma it suffices to prove that for each semi-simple $\epsilon \in \GL_N(\ol F)$, there exists at most one cocharacter $\nu$ of $\GL_{N,\epsilon, \ol F}$ such that for some uniformizer $\pi \in F$ all eigenvalues of $\epsilon^{-1} \pi^\nu$ have valuation zero. 
 	
 	Without loss of generality we may assume that $\epsilon=\diag(\lambda_1 I_{N_1},\ldots, \lambda_r I_{N_r})$ with distinct $\lambda_1,\ldots ,\lambda_r \in \overline F^\times$. Then $\GL_{N,\epsilon}$ is naturally identified with $$\GL_{N_1} \times \cdots \times \GL_{N_r}.$$ Note that if a semi-simple element $k \in \GL_{N,\epsilon}(\overline F)$ is such that all its eigenvalues have valuation zero, then the projection of $k$ in each $\GL_{N_i}(\overline F)$ satisfies the analogous condition. We have thus reduced to the case where $\epsilon$ is central in $\GL_N(\overline F)$. In this case, if 
 	$\nu$ is a cocharacter of $\GL_{N, \ol F}$ such that for some uniformizer $\pi$ all eigenvalues of $\epsilon^{-1} \pi^{\nu}$ have valuation zero,  then $\nu$ must be given by $z \mapsto \diag(z^m,\ldots, z^m)$ where $m$ is the valuation of the scalar $\epsilon$. This proves the uniqueness of $\nu$. \end{proof}

 \begin{cor}\label{b is basic}
 	Keep the setting of \S \ref{para:setting KP1}. Then $[b]$ is basic in $\B(I_{0})$. If $[b'] \in \B(I_{0})$ is another class satisfying \textbf{KP1} with respect to $\gamma_0$, then $\nu_{b} = \nu_{b'}$.  
 \end{cor}
 \begin{proof} Let $b$ be a decent representative of $[b]$. By Lemmas \ref{polar decomp'} and  \ref{uniqueness of polar decomp} (the latter applied to $F  = \LL$, $\pi = p$, $H =G_{\LL}$,  $\epsilon = \gamma_0^t$), any element of $G(\LL)$ that centralizes $\gamma_0$ has to centralize $\nu_b$. Therefore $\nu_b $ factors through the center of $I_0$, and $[b]$ is basic in $\B(I_0)$. The second statement also follows from these two lemmas. 
 \end{proof}
 \begin{cor}\label{cor:basic for adm Kott para}
 	Let $\fkc = (\gamma_0, a, [b])\in \KP_{\adm}(p^n)$. Then $[b]$ is basic in $\B(I_0(\fkc)_{\Q_p})$. If $(\gamma_0, a',[b'])$ is another element of $\KP_{\adm}(p^n)$, then $\nu_b = \nu_{b'}$. 
 \end{cor}
 \begin{proof}
 	This follows from Corollary \ref{b is basic}. 
 \end{proof}
 
 \subsection{}\label{para:comparing Kott invt}
 Our next goal is to define the notion of an isomorphism between Kottwitz parameters. We first consider a general construction. Let $\gamma_0 \in G(\QQ)_{\semi}$, and let $I_0 : = G_{\gamma_0}^0$. Let $u \in G(\Qbar)$ be an element satisfying
 $$ \gamma_0' : = u \gamma_0 u^{-1} \in G(\Q)$$ and
 $$ u^{-1} \lix^{\rho} u \in I_0, \quad \forall \rho \in \Gamma = \Gal(\Qbar/\QQ). $$ Let $I_0' : = G_{\gamma_0'}^0$.
 We have a bijection
 \begin{align}\label{eq:u_* away from p}
 	u_*: \D(I_0, G ; \adele_f^p) \isom \D(I_0', G ; \adele_f^p),
 \end{align}
 sending the class of a cocycle $(a_{\rho})_\rho \in Z^1(\Gamma, I_0(\bar \A_f^p))$ to the class of $(u a_{\rho} \lix^{\rho} u^{-1})_{\rho}$.
 Next note that the cocycle $( u^{-1}\lix{^\rho} u  )_{\rho}  \in Z^1(\QQ_p, I_0)$ becomes trivial in $\coh^1( \LL, I_0)$ by Steinberg's theorem. Hence there exists $d\in I_0(\overline {\breve \QQ}_p)$ such that
 \begin{align*}
 	u^{-1}\lix{^\rho} u  = d^{-1}\lix{^\rho} d, \quad  \forall \rho \in  \Gamma_{p,0}.
 \end{align*}
 Then $u_0:= ud^{-1} $ lies in $G(\LL)$. We have
 \begin{align*}
 	u_0 \gamma_0 u_0^{-1}= \gamma'_0,
 \end{align*}
 and
 \begin{align*}
 	u_0^{-1} \lix{^\sigma} u_0  \in  I_0 (\LL).
 \end{align*}
 By the previous two properties of $u_0$, we have a bijection
 \begin{align}\label{eq:u_* at p}
 	u_* : \B(I_{0, \Q_p}) & \isom \B(I'_{0,\Q_p}) \\ \nonumber [b] &\longmapsto [ u_0 b \sigma ( u_0)^{-1}   ] ,
 \end{align}
 which depends only on $u$, not on the choice of $d$.

 \begin{defn}\label{defn of equivalence between coh kott trip}
 	Let $ \fkc = (\gamma_0 , a ,[b]), \fkc' =  ( \gamma'_0, a' , [b'])  \in \KP$. By an \textit{isomorphism from $\fkc$ to $\fkc'$}, we mean an element $u\in G(\Qbar)$,  satisfying the following conditions.
 	\begin{enumerate}
 		\item We have $\Int(u) \gamma_0 = \gamma'_0$, and $u^{-1}\lix{^\rho} u \in I_0(\fkc)$ for all $\rho \in \Gamma$.
 		\item The bijection $u_*: \D(I_0(\fkc), G ; \adele_f^p) \isom \D(I_0(\fkc'), G ; \adele_f^p)$ as in (\ref{eq:u_* away from p}) sends $a$ to $a'$.
 		\item The bijection $u_* : \B(I_0(\fkc)_{\Q_p})\isom \B(I_0(\fkc')_{\Q_p})$ as in (\ref{eq:u_* at p}) sends $[b]$ to $[b']$.
 	\end{enumerate}
 	In such a situation we write $ u : \fkc \isom \fkc'.$ \end{defn}
 
 \subsection{}
 If $u: \fkc \isom \fkc'$ and $v: \fkc' \isom \fkc''$ are two isomorphisms between Kottwitz parameters, then $vu \in G(\Qbar)$ is an isomorphism $\fkc \isom \fkc ''$, and $u^{-1} \in G(\Qbar)$ is the isomorphism $\fkc' \isom \fkc$ inverse to $u$. Moreover, one checks that $p^n$-admissibility of Kottwitz parameters (Definition \ref{pn-adm}) is preserved under isomorphisms.
 
 We denote by $\KP/{\cong}$\index[Oct2026-sorted-notation]{kp 0 04 00 kp cong@$\mathfrak{KP}/{\cong}$} the set of isomorphism classes of Kottwitz parameters, and by ${\KP_{\adm}(p^n)}/{\cong}$\index[Oct2026-sorted-notation]{kp 0 04 00 kp subscript a p superscript n cong@$\mathfrak{KP}_{\mathrm{a}}(p^n)/{\cong}$} the set of isomorphism classes of $p^n$-admissible Kottwitz parameters.
 
 \subsection{}
 \label{from CK to K} 
 Let $\fkc = (\gamma_0, a, [b]) \in \KP_{\adm}(p^n)$. As usual, $a$ determines a conjugacy class in $G(\A_f^p)$ which is stably conjugate to $\gamma_0$. By Lemma \ref{b decides delta}, $[b]$ determines a $\sigma$-conjugacy class in $G(\Q_{p^n})$, of which $\gamma_0$ is a degree $n$ norm. Let $\gamma \in G(\A_f^p)$ and $\delta \in G(\Q_{p^n})$ be representatives of these classes. Then $(\gamma_0, \gamma, \delta)$ is an element of $\KP_{\clsc}(p^n)$, and its equivalence class depends only on $\fkc$. We thus obtain a map $$
 \KP_{\adm}(p^n) \To \KP_{\clsc}(p^n)/{\sim}, $$
 which is easily seen to factor through $\KP_{\adm}(p^n)/{\cong}$. Moreover, when $G_{\der}$ is simply connected, the induced map $\KP_{\adm}(p^n)/{\cong} ~ \To  \KP_{\clsc}(p^n)/{\sim}$ is a bijection. We will not need this fact in the sequel, but we outline here how to recover $[b] \in \B(I_{0}(\fkc)_{\Q_p})$ from $\delta$.  Since $\gamma_0$ is a degree $n$ norm of $\delta$ and since $\coh^1(\LL, I_0(\fkc))$ is trivial by Steinberg's theorem, there exists $c \in G(\LL)$ satisfying (\ref{eqn in KP1}). Define $b : = c \delta \sigma(c)^{-1} \in G(\LL)$. Since $(\delta\rtimes \sigma)^n = c^{-1} \gamma_0 c \rtimes \sigma^n$ and since $\delta \rtimes \sigma$ commutes with $\sigma^n$, we know that $\delta \rtimes \sigma$ commutes with $c^{-1} \gamma_0 c$. Hence  $b\rtimes \sigma$ commutes with $\gamma_0$. Thus $b$ commutes with $\gamma_0$, so $b \in I_0(\fkc)(\LL)$ as $G_{\der}$ is simply connected. In this way we have recovered $[b] \in \B(I_{0}(\fkc)_{\Q_p})$ from $\delta$.  
 
 \section{The Kottwitz invariant}\label{Kottinv}
 The goal of this section is to attach an \emph{invariant} to each Kottwitz parameter; our eventual point counting formula will involve a summation over parameters for which this invariant vanishes.  
 
 \subsection{}
 We first construct the abelian group in which the invariant lies. Let $G$ be a reductive group over $\Q$, and   $I$   a reductive subgroup. We have an infinite commutative diagram with exact rows and columns
 $$\xymatrix@C=1.0em@R=1.5em { \vdots \ar[d] & \vdots \ar[d] & \vdots \ar[d] & \vdots \ar[d] \\ \cdots \to  \coh^i_{\ab}(\QQ,I) \ar[r] \ar[d] & \coh^i_{\ab}(\QQ, G) \ar[r]\ar[d]  & \coh^i_{\ab}(\QQ, I \to G) \ar[r]  \ar[d] & \coh^{i+1}_{\ab}(\QQ, I) \ar[d]  \to \cdots \\   \cdots \to  \coh^i_{\ab}(\adele,I) \ar[r] \ar[d] & \coh^i_{\ab}(\adele, G) \ar[r]\ar[d] \ar@{-->}[rd] & \coh^i_{\ab}(\adele, I \to G) \ar[r]  \ar[d] & \coh^{i+1}_{\ab}(\adele, I) \ar[d] \to \cdots  \\ \cdots \to \coh^i_{\ab}(\adele/\QQ,I) \ar[r] \ar[d] & \coh^i_{\ab}(\adele/\QQ, G) \ar[r]\ar[d]  & \coh^i_{\ab}(\adele/\QQ, I \to G) \ar[r]  \ar[d] & \coh^{i+1}_{\ab}(\adele/\QQ, I) \ar[d] \to \cdots  \\ \cdots \to \coh^{i+1}_{\ab}(\QQ,I) \ar[r] \ar[d] & \coh^{i+1}_{\ab}(\QQ, G) \ar[r]  \ar[d]& \coh^{i+1}_{\ab}(\QQ, I \to G) \ar[r]  \ar[d] & \coh^{i+2}_{\ab}(\QQ,I) \ar[d] \to \cdots \\ \vdots & \vdots & \vdots &\vdots  }$$ We define $
 \E( I,G ; \adele/\QQ) $\index[Oct2026-sorted-notation]{e 0 04 00 e i g a q@$\mathfrak E( I,G ; \mathbb A/\QQ)$} to be the cokernel of the map $$\coh^0_{\ab}(\adele, G)  \To \coh^0_{\ab} (\adele/\QQ, I \to G )  $$ given by the dashed arrow in the above diagram for $i=0$. We have a natural map $$\E(I, G ; \adele) \To \E(I, G ; \A/\Q)$$ defined by first lifting an element of $\E(I, G ; \A)$ to $\coh^0_{\ab}(\A, I\to G)$, and then mapping the lift to $\coh^0_{\ab}(\A/\Q, I \to G)$ and then to $\E(I, G; \A/\Q)$. We know that the sequence \begin{align}  \label{eq:how to show zero invariant}
 	\E( I,G ; \Q ) \To \E( I,G ; \A) \To \E( I,G ; \A/\Q) \end{align} is exact; see \cite[Prop.~1.8.4]{Lab99}.
 
 In the following, we write $\Gamma$ for $\Gal(\Qbar/\Q)$ and write $\Gamma_v $ for $\Gal(\ol \QQ_v/ \QQ_v)$ for every place $v$ of $\QQ$. 
 \begin{prop}\label{computing E} Assume that $I$ contains a maximal torus in $G$. Let $K$ be the kernel of $\pi_1(I) \to \pi_1(G)$. Then the abelian group $ \E(I, G ; \adele/\Q)$ is canonically identified with $$ \frac{  K_{\Gamma ,\tors}               }{ \bigoplus_v \ker( K_{\Gamma_v,\tors} \to \pi_1(I)_{\Gamma_v               })}.$$ Here $v$ runs through all the places of $\Q$, and the quotient is with respect to the natural maps $K_{\Gamma_v,\tors} \to K_{\Gamma, \tors}$. In particular, $\E(I, G ; \adele/\Q)$ is finite.
 \end{prop}
 \begin{proof} Under our assumption, the map $\pi_1(I) \to \pi_1(G)$ is surjective. Applying Proposition \ref{prop:apply TN} (ii), we know that $\E(I, G ; \adele/\Q)$ is canonically identified with the cokernel of
 	$$ \mathscr P(K): \ker \big(\mathcal B_\Q(K) \to \mathcal B_\Q(\pi_1(I)) \big ) \To \mathcal A_\Q(K).$$
 	By Lemma \ref{lem:eq of kernels}, we have
 	$$\ker \big(\mathcal B_\Q(K) \to \mathcal B_\Q(\pi_1(I)) \big ) = \bigoplus_v \ker( K_{\Gamma_v,\tors} \to \pi_1(I)_{\Gamma_v               }). $$ The proposition follows.
 \end{proof}
 For any Hausdorff locally compact abelian group $H$, we denote by $H^D$ the Pontryagin dual. The following result is well known to experts, cf.~\cite[p.~43, Remarque]{Lab99}.
 \begin{cor}\label{cor:compare with Kottwitz}Under the assumption in Proposition \ref{computing E}, the abelian group $\E(I, G ; \adele/\Q)$ is canonically identified with $\mathfrak K (I/\Q)^D$, where $\mathfrak K (I/\Q)$ is defined in \cite[\S 4.6]{Kot86}. \end{cor}
 \begin{proof}
 	We shall freely use the definitions and results in \cite{Kot86}. Recall that $\mathfrak K(I/\Q)$\index[Oct2026-sorted-notation]{k 0 04 00 k i q@$\mathfrak K(I/\Q)$} is defined as the subgroup of $\pi_0 ([Z(\hat I) / Z(\hat G)] ^{\Gamma})$ consisting of those elements whose images in $\coh^1 (\Q, Z(\hat G)) $ are locally trivial.
 	
 	Since $X^* (Z(\hat I)) \cong \pi_1(I)$ and $X^*(Z(\hat G)) \cong \pi_1(G)$, we have $ X^*( Z(\hat I) / Z(\hat G))  \cong  K$. Hence $$X^*\bigg(\pi_0([Z(\hat I) / Z(\hat G)]^{\Gamma}) \bigg) \cong  K_{\Gamma,\tors}. $$
 	
 	Also, for each place $v$, we have $$X^*\bigg(\pi_0 (Z(\hat I) ^{\Gamma_v} ) \bigg) \cong \pi_1 (I) _{\Gamma_v,\tors}$$ and $$  X^*\bigg(\pi_0 ([Z(\hat I) / Z(\hat G)] ^{\Gamma_v}) \bigg) \cong K_{\Gamma_v,\tors}. $$
 	By the exact sequence
 	$$ \pi_0 (Z(\hat I) ^{\Gamma_v}) \to \pi_0 ([Z(\hat I) / Z(\hat G)] ^{\Gamma_v}) \to \coh^1(\Q_v, Z (\hat G)),$$
 	we can identify $\mathfrak K (I/\Q)$ with the set of $x \in (K_{\Gamma , \tors})^D$ such that for each place $v$ the composite map $x_v: K_{\Gamma_v,\tors} \to K_{\Gamma,\tors}\xrightarrow{x} \CC^{\times}$ equals the composite map $K_{\Gamma_v, \tors} \to \pi_1(I)_{\Gamma_v,\tors} \xrightarrow{y_v} \CC^{\times}$, for some  $y_v \in (\pi_1(I)_{\Gamma_v, \tors})^D$.
 	By the exact sequence
 	$$ (\pi_1(I)_{\Gamma_v, \tors})^D \To (K_{\Gamma_v, \tors})^D \To \bigg(\ker( K_{\Gamma_v,\tors} \to \pi_1(I)_{\Gamma_v, \tors})\bigg)^D, $$
 	the last condition on $x_v$ is equivalent to requiring that $x$ kills the image of $$\ker( K_{\Gamma_v,\tors} \to \pi_1(I)_{\Gamma_v, \tors})$$ in $K_{\Gamma,\tors}$. Comparing this description of $\mathfrak K (I/\Q)$ and Proposition \ref{computing E}, we see that $\mathfrak K(I/\Q) \cong \E(I, G; \adele/\Q)^D$.
 \end{proof}

 \subsection{}
 \label{defn of Kottwitz invariant} We now keep the setting of \S \ref{para:fkp}. 
 Let $\fkc = (\gamma_0, a ,[b] )\in \KP$, and write $I_0$ for $I_0(\fkc) = G_{\gamma_0}^0$.  We now construct an element \index[Oct2026-sorted-notation]{zz01 alpha 1 00 00 alpha c@$\alpha(\mathfrak c)$} $$\alpha(\fkc)  \in  \E(I_0, G; \adele/\QQ),$$ called the\index{Kottwitz invariant} \emph{Kottwitz invariant} of $\fkc$. This generalizes the construction in \cite[\S 2]{Kot90}.
 
 We write\index[Oct2026-sorted-notation]{zz02 beta 1 00 00 beta superscript p infty c@$\beta^{p,\infty}(\fkc)$} $\beta^{p,\infty}(\fkc)$ for the element $a \in \D(I_0, G; \adele_f^p)$. As discussed in \S \ref{para:two groups}, the abelianization map induces an isomorphism $\D(I_0, G ;\adele_f^p) \cong \E(I_0, G;\adele_f^p)$. By Proposition \ref{prop:apply TN} (i), we have a canonical isomorphism $$\coh^1_{\ab}(\A_f^p, I_0) \cong \bigoplus_{v \neq p,\infty} \pi_1(I_0)_{\Gamma_v, \tors}. $$ Hence we also view $\beta^{p,\infty}(\fkc)$ as an element of $$\bigoplus_{v \neq p,\infty} \pi_1(I_0)_{\Gamma_v, \tors}. $$ For each place $v \notin \set{p, \infty}$, we write $\beta_v(\fkc) \in \pi_1(I_0)_{\Gamma_v, \tors }$ for the component of $\beta^{p,\infty}(\fkc)$ at $v$.
 We pick a lift $\tilde \beta_v (\fkc) \in \pi_1(I_0)$ of $\beta_v(\fkc)$ that maps to zero in $\pi_1(G)$. Such a lift exists, since $\beta_v(\fkc)$ maps to zero in $\pi_1(G)_{\Gamma_v}$ and since the map $\pi_1(I_0) \to \pi_1(G)$ is surjective.
 Since $\beta_v (\fkc) =0$ for almost all $v$, we may and shall assume that $\tilde \beta_v(\fkc) =0$ for almost all $v$.
 
 Let\index[Oct2026-sorted-notation]{zz02 beta 1 00 00 beta subscript p c@$\beta_p (\fkc)$} $ \beta_p (\fkc) := \kappa_{I_{0}} ([b]) \in \pi_1(I_0)_{\Gamma_p}$. We pick a lift $\tilde\beta_p (\fkc) \in \pi_1( {I_0} )$ of $\beta_p(\fkc)$ that maps to $ - [\mu]_X\in \pi_1(G)$. Such a lift exists by condition \textbf{KP0} in Definition \ref{defn of coh Kott trip} and the surjectivity of the map $\pi_1(I_0) \to \pi_1(G)$.
 
 Now we take an elliptic maximal torus $T$ in $G_{\RR}$ such that $\gamma_0\in T(\RR)$. Then $T\subset I_{0,\RR}$. Since $T$ is elliptic, there exists $h\in X$ that factors through $T$. Let\index[Oct2026-sorted-notation]{zz02 beta 1 00 00 beta subscript infty c@$\beta_{\infty}(\fkc)$} $\beta_\infty(\fkc)$ be the image of $\mu_h \in X_*(T)$ in $\pi_1(I_0)_{\Gamma_\infty}$. This is independent of $T$ and $h$ by \cite[Lem.~5.1]{Kot90} and the fact that all elliptic maximal tori in $I_{0,\RR}$ are conjugate under $I_0(\RR)$. The image of $\beta_{\infty}(\fkc)$ under $\pi_1(I_0)_{\Gamma_{\infty}} \to \pi_1(G)_{\Gamma_{\infty}}$ equals the image of $[\mu]_X$. Since $\pi_1(I_0) \to \pi_1(G)$ is surjective, we can pick a lift $\tilde\beta_\infty (\fkc) \in \pi_1( I_0)$ of $\beta_\infty(\fkc)$ that maps to $[\mu]_X \in \pi_1(G)$. 
 
 Write 
 $K$ for $\ker(\pi_1(I_0) \to \pi_1(G))$. By the above construction, we have an element
 $\tilde \beta_v(\fkc) \in K$ for each place $v \notin \set{p, \infty}$, as well as an element $ \tilde \beta_p (\fkc)+ \tilde \beta_\infty (\fkc) \in K$. 
 We define $$ \tilde \beta(\fkc) : = \sum_{v} \tilde \beta_v (\fkc) \in K. $$ Here the summation is over all places $v$ of $\QQ$, and only finitely many terms are non-zero. 
 
 Note that $K_{\Gamma_{\infty}}$ is torsion. Indeed, if we take an $\RR$-elliptic maximal torus $T$ in $I_{0,\RR}$ and let $\widetilde T$ and $\widetilde S$ be the inverse images of $T$ in $G_{\sconn,\RR}$ and $I_{0, \sconn, \RR}$ respectively, then $K \cong X_*(\widetilde T)/X_*(\widetilde S)$. Since $\widetilde T$ is anisotropic over $\RR$, we know that $X_*(\widetilde T) _{\Gamma_{\infty}}$ is torsion. It follows that $K_{\Gamma_{\infty}}$ is torsion. In particular,  $K_{\Gamma}$ is torsion. Hence $\tilde \beta(\fkc)$ gives rise to an element $\alpha(\fkc) \in \E(I_0, G ; \A/\QQ)$, by Proposition \ref{computing E}. Note that the ambiguity in $\tilde\beta(\fkc)$ caused by the choices of $\tilde\beta_v(\fkc)$ always comes from $ \bigoplus_{v} \ker(K \to \pi_1(I_0)_{\Gamma_v}).$ Hence $\alpha(\fkc)$ is well defined.

 \subsection{}
 \label{subsub:local_Kottwitz_inv} 		It is convenient to have the definition of local components of Kottwitz invariants when $\gamma_0$ is not $\Q$-rational for the stabilization of the trace formula.
 Suppose $\gamma_0$ is a semi-simple element of $G(\A_f^p)$. We have $I_0 = G_{\gamma_0}^0$ over $\A_f^p$, and
 we define the pointed set $\fkD(I_0,G;\A_f^p)$ to be the restricted product $\prod'_v \D(I_0, G; \Q_v)$ with respect to the trivial elements.   If we are given an element $a \in \D(I_0, G; \A_f^p)$, we will also write $\beta^{p,\infty}(\gamma_0,a)$ for $a$, generalizing the notation $\beta^{p,\infty}(\fkc)$ in \S \ref{defn of Kottwitz invariant}. Similarly, given $\gamma_0\in G(\Q_p)_{\semi}$ and given $[b]\in \B((G_{\Q_p})_{\gamma_0}^0)$ satisfying \textbf{KP0} in Definition \ref{defn of coh Kott trip}, we have $\beta_p(\gamma_0,[b])$ and $\tilde \beta_p(\gamma_0,[b])$ (the latter involving an extra choice), generalizing $\beta_p(\fkc)$ and $\tilde \beta_p(\fkc)$ in \S \ref{defn of Kottwitz invariant}. Finally, starting from a semi-simple $\RR$-elliptic element $\gamma_0\in G(\R)$,  we can define $\beta_\infty(\gamma_0)$ and $\tilde \beta_\infty(\gamma_0)$ (the latter involving an extra choice) in the same way as in \S \ref{defn of Kottwitz invariant}, generalizing $\beta_{\infty}(\fkc)$ and $\tilde \beta_\infty (\fkc)$.

 \subsection{}
 \label{subsub:comparison_Kottwitz_triples} 
 Let us check that our definition of Kottwitz invariants coincides with Kottwitz's definition in \cite[\S 2]{Kot90}, when $G_{\der}$ is simply connected. This verification allows us to freely import results from \textit{op.~cit.}~during the stabilization process.
 
 Under the identification $\pi_1(G) \cong X^*(Z(\widehat G))$, we view $[\mu]_X \in \pi_1(G)$ as a character on $Z(\widehat{G})$. Under the current assumption on $G_{\der}$, recall that Kottwitz attaches an invariant $\alpha(\fkk)\in\fkK(I_0/\Q)^D$ to a classical Kottwitz parameter $\fkk=(\gamma_0,\gamma,\delta)$ of degree $n$. The outline is as follows. Let $I_0 : = G_{\gamma_0}$, which is connected by the assumption on $G_{\der}$. Kottwitz first defines $\alpha_v(\fkk)\in \pi_1(I_0)_{\Gamma_v} \cong X^*(Z(\widehat I_0)^{\Gamma_v})$ at every place $v$. Then the character $\alpha_v(\fkk)$ on $Z(\hat{I_0})^{\Gamma_v}$ can be extended to a character $\beta'_v(\fkk)$ on $Z(\hat{I_0})^{\Gamma_v} Z(\hat G)$, uniquely by the requirement that $\beta'_v(\fkk)$ is either trivial or equal to $-[\mu]$ or 
 $[\mu]$ on $Z(\hat G),$ according as $v\notin \set{p,\infty}$ or $v=p$ or $v=\infty$, respectively.
 \footnote{We write $\beta'$ for Kottwitz's $\beta$ to avoid conflict with our own notation.} Thereby one obtains a character $\beta(\fkk):=\prod_v \beta'_v(\fkk)$ on $\bigcap_v Z(\hat{I_0})^{\Gamma_v} Z(\hat G)$. Since $\beta(\fkk)$ is trivial on $Z(\hat G)$ by construction, it gives rise to a character $\alpha(\fkk)$ on $\left(\bigcap_v Z(\hat{I_0})^{\Gamma_v} Z(\hat G)\right)/Z(\hat G)$. The last group is canonically isomorphic to $\fkK(I_0/\Q)$ if $\gamma_0$ is elliptic in $G(\Q)$, which is true since $\gamma_0$ is $\RR$-elliptic. We note that the canonical map from $K=\ker(\pi_1(I_0)\ra \pi_1(G))$ to the character group of $\cap_v Z(\hat{I_0})^{\Gamma_v} Z(\hat G)$ is compatible with the canonical map from $K$ to $\fkE(I_0,G,\A/\Q)$ via the isomorphisms
 $$ \fkE(I_0,G;\A/\Q) \cong \fkK(I_0/\Q)^D \cong \left(\frac{\cap_v Z(\hat{I_0})^{\Gamma_v} Z(\hat G)}{Z(\hat G)}\right)^D$$
 as can be seen from the proof of Corollary \ref{cor:compare with Kottwitz}.
 
 When $G_{\der}$ is simply connected, we remarked in \S \ref{from CK to K} that we have a bijection  $\KP_{\adm}(p^n)/{\cong} \to \KP_{\clsc}(p^n)/{\sim}$. Suppose $\fkc=(\gamma_0,a,[b])\in \KP_{\adm}(p^n) $ corresponds to $\fkk=(\gamma_0,\gamma,\delta) \in \KP_{\clsc}(p^n)$. We defined $\beta_v(\fkc)$ and $\tilde \beta_v(\fkc)$ in \S\ref{defn of Kottwitz invariant}, for each place $v$. For $v\in \{p,\infty\}$, note that $\tilde \beta_v(\fkc)$ is a character on $Z(\hat I_0)$, extending the character $\beta_{v}(\fkc)$ on $Z(\hat I_0)^{\Gamma_v}$. Inspecting the definition we see that
 $$\beta_v(\fkc) =  \alpha_v(\fkk),\quad \forall v,$$
 and $$
 \left.\tilde \beta_v(\fkc)\right|_{Z(\hat{I_0})^{\Gamma_v} Z(\hat G)} = \beta'_v(\fkk),\quad  \forall v.$$
 Therefore the product $ \prod_v \tilde\beta_v(\fkc)$ gives the element $\alpha(\fkk)\in \fkK(I_0/\Q)^D$. Comparing with the definition of $\alpha(\fkc)$ in \S \ref{defn of Kottwitz invariant}, we conclude that
 $$\alpha(\fkc)=\alpha(\fkk).$$
 
 From now on we return to the general setting, i.e., we do not assume that $G_{\der}$ is simply connected.
 
 \begin{prop}\label{diff to diff}
 	Let $\fkc = (\gamma_0,  a  , [b])$, $\fkc' =  (\gamma_0, a' , [b]) \in \KP$. Write $I_0$ for the group $I_0(\fkc)  =I_0(\fkc')$. The difference
 	$\alpha(\fkc)-\alpha(\fkc') \in \E(I_0, G ;\adele/\Q)$ is equal to the image of
 	$a- a'  \in \D (I_0, G ; \A_f^p)$
 	under the composite map
 	\begin{align}\label{eq:D to E(fkc)}
 		\D (I_0 , G ; \A_f^p) \cong \E (I_0 , G ; \A_f^p)  \hookrightarrow \E (I_0, G; \adele ) \To \E(I_0, G ;\adele/\QQ) .\end{align}
 \end{prop}
 \begin{proof} Write $K$ for $\ker(\pi_1(I_0) \to \pi_1(G))$. By Proposition \ref{prop:apply TN} (ii), the diagram
 	$$\xymatrix{ \coh^0_{\ab}(\adele, I_0 \to G) \ar[r] \ar[d] &  \coh^1_{\ab}(\A, I_0) \\
 		\coh^0_{\ab}(\A/\Q, I_0 \to G)} $$ is canonically isomorphic to the diagram
 	\begin{align}\label{diag:K} \xymatrix{\tatecoh^{-1}(\Gamma_\infty, K) \oplus  \bigoplus_{v \neq \infty} K_{\Gamma_v,\tors} \ar[r] \ar[d] & \tatecoh^{-1}(\Gamma_\infty, \pi_1(I_0))\oplus  \bigoplus_{v \neq \infty} \pi_1(I_0)_{\Gamma_v ,\tors } \\
 			K _{\Gamma, \tors} } 	\end{align}

 	For each place $v \notin \set{ p, \infty}$, we choose $\tilde \beta_v(\fkc)$ and $\tilde \beta_v (\fkc')$ in $K$ as in \S \ref{defn of Kottwitz invariant}. The $\Gamma_v$-actions on $\pi_1(G)$ and  $\pi_1(I_0)$ both factor through some finite quotient $\Gamma_v'$ of $\Gamma_v$, and we have an exact sequence $$\coh_1(\Gamma_v', \pi_1(G)) \to K_{\Gamma_v} \to \pi_1(I_0)_{\Gamma_v}. $$
 	Since the image of $\tilde \beta_v(\fkc)$ in $\pi_1(I_0)_{\Gamma_v}$ is the torsion element $\beta_v(\fkc)$, and since $\coh_1(\Gamma_v', \pi_1(G))$ is torsion, the image of $\tilde \beta_v(\fkc)$ in $K_{\Gamma_v}$ is torsion. We denote this image by $\bar \beta_v(\fkc) \in K_{\Gamma_v,\tors} .$ Similarly we define $\bar \beta_v(\fkc') \in K_{\Gamma_v, \tors} .$ Let $$\Delta: = (\bar \beta_v(\fkc) - \bar \beta_v(\fkc'))_{v \neq p,\infty} \in  \bigoplus_{v \neq p, \infty} K_{\Gamma_v,\tors} . $$ Then $\Delta$ is sent to $a- a'$ by the horizontal map in (\ref{diag:K}). By the above discussion, the image of $a- a'$ under (\ref{eq:D to E(fkc)}) is equal to the image of $\Delta$ under the composite map
 	$$ \bigoplus_{v\neq p,\infty} K_{\Gamma_v, \tors }  \to K_{\Gamma,\tors} \to \E(I_0, G ; \A/\Q),$$ where the last map is the quotient map as in Proposition \ref{computing E}. On the other hand,  the image of $\Delta$ in $\E(I_0, G ; \A/\Q)$ is $\alpha(\fkc)  - \alpha(\fkc')$ by the construction of the Kottwitz invariant. Hence the image of $a- a'$ under (\ref{eq:D to E(fkc)}) is $\alpha(\fkc) - \alpha(\fkc')$.
 \end{proof}
 \subsection{}\label{para:setting invt}
 Let $u : \fkc \isom \fkc'$ be an isomorphism between Kottwitz parameters. Then the isomorphism $\Int(u): I_0(\fkc)_{\Qbar} \isom I_0(\fkc')_{\Qbar}$ is an inner twisting, and in particular it induces an isomorphism of abelian groups 
 \begin{align}\label{eq:comparing Kott inv}
 	\E(I_0(\fkc), G ; \A/\Q) \isom \E(I_0(\fkc'), G; \A/\Q)
 \end{align}The following result justifies that the Kottwitz invariant is indeed an ``invariant''.
 
 \begin{prop}\label{funct of alpha conj}
 	The isomorphism (\ref{eq:comparing Kott inv}) takes $\alpha(\fkc )$ to $\alpha(\fkc')$.  \end{prop}
 \begin{proof}We write $I_0$ and $I_0'$ for $I_0(\fkc)$ and $I_0(\fkc')$. Let $K = \ker(\pi_1(I_0) \to \pi_1(G))$ and $K' = \ker(\pi_1(I_0') \to \pi_1(G))$.  The inner twisting $\Int(u): I_{0,\Qbar} \isom I_{0,\Qbar}'$ induces a $\Gamma$-equivariant isomorphism $\pi_1(I_0) \isom \pi_1(I_0')$, which we denote by $f$. Then $f$ restricts to an isomorphism $K\isom K'$. Moreover, if we identify the two sides of (\ref{eq:comparing Kott inv}) with quotients of  $K_{\Gamma ,\tors}$ and $K'_{\Gamma ,\tors}$ respectively as in Proposition \ref{computing E}, then (\ref{eq:comparing Kott inv}) is induced by $f: K\isom K'$.
 	
 	Let $\omega \in \D(I_0, G; \Q)$ be the class of the cocycle $(u^{-1} \lix^{\rho} u)_{\rho \in \Gamma}$. For each place $v$ of $\Q$, we denote by $\beta_v(\omega)$ the image of $\omega$ under the composite map
 	$$\coh^1(\Q, I_0) \to \coh^1(\Q_v, I_0)  \xrightarrow{\ab^1_{\Q_v}} \coh^1_{\ab}(\Q_v, I_0).$$ Using the isomorphism $\coh^1_{\ab}(\Q_v, I_0) \cong \cA_{\Q_v}(\pi_1(I_0)) \subset  \pi_1(I_0)_{\Gamma_v, \tors}$ in Proposition \ref{prop:apply TN}, we also view $\beta_v(\omega)$ as an element of  $\pi_1(I_0)_{\Gamma_v, \tors}$. We claim that for each place $v$, the image of $\beta_v(\fkc) - \beta_v(\omega)$ under the isomorphism $\pi_1(I_0)_{\Gamma_v} \isom \pi_1(I_0')_{\Gamma_v}$ induced by $f$ equals $\beta_v(\fkc')$.
 	
 	Our claim for $v \notin \set{p,\infty}$ follows from the following commutative diagram, which is a special case of \cite[Prop.~3.16]{borovoi} (with horizontal arrows reversed):
 	$$ \xymatrixrowsep{3.5pc} \xymatrix { \D(I_0, G; \Q_v) \ar[rr]^{u_*}  \ar[d]^{\sim} _{\ab^1}& &  \D(I_0', G; \Q_v) \ar[d]^{\sim}_{\ab^1} \\ \E(I_0, G;\Q_v) \ar[r]^{c} & \E(I_0, G;\Q_v) \ar[r]^d_{\sim} &  \E(I_0', G;\Q_v) } $$ 
 	Here $u_*$ is the component at $v$ of the bijection (\ref{eq:u_* away from p}), $c$ is the translation map $x \mapsto x - \beta_v(\omega)$, and $d$ is the group isomorphism induced by the inner twisting $\Int(u) : I_{0, \Qbar} \isom I'_{0, \Qbar}$.  
 	
 	Similarly, the bijection $u_*: \B(I_{0,\Q_p}) \isom \B(I_{0,\Q_p}')$ as in (\ref{eq:u_* at p}) fits in the commutative diagram
 	$$ \xymatrix { \B(I_{0,\Q_p}) \ar[rr]^{u_*}  \ar[d]_{\kappa_{I_0}}&&  \B(I_{0,\Q_p}') \ar[d]_{\kappa_{I_0'}} \\ \pi_1(I_0)_{\Gamma_p} \ar[r]^{c} &  \pi_1(I_0)_{\Gamma_p}  \ar[r]^f_{\sim} &  \pi_1(I_0')_{\Gamma_p} } $$
 	where $c$ is the translation map $x \mapsto x - \beta_p(\omega)$. To see this, we use the fact that if we choose $u_0$ as in \S \ref{para:comparing Kott invt}, then the image of $u_0^{-1} \lix^\sigma u_0$ under $ \kappa_{I_0} :  I_0(\LL) \to \pi_1(I_0)_{\Gamma_p}$ equals $\beta_p(\omega)$; see \cite[Rmk.~5.7]{kottwitzisocrystal} and \cite[Rmk.~2.2 (iv)]{RV14}. Our claim for $v=p$ follows from the above commutative diagram.
 	
 	Now we prove the claim for $v = \infty$. As in \S \ref{defn of Kottwitz invariant}, we choose $(T,h)$ to define $\beta_\infty(\fkc)$ and choose $(T',h')$ to define $\beta_{\infty}(\fkc')$. Without loss of generality we may assume that $T' = \Int(k)(T)$ and $h'= \Int(k) \circ h$ for some $k \in G(\RR)$. Since elliptic maximal tori transfer between inner forms, there exists $j \in I_0(\CC)$ such that the map $\Int(uj): T_{\CC} \to I'_{0,\CC}$ is defined over $\RR$ and has image $T'$. We have a commutative diagram
 	$$ \xymatrixcolsep{5pc} \xymatrix{ X_*(T) \ar[r]^{\Int(uj)}_{\sim } \ar[d] & X_*(T') \ar[d] \\ \pi_1(I_0) \ar[r]_{\sim}^f & \pi_1(I_0') , }$$
 	where the vertical maps are the natural quotient maps. Let $n : = (uj)^{-1} k  \in G(\CC)$. Then $n$ is in the normalizer of $T$ in $G$, and $\Int(uj)^{-1} (\mu_{h'}) = \Int(n) (\mu_h)$. Let $\Delta_{\infty}$ be the image of $ \Int(n)(\mu_h) - \mu_h$ under $X_*(T) \to \pi_1(I_0)_{\Gamma_\infty}$. Then by the above discussion the image of $\beta_\infty(\fkc) + \Delta_{\infty}$ under $f: \pi_1(I_0)_{\Gamma_\infty} \isom \pi_1(I_0')_{\Gamma_\infty}$ equals $\beta_{\infty} (\fkc')$. On the other hand, by \cite[Lem.~5.1]{Kot90}, we have $\Delta_\infty = \beta_{\infty}(\omega)$. This gives our claim up to the sign of $\beta_\infty (\omega)$, and it remains to check that $\Delta_\infty$ is $2$-torsion. Let $\tilde T$ denote the inverse image of $T$ in $G_{\sconn}$. Then $\Int(n) (\mu_h) - \mu_h$ lies in the image of $X_*(\tilde T) \hookrightarrow X_*(T)$, and so $\Delta_\infty$ lies in the image of $X_*(\tilde T)_{\Gamma_\infty}$. The last group is $2$-torsion, since $\tilde T$ is anisotropic over $\RR$. This finishes the proof of our claim for $v=\infty$. 
 	
 	We have proved the claim. Now for each place $v$, since $\beta_v(\omega)$ maps to zero in $\pi_1(G)_{\Gamma_v}$, there exists $\tilde \beta_v(\omega) \in K$ lifting  $\beta_v(\omega)$. By the claim, the choices of lifts $\tilde \beta_v(\fkc)$ and $\tilde \beta_v(\fkc')$ in \S \ref{defn of Kottwitz invariant} can be arranged to satisfy that $\tilde \beta_v(\fkc) - \tilde \beta_v(\omega)$ maps to $\tilde \beta_v(\fkc')$ under $f: \pi_1(I_0) \isom \pi_1(I_0')$. To complete the proof, it remains to show that the element
 	$$\Omega: =  \sum_v \tilde \beta_v(\omega) \in K$$ is sent to zero under
 	$ K \to  K_{\Gamma} = K_{\Gamma,\tors} \to \E(I_0, G ; \A/\Q).$
 	
 	In fact, we show that there is a way to choose $\tilde \beta_v(\omega)$ such that the image of $\Omega$ in $K_{\Gamma}$ is already zero. Pick an element $\Omega' \in \coh^0_{\ab}(\Q, I_0 \to G)$ whose image along the surjection $\coh^0_{\ab}(\Q, I_0 \to G) \to \E(I_0, G ; \Q)$ is $\ab^1_{\QQ}(\omega)$. Write $(\Omega'_v)_v$ for the image of $\Omega'$ under the composite map (see Proposition \ref{prop:apply TN})
 	\begin{align}\label{eq:ab Q}
 		\coh^0_{\ab}(\Q, I_0 \to G) \to \coh^0_{\ab}(\A, I_0\to G) \cong \bigoplus_v \cA_{\Q_v} (K) \subset \bigoplus_v K_{\Gamma_v, \tors}.
 	\end{align} Then for each $v$, the element $\Omega'_v \in K_{\Gamma_v, \tors}$ is a lift of $\beta_v(\omega) \in \pi_1(I_0)_{\Gamma_v}$. Thus we may and shall choose the lift $\tilde \beta_v(\omega) \in K$ of $\beta_v(\omega)$ such that $\tilde \beta_v(\omega)$ is a lift of $\Omega_v'$. In this case, to show that $\Omega$ is sent to zero in $K_{\Gamma}$, it suffices to note that the composition of (\ref{eq:ab Q}) with the natural map
 	$  \bigoplus_v K_{\Gamma_v, \tors} \to K_{\Gamma,\tors}$ is zero. Indeed, by Proposition \ref{prop:apply TN},  this composition is identified with the composition
 	$$  \coh^0_{\ab}(\Q, I_0 \to G) \To \coh^0_{\ab}(\A, I_0\to G)  \To \coh^0_{\ab} (\A/\Q, I_0 \to G),$$ which is zero as desired.
 \end{proof}
 
 \subsection{} \label{local inner forms}  Let $\fkc = (\gamma_0, a , [b])\in \KP$, and write $I_0$ for $I_0(\fkc)$. Assume that $[b] \in \B(I_{0,\Q_p})$ is basic. This is the case, for example, when $\fkc$ is $p^n$-admissible, by Corollary \ref{cor:basic for adm Kott para}. Recall that the notion of inner forms and isomorphisms between inner forms are given in Definition \ref{defn:inner form}. The image of $a$ in $\coh^1(\A_f^p, I_0^{\ad})$ determines an inner form $I_v$ of $I_{0}$ over $\Q_v$ up to isomorphism for each place $v \neq p, \infty$. The basic element $[b] \in \B(I_{0,\Q_p})$ determines an inner form $I_p$ of $I_0$ over $\Q_p$ up to isomorphism, namely $I_p : = J_b^{I_0}$ for a decent representative $b$ of $[b]$ (see \S \ref{para:J_b}). Finally, let $(T,h)$ be as in \S \ref{defn of Kottwitz invariant}. Then $\Int(h(i)) $ induces a Cartan involution on $(I_0 / Z_G)_{\R}$, from which we obtain an inner form $I_{\infty}$ of $I_{0}$ over $\RR$ that is anisotropic modulo $Z_{G,\R}$. The isomorphism class of this inner form $I_\infty$ depends only on $\fkc$.
 
 \begin{prop}\label{existence of I}
 	In the situation of \S \ref{local inner forms}, assume that the Kottwitz invariant $\alpha(\fkc)$ is zero. Then there exists an inner form $I = I(\fkc)$\index[Oct2026-sorted-notation]{i 0 00 00 i c@$I(\mathfrak c)$} of $I_0$ over $\QQ$, unique up to isomorphism between inner forms, such that its localization over $\QQ_v$ is isomorphic to $I_v$ as inner forms of $I_{0,\QQ_v}$ for each place $v$ . 
 \end{prop}

 \begin{proof}
 	The uniqueness follows from the Hasse principle for $I_0^{\ad}$. To prove the existence, for each place $v$ we denote by ${\eta}_v$ the cohomology class in $\coh^1(\QQ_v, I_0^{\ad})$ corresponding to the inner form $I_v$ (cf.~Remark \ref{rem:inner form}). By \cite[Prop.~2.6]{Kot86} (cf.~\cite[Thm.~5.16]{borovoi}) we have an exact sequence of pointed sets
 	\begin{align*}
 		\coh^1(\QQ, I_0^{\ad}) \to \bigoplus_v \coh^1(\QQ_v, I_0^{\ad}) \xrightarrow{m} \pi_1(I_0^{\ad})_{\Gamma, \tors}.
 	\end{align*}
 	Here $m$ is defined as follows. For each place $v$, let $m_v$ be the composite
 	\begin{align*}
 		\coh^1(\QQ_v, I_0^{\ad}) \xrightarrow{\ab^1} \coh^1_{\ab}(\QQ_v, I_0^{\ad}) \cong \cA_{\Q_v}(\pi_1(I_0^{\ad})) \hookrightarrow \pi_1(I_0^{\ad})_{\Gamma_v, \tors},
 	\end{align*} (note that these maps are all isomorphisms for $v$ finite) 
 	and let $i_v$ be the natural map
 	$ \pi_1(I_0^{\ad})_{\Gamma_v, \tors} \to \pi_1(I_0^{\ad})_{\Gamma, \tors} .$ Then $m : = \sum_v i_v \circ m_v$.
 	
 	We only need to prove that $$\sum_v i_v \circ m_v( \eta_v) = 0. $$  For each $v$, we claim that $m_v(\eta_v)$ equals the image of $\beta_v(\fkc)$ under $\pi_1(I_0)_{\Gamma_v} \to \pi_1(I_0^{\ad}) _{\Gamma_v} = \pi_1(I_0^{\ad})_{\Gamma_v, \tors}$. Indeed, this statement is non-trivial only for $v \in \set{p,\infty}$. For $v =p$, there is a canonical bijection between $\coh^1(\Q_p, I_0^{\ad})$ and the set of basic elements of $\B(I_{0,\Q_p}^{\ad})$. If we identify $\coh^1(\Q_p, I_0^{\ad})$ with $\pi_1(I_0^{\ad})_{\Gamma_p} = \pi_1(I_0^{\ad})_{\Gamma_p, \tors}$, then this bijection is a section of the Kottwitz map $ \B(I_{0,\Q_p}^{\ad}) \to \pi_1(I_0^{\ad})_{\Gamma_p}$. Moreover, this bijection sends $\eta_p$ to the image of $[b]$ in $\B(I_{0,\Q_p}^{\ad})$. For more details see the end of \cite[\S 2.1]{RV14}. The claim for $v = p$ follows. For $v = \infty$, let $(T,h)$ be as in \S \ref{defn of Kottwitz invariant}, and let $\bar T: = T/Z_{G,\RR}$. Let $\bar h$ (resp.~$\bar \mu_{h}$) be the composition of $h: \mathbb S \to T$ (resp.~$\mu_h: \GG_m \to T_{\CC}$) with $T \to \bar T$. Since $\bar T$ is anisotropic, we have $\tatecoh^{-1}(\Gamma_\infty, X_*(\bar T)) = X_*(\bar T)_{\Gamma_\infty}$, and the Tate--Nakayama isomorphism $\tatecoh^{-1}(\Gamma_\infty, X_*(\bar T)) \isom \coh^1(\RR, \bar T)$ is induced by the map
 	\begin{align*}
 		X_*(\bar T) & \To Z^1(\RR, \bar T) \\
 		\lambda & \longmapsto (1\mapsto 1, \tau \mapsto \lambda(-1)),
 	\end{align*}
 	where $\tau$ denotes the complex conjugation. By definition, $\eta_\infty$ is represented by the cocycle $(1\mapsto 1, \tau \mapsto \bar h(i))$, whereas $\beta_{\infty}(\fkc)$ is represented by $\mu_h \in X_*(T)$. Thus to verify our claim it suffices to check that $ \bar h(i) = \bar \mu_{h}(-1)$. This follows from the equality $h(i) = \mu_h(-1) w_h(i)$, where $w_h$ is the weight cocharacter of $h$ and factors through $Z_{G,\RR}$.
 	
 	Write $K$ for $\ker( \pi_1(I_0) \to \pi_1(G ))$. 
 	By the above claim, $\sum_v i_v \circ m_v( \eta_v)$ is equal to the image of $\alpha(\fkc)$ under the composite map
 	\begin{align*}
 		\E(I_0,G;
 		\A/\Q)  \cong &\frac{K_{\Gamma,\tors}                }{ \bigoplus_v \ker( K_{\Gamma_v,\tors} \to \pi_1(I_0)_{\Gamma_v               })} \to \pi_1(I_0)_{\Gamma, \tors}  \to \pi_1(I_0^{\ad})_{\Gamma,\tors},
 	\end{align*}
 	where the first isomorphism is as in Proposition \ref{computing E}, and the second map is induced by the inclusion $K \hookrightarrow \pi_1(I_0)$. Since $\alpha(\fkc) =0$, we have $ \sum_v i_v \circ m_v( \eta_v) =0$, as desired.
 \end{proof}

 \section{Stating the point counting formula} 
 \label{subsec:stating pcf}  
 In this section we state the point counting formula. This is a formula expressing the quantity $T(\Phi_{\fkp}^m, f^pdg^p)$ defined in \S \ref{sec:SV} in terms of Kottwitz parameters and orbital and twisted orbital integrals. 
 \subsection{}
 \label{subsubsec:setting for point count}
 Let $m \in \ZZ_{ \geq 1} $. We now state a conjectural formula for $T(\Phi_{\fkp}^m, f^pdg^p)$ defined in \S \ref{para:unramified SD}. Let $ n=mr$. Fix $\mathfrak c = (\gamma_0, a, [b]) \in \KP_{\adm}(p^n)$, satisfying $\alpha(\mathfrak c) = 0 $. As in \S \ref{from CK to K}, $\fkc$ gives rise to a classical Kottwitz parameter $(\gamma_0, \gamma, \delta) \in \KP_{\clsc}(p^n)$ of degree $n$, well defined up to equivalence. Let $I(\fkc)$ be the global inner form of $I_0(\fkc)$ as in Proposition \ref{existence of I}. Let
 $ R : = \Res_{\QQ_{p^n}/\QQ_p} G$, and we view $\delta$ as an element of $R(\Q_p)$. Let $\theta$ be the $\Q_p$-automorphism of $R$ corresponding to the arithmetic $p$-Frobenius $\sigma \in \Gal(\QQ_{p^n}/\QQ_p)$. Let $R_{\delta\rtimes \theta}$\index[Oct2026-sorted-notation]{r 0 00 00 r subscript delta rtimes theta@$R_{\delta\rtimes \theta}$} denote the fixed subgroup of $R$ under the automorphism $\Int(\delta) \circ \theta$.
 
 Note that the $\A_f^p$-group $G_{\gamma}^0$ is isomorphic to $I(\fkc)_{\A_f^p}$, and the $\Q_p$-group $R_{\delta\rtimes \theta}^0$ is isomorphic to $I(\fkc)_{\Q_p}$. Moreover, these isomorphisms are canonical up to inner automorphisms defined over $\A_f^p$ and $\Q_p$ respectively. Choose Haar measures $di^p$ on $I(\fkc)(\A_f^p)$ and $di_p$ on $I(\fkc)(\QQ_p)$. They can be transported to  $G_{\gamma}^0 (\A_f^p)$ and  $R_{\delta\rtimes \theta}^0 (\QQ_p)$ respectively in an unambiguous way. We denote the resulting Haar measures on $G_{\gamma}^0 (\A_f^p)$ and  $R_{\delta\rtimes \theta}^0 (\QQ_p)$ still by $di^p$ and $di_p$.

 Since $r |n$ and $G$ is quasi-split over $\Q_p$, the Hodge cocharacters $\mu_h$ of $h\in X$ determine a $G(\Q_{p^n})$-conjugacy class of cocharacters of $G_{\Q_{p^n}}$, cf.~\cite[\S 1.3]{kottwitztwisted}. The \emph{negative} of this conjugacy class (i.e., with all members replaced by their inverses) further determines a $\cG(\Z_{p^n})$-double coset in $G(\Q_{p^n})$ via the Cartan decomposition, and we denote the characteristic function of this double coset by\index[Oct2026-sorted-notation]{zz21 phi 1 00 00 phi subscript n@$\phi_n$} $\phi_n : G(\Q_{p^n}) \to \set{0,1}$ (cf.~\cite[\S 2.1]{kottwitztwisted}). 
 We define \index[Oct2026-sorted-notation]{o 0 00 00 o c m f superscript p dg superscript p di subscript p di superscript p@$O(\mathfrak c, m , f^p dg^p , di_p di^p)$}
 \begin{align*}
 	O (\mathfrak c, m , f^p dg^p , di_p di^p) :=  O_{\gamma} (f^pdg^p)  TO_{\delta} (\phi_n) \in \CC, \end{align*}
 where $O_{\gamma} (f^p dg^p)$ is the orbital integral\index{orbital integral}
 $$ \int\limits_{G_{\gamma} ^0 (\A_f^p) \backslash G(\A_f^p) } f^p(g^{-1} \gamma g ) \frac{dg^p}{di^p}, $$  and $TO_{\delta} (\phi_n)$ is the twisted orbital integral\index{orbital integral!twisted} $$  \int\limits_{R_{\delta\rtimes \theta}^0 (\QQ_p)  \backslash R(\QQ_p)} \phi_n( r^{-1 } \delta \theta (r) ) \frac{ dr_p} {di_p},$$ with $dr_p$ the Haar measure on $R(\Q_p) = G(\Q_{p^n})$ giving volume $1$ to $\cG(\ZZ_{p^n})$. 
 
 \begin{rem}\label{rem:rational orbital integral} 
 	As the notation suggests, the dependence of $O(\fkc, m , f^p dg^p, di_p di^p)$ on the two Haar measures $di_p$ and $di^p$ is only via the product measure $di_p di^p$ on $I(\fkc)(\A_f)$. Moreover, if $di_pdi^p$ is rational on compact open subgroups, then $O(\fkc, m , f^p dg^p, di_p di^p)$ lies in $\QQ$. To see this, we may assume that $f^p = 1_{K^p a K^p}$ for a compact open subgroup $K^p \subset G(\A_f^p)$ and some $a\in G(\A_f^p)$, and that each of $dg^p, di_p, di^p$ is rational on compact open subgroups. We then know that $O_\gamma(f^pdg^p)$ lies in $\QQ$ by \cite[(3.4.1)]{Kot05} and the variant discussed below it (noting that $G_{\gamma}^0(\A_f^p) \backslash G_{\gamma}(\A_f^p)$ is compact). We also have $TO_{\delta}(\phi_n) \in \QQ$ by a formula similar to \textit{loc.~cit.}, cf.~the proof of \cite[Lem.~4.2.3]{ZZ}.
 \end{rem}
 \begin{rem}
 	Up to normalizations of the Haar measures on $G_{\gamma}^0(\A_f^p)$ and $R_{\delta \rtimes \theta} ^0 (\Q_p)$, the dependence of $O(\fkc, m , f^pdg^p, di_p di^p)$ on $\fkc$ is only via $(\gamma,\delta)$. However, for later purposes, it is important to normalize these Haar measures by choosing a Haar measure on $I(\fkc)(\A_f)$. Note that the classical Kottwitz parameter $(\gamma_0,\gamma,\delta)$ alone does not determine the global inner form $I(\fkc)$ of $G_{\gamma_0}^0$ (unless $G_{\der}$ is simply connected), so this way of normalization only makes sense with the presence of $\fkc$. 
 \end{rem}
 \begin{lem}\label{lem:closed subgp}
 	Let $I$ be a reductive  group over $\Q$, and let $Z$ be a $\Q$-subgroup of $Z_I$. Assume that $I/Z$ is anisotropic over $\RR$. Then for any open subgroup $U \subset Z(\A_f)$, the subgroup $I(\Q)U$ of $I(\A_f)$ is closed. 
 \end{lem}
 \begin{proof} We write $\bar I$ for $I/Z$.
 	Since $\bar I (\RR)$ is compact, $\bar I(\Q)$ is discrete in $\bar I(\A_f)$. Let $f$ denote the map $I(\A_f) \to \bar I(\A_f)$, and set $V : = f^{-1} (\bar I(\Q) )$. Then $Z(\A_f) = f^{-1}(\set{1})$ is open in $V$. Hence $I(\Q)U$ is an open subgroup of $V$, and  in particular closed in $V$. But $V$ is closed in $I(\A_f)$ since the discrete subgroup $\bar I(\QQ) \subset \bar I(\A_f)$ is closed. Hence $I(\Q) U$ is closed in $I(\A_f)$.
 \end{proof}
 
 \subsection{}\label{para:finite c_1} Keep the setting of \S \ref{subsubsec:setting for point count}. For each compact open subgroup $K^p\subset G(\A_f^p)$ we write $Z_{K_pK^p}$\index[Oct2026-sorted-notation]{z 0 00 00 z subscript k@$Z_K$} for $Z_G(\A_f) \cap K_p K^p$, and write $Z(\Q)_{K_pK^p}$\index[Oct2026-sorted-notation]{z 0 00 00 z q subscript k@$Z(\mathbb Q)_K$} for $Z_G(\Q) \cap K_pK^p$. Since $Z: = Z_G$ and $I : = I(\fkc)$ satisfy the assumptions in Lemma \ref{lem:closed subgp}, we know that $I(\fkc)(\Q)Z_{K_pK^p}$ is a closed subgroup of $I(\fkc)(\A_f)$. Hence  $I(\fkc)(\Q) Z_{K_pK^p} \backslash I(\fkc)(\A_f)$ is Hausdorff. It is also compact, since $I(\fkc)(\Q) \backslash I(\fkc)(\A_f)/U$ is finite for every compact open subgroup $U \subset I(\A_f)$ by \cite[Thm.~5.1]{Bor63}.  We equip $I(\fkc)(\A_f)$ with the Haar measure $di_p di^p$, and equip $I(\fkc)(\Q)  Z_{K_pK^p}$ with the Haar measure that gives volume $1$ to its open subgroup $Z_{K_pK^p}$. Then $I(\fkc)(\Q) Z_{K_pK^p}\backslash I(\fkc)(\A_f)$ has finite volume under the quotient measure, and we denote this volume by\index[Oct2026-sorted-notation]{c 1 00 00 c subscript 1 c k superscript p di subscript p di superscript p@$c_1(\mathfrak c, K^p, di_p di^p)$} $$c_1(\fkc, K^p, di_p di^p). $$
 
 Note that the product $$ c_1(\fkc , K^p, di_p di^p) O(\fkc, m , f^p dg^p, di_p di^p)$$ is independent of $di_p di^p$; this is because if we multiply $di_p di^p$ by a scalar $\lambda$ then the $c_1$ term gets multiplied by $\lambda$ and the $O$ term gets multiplied by $\lambda^{-1}$. Combined with Remark \ref{rem:rational orbital integral} and the observation that $c_1(\fkc , K^p, di_p di^p)  \in \QQ$ if $di_p di^p$ is rational on compact open subgroups, this implies that the above product lies in $\QQ$. In the sequel, we shall denote this product simply by 
 $$ c_1(\fkc , K^p ) O(\fkc, m , f^p dg^p) \in \QQ .$$

 We also define\index[Oct2026-sorted-notation]{c 1 00 00 c subscript 2 c c subscript 2 gamma subscript 0@$c_2(\mathfrak c) = c_2 (\gamma_0)$}
 $$ c_2(\fkc) = c_2 (\gamma_0):= \abs{\Sha_G(\QQ, G_{\gamma_0}^0)}. $$
 \subsection{}\label{para:Frob KP} Let $\SigmaREllip(G)$\index[Oct2026-sorted-notation]{zz18 sigma 0 00 00 sigma subscript r minus ell@$\Sigma_{\RR \text{-} \el}$} be the set of stable conjugacy classes of semi-simple, $\RR$-elliptic elements of $G(\QQ)$. (This is well defined, since $\RR$-elliptic maximal tori transfer between inner forms of reductive groups over $\RR$.) We fix a compact open subgroup $K^p \subset G(\A_f^p)$ such that $K_p K^p$ is neat and such that $f^p$ is $K^p$-bi-invariant. We fix a subset $\Sigma_{K^p}$\index[Oct2026-sorted-notation]{zz18 sigma 0 00 00 sigma subscript k superscript p@$\Sigma_{K^p}$} of $G(\Q)$ such that each $Z(\QQ) _{K_pK^p}$-translation-orbit in $\Sigma_{\RR \text{-} \el}(G)$ is represented by exactly one element of $\Sigma_{K^p}$.
 
 For each $\gamma_0 \in \Sigma_{K^p}$, we write $\KP(\gamma_0)$ for the set of $\fkc \in \KP$ whose first component is $\gamma_0$. 
 
 For any reductive group $H$ over $\QQ$ and any $\epsilon \in H(\QQ)_{\semi}$, we know that $(H_{\epsilon}/ H_{\epsilon}^0)(\Qbar)$ is isomorphic to a subgroup of the abelian group $\pi_1(H_{\der})$ by \cite[Cor.~2.16 (a)]{steinbergtorsion}. It follows that $H_{\epsilon}/ H_{\epsilon}^0$ is a finite commutative algebraic group over $\QQ$. We define \index[Oct2026-sorted-notation]{zz09 iota 1 00 00 iota subscript h epsilon@$\iota_H(\epsilon)$}\index[Oct2026-sorted-notation]{zz09 iota 1 00 20 bar iota subscript h epsilon@$\bar\iota_H(\epsilon)$} $$
 \iota_H(\epsilon)  : = [H_{\epsilon} (\QQ) : H_{\epsilon} ^0 (\QQ)], \quad 
 \bar \iota _H (\epsilon) : = \abs{  (H_{\epsilon} / H_{\epsilon} ^0 )(\QQ)  }. $$

 \begin{conj}\label{conj:point counting formula} For all sufficiently large integers $m$ (in a way depending on $f^pdg^p$), we have
 	\begin{multline}\label{pcf}
 		T(\Phi_{\fkp}^m, f^p dg^p) \\  =  \sum_{\gamma_0\in \Sigma_{K^p}}  \bar \iota_G (\gamma_0) ^{-1}  c_2(\gamma_0)  \tr \xi (\gamma_0)  \sum_{\substack{\fkc \in \KP(\gamma_0) \cap \KP_{\adm}(p^n) \\ \alpha(\fkc) = 0 }} c_1 (\fkc,K^p) O(\fkc, m , f^pdg^p). \end{multline}
 \end{conj}
 \begin{rem}
 	On the right hand side of (\ref{pcf}), each summand indexed by $\gamma_0$ is of the form $\tr \xi (\gamma_0)$ (which lies in $\Qbar$) times a rational number, by the discussion at the end of \S \ref{para:finite c_1}. Hence the right hand side of (\ref{pcf}) in fact lies in the smallest number field containing $\set{\tr \xi (\gamma_0)\mid \gamma_0 \in \Sigma_{K^p}}$. 	\end{rem}
 \begin{rem}\label{rem:Deligne's conj} When $G_{\der}$ is simply connected and $Z_G$ is cuspidal, the right hand side of (\ref{pcf}) recovers the formula conjectured by Kottwitz in \cite[\S 3]{Kot90}. We have formulated the conjecture only for $m$ sufficiently large since this is anyway needed for the local terms in the Grothendieck--Lefschetz--Verdier formula to be the naive local terms (i.e., Deligne's conjecture).\index{Deligne's conjecture}  For applications, it is important (and usually sufficient) to know that (\ref{pcf}) holds for all sufficiently large $m$, not just all sufficiently divisible $m$.
 \end{rem}
 
 \cleardoublepage
 \chapter{Variants of the Langlands--Rapoport Conjecture} 
 \label{sec:LRtau}
 This chapter aims to formulate the Langlands--Rapoport--$\tau$ Conjecture (Conjecture \ref{hypo about LR}), a weaker version of the original Langlands--Rapoport Conjecture (Conjecture \ref{conj:LR}). We begin in \S \ref{subsec:gerbs}--\S \ref{sec:adm morph} by developing the necessary formalism for Galois gerbs (much of which is a review). In \S \ref{sec:int mod}, we review the integral models of Shimura varieties and the statement of the original conjecture. Finally, after providing specific foundations, we state the Langlands--Rapoport--$\tau$ Conjecture in \S \ref{subsec:preparation LR}. 

 \section{The formalism of Galois gerbs}\label{subsec:gerbs}
 The Langlands--Rapoport Conjecture, as originally formulated in \cite{langlands1987gerben}, is stated in the language of \emph{Galois gerbs}. In this section, we recall the foundational definitions of this formalism. Our treatment largely follows \cite[\S 3.1]{kisin2012modp}, though we incorporate certain corrections (see especially Remark \ref{rem:conn cond}) and provide additional clarifications.
 
 In the following, let $k'/k$ be a Galois extension of fields of characteristic zero.

 \begin{defn}\label{defn:gerb} By a\index{Galois gerb} \emph{$k'/k$-Galois gerb}, we mean a pair $(G, \G)$, where $G$ is a connected linear algebraic group over $k'$, and $\G$ is an extension of topological groups $$1 \To  (G(k'),\text{discrete topology}) \To  \G \To \Gal(k'/k)  \To 1  $$ satisfying the following conditions.
 	\begin{enumerate}
 		\item For each $g \in \G$, there is a  $k'$-group isomorphism $g^{\alg} : \tau^* G \to G$, where $\tau$ is the image of $g$ in  $\Gal(k'/k)$, such that the conjugation action of $g$ on $G(k')$ is given by $G(k') \xrightarrow{\tau} (\tau^*G)(k')\xrightarrow{g^{\alg}} G(k')$.
 		\item There exists a continuous group-theoretic section of $\G \to \Gal(k'/k)$ defined on an open subgroup of $\Gal(k'/k)$.
 	\end{enumerate} We often write $\G$ for a $k'/k$-Galois gerb $(G,\G)$, and write $\G^{\Delta}$\index[Oct2026-sorted-notation]{g 0 04 00 g superscript delta@$\mathfrak G^{\Delta}$} for $G$, called the\index{Galois gerb!kernel} \emph{kernel} of $\G$. 	A \emph{morphism} between two $k'/k$-Galois gerbs $\G_1$ and $\G_2$ is a pair\index[Oct2026-sorted-notation]{zz21 phi 1 00 00 phi superscript delta@$\phi^{\Delta}$} $(\phi^{\Delta}, \phi)$ consisting of a $k'$-homomorphism $\phi^{\Delta}: \G_1^{\Delta} \to \G_2^{\Delta}$ and a continuous homomorphism $\phi: \G_1 \to \G_2$ satisfying the following conditions. 
 	\begin{itemize}
 		\item[(a)]  $\phi$ commutes with the maps $\G_i \to \Gal(k'/k)$. 
 		\item[(b)] The restriction of $\phi$ to $\G_1^{\Delta}(k')$ is given by  $\phi^{\Delta}$.
 	\end{itemize} We write $\redgb(k'/k)$\index[Oct2026-sorted-notation]{g 0 03 00 g rb k prime k@$\mathcal{G}rb(k'/k)$} for the category of $k'/k$-Galois gerbs.
 \end{defn}
 
 \begin{rem}\label{rem:conn cond} The assumption in Definition \ref{defn:gerb} that $G$ is connected is missing in \cite[\S 3.1]{kisin2012modp}, and should be added. This assumption implies that $G(k')$ is Zariski dense in $G$ (see \cite[Cor.~18.3]{borel1991}). In particular, each $g\in \G$ uniquely determines the isomorphism $g^{\alg}$, which is vital for various constructions. Another consequence is that for a morphism  $(\phi^{\Delta}, \phi): \G_1 \to \G_2$ between $k'/k$-Galois gerbs, $\phi^{\Delta}$ is uniquely determined by $\phi$. For this reason we can view $\phi$ alone as a morphism $\G_1 \to \G_2$.  \end{rem}
 \begin{rem}\label{rem:topology}
 	Let $\G \in \redgb(k'/k)$. Then all continuous group-theoretic sections of $\G \to \Gal(k'/k)$ defined on open subgroups of $\Gal(k'/k)$ belong to the same germ, i.e., any two such sections agree with each other on a sufficiently small open subgroup. Moreover, this germ can be extended to a continuous set-theoretic section of $\G \to \Gal(k'/k)$ defined on the whole $\Gal(k'/k)$. 
 \end{rem}
 \begin{defn}\label{defn:neutral gerb}
 	Let $G$ be a connected linear algebraic group over $k$. Let $\G_G$\index[Oct2026-sorted-notation]{g 0 04 00 g subscript g@$\mathfrak G_G$} be the split extension $G(k')\rtimes \Gal(k'/k)$, where $\Gal(k'/k)$ acts naturally on $G(k')$. Then $(G_{k'}, \G_G) \in \redgb(k'/k)$, and it is called the\index{Galois gerb!neutral} \emph{neutral $k'/k$-Galois gerb} associated with $G$.
 \end{defn}
 \begin{defn}\label{defn:conj of morph}
 	Let $\G \in \redgb(k'/k)$. For each $g \in \G^{\Delta}(k')$, conjugation by $g$ induces an automorphism  $\Int(g)$ of $\G$. Let $\phi, \psi: \gH \to  \G$ be two morphisms in $\redgb(k'/k)$. We say that $\phi$ and $\psi$ are \emph{conjugate} (or \emph{$\G^{\Delta}(k')$-conjugate}, for clarity), if there exists $g\in \G^{\Delta}(k')$ such that $\phi = \Int(g) \circ \psi$.
 \end{defn}
 
 \subsection{}\label{para:underline isom}
 Let $\phi: \gH \to \G$ be a morphism in $\redgb(k'/k)$. As in \cite[\S 3.1.1, Lem.~3.1.2]{kisin2012modp}, associated with $\phi$ we have  a $k$-algebraic group\footnote{The assumption in \cite[Lem.~3.1.2]{kisin2012modp} that the target of $\phi$ is a neutral gerb is not needed.}  $I_{\phi}$\index[Oct2026-sorted-notation]{i 0 00 00 i subscript phi@$I_{\phi}$} such that $I_{\phi}(k)$ is identified with $\{g \in \G^{\Delta}(k') \mid \Int(g) \circ \phi = \phi\}$ and $I_{\phi,k'}$ is identified with the centralizer $\G^{\Delta}_{\phi^{\Delta}}$\index[Oct2026-sorted-notation]{g 0 04 00 g superscript delta subscript phi superscript delta@$\mathfrak G^{\Delta}_{\phi^{\Delta}}$} of $\im(\phi^{\Delta})$ in $\G^{\Delta}$. If we choose a continuous set-theoretic section $\Gal(k'/k) \to \gH, \tau \mapsto q_{\tau}$ of $\gH \to \Gal(k'/k)$, then the action of $\tau \in \Gal(k'/k)$ on $I_{\phi}(k') \cong \G^{\Delta}_{\phi^{\Delta}}(k')$ with respect to the $k$-form $I_{\phi}$ is induced by conjugation by $\phi(q_{\tau})$ inside $\gG$.

 If $\zeta: \mathfrak K \to \mathfrak H$ and $\phi : \gH \to \gG$ are two morphisms in $\redgb(k'/k)$, then the inclusion $ \G^{\Delta}_{\phi^{\Delta}} \hookrightarrow \G^{\Delta}_{(\phi \circ \zeta)^{\Delta}}$ induces an injective $k$-homomorphism
 \begin{align}\label{eq:enlarge I}
 	I_{\phi} \hookrightarrow I_{\phi \circ \zeta} .
 \end{align}
 
 Let $\phi:\gH \to \G_G$ be a morphism in $\redgb(k'/k)$ with $G$ a reductive group over $k$ and $\gH^{\Delta}$ a torus. Then $I_{\phi,k'} \hookrightarrow G_{k'}$ is a pre-inner twisting (Definition \ref{defn:inner transfer}) between the $k$-groups $I_\phi$ and $G$.  If $\gH^\Delta$ is a torus, the induced map $\pi_1(I_{\phi}) \to \pi_1(G)$ is surjective since $I_{\phi, k'} \cong G_{k', \phi^{\Delta}}$ contains a maximal torus in $G_{k'}$. In this case we can apply the construction in \S \ref{para:construction of tilde}  to obtain reductive groups\index[Oct2026-sorted-notation]{i 0 00 00 i subscript phi superscript dagger@$I_{\phi}^{\dagger}$}\index[Oct2026-sorted-notation]{i 0 00 23 tilde i subscript phi@$\tilde I_{\phi}$} $$\tilde I_\phi \quad  \text{and} \quad I_\phi^{\dagger}$$ over $k$ equipped with pre-inner twistings into $G_{\sconn}$ and $G_{\der}$ respectively, as well as canonical $k$-maps $\tilde I_\phi \to I_{\phi}^{\dagger} \hookrightarrow I_\phi$. 
 \begin{rem}\label{rem:kernel torus}
 	Let $\G \in \redgb(k'/k)$, and let $\phi$ be the identity $\G \to \G$. Then $I_{\phi}$ is a canonical $k$-form of $\G^{\Delta}_{\phi^{\Delta}} = Z_{\G^{\Delta}}$. In particular, if $\G^{\Delta}$ is a torus, then we have a canonical $k$-form of $\G^{\Delta}$.
 \end{rem}
 
 \subsection{} \label{subsubsec:base change of gerbs}
 Let $l'/l$ be another Galois extension of fields of characteristic zero, equipped with compatible embeddings $k\hookrightarrow l$ and $k'\hookrightarrow l'$. In this situation we have the\index{Galois gerb!pull-back} \emph{pull-back} functor\index[Oct2026-sorted-notation]{pb 0 01 00 pb@$\mathrm{PB}$}
 \begin{align}\label{eq:gbBC}
 	\gbBC: \redgb(k'/k)  & \To \redgb (l'/l).
 \end{align}
 We explain its definition.	
 
 We first define $\gbBC$ of an object. Let $\G \in \redgb(k'/k)$, with kernel $G$. We have a short exact sequence
 \begin{align*}
 	1 \To G(k') \To \G_{l'/l}^0 \To \Gal(l'/l) \To 1,\end{align*} where $\G_{l'/l}^0$ is the fiber product $\G \times _{\Gal(k'/k)} \Gal(l'/l)$ in the category of topological groups. The above short exact sequence is in fact an extension of topological groups, which follows easily from Remark \ref{rem:topology}.
 
 Given any $g \in \G_{l'/l}^0$ with image $\tau \in \Gal(l'/l)$, the conjugation action of $g$ on $G(k')$ is induced by a $k'$-isomorphism $(\tau|_{k'})^* G \isom G$ uniquely determined by $g$ (i.e., the isomorphism $h^{\alg}$, where $h$ is the image of $g$ in $\G$). The last $k'$-isomorphism induces an $l'$-isomorphism $u_g: \tau^* (G_{l'}) \isom G_{l'}$, and in particular an automorphism of $G(l')$ given by $$G(l') = G_{l'}(l') \overset{\tau}{\To} (\tau^* G_{l'}) (l') \overset{u_g}{\To} G_{l'}(l') = G(l') . $$ In this way we obtain an action of $\G_{l'/l}^0$ on $G(l')$ via group automorphisms.
 
 Let $\G_{l'/l}$ be the quotient group $(G(l') \rtimes \G_{l'/l}^0 )/  \set{ x \rtimes x^{-1} \mid x \in G(k')}$. We shall denote elements of $\G_{l'/l}$ by $[x, (u,\tau)]$, for $x \in G(l'), (u,\tau) \in \G_{l'/l}^0$. We have a short exact sequence
 \begin{align}\label{eq:new ext}
 	1 \To G(l') \xrightarrow{x \mapsto [x, (1,1)]}  \G_{l'/l} \xrightarrow{[x,(u,\tau)] \mapsto \tau} \Gal(l'/l) \To 1.\end{align}
 
 We equip $\G_{l'/l}$ with the quotient topology of the product topology on $G(l') \rtimes \G^0_{l'/l}$. (As always, $G(l')$ has the discrete topology.) Then $\G_{l'/l}$ is a topological group, and (\ref{eq:new ext}) is an extension of topological groups. One checks that  $(G_{l'},\G_{l'/l}) \in \redgb(l'/l)$. We define $\gbBC(\G)$ to be $(G_{l'}, \G_{l'/l})$.
 
 We now define $\gbBC$ of a morphism. Given any morphism $\phi: \G \to \gH$ in $\redgb(k'/k)$, we define $\gbBC(\phi)$ to be $ (\phi^{\Delta}, \phi_{l'/l}) : \gbBC(\gG) \to \gbBC(\gH)$, where $\phi_{l'/l}$ is given by 
 \begin{align*}
 	\phi_{l'/l}: \G_{l'/l} & \To \gH_{l'/l} \\
 	[x, (u,\tau)] & \longmapsto [\phi^{\Delta}(x), (\phi(u), \tau)], \quad x \in G(l'), (u,\tau) \in \G^0_{l'/l}.
 \end{align*} This concludes the definition of the functor $\gbBC$. 
 
 Note that for $\G \in \redgb(k'/k)$, there is a canonical group homomorphism \index[Oct2026-sorted-notation]{zz18 sigma 1 00 00 varsigma subscript can superscript g@$\varsigma_{\can}^{\G}$}
 \begin{align*}
 	\varsigma_{\can}^{\G} : \Gal(l'/lk') \To \G_{l'/l}  , \quad
 	\tau \longmapsto [1, (1, \tau)],
 \end{align*}
 which is a section of $\G_{l'/l} \to \Gal(l'/l)$ over $\Gal(l'/lk')$. 
 
 \begin{lem}\label{lem:condition for PB} Keep the setting of \S \ref{subsubsec:base change of gerbs}, and assume in addition that $l = k$. Let $\G, \mathfrak H \in \redgb(k'/k)$.
 	A morphism $\psi : \gbBC(\gG) \to \gbBC(\gH)$ is of the form $\gbBC(\phi)$ for some morphism $\phi : \G \to \mathfrak H$ if and only if $\psi \circ \varsigma_{\can}^{\G} = \varsigma_{\can}^{\mathfrak H}$. Moreover, when this is the case, $\phi$ is unique.
 \end{lem}
 \begin{proof}
 	The ``only if'' part is trivial. We show the ``if'' part.
 	
 	From the hypothesis on $\psi$, it follows that $\psi^{\Delta}$ is the base change to $l'$ of a $k'$-homomorphism $\phi^{\Delta}: \G^{\Delta} \to \gH^{\Delta}$, and that $\psi [1, (u,\tau)]$ is of the form
 	$[1, (\phi(u), \tau)]$ for some function $\phi : \G \to \mathfrak H$. We now check that $(\phi^{\Delta}, \phi)$ is a morphism in $\redgb(k'/k)$. In fact, only the continuity of $\phi$ is non-obvious. For this, we observe that the map $\gH^0_{l'/l} \to \gH^{\Delta}(l') \rtimes \gH^0_{l'/l}, (u,\tau) \mapsto (1, u,\tau)$ is continuous and open, since $\gH^{\Delta}(l')$ is discrete. Hence the induced injective map $\gH^0_{l'/l} \to \gH_{l'/l}$ is also continuous and open, i.e., a homeomorphism onto its image. It follows that the map $\G^0_{l'/l} \to \gH^0_{l'/l}, (u,\tau) \to (\phi(u),\tau)$ is continuous. The continuity of $\phi$ then follows from the openness of the map $\Gal(l'/l)  \to \Gal(k'/k)$.
 	
 	Given that $(\phi^{\Delta}, \phi)$ is a morphism, it is clear that $\psi=\gbBC(\phi)$.
 	
 	Finally, we show that if $\psi=\gbBC(\phi)$ then $\phi$ is unique. Suppose $\phi_1$ also satisfies the condition. Then
 	$$ \psi[1,(u,\tau)]  = [1, (\phi(u), \tau)]  = [1, (\phi_1(u),\tau)], $$ and in particular $\phi(u) = \phi_1(u)$, for all $(u,\tau) \in \G_{l'/l}^0$. Since the projection $\G_{l'/l}^0 \to \G, (u,\tau) \mapsto u$ is surjective, we have $\phi = \phi_1$.
 \end{proof}
 
 \begin{defn}\label{defn:progerb}	
 	Let $\proredgb(k'/k)$\index[Oct2026-sorted-notation]{pro 1 01 00 pro minus g rb k prime k@$\mathrm{pro}\text{-}\mathcal{G}rb(k'/k)$} be the category whose objects are projective systems $(\G_i)_{i\in I}$ in $\redgb(k'/k)$ indexed by directed sets $(I,\leq)$ and whose morphisms are given by
 	$$\Hom_{\proredgb(k'/k)}((\G_i)_{i\in I}, (\mathfrak H_j)_{j\in J}) : = \varprojlim_{j \in J} \varinjlim_{i\in I}  \Hom_{\redgb(k'/k)} (\G_i, \mathfrak H_j). $$ Objects of $\proredgb(k'/k)$ are called\index{Galois gerb!pro-} \emph{pro-$k'/k$-Galois gerbs}. We view $\redgb(k'/k)$ naturally as a full subcategory of $\proredgb(k'/k)$. When we are in the situation of \S \ref{subsubsec:base change of gerbs}, the pull-back functor (\ref{eq:gbBC}) naturally extends to a functor $\proredgb(k'/k) \to \proredgb(l'/l)$, which we still call \emph{pull-back}.
 \end{defn}

 \subsection{}\label{para:notation for progerb}
 Let $\G = (\G_i)_{i\in I} \in \proredgb(k'/k)$. Then we can take the projective limit\index[Oct2026-sorted-notation]{g 0 04 00 g superscript delta@$\mathfrak G^{\Delta}$} $$\G^{\Delta} : = \varprojlim_{i} \G_i^{\Delta}$$ in the category of affine $k'$-group schemes, and take the projective limit\index[Oct2026-sorted-notation]{g 0 04 00 g superscript top@$\mathfrak G^{\mathrm{top}}$} $$\G^{\topo}: = \varprojlim_{i} \G_i$$ in the category of topological groups.
 If $\phi : \G \to \mathfrak H$ is a morphism in $\proredgb(k'/k)$, then $\phi$ naturally induces a homomorphism of affine $k'$-group schemes\index[Oct2026-sorted-notation]{zz21 phi 1 00 00 phi superscript delta@$\phi^{\Delta}$} $\phi^{\Delta}: \G^{\Delta} \to \mathfrak H^{\Delta}$, and a continuous homomorphism\index[Oct2026-sorted-notation]{zz21 phi 1 00 00 phi superscript top@$\phi^{\mathrm{top}}$} $\phi^{\topo} : \G^{\topo} \to \mathfrak H^{\topo}$. In the sequel, if $d$ is an element of $\G^{\topo}$, we shall simply write $d\in \G$. Also, we shall simply write $\phi(d) \in \mathfrak H$ for $\phi^{\topo}(d) \in \mathfrak H^{\topo}$.

 The following definition generalizes Definition \ref{defn:conj of morph}.
 \begin{defn} Let $\gH = (\gH_i)_{i\in I} \in \proredgb(k'/k), \gG \in \redgb(k'/k)$, and let $\phi ,\psi  : \gH \to \G$ be two morphisms in $\proredgb(k'/k)$. We say that $\phi$ and $\psi$ are \emph{conjugate} (or \emph{$\G^{\Delta}(k')$-conjugate}) if there exists $g \in  \G^{\Delta}(k')$ such that $\psi = \Int(g) \circ \phi$ as morphisms in $\proredgb(k'/k)$.
 \end{defn}
 
 \subsection{}\label{para:underline isom in pro case} Let $\gH= (\gH_i)_{i\in I} \in \proredgb(k'/k), \G \in \redgb(k'/k)$, and let $\phi : \gH \to \G$ be a morphism in $\proredgb(k'/k)$. We now define $I_{\phi}$,\index[Oct2026-sorted-notation]{i 0 00 00 i subscript phi@$I_{\phi}$} generalizing the definition in \S \ref{para:underline isom}.
 
 Choose $i_0 \in I$ such that $\phi$ is induced by a morphism $\phi_{i_0} : \gH_{i_0} \to \G$. For each $i \in I$ with $i \geq i_0$, let $\phi_i$ be the composition $\gH_i \to \gH_{i_0} \xrightarrow{\phi_{i_0}} \G$. For $ j  \geq i \geq i_0$, there is a natural $k$-homomorphism $I_{\phi_i} \to I_{\phi_j}$ as in (\ref{eq:enlarge I}), whose base change to $k'$ is identified with the inclusion map $\G^{\Delta}_{\phi_i^{\Delta}} \to \G^{\Delta}_{\phi_j^{\Delta}}$. For sufficiently large $i$, the decreasing subgroups $\im(\phi_i^{\Delta})$ of $\G^{\Delta}$ stabilize, since $\G^{\Delta}$ is noetherian. Hence $\G^{\Delta}_{\phi_i^{\Delta}}$ and $I_{\phi_i}$ also stabilize. We can thus define
 \begin{align*}
 	\im(\phi^{\Delta}) & : =  \im(\phi_i^\Delta) , \\
 	\G^{\Delta}_{\phi^{\Delta}} & : = \G^{\Delta}_{\phi^{\Delta}_i} , \\
 	I_{\phi} & : = I_{\phi_i},
 \end{align*} for $i\in I$ sufficiently large. Clearly these definitions are independent of the initial choices of $i_0$ and $\phi_{i_0}$.
 
 By construction, $I_{\phi, k'}$ is canonically identified with $\G^{\Delta}_{\phi^{\Delta}}$. It is also easy to see that $I_{\phi}(k)$ precisely consists of those $g\in \G^{\Delta}(k')$ such that $\Int(g) \circ \phi = \phi$ (as morphisms in $\proredgb(k'/k)$).
 
 If $\G = \G_G$ for some reductive group $G$ over $k$, and if $\gH^{\Delta}_i$ are tori for all $i\in I$, then we also extend the definitions of $I_{\phi}^{\dagger}$\index[Oct2026-sorted-notation]{i 0 00 00 i subscript phi superscript dagger@$I_{\phi}^{\dagger}$}\index[Oct2026-sorted-notation]{i 0 00 23 tilde i subscript phi@$\tilde I_{\phi}$} and $\tilde I_{\phi}$ in \S \ref{para:underline isom} to the present case in the obvious way.
 \subsection{}\label{para:twist morphism} Let $\gH, \G, \phi$ be as at the beginning of \S \ref{para:underline isom in pro case}. Given a continuous 1-cocycle $a = (a_{\rho}) \in Z^1(k'/k, I_{\phi}(k'))$, there is a morphism\index[Oct2026-sorted-notation]{a 1 00 00 a phi@$a\phi$} $a\phi : \gH \to \G$ defined as follows. Choose $i\in I$ such that $\phi$ is induced by some $\phi_i : \gH_i \to \G$ and such that $I_{\phi} = I_{\phi_i}$. Denote by $\pi$ the structural map $\gH_i \to \Gal(k'/k)$. We define
 $$ a\phi_i : \gH_i \To \G, \quad g \longmapsto a_{\pi(g)} \phi_i(g). $$
 Here, we view $a_{\pi(g)} \in I_{\phi}(k')$ as an element of $\G^{\Delta}(k') \subset \G$ via the canonical embedding $I_{\phi,k'} \hookrightarrow  \G^{\Delta}$. Then $a\phi_i$ is a morphism in $\redgb(k'/k)$, and we define $a\phi$ to be the morphism induced by $a\phi_i$. This definition is independent of choices.
 
 \begin{lem}\label{lem:twist by cocycle}
 	In the setting of \S \ref{para:twist morphism}, the map $a \mapsto a \phi$ is a bijection from $Z^1(k'/k, I_{\phi}(k'))$ to the set of morphisms $\phi': \gH \to \G$ such that $\phi^{\prime, \Delta} = \phi^{\Delta}$. Moreover, for $a, a' \in Z^1(k'/k, I_{\phi}(k'))$, we have that $a\phi$ is conjugate to $a'\phi$ if and only if $a $ is cohomologous to $a'$.
 \end{lem}
 \begin{proof}
 	Since $\G^{\Delta}$ is finitely presented over $k'$, the natural map
 	$$ \varinjlim_{i} \Hom_{k'}(\gH_i^{\Delta} , \G^{\Delta}) \To \Hom _{k'} (\gH^{\Delta} , \G^{\Delta})$$ is a bijection. Here $\Hom_{k'}$ denotes the set of homomorphisms of $k'$-group schemes. The lemma then reduces to the case where $\gH$ is in $\redgb(k'/k)$. In this case the proof is exactly the same as the proof of part (2) of  \cite[Lem.~3.1.2]{kisin2012modp}. 
 \end{proof}
 
 \begin{defn}\label{defn:twist morphism}
 	For $a$ and $\phi$ as in \S \ref{para:twist morphism}, we call $a\phi$ the\index{twist of a gerb morphism} \emph{twist of $\phi$ by $a$}. Similarly, for a class $\beta \in \coh^1(k'/k, I_{\phi}(k'))$, we define the \emph{twist of $\phi$ by $\beta$}, denoted by $\phi^{\beta}$,\index[Oct2026-sorted-notation]{zz21 phi 1 00 00 phi superscript beta@$\phi^{\beta}$} to be the conjugacy class of $a \phi$ where $a$ is any cocycle representing $\beta$. This is well defined by Lemma \ref{lem:twist by cocycle}.\footnote{Since we define $\phi^{\beta}$ to be the whole conjugacy class, its members are not necessarily of the form $a\phi$ for any $a \in Z^1(k'/k, I_{\phi}(k'))$.} We shall sometimes also write $\phi^{\beta}$ for an unspecified member of this conjugacy class.
 \end{defn}
 
 \subsection{} 
 We explain how to obtain Galois gerbs from Reimann's explicit cocycle construction in \cite[App.~B]{reimann1997zeta}. We first sketch the idea behind the construction informally. Let $\G \in \redgb(k'/k)$. Suppose that there is a continuous set-theoretic section $\varsigma : \Gal(k'/k) \to \G$ of $\G \to \Gal(k'/k)$ such that $\varsigma(\rho) \varsigma(\tau) \varsigma(\rho \tau)^{-1}$ lies in $Z_{\G^{\Delta}}(k')$ for all $\rho ,\tau \in \Gal(k'/k)$. Then the isomorphisms $(\varsigma(\tau))^{\alg} : \tau^* \G^{\Delta} \isom \G^{\Delta}$ for $\tau \in \Gal(k'/k)$ form a $k'/k$-Galois descent datum. Let $G$ be the corresponding $k$-form of $\G^{\Delta}$. Then the isomorphism class of $\G$ can be recovered from $G$ and the map $\Gal(k'/k) \times \Gal(k'/k) \to Z_G(k'), (\rho, \tau) \mapsto \varsigma(\rho) \varsigma(\tau) \varsigma(\rho \tau)^{-1}$, which is in fact a continuous $2$-cocycle.
 
 We now give the formal construction. First we define a category $ \cR(k'/k)$\index[Oct2026-sorted-notation]{r 0 03 00 r k prime k@$\mathcal R(k'/k)$} as in \cite[pp.~106--107]{reimann1997zeta}. The objects are pairs $(G,z)$, where $G$ is a connected linear algebraic group over $k$, and $z = (z_{\rho,\tau})$ is a continuous $2$-cocycle $\Gal(k'/k) \times \Gal(k'/k) \to  Z_G(k')$ satisfying $z_{1,1} = 1$. A morphism $(G',z') \to (G,z)$ is a pair $(\phi^{\Delta}, f)$, where $\phi^{\Delta}: G'_{k'} \to G_{k'}$ is a homomorphism of $k'$-groups, and $f = (f_{\rho})$ is a continuous $1$-cochain $\Gal(k'/k) \to  G(k')$ satisfying
 \begin{align*}
 	z_{\rho,\tau} f_\rho \lix^{\rho} f_{\tau} f_{\rho\tau}^{-1} &  = \phi^{\Delta}(z'_{\rho,\tau}), \\
 	\Int(f_{\rho}) \circ \rho^* (\phi^{\Delta}) & = \phi^{\Delta},
 \end{align*}
 for all $\rho ,\tau \in \Gal(k'/k)$. (In particular, $f_1 =1$.) The composition of morphisms is given by
 $$(\phi^{\Delta}, f) \circ (\psi^{\Delta}, h) : =  (\phi^{\Delta} \circ \psi^{\Delta}, ( \phi^{\Delta}(h_{\rho}) f_{\rho}) _{\rho}). $$
 As explained in \textit{loc.~cit.}, $ \cR(k'/k)$ is equivalent to the category of ``affine algebraic $k'/k$-groupoids which admit a special section''.
 
 We define a fully faithful functor\index[Oct2026-sorted-notation]{e 0 03 00 e@$\mathcal E$}
 \begin{align}\label{eq:functor E}
 	\cE: \cR(k'/k) \To  \redgb(k'/k).	\end{align}
 Let $(G,z) \in  \cR(k'/k)$. To define $\cE(G,z)$, we let $\G : = G(k') \times \Gal(k'/k)$, equipped with the product topology. Define a binary operation on $\G$ by
 $$ (g ,\rho) \cdot (h, \tau) : = (g \rho(h) z_{\rho,\tau}, \rho \tau).$$ Then $(\G,\cdot)$ is a topological group, and the maps $G(k') \to \G, g \mapsto (g,1)$ and $\G \to \Gal(k'/k), (g,\rho) \mapsto \rho$ make $\G$ a topological group extension of $\Gal(k'/k)$ by $G(k')$. Then $(G_{k'}, \G) \in \redgb(k'/k)$. (To check condition (ii) in Definition \ref{defn:gerb}, use that $z_{\rho,\tau} =1$ for all $\rho,\tau$ sufficiently close to $1$.) We define $\cE(G,z)$ to be $(G_{k'}, \G)$.
 
 For a morphism $(\phi^{\Delta}, f): (G',z') \to (G,z)$ in $\cR(k'/k)$, we define
 \begin{align*}
 	\cE(\phi^{\Delta}, f) : \cE(G',z')  = G'(k') \times \Gal(k'/k) & \To \cE(G,z) = G(k') \times \Gal(k'/k) \\
 	(g',\rho) & \longmapsto (\phi^{\Delta}(g') f_\rho, \rho).
 \end{align*} This completes the definition of the functor $\cE$. We omit the proof that $\cE$ is fully faithful, since this fact will not be used.
 
 Analogously to Definition \ref{defn:progerb}, we consider the category $\proR(k'/k)$\index[Oct2026-sorted-notation]{pro 1 01 00 pro minus r k prime k@$\mathrm{pro}\text{-}\mathcal{R}(k'/k)$} of pro-objects in $\cR(k'/k)$ indexed by directed sets. The functor (\ref{eq:functor E}) naturally extends to a functor
 \begin{align}\label{eq:functor E pro}
 	\cE : \proR(k'/k) \To \proredgb(k'/k)\end{align}
 which is also fully faithful.
 
 In \cite[App.~B]{reimann1997zeta}, various \emph{affine groupoids} are defined, which are needed for the correct formulation of the Langlands--Rapoport Conjecture. These are constructed by constructing explicit objects in $\proR(k'/k)$ (for suitable $k'/k$) and applying a natural functor from $\proR(k'/k)$ to the category of affine $k'/k$-groupoids. In the present work, we shall not use the language of affine groupoids, but we shall import Reimann's explicit constructions and obtain pro-Galois gerbs via the functor (\ref{eq:functor E pro}). This point of view is also taken in \cite{kisin2012modp}.

 \section{The Dieudonn\'e gerb, quasi-motivic gerb, and pseudo-motivic gerb} The \emph{quasi-motivic} Galois gerb $\mathfrak Q$ plays a central role in the formulation of the Langlands--Rapoport Conjecture. It is equipped with a localization map $\mathfrak G_v \to \mathfrak Q(v)$ for each place $v$ of $\QQ$. The definition of $\mathfrak G_v$ is most involved for $v=p$, when it is closely related to the \emph{Dieudonn\'e gerb} $\mathfrak D$. In this section, we first discuss $\mathfrak D$ and $\mathfrak G_p$, as well as a construction that attaches to every morphism from $\mathfrak G_p$ to a neutral Galois gerb an \emph{affine Deligne--Lusztig set}. We then introduce $\mathfrak Q$ and various morphisms into and out of it. We conclude by introducing a variant of $\mathfrak Q$, the \emph{pseudo-motivic Galois gerb} $\mathfrak P$, which will play a special role in the study of Shimura varieties of Hodge type in \S \ref{sec:Shimura Hodge type} (see especially \S \ref{subsec:parallel}). 
 
 \subsection{}
 Fix a prime $p$. We recall the definition of the\index{Dieudonne gerb@Dieudonn\'e gerb} \emph{Dieudonn\'e gerb} in terms of the functor (\ref{eq:functor E pro}), cf.~\cite[pp.~109--110]{reimann1997zeta}. For each $n \in \Z_{\geq 1}$, let\index[Oct2026-sorted-notation]{zz10 kappa 1 00 00 kappa subscript n@$\kappa_n$} $\kappa_n : \Gal(\Q_p^{\ur}/\Q_p) \times \Gal(\Q_p^{\ur}/\Q_p) \to \ZZ$ be the unique continuous function satisfying
 $$ \kappa_n (\sigma^i , \sigma^j) =
 \lfloor i/n \rfloor + \lfloor j/n\rfloor - \lfloor (i+j)/n\rfloor, \quad \forall i,j \in \ZZ. $$ Here $\sigma$ denotes the arithmetic $p$-Frobenius as usual. Let $\mathcal D_n$\index[Oct2026-sorted-notation]{d 0 03 00 d subscript n@$\mathcal D_n$} be the object in $\cR(\Q_p^{\ur}/\Q_p)$ given by
 $(\GG_m, ( p^{\kappa_n(\rho,\tau)})_{\rho ,\tau } ).$ For $n , n
 ' \in \ZZ_{\geq 1}$ with $n|n'$, let\index[Oct2026-sorted-notation]{zz11 lambda 1 00 00 lambda subscript n n prime@$\lambda_{n,n'}$} $\lambda_{n,n'} : \Gal(\Q_p^{\ur}/\Q_p) \to \ZZ$ be the unique continuous function satisfying
 $$\lambda_{n,n'} (\sigma ^i) = \lfloor i/n' \rfloor n'/n - \lfloor i/n \rfloor, \quad \forall i \in \ZZ. $$ Let\index[Oct2026-sorted-notation]{zz04 delta 1 00 00 delta subscript n prime n@$\delta_{n',n}$} $\delta_{n',n} : \mathcal D_{n'} \to \mathcal D_n$ be the morphism in $\cR(\Q_p^{\ur} / \Q_p)$ given by $$(x \mapsto x^{n'/n} , (p^{\lambda_{n,n'}(\rho)})_{\rho}). $$ Then $(\mathcal D_n)_{n\in \ZZ_{\geq 1}}$ equipped with the transition morphisms $\delta_{n',n}$ is an object in $\proR(\Q_p^{\ur}/\Q_p)$.
 
 Applying the functor (\ref{eq:functor E pro}) to $(\mathcal D_n)_{n}$, we obtain an object\index[Oct2026-sorted-notation]{d 0 04 00 d d subscript n subscript n@$\mathfrak D = (\mathfrak D_n)_n$} $\D= (\D_n)_n = (\mathcal E(\mathcal D_n))_n $ in $\proredgb(\Q_p^{\ur}/\Q_p)$. This is called the \emph{Dieudonn\'e gerb}.

 We denote by\index[Oct2026-sorted-notation]{d 0 06 00 d@$\mathbb D$} $\mathbb D$ the pro-torus $\varprojlim_{n \in \ZZ_{\geq 1}} \GG_{m}$ over $\Spec \ZZ$, where  for $n|n'$ the transition map from the $n'$-th $\GG_m$ to the $n$-th $\GG_m$ is $z \mapsto z^{n'/n}$. We have $$\D^{\Delta} = \mathbb D_{\Q_p^{\ur}}. $$
 
 Let $n \in \ZZ_{\geq 1}.$ By construction, the set underlying $\D_n$ is $  (\Q_p^{\ur}) ^\times  \times \Gal(\Q_p^{\ur}/\Q_p)$. Using this we define a canonical element \index[Oct2026-sorted-notation]{d 1 00 00 d subscript sigma n@$d_{\sigma,n}$} $$d_{\sigma,n} := (p^{-\lfloor 1/n \rfloor},\sigma) \in \D_n,$$ and a canonical map\index[Oct2026-sorted-notation]{zz18 sigma 1 00 00 varsigma subscript n@$\varsigma_n$} $$\varsigma_{n} : \Gal(\Q_p^{\ur}/\Q_{p^n}) \To \D_n , \quad \rho \longmapsto (1,\rho), $$ which is a continuous section of $\D_n \to \Gal(\Q_p^{\ur}/\Q_p)$ over $\Gal(\Q_p^{\ur}/\Q_{p^n})$. Clearly $\kappa_n(\sigma^i , \sigma^j) = 0$ when $i,j$ are divisible by $n$, so  $\kappa_n$ vanishes on $\Gal(\Q_p^{\ur}/\Q_{p^n}) \times \Gal(\Q_p^{\ur}/\Q_{p^n})$ by continuity. It follows that $\varsigma_n$ is a group homomorphism. For future use, we compute:
 \begin{align}
 	\label{eq:power of d_sigma}
 	& d_{\sigma, n}^n 
 	= \Bigl( p^{-\lfloor 1/n \rfloor} \prod_{i = 1}^{n-1} p^{\kappa_{n}(\sigma^i, \sigma)}, \sigma^n \Bigr) = (p^{-1}, \sigma^n) = p^{-1} \varsigma_{n}(\sigma^n), \\
 	\label{eq:comp of d_sigma}
 	& \mathcal{E}(\delta_{n',n})(d_{\sigma, n'}) 
 	= \bigl( p^{-\frac{n'}{n}\lfloor 1/n' \rfloor +  \lambda_{n,n'}(\sigma)}, \sigma \bigr) = (p^{-\lfloor 1/n \rfloor}, \sigma) = d_{\sigma, n}, \\
 	\label{eq:comp of varsigma_n}
 	& \mathcal{E}(\delta_{n',n})(\varsigma_{n'}(\rho)) 
 	= (p^{\lambda_{n,n'}(\rho)}, \rho) = (1, \rho) = \varsigma_{n}(\rho),
 \end{align}
 for all $n \in \NN$, all multiples $n'$ of $n$, and all $\rho \in \Gal(\Q_p^{\ur}/\Q_{p^{n'}})$. By (\ref{eq:comp of d_sigma}), the system $(d_{\sigma , n})_n$ defines an element\index[Oct2026-sorted-notation]{d 1 00 00 d subscript sigma@$d_{\sigma}$} $d_{\sigma} \in \D^{\topo}$.
 
 \begin{defn}\label{defn:G_p}
 	Let $\G_p \in \proredgb(\Qpbar/\Q_p)$\index[Oct2026-sorted-notation]{g 0 04 00 g subscript p g subscript p n subscript n@$\mathfrak G_p = (\mathfrak G_{p,n})_n$} be the pull back of $\D$. For each $n \in \ZZ_{\geq 1}$, let $\G_{p,n} \in \proredgb(\Qpbar/\Q_p)$ be the pull back of $\D_n$.\footnote{In \cite[\S 3.1.6]{kisin2012modp}, our $\G_{p,n}$ is denoted by $\tilde {\G}_p^{\QQ_{p^n}}$.} Thus $\G_p = (\G_{p,n})_n$.   
 \end{defn}
 \begin{defn}\label{defn:unram morph}
 	Let $G$ be a connected linear algebraic group over $\Q_p$. Let $\G_G \in \redgb(\Qpbar/\Q_p)$ and $\G_G^{\ur} \in \redgb(\Q_p^{\ur}/\Q_p)$\index[Oct2026-sorted-notation]{g 0 04 00 g subscript g superscript ur@$\mathfrak G_G^{\ur}$} be the associated neutral gerbs. A morphism $\theta: \G_{p} \to \G_G$ in $\proredgb(\Qpbar/\Q_p)$ is called\index{unramified!morphism of Galois gerbs} \emph{unramified}, if it is the pull-back of a morphism $\theta^{\ur}: \D \to \G_G^{\ur}$ in $\proredgb(\QQ_p^{\ur}/\QQ_p)$. By the obvious generalization of Lemma \ref{lem:condition for PB}, $\theta^{\ur}$ is uniquely determined by $\theta$. For general $\theta$, we write $\mathcal{UR}(\theta)$\index[Oct2026-sorted-notation]{ur 0 03 00 ur theta@$\mathcal{UR}(\theta)$} for the set of $g\in G(\Qpbar)$ such that $\Int(g^{-1}) \circ \theta$ is unramified.
 \end{defn}
 \begin{lem}\label{lem:unram criterion}
 	Keep the notation of Definition \ref{defn:unram morph}. The following statements hold.
 	\begin{enumerate}
 		\item For any morphism $\theta: \G_p \to \G_G$, the set $\mathcal{UR}(\theta)$ is a $G(\QQ_p^{\ur})$-torsor, where $G(\QQ_p^{\ur})$ multiplies on the right.
 		\item Let $\phi: \D \to \G_G^{\ur}$ be a morphism. For sufficiently divisible $n$, we have $\phi(d_{\sigma})^n = \phi_n^{\Delta}(p^{-1}) \rtimes \sigma^n \in \G_G^{\ur}$, where $\phi_n$ is a morphism $\D_n\to \G_G^{\ur}$ inducing $\phi$.
 	\end{enumerate}
 \end{lem}
 \begin{proof}By the discussion in \S \ref{subsubsec:base change of gerbs}, the fact that $\G_p$ is the pull-back of $\D$ gives rise to a canonical homomorphism $\varsigma: \Gamma_{p,0} = \Gal(\Qpbar/\Qpur) \to \G_p^{\topo}$, which is a section of $\G_p^{\topo} \to \Gamma_p$ over $\Gamma_{p,0}$. By (the obvious generalization of) Lemma \ref{lem:condition for PB}, a morphism $\theta : \G_p \to \G_G$ is unramified if and only if $\theta(\varsigma(\rho)) = 1 \rtimes \rho$ for all $\rho \in \Gamma_{p,0}$.
 	
 	For part (i), we write $\theta(\varsigma(\rho)) = a_{\rho} \rtimes \rho \in \G_G$, for $\rho \in \Gamma_{p,0}$. Then $(a_{\rho})_{\rho}$ is a cocycle in $ Z^1(\Qpbar/\Q_p^{\ur}, G(\Qpbar))$. By the previous paragraph, an element $g\in G(\Qpbar)$ lies in $\mathcal{UR}(\theta)$ if and only if $g^{-1} a_{\rho} \rho(g) = 1$ for all $\rho \in \Gamma_{p,0}$. By Steinberg's theorem, such a $g$ exists since $G$ is connected. It is then clear that $\mathcal{UR}(\theta)$ is a $G(\Q_p^{\ur})$-torsor.
 	
 	We now show part (ii). First pick $n$ such that $\phi$ is induced by a morphism $\phi_n : \D_n \to \G_G^{\ur}$. Since $\phi_n$ is continuous, there exists an open subgroup $U$ of $\Gal(\Q_p^{\ur}/\Q_{p^n})$ such that $\phi_n (\varsigma_n (\rho)) = 1 \rtimes \rho$ for all $\rho \in U$. We may assume that $U = \Gal(\Q_p^{\ur}/\Q_{p^{n'}})$ for some $n'$ divisible by $n$. Using (\ref{eq:comp of varsigma_n}), we may replace $n$ by $n'$. We then have
 	$$ \phi(d_{\sigma})^n = \phi_n (d_{\sigma,n}^n) = \phi_n(p^{-1} \varsigma_n(\sigma^n)) = \phi_n^{\Delta}(p^{-1}) \rtimes \sigma^n,$$ where the second equality is by (\ref{eq:power of d_sigma}).
 \end{proof}
 \begin{defn}\label{defn:b_theta} Let $G$ be a connected linear algebraic group over $\Q_p$. For any unramified morphism $\theta: \G_{p} \to \G_G$, we define\index[Oct2026-sorted-notation]{b 1 00 00 b subscript theta@$b_{\theta}$} $b_{\theta} \in G(\Q_p^{\ur})$ by the formula $\theta^{\ur} (d_{\sigma}) = b_{\theta} \rtimes \sigma$.
 \end{defn}
 \begin{prop}\label{prop:twist at p} Let $G$ be a reductive group over $\Q_p$, and let $\theta : \G_p \to \G_G$ be an unramified morphism. The following statements hold.
 	\begin{enumerate}
 		\item Viewing $\theta ^{\ur, \Delta}: \mathbb D_{\Q_p^{\ur}} \to G_{\Q_p^{\ur}}$ as a fractional cocharacter of $G_{\Qpur}$, we have $\theta ^{\ur, \Delta} = - \nu_{b_{\theta}}$. Moreover $b_{\theta}$ is decent (see \S \ref{isocdefns}).
 		\item There are natural $\QQ_p$-isomorphisms $J_{b_{\theta}} \cong I_{\theta^{\ur}} \cong I_{\theta}$.
 		\item Let $\beta\in \coh^1(\QQ_p, I_{\theta})$. Then there is a member $\theta'$ of the conjugacy class $\theta^{\beta}$ (see Definition \ref{defn:twist morphism}) satisfying the following conditions:
 		\begin{enumerate}
 			\item The morphism $\theta' : \G_p \to \G_G$ is unramified.
 			\item By part (ii), we view $\beta$ as an element of $\coh^1(\Q_p, J_{b_{\theta}})$. Then the $\sigma$-conjugacy class of $b_{\theta'}$ in $G(\Q_p^{\ur})$ is given by the twist of $b_{\theta}$ by $\beta$, as in \S \ref{para:setting for delta(b,beta)}.
 	\end{enumerate}	\end{enumerate}
 \end{prop}
 \begin{proof}
 	\textbf{(i)} Choose $n\in \ZZ_{ \geq 1}$ such that $\theta^{\ur}$ is induced by a morphism  $\theta_n: \D_n \to \G_G^{\ur}$. Let $\nu_n : = \theta_n ^{\Delta}$. Then $\nu_n$ is a cocharacter of $G_{\Q_p^{\ur}}$, and $\theta^{\ur,\Delta} = n^{-1}\nu_n$. Up to enlarging $n$, we may assume that $\nu_n$ is defined over $\QQ_{p^n}$.
 	
 	By Lemma \ref{lem:unram criterion} (ii), up to enlarging $n$ we have
 	\begin{align*}
 		b_{\theta} \sigma(b_{\theta} ) \cdots \sigma^{n-1} (b_{\theta}) = \nu_n (p^{-1}) \in G(\QQ_p^{\ur}).
 	\end{align*}
 	We conclude that $\nu_{b_{\theta}} = -n^{-1}\nu_n = - \theta^{\ur, \Delta}$, and that $b_{\theta}$ is $n$-decent.
 	
 	\textbf{(ii)}  As in \cite[Cor.~1.14]{RZ96}, the fact that $b_\theta$ is decent implies that
 	$$J_{b_{\theta}}(R) = \set { g \in G_{\Q_p^{\ur}, \nu_{b_{\theta}}}(R\otimes _{\QQ_p} \QQ_p ^{\ur} )\mid g b_\theta = b _{\theta} \sigma(g)},$$ for any $\QQ_p$-algebra $R$. In view of $G_{\Qpur , \nu _{b_\theta}}  =  G_{\Q_p^{\ur}, \theta^{\ur, \Delta}}$, the above description of $J_{b_{\theta}}$  agrees with the explicit description of $I_{\theta^{\ur}}$ as in \S \ref{para:underline isom}. The natural isomorphism $I_{\theta^{\ur}} \cong I_{\theta}$ arises from the fact that they are the same $\Q_p$-form of $G_{\Q_p^{\ur}, \theta^{\ur, \Delta}}$. 
 	
 	\textbf{(iii)} By Steinberg's theorem, $\beta$ is represented by a cocycle $$a =(a_{\rho})\in Z^1(\QQ_p^{\ur} / \QQ_p, J_{b_{\theta}} (\QQ_{p}^{\ur})).$$ Viewing $a$ as an element in $Z^1(\Q_p, I_{\theta})$, we define $\theta': = a \theta$, which is in the conjugacy class $\theta^{\beta}$. Then $\theta'$ is unramified, and we have  $b_{\theta'} = a_{\sigma} b_{\theta}$, which implies condition (b).
 \end{proof}

 \subsection{} \label{subsubsec:recall of defns in LR} Let $\cG$ be a reductive group scheme over $\ZZ_p$, with generic fiber $G$. Let $r \in \ZZ_{\geq 1}$, and let $\mu$ be a cocharacter of $\cG_{\ZZ_{p^r}}.$
 
 For each morphism $\theta: \G_p \to \G_G$ in $\proredgb(\Qpbar/\Q_p)$, we recall the definition of a set $X_\mu(\theta)$ and a bijection $\Phi : X_\mu (\theta) \to X_\mu (\theta)$ as in \cite[\S 3.3.3]{kisin2012modp}. For each $g\in \mathcal {UR}(\theta)$, we write $b_g$ for $b_{\Int(g^{-1}) \circ \theta} \in G(\Qpur)$. Define\index[Oct2026-sorted-notation]{x 0 00 23 tilde x subscript mu theta@$\tilde X_\mu(\theta)$} $$\tilde X_\mu(\theta) : = \set{  g\in \mathcal {UR}(\theta) \mid b_g\in \mathcal G(\ZZ_p^{\ur}) p^{\mu} \mathcal G(\ZZ_p^{\ur})  \subset G(\QQ_p^{\ur})}.$$ The group $ \mathcal G(\ZZ_p^{\ur})$ acts on $\tilde X_\mu(\theta)$ via right multiplication, and we define\index[Oct2026-sorted-notation]{x 0 00 00 x subscript mu theta@$X_\mu(\theta)$}
 $$X_\mu(\theta) : = \tilde X_\mu(\theta)/ \mathcal G(\ZZ_p^{\ur}). $$ 
 By the Cartan decomposition, the subset $\cG(\Z_p^{\ur}) p^{\mu} \cG(\Z_p^{\ur}) \subset G(\Qpur)$ depends on $\mu$ only via the $\cG(\Z_{p^r})$-conjugacy class of $\mu$, and so the same holds for $\tilde X_{\mu}(\theta)$ and $X_\mu(\theta)$.
 
 Consider the map 
 \begin{align*}
 	\tilde \Phi : \mathcal {UR}(\theta) & \To \mathcal {UR}(\theta) \\
 	g & \longmapsto gb_{g} \sigma(b_{g}) \cdots \sigma^{r-1}(b_{g}).
 \end{align*}
 For $g\in \mathcal {UR} (\theta)$ we have $b_{\tilde \Phi(g)} = \sigma^r(b_g)$. It follows that $\tilde \Phi$ restricts to a map $\tilde X_\mu (\theta) \to \tilde X_{\mu}(\theta)$. 
 If we fix an element $g_0 \in \mathcal {UR}(\theta)$ and use it to identify the $G(\Qpur)$-torsor $\mathcal {UR}(\theta)$ with $G(\Qpur)$, then  $\tilde \Phi$ becomes the map $G(\Qpur) \to G(\Qpur), g \mapsto b_{g_0} \sigma(b_{g_0}) \cdots \sigma^{r-1}(b_{g_0}) \sigma^r(g)$. This shows that $\tilde \Phi : \tilde X_\mu (\theta) \to \tilde X_{\mu}(\theta)$ descends to a bijection\index[Oct2026-sorted-notation]{zz21 phi 0 00 00 phi frobenius@$\Phi$ (Frobenius)} $$\Phi: X_\mu(\theta) \isom X_\mu(\theta),$$ which we call the  \emph{$p^r$-Frobenius}.
 
 The isomorphism class of the $\Phi^{\ZZ}$-set $X_\mu(\theta)$ depends on $\theta$ only through the conjugacy class of $\theta$. Moreover, after fixing an element $g_0\in \mathcal{UR}(\theta)$, from the previous paragraph we know that the map $G(\Qpbar) \to G(\Qpbar), g \mapsto g_0^{-1} g$ induces a bijection from $X_\mu (\theta)$ to the\index{affine Deligne--Lusztig set} \emph{affine Deligne--Lusztig set}\index[Oct2026-sorted-notation]{x 0 00 00 x subscript mu b@$X_\mu(b)$}
 \begin{align*}  
 	X_{\mu}(b_{g_0}) & : = 
 	\set{g\in G(\Qpur)/ \mathcal G(\ZZ_p^{\ur}) \mid  g^{-1} b_{g_0} \sigma(g) \in \mathcal G(\ZZ_p^{\ur}) p^\mu \mathcal G(\ZZ_p^{\ur})}
 	\\  & \isom  \set{g\in G(\LL)/ \mathcal G(\breve \ZZ_p) \mid  g^{-1} b_{g_0} \sigma(g) \in \mathcal G(\breve \ZZ_p) p^\mu \mathcal G(\breve \ZZ_p)} ,
 \end{align*} on which $\Phi$ acts by $g \mapsto b_{g_0} \sigma(b_{g_0}) \cdots \sigma^{r-1}(b_{g_0}) \sigma^r(g)$. Here note that we have a natural bijection  $G(\Qpur)/ \mathcal G(\ZZ_p^{\ur}) \cong G(\LL)/ \mathcal G(\breve \ZZ_p) $ by the density of $G(\Qpur)$ in $G(\LL)$ for the $p$-adic topology, and that the Cartan decomposition is functorial when passing from $\Qpur$ to $\LL$.

 \subsection{}\label{subsubsec:unique-quasi-motivic}	We keep fixing a prime $p$. For each finite prime $v \neq p$, let\index[Oct2026-sorted-notation]{g 0 04 00 g subscript v v neq p@$\mathfrak G_{v}, ~v\neq p$} $$\G_v : = \Gamma_v \in \redgb(\overline \Q_v/\Q_v). $$ Let $$\G_{\infty} \in \redgb(\CC/\RR)$$ be the Weil group of $\RR$, and let $$\G_p \in \proredgb(\Qpbar/\Q_p)$$ be as in Definition \ref{defn:G_p}.

 Consider the pro-torus\index[Oct2026-sorted-notation]{r 0 05 00 r@$\mathscr R$} $\mR : = \varprojlim_{L}\Res_{L/\QQ} \GG_m$ over $\QQ$, where $L$ runs through the set of finite Galois extensions of $\QQ$ contained in $\overline \QQ$, ordered by inclusion, and the transition maps are the norm maps. The corresponding projective system of neutral gerbs $\G _{\Res_{L/\QQ} \GG_m } \in \redgb(\Qbar/\QQ)$ forms an object in $\proredgb(\Qbar/\QQ)$, which we denote by\index[Oct2026-sorted-notation]{g 0 04 00 g subscript r@$\mathfrak G_{\mathscr R}$} $\G_{\mR}$.\footnote{\label{ft:mR}In \cite{kisin2012modp}, the objects denoted here by $\mR$ and $\G_{\mR}$ appear 
 	as $R_{\ol\QQ/\QQ} \GG_m$ and $\G_{R_{\ol\QQ/\QQ} \GG_m}$, respectively. }
 
 For each object $\G$ (resp.~morphism $\phi$) in $\proredgb(\Qbar/\QQ)$, we denote its pull-back in $\proredgb(\ol \Q_v/ \Q_v)$ by $\G(v)$\index[Oct2026-sorted-notation]{g 0 04 00 g v@$\mathfrak G(v)$} (resp.~$\phi(v)$\index[Oct2026-sorted-notation]{zz21 phi 1 00 00 phi v@$\phi(v)$}), for each place $v$ of $\Q$. Here the pull-back functor is defined with respect to our fixed embedding $\overline \Q\to \overline \Q_v$.  Reimann \cite[\S B.2]{reimann1997zeta} has constructed a\index{quasi-motivic Galois gerb} \emph{quasi-motivic Galois gerb}, which is an object\index[Oct2026-sorted-notation]{q 0 04 00 q@$\mathfrak Q$} $$\gQ  \in \proredgb(\Qbar/\Q), $$ equipped with a morphism\index[Oct2026-sorted-notation]{zz06 zeta 1 00 00 zeta subscript v g subscript v to q v@$\zeta_v:\mathfrak G_v \to \mathfrak Q(v)$} $\zeta_v : \G_v \to \Qf(v)$ in $\proredgb(\ol \Q_v/\Q_v)$ for each place $v$ of $\QQ$, and a morphism\index[Oct2026-sorted-notation]{zz23 psi 1 00 00 psi q to g subscript r@$\psi:\mathfrak Q \to \mathfrak G_{\mathscr R}$} $\psi : \gQ \to \G_{\mR}$ in $\proredgb(\Qbar/\Q)$. More precisely, Reimann constructs in the proof of \cite[Thm.~B.2.8]{reimann1997zeta} versions of $\gQ$, $\zeta_v$, and $\psi$ in the categories $\proR(\overline \Q/\Q)$ and $\proR(\overline \Q_v/ \Q_v)$. We transport his constructions via the functors (\ref{eq:functor E pro}). The unique characterization of the tuple $(\gQ , (\zeta_v)_v, \psi )$ is delicate to state. We omit this and refer the reader to \textit{loc.~cit.~}and \cite[Thm.~3.1.9]{kisin2012modp}.
 
 By construction, $\gQ$ is given by a projective system $(\gQ^L)_L$\index[Oct2026-sorted-notation]{q 0 04 00 q superscript l@$\mathfrak Q^L$} in $\redgb(\Qbar/\QQ)$, indexed by the set of finite Galois extensions $L/\QQ$ contained in $\overline {\QQ}$ ordered by inclusion. For each $L$, we have $\gQ^{L,\Delta} = Q^L_{\overline \QQ}$, where $Q^L$\index[Oct2026-sorted-notation]{q 0 00 00 q superscript l@$Q^L$} is a $\QQ$-torus explicitly described in \cite[\S B.2]{reimann1997zeta}. (Here $Q^L$ agrees with the canonical $\QQ$-form of $\gQ^{L,\Delta}$ as in Remark \ref{rem:kernel torus}.) If $L\subset L' \subset \overline \QQ$, then the transition morphism $\gQ^{L'} \to \gQ^L$ is surjective, and the induced map on kernels $ Q^{L'}_{\overline \Q} \to Q^L_{\overline \Q}$ is defined over $\Q$. We write $Q$ for the pro-torus $
 (\varprojlim_{L} Q^L)$ over $\QQ$. Since the projective system $(\gQ^L)_L$ is indexed by a countable set and since the transition morphisms are surjective, we conclude that the projections $\gQ^{\topo} \to \gQ^L$ are surjective. (See \S \ref{para:notation for progerb} for the notation $\Qf^{\topo}$.)
 
 For $v \in \set{p,\infty},$ we denote the group scheme homomorphism $\zeta_v^{\Delta} : \G_v^{\Delta} \to \Qf(v)^{\Delta}$ (see \S \ref{para:notation for progerb}) by $\nu(v)$.\index[Oct2026-sorted-notation]{zz13 nu 1 00 00 nu v@$\nu(v)$} By construction, $\nu(v)$ is defined over $\Q_v$. Thus we have $\nu(p): \mathbb D_{\Q_p} \to  Q_{\QQ_p}$ and $\nu(\infty): \GG_{m,\RR} \to  Q_{\RR}$.

 \subsection{}\label{para:Psi_mu}
 Let $T$ be a torus over $\Q$, and let $\mu \in X_*(T)$. Let $L/\Q$ be a finite Galois extension contained in $\Qbar$ such that $\mu$ is defined over $L$. Then $\mu$ induces a $\Q$-homomorphism $$\mu_*: \Res_{L/\Q} \GG_m \xrightarrow{\Res_{L/\Q} \mu} \Res_{L/\Q} T \xrightarrow{\N_{L/\Q}}  T . $$ We obtain a morphism\index[Oct2026-sorted-notation]{zz23 psi 0 00 00 psi subscript t mu@$\Psi_{T,\mu}$} $\Psi_{T,\mu} : \Qf \to \G_T$ in $\proredgb(\Qbar/\Q)$ via the composition
 $$ \Qf \overset{\psi}{\To} \G_{\mR} \To \G_{\Res_{L/\QQ} \GG_m} \overset{\mu_*}{\To} \G_T. $$ This is independent of the choice of $L$.
 
 \begin{lem}\label{lem:b in torus}Let $\theta = \Psi_{T,\mu}(p) \circ \zeta_p : \G_p \to \G_T(p)$. Choose $g \in \mathcal {UR}(\theta)$, and let $[b] \in \B(T_{\Q_p})$ be the $\sigma$-conjugacy class of $b_{\Int(g^{-1}) \circ \theta} \in T(\Qpur)$ in $T(\LL)$ (which is well defined, by Lemma \ref{lem:unram criterion} (i)). Then $\kappa_T([b])\in X_*(T)_{\Gamma_p}$ is equal to the image of $-\mu\in X_*(T)$.
 \end{lem}
 \begin{proof} The proof reduces to the ``universal case'', where $T = \Res_{L/\Q}\GG_m$, and the map $\mu_*: \Res_{L/\Q}\GG_m \to T$ is the identity. In this case, $X_*(T)_{\Gamma_p} $ is torsion free, and therefore the homomorphism $\iota: X_*(T)_{\Gamma_p} \to X_*(T) \otimes \QQ$ induced by taking averages of $\Gamma_p$-orbits in $X_*(T)$ is injective. For each $[b'] \in \B(T_{\Q_p})$, we have $\iota(\kappa_T([b'])) =\nu_{b'}$ by \cite[\S 2.8]{kottwitzisocrystal}. Hence to prove the lemma it suffices to prove that $\iota(\mu) = -\nu_b$. By Proposition \ref{prop:twist at p} (i), we have $-\nu_b = (\Int(g^{-1})\circ \theta)^{\ur, \Delta} =  \theta^{\Delta}$. By  \cite[(3.1.11)]{kisin2012modp}, $\theta^{\Delta}$ is indeed equal to $\iota(\mu)$, as desired.\footnote{A similar argument is made in the proof of \cite[Lem.~3.4.2]{kisin2012modp} in order to determine $\kappa_T([b])$ for $T = \Res_{L/\Q} \GG_m$. There $T$ is unramified over $\Q_p$, so the cited result \cite[Thm.~4.2 (ii)]{rapoportrichartz} is valid. Our present argument does not need this unramifiedness assumption.}
 \end{proof}
 
 \subsection{}\label{subsec:pseudomotivic} We recall from \cite[\S 3]{langlands1987gerben} and \cite[\S 3.5.5]{kisin2012modp} the construction 
 of the Weil pro-torus $P,$ which is a quotient of $Q.$ It has the following description. 
 Let $L$ be a CM field (contained in $\Qbar$) which is Galois over $\Q$, $m \geq 1,$ and $q = p^m.$ A\index{Weil number} \emph{$q$-Weil number} is a $p$-unit $\pi \in L$ such that 
 \begin{itemize}
 	\item For some $\iota \in \mathbb Z$, we have $|\pi|_v = q^{\iota/2}$ for all archimedean (i.e.~complex) places of $L$. 
 	\item For every place $v|p$ of $L,$ $\log_q |\prod_{\tau \in \Gal(L_v/\Q_p)} \tau(\pi)|_p \in \mathbb Z.$
 \end{itemize}
 Let $P^L(m)$ denote the torus over $\QQ$ whose character group consists of $q$-Weil numbers in $L$ modulo the subgroup 
 of roots of unity in $L.$ For a $q$-Weil number $\pi \in L$, we denote by\index[Oct2026-sorted-notation]{zz22 chi 1 00 00 chi subscript pi@$\chi_{\pi}$} $\chi_{\pi} \in X^*(P^L(m))$ the corresponding character. 
 For $m | m'$ there is a map $P^L(m') \rightarrow P^L(m)$ dual to $\pi \mapsto \pi^{m'/m}.$
 By \cite[Cor.~3.8]{langlands1987gerben}, this map is an isomorphism when $m$ is sufficiently divisible and $m'$ is an arbitrary multiple of $m$. We define $P^L$\index[Oct2026-sorted-notation]{p 0 00 00 p superscript l@$P^L$} to be $P^L(m)$ for such $m$.\footnote{This result is a consequence of   \cite[Lem.~3.7]{langlands1987gerben}, and the proof of the latter relies on a construction assigning to any torus $T$ split over $L$ and any cocharacter $\mu$ a homomorphism $\psi_{T,\mu}: Q(L,m) \to T$, or dually, a map $X^*(T) \to X^*(Q(L,m))$. Here $X^*(Q(L,m))$ is defined using ``quasi-Weil numbers'', and $m$ is a sufficiently divisible integer. This construction is incorrect, because the given formula for the map $X^*(T) \to X^*(Q(L,m))$ does not land in $X^*(Q(L,m))$ in general---this is exactly the well-known mistake which Reimann \cite{reimann1997zeta} corrected by updating the definition of $Q(L,m)$. Nevertheless, when $L$ is CM and $\mu$ satisfies the Serre condition, the formula of Langlands--Rapoport correctly defines a map $X^*(T) \to X^*(Q(L,m))$ and it factors through $X^*(P(L,m))  \subset X^*(Q(L,m))$ (as claimed by them in the Serre case). In the proof of \cite[Lem.~3.7]{langlands1987gerben}, only the Serre case of the construction is used, so the use of $Q$ can be deleted throughout, and the proof is valid.} 
 We then define\index[Oct2026-sorted-notation]{p 0 00 00 p varprojlim subscript l p superscript l@$P=\varprojlim_L P^L$}
 \[
 P = \varprojlim_L P^L,
 \]
 where $L$ runs over CM subfields of $\Qbar$ that are Galois over $\Q$, and the transition maps are induced by the natural inclusions of $q$-Weil numbers for $L \subset L'$. Because these inclusions yield injective homomorphisms 
 $ X^*(P^L(m)) \to  X^*(P^{L'}(m))$, 
 the transition morphisms $P^{L'} \to P^L$ are all surjective.
 \subsection{}
 
 For each finite Galois $L/\Q$ in $\Qbar$ which contains a CM field, there is a natural surjective $\Q$-morphism $$Q^L \To P^{L_{\mathrm{CM}}},$$ where $L_{\mathrm{CM}}$ is the maximal CM subfield of $L$; see \cite[Lem.~B2.2, B2.3]{reimann1997zeta}, \cite[Cor.~3.4]{langlands1987gerben}. These morphisms are compatible for varying $L$ and define a morphism of pro-tori $Q \to P$. Thus we can push out $\gQ$ along $Q \to P$ and obtain the\index{pseudo-motivic Galois gerb} \emph{pseudo-motivic Galois gerb}\index[Oct2026-sorted-notation]{p 0 04 00 p@$\mathfrak P$} $$\gP : =  \big (\text{push-out of } \gQ^L \text{ along } Q^L \to P^{L_{\mathrm{CM}}} \big)_L  \in \proredgb(\Qbar/\Q). $$ We have $\gP^{\Delta} = P_{\Qbar}$. By construction, if $\gG \in \redgb(\Qbar/\Q)$ and $\phi : \gQ \to \gG$ is a morphism such that $\phi^{\Delta}: Q_{\Qbar} \to \gG^{\Delta}$ factors through $P_{\Qbar}$, then $\phi$ factors uniquely through $\gP$. Comparing with \cite[Rmk.~B2.11]{reimann1997zeta}, we see that the above construction of $\gP$ is equivalent to that considered in \textit{loc.~cit.~}and \cite{langlands1987gerben}.   
 
 \subsection{}\label{subsec:deltan} 
 We again denote by $\nu(p)$ the morphism $\DD_{\Q_p} \rightarrow P_{\Q_p}$ induced by $\nu(p) : \DD_{\QQ_p} \to Q_{\QQ_p}.$  
 For $L$ a CM field Galois over $\QQ$ contained in $\Qbar$, we write $\nu(p)^L$ for the induced morphism $\DD_{\Q_p} \rightarrow P^L_{\Q_p}.$  
 We will write\index[Oct2026-sorted-notation]{zz13 nu 1 00 00 nu p superscript l subscript n@$\nu (p)^L_n$} $\nu(p)^L_n \in X_*(P^L)\otimes \Q$ for $n \cdot \nu(p)^L$ where $\nu(p)^L$ is thought of as an element of $ X_*(P^L)\otimes \Q$. If $v$ denotes the place of $L$ above $p$ determined by the fixed embedding $\Qbar \to \Qpbar$ and $n=[L_v:\Q_p]$, then we have $\nu(p)^L_n \in X_*(P^L)$ and moreover it has the following characterization (see \cite[(3.5.2)]{kisin2012modp}): For every Weil $q$-number $\pi \in L$ defining $\chi_\pi \in X^*(P^L)$, we have\footnote{Note that the right hand side can in general change if $v$ is replaced by a different place of $L$ above $p$. Indeed the definition of the quasi-motivic gerb relies on fixing the embedding $\Qbar \to \Qpbar$.}
 $$ \langle \nu(p)^L_n, \chi_\pi \rangle = \log_q \big|\prod_{\tau \in \Gal(L_v/\Q_p)} \tau(\pi)\big|_p .$$   
 
 Recall \cite[p.~142]{langlands1987gerben}, \cite[Cor.~3.5.4, \S 4.3.9]{kisin2012modp} that for fixed $L$ as above and for each positive integer $n$ divisible by $m_L \cdot |\mu(L)|$, where $m_L$ is such that $P^L = P^L(m_L)$ and $\mu(L)$ is the group of roots of unity in $L$, there exists a unique element\index[Oct2026-sorted-notation]{zz04 delta 1 00 00 delta subscript n@$\delta_n$} $\delta_n \in P^L(\Q)$ 
 with the following properties:
 \begin{enumerate} 
 	\item If $\pi \in L$ is a $q = p^m$-Weil number with $m \cdot | \mu(L)|$ dividing $n$ and $P^L = P^L(m)$, then $\chi_{\pi}(\delta_n) = \pi^{n/m}.$
 	\item For $n|n',$ we have $\delta_{n'} = \delta_{n}^{n'/n}$. 
 	\item For CM fields $L \subset L' \subset \ol \Q$ Galois over $\QQ$, the map $P^{L'} \rightarrow P^{L}$ preserves the elements $\delta_n$ (for sufficiently divisible $n$).
 	\item For $n$ sufficiently divisible, we have $\nu_{\delta_n} = -\nu(p)^L_n \in X_*(P^L),$ 
 	where $\nu_{\delta_n}$ denotes the Newton cocharacter of $\delta_n.$
 \end{enumerate}
 \section{Strictly monoidal categories}
 In order to extend the definition of admissible morphisms $\phi : \mathfrak Q \to \mathfrak G_G$ in \cite{langlands1987gerben} to groups $G$ without assuming that $G_{\der}$ is simply connected, one needs the formalism of strictly monoidal categories as in \cite{kisin2012modp}. In the more general setting, the standard morphism $\mathfrak G_G \to \mathfrak G_{G^{\ab}}$ is replaced by a strict monoidal functor $\G_G \to \G_{G/G_{\sconn}}$. In this section, we review this formalism and introduce the notion of a \emph{continuous} monoidal isomorphism. The continuity condition provides a necessary correction to the definition of admissible morphisms in \textit{loc.~cit.~}(see Definition \ref{defn of adm morphism} below). 
 \subsection{}
 Let $G, H$ be two strictly monoidal  categories (which we always assume to be small). By a\index{strictly monoidal functor} \emph{strictly monoidal functor} $G \to H$, we mean a functor that strictly respects the monoidal structures. By a\index{monoidal isomorphism} \emph{monoidal isomorphism} between two strictly monoidal functors $\phi, \psi : G \to H$, we mean an isomorphism of functors $\cA : \phi \isom \psi$ such that for any two objects $g_1, g_2 \in G$ the following diagram commutes:
 $$ \xymatrixcolsep{4pc} \xymatrix{ \phi(g_1 \otimes g_2) \ar[r]^{\cA (g_1 \otimes g_2)} \ar@{=}[d] & \psi(g_1 \otimes g_2) \ar@{=}[d]  \\ \phi(g_1) \otimes \phi(g_2) \ar[r]^{\cA(g_1) \otimes \cA(g_2)} & \psi(g_1 ) \otimes \psi(g_2) } $$
 
 Every group can be naturally viewed as a strictly monoidal category, where the only morphisms are the identities. For two groups $G$ and $H$, the set of strictly monoidal functors $G\to H$ is the same as the set of group homomorphisms $G \to H$. More generally, each crossed module $(\tilde H \to H)$ also determines a strictly monoidal category denoted by $H/\tilde H$; see \cite[\S 3.2.1]{kisin2012modp}.  When $(\tilde H \to H)$ is a crossed module, each element $h \in H$ induces via conjugation a strictly monoidal functor $\Int(h): H/\tilde H \to H/\tilde H$.
 
 \subsection{} \label{para:cocycle interpretation}
 Now let $G$ be a group and $(\tilde H \xrightarrow{\varrho} H)$ be a crossed module. We denote the structural action of $H$ on $\tilde H$ by $\Int$. Consider two strictly monoidal functors $\phi, \psi  : G \to H/\tilde H$. If $\cA: \phi  \isom  \psi$ is an isomorphism of functors, then for each $q \in G$ the isomorphism $\cA(q) : \phi(q) \isom \psi(q)$ corresponds to an element $\cA(q) \in \tilde H$. Thus we may view $\cA$ as a function $G \to \tilde H$. In this way, there is a one-to-one correspondence between monoidal isomorphisms $\cA: \phi  \isom  \psi$ and functions $\cA : G \to \tilde H$ satisfying
 \begin{align*}
 	\cA(qr) & = \cA(q) \cdot \Int(\phi(q)) (\cA(r)) , \\
 	\varrho(\cA(q) ) \cdot  \phi(q) & = \psi(q)  ,\quad \forall q,r \in G.
 \end{align*}
 When $G$ is equipped with a topology, we shall call a monoidal isomorphism $\cA : \phi \isom \psi$ \emph{continuous},\index{continuous (monoidal isomorphism)} if the corresponding function $\cA: G \to \tilde H$ is continuous with respect to the given topology on $G$ and the discrete topology on $\tilde H$.
 
 \subsection{}\label{para:phi tilde ab}  Let $k$ be a field of characteristic zero, and let $\bar k$ be a fixed algebraic closure. Let $G$ be a reductive group over $k$, and let $\G_G$ be the associated neutral gerb in $\redgb(\bar k/k)$. As in \cite[\S  3.2.2]{kisin2012modp}, we have a crossed module $G_{\sconn} (\bar k) \to \G_G$, and we denote the corresponding strictly monoidal category $\G_G/ G_{\sconn} (\bar k)$ by $\G_{G/G_{\sconn}}$.\index[Oct2026-sorted-notation]{g 0 04 00 g subscript g g subscript sconn@$\mathfrak G_{G/G_{\sconn}}$} We have a canonical strictly monoidal functor $\G_G\to  \G_{G/G_{\sconn}}$.
 
 \begin{lem}\label{lem:conj isom}
 	Keep the setting and notation of
 	\S \ref{para:phi tilde ab}. Let $\gH  = (\gH_i)_{i\in I} $ be an object in $\proredgb(\bar k/k)$ such that the projections $\gH^{\topo} \to \gH_i$ are surjective and such that $\gH^{\Delta}_i$ are tori for all $i \in I$. Let $\phi : \gH \to \G_G$ be a morphism in $\proredgb(\bar k/ k)$. Let $a \in Z^1(k, I_{\phi})$ and let $\phi' = a \phi$ (see \S \ref{para:twist morphism}). Let $\phi_{\widetilde {\ab}}$\index[Oct2026-sorted-notation]{zz21 phi 1 00 00 phi subscript widetilde ab@$\phi_{\widetilde {\ab}}$} (resp.~$\phi'_{\widetilde {\ab}}$) be the composite strictly monoidal functor
 	$$\gH^{\topo} \xrightarrow{\phi^{\topo} \text{ (resp.~$\phi^{\prime, \topo}$)}} \G_G \To \G_{G/G_{\sconn}}. $$
 	Then there exist $g \in G(\bar k)$ and a continuous monoidal isomorphism $\phi_{\widetilde {\ab}} \isom \Int(g) \circ \phi'_{\widetilde {\ab}}$, if and only if the class of $a$ in $\coh^1(k, I_{\phi})$ lies in the image of $ \coh^1(k, \tilde I_{\phi} ) \to \coh^1 (k , I_{\phi})$. Here $\tilde I_{\phi}$ is defined in \S \ref{para:underline isom} and \S \ref{para:underline isom in pro case}.
 \end{lem}
 \begin{proof} By the assumption that $\gH^{\topo} \to \gH_i$ are surjective, the lemma reduces to the case where $\gH \in \redgb(\bar k /k )$, which we now assume. Write $M$ for $G_{\bar k, \phi^{\Delta}}$. Let $\varrho$ denote the natural map $G_{\sconn} \to G$, and let $\tilde M := \varrho^{-1} (M) \subset G_{\sconn,\bar k}$. Write $\pi$ for the structural map $\gH \to \Gal(\bar k/k)$. For each $q \in \gH$, we write $\phi(q) = g_q \rtimes \pi(q)$, with $g_q \in G(\bar k)$.
 	
 	Assume there exist $g\in G(\bar k)$ and a continuous monoidal isomorphism $\cA: \phi_{ \widetilde{\ab}} \isom \Int(g) \circ \phi_{ \widetilde{\ab}}'$.
 	As discussed in \S \ref{para:cocycle interpretation}, $\cA$ can be viewed as a continuous function $\gH \to G_{\sconn} (\bar k)$
 	satisfying
 	\begin{align}\label{eq:cA1}
 		\cA(qr) & = \cA(q) \cdot \Int(\phi(q)) (\cA(r)) , \\ \label{eq:cA2}
 		\varrho(\cA(q) ) \cdot  \phi(q) & = \Int(g) [ a_{\pi(q)}\phi(q) ] ,\quad \forall q,r \in \gH.
 	\end{align}
 	Decompose $g$ as $g= g' z$, with $g' \in G_{\der}(\bar k), z \in Z_G(\bar k)$.
 	Fix a lift $\tilde g \in G_{\sconn}(\bar k)$ of $g'$, and define a continuous map $\cB : \gH \to G_{\sconn} (\bar k)$ by $$\cB(q) : = \tilde g^{-1} \cA(q) \cdot \Int(\phi(q)) (\tilde g)  \in G_{\sconn} (\bar k), \quad \forall q\in \gH. $$  Then by (\ref{eq:cA1}), we have
 	\begin{align} \label{eq:cB1}
 		\cB(qr)  = \cB(q) \cdot \Int(\phi(q)) (\cB(r)) , \quad \forall q, r \in \gH.
 	\end{align} By (\ref{eq:cA2}), we have
 	\begin{multline*}
 		\varrho(\cB(q)) = (g')^{-1} \cdot  [\varrho(\cA(q)) \cdot \phi(q)] \cdot g' \phi(q)^{-1} \\  = z a_{\pi(q)} \phi(q) g^{-1} \cdot g' \phi(q)^{-1}  = z a_{\pi(q)} g_{q} \pi(q) z^{-1} \pi(q)^{-1} g_q^{-1},
 	\end{multline*}
 	i.e.,
 	\begin{align} \label{eq:cB2}
 		\varrho(\cB(q))  = z \lix^{{\pi(q)}} z^{-1} a_{\pi(q)}, \quad \forall q \in \gH.
 	\end{align}
 	From (\ref{eq:cB2}), we see that
 	\begin{align}\label{eq:cB3}
 		\cB(q) \in \tilde M(\bar k), ~\forall q\in \gH  ;\\
 		\cB(q) \in Z_{G_{\sconn}} (\bar k), ~\forall q\in \gH^{\Delta}(\bar k)
 		. \end{align} Using (\ref{eq:cB1}) and (\ref{eq:cB3}), we see that the map $\cB|_{\gH^{\Delta} (\bar k) }: \gH^{\Delta} (\bar k) \to Z_{G_{\sconn}} (\bar k)$ is a group homomorphism. Since $\gH^{\Delta} (\bar k)$ is a divisible abelian group and $Z_{G_{\sconn}} (\bar k)$ is finite, we have $\cB|_{\gH^{\Delta} (\bar k) } \equiv 1$. Combining the last fact with (\ref{eq:cB1}), we see that $\cB(q)$ depends only on $\pi(q)$, i.e., $\cB = B \circ \pi$ for a continuous map $B : \Gal(\bar k/k)  \to \tilde M(\bar k)$. Since $\tilde I_{\phi}$ is the unique $k$-form of $\tilde M$ which is compatible with the $k$-form $I_\phi$ of $M$ (see Lemma \ref{lem:tilde rat}), the relation (\ref{eq:cB1}) precisely means that $B \in Z^1(k, \tilde I_{\phi})$. Note that the inclusion of $\bar k$-groups $Z_{G, \bar k} \hookrightarrow M$ induces an inclusion of $k$-groups $Z_G \to I_{\phi}$. Hence (\ref{eq:cB2}) implies that the class of $a$ in $\coh^1(k, I_{\phi})$ equals the image of the class of $B$ in $\coh^1(k, \tilde I_{\phi})$.
 	
 	Conversely, assume the class of $a$ in $\coh^1(k,I_{\phi})$ lies in the image of $\coh^1(k,\tilde I_{\phi})$. For each $\tau \in \Gal(\bar k/k)$ choose a lift $q_\tau \in \mathfrak H$, and write $g_\tau$ for $g_{q_\tau}$. Then there exist $g \in I_{\phi} (\bar k)$ and $B\in Z^1(k, \tilde I_{\phi})$ such that
 	$$ a_{\tau} = g \cdot \varrho(B(\tau))\cdot \Int(g_{\tau})(\lix^{\tau} g^{-1}), \quad \forall \tau \in \Gal(\bar k/k).$$ Decompose $g$ as $g = \varrho(\tilde g) z^{-1}$, with $\tilde g\in \tilde M(\bar k)$ and $z \in Z_G(\bar k)$. Then we can absorb $\tilde g$ into $B$ by replacing each $B(\tau)$ with $\tilde g B(\tau) \Int(g_{\tau}) (\lix^{\tau} \tilde g^{-1})$. Hence we may assume that $g= z^{-1}$. Let $\cB : = B \circ \pi  : \gH \to \tilde M(\bar k)$. Then $\cB$ is continuous and satisfies (\ref{eq:cB1}) and (\ref{eq:cB2}). Note that the relation (\ref{eq:cB2}) can also be written as
 	\begin{align}\label{eq:cB4}
 		\varrho(\cB(q)) \cdot \phi(q) = \Int(z)[a_{\pi(q)} \phi(q)]
 	\end{align} (since $\lix^{{\pi(q)}} z = \Int(\phi(q)) (z)$). By (\ref{eq:cB1}) and (\ref{eq:cB4}), $\cB$ is a continuous monoidal isomorphism $\phi_{ \widetilde{\ab}} \isom \Int(z) \circ \phi'_{\widetilde{\ab}}$.
 \end{proof}
 \begin{rem}\label{rem:Qf satisfies hypo}
 	By the discussion in \S \ref{subsubsec:unique-quasi-motivic}, the assumptions on $\mathfrak H$ in Lemma \ref{lem:conj isom} are satisfied by the quasi-motivic Galois gerb $\Qf \in \proredgb(\Qbar/\Q)$. 
 \end{rem}
 
 \section{Admissible morphisms for an unramified Shimura datum}\label{sec:adm morph}
 In this section, we define \emph{admissible morphisms} $\phi$ and the associated sets $X(\phi)$ and $S(\phi)$, which constitute the core components in the Langlands--Rapoport Conjecture. We also introduce a twisted variant $S_{\tau}(\phi)$ of $S(\phi)$, allowing for a  suitable weakening of the conjecture. 
 \subsection{} \label{subsubsec:SD}
 Let $(G,X, p , \cG)$ be an unramified Shimura datum, with reflex field $E$, and let $\fkp$ and $q= p^r$ be as in \S \ref{para:unramified SD}. We will use the following notation throughout the paper. For each field extension $F/E$ such that $G_F$ is quasi-split, the Hodge cocharacters\index[Oct2026-sorted-notation]{zz12 mu 1 00 00 mu subscript h@$\mu_h$} $\mu_h$ attached to $h \in X$ determine a $G(F)$-conjugacy class of cocharacters of $G_F$. We denote this conjugacy class by $\dmu_X(F)$.\index[Oct2026-sorted-notation]{zz12 mu 1 06 00 bbmu subscript x f@$\bbmu_X(F)$} Now inside $\dmu_X(\QQ_{p^r})$, there is a canonical $\cG(\Z_{p^r})$-conjugacy class consisting of those cocharacters in $\dmu_X(\QQ_{p^r})$ that extend to cocharacters of $\cG_{\Z_{p^r}}$. We denote this $\cG(\Z_{p^r})$-conjugacy class by $\dmu_X^{\cG}$.\index[Oct2026-sorted-notation]{zz12 mu 1 06 00 bbmu subscript x superscript g@$\bbmu_X^{\mathcal G}$} 
 
 A choice of $x\in X$ gives rise to a morphism\index[Oct2026-sorted-notation]{zz14 xi 1 00 00 xi subscript v@$\xi_v$} $$\xi_\infty: \G_\infty \To \G_{G}(\infty)$$ in $\redgb(\CC/\RR)$, whose conjugacy class depends only on $X$. See \cite[\S 3.3.5]{kisin2012modp} for the explicit construction. For a finite prime $v$ unequal to $p$, let $$\xi_v : \G_v = \Gamma_v \To \G_G(v) = G(\overline {\Q}_v) \rtimes \Gamma_v $$ be the natural section. Then $\xi_v$ is a morphism in $\redgb(\overline \Q_v /\Q_v)$.
 
 Let $\phi : \Qf \to \G_G$ be a morphism in $\proredgb(\Qbar/\QQ)$. (Here $\Qf$ is defined with respect to the fixed $p$.) For each place $v$ of $\Q$ except $p$, we define\index[Oct2026-sorted-notation]{x 0 00 00 x subscript v phi@$X_v(\phi)$}
 $$ X_v(\phi) : = \set{ g\in G( \overline \Q_v) \mid  \Int(g)\circ \xi_v = \phi(v)\circ \zeta_v } . $$ We define\index[Oct2026-sorted-notation]{x 0 00 00 x subscript p phi@$X_p(\phi)$}
 $$X_p(\phi) := X_{-\mu}(\phi(p) \circ \zeta_p) $$ for $\mu \in \dmu_X^{\cG}$. Here the right hand side is as in \S \ref{subsubsec:recall of defns in LR}, and it is independent of the choice of $\mu$. We have the $p^r$-Frobenius\index[Oct2026-sorted-notation]{zz21 phi 0 00 00 phi frobenius@$\Phi$ (Frobenius)} $\Phi: X_p(\phi) \isom X_p(\phi)$.
 
 \begin{defn}[cf.~{\cite[\S 3.3.6]{kisin2012modp}}] \label{defn of adm morphism} Keep the setting of \S \ref{subsubsec:SD}.	A morphism $\phi: \Qf \to \G_G$ in $\proredgb(\Qbar/\Q)$ is called\index{admissible!morphism} \emph{admissible}, if the following conditions are satisfied. \begin{enumerate}
 		\item Let $\mu\in \dmu_X(\Qbar)$. Let $\psi_{\mu_{\widetilde{\ab}}}:  \Qf \To \G_{G/G_{\sconn}}$ be the associated strictly monoidal functor, defined in \cite[\S 3.3.1]{kisin2012modp}. Let $\phi_{ \widetilde{\ab}}$ be the composite strictly monoidal functor $\Qf \xrightarrow{\phi} \G_G \to \G_{G/G_{\sconn}}$. We require that there exist $g \in G(\overline \QQ)$ and a continuous monoidal isomorphism $
 		\cA: \Int(g) \circ \psi_{\mu_{\widetilde{\ab}}} \isom \phi_{ \widetilde{\ab}}.$
 		\item For each place $v$ of $\Q$, $X_v(\phi) \neq \emptyset$.
 	\end{enumerate}
 \end{defn}
 \begin{rem}\label{rem:continuous conj isom} Condition (i) in Definition \ref{defn of adm morphism} is a correction of condition (1) in \cite[\S 3.3.6]{kisin2012modp} in that we add the requirement that $\cA$ should be continuous. The results on admissible morphisms in \cite[\S 3.4]{kisin2012modp}, especially the statement and proof of \cite[Prop.~3.4.11]{kisin2012modp}, are only correct with the present definition.
 \end{rem}
 \begin{rem} Given a morphism $\Qf \to \G_G$ in $\proredgb(\Qbar/\Q)$, whether it is admissible depends only on its conjugacy class. We may thus speak of admissible conjugacy classes of morphisms $\Qf \to \G_G$.
 \end{rem}
 
 \begin{rem}
 	As is pointed out by Reimann \cite[App.~B]{reimann1997zeta}, the original construction of the quasi-motivic gerb $\gQ$ in \cite{langlands1987gerben} is incorrect. As such, the definitions and results in \cite[\S 5]{langlands1987gerben} about ``admissible morphisms'' do not directly apply to the objects defined in Definition \ref{defn of adm morphism}. Nevertheless, most of the content of \textit{loc.~cit.~}can be modified to suit the corrected definition of $\Qf$. In the sequel, we shall cite \textit{loc.~cit.~}only for those technical results that are essentially independent of the actual construction of $\Qf$.
 \end{rem}
 \subsection{}\label{para:integral points in X_l}
 Let $\phi: \Qf \to \G_G$ be an arbitrary morphism in $\proredgb(\Qbar/\Q)$. By condition (ii) in \cite[Def.~B2.7]{reimann1997zeta}, there exists a continuous cocycle\index[Oct2026-sorted-notation]{zz06 zeta 1 00 00 zeta subscript phi superscript p infty@$\zeta_{\phi}^{p,\infty}$} $\zzeta_\phi: \Gamma \to G(\bar \A_f^p)$ that induces the morphisms $\phi(l) \circ \zeta_l$ for all finite primes $l \neq p$, in the following sense: We have a canonical map $\bar \A_f^p = \A_f^p \otimes _{\QQ} \Qbar \to \Qlbar$ given by the projection $\A_f^p \to \Q_l$ and the fixed embedding $\Qbar \to \Qlbar$. Denote the composition $$\Gamma_l \hookrightarrow  \Gamma \xrightarrow{\zzeta_{\phi}}   G(\bar \A_f^p) \To  G(\Qlbar)  $$ by $\zeta_{\phi,l}$.\index[Oct2026-sorted-notation]{zz06 zeta 1 00 00 zeta subscript phi l@$\zeta_{\phi,l}$} Then for each $\tau \in \G_l = \Gamma_l$ we have $$(\phi(l) \circ \zeta_l) (\tau) = \zeta_{\phi,l}(\tau) \rtimes \tau  \in G(\overline {\Q}_l) \rtimes \Gamma_l.$$
 
 If we choose an arbitrary $\ZZ$-structure on $G$, then for almost all primes $l \neq p$ the set $X_l(\phi)$ contains integral points in $G(\overline \Q_l)$ (and is \textit{a fortiori} non-empty).	Indeed, for almost all $l$, the chosen $\Z$-structure on $G$ has connected smooth reduction at $l$, and $\zeta_{\phi,l}$ is induced by a continuous unramified cocycle $\Gal(\Q_l^{\ur}/\Q_l) \to G(\Z_l^{\ur})$. It is a standard result (see for instance \cite[p.~294, Thm.~6.8']{plantonov-rapinchuk}) that any such cocycle is a coboundary, and this precisely means that $X_l(\phi)$ contains integral points.
 
 \subsection{}
 \label{subsubsec:defn of S(phi)}
 Let $\phi : \Qf \to \G_G$ be an admissible morphism. On choosing a $\ZZ$-structure on $G$, we form the restricted product\index[Oct2026-sorted-notation]{x 0 00 00 x superscript p phi@$X^p(\phi)$} $$ X^p(\phi) : = \prod'_{l\notin \set{p,\infty}} X_l(\phi)$$   with respect to the subsets of integral elements of $X_l(\phi)$ discussed in \S \ref{para:integral points in X_l}. Clearly $X^p(\phi)$ is independent of the choice of $\ZZ$-structure, and is a $G(\mathbb A^p_f)$-torsor under right multiplication. (It is non-empty since $\phi$ is admissible.) Equivalently, with the notation in \S \ref{para:integral points in X_l}, $X^p(\phi)$ is the right $G(\A_f^p)$-torsor consisting of $x \in G(\bar \A_f^p)$ such that $$x^{-1} \cdot  \zeta_{\phi}^{p,\infty} (\tau) \cdot \lix^{\tau}x = 1, \quad  \forall \tau \in \Gamma. $$
 
 We now define\index[Oct2026-sorted-notation]{x 0 00 00 x phi@$X(\phi)$}
 $$ X(\phi) : = X_p (\phi) \times X^p(\phi),$$ which is equipped with the action of $\Phi^{\ZZ} \times G(\adele_f^p)$. We still call $\Phi$ the \emph{$q$-Frobenius} on $X(\phi)$. 
 
 By definition, $X(\phi)$ is a subset of $G(\Qpbar)/ \cG(\Z_p^{\ur}) \times G(\bar \A_f^p)$. Under the canonical embedding $I_{\phi, \Qbar} \hookrightarrow G_{\Qbar}$, we let the group $I_{\phi}(\A_f)$  act on $G(\Qpbar)/ \cG(\Z_p^{\ur}) \times G(\bar \A_f^p)$ via left multiplication. This induces a left action of $I_{\phi}(\A_f)$ on $X(\phi)$. For each $\tau \in I_{\phi}^{\ad}(\A_f)$, we set\index[Oct2026-sorted-notation]{s 0 00 00 s subscript tau phi@$S_{\tau}(\phi)$} $$ S_{\tau}(\phi) : = \varprojlim_{K^p} I_\phi(\QQ)_{\tau}\backslash X(\phi)/ K^p,$$ where $K^p$ runs through the compact open subgroups of $G(\mathbb A^p_f)$, and $I_{\phi}(\QQ)_{\tau}$\index[Oct2026-sorted-notation]{i 0 00 00 i subscript phi q subscript tau@$I_{\phi}(\QQ)_{\tau}$} denotes the image of the composite map 
 $$ I_{\phi}(\QQ) \hookrightarrow I_{\phi}(\A_f) \xrightarrow{\Int(\tau)} I_{\phi}(\A_f).$$ Then $S_{\tau}(\phi)$ inherits the action of $\Phi^{\ZZ} \times G(\adele_f^p)$. When $\tau =1$, we write\index[Oct2026-sorted-notation]{s 0 00 00 s phi@$S(\phi)$} $S(\phi)$ for $S_{\tau}(\phi)$.

 \section{Integral models and the Langlands--Rapoport Conjecture}
 \label{sec:int mod}
 In this section we review the theory of smooth integral models of Shimura varieties and their cohomology, and then formulate the Langlands--Rapoport Conjecture. We incorporate the requirement that the integral model has \emph{well-behaved compactly supported cohomology} into our formulation, so that the description of points on the special fiber of the model provided by the conjecture reflects information about the cohomology of the generic fiber. 
 \subsection{}\label{para:model} Keep the setting and notation of \S \ref{subsubsec:SD}. Let $K_p = \mathcal G(\ZZ_p)$,\index[Oct2026-sorted-notation]{k 0 00 00 k subscript p@$K_p$} and let $\Sh_{K_p} =\Sh_{K_p}(G,X)$\index[Oct2026-sorted-notation]{sh 0 01 00 sh subscript k subscript p@$\mathrm{Sh}_{K_p}$} be the inverse limit $$\varprojlim_{K^p}\Sh_{K_pK^p}(G,X), $$ where $K^p$ runs through compact open subgroups of $G(\A_f^p)$. This inverse limit exists as an $E$-scheme, since the transition maps are affine. The right $G(\A_f^p)$-action on $\Sh_{K_p}$ induced by the $G(\A_f)$-action on $\Sh$ is admissible in the sense of Definition \ref{defn:adm action}.

 \begin{defn}\label{defn:integral model} By a\index{smooth integral model} \emph{smooth integral model} of $\Sh_{K_p}$, we mean a scheme $\Shh_{K_p}$\index[Oct2026-sorted-notation]{s 0 05 00 s subscript k subscript p@$\mathscr S_{K_p}$} over $\oo_{E,(\fkp)}$ extending $\Sh_{K_p}$, equipped with an admissible right $G(\A_f^p)$-action extending the $G(\A_f^p)$-action on $\Sh_{K_p}$ such that $\Shh_{K_p K^p} : = \Shh_{K_p}/K^p$\index[Oct2026-sorted-notation]{s 0 05 00 s subscript k subscript pk superscript p@$\mathscr S_{K_pK^p}$} is smooth, separated, and of finite type over $\oo_{E,(\fkp)}$ for all sufficiently small compact open subgroups $K^p \subset G(\adele_f^p)$.
 \end{defn}
 The following theorem is proved in \cite{kisin2010integral} for $p>2$, and in \cite{KMP16} for $p=2$. 
 \begin{thm}[{\cite{kisin2010integral,KMP16}}] \label{thm:integral model}
 	If $(G,X)$ is of abelian type, then there exists a  smooth integral model of $\Sh_{K_p}$. This model is uniquely characterized by the extension property as detailed in \cite[\S 2.3.7]{kisin2010integral}.
 \end{thm} 
 \begin{rem}In \cite{kisin2010integral} and  \cite{KMP16}, the authors do not explicitly verify that $\Shh_{K_pK^p}$ is separated over $\oo_{E,(\fkp)}$, nor that the $G(\A_f^p)$-action on $\Shh_{K_p}$ satisfies condition (iii) in Definition \ref{defn:adm action}. The former follows from the facts that the Siegel modular schemes at finite levels are separated over $\ZZ_{(p)}$, that normalization maps and closed immersions are separated, and that taking finite free quotients preserves separatedness. For the latter, see \cite[\S 3]{lanstrohII} for an explanation.
 \end{rem}
 \subsection{}\label{para:integral Galois}
 Suppose a smooth integral model $\Shh_{K_p}$ is given. Let $K^p \subset G(\A_f^p)$ be a sufficiently small compact open subgroup, and $U^p \subset K^p$ an open normal subgroup.  Since $\Shh_{K_pU^p}$ is flat, separated, and of finite type over $\oo_{E, (\fkp)}$, any $\oo_{E, (\fkp)}$-automorphism of it is uniquely determined by the induced $E$-automorphism of $\Sh_{K_p U^p}$. Moreover, if an $\oo_{E, (\fkp)}$-automorphism acts as the identity on a connected component of $\Shh_{K_pU^p}$, then it fixes a $\CC$-point. Hence it immediately follows from  (\ref{eq:(K/U)'}) that condition ($\dagger$) in \S \ref{para:sheaf and Hecke} is satisfied by the  $G(\adele_f^p)$-action on $\Shh_{K_p}$, and that we have  $$(K^p/U^p)^{\natural, \Shh_{K_p}} = (K^p/U^p)^{\natural,\Sh_{K_p}} = (K_pK^p / K_pU^p)^{\natural, \Sh} =  K^p / U^p Z(\Q)_{K_p K^p}^{-, (p)},  $$ where $Z(\Q)_{K_pK^p}^{-, (p)}$ is the image of $Z(\Q)_{K_pK^p}^-$ (notation as in \S \ref{para:aut sheaf}) under the projection $G(\A_f) \to G(\A_f^p)$. From this we also see that condition ($\ddagger$) in \S \ref{para:sheaf and Hecke} is satisfied by the $G(\adele_f^p)$-action on $\Shh_{K_p}$. Hence $ (K^p)^{\natural, \Shh_{K_p}}$ as defined in (\ref{eq:profinite Galois}) is the quotient of $K^p$ by the closure of $Z(\Q)_{K_pK^p}^{-, (p)}$ in $G(\A_f^p)$. (In fact $Z(\Q)_{K_pK^p}^{-, (p)}$ is already closed in $G(\A_f^p)$, because $Z(\Q)_{K_p K^p}^-$ is compact.) By Lemma \ref{lem:Z_s}, (the closure of) $Z(\Q)_{K_pK^p}^{-, (p)}$ is contained in $(Z_G^0)_{ac}(\A_f^p)$ when $K^p$ is sufficiently small.

 Now fix a prime $\ell \neq p$ and fix an irreducible $\overline \QQ_{\ell}$-representation $\xi$ of $G$ that factors through $G^c = G/(Z_G^0)_{ac}$, as in \S\ref{para:aut sheaf}. Via the projection $G(\A_f^p) \to G(\QQ_{\ell})$, $\xi$ induces a continuous representation of $G(\A_f^p)$.
 Then for all sufficiently small $K^p$, the restriction $\xi|_{K^p}$ factors through $(K^p)^{\natural, \Shh_{K_p}}$ by the above discussion. Thus as in \S \ref{para:sheaf and Hecke}, for each sufficiently small $K^p$ we obtain a lisse $\ol \QQ_{\ell}$-sheaf $\mL_{\xi, K^p}$ on $\Shh_{K_pK^p}$, and for each geometric point $x$ of $\Spec \oo_{E, (\fkp)}$ we have the admissible $G(\A_f^p)$-module\index[Oct2026-sorted-notation]{h 0 02 00 h superscript i subscript c s subscript k subscript p x xi@$\mathbf H^i_c(\mathscr S_{K_p, x}, \xi)$}
 $$ \coh^i_c(\Shh_{K_p, x}, \xi) : = \varinjlim_{K^p} \coh^i_c(\Shh_{K_pK^p, x}, \mL_{\xi, K^p} ). $$
 When $x = \Spec \ol E$, the above is identified with $\coh^i_c (\Sh_{\ol E}, \xi)^{K_p}$, where $\coh^i_c (\Sh_{\ol E}, \xi)$ is as in \S \ref{para:aut sheaf}. Moreover, we have a canonical adjunction morphism
 \begin{align}\label{eq:adjunction for coh}
 	\coh_c^i(\Shh_{K_p,\ol{\FF}_q},\xi) \To \coh_c^i(\Shh_{K_p,\ol E},\xi)  \cong \coh^i_c(\Sh_{\ol E}, \xi)^{K_p},
 \end{align} which is $\Gal(\ol {E_{\fkp}}/E_{\fkp}) \times G(\A_f^p)$-equivariant. Here  $\Gal(\ol {E_{\fkp}}/E_{\fkp})$ acts on the left via the quotient $\Gal(\ol \FF_q/ \FF_q)$, and acts on the right via the embedding into $\Gal(\ol E/E)$. 
 
 \begin{defn}\label{defn:coh can} We say that $\Shh_{K_p}$ has\index{well-behaved compactly supported cohomology} \emph{well-behaved $\coh^*_c$}, if (\ref{eq:adjunction for coh}) is an isomorphism for all choices of $\ell \neq p$, $\xi$, and $i$.
 \end{defn}
 \begin{thm}[{\cite[Cor.~4.6]{lanstrohII}}]\label{thm:lan-stroh}
 	The canonical smooth integral model in Theorem \ref{thm:integral model} has well-behaved $\coh^*_c$.
 \end{thm}
 
 Recall that for each admissible morphism $\phi : \Qf \to \G_G$, we have defined in \S \ref{subsubsec:defn of S(phi)} a set $S(\phi)$ equipped with an action of $G(\A_f^p)$ and a $q$-Frobenius $\Phi$, where $q$ is the residue cardinality of $\fkp$. 
 \begin{conj}[Langlands--Rapoport]\label{conj:LR} There\index{Langlands--Rapoport Conjecture} exists a smooth integral model $\Shh_{K_p}$ of $\Sh_{K_p}$ over $\oo_{E,(\fkp)}$ which has well-behaved $\coh^*_c$ and for which there is a bijection
 	$$
 	\Shh_{K_p} (\overline \FF_q)  {} \isom \coprod_{\phi }   S(\phi)  $$
 	compatible with the actions of $G(\A_f^p)$ and the $q$-Frobenius $\Phi$. Here $\phi$ runs through a set of representatives for the conjugacy classes of admissible morphisms $\Qf \to \G_G$.
 \end{conj}
 In the rest of this chapter, we formulate a variant of the above conjecture, which we call the Langlands--Rapoport--$\tau$ Conjecture.  
 
 \section{The Langlands--Rapoport--\texorpdfstring{$\tau$}{tau} Conjecture} \label{subsec:preparation LR}
 In this section, we formulate the Langlands–Rapoport–$\tau$ Conjecture by establishing the necessary foundational framework. We first introduce the notion of a \emph{g\"unstig gelegen} morphism, which provides the technical basis for studying the cohomological twisting of admissible morphisms (cf.~Proposition \ref{prop:twist adm morph}); this notion remains central to our later analysis in \S \ref{sec:adm morphisms and maximal tori}. The core strategy in weakening the original Langlands–Rapoport Conjecture lies in the introduction and \emph{precise control} of a family of twisting elements $\dertau \in \Gamma(\mathcal H)$ (see Definition \ref{defn:system of obstructions}). These elements allow us to modify the standard sets $S(\phi)$ into their twisted counterparts. By imposing two compatibility conditions on this family with respect to cohomological twisting, we define a hierarchy of subsets $\Gamma(\mathcal H)_0 \subset \Gamma(\mathcal H)_1 \subset \Gamma(\mathcal H)$. Finally, we introduce the property of \emph{tori-rationality}, which—together with the condition $\Gamma(\mathcal H)_0$—forms the basis for the statement of the conjecture. 
 \subsection{}\label{para:choice of q_tau}
 Since the projections $\gQ^{\topo} \to \gQ^L$ are surjective (see \S \ref{subsubsec:unique-quasi-motivic}),  we may choose a lift $q_{\rho} \in \gQ^{\topo}$ for each $\rho \in \Gamma$. We fix these choices in the sequel. We first study a ``well-positioned'' condition for morphisms from $\gQ$ to neutral gerbs. Let $G$ be an arbitrary reductive group over $\Q$. Recall from \S \ref{subsubsec:unique-quasi-motivic} that $\gQ^{\Delta}$ has the canonical $\Q$-structure $Q$, which is a pro-torus. Hence by the discussions in \S \ref{para:underline isom} and \S \ref{para:underline isom in pro case}, for each morphism 	$\phi : \gQ \to \G_G$ in $\proredgb(\Qbar/\Q)$  the group $I_{\phi}$ is reductive  and a $\QQ$-form of the  Levi subgroup $G_{\Qbar, \phi^{\Delta}}$ of $G_{\Qbar}$. Moreover we have $I_\phi^{\dagger}$ and $\tilde I_{\phi}$ defined. 
 
 \begin{defn}\label{defn:gg morph}	A morphism $\phi : \gQ \to \G_G$ in $\proredgb(\Qbar/\Q)$ is  called\index{gunstig gelegen@g\"unstig gelegen (gg)!morphism} \emph{g\"unstig gelegen}  (abbreviated as \emph{gg}) if the morphism $\phi^{\Delta}: \gQ^{\Delta}  \to G_{\overline \Q}$ is defined over $\Q$.\footnote{The term originates from Langlands and Rapoport \cite[\S 5]{langlands1987gerben}. Note that their original definition utilizes the elements $\delta_n$, which are not directly compatible with the corrected definition of $\gQ$ used here.}
 \end{defn}
 
 \begin{lem}\label{lem:gg morph}Let $\phi : \gQ \to \G_G$ be a morphism in $\proredgb(\Qbar/\Q)$. For each $\rho \in \Gamma = \Gal(\Qbar /\QQ )$, write $\phi(q_{\rho}) = g_{\rho} \rtimes \rho$, with $g_{\rho} \in G(\overline \Q)$. Then $\phi$ is gg if and only if $g_{\rho}$ lies in $G_{\overline \Q, \phi^{\Delta}} (\overline \Q)$ for all $\rho \in \Gamma$. If $\phi$ is gg, then the canonical $\overline \Q$-isomorphism $I_{\phi,\overline \Q} \isom G_{\phi^{\Delta}, \Qbar}$ is an inner twisting between the underlying $\Q$-groups.
 \end{lem}
 \begin{proof} Assume that $\phi$ is induced by a morphism $\phi_0: \gQ^L \to \G_G$ in $\redgb(\ol \Q/\Q)$. Let $y \in Q^L(\Qbar) =  \gQ^{L,\Delta}(\overline \Q)$ and $\rho \in \Gamma$ be arbitrary. We denote the image of $q_{\rho} \in \gQ^{\topo}$ in $\gQ^L$ still by $q_{\rho}$. Then
 	\begin{align}\label{eq:phi_0 tau t}
 		\phi_0^{\Delta}(\lix^{\rho} y) = \phi_0(q_\rho y  q_{\rho}^{-1}) =  (g_{\rho} \rtimes \rho)  \phi_0^{\Delta}(y) (g_{\rho} \rtimes \rho) ^{-1} = g_{\rho} \lix^{\rho}[\phi_0^{\Delta}(y)] g_{\rho}^{-1}.\end{align} Now $\phi^{\Delta}$ is defined over $\Q$ if and only if $\phi_0^{\Delta}$ is defined over $\Q$ (since the transition maps in the pro-torus $Q = \varprojlim_L Q^L$ are all surjective). By (\ref{eq:phi_0 tau t}), this is equivalent to the condition that each $g_{\rho}$ centralizes $\im (\phi_0^{\Delta})$. Since $\im(\phi_0^{\Delta}) = \im (\phi^{\Delta})$, the last condition is equivalent to the condition that each $g_{\rho}$ lies in $G _{\overline \Q, \phi^{\Delta}} (\overline \Q)$.
 	
 	Now if $\phi$ is gg, then for each $\rho \in \Gamma$ the $\overline \Q$-isomorphism $I_{\phi,\overline \Q} \isom G_{\phi^{\Delta}, \Qbar}$ differs from its $\rho$-twist by composition with $\Int(g_{\rho})$. Since $g_{\rho} \in G_{\phi^{\Delta}}(\Qbar)$, this means that the $\Qbar$-isomorphism is an inner twisting. 
 \end{proof}
 \begin{defn}\label{defn:G-rational tori}
 	Let $\phi : \gQ \to \G_G$ be a morphism in $\proredgb(\Qbar/\Q)$. By a\index{G-rational maximal torus@$G$-rational maximal torus} \emph{$G$-rational maximal torus in $I_\phi$}, we mean a maximal torus $T \subset I_{\phi}$ (defined over $\Q$) such that the composite embedding $T_{\overline \Q} \hookrightarrow I_{\phi,\overline \Q} \cong  G_{\overline \Q, \phi^{\Delta}}\ \hookrightarrow G_{\overline \Q}$ is  defined over $\Q$.
 \end{defn}
 \begin{rem}\label{rem:G-rational torus}
 	In Definition \ref{defn:G-rational tori}, $T$ is necessarily a maximal torus in $G$, since  $I_{\phi}$ has  the same absolute rank as $G$, cf.~\S \ref{para:choice of q_tau}.
 \end{rem}
 \begin{lem}\label{lem:suff gg}
 	Let $\phi : \gQ \to \G_G$ be a morphism in $\proredgb(\Qbar/\Q)$ such that $I_{\phi}$ contains a $G$-rational maximal torus $T$. Then $\phi$ is gg. Moreover, let $f$ denote the $\Q$-embedding underlying $T_{\Qbar} \hookrightarrow I_{\phi,\Qbar} \cong G_{\Qbar, \phi^{\Delta}} \hookrightarrow G_{\Qbar}$. Then $\phi$ factors as $\Qf \xrightarrow{\phi_T} \G_T \xrightarrow{f} \G_G$. 
 \end{lem}
 \begin{proof}For each $\rho \in \Gamma$, define $g_{\rho} \in G(\Qbar)$ as in Lemma \ref{lem:gg morph}. For   $t\in T(\overline \Q)$, we have
 	$$ f(\lix^{\rho} t ) = g_{\rho} \lix^{\rho} [f(t)] g_{\rho}^{-1},$$
 	by the definition of the $\Q$-structure of $I_\phi$; see \S \ref{para:underline isom}. Since $f$ is defined over $\Q$, we have $f(\lix^{\rho} t) = \lix^{\rho}[f(t)]$. Hence $g_{\rho}$ commutes with $f(T(\overline \Q))$. Since $f(T)$ is a maximal torus in $G$, we have $g_{\rho} \in f(T)(\Qbar)$. We conclude that $\phi$ is gg by Lemma \ref{lem:gg morph}. Moreover, since $Z_{I_{\phi}} \subset T$ and since $\phi^{\Delta}$ factors through the center of $G_{\phi^{\Delta}}$, we know that $\phi^{\Delta}$ factors through $f(T) \subset G$. We have already seen that each $g_{\rho}$ lies in $f(T)(\Qbar)$. It follows that $\phi$ factors through $f: \G_T \to \G_G$.
 \end{proof}
 
 \begin{lem}\label{lem:I_phi at infty} Let $\phi : \gQ \to \G_G$ be a morphism in $\proredgb(\Qbar/\Q)$ such that $X_{\infty}(\phi) \neq \emptyset$. Then the $\RR$-group $I_{\phi(\infty) \circ \zeta_{\infty}}$ is an inner form of $G_{\RR}$. Moreover, the $\RR$-groups $(I_{\phi}/Z_G)_{\RR} = I_{\phi(\infty)}/ Z_{G,\RR}$ and $I_{\phi(\infty) \circ \zeta_{\infty}} / Z_{G,\RR}$ are both anisotropic.
 \end{lem}
 \begin{proof} As in (\ref{eq:enlarge I}) we have an $\RR$-embedding $I_{\phi(\infty)} \hookrightarrow I_{\phi(\infty)\circ \zeta_{\infty}}$. By our assumption, $\phi(\infty)\circ \zeta_{\infty}$ is conjugate to $\xi_{\infty}$. Hence there is an $\RR$-isomorphism $I_{\phi(\infty) \circ \zeta_{\infty}} \isom I_{\xi_{\infty}}$ induced by $\Int(g)$ for some $g \in G(\CC)$. As discussed in \cite[\S 3.3.5]{kisin2012modp}, the group  $I_{\xi_{\infty}}$ is the inner form of $G_{\RR}$ with anisotropic adjoint group. The lemma follows. \end{proof}
 
 We now return to the setting of \S \ref{subsubsec:SD}. Thus we have an unramified Shimura datum $(G,X, p ,\cG)$ and the notion of admissible morphisms $\Qf \to \G_G$. 
 
 \begin{prop}\label{prop:gg morph}
 	Let $\phi : \gQ \to \G_G$ be an admissible morphism. For each maximal torus $T \subset I_\phi$ defined over $\Q$, there exists $g \in G(\overline \Q)$ such that $\Int(g) (T) $ is a $G$-rational maximal torus in $I_{\Int(g) \circ \phi}$.
 \end{prop}
 \begin{proof}This is essentially proved by Langlands--Rapoport, when they prove \cite[Lem.~5.23]{langlands1987gerben}. We sketch the argument, as the precise statement of the proposition is not explicit in \textit{op.~cit.}\footnote{The only information about $\Qf$ used in this argument is the fact that $\Qf^{\Delta}$ is a pro-torus with surjective transition maps. Hence the validity of this argument is unaffected by Reimann's correction of the definition of $\Qf$.}
 	
 	Let	$\psi: G_{\overline \Q} \isom G^*_{\overline \Q}$ be an inner twisting from $G$ to a fixed quasi-split reductive group $G^*$ over $\Q$. Then the $G^*(\Qbar)$-conjugacy class of the composite embedding $\iota: T_{\Qbar} \hookrightarrow G_{\Qbar} \xrightarrow{\psi} G^*_{\Qbar}$ is stable under $\Gamma$. By \cite[Cor.~2.2]{Kot82}, we can modify $\psi$ by an inner automorphism to arrange that $\iota$ is defined over $\QQ$.
 	
 	Let $T^* = \iota(T)$. Then $T^*$ is a maximal torus in $G^*$ defined over $\Q$, and we have a $\Q$-isomorphism $\iota: T\isom T^*$. We now check that $T^*$ transfers to $G$ locally at all places $v$, or equivalently, that some $G(\ol \QQ_v)$-conjugate of $T$ is a $\Q_v$-torus in $G_{\QQ_v}$. For $v =p$, this follows from the assumption that $G$ is unramified, and hence quasi-split, over $\QQ_p$. For $v = \infty$, this follows from the fact that $T^*/Z_{G^*} \cong T/Z_G$ is anisotropic over $\RR$ (Lemma \ref{lem:I_phi at infty}). For $v \notin \set{\infty, p}$, pick $u_v \in G(\ol \Q_v)$ such that $\Int(u_v^{-1}) \circ \phi(v) \circ \zeta_v = \xi_v$, which exists since $\phi$ is admissible. Then the canonical embedding $I_{\Int(u_v^{-1}) \circ \phi(v), \overline \Q_v} \hookrightarrow G_{\overline \Q_v}$ is defined over $\Q_v$, and hence $ \Int(u_v^{-1})(T_{\Q_v})$ is a $\Q_v$-maximal torus in $G_{\Q_v}$, as desired.
 	
 	Since $T^*$ transfers to $G$ locally and is elliptic over $\RR$, it transfers globally to $G$ by \cite[Lem.~5.6]{langlands1987gerben}.   This means there exists $g\in G(\overline \Q)$ such that $\Int (g) (T)$ is a $\Q$-maximal torus in $G$ and such that the isomorphism $\Int(g) \circ \psi^{-1} : T^* \to \Int (g) (T)$ is defined over $\Q$. It follows that $\Int(g) (T) $ is a $G$-rational maximal torus in $I_{\Int(g) \circ \phi}$. 
 \end{proof}
 \begin{cor}\label{cor:gg morph}
 	Every admissible morphism $\phi : \gQ \to \G_G$ is conjugate to a gg morphism.
 \end{cor}
 \begin{proof}By Proposition \ref{prop:gg morph}, $\phi$ is conjugate to a morphism $\phi' : \Qf \to \G_G$ such that $I_{\phi'}$ contains a $G$-rational maximal torus.
 	By Lemma \ref{lem:suff gg}, $\phi'$ is gg.
 \end{proof}
 \subsection{}\label{para:transfer I_phi}
 Let $\phi : \gQ \to \G_G$ be an admissible morphism. By  Corollary \ref{cor:gg morph} and the last statement in Lemma \ref{lem:gg morph}, the canonical embedding $I_{\phi,\overline \Q} \to G_{\overline \Q}$ is a pre-inner twisting from $I_{\phi}$ to $G$ (Definition \ref{defn:inner transfer}). We thus obtain a canonical map
 \begin{align}\label{eq:transfer I_phi}
 	\coh^1_{\ab}(\Q, I_{\phi}) \To \coh^1_{\ab}(\Q, G) ,\end{align} and we define\index[Oct2026-sorted-notation]{zzz sha 0 00 00 sha subscript g superscript infty q i subscript phi@$\Sha_G^{\infty}(\mathbb Q, I_{\phi})$} $\Sha_G^{\infty}(\Q, I_{\phi})\subset \coh^1(\Q, I_{\phi})$ as in \S \ref{subsubsec:notation of Sha}.

 \begin{prop}\label{prop:twist adm morph}
 	Let $\phi : \gQ\to \G_G$ be an admissible morphism, and $\beta \in \coh^1(\Q, I_{\phi}).$ Then $\phi^{\beta}$ (see Definition \ref{defn:twist morphism}) is admissible if and only if $\beta$ belongs to $\Sha_G^{\infty}(\Q, I_{\phi})$.
 \end{prop}
 \begin{proof}
 	The ``if'' part is proved in \cite[Lem.~4.5.6]{kisin2012modp}, under the assumption that $Z_G^0$ is cuspidal. Below we give complete proofs of both directions of the implication, taking into account the correction of the definition of admissible morphisms mentioned in Remark \ref{rem:continuous conj isom}. Since the admissibility condition is invariant under conjugacy, we use Corollary \ref{cor:gg morph} to reduce to the case where $\phi$ is gg. We now assume that $\phi$ is gg and fix a cocycle $a \in Z^1(\Q, I_{\phi})$ representing $\beta$.
 	
 	\smallskip
 	
 	\textbf{Step 1.} We show that $\beta$ has zero image under the composite
 	\begin{align} \label{eq:coh pol}
 		\coh^1(\Q, I_\phi) \xrightarrow{\ab^1_\Q} \coh^1_{\ab} (\Q, I_\phi) \xrightarrow{(\ref{eq:transfer I_phi})} \coh^1_{\ab} (\Q, G)
 	\end{align}  if and only if there exist $h \in G(\overline \Q)$ and a continuous monoidal isomorphism $\Int(h) \circ \phi_{\tilde \ab} \isom (a\phi)_{\tilde \ab}$. (See Definition \ref{defn of adm morphism} for the notation.) By Lemma \ref{lem:conj isom} and Remark \ref{rem:Qf satisfies hypo}, the latter condition is equivalent to asking that $\beta$ comes from $\coh^1(\Q, \tilde I_{\phi}).$ We have a natural exact sequence of pointed sets
 	$$ \coh^1(\Q, \tilde I_{\phi}) \To \coh^1(\Q, I_{\phi}) \To \coh^1(\Q, \tilde I_{\phi} \to I_{\phi}),$$ where $\coh^1(\Q, \tilde I_{\phi} \to I_{\phi})$ is the Galois cohomology of the crossed module $(\tilde I_{\phi} \to I_{\phi})$ of $\Q$-groups; see \cite[\S 3]{borovoi}. (The crossed module structure is the one inherited from the crossed module $G_{\sconn } \to G$.) The natural map $(Z_{G_{\sconn}} \to Z_G) \to (\tilde I_{\phi} \to I_{\phi})$ is a quasi-isomorphism of crossed modules, and therefore $\coh^1(\Q, \tilde I_{\phi} \to I_{\phi})$ is naturally isomorphic to $\coh^1_{\ab}(\Q, G)$. The composition $\coh^1(\Q, I_{\phi}) \to \coh^1(\Q, \tilde I_{\phi} \to I_\phi) \cong \coh^1_{\ab}(\Q, G)$ is equal to (\ref{eq:coh pol}). This proves the desired statement. 
 	
 	\smallskip
 	
 	\textbf{Step 2.} We show that $\beta$ has trivial image in $\coh^1 (\RR, I_{\phi})$ if and only if $\phi(\infty)\circ \zeta_{\infty}$ is conjugate to $(a\phi)(\infty) \circ \zeta_{\infty}$. By Lemma \ref{lem:twist by cocycle}, the latter condition is equivalent to the vanishing of the image of $\beta$ under $$\coh^1(\QQ, I_{\phi}) \To \coh^1(\RR, I_{\phi})  = \coh^1(\RR, I_{\phi(\infty)}) \overset{\dagger}{\To} \coh^1(\RR, I_{\phi(\infty)\circ \zeta_{\infty}}),$$ where $\dagger$ is induced by the $\RR$-inclusion $I_{\phi(\infty)} \hookrightarrow I_{\phi(\infty) \circ \zeta_{\infty}}$. Thus the desired statement boils down to $\dagger$ having trivial kernel, which follows from Lemma \ref{lem:I_phi at infty} and \cite[Lem.~4.4.5]{kisin2012modp}.
 	
 	\smallskip
 	
 	\textbf{Step 3.} We show that if $\beta$ has zero image under (\ref{eq:coh pol}), then $\phi(l) \circ \zeta_l$ is conjugate to $(a\phi)(l) \circ \zeta_l$ for all finite primes $l \neq p$. For this, it suffices to show that $\beta$ has trivial image under the composite map
 	\begin{align}\label{eq:enlarge v}
 		\coh^1(\QQ, I_{\phi}) \To \coh^1(\Q_l, I_{\phi}) = \coh^1(\Q_l, I_{\phi (l)}) \To \coh^1(\Q_l, I_{\phi(l) \circ \zeta_l}).
 	\end{align} Since $\phi$ is admissible, ${\phi(l) \circ \zeta_l}$ is conjugate to $\xi_l$. It follows that the canonical $\overline \Q_l$-embedding $I_{\phi(l) \circ\zeta_l , \overline \Q_l} \to G_{\overline \Q_l}$ is an inner twisting between $\Q_l$-groups. This induces a canonical isomorphism
 	\begin{align}\label{eq:v}
 		\coh^1_{\ab}(\Q_l, I_{\phi(l) \circ \zeta_l}) \isom \coh^1_{\ab}(\Q_l, G).
 	\end{align}
 	Now we have a commutative diagram
 	$$ \xymatrixcolsep{5pc} \xymatrix{\coh^1(\Q, I_{\phi}) \ar[r]^-{(\ref{eq:enlarge v})} \ar[d]^{(\ref{eq:coh pol})} &  \coh^1(\Q_l, I_{\phi(l) \circ \zeta_l}) \ar[r]^{\ab^1_{\Q_l}}_{\sim} &  \coh^1_{\ab}(\Q_l, I_{\phi(l) \circ \zeta_l}) \ar[d] ^{(\ref{eq:v})} _{\sim} \\ \coh^1_{\ab}(\Q, G) \ar[rr]^{\text{localization}} &  & \coh^1_{\ab}(\Q_l, G) } $$ which implies the desired statement.
 	
 	\smallskip
 	
 	\textbf{Step 4.}  We show that if $\beta$ has zero image under (\ref{eq:coh pol}), then $X_p(a\phi) \neq \emptyset$. We fix $\mu \in \dmu_X^{\cG}$ as in \S \ref{subsubsec:SD}. Let $\theta : = \phi(p) \circ \zeta_p$, and let $\theta ' : = (a\phi) (p) \circ \zeta_p$. By definition $X_p(\phi) = X_{-\mu}(\theta)$ and $X_p(a\phi) = X_{-\mu}(\theta')$.
 	
 	Fix an arbitrary $g_0 \in \mathcal{UR}(\theta)$, and write $\theta_0$ for $\Int(g_0^{-1}) \theta$. Thus $\theta_0$ is unramified. Now $\Int(g_0^{-1})$ induces a $\Q_p$-isomorphism $I_{\theta} \isom I_{\theta_0}$. Let $\beta_0$ denote the image of $\beta$ under the composite
 	$$ \coh^1(\Q, I_{\phi}) \To \coh^1(\Q_p, I_{\phi})   \To \coh^1(\Q_p, I_{\theta}) \xrightarrow{\Int(g_0^{-1})} \coh^1(\Q_p, I_{\theta_0}). $$ Then $\theta'$ belongs to the conjugacy class of $\theta_0^{\beta_0}$.
 	
 	By Proposition \ref{prop:twist at p} (iii), the conjugacy class $\theta_0^{\beta_0}$ contains an unramified member $\theta_0'$ such that $b' : = b_{\theta_0'}$ is obtained from $b: = b_{\theta_0}$ by twisting by $\beta_0$. Since $G$ is quasi-split over $\Q_p$, we can apply Proposition \ref{prop:twist p} to conclude that $\nu_{b} = \nu_{b'}$, and that $\kappa_G(b') - \kappa_G(b)$ is the image of $\beta_0$ in $\pi_1(G)_{\Gamma_p, \tors}$. Our assumption on $\beta$ implies that the last image is zero.
 	Hence we have $[b] = [b'] $ in $\B(G_{\Q_p})$, by Kottwitz's classification  (see \S \ref{para:B(G)}).
 	
 	Since $\theta_0$ (resp.~$\theta'_0$) is an unramified member in the conjugacy class of $\theta$ (resp.~$\theta'$), by the discussion in \S \ref{subsubsec:recall of defns in LR}, we have $X_{-\mu}(\theta)\cong X_{-\mu} (b)$, and $X_{-\mu}(\theta')\cong X_{-\mu}(b')$. Since $[b] = [b']$, we have $X_{-\mu} (b) \cong X_{-\mu} (b')$. Thus the non-emptiness of $X_{-\mu}(\theta)$ implies the non-emptiness of $X_{-\mu}(\theta')$.
 	
 	The proof of the proposition is completed by combining the above four steps. (The ``only if'' part follows from Steps 1 and 2 alone.)
 \end{proof}
 \subsection{}\label{para:defn H}
 Let $\phi: \Qf \to \G_G$ be an admissible morphism. We define\index[Oct2026-sorted-notation]{h 0 03 00 h phi@$\mathcal H(\phi)$} \index[Oct2026-sorted-notation]{e 0 04 00 e superscript p phi@$\mathfrak E^p(\phi)$}
 \begin{align*}
 	\mathcal H(\phi) & : = \tauhome{\phi}, \\
 	\E^p (\phi) & : =  I_{\phi} (\A_f^p) \backslash I_{\phi} ^{\ad} (\A_f^p) .
 \end{align*} 
 We have a natural map 
 $$\E^p (\phi) \To \mathcal H(\phi), $$ and it is surjective by weak approximation (see \cite[Thm.~7.8]{plantonov-rapinchuk}) applied to $I_{\phi} ^{\ad}$. 
 
 The coboundary maps arising from the short exact sequence $1 \to Z_{I_{\phi} } \to I_{\phi} \to I_{\phi}^{\ad} \to 1$ induce isomorphisms of pointed sets $\E^p(\phi) \isom \D (Z_{I_{\phi}}, I_{\phi} ;\A_f^p)$ and $I_{\phi}(\A_f) \backslash I_{\phi}^{\ad}(\A_f) \isom \D (Z_{I_{\phi}}, I_{\phi} ;\A_f)$. Note that $\D (Z_{I_{\phi}}, I_{\phi} ;\A_f^p)$ and  $\D (Z_{I_{\phi}}, I_{\phi} ;\A_f)$ are both abelian groups. In particular we have a canonical structure of abelian group on $\E^p(\phi)$.  
 \begin{lem} \label{lem:cohomological H(phi)}
 	The bijection $I_{\phi}(\A_f) \backslash I_{\phi}^{\ad}(\A_f) \isom \D (Z_{I_{\phi}}, I_{\phi} ;\A_f)$ induces an isomorphism of pointed sets  $$ \cH(\phi) \isom \D (Z_{I_{\phi}}, I_{\phi} ;\A_f) / \mathfrak D(Z_{I_{\phi,\sconn}} , I_{\phi,\sconn} ; \QQ), $$
 	defining a canonical structure of abelian group on $\cH(\phi)$. The map $\E^p(\phi) \to \mathcal H(\phi)$ is a homomorphism, and  we have a commutative diagram 
 	\begin{align}
 		\xymatrix{  \E^p (\phi) \ar[rr]^{\sim} \ar[d] && \D(Z_{I_\phi}, I_\phi ;\A_f^p) \ar[d] \\    \mathcal H(\phi) \ar[rr]^-{\sim} &&  \D (Z_{I_\phi}, I_\phi; \A_f) /\D (Z_{I_{\phi,\sconn}} , I_{\phi, \sconn} ;\Q) }
 	\end{align}
 	where the vertical arrow on the right is induced by the inclusion $\D(Z_{I_\phi}, I_\phi ;\A_f^p)  \to \D(Z_{I_\phi}, I_\phi ;\A_f)$.  
 \end{lem}
 \begin{proof}First note that $\D (Z_{I_{\phi,\sconn}} , I_{\phi, \sconn} ;\Q)$ is a subgroup of $\coh^1(\Q, Z_{I_{\phi,\sconn}} )$, so the quotient by it is defined. The composite map $I_\phi^{\ad}(\Q) \to I_{\phi}^{\ad} (\A_f) \to \coh^1(\A_f, Z_{I_{\phi}}) $ factors through the coboundary map $I_\phi^{\ad}(\Q) \to \coh^1(\Q, Z_{I_{\phi,\sconn}})$ associated with the short exact sequence $1 \to Z_{I_{\phi,\sconn}}\to I_{\phi,\sconn} \to I_{\phi}^{\ad} \to 1$, whose image is $\D (Z_{I_{\phi,\sconn}} , I_{\phi, \sconn} ;\Q)$. The first statement follows. The other statements are all obvious.  
 \end{proof}
 
 \begin{lem}\label{lem:KHC theorem} We have $\D (Z_{I_{\phi,\sconn}} , I_{\phi, \sconn} ;\Q) = \Sha^{\infty} (\Q, Z_{I_{\phi,\sconn}})$.  
 \end{lem}
 \begin{proof}
 	By Lemma \ref{lem:I_phi at infty}, $I_{\phi}^{\ad} (\RR)$ is compact. Hence it is connected (see \cite[\S 24.6]{borel1991}), and the map $I_{\phi,\sconn}(\RR) \to I_{\phi}^{\ad}(\RR)$ is surjective. Hence $\D (Z_{I_{\phi,\sconn}} , I_{\phi, \sconn} ;\RR) $ is trivial, and therefore $\D (Z_{I_{\phi,\sconn}} , I_{\phi, \sconn} ;\Q) \subset  \Sha^{\infty} (\Q, Z_{I_{\phi,\sconn}})$.
 	
 	The reverse containment follows from the fact that $\Sha^\infty(\QQ, I_{\phi, \sconn})$ is identified with $\Sha^{\infty}_{\ab} (\QQ, I_{\phi,\sconn})$, which is trivial. 
 \end{proof}
 
 \subsection{} \label{subsubsec:cocycle relation}
 Let $\AM = \AM(G,X, p ,\cG)$\index[Oct2026-sorted-notation]{am 0 03 00 am g x p g@$\mathcal{AM}(G,X,p,\mathcal G)$} be the set of all admissible morphisms $\Qf \to \G_G$. 
 On this set we define an equivalence relation\index[Oct2026-sorted-notation]{0002 approx 1 00 00 approx@$\approx$} $\approx$ by declaring $\phi_1 \approx \phi_2$ if and only if $\phi_1^{\Delta}$ is $G(\ol \QQ)$-conjugate to $\phi_2 ^{\Delta}$. Clearly $\approx$ is coarser than the equivalence relation defined by conjugacy among admissible morphisms. By Lemma \ref{lem:twist by cocycle} and Proposition \ref{prop:twist adm morph}, we have $\phi_1 \approx \phi_2$ if and only if there exists a (necessarily unique) $\beta \in \Sha_G^{\infty}(\Q, I_{\phi_1})$ such that $\phi_2$ belongs to the conjugacy class $\phi_1^{\beta}$.
 
 For any $\phi \in \AM$, by Lemmas \ref{lem:cohomological H(phi)} and \ref{lem:KHC theorem}, and the observation that for $? \in \{\A_f, \A_f^p\}$ the subgroup $\D(Z_{I_\phi}, I_\phi; ?) \subset \coh^1(?, Z_{I_\phi})$ is the kernel of the composition $$\coh^1(?, Z_{I_\phi}) \To \coh^1(?, I_\phi) \cong \coh^1_{\ab}(?, I_\phi) = \coh^1(?, Z_{I_{\phi,\sconn}} \to Z_{I_\phi})$$ which depends only on the map $Z_{I_{\phi,\sconn}} \to Z_{I_\phi}$, we see that the abelian groups $\mathcal H(\phi)$ and $\E^p (\phi)$  depend on $\phi$ only via the complex $Z_{I_{\phi, \sconn}} \to Z_{I_{\phi}}$.  Now if $\phi_1, \phi_2 \in \AM$ are such that $\phi_1 \approx \phi_2$, then for any $g \in G(\overline \Q)$ such that $\Int (g) \circ \phi_1^{\Delta} = \phi_2^{\Delta}$, the $\overline \Q$-isomorphism $\Int(g) : G_{\overline \Q, \phi_1^{\Delta}} \isom G_{\overline \Q, \phi_2^{\Delta}}$ induces inner twistings $I_{\phi_1, \overline \Q} \isom I_{\phi_2, \overline \Q}$ and $I_{\phi_1, \sconn,\overline \Q} \isom I_{\phi_2, \sconn, \overline \Q}$ between $\Q$-groups. Clearly the equivalence classes (Definition \ref{defn:inner twistings}) of these inner twistings are independent of the choice of $g$. It follows that $Z_{I_{\phi_1}}, Z_{I_{\phi_1, \sconn}}$ are canonically identified with $Z_{I_{\phi_2}}, Z_{I_{\phi_2, \sconn}}$ respectively. We conclude that there are canonical isomorphisms of abelian groups \begin{align*}
 	\mathcal H(\phi_1) \isom \mathcal H(\phi_2) \quad \text{and} \quad  
 	\E^p(\phi_1) \isom \E^p (\phi_2)
 \end{align*} which depend only on $\phi_1$ and $\phi_2$. These isomorphisms commute with the natural surjections $\E^p(\phi_i) \to \cH(\phi_i)$, and satisfy the obvious cocycle relations for any triple $\phi_1, \phi _2, \phi_3 \in \AM$ such that $\phi_1 \approx \phi_2 \approx \phi_3$.
 
 We now view the assignments $\phi \mapsto \fkE^p(\phi), \phi \mapsto \mathcal H(\phi)$ as sheaves of abelian groups $\E^p, \cH$ on $\AM$ (as a discrete topological space). We have a surjective homomorphism $\E^p \to \cH$. By the above discussion, $\E^p$ and $ \cH$ descend to sheaves of abelian groups $\E^p_{\approx}$\index[Oct2026-sorted-notation]{e 0 04 00 e superscript p subscript approx@$\mathfrak E^p_{\approx}$} and $ \cH_{\approx}$\index[Oct2026-sorted-notation]{h 0 03 00 h subscript approx@$\mathcal H_{\approx}$} on the quotient space $\AM/{\approx}$. The surjective homomorphism $\E^p \to \cH$ also descends to a  surjective homomorphism $\E^p_{\approx}\to \cH_{\approx}$. 
 
 The quotient map $\AM \to \AM/{\approx}$ factors through $ \AM/\mathrm{conj}$, the set of conjugacy classes of admissible morphisms.  We let $\cH_\mathrm{conj}$\index[Oct2026-sorted-notation]{h 0 03 00 h subscript conj@$\mathcal H_{\mathrm{conj}}$} (resp.~$\E^p_{\mathrm{conj}}$)\index[Oct2026-sorted-notation]{e 0 04 00 e superscript p subscript conj@$\mathfrak E^p_{\mathrm{conj}}$} be the pull-back of $\cH_{\approx}$ (resp.~$\E^p_{\approx}$) to $\AM/\mathrm{conj}$.  
 
 \begin{defn}\label{defn:system of obstructions} For $\cF \in \set{\E^p, \cH}$, we denote by  $\Gamma(\cF)$ the group of global sections of the sheaf $\cF$ on $\AM$. For $\dertau \in \Gamma(\cF)$, we write $\dertau(\phi) \in \cF(\phi)$ for the value of $\dertau$ at each $\phi \in \AM$. We denote by $\Gamma(\cF)_0$ (resp.~$\Gamma(\cF)_1$) the subgroup of $\Gamma(\cF)$  consisting of those global sections that descend to global sections of $\cF_{\approx}$ over $\AM/{\approx}$ (resp.~global sections of $\cF_\mathrm{conj}$ over $\AM/\mathrm{conj}$). Thus we have $$\Gamma(\cF)_0 \subset \Gamma(\cF)_1 \subset \Gamma(\cF).$$\index[Oct2026-sorted-notation]{zz03 gamma 0 00 00 gamma f@$\Gamma(\cF)$}\index[Oct2026-sorted-notation]{zz03 gamma 0 00 00 gamma f subscript 0@$\Gamma(\cF)_0$}\index[Oct2026-sorted-notation]{zz03 gamma 0 00 00 gamma f subscript 1@$\Gamma(\cF)_1$} 
 \end{defn}
 
 \subsection{}\label{para:H(phi)}
 Let $\phi\in \AM$. By Lemmas \ref{lem:cohomological H(phi)} and \ref{lem:KHC theorem}, we have a natural map 
 \begin{multline}\label{eq:second coh real'n of H(phi)}
 	\mathcal H(\phi) \cong \D (Z_{I_{\phi}}, I_{\phi} ;\A_f) / \mathfrak D(Z_{I_{\phi,\sconn}} , I_{\phi,\sconn} ; \QQ)  \\ \To  \coh^1(\A_f, Z_{I_\phi}) / \ker \big(\Sha^{\infty}(\Q, Z_{I_{\phi}}) \to \Sha^{\infty}(\Q, I_\phi) \big),
 \end{multline} fitting into the  commutative diagram
 \begin{align}\label{commdiag:1}\xymatrixcolsep{5pc}
 	\xymatrix{  \E^p (\phi) \ar[r] ^{\sim} \ar[d] & \D(Z_{I_{\phi}}, I_{\phi} ;\A_f^p) \ar[d] \\    \mathcal H(\phi) \ar[r]^-{(\ref{eq:second coh real'n of H(phi)})} &   \coh^1(\A_f, Z_{I_\phi}) / \ker \big(\Sha^{\infty}(\Q, Z_{I_{\phi}}) \to \Sha^{\infty}(\Q, I_\phi) \big) }
 \end{align}
 where the left vertical arrow is the natural surjection and the right vertical arrow is induced by the inclusion $\coh^1(\A_f^p , Z_{I_{\phi}}) \to \coh ^1 (\A_f, Z_{I_{\phi}})$.
 
 \begin{defn}\label{defn:rectification} Let $\dertau  \in \Gamma(\cH)$, and $\underline \sigma \in \Gamma(\E^p)$.
 	\begin{enumerate}
 		\item We say that $\dertau$ is\index{tori-rational} \emph{tori-rational}, if for each $\phi \in \AM$ and for each maximal torus $T\subset I_\phi$, the image of $\dertau (\phi)$ is trivial under the composite map
 		\begin{multline}\label{eq:defn of tori rational}
 			\mathcal H( \phi) \xrightarrow{(\ref{eq:second coh real'n of H(phi)})} \coh^1(\A_f, Z_{I_\phi}) / \ker \big(\Sha^{\infty}(\Q, Z_{I_{\phi}}) \to \Sha^{\infty}(\Q, I_\phi) \big) \\  \To \coh ^1 (\A_f, T) / \Sha _{G} ^{\infty} (\QQ, T).	\end{multline}
 		\item  We say that $\underline \sigma$ is \emph{tori-rational}, if for each $\phi \in \AM $ and for each maximal torus $T \subset I_{\phi}$, the image of $\underline \sigma  (\phi)$ is trivial under the composite map
 		\begin{multline*}
 			\E^p (\phi)\cong  \D (Z_{I_{\phi}}, I_{\phi} ; \A_f^p) \subset \coh ^1 (\A_f^p, Z_{I_{\phi}} ) \\ \To \coh ^1 (\A_f^p , T)  \To \coh ^1 (\A_f^p , T) / \Sha_G^{\infty, p} (\Q, T). 
 		\end{multline*} 
 	\end{enumerate}
 	In the above, $\Sha_G ^{\infty} (\QQ, T)$ and $\Sha_G^{\infty, p} (\QQ, T)$ are defined as in \S \ref{subsubsec:notation of Sha} with respect to the composite embedding $T_{\Qbar} \hookrightarrow I_{\phi,\Qbar} \hookrightarrow G_{\Qbar}$, which is $\QQ$-rational at $\pi_1$-level by \S  \ref{para:transfer I_phi} and Lemma \ref{lem:rational pi_1}.
 \end{defn}
 The next lemma relates the two notions of tori-rationality.
 \begin{lem}\label{lem:about rect}
 	Let $\dertau\in \Gamma(\cH)$. The following statements are equivalent.
 	\begin{enumerate}
 		\item  $\dertau$ is tori-rational.
 		\item $\dertau$ is the image of a tori-rational element $\underline \sigma \in \Gamma(\E^p)$ under the natural  surjection  $\Gamma(\E^p) \to \Gamma(\cH)$.
 		\item Every $\underline \sigma \in \Gamma(\E^p)$ lifting $\dertau $ is tori-rational.
 	\end{enumerate}
 \end{lem}
 \begin{proof}
 	By the commutative diagram (\ref{commdiag:1}), (ii) implies (i). Obviously (iii) implies (ii). It remains to show that (i) implies (iii).
 	
 	Let $\underline \sigma \in \Gamma(\E^p)$ be a lift of $\dertau$. For each $\phi \in \AM$, we have a natural surjection $ \D (Z_{I_\phi}, I_{\phi} ; \A_f) \to \cH(\phi)$ as in Lemma \ref{lem:cohomological H(phi)}.  Fix an element $\epsilon_{\phi} \in \D (Z_{I_\phi}, I_{\phi} ; \A_f)$ lifting $\dertau(\phi)$. Then the image of $\dertau(\phi)$ under (\ref{eq:second coh real'n of H(phi)}) is represented by $\epsilon_{\phi}$.
 	
 	By the commutative diagram (\ref{commdiag:1}), the image of $\underline \sigma (\phi)$ in $ \coh ^1 (\A_f, Z_{I_{\phi}})$
 	equals the sum of $\epsilon_{\phi}$ and the image of some $\upsilon_{\phi}\in \Sha ^{\infty}_{I_{\phi}} (\QQ, Z_{I_{\phi}})$. For each maximal torus $T \subset I_{\phi}$, by tori-rationality of $\dertau$ there exists an element $\beta _{\phi ,T} \in \Sha ^{\infty} _G (\QQ, T)$ whose image in $\coh^1(\A_f, T)$ equals that of $\epsilon_{\phi}$. Let $\beta_{ \phi ,T}' \in \Sha ^{\infty} _G(\QQ, T)$ be the sum of $\beta _{\phi ,T}$ and the image of $\upsilon_{\phi}$ in $ \Sha ^{\infty} _{I_\phi} (\QQ, T) \subset \Sha _G^{\infty} (\QQ, T)$. Then the image of $\underline \sigma (\phi)$ under
 	$$\D (Z_{I_{\phi}}, I_\phi ; \A_f^p) \To \D (Z_{I_{\phi}}, I_\phi ; \A_f)   \To \coh ^1 (\A_f, T)$$ equals that of $\beta' _{\phi ,T}$. It follows that $\beta' _{\phi ,T} $ lies in $\Sha _G^{\infty, p} (\QQ, T)$, and that the image of $\underline \sigma (\phi)$ in $\coh^1(\A_f^p, T)/ \Sha_G^{\infty, p} (\Q, T)$ is trivial, as desired.
 \end{proof}
 
 \subsection{}\label{subsubsec:defn of S_tau(phi)} Let $\dertau \in \Gamma(\cH)_1$. For each admissible morphism $\phi: \Qf \to \G_G$, we set\index[Oct2026-sorted-notation]{s 0 00 00 s subscript underline tau phi@$S_{\underline {\tau}}(\phi)$}
 $$ S_{\dertau}(\phi) : = S_{\tau_\phi} (\phi), $$ where $\tau_\phi \in I_{\phi}^{\ad}(\A_f)$ is any lift of $\dertau(\phi) \in \cH(\phi)$, and $S_{\tau_\phi} (\phi)$ is defined as in \S \ref{subsubsec:defn of S(phi)}. The isomorphism class of the $\Phi^\ZZ \times G(\adele_f^p)$-set $S_{\dertau}(\phi)$ depends only on $\dertau(\phi)$, not on the choice of $ \tau_\phi$. Moreover, from the assumption that $\dertau \in \Gamma(
 \cH)_1$, it follows that this isomorphism class depends on $\phi$ only via its conjugacy class.
 
 We write
 $\mathsf{LR} (G,X, p, \cG, \dertau) $\index[Oct2026-sorted-notation]{lr 0 07 00 lr g x p g tau underline@$\mathsf{LR} (G,X, p, \cG, \dertau)$} for the modification of Conjecture \ref{conj:LR} where each $S(\phi)$ is replaced by $S_{\dertau}(\phi)$.
 
 Combined with Theorem \ref{thm:integral model} and Theorem \ref{thm:lan-stroh}, the main result of \cite{kisin2012modp} can be stated as follows.
 \begin{thm}
 	\label{thm:kisin}
 	Let $(G,X, p ,\cG)$ be an unramified Shimura datum such that $(G,X)$ is of abelian type. Assume $p>2$. Then there exists $\dertau \in \Gamma(\cH)_1$ such that the statement $\mathsf{LR} (G,X, p, \cG, \dertau)$ holds.
 \end{thm} 
 In the following conjecture, we impose better control of $\dertau$ than merely asking that it belongs to $\Gamma(\cH)_1$. 
 \begin{conj}[Langlands--Rapoport--$\tau$]\label{hypo about LR} For\index{Langlands--Rapoport tau Conjecture@Langlands--Rapoport--$\tau$ Conjecture} each unramified Shimura datum $(G,X, p, \cG)$, there exists a tori-rational element $\dertau \in \Gamma(\cH)_0$ such that the statement $\mathsf{LR} (G,X, p, \cG, \dertau)$ holds. Namely, there exists a smooth integral model $\Shh_{K_p}$ of $\Sh_{K_p}$ over $\oo_{E,(\fkp)}$ which has well-behaved $\coh^*_c$ and for which there is a bijection
 	$$
 	\Shh_{K_p} (\overline \FF_q)  {} \isom \coprod_{\phi }   S_{\dertau}(\phi)  $$
 	compatible with the actions of $G(\A_f^p)$ and the $q$-Frobenius $\Phi$. Here $\phi$ runs through a set of representatives for the conjugacy classes of admissible morphisms $\Qf \to \G_G$.
 \end{conj}
 Note that Conjecture \ref{hypo about LR} is weaker than Conjecture \ref{conj:LR}, as the latter asserts that $\dertau$ can be taken to be trivial. However, the former is a viable substitute for the latter for  applications to computing zeta functions and $\ell$-adic cohomology, in view of the following theorem. 
 
 \begin{thm}[cf.~Theorem \ref{introthm:one conj implies the other} in the Introduction]\label{thm:main thm in point counting}
 	Conjecture \ref{hypo about LR} implies Conjecture \ref{conj:point counting formula}. 
 \end{thm}
 We devote \S \ref{sec:study of pairs} to the proof of Theorem \ref{thm:main thm in point counting}. In Part \ref{part:3},  we shall prove Conjecture \ref{hypo about LR} for Shimura varieties of abelian type.
 
 \cleardoublepage
 \chapter{Langlands--Rapoport--\texorpdfstring{$\tau$}{tau} implies point counting} \label{sec:study of pairs}
 	In this chapter, we prove Theorem \ref{thm:main thm in point counting}, namely that the Langlands--Rapoport--$\tau$ Conjecture implies the point counting formula. The proof utilizes group-theoretic arguments modeled after \cite{langlands1987gerben}, but with a necessary modification to incorporate the family of twisting elements as in the Langlands--Rapoport--$\tau$ Conjecture.

 	Throughout this chapter we fix an unramified Shimura datum $(G,X,p,\cG)$, and keep the notations $E, \fkp, q=p^r$ as in \S \ref{para:unramified SD} and \S \ref{subsubsec:SD}. 

 \section{Semi-admissible and admissible Langlands--Rapoport pairs}
 The enumeration of fixed points via the Langlands--Rapoport(--$\tau$) Conjecture naturally leads to pairs $(\phi,\epsilon)$, where $\phi$ is an admissible morphism and $\epsilon \in I_{\phi}(\QQ)$. We define these as \emph{semi-admissible Langlands--Rapoport pairs}. We then introduce a \emph{$q^m$-admissibility}  condition for such a pair, which requires that a certain operator on the affine Grassmannian constructed from $(\phi, \epsilon)$ and the $q^m$-Frobenius possesses a fixed point. Finally, we derive some consequences of the admissibility condition. 
 \begin{defn}\label{defn:LR pair} By a\index{Langlands--Rapoport pair} \emph{Langlands--Rapoport pair} (\emph{LR pair}), we mean a pair $(\phi,\epsilon)$, where $ \phi: \Qf \to \G_G$ is a morphism in $\proredgb(\Qbar/\QQ)$, and $\epsilon$ is an element of $I_{\phi}(\QQ).$
 	We call such a pair\index{semi-admissible} $(\phi,\epsilon)$ \textit{semi-admissible}, if $\phi$ is admissible. We denote by $\wgp$\index[Oct2026-sorted-notation]{lrp 0 03 00 lrp@$\mathcal{LRP}$} the set of all LR pairs, and by $\wgp_{\sa}$\index[Oct2026-sorted-notation]{lrp 0 03 00 lrp subscript sa@$\mathcal{LRP}_{\mathrm{sa}}$} the subset of semi-admissible LR pairs.
 \end{defn}
 
 \begin{rem}\label{rem:epsilon semi-simple}If $(\phi,\epsilon) \in \wgp_{\sa}$, then $\epsilon$ is semi-simple, since $(I_\phi/Z_G) (\RR)$ is anisotropic by Lemma \ref{lem:I_phi at infty}.
 \end{rem}
 \subsection{}
 The group $G({\overline {\QQ}})$ acts on the set $\wgp$ by \emph{conjugation} in the following sense. If $(\phi, \epsilon) \in \wgp$ and $g\in G({\overline {\QQ}})$, then $\Int(g) (\phi,\epsilon) := (\Int( g) \circ \phi, \Int(g)\epsilon)$ is also an element of $\wgp$. We write\index[Oct2026-sorted-notation]{lrp 0 03 11 langle lrp rangle@$\langle \mathcal{LRP} \rangle$} $$ \langle \wgp \rangle : = \wgp/G(\overline \Q)\text{-conjugacy}.$$ For $(\phi,\epsilon) \in \wgp$, we denote by
 $$\langle \phi,\epsilon \rangle \in \langle \wgp \rangle $$ the $G(\overline \QQ)$-conjugacy class of $(\phi,\epsilon)$. The subset $\wgp_{\sa} \subset \wgp$ is stable under $G({\overline {\QQ}})$-conjugacy, and we write\index[Oct2026-sorted-notation]{lrp 0 03 11 langle lrp subscript sa rangle@$\langle \mathcal{LRP}_{\mathrm{sa}} \rangle$} $$\langle \wgp_{\sa} \rangle : =  \wgp_{\sa}/G(\overline \Q)\text{-conjugacy}. $$

 \subsection{}\label{p-adic realization}
 Let $(\phi, \epsilon)\in \wgp_{\sa}$. Let $\theta = \phi(p)\circ \zeta_p  : \G_p\to \G_G(p)$, and choose $g\in\mathcal{UR} (\theta)$ (see Definition \ref{defn:unram morph} and Lemma \ref{lem:unram criterion}). Let $b_g : = b_{\Int(g^{-1}) \circ \theta} \in G(\Q_p^{\ur})$ (see Definition \ref{defn:b_theta}), and let $\epsilon_g : = \Int(g^{-1}) (\epsilon) \in G(\Qpbar)$. Since $\epsilon$ is semi-simple (Remark \ref{rem:epsilon semi-simple}), so is $\epsilon_g$.
 
 We have $\epsilon_g \in I_{\Int(g)^{-1} \circ \theta} (\Q_p)$, and hence $\epsilon_g \in J_{b_g}(\Q_p)$ by Proposition \ref{prop:twist at p}. Also, by the same proposition, $b_g$ is decent. It then also follows that $\epsilon_g \in G(\Q_p^{\ur})$, as $J_{b_g}(\Q_p) \subset G(\Q_p^{\ur})$ (see \S \ref{para:J_b}). We let\index[Oct2026-sorted-notation]{cls 1 01 00 cls subscript p phi epsilon@$\mathrm{cls}_p(\phi,\epsilon)$} $$\cls_p(\phi,\epsilon) : = \set{ (b_g, \epsilon_g) \mid g \in \mathcal{UR}(\theta)}. $$
 Let $G(\Q_p^{\ur})$ act on $\cls_p(\phi,\epsilon)$ on the left by $h \cdot (b,\epsilon') : =  (h b \sigma (h)^{-1}, h \epsilon' h^{-1})$. Since $\mathcal {UR} (\theta)$ is a $G(\Q_p^{\ur})$-torsor, the  $G(\Q_p^{\ur})$-action on $\cls_p(\phi,\epsilon)$ is transitive.

 \subsection{} \label{admissible pairs} Fix a positive integer $m$, and let $n=mr$. (Recall that $q=p^r$ is the residue cardinality of $\fkp$.) We will define \emph{$q^m$-admissible LR pairs}, which will serve to describe $\FF_{q^m}$-points of the Shimura variety.
 
 Let $(\phi,\epsilon) \in \wgp_{\sa}$, and let $(b,  \epsilon')\in \cls_p(\phi,\epsilon)$. Recall from \S \ref{subsubsec:recall of defns in LR} and \S \ref{subsubsec:SD} that the set $X_p(\phi) $ is identified with the set $ X_{-\mu}(b )$, where $\mu \in \dmu_X^{\cG}$. The action of $\Phi$ on $ X_p(\phi)$ corresponds to the left multiplication by $\Phi_b:= (b\rtimes \sigma)^r$ on $X_{-\mu}(b)$. The action of $\epsilon$ on $X_p(\phi)$ corresponds to the left multiplication by $\epsilon'$ on $X_{-\mu}(b)$. Let\index[Oct2026-sorted-notation]{x 0 00 00 x subscript p phi epsilon q superscript m@$X_p(\phi,\epsilon, q^m)$} $$ X_p(\phi,\epsilon, q^m) : = \set{x\in X_p(\phi) \mid \epsilon x = \Phi^m x}.$$
 Then we have an identification\index[Oct2026-sorted-notation]{x 0 00 00 x subscript minus mu b epsilon prime q superscript m@$X_{-\mu}(b,\epsilon', q^m)$}
 \begin{align*}
 	X_p(\phi,\epsilon, q^m) \cong X_{-\mu}(b,\epsilon' ,q^m) := \set{x\in X_{-\mu}(b)\mid\epsilon' x = \Phi_b^m x   }.
 \end{align*}
 This motivates the following definition.

 \begin{defn}\label{Defn of Adm}
 	We say that an element $(\phi,\epsilon) \in \wgp_{\sa}$ is\index{admissible!LR pair} \emph{$q^m$-admissible}, if for one (and hence every) element $(b,\epsilon')$ of $\cls_p(\phi,\epsilon)$, we have
 	$$ \set{x\in G(\LL)/ \cG(\breve \ZZ_p) \mid\epsilon' x = \Phi_b ^m x} \neq \emptyset.  $$ We denote by\index[Oct2026-sorted-notation]{lrp 0 03 00 lrp subscript a q superscript m@$\mathcal{LRP}_{\mathrm{a}}(q^m)$} $\wgp_{\adm}(q^m)$ the set of $q^m$-admissible elements of $\wgp_{\sa}$. This subset is stable under $G(\overline \Q)$-conjugacy, and we write \index[Oct2026-sorted-notation]{lrp 0 03 11 langle lrp subscript a q superscript m rangle@$\langle \mathcal{LRP}_{\mathrm{a}}(q^m) \rangle$} $$ \lprod{\wgp_{\adm}(q^m)} : =  \wgp_{\adm}(q^m)/G(\overline \Q)\text{-conjugacy}. $$
 \end{defn}

 \begin{lem}\label{polar decomp of epsilon}
 	Let $(\phi,\epsilon) \in \wgp_{\adm}(q^m)$, and let $(b, \epsilon') \in \cls_p(\phi,\epsilon)$. Then there exists $t \in \ZZ_{\geq 1}$ satisfying the following conditions. \begin{enumerate}
 		\item The fractional cocharacter $t\nu_b$ is a cocharacter of $G$ defined over $\QQ_{p^t}$.
 		\item We have $\epsilon'^{t}= p^{nt\nu_b} k $, for $k$ lying in some conjugate of $\cG(\breve \ZZ_p)$ in $G(\LL)$.
 	\end{enumerate}
 \end{lem}
 \begin{proof}We have seen in \S \ref{p-adic realization} that $b$ is decent. Take $t$ such that $b$ is $t$-decent. Then condition (i) is already satisfied. Also by assumption there exists $x\in G(\LL)/ \cG(\breve \ZZ_p)$ such that $$ \epsilon' x = \Phi_b^m x .  $$ Since $\epsilon'$ commutes with $\Phi_b = (b\rtimes \sigma)^r$, we have $$ \epsilon'^{t} x = \Phi_b^{mt} x .$$ By the density of $G(\Qpur)$ in $G(\LL)$ for the $p$-adic topology, we can replace $t$ by a multiple, and assume that $ \sigma^t x = x.$ Let $s = tn =tmr$. Then $$ \epsilon'^{t} x = \Phi_b^{mt} x = (b\rtimes \sigma)^{s} x =  b\leftidx{^\sigma} b\leftidx{^{\sigma^2}}b\cdots \leftidx{^{\sigma^{s-1} }} b x. $$
 	Therefore
 	\begin{align}\label{epsilon and newton 2}
 		k := (  b\leftidx{^\sigma} b\leftidx{^{\sigma^2}}b\cdots \leftidx{^{\sigma^{s-1} }} b)^{-1 } \epsilon'^{t}  \in  x \cG(\breve \ZZ_p) x^{-1}.
 	\end{align} Finally, since $b$ is $t$-decent we have $k = p^{-nt\nu_b} {\epsilon'}^t$. This proves condition (ii).
 \end{proof}
 Recall from \S \ref{subsubsec:unique-quasi-motivic} that we have a homomorphism of $\Q_p$-group schemes $\nu(p): \mathbb D_{\QQ_p} \to Q_{\QQ_p}$.
 \begin{prop}\label{prop:equalities at p and infty}
 	Let $(\phi_1,\epsilon_1), (\phi_2, \epsilon_2) \in \wgp_{\adm}(q^m)$. Suppose $\epsilon_1 = \epsilon_2 $ as elements of $G({\overline {\QQ}})$.  Then $\phi_1^{\Delta}\circ \nu(p)= \phi_2^{\Delta}\circ \nu(p)$ as homomorphisms $\mathbb D_{\Qpbar} \to G_{\Qpbar}$.
 \end{prop}
 \begin{proof}
 	As in \S \ref{p-adic realization}, for $i=1,2$, we choose $g_i \in \mathcal {UR} (\phi_i(p)\circ \zeta_p )$, and let $(b_i, \epsilon_i')$ be the element of $\cls_p(\phi_i,\epsilon_i)$ associated with $g_i$. By Lemma \ref{polar decomp of epsilon}, we can find $t \in \NN$ such that $$  {\epsilon'_i}^{t}  =   p^{nt\nu_{b_i}} k_i, $$ where $k_i$ lies in some $G(\LL)$-conjugate of $\cG(\breve \ZZ_p)$, for $i = 1,2$. 
 	Since ${\epsilon'_i}$ commutes with $b_i \rtimes \sigma$, we know that ${\epsilon'_i}^t$ commutes with $nt\nu_{b_i}$. Also $k_i$ lies in a bounded subgroup of $G(\LL)$. Applying Lemma \ref{uniqueness of polar decomp} to $F = \LL$, we see that $nt\nu_{b_i}$ is the unique cocharacter $\nu$ of $G$ over $\overline {\breve \QQ}_p$ commuting with ${\epsilon_i'}^t$ such that ${\epsilon_i'}^{-t} p^{\nu}$ lies in a $G(\overline {\breve \QQ}_p)$-conjugate of a bounded subgroup of $G(\LL)$. 
 	Let $g =g_2^{-1} g_1 \in G (\Qpbar)$. Then we have $\epsilon'_2 = \Int(g) \epsilon'_1$. By the above-mentioned uniqueness of $nt\nu_{b_i}$ with respect to ${\epsilon'_i}^t$, we have  
 	\begin{align}\label{eq:conjugacy of two Newton cochars}
 		\nu_{b_2} = \Int(g) \circ \nu_{b_1}.
 	\end{align}
 	By Proposition \ref{prop:twist at p} (i) we have
 	\begin{align}\label{eq:apply newton=theta}
 		-\nu_{b_i} =  \left( \Int(g_i^{-1}) \circ \phi_i(p) \circ \zeta_p   \right)^{\Delta}    = \Int(g_i^{-1}) \circ \phi_i^\Delta \circ \nu(p).
 	\end{align}
 	Comparing (\ref{eq:conjugacy of two Newton cochars}) with (\ref{eq:apply newton=theta}), we have $\phi_1^{\Delta}\circ \nu(p)=  \phi_2^{\Delta}\circ \nu(p)$ as desired.
 \end{proof}
 
 Recall from \S \ref{subsubsec:unique-quasi-motivic} that we have a homomorphism of $\RR$-group schemes $\nu(\infty): \GG_m \to Q_{\RR}$. 
 \begin{lem}\label{lem:Newton at infty}
 	Let $\phi_1, \phi_2 : \gQ \to \gG_G$ be two admissible morphisms. Then $\phi_1 ^{\Delta} \circ \nu (\infty) = \phi_2 ^{\Delta} \circ \nu(\infty)$.
 \end{lem}
 \begin{proof}
 	By condition (ii) in Definition \ref{defn of adm morphism} with $v =\infty$, we know that $\phi_i (\infty) \circ \zeta_{\infty} $ is conjugate to $\xi_\infty$. Hence $(\phi_i(\infty) \circ \zeta_{\infty} )^\Delta = \phi_i^\Delta \circ \nu(\infty) $ is conjugate to ${\xi}_{\infty}^{\Delta}$. By the definition of $\xi_{\infty}$ in \cite[\S 3.3.5]{kisin2012modp},  ${\xi}_{\infty}^{\Delta}$ is equal to the weight cocharacter for the Shimura datum $(G,X)$, which is central in $G$ (see \cite[\S 2.1.1]{deligne1979varietes}). The lemma follows.
 \end{proof}
 
 \begin{cor}\label{transporter}Let $(\phi_1,\epsilon_1), (\phi_2, \epsilon_2) \in \wgp_{\adm}(q^m)$. For $i=1,2$, assume that $\phi_i$ is gg (Definition \ref{defn:gg morph}), and that $\epsilon_i$ lies in $G(\QQ)$. Then for all $g\in G({\overline {\QQ}})$ such that $\Int(g) (\epsilon_1 )= \epsilon_2$, we have $\Int(g)\circ \phi_1^\Delta = \phi_2^\Delta.  $
 \end{cor}
 \begin{proof} Write $H := \set{   g\in G({\overline {\QQ}}) \mid\Int(g) \epsilon_1 = \epsilon_2 }$. Let $g\in H$.
 	For all $\tau \in \Gamma$, we have $\leftidx{^\tau} g\in H$, since $\epsilon_1, \epsilon_2 \in G(\QQ)$.
 	Applying Proposition \ref{prop:equalities at p and infty} and Lemma \ref{lem:Newton at infty} to $(\phi_2,\epsilon_2)$ and $\Int(\leftidx{^\tau} g) (\phi_1,\epsilon_1)$, we get
 	\begin{align}\label{gripe session 2}
 		\phi_2^\Delta \circ \nu(v) = \Int(\leftidx{^\tau} g) \circ \phi_1^\Delta \circ \nu(v), \quad \text{for }v = p,\infty.
 	\end{align}
 	Now we apply $\tau^{-1} $ to both sides of (\ref{gripe session 2}). As $\phi_1^\Delta$ and $\phi_2^\Delta$ are defined over $\QQ$ (by the gg assumption), we get
 	\begin{align}\label{gripe session}
 		\phi_2^\Delta \circ \leftidx{^{\tau^{-1}}}\nu(v) = \Int(g) \circ \phi_1^\Delta\circ  \leftidx{^{\tau^{-1}}} \nu(v), \quad \text{for } v = p,\infty.
 	\end{align}
 	
 	By construction, for each finite Galois extension $L/\QQ$ contained in $\overline \Q$, the $\Q$-torus $Q^L$ is split over $L$, and the $\QQ$-vector space $X_*(Q^L) \otimes_{\ZZ} \QQ$ is generated by the $\Gal(L/\QQ)$-conjugates of the fractional cocharacters $$\mathbb D_{\Qpbar} \xrightarrow{\nu(p)}  Q_{\Qpbar} \To Q^L_{\Qpbar} $$ and $$\GG_{m,\CC} \xrightarrow{\nu(\infty)} Q_\CC \To Q^L_{\CC}. $$ See \cite[\S 3.1]{kisin2012modp} for details. Hence (\ref{gripe session}) for all $\tau$ implies that $\phi_2^\Delta =  \Int(g) \circ \phi_1^\Delta $.
 \end{proof}
 
 \section{The g\"unstig gelegen condition}\label{certain nice}
 In \S \ref{subsec:preparation LR} we introduced g\"unstig gelegen morphisms $\phi: \mathfrak Q \to \mathfrak G_G$. In this section, we define an analogous condition for a Langlands--Rapoport pair $(\phi,\epsilon)$. For g\"unstig gelegen pairs $(\phi,\epsilon)$, the natural setting for cohomological twisting involves the \emph{connected} centralizer $I_{\phi,\epsilon}^0$ rather than $I_{\phi,\epsilon}$. We show that every $G(\Qbar)$-conjugacy class of semi-admissible or admissible Langlands--Rapoport pairs contains a g\"unstig gelegen pair. Finally, we examine how the cohomological twisting of g\"unstig gelegen pairs interacts with the semi-admissibility and admissibility conditions.

 As in \S \ref{para:choice of q_tau}, for each $\tau \in \Gamma$ we choose a lift $q_{\tau}  \in \Qf^{\topo}$. For $(\phi,\epsilon) \in \wgp$, we write $I_{\phi,\epsilon}^0$ for the connected centralizer $(I_{\phi})_{\epsilon}^0$.
 
 \begin{defn}\label{defn:gg}
 	We say that an LR pair $(\phi,\epsilon) \in \wgp$ is\index{gunstig gelegen@g\"unstig gelegen (gg)!LR pair} \textit{g\"unstig gelegen} (\textit{to be abbreviated as gg})\footnote{As in Definition \ref{defn:gg morph}, this terminology comes from \cite[\S 5]{langlands1987gerben}, but our definition is modified to suit the corrected definition of $\Qf$.}, if the following conditions are satisfied:
 	\begin{itemize}
 		\item The embedding $I_{\phi,\overline \Q} \hookrightarrow G_{\overline \Q}$ maps $\epsilon \in I_{\phi}(\Q)$ into $G(\Q)$. Moreover, $\epsilon$ is $\RR$-elliptic in $G$, cf.~\S \ref{para:fkp}.  
 		\item For each $\tau \in \Gamma, $ write $\phi(q_\tau) = g_\tau \rtimes \tau,$ with $g_{\tau} \in G(\overline \Q)$. Then $g_\tau $ lies in $(G_{\overline \Q, \phi^{\Delta}})_{\epsilon}^0 (\overline \Q) = I_{\phi, \epsilon}^0 (\overline \Q)$.
 	\end{itemize}
 	We denote by $\ggp$\index[Oct2026-sorted-notation]{lrp 0 03 00 lrp superscript gg@$\mathcal{LRP}^{\mathrm{gg}}$} the set of gg LR pairs.
 \end{defn}
 
 \begin{rem} Let $(\phi,\epsilon)\in \wgp$, and let $\phi(q_\tau) = g_\tau \rtimes \tau$, for $\tau \in \Gamma$. If one changes the choice of $q_{\tau}$, then $g_\tau$ is left multiplied by a $\overline \Q$-point of $\im(\phi^\Delta)$, which is a central torus in $G_{\overline \Q, \phi^{\Delta}}$. Hence the second condition in Definition \ref{defn:gg} is independent of the choice of $q_\tau$.
 \end{rem}
 \subsection{}
 \label{para:gg and gg}	Let $(\phi, \epsilon) \in \ggp$. By Lemma \ref{lem:gg morph},  $\phi$ is gg, and the canonical $\overline \Q$-isomorphism $I_{\phi,\overline \Q} \isom G_{\phi^{\Delta}, \overline \Q}$ is an inner twisting between the $\Q$-groups $I_{\phi}$ and $G_{\phi^{\Delta}}$. Moreover, in the current case this inner twisting restricts to an inner twisting between the $\Q$-groups $I_{\phi,\epsilon}$ and $(G_{\phi^{\Delta}})_{\epsilon}$, and an inner twisting between the $\Q$-groups $I_{\phi,\epsilon}^0$ and $(G_{\phi^{\Delta}})_{\epsilon}^0$.
 
 In fact, the inner twisting between $I_{\phi,\epsilon}^0$ and $(G_{\phi^{\Delta}})_{\epsilon}^0$ can be interpreted as follows. Let $I_0 = G_{\epsilon}^0 \subset G$. Since $(\phi,\epsilon)$ is gg, $\phi$ factors as $\Qf \to \G_{I_0} \to \G_G$. We write $\phi_{I_0}$ for $\phi$ when we view it as a morphism $\Qf \to \G_{I_0}$. Then $\phi_{I_0}$ is itself gg. Hence by Lemma \ref{lem:gg morph} applied to $\phi_{I_0}$, we obtain an inner twisting between $I_{\phi_{I_0}}$ and $(I_0)_{(\phi_{I_0})^{\Delta}} = (I_0)_{\phi^{\Delta}}$. It is easy to see that as $\Q$-groups we have $I_{\phi_{I_0}} = I_{\phi ,\epsilon}^0$ and $(I_0)_{\phi^{\Delta}} = (G_{\phi^{\Delta}})_{\epsilon}^0$.
 
 \begin{defn}\label{defn:vergessen} We write
 	$ [\ggp] $\index[Oct2026-sorted-notation]{lrp 0 03 10 lrp superscript gg leftbracket rightbracket@$[\mathcal{LRP}^{\mathrm{gg}}]$} for the quotient set of $\ggp$ divided by the equivalence relation of $G(\Qbar)$-conjugacy.\footnote{We caution the reader that the subset $\ggp \subset \wgp$ is not stable under $G(\overline \Q)$-conjugacy.} For $(\phi,\epsilon) \in \ggp$, we denote its image in $[\ggp]$ by $[\phi,\epsilon]$.  We denote the natural injection $[\ggp] \to \lprod{\wgp}$ by $\mathfrak v$\index[Oct2026-sorted-notation]{v 1 04 00 v@$\mathfrak v$} (standing for \textit{vergessen}).
 \end{defn}
 
 \begin{lem}
 	\label{lem:twisting gg pair}
 	Let $(\phi,\epsilon) \in \wgp$, and let $a = (a_{\tau})_{\tau} \in Z^1(\Q, I_{\phi,\epsilon})$. Then $(a\phi ,\epsilon) \in \wgp$ (i.e., the element $\epsilon \in G(\Qbar)$ lies in the image of $I_{a\phi}(\Q) \hookrightarrow I_{a\phi}(\Qbar) \hookrightarrow G(\Qbar)$). Here $a\phi $ is the twist of $\phi$ by the cocycle $a$ as in \S \ref{para:twist morphism}. If in addition we have $(\phi,\epsilon) \in \ggp$ and $a \in Z^1(\Q, I_{\phi,\epsilon}^0)$, then $(a\phi ,\epsilon) \in \ggp$.
 \end{lem}
 \begin{proof} First recall that $(a\phi)^{\Delta} = \phi^{\Delta}$, and $I_{\phi, \Qbar }$ and $I_{a\phi , \Qbar}$ are equal as $\Qbar$-subgroups of $G_{\Qbar}$. Write $\phi(q_{\tau}) = g_{\tau} \rtimes \tau$ for each $\tau \in \Gamma$. Then we have $(a\phi) (q_{\tau}) = a_{\tau} g_{\tau} \rtimes \tau$. Since $\epsilon \in I_{\phi}(\Q)$, we have $g_{\tau} \lix^{\tau} \epsilon g_{\tau}^{-1} = \epsilon.$ Since $\epsilon$ commutes with $a_{\tau}$, we have
 	$$ a_{\tau} g_{\tau}  \lix^{\tau} \epsilon
 	g_{\tau}^{-1} a_{\tau}^{-1} = \epsilon, $$ which means that $\epsilon \in I_{a\phi}(\Q)$. Thus we have shown that $(a\phi, \epsilon) \in \wgp$.
 	
 	To show the second statement, we need to check that $a_{\tau} g_{\tau} \in I_{a\phi, \epsilon}^0 (\Qbar)$ for all $\tau$. But both $g_{\tau}$ and $a_{\tau}$ lie in $I_{\phi, \epsilon}^0 (\Qbar)$, and we have $I_{a\phi, \epsilon}^0 (\Qbar) = I_{\phi, \epsilon}^0 (\Qbar)$. The desired statement follows.
 \end{proof}

 \begin{lem}\label{lem:factos thru coh}
 	Let $(\phi,\epsilon) \in \wgp$. Let $a, b \in Z^1(\QQ, I_{\phi, \epsilon})$. Then $a,b$ are cohomologous in $I_{\phi,\epsilon}$ if and only if $\lprod{a \phi,\epsilon} = \lprod{b\phi,\epsilon}$.
 \end{lem}
 \begin{proof}
 	Write $\phi(q_\tau) =g_\tau \rtimes \tau$, for each $\tau \in \Gamma$. If $a,b$ are cohomologous in $I_{\phi,\epsilon}$, then there exists $u\in I_{\phi, \epsilon}({\overline {\QQ}}) $ such that
 	\begin{align}\label{eq:abequiv}
 		a_\tau  = u^{-1} b_\tau g_\tau \leftidx{^\tau} u g_{\tau}^{-1}, \quad \forall \tau \in \Gamma.\end{align}
 	Here $\lix^{\tau} u$ denotes the action of $\tau$ on $u$ viewed as in $G(\overline \Q)$. Then $$ (a \phi)(q_\tau ) =  a_\tau g_\tau \rtimes \tau =  u^{-1} b_\tau g_\tau \leftidx{^\tau}u \rtimes \tau  =\Int(u^{-1}) \circ (b \phi)(q_\tau) .  $$  Hence we have $\Int(u ^{-1}) \circ (b\phi) = a \phi$ as they already agree on the kernel. Also $  \Int(u^{-1})\epsilon =\epsilon. $ Therefore $(a\phi, \epsilon) = \Int(u^{-1}) (b\phi, \epsilon)$, and so $\lprod{a \phi,\epsilon} = \lprod{b\phi,\epsilon}$.
 	
 	Conversely, assume that $\lprod{a \phi,\epsilon} = \lprod{b\phi,\epsilon}$. Then there exists $u \in G(\overline \Q)$ such that $(a\phi, \epsilon) = \Int(u^{-1}) (b\phi, \epsilon)$. Since $(a\phi)^{\Delta} = (b\phi)^{\Delta}$, we have $u \in I_{\phi,\epsilon}(\overline \Q)$. Now the relation $\Int(u^{-1}) \circ (b\phi) = a\phi$ is equivalent to (\ref{eq:abequiv}), which shows that $a$ and $b$ are cohomologous.
 \end{proof}
 \subsection{} \label{para:iota and eta}
 Let $(\psi ,\delta) \in \wgp$ and $(\phi,\epsilon) \in \ggp$.
 In view of Lemma \ref{lem:twisting gg pair} and Lemma \ref{lem:factos thru coh}, we have well-defined maps
 \begin{align*}
 	Z^1(\QQ, I_{\psi, \delta}  ) \To \wgp, \quad a \longmapsto (a\psi, \delta)
 \end{align*} and
 \begin{align*}
 	Z^1(\QQ, I_{\phi, \epsilon}^0  ) \To \ggp , \quad   a \longmapsto (a\phi, \epsilon),
 \end{align*} which induce maps\index[Oct2026-sorted-notation]{zz09 iota 1 00 00 iota subscript psi delta@$\iota_{\psi,\delta}$}
 $$\iota_{\psi,\delta} :  \coh^1(\QQ, I_{\psi , \delta} ) \To  \lprod{\wgp} $$
 and\index[Oct2026-sorted-notation]{zz07 eta 1 00 00 eta subscript phi epsilon@$\eta_{\phi,\epsilon}$} $$ \eta_{\phi,\epsilon}:
 \coh^1(\QQ, I_{\phi,\epsilon}^0 ) \To  [\ggp]
 $$ respectively. By Lemma \ref{lem:factos thru coh} the map $\iota_{\psi,\delta}$ is injective.
 
 The image of $\eta_{\phi,\epsilon}$ depends on $(\phi,\epsilon)$ only via $[\phi,\epsilon] \in [\ggp]$, and we shall denote it by $\mathcal C_{[\phi,\epsilon]}$.\index[Oct2026-sorted-notation]{c 0 03 00 c subscript phi epsilon leftbracket rightbracket@$\mathcal C_{[\phi,\epsilon]}$} 
 A straightforward computation with cocycles shows that subsets of this form  constitute a partition of $ [\ggp]$. 
 \begin{lem}\label{lem:strong conj automatic}
 	Let $(\phi_1,\epsilon_1), (\phi_2,\epsilon_2)\in \ggp$ be elements having the same image in $[\ggp],$ and let $g\in G({\overline {\QQ}})$ be such that $\Int(g)(\phi_1,\epsilon_1) = (\phi_2,\epsilon_2)$. Then we have $g~ \leftidx{^\tau}g^{-1} \in G_{\epsilon_2}^0 (\overline \Q) $ for all $\tau \in \Gamma$. In particular, $\epsilon_1, \epsilon_2 \in G(\Q)$ are stably conjugate.
 \end{lem}
 \begin{proof}
 	For $i = 1,2$, we write $\phi_i (q_\tau) = g_{i,\tau} \rtimes \tau$. Then $$ g_{2,\tau} = g g_{1,\tau} \lix^{\tau} g^{-1} = g g_{1,\tau} g^{-1} g~ \lix^{\tau} g^{-1},$$ and hence $$g ~\lix^{\tau} g^{-1} = (g g_{1,\tau} g^{-1})^{-1} g_{2,\tau}. $$ 
 	Since $(\phi_1,\epsilon_1), (\phi_2,\epsilon_2)\in \ggp$, we have $g_{1,\tau} \in G_{\epsilon_1}^0(\Qbar)$ and $g_{2,\tau} \in G_{\epsilon_2}^0(\Qbar)$. It follows that $g g_{1,\tau} g^{-1} \in G_{\epsilon_2}^0(\Qbar)$. Hence $g ~\lix^{\tau} g^{-1}\in G_{\epsilon_2}^0(\Qbar)$, as desired.
 \end{proof} 
 
 \begin{defn}\label{defn:stc} As in \S \ref{para:Frob KP}, let $\SigmaREllip(G)$ be the set of stable conjugacy classes of semi-simple, $\RR$-elliptic elements of $G(\QQ)$. We define the \emph{stable conjugacy class map}\index[Oct2026-sorted-notation]{scc 1 04 00 scc@$\mathfrak{scc}$} \begin{align*}
 		\F: [\ggp] \To \SigmaREllip (G),
 	\end{align*}
 	sending $[\phi,\epsilon]$ to the stable conjugacy class of $\epsilon$. This is well defined by Lemma \ref{lem:strong conj automatic}.
 \end{defn}
 
 \begin{defn} Fix $m$ as in \S \ref{admissible pairs}.	We set\index[Oct2026-sorted-notation]{lrp 0 03 00 lrp superscript gg subscript sa@$\mathcal{LRP}^{\mathrm{gg}}_{\mathrm{sa}}$} \index[Oct2026-sorted-notation]{lrp 0 03 00 lrp superscript gg subscript a q superscript m@$\mathcal{LRP}^{\mathrm{gg}}_{\mathrm{a}}(q^m)$} \index[Oct2026-sorted-notation]{lrp 0 03 10 lrp superscript gg subscript sa leftbracket rightbracket@$[\mathcal{LRP}^{\mathrm{gg}}_{\mathrm{sa}}]$}\index[Oct2026-sorted-notation]{lrp 0 03 10 lrp superscript gg subscript a q superscript m leftbracket rightbracket@$[\mathcal{LRP}^{\mathrm{gg}}_{\mathrm{a}}(q^m)]$} $$\ggp_{\sa} : = \ggp \cap \wgp_{\sa}, \quad \ggp_{\adm}(q^m): = \ggp \cap \wgp_{\adm}(q^m),  $$ $$ [\ggp_{\sa}] : = \mathfrak v^{-1} (\lprod{\wgp_{\sa}}), \quad [\ggp_{\adm} (q^m)] : = \mathfrak v^{-1} (\lprod{\wgp_{\adm} (q^m)}). $$ See Definition \ref{defn:vergessen} for $\mathfrak v$. 
 \end{defn}
 
 \subsection{}\label{useful remark}
 Let $(\phi,\epsilon) \in \ggp_{\adm}(q^m)$. Then $\phi$ is gg (see \S \ref{para:gg and gg}), and $\epsilon \in G(\QQ)$. Hence by Corollary \ref{transporter} the $\overline \Q$-inclusion $I_{\phi,\overline \Q} \hookrightarrow G_{\overline \Q}$ induces a $\overline \Q$-isomorphism $I_{\phi,\epsilon , \overline \Q} \isom G_{\epsilon, \overline \Q}$. By this fact and by the discussion in \S \ref{para:gg and gg}, we obtain canonical inner twistings $$
 (I^0_{\phi,\epsilon})_{{\overline {\QQ}}} \isom (G^0_{\epsilon})_{\overline{\QQ}}  \quad \text{and} \quad
 (I_{\phi,\epsilon})_{{\overline {\QQ}}} \isom (G_{\epsilon})_{\overline{\QQ}}. $$
 
 \begin{lem}\label{lem:suff gg pair}
 	Let $(\phi,\epsilon) \in \wgp$. Assume that $I_{\phi}$ contains a $G$-rational maximal torus $T$ (see Definition \ref{defn:G-rational tori}) such that $\epsilon \in T(\Q)$ and such that $T/Z_G$ is anisotropic over $\RR$. Then  $(\phi,\epsilon) \in \ggp$.
 \end{lem}
 \begin{proof} The assumptions imply that $\epsilon$ is a semi-simple $\RR$-elliptic element of $G(\QQ)$. By the proof of Lemma \ref{lem:suff gg}, if we write $\phi(q_{\tau}) = g_{\tau} \rtimes \tau$, then $g_{\tau}$ lies in $T(\Qbar)$ for each $\tau \in \Gamma$. Note that $T \subset I_{\phi, \epsilon} ^0$. Hence $g_{\tau}$ lies in $I_{\phi, \epsilon}^0 (\Qbar)$ as desired. \end{proof}
 \begin{lem}\label{lem:G-rational torus for pair}Let $(\phi,\epsilon) \in \wgp_{\sa}$. Then there exists $g\in G(\Qbar)$ such that $(\phi',\epsilon') : = \Int(g) (\phi,\epsilon) \in \wgp_{\sa}$ satisfies the following condition: The group $I_{\phi'}$ contains a $G$-rational maximal torus $T'$ such that $\epsilon' \in T'(\Q)$ and such that $T'/Z_G$ is anisotropic over $\RR$. Moreover, in this case we have $(\phi',\epsilon') \in \ggp_{\sa}$.
 \end{lem}
 \begin{proof} Choose a maximal torus $T \subset I_{\phi}$ defined over $\Q$ such that $\epsilon \in T(\Q)$. (This is possible since $\epsilon$ is semi-simple; see Remark \ref{rem:epsilon semi-simple}.)  Since $\phi$ is admissible, Proposition \ref{prop:gg morph} implies that there exists $g\in G(\overline \Q)$ such that $\Int(g)(T)$ is a $G$-rational maximal torus in $I_{\Int(g) \circ \phi}$. We write $(\phi',\epsilon') $ for $\Int(g) (\phi,\epsilon)$, and write $T'$ for $\Int(g) (T)$. Then $T'$ is a $G$-rational maximal torus in $I_{\phi'}$, and $\epsilon' \in T'(\Q)$. Since $\phi'$ is admissible, we know that $T'/Z_G$ is anisotropic over $\RR$ by Lemma \ref{lem:I_phi at infty}. Finally, we have $(\phi',\epsilon')\in \ggp_{\sa}$ by Lemma \ref{lem:suff gg pair}.
 \end{proof}
 
 \begin{cor}\label{cor:sa conj to gg}
 	The injection $\mathfrak v : [\ggp] \to \lprod{\wgp}$ restricts to bijections $[\ggp_{\sa}] \isom \lprod{\wgp_{\sa}}$ and $[\ggp_{\adm}(q^m)] \isom \lprod{\wgp_{\adm}(q^m)}$.
 \end{cor}
 \begin{proof}
 	By Lemma \ref{lem:G-rational torus for pair}, every $G(\Qbar)$-conjugacy class in ${\wgp_{\sa}}$ contains an element of $\ggp_{\sa}$. Hence $\mathfrak v$ restricts to a bijection $[\ggp_{\sa}] \isom \lprod{\wgp_{\sa}}$. Since $[\ggp_{\adm}(q^m)]$ is by definition $ \mathfrak v^{-1} (\lprod{\wgp_{\adm} (q^m)})$, we see that $\mathfrak v$ also restricts to a bijection $[\ggp_{\adm}(q^m)] \isom \lprod{\wgp_{\adm}(q^m)}$. \end{proof}
 \subsection{}\label{para:defn of Sha_G^infty for gg pair}
 Let $(\phi,\epsilon) \in \ggp_{\sa}.$ As $I_{\phi,\epsilon}^0$ is a $\Q$-subgroup of $I_{\phi}$, by the discussion in \S \ref{para:transfer I_phi} and Lemma \ref{lem:rational pi_1} the canonical embedding $(I_{\phi,\epsilon}^0)_{\Qbar} \hookrightarrow G_{\Qbar}$ is $\QQ$-rational at $\pi_1$-level. We then define $\Sha_G^{\infty}(\Q, I_{\phi,\epsilon}^0)$ as in \S \ref{subsubsec:notation of Sha}.
 
 \begin{prop}\label{semi-adm locus} For $(\phi,\epsilon) \in \ggp_{\sa}$, the inverse image of $ [\ggp_{\sa}]$ under the map $\eta_{\phi,\epsilon}$  in \S \ref{para:iota and eta} is $ \Sha^{\infty} _G(\QQ,I_{\phi,\epsilon}^0). $  
 \end{prop}
 \begin{proof}
 	By Proposition \ref{prop:twist adm morph}, an element $\beta \in \coh^1(\Q,I_{\phi,\epsilon}^0)$ lies in $\eta_{\phi,\epsilon}^{-1} ([\ggp_{\sa}])$ if and only if the image of $\beta$ in $\coh^1(\Q, I_{\phi})$ lies in $\Sha_G^{\infty} (\Q, I_{\phi})$. Thus we only need to check that the map $\coh^1(\R, I_{\phi,\epsilon}^0) \to \coh^1(\R, I_\phi)$ has trivial kernel. Note that $(I_{\phi,\epsilon}^0)_\RR$ and $(I_{\phi})_{\RR}$ are both reductive groups over $\RR$, and their centers both contain $Z_{G,\RR}$. The desired statement then follows from Lemma \ref{lem:I_phi at infty} and \cite[Lem.~4.4.5]{kisin2012modp}.	
 \end{proof}

 \begin{prop}\label{one cell}
 	Let $x\in [\ggp_{\adm}(q^m)]$. Write $e$ for $\F (x) \in \SigmaREllip(G)$, and let $\mathcal C_x \subset [\ggp]$ be as in \S \ref{para:iota and eta}. Then $$
 	\mathcal C_x  \cap [\ggp_{\adm}(q^m)] = \F^{-1} (e) \cap [\ggp_{\adm}(q^m)].$$
 \end{prop}
 \begin{proof}We only need to show the containment $\F^{-1} (e) \cap [\ggp_{\adm}(q^m)] \subset \mathcal C_x$. Let $y \in \F^{-1} (e) \cap [\ggp_{\adm}(q^m)]$ . Write $x= [ \phi_0,\epsilon_0], y= [ \phi_1, \epsilon_1]$, for some $(\phi_i, \epsilon_i) \in \ggp_{\adm}(q^m)$, $i=0,1$. Since $\F(x) =\F(y)$, there exists $g\in G({\overline {\QQ}})$ such that
 	\begin{align}\label{same cell -1} \Int(g)\epsilon_1 & = \epsilon_0 , \\
 		\label{same cell 0} g \leftidx{^\tau}g^{-1} &\in G_{\epsilon_0}^0 (\Qbar), \quad   \forall \tau\in \Gamma.
 	\end{align}
 	
 	Let $\phi' = \Int(g) \circ \phi_1$. By (\ref{same cell -1}), we have $(\phi', \epsilon_0) \in \wgp$, and $(\phi',\epsilon_0)$ is $G(\Qbar)$-conjugate to $(\phi_1,\epsilon_1)$. By Corollary \ref{transporter} applied to the relation (\ref{same cell -1}), we have $\phi'^\Delta = \phi_0^\Delta$. By Lemma \ref{lem:twist by cocycle} we have $\phi' = a \phi_0$ for some $a \in Z^1(\QQ, I_{\phi_0})$. For $i = 0,1$, we write $
 	\phi_i(q_\tau) = g^{(i)}_\tau \rtimes \tau.$
 	Then
 	\begin{align}\label{same cell 1}
 		\phi'(q_\tau) =a_\tau g^{(0)}_{\tau} \rtimes \tau.
 	\end{align}
 	By the definition of $\phi'$ we have\begin{align}\label{same cell 3}
 		\phi'(q_\tau) = \Int(g) \circ \phi_1 (q_{\tau}) = g (  g^{(1)}_\tau \rtimes \tau ) g^{-1} = g g^{(1)}_\tau \leftidx{^{\tau}} g^{-1}  \rtimes \tau.
 	\end{align}
 	Comparing (\ref{same cell 1}) and (\ref{same cell 3}), we have \begin{align}\label{eq:formula for a_rho}
 		a_\tau = g g^{(1)}_\tau \leftidx{^{\tau}} g^{-1} (g^{(0)}_\tau) ^{-1} = \Int(g)(g_{\tau}^{(1)}   ) \cdot (g\leftidx{^{\tau}}g^{-1}) \cdot  (g^{(0)}_\tau) ^{-1}.
 	\end{align}
 	
 	Since $(\phi_i,\epsilon_i)$ is gg, we have $g_{\tau}^{(i)} \in G_{\epsilon_i}^{0}(\Qbar)$. Hence $\Int(g) (g_{\tau}^{(1)})$ and $(g^{(0)}_\tau) ^{-1}$ both lie in $G_{\epsilon_0}^{0}(\Qbar)$. Thus by (\ref{same cell 0}) and (\ref{eq:formula for a_rho}), we have $a_\tau \in G_{\epsilon_0}^{0}(\Qbar)$. By the discussion in \S \ref{useful remark}, we have $(G_{\epsilon_0}^0)_{\Qbar} = (I_{\phi_0,\epsilon_0}^0)_{\Qbar}$, since $(\phi_0,\epsilon_0) \in \ggp_{\adm}(q^m)$. Hence $a= (a_\tau)$ is a cocycle in $Z^1(\Q, I_{\phi_0,\epsilon_0}^0 )$. It follows from Lemma \ref{lem:twisting gg pair} that the pair $(\phi', \epsilon_0) = (a\phi_0 , \epsilon_0)$ is gg. Since this pair is $G(\Qbar)$-conjugate to $(\phi_1,\epsilon_1)$, we have $y = [\phi',\epsilon_0] = [a\phi_0, \epsilon_0] \in \im(\eta_{\phi_0,\epsilon_0}) = \mathcal C_x $.
 \end{proof}
 
 \section{Admissible morphisms and maximal tori}\label{sec:adm morphisms and maximal tori}
 In this section, we show that every \emph{special point datum} for $(G,X)$ gives rise to an admissible morphism $\phi$ together with a maximal torus $T \subset I_\phi$ and, conversely, that every pair $(\phi, T)$ arises in this way.  
 
 \begin{defn}\label{defn:special point datum} By a\index{special point datum} \textit{special point datum} for $(G,X)$, we mean a triple $(T, i, h)$, where $T$ is a torus over $\QQ$, $i: T \to  G$ is an injective $\Q$-homomorphism whose image is a maximal torus in $G$, and $h: \mathbb S  \to T_{\RR} $ is an $\RR$-homomorphism such that $i\circ h \in X$. We denote the set of special point data by $\spd(G,X)$.\index[Oct2026-sorted-notation]{spd 0 03 00 spd g x@$\mathcal{SPD}(G,X)$}
 \end{defn}
 \begin{defn}\label{defn:special morphisms}
 	Let $(T,i,h)\in \spd(G,X)$. Let $\mu_h \in X_*(T)$ be the Hodge cocharacter associated with $h$. We denote by $\phi(T,i,h)$\index[Oct2026-sorted-notation]{zz21 phi 1 00 00 phi t i h@$\phi(T,i,h)$} the composite morphism $$ \Qf \xrightarrow{\Psi_{T, \mu_h}} \G_T \overset{i}{\To} \G_G $$ in $\proredgb(\Qbar/\Q)$. (See \S \ref{para:Psi_mu} for $\Psi_{T,\mu_h}$.)
 \end{defn}
 The following result considerably strengthens Corollary \ref{cor:gg morph}.
 \begin{thm}[{\cite[Lem.~3.5.8,  Thm.~3.5.11]{kisin2012modp}}] \label{thm:kisin special}
 	A morphism $\phi: \Qf \to \G_G$ is admissible if and only if there exists $(T,i,h) \in \spd(G,X)$ such that $\phi$ is conjugate to $\phi(T,i,h)$.
 \end{thm}
 \begin{rem}\label{rem:special morph is gg}
 	Let $(T,i,h) \in \spd(G,X)$, and let $\phi  = \phi(T,i,h)$. Then $i(T_{\Qbar})$ is contained in $I_{\phi,\Qbar}$ (when they are both viewed as subgroups of $G_{\Qbar}$), and the inclusion $T_{\Qbar} \hookrightarrow I_{\phi ,\Qbar}$ is defined over $\Q$. In other words, $T$ is naturally a $G$-rational maximal torus in $I_{\phi}$. By Lemma \ref{lem:suff gg}, $\phi$ is gg. With this observation, we see that Theorem \ref{thm:kisin special} strengthens Corollary \ref{cor:gg morph}.
 \end{rem}
 
 \subsection{}	\label{para:bar mu}
 Our next goal is to prove a more precise version of the ``only if'' direction in Theorem \ref{thm:kisin special}. We introduce some general notation.  For a $\Q_p$-torus $T$, we write $X_*(T)_{\QQ}$ for $X_*(T)\otimes _{\ZZ} \QQ$. For a cocharacter $\lambda \in X_*(T)$, we write $\bar \lambda^T$, or simply $\bar \lambda$, for $$[L:\Q_p]^{-1} \sum_{\tau \in \Gal(L/\Q_p)} \tau (
 \lambda) \in X_*(T)_{\QQ},$$ where $L/\Q_p$ is any finite Galois extension over which $\lambda$ is defined. For an unramified reductive group $M$ over $\Q_p$, we write
 \begin{align*}
 	w_M : M(\LL) \To \pi_1(M)_{\Gamma_{p,0}} = \pi_1(M)
 \end{align*}
 for the Kottwitz homomorphism associated with the $p$-adic valuation on $\LL$, as in \S \ref{subsec:Kottwitz hom} and \S \ref{para:B(G)}. (Here $\Gamma_{p,0}$ acts trivially on $\pi_1(M)$ since $M$ is unramified over $\Q_p$.)
 
 We fix $\mu \in \dmu_X^{\cG}$ as in \S \ref{subsubsec:SD}. For each $b \in G(\LL)$, we define $X_{-\mu}(b)$ as in \S \ref{subsubsec:recall of defns in LR}.

 \begin{lem}\label{lem:LR 5.11}
 	Let $T \subset G_{\Q_p}$ be a maximal torus over $\Q_p$. Let $b \in T(\LL) 
 	\subset G(\LL)$ be such that $X_{-\mu}(b)$ is non-empty. Then there exists $\mu_T \in X_*(T)$ which is $G(\Qpbar)$-conjugate to $\mu$ and such that $\nu_b$ is equal to $-\overline {\mu_T}^T$ as elements of $X_*(T)_ \QQ$. 
 \end{lem}
 \begin{proof}
 	We recall the argument from the proof of  \cite[Lem.~5.11]{langlands1987gerben}, with our notation.
 	
 	First note that $\nu_b$ is a fractional cocharacter of $T$ defined over $\Q_p$, which holds for any $b\in T(\breve \Q_p)$. Let $M$ be the centralizer in $G_{\Q_p}$ of the maximal $\Q_p$-split subtorus of $T$. Then $M$ is a $\Q_p$-Levi subgroup of $G_{\Q_p}$ containing $T$, and $\nu_b$ factors through the center of $M$. Up to conjugating $T$ and $b$ by an element of $G(\Q_p)$, we may assume that $M$ contains a maximal torus $T'$ that is the centralizer of (the generic fiber of) a maximal $\ZZ_p$-split torus in $\cG$. In particular, there is a connected reductive model $\cM$ of $M$ over $\ZZ_p$ such that the embedding $M\hookrightarrow G_{\Q_p}$ extends to $\cM \hookrightarrow \cG$. Without loss of generality, we may also assume that $\mu$ is a cocharacter of $T'$ defined over $E_{\mathfrak p} = \QQ_{p^r}$.
 	
 	Write $\Omega$ and $\Omega_M$ for  the absolute Weyl groups of $G$ and $M$. By Cartan decomposition we have maps
 	\begin{align*}
 		c_{\cG} : & ~ G(\LL) \To \cG(\breve \ZZ_p) \backslash G(\LL)/ \cG(\breve \ZZ_p) \isom \Omega \backslash X_*(T'), \\
 		c_{\cM} : & ~M(\LL) \To \cM(\breve \ZZ_p) \backslash M(\LL)/ \cM(\breve \ZZ_p) \isom \Omega_M \backslash X_*(T'),	
 	\end{align*} which lift $w_G$ and $w_M$ (cf.~Corollary \ref{cor:kottCartan}).
 	
 	By the non-emptiness of $X_{-\mu}(b)$, there exists $x \in G(\LL)$ such that $$c_{\cG} (x^{-1} b \sigma ( x)) = \Omega \cdot (-\mu). $$ It is shown on p.~178 of \cite{langlands1987gerben} that this condition implies the existence of $\mu_1 \in \Omega \cdot \mu$ and  $m \in M(\LL)$ such that
 	\begin{align}\label{eq:adm in M}
 		c_{\cM} (m^{-1} b  \sigma( m ) ) = \Omega_M \cdot  (-\mu_1).\end{align}
 	(The argument uses the Iwasawa decomposition and the fact that $\mu$ is minuscule, cf.~also the proof of \cite[Lem.~2.2.2]{kisin2012modp}.)

 	Choose $\mu_T \in X_*(T)$ conjugate to $\mu_1\in X_*(T')$ under $M(\Qpbar)$. Note that $\overline {\mu_T}^T$ factors through the maximal $\Q_p$-split subtorus of $T$, and is therefore central in $M$. We have seen that $\nu_b$ is also central in $M$. Hence in order to check $\nu_b = -\overline {\mu_T}^T$, it suffices to compare their images in $\pi_1(M)_{\Q} :  = \pi_1(M) \otimes_{\Z} \Q$. We may assume that $b$ is decent in $T(\LL)$ by $\sigma$-conjugating it in $T(\LL)$ (see \S \ref{isocdefns}); relation (\ref{eq:adm in M}) is preserved by adjusting $m$. Then for sufficiently divisible $n$ we have
 	$$
 	b  \sigma (b) \cdots {\sigma^{n-1}} (b) = p^{n \nu_b}.
 	$$
 	Write $[\cdot]$ for the map $ X_*(T')_\Q \to \pi_1(M)_{\Q}$. We have  \begin{align*} [n\nu_b] =
 		w_M (p^{n\nu_b}) & = w_M(b \sigma (b) \cdots \sigma^{n-1}(b)) \\
 		&  = w_M(m) - w_M(\sigma ^{n} (m)) - [\mu_1 + \sigma (\mu_1) + \cdots + \sigma^{n-1}(\mu_1)]. \end{align*} The last equality follows from (\ref{eq:adm in M}), and the $\sigma$-action on $X_*(T')$ is defined since $T'$ is unramified.  We can choose $n$ such that the coset $m \cM(\breve \Z_p) \in M(\LL)/\cM(\breve \ZZ_p)$ is fixed by $\sigma^n$ (by the density of $M(\Qpur)$ in $M(\LL)$), and such that $T'$ splits over $\Q_{p^n}$. Then the above relation becomes
 	$$ [\nu_b]= - [\ol {\mu_1}^{T'}].$$
 	Finally, $[\ol {\mu_1} ^{T'}]$ is equal to the image of $\overline{\mu_T}^T$ in $\pi_1(M)_{\Q}$. This is because $\mu_1$ and $\mu_T$ have the same image in $\pi_1(M)$, and both maps $X_*(T)\to \pi_1(M)$ and $X_*(T') \to \pi_1(M)$ are Galois equivariant. Thus we conclude that $\nu_b$ and $-\overline{\mu_T}^T$ have the same image in $\pi_1(M)_\Q$, as desired.
 \end{proof}
 \begin{lem}\label{lem:LR 5.12}
 	Let $T\subset G$ be a maximal torus over $\Q$ such that $T_{\RR }$ is elliptic in $G_{\RR}$. Let $i$ denote the inclusion $T \hookrightarrow G$. Let $\mu_T \in X_*(T)$ be such that $i\circ \mu_T$ lies in $\dmu_X(\Qbar)$. Then there exist $u \in G(\Qbar)$ and an $\RR$-homomorphism $h: \mathbb S \to T_{\RR}$, satisfying the following conditions. 
 	\begin{enumerate}
 		\item $i' : = \Int(u) \circ i : T_{\Qbar} \to G_{\Qbar}$ is defined over $\Q$.
 		\item $\mu_T = \mu_h$.	\item  $i' \circ h \in X$.
 	\end{enumerate}
 \end{lem}
 \begin{proof}
 	The proof of this result is given inside the proof of  \cite[Lem.~5.12]{langlands1987gerben}. In fact, in that lemma $\mu_T$ is assumed to be of the form $\omega \mu_{h_0}$, where $\omega$ is an element of the absolute Weyl group of $(G,T)$, and $h_0$ is an $\RR$-homomorphism $\mathbb S \to T_{\RR}$ such that $i \circ h_0 \in X$. We explain why our hypothesis implies that setting. Since $T_{\RR}$ is elliptic in $G_{\RR}$, there indeed exists an $\RR$-homomorphism $h_0: \mathbb S \to T_{\RR}$ such that $i\circ h_0 \in X$. Then $i\circ \mu_{h_0}$ and $i\circ \mu_T$ are conjugate by $G(\CC)$, and it follows that $\mu_T = \omega \mu_{h_0}$ for some $\omega$ in the absolute Weyl group.
 \end{proof}
 \subsection{}\label{para:admissible cochar}
 Let $\phi : \Qf \to \G_G$ be an admissible morphism. Let $T \subset I_{\phi}$ be a maximal torus over $\Q$. Since $(\phi(p) \circ \zeta_p)^{\Delta}$ is a central fractional cocharacter of $I_{\phi, \Qpbar}$, it can be viewed as an element of $X_*(T)\otimes \Q$. We say that a cocharacter $\mu_T \in X_*(T)$ is\index{admissible!phi-admissible@$\phi$-admissible}\index{phi-admissible@$\phi$-admissible|see{admissible}} \emph{$\phi$-admissible}, if the composition $\GG_{m,\Qbar} \xrightarrow{\mu_T} T_{\Qbar} \hookrightarrow I_{\phi,\Qbar} \hookrightarrow G_{\Qbar}$ is a cocharacter in $\dmu_X(\Qbar)$, and if $\overline {\mu_T}^{T_{\Q_p}} = (\phi(p) \circ \zeta_p)^{\Delta} $ as elements of $X_*(T)\otimes \Q$.

 \begin{thm}[{\cite{langlands1987gerben}}]\label{thm:precise speciality}
 	Let $\phi: \Qf \to \G_G$ be an admissible morphism, and let $T\subset I_{\phi}$ be a maximal torus over $\Q$. The following statements hold.
 	\begin{enumerate}
 		\item There exists $\mu_T \in X_*(T)$ that is $\phi$-admissible in the sense of \S \ref{para:admissible cochar}.
 		\item Let $\mu_T \in X_*(T)$ be as in (i).   Then there exists a special point datum of the form $(T,i,h) \in \spd(G,X)$, satisfying the following conditions:
 		\begin{enumerate}
 			\item We have $\mu_T = \mu_h$.
 			\item There exists $g\in G(\Qbar)$ such that $\Int(g)\circ \phi = \phi(T,i,h)$, and such that the embedding $i : T \to G$ equals the composition $$   T_{\Qbar} \hookrightarrow I_{\phi,\Qbar}  \xrightarrow{\Int(g)}    I_{\phi(T,i,h), \Qbar} \hookrightarrow G_{\Qbar}  .$$
 			(Here the first two maps are defined over $\Q$, and the second map is an isomorphism.)
 		\end{enumerate}
 	\end{enumerate}	
 \end{thm}
 \begin{proof} This essentially follows from the proof of \cite[Satz 5.3]{langlands1987gerben}. We reproduce the argument for the convenience of the reader, and we remove the assumption in \textit{loc.~cit.} that $G_{\der}$ is simply connected. 
 	
 	\textbf{(i)} By Proposition \ref{prop:gg morph}, we may assume that $T$ is a $G$-rational maximal torus in $I_{\phi}$. By Lemma \ref{lem:suff gg}, $\phi$ factors through $\G_T \subset \G_G$. Hence the $G(\Q_p^{\ur})$-torsor $\mathcal {UR}(\phi(p) \circ \zeta_p)$ contains elements of $T(\Qpbar)$. We choose $t \in \mathcal {UR}(\phi(p) \circ \zeta_p) \cap T(\Qpbar)$, and let $b = b_{\Int(t^{-1}) \circ \phi(p)\circ \zeta_p} \in G(\Qpur)$ (see Definition \ref{defn:b_theta}).  Then $b \in T(\Qpur)$. Since $\phi$ is admissible, we have $X_{-\mu}(b) \neq \emptyset$ for $\mu \in \dmu_X^{\cG}$. Hence by Lemma \ref{lem:LR 5.11} we find $\mu_T \in X_*(T)$ such that $\overline {\mu_T}^{T_{\Q_p}} = - \nu_b$. Finally, by Proposition \ref{prop:twist at p} (i), we have $ - \nu_b = (\Int(t)^{-1} \circ \phi(p) \circ \zeta_p)^{\Delta}$, which equals $(\phi(p) \circ \zeta_p)^{\Delta}$. Hence $\mu_T$ is $\phi$-admissible.
 	
 	\textbf{(ii)} Again by Proposition \ref{prop:gg morph}, we may assume that $T$ is a $G$-rational maximal torus in $I_{\phi}$. We denote by $i_0$ the inclusion $T_{\Qbar} \hookrightarrow I_{\phi, \Qbar} \hookrightarrow G_{\Qbar}$, which is defined over $\Q$. By Lemma \ref{lem:I_phi at infty}, $i_0(T_{\RR})$ is an elliptic maximal torus in $G_{\RR}$. Applying Lemma \ref{lem:LR 5.12}, we find $u \in G(\Qbar)$ and an $\RR$-homomorphism $h : \mathbb S\to T_{\RR}$ such that $\Int(u) \circ i_0 : T \to G$ is still defined over $\Q$, such that $\mu_T = \mu_h$, and such that $\Int(u) \circ i_0 \circ h \in X$.
 	
 	We may replace $\phi$ and $i_0$ by $\Int(u) \circ \phi$ and $\Int(u) \circ i_0$ respectively, and assume that $u =1$. This does not change the property that $T$ is $G$-rational, and we now have $(T,i_0,h) \in \spd(G,X)$ such that $\mu_T = \mu_h$.
 	
 	By Lemma \ref{lem:suff gg}, we may factorize $\phi$ uniquely as $\phi = i_0 \circ \phi_T$, where $\phi_T$ is a morphism $\Qf \to \G_T$. Analogously, $\phi(T,i_0,h)$ factors as $i_0 \circ \Psi_{T,\mu_T}$ by its definition. We claim that $\phi_T^{\Delta} = \Psi_{T,\mu_T}^{\Delta}$, or equivalently $\phi^{\Delta} = \phi(T,i_0,h)^{\Delta}$. In fact, the $\phi$-admissibility of $\mu_T$ implies that
 	\begin{align}\label{eq:comp Newton at p}
 		(\phi_T(p ) \circ \zeta_p)^{\Delta} = (\Psi_{T,\mu_T}(p) \circ \zeta_p)^{\Delta}
 	\end{align}
 	(see \cite[(3.1.11)]{kisin2012modp}). Also, since $\phi$ and $\phi(T,i_0,h)$ are both admissible (the former by assumption and the latter by the ``if'' direction in Theorem \ref{thm:kisin special}), we have
 	\begin{align}\label{eq:comp Newton at infty}
 		(\phi_T(\infty) \circ \zeta_\infty)^{\Delta} = (\Psi_{T,\mu_T}(\infty) \circ \zeta_\infty)^{\Delta}
 	\end{align}
 	by Lemma \ref{lem:Newton at infty}. Since $\phi_T^{\Delta}$ and $\Psi_{T,\mu_T}^{\Delta}$ are both homomorphisms $\Qf^{\Delta} \to T_{\Qbar}$ defined over $\Q$ (for instance by Lemma \ref{lem:suff gg}), the above relations (\ref{eq:comp Newton at p}) and (\ref{eq:comp Newton at infty}) imply that $\phi_T^{\Delta} = \Psi_{T,\mu_T}^{\Delta}$ as desired (cf.~the discussion on the tori $Q^L$ in the proof of Corollary \ref{transporter}).
 	
 	By the claim above and by Lemma \ref{lem:twist by cocycle}, we have $\phi_T = a \Psi_{T,\mu_T}$, for some element $a \in Z^1(\Q, I_{\Psi_{T,\mu_T}}) = Z^1(\Q, T)$. By the admissibility of $\phi$ and $\phi(T,i_0,h)$ and by Proposition \ref{prop:twist adm morph}, the image of $a$ under $$Z^1(\Q, T) \To \coh^1(\Q, T) \To \coh^1(\Q, I_{\phi(T,i_0,h)})$$ lies in $\Sha_G^{\infty}(\Q, I_{\phi(T,i_0,h)})$. By Lemma \ref{lem:I_phi at infty} and \cite[Lem.~4.4.5]{kisin2012modp}, the map $\coh^1(\R, T) \to \coh^1(\R, I_{\phi(T,i_0,h)})$ has trivial kernel. Thus the class of $a$ in $\coh^1(\Q,T)$ lies in the kernel of the map $\Sha^{\infty} (\Q, T) \to \Sha^{\infty} (\Q, G)$ induced by $i_0: T \to G$. It follows that there exists $g \in G(\Qbar)$, satisfying:
 	\begin{align} \label{eq:twist T1}
 		i_0(a_{\tau}) = g^{-1} \lix^{\tau} g, \quad \forall \tau \in \Gamma ;\\ \label{eq:twist T2}
 		g \in   G(\RR)  \cdot i_0(T(\CC)) \subset G(\CC).
 	\end{align}
 	For each $\tau \in \Gamma$, we write $\phi(T,i_0,h) (q_{\tau})= g_{\tau} \rtimes \tau $, and $\phi(q_{\tau}) = g'_{\tau} \rtimes \tau$. (Here $q_{\tau} \in \Qf$ is a lift of $\tau$, as always).
 	Because $i_0(a_{\tau})$ and $g_{\tau}$ both lie in $i_0(T(\Qbar))$ and commute with each other, we have 	$g'_{\tau} = i_0(a_\tau) g _{\tau} = g_{\tau} i_0(a_{\tau}) = g_{\tau} g^{-1} \lix^{\tau}g$. Then we have
 	\begin{align}\label{eq:twist T3}(\Int(g) \circ \phi) (q_{\tau}) = g g'_{\tau} \lix^{\tau}g^{-1} \rtimes \tau = g g_{\tau} g^{-1} \rtimes \tau.
 	\end{align}
 	
 	Now let $i : = \Int(g) \circ i_0$. By (\ref{eq:twist T1}), $i$ is a $\Q$-embedding $T \to G$. By (\ref{eq:twist T2}), $i\circ h \in X$. In particular, $(T,i,h) \in \spd(G,X)$. By (\ref{eq:twist T3}), we have $(\Int(g) \circ \phi ) (q_{\tau}) = \phi(T,i,h)(q_{\tau})$. Since we also have $(\Int(g) \circ \phi ) ^{\Delta}= \phi(T,i,h)^{\Delta}$ (because $\phi_T^{\Delta} = \Psi_{T,\mu_T}^{\Delta}$), we have $\Int(g) \circ \phi = \phi(T,i,h)$. Thus $(T,i,h)$ and $g$ are the desired elements. 	
 \end{proof}
 
 \section{Admissible stable conjugacy classes}
 
 The goal of this section is to identify sufficient conditions for a stable conjugacy class in $\SigmaREllip(G)$ to lie in the image of  $[\ggp_{\adm}(q^m)] \subset [\ggp]$ under $\F: [\ggp] \to \SigmaREllip(G)$ (see Definition \ref{defn:stc}). We show in Theorem \ref{thm:*epsilon} that  \emph{$p^n$-admissibility} is sufficient. We then show in Proposition \ref{prop:sufficient condition for *epsilon} that the non-vanishing of a twisted orbital integral is sufficient to guarantee $p^n$-admissibility.

 For $M$ a reductive group over $\QQ_p$, we write $w_M: M(\LL) \to \pi_1(M)_{\Gamma_{p,0}}$ for the Kottwitz homomorphism, as in \S \ref{para:B(G)}. For each $k \in \ZZ_{ \geq 1}$, we denote by $\B^{(k)}(M)$\index[Oct2026-sorted-notation]{b 0 01 00 b superscript k m@$\mathrm B^{(k)}(M)$} the set of $\sigma^k$-conjugacy classes in $M(\LL)$.
 \begin{lem}\label{lem:ker of lambda}
 	Let $M$ be a  reductive group over $\QQ_p$. Let $x \in M(\LL)$ be in the kernel of $w_M$. For all $k \in \ZZ_{ \geq 1}$, the class of $x$ in $\B^{(k)}(M)$ is in the image of the natural map $\B^{(k)}(M_{\sconn}) \to \B^{(k)}(M)$.
 \end{lem}
 \begin{proof} Write $\tau: =  \sigma^k$. We first assume that $M_{\der } = M_{\sconn}$. Consider $M^{\ab} = M / M_{\der}$, which is a  torus over $\QQ_p$. The image $\bar x$ of $x$ in $M^{\ab} (\LL)$ is in the kernel of $w_{M^{\ab}}$, and this kernel is the unique parahoric subgroup of $M^{\ab} (\LL)$; see \cite[Rmk.~2.2 (iii)]{rapoportguide}. By Greenberg's theorem \cite[Prop.~3]{Greenberg} applied to the parahoric group scheme (i.e.~the connected N\'eron model of $M^{\ab}$), we have $\bar x = \bar c\cdot  \leftidx^{\tau}  \bar c^{-1}$ for some $\bar c\in M^{\ab} (\LL)$. Let $c \in M(\LL)$ be a lift of $\bar c$, which exists because $\coh^1 (\LL, M_{\der})$ is trivial by Steinberg's theorem. Then $c^{-1} x \leftidx^{\tau} c \in M_{\der} (\LL)$, which means that the class of $x$ in $\B^{(k)}(M)$ comes from $\B^{(k)}(M_{\der}) = \B^{(k)}(M_{\sconn})$.
 	
 	In the general case, as in \cite[\S 1]{Kot82} we take a $z$-extension
 	\begin{align*}
 		1 \To Z \To H \To M \To 1
 	\end{align*}
 	over $\QQ_p$, where $H$ is a reductive group with simply connected derived subgroup, and $Z$ is an induced torus contained in the center of $H$. Then we have a commutative diagram with exact rows:
 	\begin{align*}
 		\xymatrix{ 1  \ar[r]  & Z(\LL) \ar[d] ^{w_Z} \ar[r] & H(\LL) \ar[r] \ar[d] ^{w_H}  & M(\LL) \ar[r] \ar[d]^{w_M} & 1 \\ 1 \ar[r] &   \pi_1(Z) _{\Gamma_{p,0}}  \ar[r] & \pi_1(H)  _{\Gamma_{p,0}} \ar[r] & \pi_1(M) _{\Gamma_{p,0}} \ar[r] & 1 } \end{align*}
 	By \cite[\S 7]{kottwitzisocrystal2}, the vertical arrows in the above diagram are surjective. Since $w_M(x)= 0$, there exists $h \in \ker (w_H) \subset H(\LL)$ that maps to $x \in M(\LL)$. Applying the first part of the proof to $H$, we know that the class of $h$ in $\B^{(k)}(H)$ comes from $\B^{(k)}(H_{\sconn})$. But the composite map  $H_{\sconn} \to H  \to M$ factors through $M_{\sconn}$. The lemma follows.
 \end{proof}
 \begin{lem}\label{lem:compact elt}
 	Let $M$ be a  reductive group over $\QQ_p$, and let $T \subset M$ be an elliptic maximal torus. Let $x \in T(\Q_p)$ be such that $w_M(x) =0 \in \pi_1(M)_{\Gamma_{p,0}}$. Then $x$ lies in a compact subgroup of $T(\Q_p)$.
 \end{lem}
 \begin{proof}
 	Let $\cT^{\circ}$ be the connected N\'eron model of $T$ over $\ZZ_p$. Then the kernel of $w_T: T(\LL) \to X_*(T)_{\Gamma_{p,0}}$ is $\cT^{\circ} (\breve \ZZ_p)$; see \cite[Rmk.~2.2 (iii)]{rapoportguide}. If $w_T(x^k) = 0$ for some positive integer $k$, then $x^k$ lies in $ \cT^{\circ } (\breve \ZZ_p) \cap T(\Q_p) = \cT^{\circ} (\Z_p)$, which is a compact subgroup of $T(\Q_p)$. It then follows that $x$ lies in a compact subgroup of $T(\Q_p)$. Hence to prove the lemma it  suffices to show that $w_T(x^k) = 0$ for some positive integer $k$. We know that $w_T$ maps $T(\Q_p) \subset T(\LL)$ into the group of $\sigma$-invariants $(X_*(T)_{\Gamma_{p,0}})^{\sigma}$. It remains to show that the natural map
 	$$(X_*(T)_{\Gamma_{p,0}})^{\sigma} \hookrightarrow X_*(T)_{\Gamma_{p,0}} \To \pi_1(M)_{\Gamma_{p,0}}$$ has torsion kernel. 
 	
 	For this, let $A$ be the maximal $\Q_p$-split subtorus of $Z_M$. Since $T$ is elliptic in $M$, we know that $A$ is also the maximal $\Q_p$-split subtorus of $T$. In particular $X_*(A)$ is identified with $X_*(T)^{\Gamma_p}$. The projection $X_*(T)^{\Gamma_{p,0}} \to X_*(T)_{\Gamma _{p,0}}$ induces an isomorphism $X_*(T)^{\Gamma_{p,0}} \otimes \QQ \isom X_*(T)_{\Gamma _{p,0}} \otimes \QQ$ (with inverse given by taking average). This isomorphism is $\sigma$-equivariant, and the group of $\sigma$-invariants of the left hand side  is naturally identified with $X_*(A) \otimes \QQ$. Hence the natural map $X_*(A)\otimes \QQ\to (X_*(T)_{\Gamma_{p,0}})^{\sigma} \otimes \QQ$ is an isomorphism. It remains to show that the map $X_*(A) \to \pi_1(M)_{\Gamma_{p,0}}$ has torsion kernel. The map $X_*(A) \to \pi_1(M)$ is injective since $A$ is central in $M$, and the image is contained in $\pi_1(M)^{\Gamma_{p,0}}$. Hence it remains to show that the map $\pi_1(M)^{\Gamma_{p,0}} 
 	\to \pi_1(M)_{\Gamma_{p,0}}$ has torsion kernel. This is true since after tensoring with $\QQ$ the map becomes an isomorphism. 
 \end{proof}
 \begin{defn}\label{defn:M(epsilon)}
 	For $\epsilon \in G(\QQ_p)_{\semi}$, we let $M(\epsilon)$\index[Oct2026-sorted-notation]{m 0 00 00 m epsilon@$M(\epsilon)$} be the Levi subgroup of $G_{\QQ_p}$ that is the centralizer of the maximal $\QQ_p$-split subtorus of the center of  $G_{\epsilon} ^0$.  Equivalently, $M(\epsilon)$ is the smallest Levi subgroup of $G_{\Q_p}$ containing $G_{\Q_p, \epsilon}^0$. \end{defn}
 \begin{defn}\label{defn:*epsilon}
 	Let $\epsilon \in G(\QQ_p)_{\semi}$, and let $M : = M(\epsilon)$. Let $n = mr$ as in \S \ref{admissible pairs}. We say that $\epsilon$ is\index{admissible!semi-simple element} \emph{$p^n$-admissible}, if there exists a cocharacter $\mu_M$ of $M_{\Q_{p^n}}$ satisfying the following conditions.
 	\begin{itemize}
 		\item $\mu_M \in \dmu_X(\Q_{p^n})$. (Here $\dmu_X(\Q_{p^n})$ is well defined since $E_{\fkp}  =\Q_{p^r} \subset \Q_{p^n}$.)
 		\item  We have\footnote{Since $M$ is a Levi subgroup of $G_{\QQ_p}$, it is unramified. Hence the target of $w_M$ is $\pi_1(M)_{\Gamma_{p,0}} = \pi_1(M)$.} \begin{align}\label{(*epsilon)}
 			w_M(\epsilon) = - \sum_{i = 0}^{n-1} \sigma ^i [\mu_M]^M \in \pi_1(M),
 		\end{align}
 		where $[\mu_M]^M$ denotes the image of $\mu_M$ in $\pi_1 (M)$.
 	\end{itemize} 
 \end{defn}
 The following lemma generalizes \cite[Lem.~5.17]{langlands1987gerben} to the case where $G_{\der}$ is not necessarily simply connected.
 \begin{lem}\label{lem:*epsilon is invariant}
 	The set of $p^n$-admissible elements in $G(\QQ_p)_{\semi}$ is invariant under stable conjugacy over $\Q_p$.
 \end{lem}
 \begin{proof}
 	Evidently this set is invariant under $G(\QQ_p)$-conjugacy. Now let $\epsilon,\epsilon'\in G(\QQ_p)_{\semi}$ be stably conjugate and suppose that $\epsilon$ is $p^n$-admissible. We show that $\epsilon'$ is $p^n$-admissible. Let $M = M(\epsilon)$. Since $M$ is a Levi subgroup of $G_{\Q_p}$, the natural map $\coh^1(\QQ_p,M) \to \coh^1(\QQ_p,G)$ is an injection, cf.~\cite[\S 4.1]{hainesbasechange}. Also note that $G_{\epsilon}^0 = M_{\epsilon}^0$. Hence we have a natural bijection $\D (M_{\epsilon}^0, M ; \QQ_p) \isom \D (G_{\epsilon}^0, G ; \QQ_p)$. It follows that we have a natural surjection from the set of $M(\QQ_p)$-conjugacy classes in the stable conjugacy class of $\epsilon$ in $M$ onto the set of $G(\QQ_p)$-conjugacy classes in the stable conjugacy class of $\epsilon$ in $G_{\QQ_p}$.\footnote{Recall that these two sets are mapped onto by $\D(M_{\epsilon}^0, M ;\Q_p)$ and $\D(G_{\epsilon}^0, G ;\Q_p)$ respectively.} Conjugating $\epsilon'$ by an element of $G(\QQ_p)$ if necessary, we may assume that $\epsilon'$ lies in $M(\QQ_p)$, and that $\epsilon'$ is stably conjugate to $\epsilon$ inside $M$. Under these assumptions, we have $M = M(\epsilon')$. To check that $\epsilon'$ is $p^n$-admissible, it suffices to show that $w_M(\epsilon)=w_M(\epsilon')$. For this, note that since $\epsilon$ and $\epsilon'$ are stably conjugate in $M$, they are conjugate in $M(\LL)$ by Steinberg's theorem. The desired statement follows from the fact that $w_M$ is a group homomorphism from $M(\LL)$ to an abelian group. \end{proof}
 
 By Lemma \ref{lem:*epsilon is invariant}, we have a well-defined notion of \emph{$p^n$-admissibility} for stable conjugacy classes in $G(\QQ_p)_{\semi}$. The following result is a generalization of one direction in \cite[Satz 5.21]{langlands1987gerben}.
 
 \begin{thm} \label{thm:*epsilon}
 	Let $\epsilon \in G(\Q)_{\semi}$ represent an element of $\Sigma_{\RR \text{-} \el} (G)$ whose localization over $\QQ_p$ is $p^n$-admissible. Then the stable conjugacy class of $\epsilon$ lies in the image of $[\ggp_{\adm}(q^m)]$ under the map
 	$\F: [\ggp] \to \SigmaREllip(G)$. (Here $n = mr$ and $p^n = q^m$.)
 \end{thm}
 \begin{proof} By assumption there exists an $\RR$-elliptic maximal torus in $(G_{\epsilon}^0)_{\RR}$. By a theorem of Kneser \cite{KneserII}, there exists a $\Q_p$-elliptic maximal torus in $(G_{\epsilon}^0)_{\QQ_p}$. It then follows from \cite[Lem.~5.10]{langlands1987gerben} that there exists a maximal torus $T$ in $G_{\epsilon}^0$ such that $T_{\Q_v}$ is elliptic in $(G_{\epsilon}^0)_{\Q_v}$ for $v = \infty$ and $p$.
 	
 	Let $M = M(\epsilon)$ (which is defined over $\Q_p$) be as in Definition \ref{defn:M(epsilon)}, and let $\mu_M$ be as in Definition \ref{defn:*epsilon}. Then $T_{\Q_p}$ is a maximal torus in $M$. Let $\mu_T \in X_*(T)$ be a conjugate of $\mu_M$ under $M(\Qpbar).$ By Lemma \ref{lem:LR 5.12}, we find $u \in G(\Qbar)$ such that $\Int(u): T \to G$ is defined over $\Q$, and such that $\Int(u) \circ \mu_T = \mu_h$ for some $h \in X$. Note that $\epsilon$ is stably conjugate to $\Int(u) (\epsilon)$, because $u^{-1} \lix^{\tau} u \in T(\Qbar) \subset G_{\epsilon}^0(\Qbar)$ for all $\tau \in \Gamma$. Hence we may replace $\epsilon, T, M$ by $\Int(u) (\epsilon), \Int(u)(T), \Int(u)( M)$ respectively, and assume that $\mu _T = \mu_h$ for some $h \in X$ (which necessarily factors through $T_{\RR}$).
 	
 	We denote the inclusion $T \hookrightarrow G$ by $i$. Then we have $(T,i,h) \in \spd(G,X).$ Let $\phi = \phi(T,i,h)$ (Definition \ref{defn:special morphisms}). Then $\epsilon \in T(\Q) \subset I_{\phi}(\Q)$, and we have $(\phi,\epsilon) \in \wgp_{\sa}$ by Theorem \ref{thm:kisin special}. Moreover, by Remark \ref{rem:special morph is gg} we have $(\phi,\epsilon) \in \ggp_{\sa}$. It remains to show that $(\phi,\epsilon) \in \ggp_{\adm}(q^m)$.
 	
 	Since $M$ is an unramified reductive group over $\Q_p$, it contains an unramified \emph{and} elliptic maximal torus $T'$ (see \cite[p.~171]{langlands1987gerben} or \cite[\S 2.4]{debacker2006}). Let $\mu'\in X_*(T')$ be a conjugate of $\mu_T$ and $\mu_M$ under $M(\Qpbar)$. Let $b' = \mu'(p^{-1}) \in T'(\Q_p^{\ur})$. Then one immediately checks the following properties:
 	\begin{enumerate}
 		\item The element $b'$ is decent in $T'(\LL)$.
 		\item We have $w_{T'}(b')  = -\mu ' \in X_*(T')$. In particular, we have
 		\begin{align}
 			\label{eq:omega_M(b')}		w_{M}(b') = - [\mu_M]^M \in \pi_1(M). \end{align}
 	\end{enumerate}
 	
 	Since $\phi  = \phi(T,i,h) = i\circ \Psi_{T,\mu_T}$, the $G(\Q_p^{\ur})$-torsor $\mathcal{UR}(\phi(p) \circ \zeta_p)$ contains the $T(\Q_p^{\ur})$-torsor $\mathcal{UR}(\Psi_{T, \mu_T} (p) \circ \zeta_p)$. Let $(b_0,\epsilon_0)\in \cls_p(\phi,\epsilon)$ be the element associated with some $g\in \mathcal{UR}(\Psi_{T, \mu_T} (p) \circ \zeta_p)$ (see \S \ref{p-adic realization}). Then $\epsilon_0 = \epsilon$, and $b_0 \in T(\Q_p^{\ur})$. Moreover, by Lemma \ref{lem:b in torus}, the element $\kappa_T(b_0) \in X_*(T)_{\Gamma_p}$ equals the image of $-\mu_T$.
 	
 	Since $T_{\Q_p}$ and $T'$ are both elliptic in $M$, we know that $b_0 \in T(\LL)$ and $b' \in T'(\LL)$ are both basic in $M(\LL)$.
 	By what we have seen about $w_{T'}(b')$ and $\kappa_T(b_0)$, we have $\kappa_{M}(b') = \kappa_M(b_0) \in \pi_1(M)_{\Gamma_p}$. It follows that $b'$ and $b_0$ represent the same (basic) class in $\B(M)$ (see \cite[Prop.~5.6]{kottwitzisocrystal}). Since $b_0$ and $b'$ are decent, there exists $s \in M(\Q_p^{\ur})$ such that $b' = s b_0 \sigma(s)^{-1}$ (see \S\ref{isocdefns}). We then have $(b', s \epsilon s^{-1} ) \in \cls_p(\phi,\epsilon)$.
 	
 	We write $\epsilon'$ for $s \epsilon s^{-1}$, which is an element of $M(\Q_p^{\ur})$. It remains to show that $\epsilon'^{-1} (b'\rtimes \sigma)^n$ has a fixed point in $G(\LL)/ \cG(\breve \ZZ_p)$. For each $k \in \ZZ_{\geq 1}$, let $$U_k: = (\epsilon')^{-k} \cdot b' \cdot \sigma(b') \cdots  \sigma^{kn-1}(b') \in M(\Q_p^{\ur}), $$ so that
 	\begin{align}\label{eq:epsilon'b'U}
 		(\epsilon')^{-k} (b'\rtimes \sigma)^{kn} = U_k \rtimes \sigma^{kn}.
 	\end{align} Since $\epsilon'$ is $M(\Q_p^{\ur})$-conjugate to $\epsilon$, we have $w_M(\epsilon' ) = w_M(\epsilon)$. Hence $w_M(U_k) =0$ by (\ref{(*epsilon)}) and (\ref{eq:omega_M(b')}). By Lemma \ref{lem:ker of lambda}, the class of $U_1$ in $\B^{(n)}(M)$ comes from $\B^{(n)}(M_{\sconn})$. We claim that the class of $U_1$ in $\B^{(n)}(M)$ is basic with trivial Newton point. Given the claim, we use the fact that the only basic class in $\B^{(n)}(M_{\sconn})$ is the trivial class (\cite[Prop.~5.4]{kottwitzisocrystal}) to deduce that $U_1 = c \sigma^n (c^{-1})$ for some $c \in M(\LL)$. It follows that $\epsilon'^{-1} (b'\rtimes \sigma)^n = U_1 \rtimes \sigma^n$ has a fixed point in $G(\LL)/\cG(\breve \ZZ_p)$, namely $c \cG(\breve \Z_p)$. The proof of the theorem will then be finished.
 	
 	To prove the claim, it suffices to find $t\in \NN$ and $e\in M(\LL)$ such that
 	\begin{align}\label{eq:desiring t,e}
 		U_1 \cdot \sigma^n (U_1) \cdot \sigma^{2n}(U_1) \cdots \sigma^{(t-1)n}(U_1) = e \cdot \sigma ^{tn} (e^{-1}).
 	\end{align}
 	Since $\epsilon'$ commutes with $b'\rtimes \sigma$, it easily follows from (\ref{eq:epsilon'b'U}) that the left hand side of (\ref{eq:desiring t,e}) is equal to $U_t$. Fix an arbitrary reductive model $\cM$ of $M$ over $\ZZ_p$. By Greenberg's theorem \cite[Prop.~3]{Greenberg} (cf.~\cite[Lem.~1.4.9]{kottwitztwisted}), in order to find $t$ and $e$ such that $U_t = e \cdot \sigma^{tn}(e^{-1})$, it suffices to find $t$ such that $U_t \in \cM(\breve \ZZ_p)$.
 	
 	Now for each $k\in \NN$ we have
 	\begin{align}\label{eq:nu_l}
 		U_k = (\epsilon')^{-k} p^{\lambda_k} ,
 	\end{align} where
 	$$ \lambda_k  : = - \sum _{j=0} ^{kn -1} \sigma^j(\mu'). $$
 	Since $T'$ is unramified, for sufficiently divisible $k$ we have $\sigma^{kn}(\mu') = \mu'$. In this case $\lambda_k$ is defined over $\Q_p$, and in particular it is central in $M$ by the ellipticity of $T'$. Hence for sufficiently divisible $k$ we have
 	$$s^{-1} U_k s  = \epsilon^{-k} p^{\lambda_k} \in T(\Q_p).$$
 	Since $w_M(s^{-1} U_k s) = w_M(U_k) = 0$, and since $T_{\Q_p}$ is elliptic in $M$, we apply Lemma \ref{lem:compact elt} to conclude that  $s^{-1}U_k s$ lies in a compact subgroup of $T(\Q_p)$. In particular, for any given neighborhood $\mathcal N$ of $1$ in $T(\Q_p)$, all sufficiently divisible powers of $s^{-1} U_k s$ lie in $\mathcal N$. Now observe that when $\lambda_k$ is defined over $\Q_p$, we have $\lambda_{kl} = l \cdot \lambda_k$ for all $l \in \NN$, and therefore $$ s^{-1} U_{kl} s = (s^{-1} U_k s) ^l, \quad \forall l \in \NN.$$ We conclude that for $\mathcal N$ as above, we have $s^{-1} U_k s \in \mathcal N$ for all sufficiently divisible $k$. If we take $\mathcal N$ to be $(s^{-1}\cM(\breve \ZZ_p)s) \cap T(\Q_p)$, then we see that $U_k \in \cM(\breve \ZZ_p)$ for all sufficiently divisible $k$, as desired.
 \end{proof}
 For a quasi-split reductive group $M$ over $\Q_p$, recall that the \emph{degree $n$ norm} is a map from the set of $\sigma$-conjugacy classes in $M(\Q_{p^n})$ to the set of stable conjugacy classes in $M(\Q_p)$; see \cite[\S 5]{Kot82}.
 \begin{lem}\label{lem:norm and adm}
 	Let $\epsilon \in G(\Q_p)_{\semi}$, and let $M = M(\epsilon)$. Assume
 	that there exists a cocharacter $\mu_M$ of $M_{\Q_{p^n}}$ satisfying the first condition in Definition \ref{defn:*epsilon} and such that $\epsilon$ is a degree $n$ norm with respect to $M$ of some $\delta \in M(\Q_{p^n})$ satisfying $w_M(\delta) = - [\mu_M]^M$. Then $\epsilon$ is $p^n$-admissible.
 \end{lem}
 \begin{proof}
 	It suffices to check (\ref{(*epsilon)}), for the given $\mu_M$. First assume that $M_{\der}$ is simply connected. Then we have $\pi_1(M) = \pi_1(M^{\ab})$. Let $\bar \epsilon$ and $\bar \delta$ be the images of $\epsilon$ and $\delta$ in $M^{\ab}(\Q_p)$ and $M^{\ab}(\Q_{p^n})$ respectively. Then we have $\bar \epsilon = \bar \delta \sigma(\bar \delta) \cdots \sigma^{n-1}(\bar \delta)$. It follows that
 	\begin{align}\label{eq:omega norm}
 		w_M(\epsilon) = w_{M^{\ab}} (\bar \epsilon) = \sum_{i=0}^{n-1} \sigma^i (w_{M^{\ab}} (\bar \delta)) = \sum_{i=0}^{n-1} \sigma^i (w_M(\delta)),
 	\end{align} which gives the desired (\ref{(*epsilon)}).
 	
 	In the general case, we take a  $z$-extension $1 \to Z \to H \to M \to 1$ as in the proof of Lemma \ref{lem:ker of lambda}. Let $\tilde \delta \in H(\Q_{p^n})$ be a lift of $\delta$, and let $\tilde \epsilon \in H(\Q_p)$ be a degree $n$ norm of $\tilde \delta $. By the identity (\ref{eq:omega norm}) applied to $H, \tilde \epsilon, \tilde \delta$, we have
 	\begin{align}\label{eq:omega H}
 		w_H (\tilde \epsilon) = \sum_{i=0}^{n-1} \sigma ^i (w_H(\tilde \delta)) .\end{align} (This is an equality in $\pi_1(H)_{\Gamma_{p,0}}$.)
 	Now the image of $\tilde \epsilon$ in $M(\Q_p)$ is stably conjugate to $\epsilon$ over $\Q_p$, and therefore $M(\LL)$-conjugate to $\epsilon$ by Steinberg's theorem. Hence the image of $w_{H} (\tilde \epsilon)$ under $\pi_1(H)_{\Gamma_{p,0}} \to \pi_1(M)_{\Gamma_{p,0}} = \pi_1(M)$ equals $w_M(\epsilon)$. Obviously the image of $w_H(\tilde \delta)$ under $\pi_1(H)_{\Gamma_{p,0}}  \to \pi_1(M)_{\Gamma_{p,0}}  = \pi_1(M) $ equals $w_M(\delta)$. Hence from (\ref{eq:omega H}) we get
 	\begin{align*}
 		w_M(\epsilon) = \sum_{i=0}^{n-1} \sigma^i (w_M(\delta)),
 	\end{align*} which gives the desired (\ref{(*epsilon)}).
 \end{proof}
 In the next proposition, let the function $\phi_n : G(\Q_{p^n}) \to \set{0,1}$ be as in \S \ref{subsubsec:setting for point count}. Thus $\phi_n$ is the characteristic function of the double coset $\cG(\Z_{p^n}) p^{-\mu} \cG(\Z_{p^n})$ for arbitrary $\mu \in \dmu_X^{\cG}$, in the notation of \S\ref{subsubsec:SD}. 
 \begin{prop}\label{prop:sufficient condition for *epsilon} 	Let $\epsilon \in G(\Q)_{\semi}$ represent an element of $\Sigma_{\RR \text{-} \el} (G)$. Assume that $\epsilon\in G(\Q_p)$ is a degree $n$ norm of some $\delta \in G(\QQ_{p^n})$ whose $\sigma$-conjugacy class in $G(\Q_{p^n})$ intersects the support of $\phi_n$ non-trivially. Then the stable conjugacy class of $\epsilon$ is in the image of $[\ggp_{\adm}(q^m)]$ under $\F$. 
 \end{prop}
 \begin{proof} 
 	Let $\epsilon'$ be a $G(\Q_p)$-conjugate of $\epsilon$ such that $M(\epsilon')$ contains the generic fiber of a maximal $\Z_p$-split torus $\cA$ in $\cG$. By Theorem \ref{thm:*epsilon}, we only need to check that $\epsilon'$ is $p^n$-admissible.
 	
 	Write $M$ for $M(\epsilon')$. Since $M \supset \cA_{\Q_p}$, we have a reductive model $\cM$ of $M$ over $\Z_p$ such that the embedding $M \hookrightarrow G_{\Q_{p}}$ extends to $\cM \hookrightarrow \cG$. Let $\cS$ be a maximal $\ZZ_{p^n}$-split torus in $\cG_{\Z_{p^n}}$ such that $\cS$ contains $\cA_{\Z_{p^n}}$. In particular, $S: = \cS_{\Q_{p^n}}$ is a maximal $\Q_{p^n}$-split torus in $G_{\Q_{p^n}}$, and it is contained in $M_{\Q_{p^n}}$. There exists $\mu' \in X_*(S)$ that is $\cG(\Z_{p^n})$-conjugate to some $\mu \in \dmu_X^{\cG}$. We fix such a $\mu'$. Note that $\phi_n$ is the characteristic function of $\cG(\Z_{p^n}) p^{-\mu'} \cG(\Z_{p^n})$.  
 	
 	Denote by $s_G^M(\phi_n)$ the function in the unramified Hecke algebra of $M(\QQ_{p^n})$ with respect to $\cM(\Z_{p^n})$, obtained from $\phi_n$ via the partial Satake transform (a.k.a.~normalized constant term). Now
 	\cite[Lem.~4.2.1]{hainesbasechange} implies that there exists $\delta_M\in M(\QQ_{p^n})$ which is $\sigma$-conjugate to $\delta$ in $G(\QQ_{p^n})$ and whose degree $n$ norm with respect to $M$ is the stable conjugacy class of $\epsilon'$ in $M(\Q_p)$. Recall the parabolic descent formula for twisted orbital integrals
 	\[
 	TO^{G(\Q_{p^n})}_{\delta}(\phi_n)=|D^G_M(\epsilon')|^{-1/2}TO^{M(\Q_{p^n})}_{\delta_M}(s_G^M(\phi_n)),
 	\] 
 	where $D^G_M$ denotes the Weyl discriminant. For our purpose, it is enough to know that $D^G_M(\epsilon')\neq 0$, which is true since $M$ contains the identity component of the centralizer of $\epsilon'$ in $G$. The above formula is a special case of \cite[(4.4.4)]{hainesbasechange} with $\phi=\phi_n$, $E=\Q_{p^n}$, and $F=\Q_p$. (The sum on the right-hand side of \emph{loc.~cit.~}has a single term over $w=1$ by the Iwasawa decomposition, cf.~the definition of $_E W(P,J)$ in that paper.) Since the left hand side in the displayed formula is nonzero by initial hypothesis, it follows that the $\sigma$-conjugacy class of $\delta_M$ in $M(\Q_{p^n})$ intersects the support of $s_G^M(\phi_n)$ non-trivially.      
 	
 	On the other hand, $s_G^M(\phi_n)$ can be computed as follows. Let $\Omega_n,\Omega_{M,n}$ denote the $\Q_{p^n}$-rational Weyl groups of $S$ in $G,M$, respectively. For each $\mu_M\in X_*(S)$, write  $f^M_{-\mu_M}$ for the characteristic function of $\mathcal M(\Z_{p^n}) p^{-\mu_M} \mathcal M(\Z_{p^n})$; this depends only on the $\Omega_{M,n}$-orbit of $\mu_M$. The Satake transform $s$ for $G$ (resp.~$s_M$ for $M$) maps the unramified Hecke algebra of $G(\Q_{p^n})$ 
 	(resp.~$M(\Q_{p^n})$) onto $\CC[X_*(S)]^{\Omega_n}$ (resp.~$\CC[X_*(S)]^{\Omega_{M,n}}$). Moreover $s=s_M\circ s_G^M$ via the inclusion of $\CC[X_*(S)]^{\Omega_n}$ into $\CC[X_*(S)]^{\Omega_{M,n}}$. (This is a special case of the second commutative diagram in \cite[3.3.1]{hainesbasechange}.)     We claim that
 	\[
 	s_G^M(\phi_n)   = \sum_{w'\in \Omega_{M,n}\backslash \Omega_n} c_{w'} f^M_{-w'\mu'}\quad \mbox{for~some~nonzero}~c_{w'}\in \CC.
 	\] 
 	We may verify the equality after applying the isomorphism $s_M$. The computation of Satake transform on p.~297 of \cite{kottwitztwisted} (for the base field $\QQ_{p^n}$) tells us that
 	\[
 	s(\phi_n)=c \sum_{w\in \Omega_n} w(-\mu'),\qquad s_M(f^M_{-w'\mu'})=c(w') \sum_{w_M\in \Omega_{M,n}} w_M(-w'\mu')
 	\]
 	with explicit non-zero constants $c$ and $c(w')$ for each $w'\in \Omega_{M,n}\backslash \Omega_n$.
 	(The computation we need from \emph{loc.~cit.~}is summarized in line 8 therein. We apply his result to $-w_0\mu'$ in place of $\mu$ there, where $w_0$ is the longest element of $\Omega_n$, and similarly for the Levi $M$. To employ results from there, $-\mu'$ and $w'\mu'$ should be minuscule for $G$ and $M$, respectively; this holds true in our case.) Now it is clear that $s_G^M(\phi_n) = \sum_{w'} c_{w'} f^M_{-w'\mu'}$ with $c_{w'}=c/c(w')$. The claim is proved.
 	
 	Going back to $\delta_M$, since the $\sigma$-conjugacy class of $\delta_M$ in $M(\Q_{p^n})$ meets the support of $s_G^M(\phi_n)$, the same $\sigma$-conjugacy class meets the support of $f^M_{-w'\mu'}$ for some $w'\in \Omega_{M,n}\backslash \Omega_n$. Thus we may assume that $\delta_M\in \mathcal M(\ZZ_{p^n}) p^{-\mu_M} \mathcal M(\ZZ_{p^n})$ for $\mu_M:=w'\mu'$.
 	Then $\mu_M$ satisfies the first condition in Definition \ref{defn:*epsilon}, and $w_M(\delta_M) = - [\mu_M]^M$. Since the stable conjugacy class of $\epsilon'$ in $M(\Q_p)$ is the degree $n$ norm of $\delta_M$ with respect to $M$, we apply Lemma \ref{lem:norm and adm} to conclude that $\epsilon'$ is $p^n$-admissible.
 \end{proof}

 \section{Constructing Kottwitz parameters} \label{subsec:attach}
 In this section we attach a Kottwitz parameter to each gg semi-admissible pair $(\phi,\epsilon)\in \ggp_{\sa}$ and twisting element $\tau \in I^{\ad} _{\phi}(\A_f^p)$. After establishing basic functorial properties of this assignment, we show that the $p^n$-admissibility of $(\phi,\epsilon)$ is equivalent to that of the associated Kottwitz parameter. When the twisting element $\tau$ is trivial, we show that the Kottwitz parameter attached to $(\phi,\epsilon)$ has trivial invariant, and we establish a weak converse to show that the vanishing of the invariant is nearly sufficient for a $p^n$-admissible Kottwitz parameter to come from some $(\phi,\epsilon)$.
 \subsection{} \label{construction of coh kott trip}
 Let $(\phi,\epsilon)\in \ggp_{\sa}$, and let $\tau \in I^{\ad} _{\phi}(\A_f^p)$ be an arbitrary element. Our first goal is to assign to the triple $(\phi,\epsilon,\tau)$ a Kottwitz parameter (see Definition \ref{defn of coh Kott trip})\index[Oct2026-sorted-notation]{t 1 02 00 t phi epsilon tau@$\mathbf t (\phi,\epsilon,\tau)$} 	$$\mathbf t (\phi,\epsilon,\tau)\in \KP.$$ Let $\gamma_0 :=\epsilon,$ and $ I_0 : = G_{\gamma_0}^0$. The Kottwitz parameter $\mathbf t (\phi,\epsilon,\tau)$ shall be of the form $( \gamma_0 ,a ,[b])$. Note that by the first condition in Definition \ref{defn:gg}, $\gamma_0$ is indeed a semi-simple and $\RR$-elliptic element of $G(\Q)$, meeting the requirement in Definition \ref{defn of coh Kott trip}.
 
 We construct $a \in \D(I_0, G ; \A_f^p)$. Fix a lift $\tilde \tau \in I_{\phi}(\bar \A_f^p)$ of $\tau$, and let $\zzeta_{\phi} : \Gamma \to G(\bar \A_f^p)$ be the cocycle as in \S \ref{para:integral points in X_l}. Since $(\phi,\epsilon)$ is gg, $\zzeta_{\phi}$ takes values in $I_0(\bar \A_f^p)$. For each $\rho \in \Gamma$, let $t_{\rho} : = \tilde \tau^{-1} \cdot \lix^{\rho} {\tilde \tau } \in Z_{I_{\phi}} (\bar \A_f^p)$. (Here $\lix^{\rho} {\tilde \tau }$ is defined with respect to the $\Q$-structure of $I_{\phi}$.) Then $(t_{\rho})_{\rho}$ is a continuous cocycle $\Gamma \to Z_{I_{\phi}} (\bar \A_f^p)$. Consider the map
 \begin{align*}
 	A: \Gamma & \To G(\bar \A_f^p) \\
 	\rho & \longmapsto t_{\rho} \zzeta _{\phi}(\rho).
 \end{align*}
 Using the fact that the natural map $Z_{I_{\phi}} \to G$ is defined over $\Q$ and factors through $I_0$, we know that $A$ is a continuous cocycle $\Gamma \to I_0(\bar \A_f^p).$ We claim that $A$ has trivial image in $\coh^1(\Gamma, G(\bar \A_f^p))$. In fact, if we denote by $\tilde \tau_G$ the image of $\tilde \tau$ under the canonical embedding $I_{\phi}(\bar \A_f^p) \to G(\bar \A_f^p)$, then
 $$ t_{\rho} =  \tilde \tau_G^{-1} \cdot \zzeta_{\phi}(\rho) \cdot \lix^{\rho} {\tilde \tau_G} \cdot \zzeta_{\phi}(\rho)^{-1} . $$
 (Here $\lix^{\rho} {\tilde \tau_G} $ is defined with respect to the $\Q$-structure of $G$.) Hence
 \begin{align}\label{eq:A(rho)}
 	A(\rho) =  \tilde \tau_G^{-1} \cdot \zzeta_{\phi}(\rho) \cdot \lix^{\rho} {\tilde \tau_G}, \quad \forall \rho \in \Gamma.\end{align}
 Since $\phi$ is admissible, the cocycle $\zzeta_{\phi}$ has trivial image in $\coh^1(\Gamma, G(\bar \A_f^p))$. It follows from (\ref{eq:A(rho)}) that $A$ also has trivial image in $\coh^1(\Gamma, G(\bar \A_f^p))$. Now from (\ref{eq:A(rho)}) it is clear that the class of $A$ in $\coh^1(\A_f^p, I_0)$ is independent of the choice of $\tilde \tau$. We define the desired element $a \in \D(I_0, G;\A_f^p)$ to be the class of $A$.

 Next we construct $[b] \in \B(I_{0, \Q_p})$. By \S \ref{para:gg and gg}, we may view $\phi$ as a morphism $\phi_{I_0}: \Qf \to \G_{I_0}.$ We choose $g \in \mathcal {UR} (\phi_{I_0}(p) \circ \zeta_p)$ and let $b = b_{\Int(g^{-1}) \circ \phi_{I_0}(p) \circ \zeta_p} \in I_0(\Qpur)$; see Definition \ref{defn:b_theta}. By Lemma \ref{lem:unram criterion} (i), the class $[b]$ of $b$ in $\B(I_{0,\Q_p})$ is independent of choices. We now check condition \textbf{KP0} in Definition \ref{defn of coh Kott trip} for $[b]$. By the admissibility of $\phi$, we have $X_p(\phi) \neq \emptyset$. Comparing \S \ref{subsubsec:recall of defns in LR} and \S \ref{subsubsec:SD}, we have $X_p(\phi) \cong X_{-\mu}(b)$ for $\mu \in \dmu_X^{\cG}$, and so $X_{-\mu}(b) \neq \emptyset$. It immediately follows that $b$ and $-\mu$ have the same image in $\pi_1(G)_{\Gamma_p}$ (cf.~Corollary \ref{cor:kottCartan}). Thus \textbf{KP0} is satisfied by $[b]$. This finishes the construction of $\mathbf t (\phi,\epsilon,\tau)\in \KP.$

 \begin{prop}\label{funct of t changing tau}
 	Let $(\phi,\epsilon)\in \ggp_{\sa}$, and let $\tau \in I^{\ad} _{\phi}(\A_f^p)$. Let
 	\begin{align*}
 		(\gamma_0,a ,[b]) & : = \mathbf t (\phi,\epsilon,1) , \\ (\gamma_0,a',[b']) & : = \mathbf t (\phi,\epsilon,\tau).
 	\end{align*} Then $[b']  = [b] $. The difference $$a' - a \in \D(I_0, G ; \A_f^p) \cong  \E (I_0, G ;\A_f^p) $$ is equal to the image of $\tau$ under the composite map
 	\begin{align}\label{eq:tau to E(I_0,G)}
 		I^{\ad}_{\phi}(\A_f^p) \To \D(Z_{I_\phi}, I_\phi; \A_f^p) \To \D(I_0, G ;\A_f^p) .
 	\end{align}
 	In particular, the dependence of $\mathbf t (\phi ,\epsilon, \tau )$ on $\tau$ is only through the image of $\tau$ in $\D (I_0 , G ; \A_f^p)$.
 \end{prop}
 \begin{proof}
 	The equality $[b]=[b']$ is clear. In the notation of \S \ref{construction of coh kott trip}, $a$ is represented by the cocycle $\rho \mapsto \zzeta_{\phi}(\rho)$ and $a'$ is represented by the cocycle $ \rho \mapsto t_{\rho} \zzeta_\phi (\rho)$, where $\rho \mapsto t_{\rho} = \tilde \tau^{-1} \cdot \lix^{\rho} {\tilde \tau }$ is the cocycle $
 	\Gamma \to Z_{I_\phi} (\bar \A_f^p)$ determined by a choice of $\tilde \tau \in I_\phi(\bar \A_f^p)$ lifting $\tau$. The gg condition implies that the cocycle $\zeta_\phi^{p,\infty}$
 	takes values in 
 	$$ H = (G_{\phi^{\Delta}})_{\epsilon}^0 \subset I_0.$$ Since $Z_{I_\phi}$ is contained in the center of $H$, the abelian group structure on $\coh^1_{\ab}(\A_f^p, H)$ is compatible with  multiplication of cocycles: $$Z^1(\Gamma, Z_{I_\phi} (\bar \A_f^p)) \times Z^1 (\Gamma, H (\bar \A_f^p)) \To Z^1 (\Gamma, H (\bar \A_f^p)) , \quad  (c_\rho, d_\rho) \longmapsto c_\rho d_\rho.$$ Functoriality for $H \hookrightarrow I_0$ then gives the desired difference in $\coh^1_{\ab} (\A_f^p, I_0)$.  
 \end{proof}
 
 \begin{prop}\label{funct t conj}

 	Let $(\phi,\epsilon) , (\phi',\epsilon') \in \ggp_{\sa}$ such that $[\phi,\epsilon] = [\phi',\epsilon']$. Let $u \in G(\Qbar)$ be an element such that $\Int(u)(\phi,\epsilon) = (\phi',\epsilon')$. Let $\tau \in I_{\phi}^{\ad}(\A_f^p)$, and let $\tau'\in I_{\phi'}^{\ad}(\A_f^p)$ be the image of $\tau$ under the $\Q$-isomorphism $I_{\phi}^{\ad} \isom I_{\phi'}^{\ad}$ induced by the $\Q$-isomorphism $\Int(u): I_{\phi} \isom I_{\phi'}$. Then $u$ is an isomorphism
 	$\mathbf t(\phi,\epsilon,\tau) \isom \mathbf t(\phi',\epsilon',\tau')$
 	in the sense of Definition \ref{defn of equivalence between coh kott trip}.
 \end{prop}
 
 \begin{proof} We write $\mathbf t (\phi,\epsilon,\tau) = (\gamma_0, a,[b])$ and $\mathbf t (\phi', \epsilon', \tau') = (\gamma_0', a', [b'])$. By definition, $\gamma_0 = \epsilon$ and $\gamma_0' = \epsilon'$. We check that $u$ satisfies the three conditions in Definition \ref{defn of equivalence between coh kott trip}. 	Condition (i) follows from  Lemma \ref{lem:strong conj automatic}.
 	
 	To check condition (ii), we define $t_{\rho}$ as in \S \ref{construction of coh kott trip}, with respect to $(\phi,\epsilon,\tau)$. Then the counterpart $t_{\rho}'$ with respect to $(\phi',\epsilon',\tau')$ can be chosen to be $u t_{\rho} u^{-1}$. Now $a$ is represented by the cocycle $A: \rho \mapsto t_{\rho}\zzeta_{\phi}(\rho)$, and $a'$ is represented by the cocycle $A': \rho \mapsto t_{\rho}'\zzeta_{\phi'}(\rho) = u t_{\rho} u^{-1} \zzeta_{\phi'}(\rho)$. Note that $\zzeta_{\phi'}(\rho) = u \cdot \zzeta_{\phi}(\rho) \cdot \lix^{\rho}u^{-1}$, since $\phi' = \Int(u ) \circ \phi$.
 	Hence $$ A'(\rho) = u A(\rho) \lix^{\rho} u^{-1}, \quad \forall \rho \in \Gamma.$$ This proves condition (ii).
 	
 	To check condition (iii), we write $I_0$ for $G_{\epsilon}^0$ and write $I_0'$ for $G_{\epsilon'}^0$. We choose $d\in I_0(\overline {\breve \Q}_p)$ such that $u_0: = ud^{-1}$ lies in $G(\LL)$, as in \S \ref{para:comparing Kott invt}. Choose $g \in \mathcal {UR}(\phi_{I_0}(p) \circ \zeta_p) \subset I_0(\Qpbar)$, and $g' \in \mathcal {UR}(\phi'_{I_0'}(p) \circ \zeta_p) \subset I_0' (\Qpbar)$. We may assume that $b = b_{\Int(g^{-1}) \circ \phi_{I_0}(p) \circ \zeta_p}$ and $b' = b_{\Int(g'^{-1}) \circ \phi'_{I_0'}(p) \circ \zeta_p}$.
 	Then we have
 	\begin{align}\label{b and b' 0} &
 		g^{-1} u^{-1} g' \in G (\QQ_p^{\ur}),
 	\end{align}
 	since this element conjugates the unramified morphism $\Int(g'^{-1}) \circ \phi'(p) \circ \zeta_p: \G_p \to \G_G(p)$ to  the unramified morphism $\Int(g^{-1}) \circ \phi(p) \circ\zeta_p: \G_p \to \G_G(p)$. For the same reason we have
 	\begin{align}
 		\label{b and b' 1}
 		b =( g^{-1} u^{-1} g') b' \sigma ( g^{-1} u^{-1} g')^{-1}.
 	\end{align}The bijection $u_*: \B(I_{0, \Q_p}) \to \B(I_{0,\Q_p}')$ as in (\ref{eq:u_* at p}) sends $[b]$ to $[u_0 b \sigma(u_0)^{-1}]$, and by (\ref{b and b' 1}) the latter element   equals $[w b' \sigma (w)^{-1}]$, where $$w = u_0 g^{-1} u^{-1} g'\in G(\overline {\breve \QQ}_p) .$$ To finish the proof it suffices to show that $ w \in I'_0(\LL).$ Since $d^{-1} , g^{-1}\in I_0 (\overline {\breve \QQ}_p), $ we have $u_0 g^{-1} u^{-1} =\Int(u) (d^{-1} g^{-1}) \in I'_0 (\overline {\breve \QQ}_p)$. Also $ g' \in I'_0(\Qpbar)$, so $w \in I'_0 (\overline {\breve \QQ}_p)$. By the fact that $u_0 \in G(\LL)$ and by (\ref{b and b' 0}), we have $w \in G(\LL)$. Hence $  w \in I_0' (\LL) $, as desired.
 \end{proof}
 
 \subsection{}\label{para:setting for funct t twist}
 Let $(\phi,\epsilon) \in \ggp_{\sa}$. Let $I_0 := G_{\epsilon}^0.$ Recall that the group $\Sha_G^\infty(\Q, I_{\phi,\epsilon}^0)$ is defined in \S \ref{para:defn of Sha_G^infty for gg pair}. Using the canonical inner twisting between $I_{\phi, \epsilon}^0 $ and $(G_{\phi^{\Delta}})_{\epsilon}^0 \subset I_0$ as in \S \ref{para:gg and gg}, we have a natural homomorphism $\Sha_G^{\infty} (\Q, I_{\phi,\epsilon}^0) \to \E(I_0, G;\Q).$ Let $e = (e_\rho)_\rho\in Z^1(\QQ, I_{\phi,\epsilon}^0)$ be a cocycle representing a class in $\Sha_G^\infty(\Q, I_{\phi,\epsilon}^0)$. By Proposition \ref{semi-adm locus}, we obtain $(e\phi, \epsilon)\in \ggp_{\sa}$.
 
 \begin{prop}\label{funct t twist}   Let $
 	(\epsilon, a , [b])  : =\mathbf t(\phi,\epsilon,1)$, and $(\epsilon, a', [b'])  : =\mathbf t(e \phi,\epsilon, 1).$  Then the difference $a' - a \in \E(I_0, G, \A_f^p)$ is equal to the natural image of $e$. The elements $ [b'] , [b]\in \B(I_{0,\Q_p})$ have conjugate Newton cocharacters. Moreover, if $\nu_b$ is central in $I_0$, then $\kappa_{I_0} (b') - \kappa_{I_0 } (b) \in  \pi_1(I_0) _{\Gamma_p}$ is equal to the image of $e$ in $ \coh^1_{\ab} (\QQ_p, I_0)$, which is identified with $\pi_1 (I_0) _{\Gamma_p ,\tors}$  as in Proposition \ref{prop:apply TN}.
 \end{prop}
 
 \begin{proof} By construction, $a$ and $a'$ are represented by the cocycles $\zzeta_{\phi}$ and $\zzeta_{e\phi}$ respectively. Clearly \begin{align}\label{eq:zeta and zeta'}
 		\zzeta_{e\phi}(\rho) = e_\rho \zzeta_{\phi}(\rho), \quad  \forall  \rho \in \Gamma
 	\end{align} The statement about $a'-a$ follows from this and \cite[Lem.~3.15.1]{borovoi}. 
 	
 	We now prove the statement about $[b]$ and $[b']$. As in \S \ref{para:gg and gg}, we may view $\phi$ as a morphism
 	$ \Qf \to \G_{I_0 },$ which we denote by $\phi_{I_0}$. We may also view $e$ as a cocycle in $Z^1(\QQ, I_{\phi_{I_0}})$. Then $e\phi$ is induced by $ e\phi_{I_0} : \Qf \to \G_{I_0}$. By construction, $b$ (resp.~$b'$) is associated with the choice of an element $g$ in $\mathcal{UR}(\phi_{I_0}(p) \circ \zeta_p)$ (resp.~$g'$ in $\mathcal{UR}( (e\phi_{I_0})(p) \circ \zeta_p)$). Let $\theta:= \Int(g)^{-1} \circ \phi_{I_0}(p) \circ \zeta_p$ and let $\theta' : = \Int(g')^{-1} \circ  (e\phi_{I_0})(p) \circ \zeta_p$, and we view both as unramified morphisms $\G_p \to \G_{I_0 }(p)$. Thus $b = b_\theta$ and $b' = b_{\theta'}$. Now $\Int(g^{-1})$ induces a $\Q_p$-map $I_{\phi_{I_0}, \Q_p} \to I_{\theta}$. We let $e_\theta$ denote the image of $e$ under
 	$$ Z^1(\Q, I_{\phi_{I_0} }) \To Z^1(\Q_p, I_{\phi_{I_0}}) \xrightarrow{\Int(g^{-1})} Z^1(\Q_p, I_\theta).$$ By Proposition \ref{prop:twist at p} (ii), we have a canonical isomorphism $I_{\theta} \cong  J_b^{I_0}$. (See \S \ref{para:J_b} for the notation $J_b^{I_0}$.) We write $e_b$ for $e_{\theta}$ when we view it as an element of $Z^1(\Q_p, J_b^{I_0})$. Now $\theta'$ is in the conjugacy class of $e_\theta \theta$ (as morphisms $\G_p \to \G_{I_0} (p)$). By Proposition \ref{prop:twist at p} (iii), the $\sigma$-conjugacy class of $b'$ in $I_0(\Q_p^{\ur})$ is the twist of $b$ by $e_{b}$. In particular, $[b]$ and $[b']$ have conjugate Newton cocharacters by Proposition \ref{prop:twist p}.
 	
 	If $\nu_b$ is central, then we can apply Proposition \ref{prop:twist p} (to the group $I_{0,\Q_p}$), and conclude that $\kappa_{I_0}(b') - \kappa_{I_0}(b)$ is equal to the image of $e_{b}$ under
 	$$ \coh^1(\QQ_p, J_b^{I_0}) \isom  \coh^1_{\ab}(\Q_p, J_b^{I_0}) \isom   \coh^1_{\ab} (\QQ_p , I_0) \isom  \pi_1(I_0) _{\Gamma_p ,\tors}. $$
 	Here the second isomorphism is induced by the canonical inner twisting $\iota$ between the $\Q_p$-groups $J_b^{I_0}$ and $I_{0,\Q_p}$, as in case (ii) of \S \ref{para:J_b}. We are left to check that the above image of $e_b$ in $\pi_1(I_0)_{\Gamma_p, \tors}$ is equal to the image of $e$ in $\pi_1(I_0)_{\Gamma_p, \tors}$ as in the proposition. This boils down to the commutativity of the following diagram up to inner automorphisms:
 	$$  \xymatrix{(I_{\phi_{I_0}})_{\Qpbar} \ar[r] \ar[d]_{\sim}^{\Int(g^{-1})}& I_{0,\Qpbar} \\  I_{\theta , \Qpbar} \ar[r]^-{\sim} & (J_b^{I_0} )_{\Qpbar} \ar[u]^{\sim}_{\iota} } $$
 	Here the top arrow is given by the canonical inner twisting between $I_{\phi_{I_0}}$ and $(I_0)_{\phi^{\Delta}} \subset I_0$. The bottom arrow is the canonical isomorphism in Proposition \ref{prop:twist at p} (ii). It suffices to check that the composition $I_{\theta, \Qpbar} \isom ( J_b^{I_0})_{\Qpbar} \xrightarrow{\iota} I_{0,\Qpbar}$ is equal to the canonical embedding $I_{\theta, \Qpbar} \to I_{0,\Qpbar}$ attached to $\theta: \G_p \to \G_{I_0}(p)$. This is straightforward by the proof of Proposition \ref{prop:twist at p} (ii).
 \end{proof}

 \begin{prop}\label{vanishing of kott inv}
 	Let $(\phi,\epsilon) \in \ggp_{\sa}$, and let $\fkc = \mathbf t (\phi,\epsilon , 1) \in \KP$. Then the Kottwitz invariant $\alpha(\fkc)$ is trivial.
 \end{prop}
 \begin{proof}
 	By Proposition \ref{funct t conj} and Proposition \ref{funct of alpha conj}, we may replace $(\phi,\epsilon)$ by any other element $(\phi', \epsilon') \in \ggp_{\sa}$ such that $[\phi,\epsilon] = [\phi',\epsilon']$, in the course of the proof. Choose a maximal torus $T \subset I_\phi$ defined over $\Q$ such that $\epsilon \in T(\Q)$. (This is possible since $\epsilon$ is semi-simple.) By Theorem \ref{thm:precise speciality}, there exists $(T,i,h) \in \spd(G,X)$ and $g \in G(\Qbar)$ such that $\Int(g) \circ \phi = \phi(T,i,h)$, and such that $\Int(g) (\epsilon) = i(\epsilon) \in G(\Q).$ Using Lemma \ref{lem:suff gg pair} one checks that $(\Int(g ) \circ \phi, \Int(g)(\epsilon)  ) \in \ggp$. Hence up to replacing $(\phi, \epsilon)$ by $(\Int(g ) \circ \phi, \Int(g)(\epsilon)  )$, we may assume that $\phi = \phi(T,i,h)$, that $\epsilon \in T(\Q)$, and that the embedding $T_{\Qbar} \hookrightarrow I_{\phi, \Qbar} \hookrightarrow G_{\Qbar}$ coincides with $i$. 
 	
 	Write $\fkc  = (\gamma_0 = \epsilon, a, [b])$. By definition (see \S \ref{para:Psi_mu} and Definition \ref{defn:special morphisms}), $\phi$ factors as $\Qf \to \G_{\Res_{L/\QQ} \GG_m} \to \G_G$, where $L/\Q$ is a finite Galois extension contained in $\Qbar$ splitting $T$, and $\G_{\Res_{L/\QQ} \GG_m} \to \G_G$ is induced by a  $\Q$-homomorphism  $f: \Res_{L/\QQ} \GG_m \to G$. In particular,  $\zzeta_\phi$ is induced by a cocycle $\Gamma \to (\Res_{L/\Q} \GG_m)(\bar \A_f^p)$ and $f$. By Shapiro's lemma and Hilbert 90, the class $a$ is trivial.
 	
 	Recall that $b = b_{\theta},$ where $\theta = \Int(g^{-1}) \circ \phi(p) \circ \zeta_p$ for some $g \in \mathcal {UR}(  \phi(p) \circ \zeta_p) \cap I_0(\Qpbar)$. Since $\phi$ factors through $i:\G_T \to \G_G$, we may take $g$ to be inside $i(T)(\Qpbar)$. Then $\theta$ factors through $i: \G_T(p) \to \G_G(p)$, and $b \in i(T)(\Qpur)$. Write $b = i(b_T)$. By Lemma \ref{lem:b in torus}, the element $\kappa_T(b_T) \in X_*(T)_{\Gamma_p}$ is equal to the image of $-\mu_h \in X_*(T)$. Therefore, keeping the notation in \S \ref{defn of Kottwitz invariant}, we may choose $\tilde \beta_p(\fkc)$ to be the image of $-\mu_h\in X_*(T)$ in $\pi_1(I_0(\fkc))$ (with respect to $i: T \to I_0(\fkc)$). Also, we may choose $ \tilde \beta_{\infty}(\fkc)$ to be the image of $\mu_h\in X_*(T)$ in $\pi_1(I_0(\fkc))$. We have seen that $a$ is trivial, so we may choose $\tilde \beta_l(\fkc)$ to be zero for all $l \notin \set{p ,\infty }$. Then we have $\alpha(\fkc) =0$.
 \end{proof}
 
 \begin{cor}\label{value of kott inv}
 	Let $(\phi,\epsilon) \in \ggp_{\sa}$ and let $\tau \in I_{\phi}^{\ad}(\A_f^p).$ Let $\fkc  = \mathbf t (\phi,\epsilon ,\tau) \in \KP$. Then $\alpha(\fkc)$ is equal to the image of $\tau$ under
 	$$I^{\ad}_{\phi}(\A_f^p) \xrightarrow{(\ref{eq:tau to E(I_0,G)})} \E(I_0(\fkc), G ;\A_f^p) \hookrightarrow \E(I_0 (\fkc),G ; \A) \To \E( I_0 (\fkc) , G ; \adele/\QQ) .$$
 \end{cor}
 \begin{proof}
 	This follows from Proposition \ref{vanishing of kott inv}, Proposition \ref{funct of t changing tau}, and Proposition \ref{diff to diff}.
 \end{proof}
 
 As in \S\ref{admissible pairs}, let $n$ be a positive multiple of $r$. 
 \begin{prop}\label{admissible implies admissible} Let $(\phi,\epsilon) \in \ggp_{\sa}$ and $\tau \in I_{\phi}^{\ad} (\A_f^p)$. Then $(\phi,\epsilon)$ is $p^n$-admissible (Definition \ref{Defn of Adm}) if and only if $\mathbf t (\phi,\epsilon,\tau)$ is $p^n$-admissible (Definition \ref{pn-adm}).
 \end{prop}
 \begin{proof}
 	Write $\mathbf t(\phi,\epsilon,\tau) = (\epsilon, a, [b])$ and write $I_0$ for $ G_{\epsilon } ^0$. By construction, $[b] \in \B(I_{0,\Q_p})$ has a representative $b \in I_0(\Q_p^{\ur})$ such that $b = b_{\theta}$, where $\theta = \Int(g^{-1}) \circ \phi(p) \circ \zeta_p$ for some $g \in \mathcal{UR} (\phi(p)\circ \zeta_p)\cap I_0(\Qpbar)$. Since $g$ commutes with $\epsilon$, we see that $(b,\epsilon) \in \cls_p(\phi,\epsilon)$. By \cite[Lem.~1.4.9]{kottwitztwisted}, $\mathbf t(\phi,\epsilon,\tau) $ is $p^n$-admissible if and only if $\epsilon^{-1} (b \rtimes \sigma)^n$ has a fixed point in $G(\LL)/ \cG(\breve \Z_p)$. The latter condition is precisely the definition that $(\phi,\epsilon)$ is $p^n$-admissible.
 \end{proof}
 \begin{prop}\label{image of t} Let $\fkc = (\gamma_0, a, [b]) \in \KP_{\adm}(p^n)$ with $\alpha(\fkc)=0$. Assume that there exists an element of $\ggp_{\adm}(p^n)$ of the form $(\phi_0, \gamma_0)$. Then there exists an element of $\ggp_{\adm}(p^n)$ of the form $(\phi,\gamma_0)$ such that $\fkc = \mathbf t ( \phi,\gamma_0 ,  1). $
 \end{prop}
 
 \begin{proof}
 	Let $I_0 : = G_{\gamma_0}^0$. Write $\mathbf t( \phi_0, \gamma_0  ,1) = \fkc_0 = (\gamma_0, a_0,[b_0]) \in \KP.$ By Proposition \ref{admissible implies admissible}, we have $\fkc_0 \in \KP_{\adm}(p^n)$. Hence by Corollary \ref{cor:basic for adm Kott para} we know that $\nu_{b_0}$ is central in $I_0$ and that $\nu_b = \nu_{b_0}$. In particular, $\kappa_{I_0}(b)- \kappa_{I_0}(b_0)$ lies in $ \pi_1(I_0)_{\Gamma_p, \tors}$. Consider the element
 	$$ 		e_0 :=  
 	\big( a-a_0 ,   \kappa_{I_{0}}(b)  -\kappa_{I_{0}}(b_0)   \big)  \in  \coh^1_{\ab}(\A_f^p , I_0) \oplus \pi_1(I_0)_{\Gamma_p, \tors} \cong   \coh^1_{\ab}(\A_f, I_0).
 	$$ It directly follows from Definition \ref{defn of coh Kott trip} that  $e_0 \in \E(I_0, G ; \adele_f).$
 	
 	Note that the image of $e_0$ in $\E(I_0, G ; \adele/\QQ) $ is just $\alpha(\fkc) - \alpha (\fkc_0)$. We have $\alpha(\fkc) = 0$ by hypothesis, and $\alpha(\fkc_0)=0 $ by Proposition \ref{vanishing of kott inv}. By the exact sequence (\ref{eq:how to show zero invariant}), there exists a lift $e_1  \in \Sha_G^{\infty}(\Q, I_0)$ of $e_0$. Let $I:= I_{\phi_0,\gamma_0}^0$. Then we have a canonical inner twisting $I_{\Qbar} \isom I_{0,\Qbar}$ as in \S \ref{useful remark}, and hence a canonical isomorphism $\Sha_G^{\infty} (\Q,I) \cong \Sha_G^{\infty}(\Q, I_0)$. Let $e_2 \in \Sha_G^{\infty}(\QQ , I)$ be the element corresponding to 
 	$e_1 \in \Sha_G^{\infty} (\Q, I_0)$.
 	
 	By Proposition \ref{semi-adm locus}, we obtain an element  $(e_2\phi_0, \gamma_0 ) \in \ggp_{\sa}.$	We claim that $\mathbf t (e_2\phi_0, \gamma_0 , 1) = \fkc$. Write $\mathbf t (e_2\phi_0, \gamma_0 ,1) = (\gamma_0 , a', [b'])$. By Proposition \ref{funct t twist}, we have $a'- a_0 = a -a_0$, and so $a' = a$. As $\nu_{b_0}$ is central in $I_0$, we can apply Proposition \ref{funct t twist} to conclude that $\nu_{b'} = \nu_{b_0}$, and that $\kappa_{I_0}(b') - \kappa_{I_0}(b_0) = \kappa_{I_0}(b) - \kappa_{I_0}(b_0)$. Thus we have $[b] = [b']$ by the classification of $\B(I_{0,\Q_p})$ (see \S \ref{para:B(G)}). Having checked that $\mathbf t(e_2 \phi_0 , \gamma_0,1) = \fkc$, we deduce that $(e_2 \phi_0 , \gamma_0)\in \ggp_{\adm}(p^n)$ by Proposition \ref{admissible implies admissible}.
 \end{proof}
 
 \section{The effect of a controlled twist}
 \label{subsec:axiomatic setting}
 This section examines the Kottwitz parameters constructed in \S \ref{subsec:attach}, focusing on the specific case where the twisting elements are derived from the Langlands--Rapoport--$\tau$ Conjecture. 
 
 Specifically, Conjecture \ref{hypo about LR} predicts the existence of a 
 tori-rational element $\dertau \in \Gamma(\cH)_0$ such that $\mathsf{LR} (G,X, p, \cG, \dertau)$ holds. Fixing such a $\dertau$, we apply Lemma \ref{lem:about rect} to obtain a tori-rational $\dersigma \in \Gamma(\E^p)_0$ lifting $\dertau$. By the last assertion in Proposition \ref{funct of t changing tau}, we have a well-defined map\index[Oct2026-sorted-notation]{t 1 02 00 t subscript underline sigma phi epsilon@$\mathbf t_{\underline{\sigma}}(\phi,\epsilon)$}  $$\mathbf t_{\dersigma} : \ggp_{\sa} \To \KP, \quad (\phi,\epsilon) \longmapsto \mathbf{t} (\phi,\epsilon, \dersigma(\phi)).$$
 We proceed to establish properties of this assignment. 
 
 As in \S\ref{admissible pairs}, let $n$ be a positive multiple of $r$. 
 
 \begin{prop}\label{prop:funct t conj}
 	Let $(\phi,\epsilon)$ and $(\phi',\epsilon')$ be elements of $\ggp_{\sa}$ such that $[\phi,\epsilon] = [\phi',\epsilon']$. Let $u \in G(\Qbar)$ be an element such that $\Int(u)(\phi,\epsilon) = (\phi',\epsilon')$. Then $u$ is an isomorphism $ \mathbf t_{\dersigma}(\phi,\epsilon) \isom \mathbf t_{\dersigma}(\phi',\epsilon').$ In particular, the isomorphism class of $\mathbf{t}_{\dersigma} (\phi ,\epsilon)$ depends on $(\phi,\epsilon)$ only through $[\phi,\epsilon] \in [\ggp_{\sa}]$. \end{prop}
 \begin{proof}
 	Since $\dersigma \in \Gamma(\E^p)_0$, we may find a lift $\tau \in I_{\phi}^{\ad}(\A_f^p)$ of $\dersigma(\phi)$ and a lift $\tau ' \in I_{\phi'}^{\ad}(\A_f^p)$ of $\dersigma(\phi')$ such that $\tau$ maps to $\tau'$ under the isomorphism $I_{\phi}^{\ad} \isom I_{\phi'}^{\ad}$ induced by $\Int (u): I_{\phi} \isom I_{\phi'}$. The proposition then follows from Proposition \ref{funct t conj}.
 \end{proof}
 
 \begin{prop}\label{prop:funct t twist}
 	Keep the setting and notation of \S \ref{para:setting for funct t twist}, and assume in addition that $(\phi,\epsilon) \in \ggp_{\adm}(p^n)$. Let $
 	(\epsilon, a , [b])  : =\mathbf t_{\dersigma}(\phi,\epsilon)$, and $ (\epsilon, a', [b']) : =\mathbf t_{\dersigma}(e \phi,\epsilon)$. The following statements hold. 
 	\begin{enumerate}
 		\item 	The difference $a' - a \in \E(I_0, G; \A_f^p)$ is equal to the natural image of $e$. 
 		\item 
 		The elements $ [b'] , [b]\in \B(I_{0,\Q_p})$ are basic and have equal Newton cocharacter. The difference $\kappa_{I_0} (b') - \kappa_{I_0 } (b) \in  \pi_1(I_0) _{\Gamma_p}$ is equal to the image of $e$ in $\pi_1 (I_0) _{\Gamma_p ,\tors} $ as in Proposition \ref{funct t twist}.
 		\item Let $e' \in Z^1(\QQ, I_{\phi,\epsilon}^0)$ be another cocycle representing the same cohomology class as $e$. Then we have $\mathbf t_{\dersigma} (e \phi , \epsilon) = \mathbf t_{\dersigma} (e' \phi , \epsilon)$. (These two Kottwitz parameters are equal, not just isomorphic.) 
 	\end{enumerate}
 	
 \end{prop}
 \begin{proof}For (i), in view of Proposition \ref{funct of t changing tau} and the statement about $a'-a$ in Proposition \ref{funct t twist}, it suffices to check that $\dersigma (\phi)$ and $\dersigma(e\phi)$ have the same image in $\E(I_0, G ; \A_f^p)$. But this follows easily from the fact that $\dersigma \in \Gamma(\E^p)_0$.
 	
 	By Proposition \ref{funct of t changing tau}, the components $[b]$ and $[b']$ are unaffected by $\dersigma$. Hence to prove (ii) we may assume that $\dersigma$ is trivial. By Proposition \ref{admissible implies admissible} and Corollary \ref{cor:basic for adm Kott para}, we know that $[b]$ is basic in $\B(I_0)$. The remaining statements in (ii) follow from Proposition \ref{funct t twist}.
 	
 	Part (iii) follows from the previous two parts, and Kottwitz's classification of $\B(I_{0,\Q_p})$ (see \S \ref{para:B(G)}).
 \end{proof}
 For each $\fkc \in \KP$, recall that we have defined the Kottwitz invariant $\alpha(\fkc)$ in \S \ref{defn of Kottwitz invariant}. 
 \begin{prop}\label{prop:practical vanishing} 
 	For each $(\phi,\epsilon) \in \ggp_{\sa}$, we have $\alpha (\mathbf{t}_{\dersigma} (\phi,\epsilon)) = 0$. 
 \end{prop}
 \begin{proof} Choose a maximal torus $T \subset I_{\phi}$ defined over $\QQ$ such that $\epsilon \in T(\Q)$. (This is possible since $\epsilon$ is semi-simple.) By Corollary \ref{value of kott inv} and the exact sequence (\ref{eq:how to show zero invariant}), it suffices to show that the image of $\dersigma(\phi)$ in $\coh^1(\A_f^p, T)$ comes from $\Sha^{\infty, p}_G(\QQ, T)$. This follows from the fact that $\dersigma$ is tori-rational.
 \end{proof}
 Let $n$ be a positive multiple of $r$.  
 
 \begin{lem}\label{lem:converse stabl conj}
 	Let $(\phi_1,\epsilon_1) \in \ggp_{\adm}(p^n)$. Let $g\in G({\overline {\QQ}})$ be an element that stably conjugates $\epsilon_1$ to some $\epsilon_2 \in G(\QQ)$, i.e., $\Int(g)(\epsilon_1) = \epsilon_2$, and $g~ \leftidx{^\tau}g^{-1} \in G_{\epsilon_2}^0 (\overline \Q) $ for all $\tau \in \Gamma$. Define $\phi_2 : = \Int (g) \circ \phi_1$. Then $(\phi_2, \epsilon_2) \in \ggp_{\adm}(p^n)$.
 \end{lem}
 \begin{proof}We only need to show that $(\phi_2, \epsilon_2) \in \ggp$. Since $\epsilon_2 \in G(\Q)$ is stably conjugate to $\epsilon_1$, it is $\RR$-elliptic in $G$ as $\epsilon_1$ is (cf.~\S \ref{para:Frob KP}). Now write $\phi_i (q_\tau) = g_{i,\tau} \rtimes \tau$, for $i = 1,2$. It remains to show that $g_{2,\tau} \in I_{\phi_2, \epsilon_2}^0(\Qbar)$, for each $\tau \in \Gamma$. As in the proof of Lemma \ref{lem:strong conj automatic}, we have
 	$$ g ~\lix^{\tau} g^{-1} = [g g_{1,\tau} g^{-1}]^{-1} g_{2,\tau}.$$ The left hand side lies in $G_{\epsilon_2}^0(\Qbar)$ by hypothesis, and the term $g g_{1,\tau} g^{-1}$ lies in $G_{\epsilon_2}^0(\Qbar)$ since $g_{1,\tau} \in G_{\epsilon_1}^0(\Qbar)$. Hence we have $g_{2,\tau} \in G_{\epsilon_2}^0(\Qbar)$. By \S \ref{useful remark} we have $G_{\epsilon_1 ,\Qbar}^0 = I_{\phi_1, \epsilon_1, \Qbar}^0$. Hence $G_{\epsilon_2,\Qbar}^0 = \Int(g) (G_{\epsilon_1 ,\Qbar}^0 )= \Int(g)(I_{\phi_1, \epsilon_1, \Qbar}^0) = I_{\phi_2, \epsilon_2,\Qbar}^0$, and we have $g_{2,\tau} \in I_{\phi_2, \epsilon_2}^0(\Qbar)$ as desired. 
 \end{proof}
 
 \begin{prop}\label{prop:practical effective}
 	Let $\fkc = (\gamma_0 , a , [b])\in \KP_{\adm}(p^n)$. Let $(\gamma_0 , \gamma,\delta) \in \KP_{\clsc}(p^n)$ be the classical Kottwitz parameter of degree $n$ (up to equivalence) assigned to $\fkc$ as in \S \ref{from CK to K}. Assume that the $\sigma$-conjugacy class of $\delta$ in $G(\Q_{p^n})$ intersects the support of $\phi_n$ non-trivially. (Here $\phi_n$ is as in \S \ref{subsubsec:setting for point count}.) Assume that $\alpha (\fkc) = 0$. Then there exists an element of $\ggp_{\adm}(p^n)$ of the form $(\phi,\gamma_0)$ such that $\mathbf t_{\dersigma} (\phi,\gamma_0) = \fkc$.
 \end{prop}	
 \begin{proof} Let $I_0: = G_{\gamma_0} ^0$. Since $\gamma_0 \in G(\Q_p)$ is a degree $n$ norm of $\delta \in G(\Q_{p^n})$, by Proposition \ref{prop:sufficient condition for *epsilon} and Lemma \ref{lem:converse stabl conj} there exists an element of $\ggp_{\adm}(p^n)$ of the form $(\phi_0, \gamma_0)$. By Proposition \ref{image of t}, we find $(\phi_1,\gamma_0) \in \ggp_{\adm}(p^n)$ such that $\mathbf t(\phi_1, \gamma_0, 1) = \fkc$.
 	
 	Consider the Kottwitz parameter
 	\begin{align}\label{eq:KP1}
 		(\gamma_0 ,a_1, [b_1]) : = \mathbf{t}_{\dersigma} (\phi_1 ,\gamma_0).\end{align} By Proposition \ref{funct of t changing tau}, we know that $[b_1] = [b]$, and that $a_1 -a$ is equal to the image of $\dersigma (\phi_1)$ in $\E(I_0, G;\A_f^p)$. Fix a maximal torus $T \subset I_{\phi_1}$ such that $\gamma_0 \in T(\Q)$. Then by tori-rationality of $\dersigma$, the image of $\dersigma(\phi_1)$ in $\coh ^1 (\A_f^p, T) $ is equal to the image of some $\beta \in \Sha_{G}^{\infty,p}(\QQ, T)$.
 	
 	Under $T\hookrightarrow I_{\phi_1,\gamma_0} ^0$ the class $- \beta$ determines a class in $\Sha^{\infty,p}_{G}(\QQ, I_{\phi _1 ,\gamma_0} ^0)$. Fix a cocycle $ e \in Z^1(\QQ, I_{\phi_1,\gamma_0} ^0)$ representing the latter class. By Proposition \ref{semi-adm locus}, we obtain $(e\phi_1,\gamma_0) \in \ggp_{\sa}$. Write $\phi$ for $e\phi_1$, and let
 	\begin{align}\label{eq:KP2}
 		(\gamma_0, a', [b'] ) : = \mathbf{t}_{\dersigma} (\phi,\gamma_0).	\end{align}
 	Comparing (\ref{eq:KP1}) and (\ref{eq:KP2}), we see from Proposition \ref{prop:funct t twist} that $a'-a_1$ equals the image of $-\beta$
 	in $\E(I_0, G; \A_f^p)$, and that $[b'] = [b_1]$ (since $-\beta$ is trivial at $p$). Thus we have $a' = a$ and $[b'] = [b]$. Hence $\mathbf t_{\dersigma} (\phi ,\gamma_0) = \fkc$. By Proposition \ref{admissible implies admissible}, we have $(\phi,\gamma_0) \in \ggp_{\adm}(p^n)$. 
 \end{proof}
 
 \section{Proof of the point counting} In this section we prove Theorem \ref{thm:main thm in point counting}. 
 \subsection{}\label{para:correspondence} Throughout we fix a tori-rational element $\dertau \in \Gamma(\cH)_0$ such that the statement $\mathsf{LR} (G,X, p, \cG, \dertau)$ holds. Namely, 
 we have	a smooth integral model $\Shh_{K_p}$ of $\Sh_{K_p}$ over $\oo_{E,(\fkp)}$ which has well-behaved $\coh^*_c$, and we have a bijection
 \begin{align}\label{eq:LR twisted}
 	\Shh_{K_p} (\overline \FF_q)  {} \isom \coprod_{\phi }   S_{\dertau}(\phi) 
 \end{align}
 compatible with the actions of $G(\A_f^p)$ and the $q$-Frobenius $\Phi$.

 Our goal is to prove (\ref{pcf}) for all sufficiently large $m$.  First observe that in the proof we may arbitrarily replace $K^p$ by an open subgroup. In particular, we may and shall assume that $K^p$ is neat (see \cite[Def.~1.4.1.8, Rmk.~1.4.1.9]{lan2013arithmetic}), and that $\Shh_{K_pU^p}$ is defined for each open subgroup $U^p \subset K^p$. It follows from the neatness of $K^p$ that $K = K_p K^p$ is neat in the sense of Pink \cite[\S 0.6]{pink1989compactification}.
 In the sequel we write\index[Oct2026-sorted-notation]{z 0 00 00 z q subscript k@$Z(\mathbb Q)_K$} $Z(\QQ)_K$ for $Z_G(\QQ) \cap K$ and write\index[Oct2026-sorted-notation]{z 0 00 00 z subscript k@$Z_K$} $Z_K$ for $Z_G(\A_f) \cap K$, as in \S \ref{para:finite c_1}.

 By linearity, we may assume that $f^p = 1_{K^p g^{-1} K^p}$ for some $g \in G(\A_f^p)$, and that $dg^p$ assigns volume $1$ to $K^p$. 
 In the following, we write 
 $$ K^p_g : = K^p \cap g K^p g^{-1}. $$
 
 By Lemma \ref{lem:about rect} we fix a tori-rational $\dersigma \in \Gamma(\E^p)_0$ lifting $\dertau$. For each admissible morphism $\phi$, we fix a lift $\tau_\phi \in I_{\phi}^{\ad} (\A_f^p)$ of $\dersigma(\phi) \in  \E^p(\phi) = I_{\phi} (\A_f^p) \backslash I_{\phi} ^{\ad} (\A_f^p)$. We let $X^p (\phi)'$ be the $I_{\phi}(\A_f^p)$-set whose underlying set is $X^p(\phi)$, but the $ I_{\phi } (\A_f^p)$-action is given by the natural action pre-composed with conjugation by $\tau_\phi$. Thus $S_{\dertau} (\phi)$ is isomorphic to $$\varprojlim_{U^p} I_{\phi}(\Q) \backslash (X_p(\phi) \times  X^p(\phi)'/U^p).$$ 
 
 \begin{lem}\label{lem:check app licable}
 	For each admissible morphism $\phi$, the following statements hold.\begin{enumerate}
 		\item Let $\epsilon'$ be an element of $I_{\phi} (\A_f)$ which is conjugate under $I_{\phi}^{\ad}(\A_f)$ to some $\epsilon \in I_{\phi}(\QQ)$. If there exist $y ^p \in  X^p(\phi)/ K^p _g$ and $y_p\in X_p (\phi)$ such that $$ y ^p g \equiv \epsilon'  y^p g  \mod K^p, \quad \text{and } \quad \epsilon' y_p = y_p , $$ then $\epsilon'\in Z(\QQ) _K$.
 		\item We have $I_{\phi} (\QQ)_{\der} \cap Z(\QQ) _K =  \set{1}$.
 	\end{enumerate}
 \end{lem}
 \begin{proof}
 	\textbf{(i)} The proof follows the same idea as  \cite[Lem.~5.5]{milne92}.
 	We view $\epsilon$ as an element of $G({\overline {\QQ}})$ and let $\bar \epsilon$ be its image in $G^{\ad} ({\overline {\QQ}})$. By Lemma \ref{lem:I_phi at infty}, $\epsilon$ is semi-simple, and $\bar \epsilon$ lies in the $\RR$-points of a compact form of $G^{\ad}_{\RR}$. By the existence of $y^p$, the image of $\epsilon$ under $G(\Qbar) \to G(\bar \A_f^p)$ has a conjugate $u$ that lies in $K^p \subset G(\A_f^p) \subset G(\bar \A_f^p)$. By the existence of $y_p$, the image of $\epsilon$ under $G(\Qbar) \to G(\Qpbar)$ has a conjugate $v$ that lies in $\cG(\Z_p^{\ur})$. It follows that  $\bar \epsilon$ is torsion. But $\bar \epsilon$ can be  conjugated under $G^{\ad}(\bar \A_f^p)$ into the image of $K^p$ under $G(\A_f^p) \to G^{\ad}(\A_f^p)$, which is neat. Hence $\bar \epsilon =1$. We have thus shown that $\epsilon \in Z(\Qbar)$, and in particular $\epsilon = \epsilon'$.
 	
 	Since the natural embedding $Z \to I_{\phi}$ is defined over $\Q$ and since $\epsilon \in I_{\phi}(\Q)$, we have $\epsilon \in Z(\Q)$. Now using the existence of $y^p$ and $y_p$ it is easy to see that $\epsilon \in Z(\QQ)_K$.

 	\textbf{(ii)} Clearly $I_{\phi}(\QQ)_{\der} \cap Z(\QQ)_K$ is contained in the $\Q$-points of the center of $I_{\phi,\der}$, and is hence finite. Since $K$ is neat, $Z(\Q)_K$ is torsion free, and so $I_{\phi}(\QQ)_{\der} \cap Z(\QQ)_K = \set{1}$.
 \end{proof}
 \begin{lem} \label{lem:apply Deligne} We keep the assumptions on $f^p$ and $dg^p$ in \S \ref{para:correspondence}. When $m$ is sufficiently large (in a way depending on $K^p$ and $f^p$), we have
 	\begin{align}\label{eq:formula for T}
 		T(\Phi_{\fkp}^m , f^p dg^p)   = \sum_{\phi} \sum _{\epsilon} \# {\mathscr O(\phi,\epsilon, m , g)} \cdot \tr {\xi} (\epsilon),
 	\end{align} where
 	\begin{itemize}
 		\item
 		the first summation is over a set of representatives for the conjugacy classes of admissible morphisms $\phi$.
 		\item for each $\phi$ the second summation is over a subset of $I_{\phi} (\QQ)$ such that each conjugacy class in $I_{\phi} (\QQ)/ Z(\QQ)_K$ is represented exactly once.
 		\item the set $ \mathscr O(\phi,\epsilon, m ,g)$\index[Oct2026-sorted-notation]{o 0 05 00 o phi epsilon m g@$\mathscr O(\phi,\epsilon,m,g)$} is defined as the quotient of $$ X_p(\phi, \epsilon, q^m)\times \set{y^p \in X^p(\phi)'/ K^p _g \mid y^p  \equiv \epsilon y^p g \mod K^p}  $$ by the diagonal left action of $I_{\phi,\epsilon}(\QQ)$. Here $X_p(\phi, \epsilon, q^m)$ is defined in \S \ref{admissible pairs}.
 	\end{itemize}
 \end{lem}
 \begin{proof}  Since  (\ref{eq:adjunction for coh}) is an isomorphism by assumption, $T(\Phi_{\fkp}^m, f^pdg^p)$ is equal to
 	\begin{align}\label{eq:T1}
 		\sum_{i} (-1)^i \tr \bigg({\Phi}_{\mathfrak p}^m \times (f^p dg^p)  \mid \coh_c^i(\Shh_{K_p, \ol{\FF}_q},\xi) \bigg) .
 	\end{align}
 	We have two maps $\pi_g, \pi_1: \Shh_{K_p K^p_g} \to \Shh_{K}$, induced by the actions of $g\in G(\A_f^p)$ and $1\in G(\A_f^p)$ on $\Shh_{K_p}$, respectively. By our specific choices of $f^p$ and $dg^p$, the endomorphism $\Phi_{\fkp}^m \times (f^p dg^p)$ of $\coh_c^i(\Shh_{K, \ol{\FF}_q},\mL_{\xi,K^p})$ is  induced by the correspondence
 	\begin{align}\label{diag:correspondence}
 		\xymatrix{ &\Shh_{K_pK^p_g,\FF_q}   \ar[dl] _{\pi_g} \ar[dr]  ^{F^{mr} \circ \pi _1}\\ \Shh_{K,\FF_q} && \Shh_{K, \FF_q} }
 	\end{align}
 	and the cohomological correspondence 
 	\begin{align}\label{eq:Frob twisted}
 		\pi_g^* \mL_{\xi,K^p} \xrightarrow{\overrightarrow{g}^*} \mL_{\xi, K^p_g } \xrightarrow{(\overrightarrow{1}^* )^{-1}} \pi_1^* \mL_{\xi, K^p} \cong \pi_1^* (F^{mr})^* \mL_{\xi, K^p}.
 	\end{align}  (See \S \ref{para:sheaf and Hecke} for $\overrightarrow{g}^*$ and $\overrightarrow{1}^*$.)
 	Here $F$ is the absolute $p$-Frobenius endomorphism, and the last isomorphism in (\ref{eq:Frob twisted}) is induced by the canonical isomorphism between any \'etale sheaf and its pull-back under $F$.
 	
 	We compute (\ref{eq:T1}) by applying the Grothendieck--Lefschetz--Verdier trace formula in conjunction with Deligne's conjecture.\index{Deligne's conjecture} The latter---proven by Fujiwara \cite{fujiwara}  and Varshavsky \cite{varshavsky,varshavsky07} (see also \cite{pinkcalc})---permits the replacement of the local terms in the trace formula by the naive local terms, provided that $m$ is sufficiently large for fixed $K^p$ and $g$. 
 	
 	Let $\Fix$ be the set of $\ol \FF_q$-valued fixed points of the correspondence (\ref{diag:correspondence}).  Using the bijection (\ref{eq:LR twisted}), we obtain a description of $\Fix$ as follows. For each admissible morphism $\phi$,  by Lemma \ref{lem:check app licable} we know that the data
 	\begin{align}\label{situation}
 		\begin{cases}
 			Y = X_p(\phi) \times  X^p(\phi)'/ K^p_g, \\
 			X =  X_p(\phi) \times X^p(\phi)'/ K^p ,\\
 			I = I_{\phi} (\QQ), \\  C = Z(\QQ) _K ,\\
 			a: Y \to X , ~ ( y_p, y^p ) \mapsto (y_p, y^p g \mod K^p),\\
 			b :  Y \to X , ~ (y_p ,y^p) \mapsto (\Phi ^m y_p , y^p \mod K^p).
 		\end{cases}
 	\end{align} satisfy the hypotheses of \cite[Lem.~5.3]{milne92}. By \textit{loc.~cit.},  (\ref{eq:LR twisted}) induces a  bijection
 	\begin{align}\label{finite LR}
 		\Fix    \cong  \coprod _{\phi} \coprod _{ \epsilon} \mathscr O(\phi ,\epsilon,m , g) ,
 	\end{align}
 	where $\phi$ and $\epsilon$ run through the same ranges as in (\ref{eq:formula for T}).
 	
 	It remains to calculate the naive local term at each point in $\Fix$. We need to show that if $x \in \Fix$ corresponds to a point in $\mathscr O(\phi ,\epsilon, m , g)$ under (\ref{finite LR}), then the naive local term at $x$ is equal to $\tr \xi (\epsilon)$. Note that $x$ only determines the conjugacy class of $\epsilon$ in $I_\phi (\Q)/Z(\Q)_K$.  By our assumption that $\xi $ factors through $G^c$ and by Lemma \ref{lem:Z_s}, we know that $\tr \xi$ is invariant under $Z(\Q)_{K}$ since $K$ is neat. Hence $\tr \xi$ defines a class function on $G(\Qbar)$ that is translation-invariant under $Z(\Q)_K$, and our desideratum makes sense.
 	
 	Now let $\tilde x$ be a lift of $x$ in $\Shh_{K_p} (\overline \FF_q)$. Since $x \in \Fix$, there exists $k_{\asim{x}} \in K^p$ such that
 	$$\Phi^m \tilde x = \tilde x g k_{\asim{x}}.$$ Let $k_{\ell} \in G(\Q_{\ell})$ be the component of $k_{\asim{x}}$ at $\ell$, and let $g_{\ell} \in G(\Q_{\ell})$ be the component of $g \in G(\A_f^p)$ at $\ell$. By the same argument as \cite[p.~433]{kottwitz1992points}, the naive local term at $x$ is given by $$ \tr \xi(k_{\ell}^{-1} g_{\ell}^{-1}).$$ If $x$ corresponds to an element of $\mathscr  O(\phi, \epsilon, m ,g)$ under (\ref{finite LR}), then $k_{\ell}^{-1} g_{\ell}^{-1}$ is conjugate to $\epsilon$ in $G(\overline \QQ_{\ell})$. Thus the naive local term at $x$ is $\tr \xi (\epsilon)$, as desired.\end{proof}
 \begin{lem}\label{lem:card and orb int}
 	Let $(\phi,\epsilon ) \in \ggp_{\sa}$, and let $\fkc = \mathbf t_{\dersigma}(\phi,\epsilon)$. The set $\mathscr O(\phi,\epsilon,m, g)$ defined in Lemma \ref{lem:apply Deligne} is empty unless $(\phi,\epsilon)$ lies in $ \ggp_{\adm}(q^m)$. In the latter case, we have $\fkc \in \KP_{\adm}(p^n)$ by Proposition \ref{admissible implies admissible}, and we define  $$c_1(\fkc, K^p) O(\fkc, m,f^pdg^p)$$ as in \S \ref{subsubsec:setting for point count} and \S \ref{para:finite c_1}. We have
 	$$ \# \mathscr O(\phi,\epsilon,m, g)=  \iota_{I_{\phi}}(\epsilon)^{-1} c_1(\fkc, K^p)  O (\fkc, m , f^pdg^p).$$ See \S \ref{para:Frob KP} for $\iota_{I_{\phi}}(\epsilon)$. 
 \end{lem}
 \begin{proof}
 	If $\mathscr O(\phi,\epsilon,m,g) \neq \emptyset$, then $X_p(\phi, \epsilon,q^m) \neq \emptyset$, and it follows that $(\phi,\epsilon)$ is $q^m$-admissible. From now on assume that $(\phi,\epsilon)$ is $q^m$-admissible. 
 	We write $\fkc = (\gamma_0, a,[b])$, where $\gamma_0 = \epsilon$, and let  $(\gamma_0, \gamma ,\delta) \in \KP_{\clsc}(p^n)$ be the  classical Kottwitz parameter associated with $\fkc$ (which is well defined up to equivalence). 
 	The computation of $\#\mathscr O (\phi,\epsilon, m ,g)$ is essentially the same as the computation by Kottwitz in \cite[\S 1.4, \S 1.5]{kottwitztwisted}. We explain how to make the transition from our setting to the setting of \textit{loc.~cit.}
 	
 	We write $Y^p$ for the set $\{y^p \in X^p(\phi)'/ K^p _g \mid y^p  \equiv \epsilon y^p g \mod K^p \}$, and write $Y_p$ for the set $X_p(\phi,\epsilon, q^m)$. Thus $\mathscr O (\phi,\epsilon, m ,g) = I_{\phi,\epsilon}(\Q) \backslash (Y_p \times Y^p)$. Note that $\# \mathscr O (\phi,\epsilon, m,g)$ is equal to the cardinality of $I_{\phi,\epsilon}^0(\Q) \backslash (Y_p \times Y^p)$ multiplied by $\iota_{I_{\phi}}(\epsilon)^{-1}$. This is because, if an element $u\in I_{\phi,\epsilon}(\Q)$ has a fixed point in $Y_p \times Y^p$, then by Lemma \ref{lem:check app licable} (applied to $g=1$), we must have $u\in Z(\Q)_K \subset I_{\phi,\epsilon}^0 (\Q)$. Thus it remains to show that
 	\begin{align}\label{eq:desideratum in card}
 		\# \bigg (I_{\phi,\epsilon}^0(\Q) \backslash (Y_p \times Y^p)  \bigg) = c_1(\fkc)  O (\fkc, m , f^pdg^p).
 	\end{align}
 	
 	As explained in \S \ref{subsubsec:defn of S(phi)}, we identify $X^p(\phi)' = X^p(\phi)$ with a right $G(\A_f^p)$-coset inside $G(\bar \A_f^p)$. Fix $y^p_0 \in X^p(\phi)' \subset G(\bar \A_f^p)$. Inspecting definitions, we see that $\gamma \in G(\A_f^p)$ is conjugate to $(\Int (y^p_0)^{-1} \circ \Int(\tau_\phi)) (\epsilon) \in G(\A_f^p)$ inside $G(\A_f^p)$. We may assume that they are equal. Now if we use the ``base point'' $y^p_0$ to identify $X^p(\phi)'$ with $G(\A_f^p)$, then we have a bijection
 	\begin{align*}
 		Y^p \isom   W^p : =  \set{y^p \in G(\A_f^p)/ K^p _g \mid y^p  \equiv \gamma y^p g \mod K^p}.
 	\end{align*}
 	Under this bijection, the action of $I_{\phi,\epsilon}^0(\Q)$ on $Y^p$ corresponds to the action of $I_{\phi,\epsilon}^0(\Q)$ on $W^p$ given by the composition of the $\A_f^p$-isomorphism $\Int (y^p_0)^{-1} \circ \Int(\tau_\phi): (I_{\phi,\epsilon}^0)_{\A_f^p} \hookrightarrow  (G_{\gamma}^0)_{\A_f^p}$ (see \S \ref{useful remark}) followed by left multiplication of $G_{\gamma}^0(\A_f^p)$ on $W^p$.  
 	
 	Let $\mu \in \dmu_X^{\cG}$. By construction, we have $(b,\epsilon) \in \cls_p(\phi ,\epsilon)$, so $X_p(\phi, \epsilon, q^m)$ is identified with $X_{-\mu}(b, \epsilon, q^m)$; see \S \ref{admissible pairs}. More precisely, write $I_0$ for $G_{\epsilon}^0$, and suppose that $$b= b _{ \Int(g^{-1}) \circ  \phi_{I_0}(p) \circ \zeta_p}$$ for $g\in \mathcal{UR}(\phi(p) \circ \zeta_p) \cap I_0(\Qpbar)$ (following the notation in \S \ref{construction of coh kott trip}). 
 	Then $y_p \mapsto g^{-1} y_p$ induces a bijection $$Y_p \isom X_{-\mu}(b, \epsilon, q^m) ;$$ see \S \ref{subsubsec:recall of defns in LR} and \S \ref{admissible pairs}. Moreover, under this bijection the action of $I_{\phi,\epsilon}^0(\Q)$ on $Y_p$ corresponds to the action of $I_{\phi,\epsilon}^0(\Q)$ on $X_{-\mu}(b, \epsilon, q^m)$ given by the composition of the  $\Q_p$-isomorphism $\Int(g^{-1}): (I_{\phi,\epsilon}^0)_{ \Q_p } \isom  J_{b}^{I_0}$ followed by left multiplication of $J_b^{I_0}(\Q_p)$ on $X_{-\mu}(b, \epsilon, q^m)$. Here, to see that $\Int(g^{-1}): (I_{\phi,\epsilon}^0)_{ \Q_p } \isom  J_{b}^{I_0}$ is an isomorphism, use the identification $I_{\phi,\epsilon}^0\cong I_{\phi_{I_0}}$ as in \S \ref{para:gg and gg} and \S \ref{useful remark}, and use the fact that $b$ is basic in $I_0$ by Corollary \ref{cor:basic for adm Kott para}.

 	If we fix $c \in G(\LL)$ such that $\delta = c^{-1} b \sigma (c)$ and such that (\ref{eqn in KP1}) holds, then the map $x \mapsto c^{-1} x$ induces a  bijection
 	\begin{multline*}
 		X_{-\mu}(b, \epsilon, q^m)\isom \\ W_p : = \set{y_p \in G(\Q_{p^n})/ \cG(\Z_{p^n}) \mid y_p^{-1} \delta \sigma (y_p) \in \cG(\Z_{p^n}) p^{-\mu} \cG(\Z_{p^n}) }.
 	\end{multline*}
 	Under this bijection, the action of $J_b^{I_0}(\Q_p)$ on the left corresponds to the following action on $W_p$: We have an injective $\Q_p$-homomorphism $J_b ^{I_0} \to  R_{\delta \rtimes \theta}^0$ induced by $\Int (c^{-1})$. (See \S \ref{subsubsec:setting for point count} for $R_{\delta \rtimes \theta}^0$.) 	
 	This homomorphism is in fact an isomorphism, because both groups are connected, and their dimensions are equal to that of $I_0$. We thus identify $J_b^{I_0}(\Q_p)$ with $R^0_{\delta \rtimes \theta}(\Q_p)$, and let the latter group act on $W_p$ by left multiplication.
 	
 	In conclusion, we have bijections
 	\begin{align}\label{eq:three sets}
 		I_{\phi,\epsilon}^0(\Q) \backslash (Y_p \times Y^p) \isom  I_{\phi,\epsilon}^0(\Q) \backslash ( X_{-\mu}(b, \epsilon, q^m)  \times W^p)  \isom  I_{\phi,\epsilon}^0(\Q) \backslash (W_p \times W^p),
 	\end{align}	where $I_{\phi,\epsilon}^0(\Q)$ acts on $X_{-\mu}(b, \epsilon, q^m)$, $W^p$, and $W_p$ in the way described above. By abuse of notation, we still write $I_{\phi,\epsilon}^0(\Q)$ for the image of $I_{\phi,\epsilon}^0(\Q)$ inside $  J_b^{I_0}(\Q_p) \times G_{\gamma}^0(\A_f^p) \cong  R_{\delta \rtimes \theta}^0(\Q_p) \times G_{\gamma}^0(\A_f^p)$, under the embeddings described above. We also identify the last two product groups with $I(\fkc)(\A_f)$, canonically up to conjugation by $I(\fkc)^{\ad}(\A_f)$. We assume for a moment that $I_{\phi,\epsilon}^0(\Q) Z_K$ is closed and has finite co-volume inside $I(\fkc)(\A_f)$. Then the computation in \cite[\S 1.5]{kottwitztwisted}\footnote{In \textit{loc.~cit.}, it is assumed that $G_{\der}$ is simply connected, so that $G_{\gamma}$ and $R_{\delta \rtimes \theta}$ (which is denoted by $G_{\delta}^{\sigma}$) are connected. To transport the computation to the current situation, one simply replaces all appearances of $G_{\gamma}$ and $R_{\delta \rtimes \theta}$ by their identity components.} shows that the cardinality of the third set in (\ref{eq:three sets}) is given by $$ \vol (I_{\phi,\epsilon}^0(\Q) Z_K \backslash I(\fkc)(\A_f))  \cdot O (\fkc, m , f^pdg^p).$$
 	Here, $I_{\phi,\epsilon}^0(\Q) Z_K$ is equipped with the Haar measure giving volume $1$ to its open subgroup $Z_K = Z_G(\A_f) \cap K $.
 	
 	To complete the proof, we need to verify our assumption on $I_{\phi,\epsilon}^0 (\Q) Z_K$, and we need to identify $\vol (I_{\phi,\epsilon}^0(\Q) Z_K \backslash I(\fkc)(\A_f))$ with $c_1(\fkc,K^p)$. For both purposes, it suffices to prove the following claim: The image of $I_{\phi,\epsilon}^0(\Q)$ inside $I(\fkc)(\A_f)$ is $I(\fkc)^{\ad}(\A_f)$-conjugate to $I(\fkc)(\Q)$. Indeed, the claim together with Lemma \ref{lem:closed subgp} and the discussion in \S \ref{para:finite c_1} will imply our assumption on $I_{\phi,\epsilon}^0 (\Q) Z_K$. The claim together with the observation that the Haar measure on $I(\fkc)(\A_f)$ is invariant under conjugation by $I(\fkc)^{\ad}(\A_f)$ will imply the identity $\vol (I_{\phi,\epsilon}^0(\Q) Z_K \backslash I(\fkc)(\A_f)) = c_1(\fkc,K^p)$.
 	
 	To prove the claim, note that we have a canonical inner twisting between the $\Q$-groups $I_{\phi,\epsilon}^0$ and $I_{0}$, as in \S \ref{useful remark}. We have described an $\A_f^p$-isomorphism $I_{\phi,\epsilon}^0 \isom G_{\gamma}^0$, and a $\Q_p$-isomorphism $I_{\phi,\epsilon}^0 \isom J_b^{I_0}$. One checks that these isomorphisms are isomorphisms between inner forms of $I_0$ (in the sense of Definition \ref{defn:inner form}). Also, $I_{\phi,\epsilon}^0/Z_G$ is anisotropic over $\R$, and up to isomorphism there is at most one inner form of $I_0$ over $\R$ which is anisotropic modulo $Z_G$. The claim then follows from the unique characterization of $I(\fkc)$ as an inner form of $I_0$ over $\Q$ (up to isomorphism), in Proposition \ref{existence of I}. \end{proof}
 
 \subsection{} We have a natural action 	\begin{align*}
 	Z(\QQ)_K \times [\ggp_{\adm}(q^m)] & \To [\ggp_{\adm}(q^m)]  \\
 	(z, [\phi,\epsilon] ) & \longmapsto [\phi, z\epsilon].
 \end{align*}
 We fix a set of representatives for the $Z(\QQ)_K$-orbits in $[\ggp_{\adm}(q^m)]$  and by abuse of notation denote this subset of $[\ggp_{\adm}(q^m)]$ by  $$[\ggp_{\adm}(q^m)] / {Z(\QQ)_K}.$$ By the first statement in Lemma \ref{lem:card and orb int} and by Corollary \ref{cor:sa conj to gg}, we may replace the two summations in (\ref{eq:formula for T}) by the summation over $[\ggp_{\adm}(q^m)] / {Z(\QQ)_K}$. Applying Lemma \ref{lem:card and orb int} to the summands, we obtain \begin{multline}\label{eq:with sa pairs}
 	T(\Phi_{\fkp}^m, f^pdg^p)  \\ = \sum _{[\phi,\epsilon]\in [\ggp_{\adm}(q^m)] / Z(\QQ)_K} \iota_{I_{\phi}}(\epsilon)^{-1} c_1 (\mathbf t_{\dersigma}(\phi,\epsilon ) ) ~ O \big (\mathbf t_{\dersigma}(\phi,\epsilon ) , m , f^pdg^p \big ) \tr {\xi} (\epsilon),
 \end{multline} for all sufficiently large $m$. Here, for $(\phi,\epsilon) \in \ggp_{\adm}(q^m)$, we know that the isomorphism class of $\mathbf t_{\dersigma} (\phi,\epsilon)$ depends only on $[\phi,\epsilon]$, by Proposition \ref{prop:funct t conj}. Moreover, the value of $$c_1 (\mathbf t_{\dersigma}(\phi,\epsilon ) ) \cdot O \big (\mathbf t_{\dersigma}(\phi,\epsilon ) , m , f^pdg^p \big ) $$ depends on $\mathbf t_{\dersigma}(\phi,\epsilon )$ only via its isomorphism class, which can be checked using the definitions. 
 
 Let $\Sigma_{K^p}$ be as in \S \ref{para:Frob KP}, and for each $\gamma_0 \in \Sigma_{K^p}$ we write $\KP(\gamma_0, q^m)_{\alpha}$ for the set of $\fkc \in \KP(\gamma_0) \cap \KP_{\adm}(q^m)$ satisfying $\alpha(\fkc) = 0$. (See \S \ref{para:Frob KP} for the notation $\KP(\gamma_0)$.) By Proposition \ref{funct of alpha conj}, $\KP(\gamma_0, q^m)_{\alpha}$ is stable under isomorphisms between Kottwitz parameters. Using Proposition \ref{prop:practical vanishing}, we rewrite (\ref{eq:with sa pairs}) as
 \begin{align}\label{eq:sum over triples with unknown multiplicities}
 	T(\Phi_{\fkp}^m, f^pdg^p)  = \sum_{\gamma_0 \in \Sigma_{K^p}} \tr {\xi} (\gamma_0) \sum_{ \fkc \in \KP(\gamma_0, q^m)_{\alpha}/{\cong}} c_1 (\fkc)  O(\fkc, m , f^pdg^p)  \mathscr A (\fkc),
 \end{align}
 where \index[Oct2026-sorted-notation]{a 0 05 00 a c@$\mathscr A (\fkc)$}	\begin{align}\label{eq:defn of A(z)}
 	\mathscr A (\fkc): = \sum_{\substack{[\phi,\epsilon] \in [\ggp_{\adm}(q^m)] \\ \mathbf t_{\dersigma} (\phi,\epsilon) \cong \fkc  }} \iota_{I_\phi}(\epsilon)^{-1}.	\end{align}

 \begin{lem}\label{lem:compute A} Let $\gamma_0 \in \Sigma_{K^p}$ and let $\fkc \in \KP(\gamma_0, q^m)_{\alpha}$. Assume that $O(\fkc, m , f^p dg^p)$ is non-zero. Let $k$ be the number of elements of $\KP(\gamma_0, q^m)_{\alpha}$ that are isomorphic to $\fkc$. Then we have
 	$$ \mathscr A(\fkc) = k \cdot  \bar \iota_G(\gamma_0) ^{-1} \cdot  c_2(\gamma_0). $$ Here the notation is as in \S \ref{para:finite c_1} and \S \ref{para:Frob KP}. 	\end{lem}
 \begin{proof} First note that $k$ is finite, since it is at most $\abs{G_{\gamma_0}/G_{\gamma_0}^0}$. From the non-vanishing $O(\fkc, m , f^p dg^p)$ (or rather just the non-vanishing of the twisted orbital integral at $p$), it follows that $\fkc$ satisfies the assumptions of  Proposition \ref{prop:practical effective} (for $n = mr$). By that proposition, there exists an element of $ \ggp_{\adm}(q^{m})$ of the form $(\phi_0,\gamma_0)$ such that $\mathbf{t}_{\dersigma} (\phi_0,\gamma_0) = \fkc$. By Proposition \ref{one cell}, all $[\phi,\epsilon]$ appearing in the summation in (\ref{eq:defn of A(z)}) necessarily lie inside $\mathcal C_{[\phi_0,\gamma_0]}$. On the other hand, for each $[\phi,\epsilon] \in \mathcal C_{[\phi_0,\gamma_0]} \cap [\ggp_{\sa}]$, if $\mathbf t_{\dersigma} (\phi,\epsilon) \cong  \fkc$, then $[\phi,\epsilon]$ is automatically $q^m$-admissible, by Proposition \ref{admissible implies admissible}. Hence we can rewrite (\ref{eq:defn of A(z)}) as
 	\begin{align}\label{eq:formula for A'}
 		\mathscr A(\fkc) = \sum_{\substack{[\phi,\epsilon] \in \mathcal C_{[\phi_0, \gamma_0]}\cap [\ggp_{\sa}] \\ \mathbf t_{\dersigma} (\phi, \epsilon) \cong  \fkc }} \iota_{I_{\phi}} (\epsilon) ^{-1}.
 	\end{align}
 	
 	Let $\tilde \Sha_{G}^{\infty} (\Q, I_{\phi_0,\gamma_0} ^0)$ denote the inverse image of $\Sha _{G}^{\infty} (\Q, I_{\phi_0,\gamma_0} ^0)$ in $Z^1(\Q, I_{\phi_0,\gamma_0} ^0)$.
 	By Proposition \ref{semi-adm locus}, we have a surjection $$
 	\tilde \Sha _G^{\infty} (\QQ, I_{\phi_0 ,\gamma_0} ^{0}) \To \mathcal C_{[\phi_0,\gamma_0]}\cap [\ggp_{\sa}], \quad  e   \longmapsto [e \phi_0, \gamma_0].
 	$$ which factors through a surjection
 	\begin{align*}
 		\eta_{\phi_0,\gamma_0}: \Sha_G^{\infty} (\QQ, I_{\phi_0 ,\gamma_0} ^{0}) \To \mathcal C_{[\phi_0,\gamma_0]} \cap [\ggp_{\sa}].
 	\end{align*}
 	Define $\tilde D$ to be the subset of $\tilde \Sha_G^{\infty} (\QQ, I_{\phi_0,\gamma_0} ^0) $ consisting of elements $e$ such that $$\mathbf t_{\dersigma} (e\phi_0 , \gamma_0) \cong  \fkc.$$
 	Now let $\fkc_1, \fkc_2,\ldots, \fkc_{k}$ be all the distinct elements of $\KP(\gamma_0, q^m)_{\alpha}$ that are isomorphic to $\fkc$.  For each $e\in \tilde D$, we have $\mathbf t_{\dersigma} (e\phi_0 , \gamma_0)  = \fkc_i$ for a unique $i \in \set{1,\ldots, k}$. We thus obtain a partition $$\tilde D = \coprod_{i =1}^k  \tilde D_i. $$ By Proposition \ref{prop:funct t twist} (iii), for each $1\leq i \leq k$, the set $\tilde D_i$ is the inverse image of a subset $D_i$ of $\Sha_G^{\infty} (\Q, I_{\phi_0, \gamma_0}^0)$. 	We can thus rewrite (\ref{eq:formula for A'}) as
 	\begin{align}\label{eq:rewriting multiplicity}
 		\mathscr A(\fkc) = \sum _{i=1}^k \sum_{\beta \in D_i} \frac{1} { \iota _{I_{e \phi_0}} (\gamma_0) \cdot \# \eta_{\phi_0,\gamma_0}^{-1} (\eta_{\phi_0,\gamma_0} (\beta))  }
 	\end{align}
 	where $e$ is a cocycle representing $\beta$. (We will soon see that the fibers of $\eta_{\phi_0,\gamma_0}$ are indeed finite.)
 	
 	By \S \ref{useful remark}, the $\Qbar$-embedding $(I_{\phi_0,\gamma_0}^0)_{\Qbar} \to (G_{\gamma_0}^0)_{\Qbar}$ is an isomorphism and is an inner twisting between the $\Q$-groups $I_{\phi_0,\gamma_0}^0$ and $G_{\gamma_0}^0$. Using this observation and Proposition \ref{prop:funct t twist}, it is easy to see that for each $1\leq i\leq k$, the set $D_i$ is either empty, or a coset of $\Sha_G(\Q, I_{\phi_0,\gamma_0}^0)$ inside $\Sha_G^{\infty}(\Q, I_{\phi_0,\gamma_0}^0)$. We claim that it is never empty. Let $u_i \in G(\Qbar)$ be an isomorphism $\fkc \isom \fkc_i$. Then $u_i \in G_{\gamma_0}(\Qbar)$, and $u_i \lix^{\tau} u_i^{-1} \in G_{\gamma_0}^0(\Qbar), \forall \tau \in \Gamma$. By Lemma \ref{lem:converse stabl conj}, we have $(\Int(u_i)\circ \phi_0, \gamma_0) \in \ggp_{\adm}(q^m)$. By Proposition \ref{prop:funct t conj}, we have $\mathbf t_{\dersigma} ( \Int(u_i) \circ \phi_0 , \gamma_0) = \fkc_i$. It remains to check that $\Int(u_i) \circ \phi_0$ is of the form $e\phi_0$ for some $e \in Z^1(\Q, I_{\phi_0 ,\gamma_0}^0)$. For this, we fix a lift $q_\tau \in \Qf$ of each $\tau \in \Gamma$, and write $\phi_0(q_\tau) = g_{\tau} \rtimes \tau$. Define $e_{\tau} : = u_i g_{\tau} \lix^{\tau} u_i^{-1} g_{\tau}^{-1}$. Then $e = (e_{\tau})_{\tau}$ is a cocycle in $Z^1(\Q, I_{\phi_0})$, and $\Int(u_i) \circ \phi_0 = e \phi_0$. It remains to check that $e_{\tau} \in I_{\phi_0, \gamma_0}^0(\Qbar)$ for each $\tau$. Let $\pi_0 : = (G_{\gamma_0}/G_{\gamma_0}^0 ) (\Qbar)$, which is an abelian group as explained in \S \ref{para:Frob KP}. We have seen that $u_i$ and $\lix^{\tau} u_i$ map to the same element of $\pi_0$. Since $(\phi_0,\gamma_0)$ is gg, $g_{\tau}$ maps to the identity in $\pi_0$. Hence $e_{\tau}$ maps to the identity in $\pi_0$, i.e., $e_{\tau} \in G_{\gamma_0}^0 (\Qbar)=  I_{\phi_0, \gamma_0}^0(\Qbar)$, as desired.
 	
 	We have proved the claim. Hence for each $1\leq i \leq k$, we have $$\abs{D_i}  = \abs{\Sha_G(\Q, I_{\phi_0,\gamma_0}^0)} . $$ Since $I_{\phi_0,\gamma_0}^0$ is an inner form of $G_{\gamma_0}^0$ over $\Q$, this number is equal to $c_2(\gamma_0) = \abs{\Sha_G(\Q, G_{\gamma_0}^0)}$.
 	
 	To complete the proof of the lemma, it suffices to check that each summand in (\ref{eq:rewriting multiplicity}) is equal to $\bar \iota _G(\gamma_0)^{-1}$.
 	Recall from \S \ref{para:iota and eta} that the composition of $\eta_{\phi_0,\gamma_0}$ with the natural injection $\mathfrak v : [\ggp] \to \lprod{\wgp}$ factors as $$\coh^1(\Q, I_{\phi_0,\gamma_0}^0) \to \coh^1(\Q, I_{\phi_0,\gamma_0})\xrightarrow{\iota_{\phi_0,\gamma_0}} \lprod{\wgp}. $$ The map $\iota_{\phi_0,\gamma_0}$ is injective. Hence for each $\beta\in \coh^1(\Q, I_{\phi_0, \gamma_0}^0)$, the set $$\eta_{\phi_0,\gamma_0}^{-1}(\eta_{\phi_0,\gamma_0}(\beta))$$ is equal to  the fiber of the map $\coh^1(\Q, I_{\phi_0,\gamma_0}^0) \to \coh^1(\Q, I_{\phi_0,\gamma_0})$ containing $\beta$. This fiber can be identified with the kernel of $\coh^1(\Q, I_{e\phi_0,\gamma_0}^0) \to \coh^1(\Q, I_{e\phi_0,\gamma_0})$, where $e$ is a cocycle representing $\beta$. By \cite[Cor.~III.1.3]{Lab04}, the cardinality of this kernel is
 	$$  \bar \iota _{I_{e\phi_0}} (\gamma_0) \cdot (\iota _{I_{e\phi_0} } (\gamma_0))^{-1}.$$
 	From this, we see that if a summand in (\ref{eq:rewriting multiplicity}) is indexed by $\beta$, then this summand is equal to $\bar \iota _{I_{e\phi_0}} (\gamma_0) ^{-1}$ for any cocycle $e$ representing $\beta$. Since $\beta \in D_i$ for some $i$ and since $\fkc_i$ is $q^m$-admissible (as $\fkc$ is), we have $(e\phi_0, \gamma_0) \in \ggp_{\adm}(q^m)$ by Proposition \ref{admissible implies admissible}. Thus by \S \ref{useful remark}, we have canonical inner twistings $(I_{e\phi_0, \gamma_0})_{\Qbar} \isom (G_{\gamma_0})_{\Qbar}$ and $(I_{e\phi_0, \gamma_0}^0)_{\Qbar} \isom (G_{\gamma_0}^0)_{\Qbar}$.  In particular, we have an inner twisting between the commutative groups $I_{e\phi_0, \gamma_0}/ I_{e\phi_0, \gamma_0}^0$ and $G_{\gamma_0}/ G_{\gamma_0}^0$ (see \S\ref{para:Frob KP}), which must be an isomorphism over $\Q$. It follows that $\bar \iota _{I_{e\phi_0}} (\gamma_0) = \bar \iota_{G} (\gamma_0)$. The proof is complete.
 \end{proof}
 
 The proof of  (\ref{pcf}) is completed by combining (\ref{eq:sum over triples with unknown multiplicities}) and Lemma \ref{lem:compute A}.

 \cleardoublepage \part{Shimura varieties of abelian type}
 \label{part:3}
 
 \cleardoublepage
 \chapter{Results on crystalline representations} \label{sec:crystalline}
Throughout this chapter we fix a prime $p$. Our goal is to develop the \emph{local theory} necessary for the study of Shimura varieties of Hodge type in \S \ref{sec:Shimura Hodge type}. While much of this material recasts the framework established in \cite{kisin2010integral, kisin2012modp}, our approach shifts the focus toward the Tannakian formalism. After reviewing generalities on fiber functors in \S \ref{subsec:fiber functor}, we lay the foundations for a Tannakian treatment of $F$-isocrystals with integral structure---specifically, \S \ref{subsec:intFisoc} introduces the notion of integral $F$-isocrystals with $\mathcal{G}$-structure, where $\mathcal{G}$ is a flat, finite-type, affine group scheme over $\mathbb{Z}_p$. The central construction, developed in \S \ref{subsec:crystalline lattices}, associates such an object to a $\mathcal{G}(\mathbb{Z}_p)$-valued crystalline representation of a $p$-adic field. This construction employs the ``Breuil–Kisin module'' functor $\mathfrak{M}$ from \cite{kisin2010integral} and relies on the assumption that $\mathcal{G}$ is smooth with connected fibers and satisfies \emph{$\mathfrak{S}$-purity} (a property established by Ansch\"utz for parahoric group schemes). Finally, \S \ref{subsec:Refineintst} examines functoriality properties when the crystalline representation factors through a maximal torus.

 \section{Generalities on fiber functors}\label{subsec:fiber functor}
 In this section, we review the general theory of fiber functors, with a particular focus on those defined over the integral base $\mathbb{Z}_p$. We compare the abstract Tannakian perspective with the more concrete approach of defining fiber functors via tensors. The section concludes by introducing the torsors of rational and integral trivializations, respectively. 
 
 \subsection{}\label{para:setting for fiber functor} Let $\cG$ be a flat, finite-type, affine group scheme over $\ZZ_p,$ and
 $\Rep_{\Z_p} \cG$\index[Oct2026-sorted-notation]{rep 0 01 00 rep subscript z subscript p g@$\mathrm{Rep}_{\mathbb Z_p}\mathcal G$} the category of representations of $\cG$ on finite free $\Z_p$-modules. For any commutative ring $S$, we write $\Modfp_S$\index[Oct2026-sorted-notation]{mod 0 01 00 mod subscript s superscript fp@$\mathrm{Mod}_S^{\mathrm{fp}}$} for the category of finite projective $S$-modules. 
 
 Let $R$ be a faithfully flat $\Z_p$-algebra. By a\index{fiber functor} \emph{fiber functor over $R$}, we mean a faithful exact $\otimes$-functor
 $\omega:\Rep_{\Z_p} \cG \rightarrow \Modfp_R$.
 We denote by\index[Oct2026-sorted-notation]{0001 one 1 06 00 1 subscript r@$\mathbbl {1}_R$} $\idfunc_R :\Rep_{\Z_p} \cG \rightarrow \Modfp_R$ the functor which takes $L$ to $L\otimes R.$ Then $\idfunc_R$ is a fiber functor, called \emph{the standard fiber functor}.
 
 If $\omega: \Rep_{\Z_p} \cG \to \Modfp_R$ is a $\otimes$-functor and $S$ is an $R$-algebra, we write $\omega_S$ for the composition of $\otimes$-functors 
 $$  \Rep_{\Z_p} \cG  \overset{\omega} {\To} \Modfp_R \xrightarrow{\cdot \otimes_R S}  \Modfp_S,$$ 
 which is called the \emph{base change} or \emph{pull-back} of $\omega$ over $S$. 
 
 For two $\otimes$-functors $\omega_1, \omega_2: \Rep_{\Z_p} \cG \to \Modfp_R$, we write\index[Oct2026-sorted-notation]{isom 0 01 24 underline isom superscript otimes omega subscript 1 omega subscript 2@$\underline{\mathrm{Isom}}^{\otimes}(\omega_1,\omega_2)$} $\SIsom^{\otimes} (\omega_1, \omega_2)$ for the $R$-scheme of $\otimes$-isomorphisms from $\omega_1 $ to $\omega_2$; see \cite[\S 1.11]{deligneCT}. Thus for any $R$-algebra $S$, $\SIsom^{\otimes} (\omega_1, \omega_2)(S)$ is the set of $\otimes$-isomorphisms $\omega_{1,S} \isom \omega_{2,S}$. We write\index[Oct2026-sorted-notation]{aut 0 01 24 underline aut superscript otimes omega@$\underline{\mathrm{Aut}}^{\otimes}(\omega)$} $\underline{\Aut}^{\otimes}(\omega_1)$ for $\SIsom^{\otimes}(\omega_1, \omega_1)$.
 
 The reconstruction theorem, which is well known over a field, is valid in our current setting, guaranteeing that the natural morphism $\cG \to \underline{\Aut}^{\otimes}(\idfunc_{\Z_p})$ is an isomorphism of $\Z_p$-group schemes; see for instance \cite[X, Thm.~1.2, Rmk.~1.6]{milneAGS}, or \cite[Thm.~5.17]{wedhorn2004tannakian}.
 
 Let $\Rep_{\Q_p} \cG$\index[Oct2026-sorted-notation]{rep 0 01 00 rep subscript q subscript p g@$\mathrm{Rep}_{\mathbb Q_p} \mathcal G$} be the category of $\cG$-representations on finite-dimensional $\Q_p$-vector spaces (which depends only on $\mathcal G_{\QQ_p}$). By \cite[\S 1.5]{Serre68}, every representation in $\Rep_{\Q_p} \cG$ contains a $\cG$-stable $\Z_p$-lattice. Using this fact, for each fiber functor $\omega:\Rep_{\Z_p} \cG \rightarrow \Modfp_R$ we can define a $\otimes$-functor
 $$\omega[1/p]: \Rep_{\Q_p} \cG \To  \Modfp_{R[1/p]}$$
 as follows: If $V$ is in $\Rep_{\Q_p} \cG,$ then we write $V = \varinjlim_i V_i$ as a direct limit of $\cG$-stable $\Z_p$-lattices. We set $\omega[1/p](V) = \varinjlim_i \omega(V_i),$ which is naturally isomorphic to $\omega(V_i)\otimes_R R[1/p]$ for any $i.$ Note that $\omega[1/p]$ is again a fiber functor, i.e., a faithful exact $\otimes$-functor.
 
 Given the above construction, it is easy to see that the category of fiber functors $\Rep_{\Z_p} \cG \to \Modfp_R$ is equivalent to the category of fiber-wise faithful exact functors between fibered categories $\mathbf{Rep}~ \cG \to \mathbf{Bun}_{\Spec R}$ (which are fibered over the small Zariski site of $\Spec \Z_p$), as considered in \cite{broshi2013}. 
 Thus we shall freely import results from \textit{loc.~cit.~}for the fiber functors in our sense. Also cf.~the last remark in the introduction of \textit{loc.~cit}.

 Let
 $\omega: \Rep_{\Z_p} \cG \to \Modfp_R$ be a fiber functor. Then $\SIsom^{\otimes}(\idfunc_R, \omega)$ has a right $\cG_R$-action via the natural homomorphism  $\cG_R \to  \underline{\Aut}^{\otimes} (\idfunc_R)$ (which we have seen is an isomorphism). By \cite[Thm.~4.8]{broshi2013}, $\SIsom^{\otimes}(\idfunc_R, \omega)$ is in fact a right $\cG_R$-torsor over $R$ (for the fppf topology). This result could be viewed as a generalization of the reconstruction theorem
 recalled above. In the sequel, we denote the $\cG_R$-torsor $\SIsom^{\otimes}(\idfunc_R, \omega)$ by $P_{\omega}$.\index[Oct2026-sorted-notation]{p 0 00 00 p subscript omega@$P_\omega$}
 
 \subsection{}\label{para:short-hand}
 Let $\cG$ and $R$ be as in \S \ref{para:setting for fiber functor}. We introduce a short-hand notation. Let $\omega: \Rep_{\Z_p} \cG \to \Modfp_R$ be a fiber functor, and let $S$ be a $R[1/p]$-algebra. Then the base change $\omega_S: \Rep_{\Z_p} \cG \to \Modfp_S$ of $\omega$ over $S$ factors as $$ \Rep_{\Z_p} \cG \xrightarrow{\cdot \otimes \Q_p}  \Rep_{\Q_p} \cG \xrightarrow{\omega[1/p]}  \Modfp_{R[1/p]} \xrightarrow{\cdot \otimes_{R[1/p]} S}    \Modfp_S.  $$ We denote the composite functor $$(\cdot \otimes_{R[1/p]} S) \circ \omega[1/p]: \Rep_{\Q_p} \cG \To \Modfp_S$$ again by $\omega_S$.

 \begin{defn}\label{defn:M^otimes}
 	For any finite free module $M$ over any commutative ring $S$, we denote by $M^{\otimes}$\index[Oct2026-sorted-notation]{m 0 00 00 m superscript otimes@$M^{\otimes}$} the direct sum of all the $S$-modules
 	which can be formed from $M$ using the operations of taking duals, tensor products,
 	symmetric powers, and exterior powers. Elements of $M^{\otimes}$ are called\index{tensor} \emph{tensors over $M$}.
 \end{defn}
 \subsection{}\label{para:symm and ext}
 Let $\omega : \Rep_{\Z_p} \cG \to \Modfp_R$ be a fiber functor. Then $\omega$ is compatible with the operations considered in Definition \ref{defn:M^otimes}, by \cite[Thm.~4.8, Rmk.~4.2]{broshi2013}. (The compatibility with taking duals follows from rigidity; see \cite[Prop.~1.9]{deligne1982tannakian}.) In this case, for each $L$ in $\Rep_{\Z_p} \cG$ and each $\cG$-invariant element $s \in L^{\otimes}$, we may think of $s$ as a morphism $\mathbf 1 \to L^{\otimes}$ (where $\mathbf 1$\index[Oct2026-sorted-notation]{0001 one 1 02 00 1@$\mathbf 1$} is the unit object) and apply $\omega$ to it. We then get a morphism $\omega(s) : \mathbf 1 \to \omega(L^{\otimes}) \cong \omega(L)^{\otimes}$, or equivalently an element $\omega (s) \in \omega(L)^{\otimes}$. Here, it is understood that $\omega$ has been extended to infinite direct sums of objects, when we apply it to $L^{\otimes}$.

 \begin{defn} We say that $\cG$ is\index{definable by tensors} \emph{definable by tensors}, if there exists $L$ in $\Rep_{\Z_p} \cG$ and a subset $(s_{\alpha})_{\alpha\in \bm \alpha} \subset L^\otimes$ such that the $\Z_p$-homomorphism $\cG \to \GL(L)$ is a closed embedding whose image is the point-wise stabilizer of $(s_{\alpha})_{\alpha\in \bm \alpha}$. When this is the case we call \emph{$(L, (s_{\alpha})_{\alpha\in \bm \alpha})$} a \emph{defining datum for $\cG$}.
 \end{defn}
 
 \begin{rem}
 	If $\cG_{\Q_p}$ is reductive, then we know that $\cG$ is definable by tensors, by combining \cite[Lem.~3.2]{broshi2013} and \cite[Prop.~1.3.2]{kisin2010integral}.\end{rem}
 \subsection{}\label{para:s_alpha} Suppose $\cG$ is definable by tensors. Fix a defining datum $(L, (s_\alpha)_{\alpha\in \bm \alpha})$ for $\cG$. Let $R$ be a faithfully flat $\Z_p$-algebra. It will be useful to give a description of more explicit data giving rise to a fiber functor over $R$.
 For this, let $D$ be a finite free $R$-module equipped with a collection of tensors
 $(s_{\alpha,0})_{\alpha\in \bm \alpha} \subset D^\otimes$ indexed by the set $\bm \alpha$, and suppose that there exists an $R$-module isomorphism
 $L\otimes_{\Z_p} R \isom D$ taking each $s_{\alpha}$ to $s_{\alpha,0}.$ (Obviously such $D$ and $(s_{\alpha,0})_{\alpha \in \bm \alpha}$ always exist.)

 \begin{lem}\label{lem:explicitfiberfunctor} In the setting of \S \ref{para:s_alpha}, there exists a fiber functor $$\omega:\Rep_{\Z_p} \cG \To \Modfp_R$$
 	equipped with an isomorphism $\iota: \omega(L) \isom D$  such that $\iota$ maps $\omega(s_{\alpha}) \in \omega(L) ^{\otimes}$ to $s_{\alpha, 0}$ for each $\alpha \in \bm \alpha$. The pair $(\omega,\iota)$ is unique up to unique isomorphism in the following sense: Given two such pairs $(\omega,\iota)$ and  $(\omega',\iota')$, there is a unique $\otimes$-isomorphism $\omega \isom \omega'$
 	which takes $\iota$ to $\iota'.$
 \end{lem}
 \begin{proof} To show the existence of $\omega,$ we set $\omega(Q) = Q\otimes_{\Z_p} R$ for $Q$ in $\Rep_{\Z_p} \cG,$
 	and we take $\iota$ to be any isomorphism $L\otimes_{\Z_p} R \isom D$ taking $s_{\alpha}$ to $s_{\alpha,0}.$ (Such an isomorphism is assumed to exist in \S \ref{para:s_alpha}.)
 	
 	For the uniqueness, consider two such pairs $(\omega,\iota),$ $(\omega',\iota').$ We may assume that $(\omega , \iota)$ is as constructed above, so in particular $\omega = \idfunc_R$. Consider the $R$-scheme
 	$\SIsom_{(s_{\alpha})}(\omega(L), \omega'(L))$ whose points valued in any $R$-algebra $S$ classify $S$-linear isomorphisms $\omega(L) \otimes_R S\isom \omega'(L) \otimes_R S $ taking
 	$\omega(s_{\alpha})$ to $\omega'(s_{\alpha})$ for each $\alpha$. Using the existence of $\iota$ and $\iota'$, one sees that $\SIsom_{(s_{\alpha})}(\omega(L), \omega'(L))$ is a trivial $\cG_R$-torsor. (Here $\cG_R$ acts on 
 	$\omega(L) = L \otimes _{\Z_p} R$.)
 	There is a natural $\cG_R$-equivariant map
 	$$P_{\omega'} = \SIsom^{\otimes}(\idfunc_R,\omega') \To \SIsom_{(s_{\alpha})}(\omega(L), \omega'(L)),$$
 	which is necessarily an isomorphism since $P_{\omega'}$ is a $\cG_R$-torsor. It follows that $P_{\omega'}$ is a trivial $\cG_R$-torsor, and there exists an isomorphism of $\otimes$-functors $\omega \isom \omega'$ which is unique up to multiplication by elements of $\cG(R).$
 	In particular, there is a unique choice of such an isomorphism which takes $\iota$ to $\iota'.$
 \end{proof}
 \begin{rem}\label{rem:omega trivial}
 	The proof of Lemma \ref{lem:explicitfiberfunctor} also shows that if $(\omega,\iota)$ is a pair as in that lemma, then $\omega$ is $\otimes$-isomorphic to $\idfunc_R$.
 \end{rem}
 \begin{lem}\label{lem:fiberfunctorisogeny}
 	Let $(\omega, \iota)$ be a pair as in Lemma \ref{lem:explicitfiberfunctor}. Let $\omega'$ be a fiber functor $\Rep_{\Q_p} \cG \to \Modfp_{R[1/p]}$, and let $\iota'$ be an isomorphism $\omega'(L[1/p]) \isom D[1/p]$ taking $\omega'(s_\alpha) \in \omega'(L[1/p])^\otimes$ to $s_{\alpha,0} \in D[1/p]^\otimes$ for each $\alpha \in \bm \alpha$. Then there is a unique $\otimes$-isomorphism $\omega[1/p] \isom \omega'$ which takes $\iota$ to $\iota'$.
 \end{lem}
 \begin{proof}
 	By Remark \ref{rem:omega trivial}, we may assume that $\omega = \idfunc_R$ without loss of generality. The proof of the lemma is then completely analogous to the proof of the uniqueness in Lemma \ref{lem:explicitfiberfunctor}.
 \end{proof}
 
 \begin{lem}\label{lem:fiberfunctortorsor} Let $\cG$ be a smooth affine group scheme over $\Z_p$ with connected fibers. (We do not need to assume that $\cG$ is definable by tensors.) Let $R$ be the ring of integers in either a finite unramified extension of $\Q_p$ or in $\Q_p^{\ur}.$  Let $\omega: \Rep_{\Z_p}\cG \to \Modfp_R$ be a fiber functor. Then $P_{\omega}$ is a trivial $\cG_R$-torsor. In particular, the set $P_{\omega}(R)$ is a $\cG(R)$-torsor.
 \end{lem}
 \begin{proof}
 	Our assumptions, together with Lang's theorem, imply that any $\cG_R$-torsor over $R$, such as $P_{\omega}$, is necessarily trivial.
 \end{proof}
 \subsection{}\label{para:defnY}
 We continue to consider $\cG$ and $R$ as in \S \ref{para:setting for fiber functor}. For each fiber functor $\omega:\Rep_{\Z_p} \cG \rightarrow \Modfp_R$ we set\index[Oct2026-sorted-notation]{y 0 00 00 y omega@$Y(\omega)$}
 $$Y(\omega) : =  P_{\omega}(R[1/p]) = \{ \text{$\otimes$-isomorphisms } \idfunc_{R[1/p]} \isom \omega_{R[1/p]} \}.$$
 This is either empty or a right $\cG(R[1/p])$-torsor (i.e., for $\eta \in Y(\omega)$ and $g \in \cG(R[1/p])$, we define $\eta \cdot g : = \eta \circ g$). When $P_\omega(R)$ is non-empty and thus a (right) $\cG(R)$-torsor, there is a canonical isomorphism 
 $$ Y(\omega) \cong P_{\omega}(R)\times^{\cG(R)}\cG(R[1/p]).$$ In this case we write\index[Oct2026-sorted-notation]{y 0 00 00 y omega superscript circ@$Y(\omega)^{\circ}$} $Y(\omega)^{\circ}$ for the subset of $Y(\omega)$ corresponding to $P_{\omega}(R)$. Its elements are called \emph{integral points} of $Y(\omega)$.

 \section{Integral \texorpdfstring{$F$}{F}-isocrystals with \texorpdfstring{$\cG$}{G}-structure}\label{subsec:intFisoc}
 In this section, we employ the Tannakian perspective to formally define \emph{integral $F$-isocrystals with $\mathcal G$-structure}.  
 \subsection{} \label{para:K_0 and sigma}
 Let $K_0$ be either a finite unramified extension of $\Q_p,$ or $\Q_p^{\ur}$. Let $\sigma \in \Aut(K_0)$ be the arithmetic $p$-Frobenius. We will apply the considerations in \S \ref{subsec:fiber functor} to $F$-isocrystals over $K_0$ with additional structures.
 
 Recall that an\index{F-isocrystal@$F$-isocrystal} \emph{$F$-isocrystal over $K_0$} is a
 finite-dimensional $K_0$-vector space $V,$ equipped with a $\sigma$-linear map $\varphi_V,$ called the \emph{Frobenius}, 	whose linearization 
 $1\otimes {\varphi_V}: \sigma^*V \isom V$ is an isomorphism. 
 We denote by $\Isoc_{K_0}$\index[Oct2026-sorted-notation]{isoc 0 01 00 isoc subscript k subscript 0@$\mathrm{Isoc}_{K_0}$} the category of $F$-isocrystals over $K_0.$
 
 Let\index[Oct2026-sorted-notation]{isoc 0 01 00 isoc subscript k subscript 0 superscript circ@$\mathrm{Isoc}_{K_0}^{\circ}$} $\Isoc^\circ_{K_0}$ be the category of pairs $(M, \varphi_M)$, where $M$ is a finite free $\O_{K_0}$-module, and $(M[1/p], \varphi_M)$ is an $F$-isocrystal. We shall often write $\varphi$ for $\varphi_M$ when no confusion can arise. Morphisms in this category are morphisms of $\oo_{K_0}$-modules that respect the Frobenii (after inverting $p$). We call an object of $\Isoc^\circ_{K_0}$ an {\em integral $F$-isocrystal over $K_0$}. Note that $\Isoc^\circ_{K_0}$ contains the category of $F$-crystals, but has the advantage of containing duals of objects.
 We write $$v: \Isoc_{K_0}^{\circ} \To \Modfp_{\O_{K_0}}$$ for the forgetful functor, taking each $(M,\varphi_M)$ to $M$. (Here $v$ stands for \textit{vergessen}, as in \cite[\S 1]{RZ96}.)
 
 Let $\cG$ be a flat, finite-type, affine group scheme over $\Z_p$. By an\index{F-isocrystal@$F$-isocrystal!integral with $\cG$-structure} \emph{integral $F$-isocrystal with $\cG$-structure over $K_0$}, we mean a faithful, exact, $\otimes$-functor
 $$\tomega: \Rep_{\Z_p} \cG \To \Isoc^\circ_{K_0}.$$ Equivalently, $\tomega$ is a $\otimes$-functor such that the composition $v\circ \tomega: \Rep_{\Z_p} \cG \rightarrow \Modfp_{\O_{K_0}}$ is a fiber functor. Note that we do not require $v\circ \tomega$ to be equal to $\idfunc_{\mathcal O_{K_0}}$.
 
 Similarly, by an\index{F-isocrystal@$F$-isocrystal!with $\cG$-structure} \emph{$F$-isocrystal with $\cG$-structure over $K_0$}, we mean a faithful, exact, $\otimes$-functor $$\tomega: \Rep_{\Q_p} \cG \To \Isoc_{K_0} .$$
 
 If $\tomega$ is an integral $F$-isocrystal with $\cG$-structure, then a construction as in \S \ref{para:setting for fiber functor} yields a natural $F$-isocrystal with $\cG$-structure $\tomega[1/p]$.
 
 We denote the categories of integral $F$-isocrystals (resp.~$F$-isocrystals) with $\cG$-structure over $K_0$
 by $\GIsoc^\circ_{K_0}$\index[Oct2026-sorted-notation]{g 0 03 00 g minus isoc subscript k subscript 0 superscript circ@$\mathcal G\textbf{-}\mathrm{Isoc}_{K_0}^{\circ}$} (resp.~$\GIsoc_{K_0}$).\index[Oct2026-sorted-notation]{g 0 03 00 g minus isoc subscript k subscript 0@$\mathcal G\textbf{-}\mathrm{Isoc}_{K_0}$} In these categories, morphisms are, by definition, $\otimes$-isomorphisms between $\otimes$-functors $\Rep_{\Z_p} \cG \to \Isoc_{K_0}^{\circ}$ (resp.~between $\otimes$-functors $\Rep_{\Q_p} \cG \to \Isoc_{K_0}$).

 \subsection{}\label{para:setting for explicit isoc} Suppose $\cG$ is definable by tensors, and fix a defining datum $(L, (s_\alpha)_{\alpha\in \bm \alpha})$ for $\cG$ (see \S \ref{para:s_alpha}).	Let $(D,\varphi_D)$ be an object in $\Isoc_{K_0}^{\circ}$, equipped with a collection of $\varphi_D$-invariant tensors $(s_{\alpha,0})_{\alpha \in \bm \alpha} \subset D^\otimes$, such that there
 exists an $\O_{K_0}$-module isomorphism $L\otimes \O_{K_0} \isom D$ taking each $s_{\alpha}$ to $s_{\alpha,0}$.  
 
 \begin{lem}\label{lem:Isocexplicitfiberfunctor} With the notation of \S \ref{para:setting for explicit isoc}, there exists an integral 
 	$F$-isocrystal with $\cG$-structure $\tomega,$ 
 	equipped with an isomorphism $\iota: \tomega(L) \isom (D,\varphi_D)$ in $\Isoc_{K_0}^\circ$ such that $\iota$ maps $\tomega(s_{\alpha}) \in \tomega(L) ^{\otimes}$ to $s_{\alpha, 0}$ for each $\alpha \in \bm \alpha$. The pair $(\tomega,\iota)$ is unique up to unique isomorphism. 
 \end{lem}
 \begin{proof}
 	Fix an $\oo_{K_0}$-module isomorphism $\iota: L\otimes \O_{K_0} \isom D$ taking each $s_{\alpha}$ to $s_{\alpha,0}.$
 	The composite map
 	$$  L\otimes K_0  \underset{\sim}{\xrightarrow{\sigma^*\iota} } \sigma^* D[1/p] \overset{\varphi_D} {\To} D[1/p] \underset{\sim} { \xrightarrow{\iota^{-1}}}  L\otimes K_0  $$
 	fixes $(s_{\alpha}),$ and is hence induced by some $\delta \in \cG(K_0).$ For $Q$ in $\Rep_{\Z_p} \cG$
 	we set $\tomega(Q) = Q \otimes \O_{K_0}$, and set $
 	\varphi_Q: = \delta \sigma: \sigma^*Q\otimes K_0 \isom Q \otimes K_0$. This shows the existence of $(\tomega, \iota).$
 	
 	To show uniqueness, let $(\tomega, \iota),$ $(\tomega', \iota')$ be two pairs as in the lemma. Let $\omega = v \circ \tomega$ and $\omega' = v\circ \tomega'$.
 	By Lemma \ref{lem:explicitfiberfunctor}, there is a unique isomorphism  $\eta: \omega \isom \omega'$ of fiber functors over $\O_{K_0}$ taking $\iota$ to $\iota'.$ In particular the isomorphism $\eta(L) : \omega(L) \isom \omega'(L)$ is compatible with the Frobenii on $\omega(L)[1/p]$ and on $\omega'(L)[1/p]$.
 	Since for any $Q$ in $\Rep_{\Z_p} \cG$, the $\Q_p$-representation $Q\otimes \Q_p$ of $\cG$ is a subquotient of $(L\otimes \Q_p)^\otimes,$ this implies that the isomorphism $\eta(Q): \omega(Q) \isom \omega'(Q)$ is compatible with
 	the Frobenii on  $ \omega(Q)[1/p]$ and on $\omega'(Q) [1/p]$. Hence the isomorphism $\eta: \omega \isom \omega'$ comes from a (necessarily unique) isomorphism $\tomega \isom \tomega'$ in $\GIsoc_{K_0}^{\circ}$.
 \end{proof}
 \begin{lem}\label{lem:Isocisogeny} Let $(\tomega, \iota)$ be a pair as in Lemma \ref{lem:Isocexplicitfiberfunctor}, 
 	$\tomega'$ an $F$-isocrystal with $\cG$-structure, and $\iota'$ be an isomorphism $ \tomega'(L[1/p]) \isom (D[1/p], \varphi_D)$ in  $\Isoc_{K_0}$ such that $\iota'$ maps $\tomega'(s_{\alpha}) \in \tomega'(L[1/p]) ^{\otimes}$ to $s_{\alpha, 0} \in D[1/p]^{\otimes}$ for each $\alpha \in \bm \alpha$. Then there is a unique isomorphism $\tomega[1/p] \isom \tomega'$ in the category $\GIsoc_{K_0}$ which takes $\iota$ to $\iota'.$
 \end{lem}
 \begin{proof} Let $\omega =  v \circ \tomega[1/p]$ and $\omega' = v \circ \tomega' $. These are fiber functors $\Rep_{\Q_p} \cG \to \Modfp_{K_0}$. By Lemma \ref{lem:fiberfunctorisogeny}, there is a unique $\otimes$-isomorphism $\omega \isom \omega'$ which takes $\iota$ to $\iota'$. By exactly the same argument as in the proof of Lemma \ref{lem:Isocexplicitfiberfunctor}, this isomorphism comes from a unique isomorphism $\tomega[1/p] \isom \tomega'$ in $\GIsoc_{K_0}$.
 \end{proof}
 \subsection{}\label{para:delta assoc to isoc}
 Let $\cG$ be a smooth affine group scheme over $\Z_p$ with connected fibers. Let $\tomega$ be an integral $F$-isocrystal with $\cG$-structure over $K_0$. By Lemma \ref{lem:fiberfunctortorsor}, the set $P_{v\circ \tomega}(\O_{K_0})$ is a non-empty $\cG(\O_{K_0})$-torsor. Then we have the $\cG(K_0)$-torsor $Y(v\circ \tomega)$, and a canonical $\cG(\oo_{K_0})$-torsor $Y(v\circ \tomega)^{\circ} \subset Y(v\circ \tomega)$, as in \S \ref{para:defnY}. In the sequel, we shall write $Y(\tomega)$ and $Y(\tomega)^\circ$ for $Y(v\circ \tomega)$ and $Y(v \circ \tomega)^\circ$ respectively.
 
 Now consider an element $y\in Y(\tomega)$. For each $L$ in $\Rep_{\Z_p} \cG,$ the $K_0$-linear
 isomorphism $L\otimes_{\ZZ_p} K_0 \isom \tomega(L)[1/p]$ induced by $y$ allows us to view the Frobenius on $\tomega(L)[1/p]$ as being given by
 $\delta_{y,L}\cdot \sigma$ for some $\delta_{y,L} \in \GL(L \otimes K_0).$ Since every representation in $\Rep_{\Q_p} \cG$ contains a $\cG$-stable $\Z_p$-lattice (see \S \ref{para:setting for fiber functor}), the reconstruction theorem over a field implies that the elements $\delta_{y,L}$ for all $L$ come from a unique element $\delta_y \in \cG(K_0)$. Thus we have obtained a map $$ Y(\tomega) \To \cG(K_0), \quad y \longmapsto \delta_y.$$
 It is clear that the set $\set{\delta_y \mid y \in Y(\tomega)}$ is a $\cG(K_0)$-orbit under $\sigma$-conjugation, and the set $\set{\delta_y \mid y \in Y(\tomega)^\circ}$ is a $\cG(\O_{K_0})$-orbit under $\sigma$-conjugation. We call the last $\cG(\O_{K_0})$-orbit the \emph{invariant} of $\tomega$, and denote it by $\inv(\tomega)$.\index[Oct2026-sorted-notation]{inv 1 01 00 inv upsilon@$\mathrm{inv}(\Upsilon)$}
 
 \begin{lem}\label{lem:isocrystalclass} Keep the assumptions in \S \ref{para:delta assoc to isoc}. The construction $\tomega \mapsto \inv(\tomega)$ induces a bijection from the set of isomorphism classes in the category
 	$\GIsoc^\circ_{K_0}$, to the set of $\cG(\O_{K_0})$-$\sigma$-conjugacy classes in $\cG(K_0)$.
 \end{lem}
 \begin{proof} It is easy to check that the map described in the lemma is well defined. The inverse of the map is induced by the following construction: Given $\delta \in \cG(K_0),$ we define a functor $\tomega_{\delta}:\Rep_{\Z_p} \cG \rightarrow \Isoc^\circ_{K_0}$ by
 	sending $L$ to $L\otimes_{\Z_p}\O_{K_0}$ and equipping $L\otimes_{\Z_p} K_0$ with $\delta\sigma$ as the Frobenius.
 \end{proof}
 \begin{defn}\label{defn:intrinsic ADLV}\index{affine Deligne--Lusztig set}
 	Let $\cG$ be a reductive group scheme over $\Z_p$, and let $\mu: \GG_m \to \cG_{\oo_{K_0}}$ be a homomorphism of $\oo_{K_0}$-group schemes. Let $\tomega$ be an integral $F$-isocrystal with $\cG$-structure over $K_0$. With the notation as in \S \ref{para:delta assoc to isoc}, we define\index[Oct2026-sorted-notation]{x 0 00 23 tilde x subscript mu upsilon@$\tilde X_{\mu}(\Upsilon)$}
 	$$ \tilde X_{\mu}(\tomega) : = \set{ y \in Y(\tomega) \mid \delta_{y} \in \cG(\oo_{K_0}) \mu(p) \cG(\oo_{K_0}) }. $$ This is a subset of $Y(\tomega)$ stable under the $\cG(\oo_{K_0})$-action. We define\index[Oct2026-sorted-notation]{x 0 00 00 x subscript mu upsilon@$X_{\mu}(\Upsilon)$}
 	$$ X_{\mu} (\tomega) : = \tilde X_{\mu} (\tomega)/ \cG(\oo_{K_0}). $$
 \end{defn}
 \begin{rem}\label{rem:intrinsic ADLV} In the setting of Definition \ref{defn:intrinsic ADLV}, if we fix an element $y_0 \in Y(\tomega)$ and use it to identify the $\cG(K_0)$-torsor $Y(\tomega)$ with $\cG(K_0)$, then $X_{\mu}(\tomega)$ is identified with
 	$$\set{g \in \cG(K_0)/ \cG(\oo_{K_0}) \mid g^{-1} \delta_{y_0} \sigma (g)  \in \cG(\oo_{K_0}) \mu(p) \cG(\oo_{K_0}) }. $$
 	When $K_0 = \Q_p^{\ur}$, the above is the\index{affine Deligne--Lusztig set} affine Deligne--Lusztig set $X_{\mu}(\delta_{y_0})$ introduced in \S \ref{subsubsec:recall of defns in LR}.
 \end{rem}
 
 \section{Crystalline representations with \texorpdfstring{$G$}{G}-structure}\label{subsec:crystalline-representations-with-g-structure}
 This section begins with a review of the pairs $(\mu, \delta)$ used by Rapoport and Zink to parametrize (admissible) filtered $\varphi$-modules with additional structures. We then use Fontaine's crystalline functor $D_{\cris}$ to assign such a pair to a $G(\mathbb{Q}_p)$-valued crystalline representation. Finally, we investigate specific torus-valued crystalline representations arising from global constructions with a close relation to the reciprocity law in the theory of canonical models of Shimura varieties. 
 
 \subsection{}\label{subsubsec:RZ} Let $G$ be a reductive  group over $\Q_p.$ Throughout the remainder of this chapter, all finite extensions of $\Q_p$ are understood 
 to be taken inside $\Qpbar$. Let $K/\Q_p$ be a finite extension, and let $K_0$ denote the maximal unramified subextension of $\Q_p$ in $K$.\footnote{The notation $K$ and $K_0$ is 
 	standard in $p$-adic Hodge theory. In the context of Shimura varieties, where $K$ 
 	usually denotes the level subgroup, we will instead write $F$ and $F_0$ for local fields.}  We denote by $\RZ_G(K)$\index[Oct2026-sorted-notation]{rz 0 01 00 rz subscript g k@$\mathrm{RZ}_G(K)$} the set of pairs $(\mu,\delta)$ consisting of a $K$-homomorphism $\mu:\GG_{m,K} \to G_K$ and an element $\delta \in G(K_0)$. These pairs are considered by Rapoport--Zink in \cite[\S 1]{RZ96}.
 
 Let $\mathrm{MF}_K^{\varphi}$\index[Oct2026-sorted-notation]{mf 0 01 00 mf subscript k superscript varphi@$\mathrm{MF}_K^{\varphi}$} be the category of filtered $\varphi$-modules over $K$, with exact sequences taken to be strict for the filtrations. Write $v : \mathrm{MF}_K^{\varphi} \to \Modfp_{K_0}$ for the forgetful functor taking each filtered $\varphi$-module to its underlying $K_0$-vector space. Associated with each $(\mu,\delta)\in \RZ_G(K)$ we have a faithful exact $\otimes$-functor\index[Oct2026-sorted-notation]{i 0 03 00 i subscript mu delta@$\mathcal I_{\mu,\delta}$} $$\mathcal I_{\mu,\delta}: \Rep_{\QQ_p}G \To \mathrm{MF}_K^{\varphi}$$
 sending $V$ to $V \otimes_{\Q_p} K_0$ equipped with Frobenius $\delta \sigma$ and filtration on $(V\otimes K_0)\otimes_{K_0} K$ given by $\mu$.  
 Conversely, for every faithful exact $\otimes$-functor $\mathcal I: \Rep_{\QQ_p}G \to \mathrm{MF}_K^{\varphi}$ such that $v\circ \mathcal I$ is $\otimes$-isomorphic to the standard functor $V\mapsto V\otimes_{\Q_p}K_0$, we can choose such an isomorphism and choose a splitting of the resulting filtration on the functor $V \mapsto V\otimes_{\Q_p} K$ (by Deligne's splitting theorem \cite[IV, 2.4]{Saavedra}) to obtain $(\mu,\delta)\in \RZ_G(K)$ such that $\mathcal I$ is $\otimes$-isomorphic to $\mathcal I_{\mu,\delta}$. 
 
 We define an equivalence relation $\sim$ on $\RZ_G(K)$ by declaring $(\mu_1,\delta_1)\sim(\mu_2,\delta_2)$ if $\mathcal I_{\mu_1,\delta_1}$ and $\mathcal I_{\mu_2,\delta_2}$ are $\otimes$-isomorphic. 
 We call an element $(\mu,\delta)$ of $\RZ_G(K)$ \emph{admissible},
 if $\mathcal I_{\mu,\delta}$ lands in
 the subcategory $\wa{K}$\index[Oct2026-sorted-notation]{mf 0 01 00 mf subscript k superscript varphi a@$\mathrm{MF}_K^{\varphi, \mathrm{a}}$} of admissible filtered $\varphi$-modules.
 We call $(\mu,\delta)$  \emph{neutral}, if $\kappa_G(\delta) = [\mu] \in \pi_1(G)_{\Gamma_p}$, where $\kappa_G$ is the Kottwitz map as in \S \ref{para:B(G)}. We denote by $\RZ_G^{\adm}(K)$\index[Oct2026-sorted-notation]{rz 0 01 00 rz subscript g superscript a k@$\mathrm{RZ}_G^{\mathrm a}(K)$} (resp.~$\RZ_G^{\neu} (K)$)\index[Oct2026-sorted-notation]{rz 0 01 00 rz subscript g superscript n k@$\mathrm{RZ}_G^{\mathrm{n}}(K)$} the subset of admissible (resp.~neutral) elements of $\RZ_G(K)$, and write $\RZ_G^{\adm, \neu}(K)$\index[Oct2026-sorted-notation]{rz 0 01 00 rz subscript g superscript a n k@$\mathrm{RZ}_G^{\mathrm{a}, \mathrm{n}}(K)$} for $\RZ_G^{\adm}(K) \cap \RZ_G^{\neu} (K)$.
 
 It is easy to see (cf.~\cite[Def.~1.23]{RZ96}) that  $(\mu_1,\delta_1), (\mu_2,\delta_2 ) \in \RZ_G(K)$ are equivalent if and only if there exists $g \in G(K_0)$ such that $g\delta_1 \sigma(g)^{-1} = \delta_2$ and such that $\Int(g) \circ \mu_1$ and $\mu_2 $ define the same filtration on the functor $V \mapsto V\otimes_{\Q_p} K$. By \cite[IV, 2.2.5 (2)]{Saavedra}, the last condition   implies that $\mu_1$ and $\mu_2$ are conjugate by $G(K)$. It follows that the subset $\RZ_G^{\adm} (K)\subset \RZ_G(K)$ is stable under $\sim$. Clearly the subset $\RZ_G^{\neu} (K)\subset \RZ_G(K)$ is also stable under $\sim$.

 \begin{prop}\label{prop:admissible pair for torus}
 	Let $T$ be a torus over $\Q_p$. The  following statements hold. 
 	\begin{enumerate}
 		\item Two elements $(\mu_1,\delta_1), (\mu_2,\delta_2) \in \RZ_T(K)$ are equivalent if and only if $\mu_1 = \mu_2$ and $\delta_1$ is $\sigma$-conjugate to $\delta _2$ in $T(K_0)$.
 		\item We have $ \RZ_T^{\adm, \neu}(K)= \RZ_T^{\neu}(K) $.
 	\end{enumerate}
 \end{prop}
 \begin{proof} Part (i) follows immediately from the preceding discussion. For part (ii), the condition of being neutral implies the second condition in \cite[Prop.~1.21]{RZ96}. The result then follows from that proposition, descent of weak admissibility along unramified extensions, and the Colmez--Fontaine Theorem \cite{CF}.
 \end{proof}
 \subsection{}
 Keep the setting of \S \ref{subsubsec:RZ}. If $K'/K$ is a finite extension , then there is a natural map $\mathrm{RZ}_G(K) \to \mathrm{RZ}_G (K')$ preserving the equivalence relation, neutrality, and admissibility. We define\index[Oct2026-sorted-notation]{rz 0 01 00 rz subscript g@$\mathrm{RZ}_G$}\index[Oct2026-sorted-notation]{rz 0 01 00 rz subscript g superscript a@$\mathrm{RZ}_G^{\mathrm{a}}$}\index[Oct2026-sorted-notation]{rz 0 01 00 rz subscript g superscript n@$\mathrm{RZ}_G^{\mathrm{n}}$}\index[Oct2026-sorted-notation]{rz 0 01 00 rz subscript g superscript a n@$\mathrm{RZ}_G^{\mathrm{a},\mathrm{n}}$}
 $$ \mathrm{RZ}_G: = \varinjlim_{K} \mathrm{RZ}_G(K) ,$$ where $K$ runs through finite extensions of $\Q_p$. Similarly, we define the subsets $\RZ_G^{\adm}$, $\RZ_G^{\neu}$, and $\RZ_G^{\adm , \neu}$ of $\RZ_G$ by taking direct limits. We write $\sim$ for the finest equivalence relation on $\RZ_G$ such that each map $\RZ_G(K) \to \RZ_G$ preserves the equivalence relation. The three subsets of $\RZ_G$ introduced above are all stable under $\sim$.
 
 \begin{cor}\label{cor:pairs for torus} Let $T$ be a torus over $\Q_p$. We have $\mathrm{RZ}^{\adm,\neu}_T = \RZ^{\neu}_T$. There is a natural bijection between $\mathrm{RZ}^{\adm,\neu}_T/{\sim}$ and the set
 	\begin{align*}
 		\set{(\mu, [\delta]) \mid  \mu \in X_*(T), [\delta] \in T(\Q_p^{\ur})/(1-\sigma), \kappa_T([\delta]) = [\mu]} .
 	\end{align*}
 \end{cor}
 \begin{proof}
 	This follows from Proposition \ref{prop:admissible pair for torus}.
 \end{proof}
 \subsection{}\label{para:rational crys rep}
 Keep the setting of \S \ref{subsubsec:RZ}. We say that a homomorphism $$\rho: \Gamma_K \To G(\Q_p)$$ is a\index{G(Qp)-valued crystalline representation@$G(\Q_p)$-valued crystalline representation} \emph{$G(\Q_p)$-valued crystalline representation of $\Gamma_K$}, if for some faithful representation $V$ in $\Rep_{\Q_p} G$, the homomorphism $\Gamma_K \to \GL(V)(\Q_p)$ arising from $\rho$ is a crystalline representation. This condition is  independent of the choice of $V$, since the category of crystalline representations is closed under direct sums, tensor products, duals, and subquotients. We denote by $\Crys_G(K)$ the set of all such $\rho$, and define
 \index[Oct2026-sorted-notation]{crys 0 01 00 crys subscript g k@$\mathrm{Crys}_G(K)$}\index[Oct2026-sorted-notation]{crys 0 01 00 crys subscript g@$\mathrm{Crys}_G$}$$ \Crys_{G} : = \varinjlim_{K} \Crys_{G}(K),$$ where $K$ runs through finite extensions of $\Q_p$, and the transition maps are given by restriction. Denote by $\sim$ the equivalence relation on $\Crys_G$ given by $G(\Q_p)$-conjugation. 
 
 For $\rho \in \Crys _{G}(K)$, we denote its image in $\Crys _{G}$ by $[\rho]$. Associated with $\rho$ we have a faithful exact $\otimes$-functor $\cI_{\rho} : \Rep_{\Q_p} G \to \mathrm{MF}^{\varphi}_{K}$ taking each $V$ to $$D_{\cris}(V) = (B_{\cris} \otimes_{\Q_p} V)^{\Gamma_K} , $$ where $V$ is viewed as a crystalline representation of $\Gamma_{K}$ via $\rho$. (The fact that $\cI_{\rho}$ is a faithful exact $\otimes$-functor follows from the fundamental properties of $D_{\cris}$ \cite[\S 3.4]{Fontaine79}.) Since $\rho$ is crystalline, the formation of $D_{\cris}$ is compatible with base change. Hence after enlarging $K$ the fiber functor $v\circ \cI_{\rho}$ is $\otimes$-isomorphic to $\idfunc_{K_0}$. Then $\cI_{\rho}$ is $\otimes$-isomorphic to $\cI_{\mu,\delta}$ for some $(\mu,\delta)\in \mathrm{RZ}_G(K)$ which is unique up to equivalence; see \S \ref{subsubsec:RZ}. The image of $(\mu,\delta)$ in $\mathrm{RZ}_G/{\sim}$ is a well-defined invariant of $[\rho] \in \Crys_G$, and it depends only on the image of $[\rho]$ in  $\Crys_G/{\sim}$. Thus we obtain a map\index[Oct2026-sorted-notation]{d 0 00 00 d subscript cris superscript g@$D_{\mathrm{cris}}^G$} $$D_{\cris}^{G}: \Crys_G /{\sim} \To  \mathrm{RZ}_G/{\sim} .$$

 \begin{prop}\label{prop:adm RZ pair}
 	The map $D_{\cris}^{G}$ is injective with image $\mathrm{RZ}^{\adm,\neu}_G /{\sim}$.
 \end{prop}
 \begin{proof} Firstly, it is clear from the definitions that $\im (D_{\cris}^{G}) \subset \RZ_G^{\adm}/{\sim}$.
 	
 	For each finite extension $K/\Q_p$,  let $\Rep_{\Gamma_{K}}^{\cris}$ be the category of crystalline representations of $\Gamma_{K}$ over $\Q_p$. Then $D_{\cris} : \Rep_{\Gamma_{K}}^{\cris} \to \mathrm{MF}_{K}^{\varphi, \adm}$ is a $\otimes$-equivalence of $\otimes$-categories, with a quasi-inverse given by the functor $$V_{\cris} : D \longmapsto \Fil^0(B_{\cris}\otimes_{K_0} D) ^{\varphi= 1} $$ (which is also a $\otimes$-functor). Let $u: \Rep_{\Gamma_{K}}^{\cris} \to \Modfp_{\Q_p}$ be the functor sending a crystalline representation to its underlying $\Q_p$-vector space.
 	
 	By a result of Wintenberger \cite{Wintenberger97},   an element $(\mu,\delta) \in \RZ_G^{\adm} (K)$ is neutral if and only if the composite $\otimes$-functor $$u \circ V_{\cris} \circ \cI_{\mu,\delta} : \Rep_{\Q_p} G \To \Modfp_{\Q_p}$$ is $\otimes$-isomorphic to the standard fiber functor $\idfunc_{\Q_p}$.\footnote{This was originally a conjecture of Rapoport--Zink; see the paragraph below \cite[Prop.~1.20]{RZ96}. Wintenberger showed in \cite{Wintenberger97} that it is a consequence of the Colmez--Fontaine Theorem \cite{CF}.} It follows that $\im (D_{\cris}^{G}) \subset \RZ_G^{\adm, \neu}/{\sim}$.
 	
 	Now we construct a map $\mathrm{RZ}^{\adm,\neu}_G/{\sim} \to \Crys_G/{\sim}$ inverse to $D_{\cris}^{G}$. Let  $(\mu,\delta) \in \mathrm{RZ}^{\adm,\neu}_G(K)$. By the aforementioned result, the $\otimes$-functor $\cF: = u \circ V_{\cris} \circ \cI_{\mu,\delta}$ is $\otimes$-isomorphic to $\idfunc_{\Q_p}$. By composing the tautological homomorphism $\Gamma_{K} \to \underline{\Aut}^{\otimes} (u) (\Q_p)$ with the natural homomorphism $\underline{\Aut}^{\otimes} (u) (\Q_p) \to \underline{\Aut}^{\otimes} (\cF) (\Q_p)$, we obtain a homomorphism $\rho': \Gamma_{K} \to \underline{\Aut}^{\otimes} (\cF) (\Q_p)$. By choosing a $\otimes$-isomorphism between $\idfunc_{\Q_p}$ and $\cF$, we identify the $\Q_p$-group $\underline{\Aut}^{\otimes} (\cF)$ with
 	$G$, and identify $\rho'$ with a homomorphism $\rho : \Gamma_{K} \to G(\Q_p)$. Clearly the $G(\Q_p)$-conjugacy class of $\rho $
 	is independent of choices. It is straightforward to check that the construction $(\mu,\delta) \mapsto \rho$ gives rise to the desired inverse to $D_{\cris}^{G}$.
 \end{proof}
 
 \begin{cor}\label{cor:classification of T-valued}
 	Let $T$ be a torus over $\Q_p$. The map $D_{\cris}^T$ induces a bijection $$\Crys_T \isom 	\set{(\mu, [\delta]) \mid  \mu \in X_*(T), [\delta] \in T(\Q_p^{\ur})/(1-\sigma), \kappa_T([\delta]) = [\mu]} .$$
 \end{cor}
 \begin{proof}
 	This follows from Corollary \ref{cor:pairs for torus} and Proposition \ref{prop:adm RZ pair}. (Note that $\Crys_T = \Crys_T /{\sim}$.)
 \end{proof}

 \subsection{}\label{para:motivic element} Let $T$ be a cuspidal torus over $\Q$ (see Definition \ref{defn:anti-cusp}). We call an element $\delta \in T(\Q_p^{\ur})$ \emph{motivic}, if\index{motivic element} for some $n \geq 1,$ the element $\gamma =  \delta\sigma(\delta) \dots \sigma^{n-1}(\delta)$ is in $T(\Q)$, and is a $p$-unit (i.e., $\gamma$ lies in a compact open subgroup of $T(\A_f^p)$). We denote by\index[Oct2026-sorted-notation]{t 0 00 00 t q subscript p superscript ur superscript mot@$T(\mathbb Q_p^{\ur})^{\mot}$} $T(\Q_p^{\ur})^{\mot} \subset T(\Q_p^{\ur})$ the subset of motivic elements. Note that $T(\Q_p^{\ur})^{\mot} $ is stable under $\sigma$-conjugation by $T(\Q_p^{\ur})$. We denote by $\cT^{\circ}$\index[Oct2026-sorted-notation]{t 0 03 00 t superscript circ@$\mathcal T^{\circ}$} the connected N\'eron model of $T_{\Q_p}$ over $\ZZ_p$. Let\index[Oct2026-sorted-notation]{0003 sim 1 00 26 overset circ sim@$\overset{\circ}{\sim}$} $\circsim$ be the equivalence relation on $T(\Q_p^{\ur})^{\mot}$ defined by $\sigma$-conjugation by $\cT^{\circ} (\ZZ_p^{\ur})$. As in \S \ref{para:B(G)}, we denote by $w_T : T(\LL) \to X_*(T)_{\Gamma_{p,0}}$ the Kottwitz homomorphism, which is surjective. By \cite[Rmk.~2.2 (iii)]{rapoportguide}, we have $\ker (w_T) = \cT^{\circ}(\breve \Z_p)$. Hence the restriction of $w_T$ to $T(\Q_p^{\ur})^{\mot}$ factors through $T(\Q_p^{\ur})^{\mot}/{\circsim}$.

 \begin{lem}\label{lem:motelts}In the setting of \S \ref{para:motivic element}, the map $w_T$ induces a bijection
 	\begin{align}\label{eq:bij for mot}
 		T(\Q_p^{\ur})^{\mot}/{\circsim} \isom  X_*(T)_{\Gamma_{p,0}}.
 	\end{align}
 \end{lem}
 \begin{proof} The surjectivity of (\ref{eq:bij for mot}) follows from the second construction in \cite[\S 4.3.9]{kisin2012modp}. We explain this in more detail. Let $\mu \in X_*(T)$, and let $L/\Q$ be a finite Galois extension inside $\Qbar$ splitting $T$. Let $L_p$ be the completion of $L$ at the place above $p$ determined by the fixed embedding $\Qbar \hookrightarrow \Qpbar$, and let $\pi$ be a uniformizer of $L_p$. Then the same argument as in \textit{loc.~cit.~}(with $\mu_{h_T}$ replaced by $\mu$) shows that there exist $s \in \NN $ and a $p$-unit $\gamma\in T(\Q)$ such that $\N_{L_p/\Q_p}(\mu(\pi))^s \gamma^{-1} \in \cT^{\circ} (\Z_p)$. By Greenberg's theorem \cite[Prop.~3]{Greenberg}, the map $\cT^{\circ}(\breve \Z_p ) \to \cT^{\circ}(\breve \Z_p ), c \mapsto c\sigma(c)^{-1}$ is surjective. Hence we can find $c\in \cT^{\circ}(\breve \Z_p )$ such that\footnote{In the last paragraph of \cite[\S 4.3.9]{kisin2012modp}, it is used that such $c$ can be found in $T(\LL)$.}
 	$$ c \N_{L_p/\Q_p} (\mu(\pi)) ^s \sigma^n (c)^{-1} = \gamma, $$
 	where $n = s [L_{p,0}: \Q_p].$
 	Let $$\delta : = c \N_{L_p/L_{p,0}} (\mu(\pi)) \sigma (c)^{-1} \in T(\LL), $$ where $L_{p,0}$ is the maximal unramified extension of $\Q_p$ inside $L_p$. As in \textit{loc.~cit.,} we have $\delta \sigma(\delta) \cdots \sigma^{n-1}(\delta) = \gamma$, and we have $\delta \in T(\Q_{p^n}).$ Hence $\delta \in T(\Q_p^{\ur})^{\mot}$. 
 	
 	We now check that $w_T(\delta)$ equals the image of $\mu$ in $X_*(T)_{\Gamma_{p,0}}$, which will prove the surjectivity of (\ref{eq:bij for mot}). For this, it suffices to show that  $w_T(\N_{L_p/L_{p,0}} (\mu(\pi) ))$  equals the image of $\mu$, since $w_T(c)$ is trivial. Writing $F$ for $L_{p,0}$, we have $X_*(T)_{\Gamma_F} = X_*(T)_{\Gamma_{p,0}}$. Let $\sigma_F = \sigma ^{[F: \Q_p]}$. By \cite[\S 2.5]{kottwitzisocrystal}, $w_T$ induces a bijection $$\B_F(T): = \set{\sigma_F\text{-conjugacy classes in }T(\LL)} \isom X_*(T)_{\Gamma_F},$$ whose inverse is induced by $\mu \mapsto \N_{L_p / F} (\mu(\pi)). $ (Here we use that $T$ splits over $L_p$ and that $L_p/F$ is totally ramified.) This gives what we want.
 	
 	For the injectivity of (\ref{eq:bij for mot}), let $\delta_1, \delta_2 \in T(\Q_p^{\ur})^{\mot}$ be such that $w_T(\delta_1) =w_T(\delta_2)$. Let $\delta = \delta_1 \delta_2^{-1}$. Then $\delta$ is also motivic, so we can choose $n$ such that $\gamma: = \delta \sigma (\delta) \cdots \sigma^{n-1}(\delta)$ is a $p$-unit in $T(\Q)$. Since $w_T(\delta) = 0$, we have $\delta \in \cT^{\circ}(\breve \Z_p)$, and in particular $\gamma \in \cT^{\circ} (\Z_p)$. Therefore $\gamma$ lies in a congruence subgroup of $T(\Q)$, and has finite order by Lemma \ref{lem:cusp TFAE}. Up to enlarging $n$, we may assume that $\gamma = 1$.
 	Again by Greenberg's theorem, we can write $\delta \in \cT^{\circ} (\breve \Z_p)$ as $c \sigma (c)^{-1}$, for some $c \in\cT^{\circ} (\breve \Z_p)$. Since $\gamma=1$, we have $c \sigma^n (c)^{-1} =1$, i.e., $c \in
 	\cT^{\circ} (\ZZ_{p^n})$. Hence $\delta_1 = \delta \delta_2 =  c \delta_2 \sigma(c)^{-1}$ with $c\in \cT^\circ (\Z_p^{\ur})$, which means $\delta_1 \circsim \delta_2.$
 \end{proof}
 \begin{lem}\label{lem:mot rational}
 	In the setting of \S \ref{para:motivic element}, the map $w_T$ induces a bijection from the set of $T(\Q_p^{\ur})$-$\sigma$-conjugacy classes in $T(\Q_p^{\ur})^{\mot}$ to $X_*(T)_{\Gamma_p}$.
 \end{lem}
 \begin{proof}
 	The map $w_T$ induces a group isomorphism $T(\LL)/\cT^{\circ}(\breve \Z_p) \isom X_*(T)_{\Gamma_{p,0}}$. In view of the density of $T(\Qpur)$ in $T(\LL)$, this induces a group isomorphism $$f: T(\Q_p^{\ur})/\cT^{\circ}(\Z_p^{\ur}) \isom X_*(T)_{\Gamma_{p,0}},$$ which is equivariant for the natural actions of $\sigma$ on the two sides. Now consider the two sides of the bijection (\ref{eq:bij for mot}). On the left hand side, the group $T(\Q_p^{\ur})/\cT^{\circ}(\Z_p^{\ur})$ acts by $\sigma$-conjugation. On the right hand side, the group $X_*(T)_{\Gamma_{p,0}}$ acts (on itself) by $\sigma$-conjugation. Moreover, the bijection (\ref{eq:bij for mot}) is equivariant with respect to these two actions if we identify the two acting groups via $f$. Taking orbits on the two sides we obtain the desired bijection.
 \end{proof}
 
 \begin{defn}\label{defn:motivic rep}
 	Let $T$ be a cuspidal torus over $\Q$. Let\index[Oct2026-sorted-notation]{mot 0 01 00 mot subscript t@$\mathrm{Mot}_T$} $\Mot_T$ be the subset of $\Crys_{T_{\Q_p}}$ consisting of those $[\rho]$ whose image under the bijection in Corollary \ref{cor:classification of T-valued} is of the form $(\mu,[\delta])$ where $[\delta]\in T(\Q_p^{\ur})/(1-\sigma)$ is in the image of $T(\Q_p^{\ur})^{\mot}$. Note that the definition of $\Mot_T$ depends on $T$ over $\Q$, not just $T_{\Q_p}$.
 \end{defn}
 
 \begin{prop}\label{prop:classifyingmotrep} In the setting of Definition \ref{defn:motivic rep}, the $\mu$-component of the map $D_{\cris}^{T_{\Q_p}}$ induces a bijection\index[Oct2026-sorted-notation]{m 0 05 00 m subscript t@$\mathscr M_T$}
 	$$\motinv_T: \Mot_T \isom X_*(T).$$
 \end{prop}
 \begin{proof}This follows from Corollary \ref{cor:classification of T-valued} and Lemma \ref{lem:mot rational}.
 \end{proof}
 \subsection{}\label{para:r(mu)_U}
 Let $T$ be a cuspidal torus over $\QQ$. We now use class field theory to construct certain $T(\QQ_p)$-valued global Galois representations, and show that their localizations at places above $p$ give rise to elements of $\Mot_T$.
 
 Let $\mu \in X_*(T)$, and let $E_\mu\subset \Qbar$ be the field of definition of $\mu$. Similar to \S \ref{para:Shimura variety}, we consider the composite  homomorphism of $\Q$-algebraic groups \index[Oct2026-sorted-notation]{r 1 00 00 r mu superscript alg@$r(\mu)^{\alg}$}
 \begin{align*}
 	r(\mu)^{\alg}: \Res_{E_\mu/\Q} \GG_m \xrightarrow{\Res_{E_\mu/\Q} \mu} \Res_{E_\mu/\Q} T \xrightarrow{\N_{E_\mu/\Q}} T.
 \end{align*} We have an induced homomorphism between topological groups
 \begin{align}\label{eq:ST rec law}
 	E_\mu^{\times}  \backslash \A_{E_\mu}^{\times} \To T(\Q) \backslash T(\A).
 \end{align}
 By Lemma \ref{lem:cusp TFAE}, $T(\Q)$ is discrete in $T(\A_f)$, and so we have $$T(\Q) \backslash T(\A_f) = \varprojlim_{U} T(\Q) \backslash T(\A_f)/ U, $$ where $U$ runs through compact open subgroups of $T(\A_f)$. For each such $U$, we have a natural map $\pi_0(T(\Q) \backslash T(\A)) \to  T(\Q) \backslash T(\A_f)/ U$, cf.~\S \ref{para:Shimura variety}. In the limit we obtain a map
 \begin{align}\label{eq:pi_0T}
 	\pi_0(T(\Q) \backslash T(\A)) \To  T(\Q) \backslash T(\A_f).
 \end{align} The composition
 \begin{align}\label{eq:ST law}   		E_\mu^{\times}  \backslash \A_{E_\mu}^{\times} \xrightarrow{(\ref{eq:ST rec law})}   T(\Q) \backslash T(\A)  \To \pi_0(T(\Q) \backslash T(\A)) \xrightarrow{(\ref{eq:pi_0T})}   T(\Q) \backslash T(\A_f)  
 \end{align}
 factors through the global Artin map $E_\mu^{\times}  \backslash \A_{E_\mu}^{\times} \to \pi_0(E_\mu^{\times}  \backslash \A_{E_\mu}^{\times} ) \cong \Gal( E^{\ab}_\mu/E_\mu)$. (Recall from \S \ref{para:Shimura variety} that we take the geometric normalization of the global Artin map.) We thus obtain a map\index[Oct2026-sorted-notation]{r 1 00 00 r mu@$r(\mu)$}
 $$ r(\mu): \Gal(E^{\ab}_\mu/E_\mu) \To T(\Q)\backslash T(\A_f).$$
 
 Let $U \subset T(\A_f)$ be a neat compact open subgroup. Since $T(\Q)$ is discrete in $T(\A_f)$, we have $T(\Q) \cap U = \set{1}$ (cf.~the proof of Lemma \ref{lem:Z_s}). The kernel of the projection $$\pi_U: T(\Q) \backslash T(\A_f) \To T(\Q) \backslash T(\A_f)/U$$ is $T(\Q) \backslash T(\Q) U$, which we identify with $U$, using that $T(\Q) \cap U = \set{ 1}$. Denote by\index[Oct2026-sorted-notation]{e 0 00 00 e subscript mu u@$E_{\mu,U}$} $E_{\mu,U}$ the fixed field of $\ker(\pi_U \circ r(\mu))$ in $E^{\ab}_\mu$; this is a finite extension of $E_\mu$. Then $r(\mu)$ induces a homomorphism
 $$ r(\mu)_U : \Gal(E^{\ab}_\mu/E_{\mu,U}) \To \ker \pi_U \cong U. $$
 We denote by $r(\mu)_{U,p}$ the composite homomorphism
 $$ \Gal (E^{\ab}_\mu/E_{\mu,U}) \xrightarrow{r(\mu)_U}   U \hookrightarrow T(\A_f) \xrightarrow{\text{projection}}   T(\Q_p).  $$
 
 The fixed embeddings $E_{\mu,U} \hookrightarrow \Qbar \hookrightarrow \Qpbar$ give rise to a place $v$ of $E_{\mu,U}$ above $p$. Let $K = E_{\mu, U,v} \subset \Qpbar$. We denote by\index[Oct2026-sorted-notation]{r 1 00 00 r mu subscript u p loc@$r(\mu)_{U,p,\mathrm{loc}}$} $r(\mu)_{U,p,\mathrm{loc}}$ the composite map
 $$ \Gamma_{K} = \Gal (\Qpbar/K) \To  \Gal (E^{\ab}_\mu/E_{\mu, U}) \xrightarrow{r(\mu)_{U,p}}   T(\Q_p).    $$
 
 \begin{prop}\label{prop:local component} In the setting of \S \ref{para:r(mu)_U}, assume in addition that $U$ is of the form $U^p U_p$, where $U^p$ is a neat compact open subgroup of $T(\A_f^p)$ and $U_p$ is a compact open subgroup of $T(\Q_p)$. Then $ r(\mu)_{ U, p , \mathrm{loc}} : \Gamma_K \to T(\Q_p)$ is a $T(\Q_p)$-valued crystalline representation. The element $[r(\mu)_{ U, p , \mathrm{loc}}] \in \Crys_{T_{\Q_p}}$ lies in $\Mot_T$. Moreover, the image of $[r(\mu)_{ U, p , \mathrm{loc}}]$ under the bijection $\motinv_T$ in Proposition \ref{prop:classifyingmotrep} is $-\mu$.
 \end{prop}
 \begin{proof} Let $f$ be the composite map of topological groups
 	$$   K^{\times} \xrightarrow{\mathrm{Art}_K}  \Gamma_K^{\ab} \To  \Gal(E^{\ab}_\mu/E_{\mu, U}) \xrightarrow{r(\mu)_U}  U, $$ where $\mathrm{Art}_K$ is the local Artin map (normalized geometrically). Let $F$ be the completion of $E_\mu$ inside $K = E_{\mu, U,v}$. Let $f_1$ be the composite map of $\Q_p$-algebraic groups
 	$$ \Res_{K/\Q_p}  \GG_m \xrightarrow{\N_{K/F}} \Res_{F/\Q_p} \GG_m \hookrightarrow (\Res_{E_\mu/\Q} \GG_m ) \otimes_{\Q} \Q_p \xrightarrow{r(\mu)^{\alg}\otimes_{\Q} \Q_p} T_{\Q_p}.$$ Then $f_1$ induces a map
 	$$   K^{\times} = (\Res_{K/\Q_p} \GG_m)(\Q_p) \overset{f_1}{\To} T(\Q_p) \hookrightarrow T(\A_f),   $$ which we again denote by $f_1$.
 	
 	We claim that $f$ and $f_1$ induce the same map $\O_K^{\times} \to T(\A_f)$. In fact, by the definition of $r(\mu)_U,$ 
 	the fact that inclusions of Galois groups correspond under the Artin map to taking norms, and 
 	the compatibility of the local and global Artin maps, we have
 	\begin{align}\label{eq:f_1 and f}
 		f_1(x) \in f(x) T(\Q) \subset T(\A_f) ,\quad \forall x \in K^\times.\end{align}
 	Take $x \in \O_K^\times$, and let $\gamma \in T(\Q)$ be such that $f_1(x) = f(x) \gamma$. Note that $f_1(x)$ lies inside the maximal compact subgroup $U_{p,\max}$ of $T(\Q_p)$, by the compactness of $\O_K^{\times}$. Hence $\gamma$ lies in $U^p U_{p,\max}$, which is a neat compact open subgroup of $T(\A_f)$ by the neatness of $U^p$. Since $T(\Q)$ is discrete in $T(\A_f)$, we have $\gamma=1$. The claim is proved.
 	
 	Now to check that $r(\mu)_{U,p , \mathrm{loc}}$ is crystalline, we take an arbitrary representation $V$ of $T_{\Q_p}$ and check that 
 	\begin{align*}  \Gamma_K \xrightarrow{r(\mu)_{U,p, \mathrm{loc}}}   T(\Q_p) \To   \GL(V)(\Q_p)  
 	\end{align*}  is crystalline. For this, it suffices to check that the composition 
 	\begin{align*} 
 		K^\times \overset{f}{\To}   U \To  T(\Q_p) \To   \GL(V)(\Q_p)  
 	\end{align*}
 	agrees with a $\Q_p$-algebraic group homomorphism $\Res_{K/\Q_p} \GG_m \to \GL(V)$ on $\O_K^{\times}$, by a well-known criterion in $p$-adic Hodge theory (see for instance \cite[Prop.~B.4 (i)]{conradlifting} and the remark following it). But this follows immediately from our claim proved above.
 	
 	We now check that $[r(\mu)_{ U, p , \mathrm{loc}}] \in \Crys_{T_{\Q_p}}$ lies in $\Mot_T$. Let $[\delta]\in T(\Q_{p}^{\ur})/(1- \sigma)$ be the element attached to  $[r(\mu)_{ U, p , \mathrm{loc}}]$ as in Corollary \ref{cor:classification of T-valued}, and let $\delta \in T(\Q_{p^n})$ be a representative of $[\delta]$. By taking a faithful representation of $T_{\Q_p}$ and applying \cite[Prop.~B.4 (ii)]{conradlifting}, we know that up to enlarging $n$ the element $$ \gamma_n: = \delta \sigma(\delta)\cdots \sigma^{n-1}(\delta) \in T(\Q_p)$$ is equal to\footnote{Note that in \cite{conradlifting} the arithmetic normalization of the local Artin map is used, which is opposite to our normalization. This results in the sign difference in the exponent in the expression below.}
 	$$[f(\pi)_p^{-1} f_1(\pi) ] ^{-n/n_K}\in T(\Q_p). $$
 	Here $\pi$ is a uniformizer of $K$, $n_K$ is the residue degree of $K$, and $f(\pi)_p$ denotes the component at $p$ of $f(\pi) \in U$. It remains to show that $f(\pi)_p^{-1} f_1(\pi)$ is a $p$-unit in $T(\Q)$.
 	By (\ref{eq:f_1 and f}) there exists $\gamma \in T(\Q)$ such that $f_1(\pi) = f(\pi) \gamma \in T(\A_f)$. Then $f(\pi)^{-1}_p f_1(\pi)$ equals the image of $\gamma$ in $T(\Q_p)$. In addition, $\gamma$ and $f(\pi)^{-1}$ have the same image in $T(\A_f^p)$, which shows that $\gamma \in U^p$. Hence $\gamma$ is a $p$-unit.
 	
 	Finally, we check that  $\motinv_T([r(\mu)_{U,p , \mathrm{loc}}])  = -\mu$. Let $\CC_p$ be the completion of $\Qpbar$. For each faithful representation $V$ of $T_{\Q_p}$, the $\Gamma_K$-representation on $V$ induced by $r(\mu)_{U,p , \mathrm{loc}}$ is Hodge--Tate, and we have the cocharacter $h_V: \GG_{m,\CC_p} \to \GL(V)_{\CC_p}$ as in \cite[\S 1.4]{Serre79}. We know that $h_V$ factors through $T_{\CC_p}$ (see \textit{loc.~cit.}), and the resulting \textit{Hodge--Tate cocharacter} $\mu_{\mathrm{HT}} \in X_*(T)$ is independent of the choice of $V$ by the functoriality of the construction. Since the filtration on $D_{\dR}(\QQ_p(1))$ jumps at $-1$, it is easy to see that $\motinv_T([r(\mu)_{U,p , \mathrm{loc}}])  = - \mu_{\mathrm{HT}}$. We are left to check that $\mu_{\mathrm{HT}} = \mu$.
 	
 	Let $T' = \Res_{K/\Q_p} \GG_m$. We identify $T'_{\Qpbar}$ with $\prod_{\tau \in \Hom_{\Q_p}(K , \Qpbar)} \GG_{m,\Qpbar}$, and define $\mu' \in X_*(T')$ by $\mu'(z) = (z,1,\ldots, 1)$, where the first entry corresponds to the canonical embedding $K \hookrightarrow \Qpbar$. Let $r' : \Gamma_{K} \to \oo_K^{\times}$ be the Lubin--Tate character for a choice of a uniformizer (cf.~\cite[\S 2.1]{Serre79}). Then $r'$ is a $T'(\Q_p)$-valued crystalline representation. Since $\mathrm{Art}_K \circ r'|_{\Gamma_{K,0}}$ is the composition of the inclusion $\Gamma_{K,0} \hookrightarrow \Gamma_K$ followed by the projection $\Gamma_K \to \Gamma_K^{\ab}$ (thanks to the geometric normalization of $\mathrm{Art}_K$), the restriction of $r(\mu)_{U,p}$ to $\Gamma_{K,0}$ equals the composition $$   \Gamma_{K,0} \xrightarrow{r'|_{\Gamma_{K,0}}}    \oo_K^{\times} \overset{f}{\To}  U \To T(\Q_p) .   $$ Moreover, by our previous claim that $f$ and $f_1$ induce the same map $\O_K^{\times} \to T(\A_f)$, we know that the above composition is equal to the composition of $r'|_{\Gamma_{K,0}}$ followed by $f_1: K^{\times } \to T(\Q_p)$. Therefore if we let $\mu'_{\mathrm{HT}} \in X_*(T')$ be the Hodge--Tate cocharacter of $r'$, then $\mu_{\mathrm{HT}}$ equals $f_1 \circ \mu'_{\mathrm{HT}}$. (Recall that $f_1$ is an algebraic homomorphism $T' \to T_{\Q_p}$.) By the last paragraph of \cite[\S 2.1]{Serre79}, we have $\mu'_{\mathrm{HT}} = \mu'$. Therefore $\mu_{\mathrm{HT}} = f_1 \circ \mu'$, which is easily seen to be equal to $\mu$.
 \end{proof}
 \section{Crystalline lattices with \texorpdfstring{$\cG$}{G}-structure}\label{subsec:crystalline lattices}
 This section constitutes the technical core of the chapter. We apply the ``Breuil--Kisin module'' functor $\mathfrak M$ to a $\mathcal G(\Z_p)$-valued crystalline representation to obtain a ``Breuil--Kisin module with $\mathcal G$-structure'', i.e., a $\otimes$-functor $\Rep_{\Z_p} \cG \to \Mod_{/\gS}^{\varphi}$. Assuming $\mathcal G$ is smooth affine  with connected fibers and satisfies a condition called \emph{$\mathfrak S$-purity}, we show that this $\otimes$-functor induces a fiber functor over $\mathfrak S$, thereby giving rise to an integral $F$-isocrystal with $\mathcal G$-structure. 
 
 Using the theory of Kottwitz homomorphisms in families developed in \S \ref{subsec:Kottwitz hom}, we achieve two further goals. First, we compute the image of the crystalline Frobenius under the Kottwitz homomorphism. Second, we compare the relative position of two lattices in the crystalline representation with the relative position of the resulting lattices in the $F$-isocrystal, by considering their respective images under the Kottwitz homomorphism.  
 \subsection{}\label{gSmodules}
 Let $K$ be a finite extension of $\Q_p$ with residue field $k$. Let $K_0$ be the maximal unramified subextension of $\Q_p$ in $K$, and let $\sigma \in \Aut(K_0)$ be the arithmetic $p$-Frobenius.
 
 We write $W$ for $W(k) = \O_{K_0}$. We set $\gS = W\lps u \rps$,\index[Oct2026-sorted-notation]{s 0 04 00 s@$\mathfrak S$} and let $\varphi$ be the endomorphism of $\gS$ that restricts to $\sigma$ on $W$ and sends $u$ to $u^p$. Fix a uniformizer $\pi$ of $K$, and let $E = E(u)\in \gS$ be its Eisenstein polynomial over $K_0$. We have a $W$-algebra isomorphism $\gS/(E) \isom \O_K$ sending $u$ to $\pi$. Thus we have specialization maps $\gS \to \oo_K$ and $\gS \to W$, sending $u$ to $\pi$ and to $0$ respectively. Using these two maps we also view $K$ and $K_0$ as $\gS$-algebras respectively. (We shall never consider the composition of $\gS \to W$ with $W \hookrightarrow \O_K$.)
 
 For any height $1$ prime ideal $\fkp \subset \gS$, the localization ${\gS}_{\fkp}$ is a DVR, and we write $\widehat {\gS}_{\fkp}$ for its completion.
 
 We denote by $\Mod_{/\gS}^{\varphi}$\index[Oct2026-sorted-notation]{mod 0 01 00 mod subscript s superscript varphi@$\mathrm{Mod}_{/\mathfrak S}^{\varphi}$} the category of pairs $(\gM, \varphi_{\gM})$, where $\gM$ is a finite free $\gS$-module, and $\varphi_{\gM}$ is a $\gS$-module
 isomorphism $\varphi^*\gM[1/E] \isom \gM[1/E].$ For such a pair $(\gM, \varphi_{\gM})$, the $\gS$-module $\varphi^* \gM$ carries a filtration, given by
 $$\Fil^i(\varphi^*\gM )= \varphi_{\gM}^{-1} (E^i\gM) \cap \varphi^*\gM \subset \varphi^*\gM[1/E]$$ for $i \in \mathbb Z$.
 
 Let\index[Oct2026-sorted-notation]{rep 0 01 00 rep subscript gamma subscript k superscript cris circ@$\mathrm{Rep}_{\Gamma_K}^{\cris \circ}$}
 $\Rep_{\Gamma_K}^{\cris \circ}$ be the category of $\Gamma_K$-stable $\Z_p$-lattices in  crystalline representations of $\Gamma_K$ over $\Q_p$.
 Recall from \cite[\S 1.2]{kisin2010integral} that there is a faithful $\otimes$-functor\index[Oct2026-sorted-notation]{m 0 04 00 m@$\mathfrak M$} $$\gM: \Rep_{\Gamma_K} ^{\cris \circ} \To  \Mod_{/\gS}^{\varphi}. $$ For $L$ in $\Rep_{\Gamma_K}^{\cris \circ}$, the following statements hold, where all isomorphisms are functorial in $L$ and compatible with tensor products. 
 \begin{enumerate}
 	\item There is a canonical isomorphism
 	$$\gM(L) \otimes_{\gS} K_0  \cong D_{\cris}(L\otimes \Q_p) $$
 	of isocrystals over $K_0$. The Frobenius on the left is induced by $\varphi_{\gM(L)}$. 
 	\item  There is a canonical isomorphism
 	\begin{align}
 		\label{eq:filtered iso general 1}
 		\varphi^*\gM (L) \otimes _{\gS} K  \cong  D_{\dR}(L\otimes \Q_p)
 	\end{align}
 	of filtered $K$-vector spaces. The filtration on the left is induced by the filtration on $\varphi^*\gM(L)$. 
 	\item There is a faithfully flat and formally \'etale $\hat \gS_{(p)}$-algebra\index[Oct2026-sorted-notation]{o 0 03 22 widehat o subscript e superscript ur@$\widehat {\mathcal O}_{\mathcal E^{\ur}}$} $\widehat \O_{\mathcal E^{\ur}},$ and a canonical $\widehat \O_{\mathcal E^{\ur}}$-linear isomorphism
 	\begin{align*}
 		\widehat \O_{\mathcal E^{\ur}} \otimes_{\Z_p} L \cong \widehat \O_{\mathcal E^{\ur}}\otimes_{\gS}\gM(L).
 	\end{align*}
 \end{enumerate}
 
 For $L$ in $\Rep_{\Gamma_K}^{\cris \circ}$, we set\index[Oct2026-sorted-notation]{m 0 00 00 m subscript cris@$M_{\cris}$} $$M_{\cris}(L) : = \gM(L) \otimes_{\gS} W.$$ Thus $M_{\cris}$ is a $\otimes$-functor
 $$ \Rep_{\Gamma_K} ^{\cris \circ}  \To \Isoc_{K_0}^{\circ}. $$ By (i) above, $M_{\cris}(L)$ is a $W$-lattice in the $K_0$-vector space $D_{\cris} (L \otimes \Q_p)$. The following property is proved by Tong Liu in \cite[\S 4]{Liucompatibility}.
 \begin{enumerate}
 	\setcounter{enumi}{3}
 	\item Inside $D_{\cris} (L \otimes \Q_p)$, the $W$-lattice $M_{\cris}(L)$ is independent of the choice of a uniformizer in $K$ (which is needed to define the functor $\gM$). Moreover, if $K'/K$ is a finite extension  with residue field $k'$ and if $L' \in \Rep_{\Gamma_{K'}}^{\cris \circ}$ is $L$ equipped with the inherited $\Gamma_{K'}$-action, then we have a canonical identification $$ M_{\cris} (L') \cong M_{\cris}(L) \otimes_{W(k)} W(k').$$ This is compatible with the usual identification
 	$$ D_{\cris} (L') \cong D_{\cris}(L) \otimes_{K_0} K_0',$$ where $K_0'$ denotes $W(k')[1/p]$.
 \end{enumerate}
 
 The filtered isomorphism (\ref{eq:filtered iso general 1}) is induced by another canonical filtered isomorphism, which we now recall. Note that $\widehat{\gS}_{(E)}$ is a complete DVR with residue field $K$, which has characteristic zero. 
 Hence there is a unique $K_0$-linear multiplicative lifting $K\to \widehat{\gS}_{(E)}$, making $\widehat{\gS}_{(E)}$ a $K$-algebra. We equip it with its $E$-adic filtration. 
 Using the proof of \cite[Thm.~1.2.1]{kisin2010integral}, and \cite[Lem.~1.2.12 (4)]{kisin2006crystalline}, 
 one obtains, functorially in $L$ and compatibly with tensor products, a canonical filtered isomorphism
 \begin{align}\label{eq:filtered iso general 2}
 	\varphi^* \gM(L) \otimes_{\gS} \widehat \gS_{(E)} \isom  D_{\dR}(L\otimes \Q_p)\otimes_K \widehat \gS_{(E)}. 
 \end{align} The isomorphism (\ref{eq:filtered iso general 1}) is induced by (\ref{eq:filtered iso general 2}) by passing to the residue field $K$ of $\widehat{\gS}_{(E)}$.

 \subsection{} \label{para:crystallinerep} Let $\cG$ be a flat, finite-type, affine group scheme over $\Z_p$. Write $G$ for $\cG_{\QQ_p}$, which we shall always assume to be reductive. We say that a homomorphism
 $$ \rho: \Gamma_K \longrightarrow \cG(\Z_p)$$ is a (\emph{$\cG(\Z_p)$-valued) crystalline representation}, if the composition of $\rho$ with the inclusion $\cG(\Z_p) \subset G(\Q_p)$ is a $G(\Q_p)$-valued crystalline representation as in \S \ref{para:rational crys rep}.
 
 Fix a $\cG(\Z_p)$-valued crystalline representation $\rho$. From it we obtain a tautological functor\index[Oct2026-sorted-notation]{0001 one 1 06 25 sideset subscript rho 1@$\sideset{_\rho}{}{\idfunc}$} $$\idrho: \Rep_{\Z_p} \cG \To \Rep_{\Gamma_K}^{\cris \circ}, $$ sending each $L$ to the $\Gamma_K$-stable lattice $L$ in the crystalline representation $L\otimes \Q_p$.
 
 We shall need a variant of this definition. 
 Let $\Gamma_K$ act on $G(\Q_p)/\cG(\ZZ_p)$ by $$\gamma(g \cG(\ZZ_p)):= \rho (\gamma)  g \cG(\ZZ_p), \quad  \forall \gamma \in \Gamma_K, g\in G(\Q_p). $$ Since $\rho(\Gamma_K) \subset \cG(\ZZ_p)$, the coset of $1$ in $G(\QQ_p)/\cG(\ZZ_p)$ is fixed by $\Gamma_K$. Now consider an arbitrary element  $\lambda \in G(\Q_p)/\cG(\Z_p)$  fixed by $\Gamma_K.$ Choose a representative $\lambda\in G(\QQ_p)$. 
 Then $\lambda^{-1}\rho(\Gamma_K)\lambda \subset \cG(\Z_p)$ and we set \index[Oct2026-sorted-notation]{zz17 rho 1 00 00 rho superscript lambda@$\rho^{\lambda}$}
 $$ \rho^{\lambda}: \Gamma_K \longrightarrow \cG(\Z_p), \quad \gamma \longmapsto \lambda^{-1}\rho(\gamma)\lambda. $$  		
 We define\index[Oct2026-sorted-notation]{0001 one 1 06 25 sideset subscript rho 1 superscript lambda@$\sideset{_\rho}{}{\idfunc}^{\lambda}$}  $$\idrholambda : \Rep_{\Z_p} \cG \To \Rep_{\Gamma_K}^{\cris \circ}$$ to be the functor sending each $L$ to the 
 $\Gamma_K$-stable lattice $\lambda \cdot L$ inside $L\otimes \Q_p$. 
 Then there is a natural isomorphism $\idlambdarho \isom \idrholambda$ given by the multiplication-by-$\lambda$ maps $L \isom  \lambda\cdot L.$
 
 \subsection{}\label{para:crystallinerepII}
 We define $\otimes$-functors \index[Oct2026-sorted-notation]{zz24 omega 1 00 00 omega subscript rho superscript lambda@$\omega_{\rho}^{(\lambda)}$}
 $$ \omega_{\rho}:   \Rep_{\Z_p} \cG \xrightarrow{\gM \circ \idrho} \Mod_{/\gS}^{\varphi} \to \Modfp_{\gS}, \quad L \longmapsto \gM(L) $$
 and \index[Oct2026-sorted-notation]{zz24 omega 1 00 00 omega subscript rho 0 superscript lambda@$\omega_{\rho,0}^{(\lambda)}$}
 $$ \omega_{\rho, 0}:  \Rep_{\Z_p} \cG  \xrightarrow{M_{\cris} \circ  \idrho } \Isoc_{K_0}^\circ \to   \Modfp_{W} , \quad L \longmapsto M_{\cris}(L).$$
 Here the last arrows in both cases are the natural forgetful functors.
 Clearly $$\omega_{\rho, 0} = (\omega_{\rho})_W .$$
 We also write\index[Oct2026-sorted-notation]{zz20 upsilon 0 00 00 upsilon subscript rho superscript lambda@$\Upsilon_{\rho}^{(\lambda)}$} $\tomega_{\rho}$ for the $\otimes$-functor
 $$ M_{\cris} \circ  \idrho  :   \Rep_{\Z_p} \cG  \To \Isoc_{K_0}^\circ .$$
 
 Note that the definition of $\omega_{\rho}$ depends on the choice of a uniformizer in $K$, but the definitions of 
 $\tomega_{\rho}$ and $\omega_{\rho, 0}$ are independent of such a choice (up to canonical $\otimes$-isomorphism), 
 by property (iv) in \S \ref{gSmodules}. We will write $\omega^\lambda_{\rho}, \omega^\lambda_{\rho, 0}$ and 
 $\tomega^\lambda_{\rho}$ when we apply these constructions to $\idlambdarho \cong  \idrholambda.$
 
 \subsection{}\label{para:crystallinerepIII}
 As in \S \ref{para:setting for fiber functor}, we denote by $P_{\omega_{\rho}}$ the $\gS$-scheme $\SIsom^{\otimes}(\idfunc_{\gS}, \omega_{\rho})$. Since we have  $$\underline{\Aut}^{\otimes}(\idfunc_{\gS}) \cong \cG_{\gS}$$ by the reconstruction theorem, we know that $P_{\omega_{\rho}}$ is a pseudo-torsor under $\cG_{\gS}$ (i.e., for each $\gS$-scheme $S$, the set $P_{\omega_{\rho}}(S)$ is either empty or a principal homogeneous space under $\cG(S)$). However, in general 
 $\omega_{\rho}$ may not be a fiber functor (as it may not be exact), and $P_{\omega_{\rho}}$ may not be a 
 $\cG_{\gS}$-torsor.  In the following, we denote by $U$ the complement of the closed point in $\Spec \gS$.  
 \begin{lem}\label{lem:torsor over U}
 	The restriction $P_{\omega_{\rho}}|_U$ is a $\cG_U$-torsor over $U$.
 \end{lem}
 \begin{proof} We write $\omega$ for $\omega_{\rho}$.
 	By \cite[Thm.~3.3.2]{KisinPappas}, the functor $\omega_U: L \mapsto \omega(L) |_{U}$ is an exact faithful $\otimes$-functor from $\Rep_{\Z_p} \cG$ to the category of vector bundles on $U$. As in \S \ref{para:setting for fiber functor}, we may also regard $\omega_U$ as a fiber-wise faithful exact functor between the fibered categories $\mathbf{Rep}~\cG$ and $\mathbf{Bun}_U$, where the fibers of $\mathbf{Bun}_U$ over $\Spec \Z_p$ and $\Spec \Q_p$ are respectively the categories of vector bundles on $U$ and on $U \times_{\Spec \Z_p} \Spec \Q_p$. Since $P_{\omega}|_U$ is identified with $P_{\omega_U}$, we know that it is a $\cG_U$-torsor by \cite[Thm.~4.8]{broshi2013}.
 \end{proof}
 \subsection{}\label{para:KL}
 Let $\cG$ be a smooth affine group scheme over $\Z_p$ with connected fibers.  We say that\index{S-purity@$\mathfrak S$-purity} \emph{$\cG$ satisfies $\gS$-purity}, if every $\cG$-torsor over $U$ (the complement of the closed point in $\Spec \gS$) extends to $\Spec \gS$. This is equivalent to the condition that all $\cG$-torsors on $U$ are trivial. Indeed, since $\cG$ is smooth with connected fibers and since the residue field of the closed point in $\Spec \gS$ is the finite field $k$, by Lang's theorem we know that all $\cG$-torsors on $\Spec \gS$ are trivial.

 By descent, $\cG$ satisfies $\gS$-purity if and only if all $\cG$-torsors on the complement of the closed point of $\Spec W(\Fpbar)\lps u \rps$ extend to $\Spec W(\Fpbar)\lps u \rps$, cf.~Step 4 in the proof of \cite[Prop.~1.4.3]{KisinPappas}. Thus $\gS$-purity is intrinsic to the group $\cG$, and is independent of the finite extension $k/\FF_p$ which determines $\gS = W(k)\lps u \rps$.

 It has been proved by Ansch\"utz \cite[Cor.~1.2]{Anschutz} that all parahoric group schemes $\cG$ over $\Z_p$ satisfy $\gS$-purity, generalizing earlier results in \cite{CS79} and \cite{KisinPappas}. We will make use of this result mainly when $\cG$ is either a connected reductive group scheme over $\Z_p$, or the connected N\'eron model of a torus.\footnote{For connected reductive group schemes this result already follows from  \cite{CS79}, as explained in Step 5 in the proof of \cite[Prop.~1.3.4]{kisin2010integral}.}  
 
 \begin{lem}\label{keylemma} Suppose that $\cG$ satisfies $\gS$-purity. Then $P_{\omega_{\rho}}$ is a $\cG$-torsor. 
 	In particular, $\omega_{\rho}$  is (non-canonically) $\otimes$-isomorphic to $\idfunc_{\gS},$  the functors 
 	$\omega_{\rho}$ and $\omega_{\rho,0}$ are fiber functors, and $\tomega_{\rho}$ is an object in $\GIsoc_{K_0}^\circ$.
 \end{lem}
 \begin{proof}  We write $\omega$ for $\omega_{\rho}$. By Lemma \ref{lem:torsor over U} and the discussion in \S \ref{para:KL},  $P_{\omega}|_U$ is a trivial $\cG_U$-torsor. Fix a section of $P_{\omega}|_U$ over $U$. Then for each $L$ in $\Rep_{\Z_p}\cG$, we obtain an isomorphism $\iota_L: \idfunc_{\gS} (L) |_U \isom \omega(L)|_U$ between vector bundles on $U$, which is functorial in $L$ and compatible with tensor products. Since $\gS$ is a noetherian normal domain and since the closed point in it has codimension $2$, the isomorphism $\iota_L$ extends uniquely to an isomorphism $\tilde \iota_L: \idfunc_{\gS}(L) \isom \omega(L)$ between finite projective $\gS$-modules. By the uniqueness, we know that $\tilde \iota _L$ is functorial in $L$ and compatible with tensor products. Thus we have constructed a $\otimes$-isomorphism $\idfunc_{\gS} \isom \omega$ between $\otimes$-functors.
 	
 	Since $\omega_{\rho}$ is $\otimes$-isomorphic to $\idfunc_{\gS}$, it is a fiber functor. Since $\omega_{\rho, 0} \cong (\omega_{\rho})_W$, it is also a fiber functor. It follows that $\tomega_{\rho}$ is in $\GIsoc_{K_0}^\circ$.
 \end{proof}
 \begin{rem}\label{rem:top crit}
 	Lemma \ref{keylemma} together with (i) in \S \ref{gSmodules} has the following consequence. Let $\rho: \Gamma_K \to G(\Q_p)$ be a crystalline representation. In the notation of \S \ref{para:rational crys rep}, a sufficient condition for $v \circ \mathcal I_{\rho}$ to be $\otimes$-isomorphic to $\idfunc_{K_0}$ (or equivalently, for $\mathcal I_\rho$ to be $\otimes$-isomorphic to some $\mathcal I_{\mu,\delta}$) is that $\rho$ factors through a subgroup of the form $\cG(\Z_p)$, for $\cG$ satisfying $\gS$-purity. 
 \end{rem}
 \subsection{}\label{para:delta and mu} For the rest of \S \ref{subsec:crystalline lattices} we assume that $\cG$ satisfies $\gS$-purity, and write $G$ for $\cG_{\QQ_p}$ (assumed to be reductive). We also assume that $K_0$ is sufficiently large such that the images of $\Gamma_{K_0}$ and $\Gamma_{p,0}$ (the inertia subgroup of $\Gamma_{\QQ_p}$) in $\Aut(\pi_1(G))$ are the same. Fix a $\otimes$-isomorphism $\eta: \idfunc_{\gS} \isom \omega_{\rho}.$ Then for $L$ in $\Rep_{\Z_p} \cG,$ the isomorphism $ L \otimes_{\Z_p} \gS \isom \omega_\rho(L)$ induced by $\eta$ carries the Frobenius $\varphi^*(\omega_\rho(L))[1/E] \isom \omega_\rho(L)[1/E] $ to an isomorphism $\varphi^*(L \otimes_{\Z_p} \gS)[1/E] \isom L \otimes_{\Z_p} \gS[1/E] $, which is of the form $\delta_{\gS, L}  \otimes \varphi$ for some $\delta_{\gS, L} \in \GL(L) (\gS[1/E])$. By the reconstruction theorem, the elements $\delta_{\gS, L}$ for all $L$ come from a common, unique element $\delta_{\gS} \in \cG (\gS[1/E])$ (cf.~the similar argument in \S \ref{para:delta assoc to isoc}). If we change the choice of $\eta$, then $\delta_{\gS}$ gets $\varphi$-conjugated by an element of $\cG(\gS)$.
 
 Let $\delta \in \cG(K_0)$ be the image of $\delta_{\gS}$ under the specialization $u \mapsto 0.$ Then the $\cG(\O_{K_0})$-$\sigma$-conjugacy class of $\delta$ is independent of the choice of $\eta$, and it coincides with $\inv(\tomega_{\rho})$ defined in \S \ref{para:delta assoc to isoc}. By applying the map $D_{\cris}^{G}$ in \S \ref{para:rational crys rep} to $\rho$ we obtain an element of $\RZ_G/{\sim}$, and in particular a cocharacter $\mu$ of $G_{\Qpbar}$, well defined up to $G(\Qpbar)$-conjugacy (cf.~\S \ref{subsubsec:RZ}).  As usual, write $w_G: G(\LL) \to \pi_1(G)_{\Gamma_{p,0}}$ for the Kottwitz homomorphism associated with the $p$-adic valuation on $\LL$.
 
 \begin{prop}\label{mark:5.1} With the above notation, $w_G(\delta)$ is equal to the image of $\mu$. 
 \end{prop}

 The proof of this proposition requires the following lemma, which we shall invoke again for Proposition \ref{prop:compbuildingpoints}.

 \begin{lem}\label{lem:Bruhatorder}   Let  $\mathcal{G}$  be a smooth affine group scheme over $\ZZ_p$ with connected fibers and reductive generic fiber $G$. Assume that $\Gamma_{K_0}$ and $\Gamma_{p,0}$ have the same image in $\Aut(\pi_1(G))$. Let $f$ be a non-zero irreducible element of $\gS$ such that $f(0) \in W$ has $p$-adic valuation $1$. Let $v_f$ be the $f$-adic valuation
 	on $\Frac \gS$, and let $$   \kappa_{G_{K_0}}^{ v_f} :  G(\Frac \gS) \To \pi_1(G)_{\Gamma_{K_0}} = \pi_1(G)_{\Gamma_{p,0}}$$ be the associated Kottwitz homomorphism as in Definition \ref{defn:Kott map}. For any $g \in \cG(\gS[1/f])$ with image $g_0 \in G(K_0)$ under the specialization $u \mapsto 0$, we have $\kappa_{G_{K_0}}^{ v_f}  (g) = w_G(g_0)$. 
 \end{lem}
 \begin{proof} 	Let $v_f$ and $v_p$ be the discrete valuations on $\gS[1/fp]$ attached to the prime ideals $f\gS$ and $p\gS$ of $\gS$. Let $v_0$ be the discrete valuation on $\gS[1/fp]$ given by  $$ \gS[1/fp] \xrightarrow{u\mapsto 0} K_0 \overset{v_p}{\To} \ZZ \cup \set{\infty}.$$  
 	Any unit $w \in \gS[1/fp]^\times$ has the form $w = f^ip^jy$ with $i,j \in \ZZ,$ and
 	$y \in \gS^\times.$ We have $$v_0 (w) = i+ j = \begin{cases}
 		v_f (w), & \text{if } f \in p\gS^\times ,\\
 		v_f (w) + v_p(w), & \text{otherwise}.
 	\end{cases} $$ By functoriality and Proposition \ref{prop:sumkottmaps}, for any $g \in G(\gS[1/fp])$ we have
 	$$ w_G(g_0) = \kappa_{G_{K_0}}^{v_0}(g) = \begin{cases}  \kappa_{G_{K_0}}^{v_f}(g), & \text{if } f \in p \gS^\times , \\ 
 		\kappa_{G_{K_0}}^{v_f}(g) + \kappa_{G_{K_0}}^{ v_p}(g), & \text{otherwise}.
 	\end{cases}  $$ It remains to show that when $f \notin p\gS^\times$ we have $\kappa^{v_p}_{G_{K_0}}(g) = 0$ for all $g \in \cG(\gS[1/f])$. This follows from Corollary \ref{cor:kottparahoric} since $g \in \cG(\gS_{(p)})$. 
 \end{proof}
 
 \begin{proof}[Proof of Proposition \ref{mark:5.1}]   For each $L \in \Rep_{\Z_p} \cG$ and $V \in \Rep_{\Q_p} G$, we understand that $\Gamma_K$ acts on $L$ and $V$ via $\rho$. By Remark \ref{rem:top crit}, the fiber functor $\Rep_{\Q_p} G \to \Modfp_K, V\mapsto D_{\cris}(V) \otimes_{K_0} K \cong D_{\dR}(V)$ is $\otimes$-isomorphic to $V \mapsto V\otimes_{\QQ_p} K$. Fix such a $\otimes$-isomorphism $j: D_{\dR}(V) \isom V\otimes_{\QQ_p} K$. By the definition of $\mu$, we may assume that it is defined over $K$ and that $j$ is a filtered isomorphism where the filtration on the right hand side is defined by $\mu$.    
 	
 	Consider $L$ in $\Rep_{\Z_p} \cG$.  By (\ref{eq:filtered iso general 2}) we have a canonical filtered isomorphism 
 	\begin{align*}
 		c:	\varphi^* \gM(L) \otimes_{\gS} \widehat \gS_{(E)} \isom  D_{\dR}(L\otimes \Q_p)\otimes_K \widehat \gS_{(E)}. 
 	\end{align*}
 	Then $(j \otimes_K   \widehat \gS_{(E)} ) \circ c \circ (\varphi^*\eta \otimes_{\gS} \widehat \gS_{(E)})$ is a  $\otimes$-automorphism of the functor $L \mapsto L \otimes_{\Z_p} \widehat \gS_{(E)} $, and is therefore given by an element $a \in G(\widehat \gS_{(E)})$. 
 	
 	Since the isomorphism $$(j \otimes_K   \widehat \gS_{(E)} ) \circ c :  \varphi^* \gM(L) \otimes_{\gS} \widehat \gS_{(E)} \isom L \otimes_{\Z_p} \widehat \gS_{(E)}   = L[1/p]\otimes _{\Q_p} \widehat \gS_{(E)} $$ is filtered (where the right hand side has the tensor product filtration given by the $\mu$-filtration on $L[1/p]$ and the $E$-adic filtration on $\widehat \gS _{(E)}$), it maps $\Fil^r$ of the left hand side onto $$\mu(E^{-1}) E^r (L \otimes_{\Z_p} \widehat \gS_{(E)} ) \cap (L \otimes_{\Z_p} \widehat \gS_{(E)} ), $$ where $\mu(E^{-1})$ acts on $L \otimes_{\Z_p} \widehat \gS_{(E)}[1/E]$. By the definition of the filtration on $\varphi^* \gM(L)$ and the definition of $\delta_{\gS}$, the isomorphism $$(\varphi^*\eta \otimes_{\gS} \widehat \gS_{(E)}) :  L \otimes_{\Z_p} \widehat \gS_{(E)} \isom   	\varphi^* \gM(L) \otimes_{\gS} \widehat \gS_{(E)} $$ pulls back $\Fil^r$ of the right hand side to a submodule of the left hand side with the property that $\delta_{\gS}$ times it is contained in $E^r (L \otimes_{\Z_p} \widehat \gS_{(E)})$. (Here $\delta_{\gS}$ acts on $L \otimes_{\Z_p} \widehat \gS_{(E)}[1/E]$.) Combining these two facts, we see that 
 	$$ \delta_{\gS} a^{-1} \mu(E^{-1}) E^r(L \otimes_{\Z_p} \widehat \gS_{(E)}) \subset E^r ( L \otimes_{\Z_p} \widehat \gS_{(E)})   $$ for all sufficiently large $r$. Canceling $E^r$ we obtain  $$ \delta_{\gS} a^{-1} \mu(E^{-1}) (L \otimes_{\Z_p} \widehat \gS_{(E)}) \subset  L \otimes_{\Z_p} \widehat \gS_{(E)} .  $$
 	This is true for all $L$, so the element $\delta_{\gS} a^{-1} \mu(E^{-1}) \in G(\widehat \gS_{(E)}[1/E])$ lies in $ G(\widehat \gS_{(E)})$. Hence we have $ \delta_{\gS} \in G(\widehat \gS_{(E)}) \mu(E) G(\widehat \gS_{(E)}).  $ 
 	
 	By Corollary \ref{cor:kottCartan} applied to $F= \widehat \gS_{(E)}[1/E]$ and $\O_F = \widehat \gS_{(E)}$, the above relation implies that the element   $\kappa_{G_{K_0}}^{v_E} (\delta_{\gS}) \in \pi_1(G)_{\Gamma_{p,0}}$ defined as in Lemma \ref{lem:Bruhatorder} is equal to the image of $\mu$. We conclude the proof by that lemma.  
 \end{proof}
 
 \subsection{}\label{fiberfunctoranew} 
 The natural base change functor 
 $\Isoc_{K_0}^{\circ} \to \Isoc_{\Q_p^{\ur}}^{\circ}, D \mapsto D\otimes_{\oo_{K_0}} \ZZ_p^{\ur}$ induces a base change functor 
 $\GIsoc_{K_0}^{\circ} \to \GIsoc_{\Q_p^{\ur}}^{\circ} $. For a crystalline representation $\rho : \Gamma_K \to \cG(\Z_p)$,  we denote by\index[Oct2026-sorted-notation]{zz20 upsilon 0 00 00 upsilon subscript rho ur@$\Upsilon_{\rho,\ur}$} $\tomega_{\rho,\ur} \in \GIsoc_{\Q_p^{\ur}}^{\circ}$ the image of  $\tomega_{\rho} \in \GIsoc_{K_0}^{\circ}$ under the  latter functor. 
 
 If $K'/K$ is a finite extension and if $\rho'$ is the restriction of $\rho$ to $\Gamma_{K'}$, then the same construction gives rise to $ \tomega_{\rho', \ur} \in \GIsoc_{\Q_p^{\ur}}^{\circ}.$ By property (iv) in \S \ref{gSmodules}, we know that $\tomega_{\rho', \ur}$ is canonically $\otimes$-isomorphic with $\tomega_{\rho, \ur}$. 
 Hence the definition of $\tomega_{\rho,\ur} \in \GIsoc_{\Q_p^{\ur}}^{\circ}$ depends only on the germ of $\rho$ at 
 $1\in \Gamma_K.$
 
 Now let $[\rho] \in \Crys_{G}$ (with notation as in \S \ref{para:rational crys rep}). By continuity, we can find a finite extension $K/\QQ_p$  and  a representative $\rho \in \Crys_{G}(K)$ of $[\rho]$ such that $\rho(\Gamma_{K}) \subset \cG(\Z_p)$. We then define\index[Oct2026-sorted-notation]{zz20 upsilon 0 00 00 upsilon subscript rho superscript lambda leftbracket rightbracket@$\Upsilon_{[\rho]}^{(\lambda)}$}
 $$ \tomega_{[\rho]} : = \tomega_{\rho, \ur} \in \GIsoc_{\Q_p^{\ur}}^{\circ},$$ which is independent of choices. 
 For $\lambda \in G(\Q_p)/\cG(\Z_p)$ we write 	 $\tomega_{[\rho]}^\lambda$ for $\tomega_{[\rho^{\lambda}]}.$ (As we can arbitrarily enlarge $K$, we may assume that $\lambda$ is fixed by $\Gamma_K$ as in \S \ref{para:crystallinerep}.)

 \begin{cor}[cf.~Proposition \ref{introprop:p-adic Hodge} in the Introduction]\label{cor:comparing delta and mu}
 	Let $[\rho] \in \Crys_G$. 
 	Let $\delta_{\ur} \in G(\Qpur)$ be a representative of the  $\cG(\Z_p^{\ur})$-$\sigma$-conjugacy class $\inv(\tomega_{[\rho]})$ (\S \ref{para:delta assoc to isoc}), and let $\mu$ be the cocharacter in $D_{\cris}^G([\rho])$ (up to $G(\Qpbar)$-conjugacy). Then $w_G(\delta_{\ur})\in \pi_1(G)_{\Gamma_{p,0}}$ is equal to the image of $\mu$.  
 \end{cor}
 \begin{proof} Pick a representative $\rho: \Gamma_K \to \cG(\Z_p)$ of $[\rho]$, with $K$ large enough such that $\Gamma_{K_0}$ and $\Gamma_{p,0}$ have the same image in $\Aut(\pi_1(G))$. Then $\delta$ as defined in \S \ref{para:delta and mu} represents $\inv(\tomega_{[\rho]})$. Since $w_G$ is constant on $\inv(\tomega_{[\rho]})$ (Corollary \ref{cor:kottparahoric}), the result follows immediately from Proposition \ref{mark:5.1}. 
 \end{proof}
 \begin{rem}
 	The above result strengthens the assertion that the image of $D_{\cris}^G$ is contained in $\mathrm{RZ}^{\neu}_G /{\sim}$ in Proposition \ref{prop:adm RZ pair}. 
 \end{rem}
 \begin{cor}\label{cor:inv of T-valued rep} 
 	Let $T$ be a cuspidal torus over $\Q$, and $\cT^\circ$ the connected N\'eron model of $T_{\Q_p}$ over $\Z_p.$ 
 	Let $[\rho_T]$ be an element of $\Mot_T,$ and let $\mu = \motinv_T([\rho_T]) \in X_*(T)$, where $\motinv_T$ is the bijection in Proposition \ref{prop:classifyingmotrep}. Let $\delta_T \in T(\Qpur) $ be a representative of  $\inv(\tomega_{[\rho_T]})$. Then $\delta_T $ lies in $T(\Qpur)^{\mot}$, and  $w_{T_{\Q_p}} (\delta_T) \in X_*(T)_{\Gamma_{p,0}}$ is equal to the image of $\mu$.
 \end{cor}
 \begin{proof} Since $[\rho_T]$ lies in $\Mot_T$, we have $\delta_T \in T(\Qpur)^{\mot}$. The claim about $w_{T_{\Q_p}} (\delta_T)$ follows from Corollary \ref{cor:comparing delta and mu} applied to $\cG = \cT^{\circ}$. Here we have used that $\cT^\circ$ satisfies $\gS$-purity; see \S \ref{para:KL}.
 \end{proof}
 \begin{rem}
 	By Lemma \ref{lem:motelts}, the two properties satisfied by $\delta_T$ asserted in Corollary \ref{cor:inv of T-valued rep} uniquely characterize the $\cT^\circ(\Z_p^{\ur})$-$\sigma$-conjugation orbit of $\delta_T$.
 \end{rem}
 \subsection{}\label{defnYanew} Let $[\rho] \in \Crys_G$ and  $\lambda \in G(\Q_p)/\cG(\Z_p).$ 
 There is a tautological $\otimes$-isomorphism $\tomega_{[\rho]}^{\lambda}[1/p] \isom \tomega_{[\rho]}[1/p]$, inducing a tautological isomorphism $Y(\tomega_{[\rho]}^{\lambda}) \isom Y(\tomega_{[\rho]})$. Under this isomorphism, the relative position of the $\cG(\Z_p^{\ur})$-torsors $Y(\tomega_{[\rho]}^{\lambda})^{\circ}$ and $Y(\tomega_{[\rho]})^{\circ}$  is measured by an element $$\lambda_0 \in \cG(\Z_p^{\ur}) \backslash G(\Qpur)/\cG(\Z_p^{\ur}). $$ Namely, after picking $y_0 ' \in Y(\tomega_{[\rho]}^{\lambda})^{\circ}$ and $y_0 \in Y(\tomega_{[\rho]})^{\circ} $ we define $\lambda_0 \in G(\Qpur)$ by  $y_0' \mapsto y_0 \lambda_0$. By Corollary \ref{cor:kottparahoric}, we obtain well-defined elements $w_G(\lambda)$ and $ w_G(\lambda_0) $ of $\pi_1(G)_{\Gamma_{p,0}}$.
 \begin{prop}\label{prop:compbuildingpoints} We have $w_G(\lambda) = w_G(\lambda_0)$.
 \end{prop}
 \begin{proof} Pick a representative $\rho: \Gamma_K \to \cG(\Z_p)$ of $[\rho]$. Up to enlarging $K$,  we may assume that the $\Gamma_K$-action on $G(\Q_p)/\cG(\Z_p)$ induced by $\rho$ fixes $\lambda$, and that  $\Gamma_{K_0}$ and $\Gamma_{p,0}$ have the same image in $\Aut(\pi_1(G))$.

 	As in \S \ref{para:crystallinerepII}, we have $\otimes$-functors $\omega_{\rho}$ and $\omega^{\lambda}_{\rho}: \Rep_{\Z_p} \cG \to \Modfp_{\gS}$. For brevity we denote them by $\omega$  and $\omega^{\lambda}$. By Lemma \ref{keylemma}, $P_{\omega} (\gS)$ and $P_{\omega^{\lambda}} (\gS)$ are non-empty. In particular, the $G(\gS[1/p])$-torsor $Y(\omega)$ (resp.~$Y(\omega^{\lambda})$) contains the  $\cG(\gS)$-torsor $Y(\omega)^{\circ}$ (resp.~$Y(\omega^{\lambda})^{\circ}$) (see \S \ref{para:defnY}). The tautological isomorphism $\omega^{\lambda}[1/p] \isom \omega[1/p]$ induces an isomorphism $Y(\omega^{\lambda}) \isom  Y(\omega)$, under which the relative position of $Y(\omega^{\lambda})^{\circ}$ and $Y(\omega)^{\circ}$ is measured by an element $$\lambda_{\gS} \in \cG(\gS) \backslash G(\gS[1/p])/\cG(\gS)$$ (in the same way as in \S \ref{defnYanew}). Clearly  $\lambda_0$ is equal to the image of $\lambda_{\gS}$ under $$\cG(\gS)\backslash G(\gS[1/p]) /\cG(\gS) \to \cG(\O_{K_0})\backslash G(K_0)  / \cG(\O_{K_0})\to \cG(\Z_p^{\ur}) \backslash  G(\Qpur)/\cG(\Z_p^{\ur}), $$ where the first map is induced by the specialization $u\mapsto 0$.
 	
 	Now  write $C$ for  $\hat \O_{\cE^{\ur}}$. By (iii) in \S \ref{gSmodules}, the base change $\omega^{\lambda}_C$ of $\omega^{\lambda}$ to $C$ is canonically $\otimes$-isomorphic to the functor $$\idfunc^{\lambda}_{C}: \Rep_{\Z_p } \cG \To \Modfp_C, \quad  L \longmapsto (\lambda \cdot L)\otimes_{\Z_p} C. $$ Similarly, $\omega_C$ is canonically $\otimes$-isomorphic to $\idfunc_C$. Moreover, the canonical $\otimes$-isomorphisms $\idfunc_{C}^{\lambda} \cong \omega^{\lambda}_C$ and $\idfunc_C \cong \omega_C$ are compatible with the tautological isomorphisms $\idfunc_C^{\lambda}[1/p] \cong \idfunc_C [1/p]$ and $\omega^{\lambda}_C[1/p] \cong \omega_C[1/p]$. It follows that the image of $\lambda_{\gS}$ in $\cG(C) \backslash G(C[1/p])/\cG(C)$ is equal to the image of $\lambda$. 
 	
 	Write $v_p$ for the $p$-adic valuation on $\Frac(\gS)$. By the construction of $C = \widehat \O_{\mathcal E^{\ur}} $  in \cite[\S 2.1.1]{kisin2006crystalline} (cf.~\cite[\S A.3.2.1]{FontaineGroth}), we know that $C$ is a complete DVR with uniformizer $p$. In particular,
 	$C[1/p]$  is a field, and it has a discrete valuation $w_p$ extending $v_p$ with valuation ring $C$.  
 	By Corollary \ref{cor:kottparahoric} and the functoriality of the Kottwitz homomorphism,  
 	we have a commutative diagram 
 	$$ \xymatrix{ G(\Q_p)/\cG(\Z_p) \ar[r] \ar[d]^{w_G}& \cG(\gS)\backslash G(\gS[1/p]) /\cG(\gS)  \ar[r] \ar[d]^{\kappa_{G_{K_0}}^{v_p}}& \cG(C) \backslash G(C[1/p])/\cG(C)  \ar[d]^{\kappa_{G_{K_0}}^{w_p}} \\ \pi_1(G)_{\Gamma_{p,0}} \ar@{=}[r] &  \pi_1(G)_{\Gamma_{K_0}} \ar@{=}[r] &  \pi_1(G)_{\Gamma_{K_0}} }$$  Since $\lambda_{\gS}$ and $\lambda$ have the same image in the top right corner, we have 
 	$ w_G(\lambda) = \kappa_{G_{K_0}}^{v_p} (\lambda_{\gS}). $
 	Now by Lemma \ref{lem:Bruhatorder} applied to $f= p$, we have $\kappa_{G_{K_0}}^{v_p} (\lambda_{\gS})  =w_G(\lambda_0)$ since $\lambda_0$ is the image of $\lambda_{\gS}$ under the specialization $u \mapsto 0$. Therefore $w_G(\lambda) = w_G(\lambda_0)$ as desired.
 \end{proof} 
 \section{Compatibility with a maximal torus}\label{subsec:Refineintst}
 In this section we establish compatibility results for a reductive group and a maximal torus in it. 
 \subsection{}\label{para:rho_T} We continue to use the notation of \S \ref{subsec:crystalline lattices}, and we assume that 
 $\cG$ is a parahoric group scheme over $\Z_p.$  Let $T \subset G = \cG_{\Q_p}$ be a maximal torus. 
 Let $\cT^\circ$ (resp.~$\cG^{\ab}$) be the connected N\'eron model of $T$ (resp.~$G^{\ab}$) over $\Z_p.$ 
 As recalled in \S \ref{para:KL}, the $\Z_p$-group schemes $\cG, \cT^\circ, \cG^{\ab}$ satisfy $\gS$-purity. 
 
 Let $[\rho_T] \in \Crys_T,$ and denote by $[\rho]$ and $[\rho^{\ab}]$ the images of $[\rho_T]$ in $\Crys_G$ and  $\Crys_{G^{\ab}}$ respectively. 
 Using the constructions of \S \ref{fiberfunctoranew} and \S \ref{para:delta assoc to isoc} we obtain the $G(\Qpur)$-torsor $Y(\tomega_{[\rho]})$, the $T(\Qpur)$-torsor $Y(\tomega_{[\rho_T]})$, and the $G^{\ab}(\Qpur)$-torsor $Y(\tomega_{[\rho^{\ab}]})$. 
 The maps $T \rightarrow G \rightarrow G^{\ab}$ induce equivariant maps 
 $$ Y(\tomega_{[\rho_T]}) \To Y(\tomega_{[\rho]}) \To  Y(\tomega_{[\rho^{\ab}]}).$$

 Since the inclusion $T \hookrightarrow G$ does not necessarily extend to a map $\cT^\circ \rightarrow \cG$ over $\Z_p,$
 one cannot expect that the map $Y(\tomega_{[\rho_T]}) \to Y(\tomega_{[\rho]})$ sends integral points to integral points. Nevertheless, we have the following compatibility result. Denote by $\varrho$ the natural map $G_{\sconn} \to G$.\index[Oct2026-sorted-notation]{zz17 rho 1 00 00 varrho g subscript sconn to g@$\varrho: G_{\sconn}\to G$}
 
 \begin{prop}\label{prop:intstrcomp} The following statements hold. 
 	\begin{enumerate}
 		\item The	natural maps $ Y(\tomega_{[\rho_T]}) \to Y(\tomega_{[\rho^{\ab}]})$  and  $Y(\tomega_{[\rho]}) \to Y(\tomega_{[\rho^{\ab}]})$  send integral points to integral points.
 		\item  The image of $Y(\tomega_{[\rho_T]})^\circ$ in $Y(\tomega_{[\rho]})$ is contained in $ Y(\tomega_{[\rho]}) ^\circ \cdot \varrho (G_{\sconn}(\Qpur))$.
 	\end{enumerate}
 \end{prop}
 For the proof we need two lemmas from group theory, both of which should be well known to experts. 
 \begin{lem}\label{lem:ab surj}
 	The image of $\cG(\ZZ_p^{\ur})$ in $G^{\ab}(\QQ_p^{\ur})$ is $\cG^{\ab} (\ZZ_p^{\ur})$.  
 \end{lem}
 \begin{proof}
 	The map $G \to G^{\ab}$ extends to a map of $\Z_p$-group schemes $\cG\to \cG^{\ab}$, so the image of $\cG(\ZZ_p^{\ur})$ in $G^{\ab}(\QQ_p^{\ur})$ is contained in $\cG^{\ab} (\ZZ_p^{\ur})$. To show that $\cG(\ZZ_p^{\ur})$ surjects onto $\cG^{\ab} (\ZZ_p^{\ur})$, we have the following simple argument communicated to us by Rong Zhou. 
 	
 	Let $S$ be the centralizer of a maximal $\Qpur$-split torus in $G_{\Qpur}$ whose apartment over $\Qpur$ contains the facet defining the parahoric $\cG_{\Zpur}$. Then $S$ is a maximal torus since $G_{\Qpur}$ is automatically quasi-split over $\Qpur$. Let $\mathcal S^\circ$ denote the connected N\'eron model of $S$ over $\Zpur$. Then $\mathcal S^{\circ}(\Z_p^{\ur})$ is contained in $\cG(\Z_p^{\ur})$, so it suffices
 	to prove that $\mathcal S ^{\circ}(\Z_p^{\ur})$ maps onto $\cG^{\ab}(\Z_p^{\ur})$.  The kernel of $S \to G^{\ab}_{\Qpur}$ is  a maximal torus in $G_{\der,\Qpur}$. Thus it suffices to show the following general fact: Let $1 \to T_1 \to T_2 \to T_3 \to 1$ be a short exact sequence of tori over $\Qpur$. Let $\mathcal T_i^{\circ}$ denote the connected N\'eron model of $T_i$ over $\Zpur$. Then the map $\mathcal T_2^{\circ}(\Zpur) \to \mathcal T_3^{\circ}(\Zpur)$ is surjective. 
 	
 	We have the commutative diagram
 	$$ 
 	\xymatrix{
 		1 \ar[r]  & T_1(\Q_p^{\ur})\ar[r] \ar[d] & T_2(\Q_p^{\ur}) \ar[r] \ar[d] & T_3(\Q_p^{\ur}) \ar[r] \ar[d] & 1\\ 
 		& X_*(T_1)_{\Gamma_{p,0}}\ar[r] & X_*(T_2)_{\Gamma_{p,0}}\ar[r]  & X_*(T_3)_{\Gamma_{p,0}} \ar[r] & 1.}
 	$$ 
 	The second row is clearly exact, and the first row is exact by Steinberg's theorem. The vertical maps are the Kottwitz homomorphisms, and they are surjective by \cite[\S 7.7]{kottwitzisocrystal2}. The kernels of the three vertical maps are  $\mathcal T_i^{\circ}(\ZZ_p^{\ur})$, by \cite[Rmk.~2.2 (iii)]{rapoportguide}. By the Snake Lemma we have the surjectivity of $\mathcal T_2 ^{\circ}(\Z_p^{\ur}) \to \mathcal T_3^\circ (\Z_p^{\ur})$. 
 \end{proof}
 
 \begin{rem}\label{rem:ab surj}
 	If $G_{\der}= G_{\sconn}$, or if $G$ splits over an unramified extension, then we also know that $\cG(\ZZ_{p^n})$ surjects onto $\cG^{\ab}(\ZZ_{p^n})$ for all $n\geq 1$. Indeed, in either case the map $\pi_1(G_{\der})_{\Gamma_{p,0}} \to \pi_1(G) _{\Gamma_{p,0}}$ is injective, so we have $G_{\der}(\Qpur) \cap \cG(\Z_p^{\ur}) = \cG_{\der}(\Z_p^{\ur})$, where $\cG_{\der}$ is the parahoric model of $G_{\der}$ over $\Z_p$ determined by $\cG$ (via the identification of the buildings of $G_{\der}$ and $G$). The desired surjectivity follows from Lemma \ref{lem:ab surj} and the surjectivity of $\cG_{\der} (\breve \Z_p)\to \cG_{\der}(\breve \Z_p), g \mapsto g^{-1} \sigma^n(g)$ (by Greenberg's theorem applied to $\cG_{\der}$). 
 \end{rem}
 
 \begin{lem}\label{lem:ker of Kott}
 	The kernel of the Kottwitz homomorphism $w_G: G(\Qpur) \to \pi_1(G)_{\Gamma_{p,0}}$ is $\cG(\Z_p^{\ur}) \varrho(G_{\sconn}(\Qpur))$. 
 \end{lem}
 \begin{proof} Note that $\cG(\ZZ_p^{\ur}) \varrho(G_{\sconn}(\QQ^{\ur}_p))$ is  a subgroup of $G(\Qpur)$, since $ \varrho(G_{\sconn}(\QQ^{\ur}_p))$ is a normal subgroup of $G(\QQ^{\ur}_p)$. By \cite[Lem.~17]{HainesRapoport}, the kernel of the Kottwitz homomorphism $G(\breve \QQ_p) \to \pi_1(G)_{\Gamma_{p,0}}$ is $\cG(\breve \ZZ_p) \varrho (G_{\sconn}(\breve \Q_p))$. It remains to show the containment 
 	$$\varrho(G_{\sconn}(\breve \QQ_p)) \subset \cG(\breve \Z_p) \varrho(G_{\sconn} (\Qpur)).$$
 	Via the canonical identification of buildings, the parahoric model $\cG$ of $G$ determines a parahoric model $\cG_{\sconn}$ of $G_{\sconn}$, and we have $\varrho (\cG_{\sconn}(\Zpur)) \subset \cG(\Zpur)$.  By the density of $G_{\sconn}(\Qpur)$ in $G_{\sconn}(\LL)$, we have $G_{\sconn}(\breve \QQ_p) =\cG_{\sconn}(\breve \ZZ_p) G_{\sconn}(\QQ_p^{\ur})$. Hence 
 	$$ \varrho(G_{\sconn}(\breve \QQ_p)) \subset  \varrho (\cG_{\sconn}(\breve \ZZ_p) ) \varrho( G_{\sconn}(\QQ_p^{\ur}) ) \subset   \cG(\breve \Z_p) \varrho(G_{\sconn} (\Qpur)),  $$  as desired.   \end{proof}
 \begin{proof}[Proof of Proposition \ref{prop:intstrcomp}] Part (i) follows immediately from the functoriality of the constructions, and the fact that the $\Q_p$-homomorphisms $T \to G^{\ab}$ and $G\to G^{\ab}$ extend to $\Z_p$-homomorphisms $\cT^{\circ} \to \cG^{\ab}$ and $\cG\to \cG^{\ab}$ respectively. For part (ii), we first consider the special case where $G_{\der} = G_{\sconn}$. Then part (ii) follows from part (i) once we know that the map $Y(\tomega_{[\rho]})^\circ \to Y(\tomega_{[\rho^{\ab}]})^\circ$ (provided by part (i)) is surjective. This surjectivity follows from Lemma \ref{lem:ab surj}.
 	
 	We now prove part (ii) in general. By Corollary \ref{cor:kottparahoric} and the functoriality of the Kottwitz homomorphism, we know that $\cT^{\circ}(\Z_p^{\ur})$ lies in the kernel of the Kottwitz homomorphism $G(\Qpur) \to \pi_1(G)_{\Gamma_{p,0}}$. Thus by Lemma \ref{lem:ker of Kott} we have $$\cT^{\circ} (\Z_p^{\ur}) \subset \cG(\Z_p^{\ur}) \varrho(G_{\sconn}(\Qpur)) . $$ 
 	By this containment and as $Y(\tomega_{[\rho_T]})^{\circ}$ is a 
 	$\cT^{\circ}(\Z_p^{\ur})$-torsor and $Y(\tomega_{[\rho]})^{\circ}$ is a 
 	$\cG(\Z_p^{\ur})$-torsor, we see that to prove part (ii) it suffices to show that the image of $Y(\tomega_{[\rho_T]})^{\circ} \to Y(\tomega_{[\rho]})$ contains one element of $ Y(\tomega_{[\rho]})^{\circ} \cdot \varrho(G_{\sconn} (\Qpur))$.  
 	
 	To show this, let $1 \to Z \to G_1 \to G \to 1$ be a $z$-extension over $\QQ_p$ as in \cite[\S 1]{Kot82}, where $G_1$ has simply connected derived group and $Z$ is a central torus in $G_1$. Let $T_1$ be the inverse image of $T$ in $G_1$. By the identification $G_1^{\ad} \isom G^{\ad}$, the parahoric model $\cG$ of $G$ determines a parahoric model $\cG_1$ of $G_1$, and the map $G_1 \to G$ extends to a map of $\Z_p$-group schemes $\cG_1 \to \cG$. We claim that $[\rho_T] \in \Crys_T$ lifts to an element $[\rho_{T_1}] \in \Crys_{T_1}$. Given the claim, we apply the special case of part (ii) already proved above  to the maximal torus $T_1$ in $G_1$, $[\rho_{T_1}]$, and the element $[\rho_1]\in \Crys_{G_1}$ induced by $[\rho_{T_1}]$. Then we know that the natural map $Y(\tomega_{[\rho_{T_1}]}) \to Y(\tomega_{[\rho_1]})$ maps $Y(\tomega_{[\rho_{T_1}]})^{\circ}$ into $ Y(\tomega_{[\rho_1]})^{\circ} \cdot \varrho_1 (G_{1,\sconn} (\Qpur)),$ where $\varrho_1$ denotes the natural map $G_{1,\sconn} \to G_1$. By functoriality we have a commutative diagram 
 	$$ \xymatrix{ Y(\tomega_{[\rho_{T_1}]}) \ar[d] \ar[r] & Y(\tomega_{[\rho_1]}) \ar[d] \\  Y(\tomega_{[\rho_{T }]})  \ar[r] & Y(\tomega_{[\rho ]})   }  $$ and the vertical maps send $Y(\tomega_{[\rho_{T_1}]})^\circ$ into $Y(\tomega_{[\rho_{T}]})^\circ$ and send $Y(\tomega_{[\rho_1]})^\circ$ into $Y(\tomega_{[\rho]})^\circ$ respectively. The right vertical map is equivariant with respect to the homomorphism $G_1(\Qpur) \to G(\Qpur)$, and the latter maps $\varrho_1(G_{1,\sconn}(\Qpur))$ into $\varrho (G_{\sconn}(\Qpur))$. Thus the image of $Y(\tomega_{[\rho_T]})^\circ$ under the bottom horizontal map indeed contains an element of $Y(\tomega_{[\rho]})^{\circ} \cdot \varrho(G_{\sconn} (\Qpur))$, as desired. 
 	
 	It remains to prove the claim. Let $(\mu,[\delta]) \in X_*(T) \times (T(\Qpur)/1-\sigma)$ correspond to $[\rho_T]$ under the bijection in Corollary \ref{cor:classification of T-valued}. It suffices to lift this to a pair $(\mu_1,[\delta_1]) \in X_*(T_1) \times (T_1(\Qpur)/1-\sigma)$ such that $\mu_1$ and $[\delta_1]$ have the same image in $X_*(T_1)_{\Gamma_p}$. Let $(\mu_1',[\delta_1])$ be an arbitrary lift, which exists since $T_1(\Qpur) \to T(\Qpur)$ is surjective (by Steinberg's theorem). Then the difference of the images of $\mu_1'$ and $[\delta_1]$ in $X_*(T_1)_{\Gamma_p}$ lies in the kernel of   $X_*(T_1)_{\Gamma_p} \to X_*(T)_{\Gamma_p}$. Since we have an exact sequence 
 	$ X_*(Z)_{\Gamma_p} \to X_*(T_1)_{\Gamma_p} \to X_*(T)_{\Gamma_p}$, we can modify $\mu_1'$ by some element of $X_*(Z)$ to obtain the desired pair $(\mu_1, [\delta_1])$. 
 \end{proof}
 \cleardoublepage
 \chapter{Shimura varieties of Hodge type}\label{sec:Shimura Hodge type}

 	The goal of this chapter is to prove a strengthening of the Langlands--Rapoport--$\tau$ Conjecture for Shimura varieties of Hodge type.

 	In \S\ref{subsec:Hodge setting}--\S\ref{recall of defns} we mostly review the constructions from \cite{kisin2012modp}. In \S \ref{sec:conn comp}, we discuss the parametrization of connected components and introduce a refinement $\fkH(\phi)$ of the group $\cH(\phi)$ containing $\tau$ in the Langlands--Rapoport--$\tau$ Conjecture, which is closely related to the connected components. 
 	
 	In \S \ref{subsec:unif}, we study the uniformizations of (certain enlargements of) the isogeny classes, in preparation for the introduction and study of integral special points in \S \ref{subsec:intsppt}, where results from \S \ref{sec:crystalline} are applied in an essential way. 
 	
 	In \S \ref{subsec:markingsamic}, we introduce the notion of \emph{amicable pairs}, which plays a central role in our argument. The main technical tools to study them are provided in \S \ref{subsec:gauges}. In \S \ref{subsec:Galcohprop}, we then extend these technical results to general statements about amicable pairs using the CM lifting theorems on both sides. 
 	
 	Finally, in \S \ref{subsec:LR str}, we establish the main result of this chapter: the strengthened version of the Langlands--Rapoport--$\tau$ Conjecture in the case of Hodge type.	This strengthening incorporates additional structures (namely, connected components and certain twistings by central cohomology classes) as well as gives more refined information about the twisting elements $\tau$. This will ultimately enable us to derive the (unstrengthened) Langlands--Rapoport--$\tau$ Conjecture for Shimura varieties of abelian type in \S \ref{sec:abelian type}.

 \section{The abelian scheme and related structures on the Shimura variety}\label{subsec:Hodge setting}
 In this section, we review the integral models of Shimura varieties of Hodge type 
 along with their universal abelian schemes and related structures. The key ingredient 
 is the construction of an integral $F$-isocrystal with $\cG$-structure associated 
 with each $\ol\FF_p$-point.
 \subsection{}\label{para:Hodge embedding}
 Throughout \S \ref{sec:Shimura Hodge type} we keep the following assumptions. Let $(G,X, p ,\cG)$ be an unramified Shimura datum as in \S \ref{para:unramified SD}. Assume that $(G,X)$ is of \emph{Hodge type};\index{Shimura datum!Hodge type} recall that this means that there exists an embedding of Shimura data \begin{align}\label{eq:Siegel emb} (G,X) \hookrightarrow (\mathrm{GSp}(V_\Q), S^{\pm}),
 \end{align} where $V_\Q$ is a symplectic space over $\Q$, and $(\mathrm{GSp}(V_\Q), S^{\pm})$ is the corresponding Siegel Shimura datum.  Let $E \subset \CC$ be the reflex field of $(G,X)$, viewed as a subfield of $\Qbar$ via our fixed embedding $\Qbar \hookrightarrow \CC$. Let $\fkp$ be the prime of $E$ determined by the embedding $E \hookrightarrow \Qbar \hookrightarrow \Qpbar$. We write $K_p$ for $\cG(\Z_p)$.\index[Oct2026-sorted-notation]{k 0 00 00 k subscript p@$K_p$}

 Since $\cG$ is smooth over $\Z_p$, its Hopf algebra $\oo_{\cG}(\cG)$ injects into $\oo_{G_{\Q_p}} (G_{\Q_p})$. The intersection $\oo_\cG(\cG) \cap \oo_G(G)$  inside $\oo_{G_{\Q_p}} (G_{\Q_p})$ has the natural structure of a Hopf algebra over $\ZZ_{ (p) }$, and defines a $\ZZ_{ (p) }$-group scheme $G_{\ZZ_{(p)}}$.\index[Oct2026-sorted-notation]{g 0 00 00 g subscript z subscript p@$G_{\ZZ_{(p)}}$} Thus $G_{\ZZ_{(p)}}$ is the unique (up to unique isomorphism) reductive group scheme over $\ZZ_{(p)}$ whose generic fiber is identified with $G$ and whose base change to $\ZZ_p$ is identified with $\cG$. 
 
 As in \cite[\S 1.3.3]{kisin2012modp}, we may choose the Siegel embedding (\ref{eq:Siegel emb}) such that $G \hookrightarrow \mathrm{GSp} (V_{\Q})$ is induced by a closed embedding of $\Z_{(p)}$-group schemes $G_{\ZZ_{(p)}}\hookrightarrow \GL (V_{\Z_{(p)}})$, for some self-dual $\Z_{(p)}$-lattice $V_{\Z_{(p)}}$ in $V_{\Q}$.\footnote{This uses \cite[Lem.~2.3.1]{kisin2010integral} and Zarhin's trick. In the former result, there is an extra assumption on $G$ when $p=2$. However this extra assumption can be removed, as explained in the proof of \cite[Lem.~4.7]{KMP16}.} We fix such choices.   We write $\mathrm{GSp}$ for the $\Q$-algebraic group $\mathrm{GSp}(V_{\Q})$, and write $\spK_p$\index[Oct2026-sorted-notation]{k 0 06 00 k subscript p@$\mathbb K_p$} for the compact open subgroup $\mathrm{GSp}(V_{\Z_{(p)}}) (\Z_p)$ of $\mathrm{GSp}(\Q_p)$.
 Thus the embedding $G(\Q_p) \hookrightarrow \mathrm{GSp}(\Q_p)$ maps $K_p$ into $\spK_p$. 
 
 As in \cite[\S 1.3.6]{kisin2012modp}, we fix once and for all a collection of tensors\index[Oct2026-sorted-notation]{s 1 00 00 s subscript alpha subscript alpha in alpha@$(s_{\alpha})_{\alpha \in \bm \alpha}$} $(s_{\alpha})_{\alpha \in \bm \alpha} \subset V_{\Z_{(p)}}^{\otimes} $ such that the image of $G_{\ZZ_{(p)}} \hookrightarrow \GL (V_{\Z_{(p)}})$ is the scheme-theoretic stabilizer of these tensors. In the sequel, we shall assume that the collection $(s_{\alpha})_{\alpha \in \bm \alpha}$ is maximal in the sense that it consists of all elements of $V_{\ZZ_{(p)}}^{\otimes}$ stabilized by $G_{\ZZ_{(p)}}$.
 
 Let $V^*_{\Z_p}$ be the $\Z_p$-linear dual of $V_{\Z_p}$. We shall view it as a representation of $\cG$ over $\Z_p$, i.e., the contragredient of $V_{\Z_p}$. Similarly, we view the $\Q$-linear dual $V^*_\Q$ of $V_{\Q}$ as a representation of $G$ over $\Q$. We view $(s_{\alpha})_{\alpha \in \bm \alpha}$ also as tensors in $V^{*\otimes}_{\Z_p}$ or $V^{*\otimes}_{\Q}$.

 \subsection{}\label{para:universal av} For each compact open subgroup $K \subset G(\A_f)$ (resp.~$\spK \subset \mathrm{GSp} (\A_f)$), we write $\Sh_K$ (resp.~$\Sh_{\spK}$) for the Shimura variety $\Sh_K(G,X)$ (resp.~$\Sh_{\spK} (\mathrm{GSp}, S^{\pm})$). Below we recall the construction of the canonical smooth integral model $\Shh_{ K_p}$ of $\Sh_{K_p} = \varprojlim_{K^p}\Sh_{K_pK^p}$ stated in Theorem \ref{thm:integral model}.
 
 For a compact open subgroup $K^p \subset G(\A_f^p)$, we call it \emph{$\mathrm{GSp}$-neat} if its image in $\mathrm{GSp}(\A_f^p)$ is contained in a neat compact open subgroup. Clearly this holds for all sufficiently small $K^p$. For any $\mathrm{GSp}$-neat $K^p$, by \cite[Lem.~2.1.2]{kisin2010integral} there exists a neat compact open subgroup $\spK^p \subset \mathrm{GSp}(\A_f^p)$ such that the natural map $$\Sh_{K_pK^p} \To \Sh_{\spK_p \spK^p}\times _{\Spec\Q} \Spec E $$ is a closed immersion of $E$-schemes.
 
 Now $\Sh_{\spK_p \spK^p} $ has a canonical model $\Shh_{\spK_p \spK^p}$ over $\Z_{(p)}$ representing the usual Siegel moduli functor. We define $\Shh_{K_pK^p}$ to be the normalization of the closure of $\Sh_{K_pK^p}$ inside $\Shh_{\spK_p \spK^p} \times_{\Z_{(p)}} \oo_{E, (\fkp)}$. (It is shown in \cite{Xu} that taking the normalization is redundant.) By the main results of \cite{kisin2010integral} and \cite{KMP16}, $\Shh_{K_pK^p}$ is smooth over $\oo_{E, (\fkp)}$, and the inverse limit $$ \Shh_{K_p} = \varprojlim_{K^p} \Shh_{K_p K^p}$$ is the canonical smooth integral model of $\Sh_{K_p}$. We have $\Shh_{K_pK^p} = \Shh_{K_p}/K^p$, so the notation here is consistent with that in Definition \ref{defn:integral model}.
 
 In the sequel, the symbol $K^p$ will always denote a $\mathrm{GSp}$-neat compact open subgroup of $G(\A_f^p)$, and we write $K$ for $K_p K^p$. Over $\Shh_{K} = \Shh_{K_p}/K^p$, we have an abelian scheme up to prime-to-$p$ isogeny $\cA$\index[Oct2026-sorted-notation]{a 0 03 00 a abelian scheme@$\mathcal A$ (abelian scheme)}\index[Oct2026-sorted-notation]{a 0 03 00 a subscript x abelian scheme@$\mathcal A_x$ (abelian scheme)} that is the pull-back of the universal abelian scheme up to prime-to-$p$ isogeny over $\Shh_{\spK_p \spK^p}$. For any morphism of $\oo_{E,(\fkp)}$-schemes $x: Y\to \Shh_{K}$ or $x : Y \to \Shh_{K_p}$, we denote by $\cA_x$ the pull-back of $\cA$ to $Y$. 
 
 Over $\Sh_{\mathbb K_p \mathbb K^p}$, there is a canonical $\mathbb K_p$-level structure of the universal abelian scheme. The functor of trivializations (i.e., isomorphisms with $V^*_{\Z_p}$) of the first relative \'etale $\Z_p$-cohomology of the universal abelian scheme which are compatible with this level structure is pro-represented over $\Sh_{\mathbb K_p \mathbb K^p}$ by the tower $\varprojlim_{\mathbb K_p' \subset \mathbb K_p} \Sh_{\mathbb K_p' \mathbb K^p}$, where $\mathbb K_p'$ runs over all open subgroups of $\mathbb K_p$. The pull-back of this to $\Sh_{K}$ admits a morphism from the tower $\varprojlim_{ K_p' \subset K_p} \Sh_{ K_p'  K^p}$, where $K_p'$ runs over all open subgroups of $K_p$, and this morphism is equivariant with respect to $K_p \to \mathbb K_p$. By \S \ref{para:sheaf and Hecke}, Lemmas \ref{lem:Z_s} and \ref{lem:G=G^c}, and \S \ref{para:aut sheaf}, we know that $\Sh_{K_p'K^p} \to \Sh_{K_p K^p}$ is a $K_p/K_p'$-torsor for any normal open subgroup $K_p'$ of $K_p$. From this we see that over $\Sh_K$ there is a canonical $K_p$-level structure of $\cA$ which induces the $\mathbb K_p$-level structure inherited from $\Sh_{\mathbb K_p \mathbb K^p}$. 
 
 Similarly, by the results cited above together with \S \ref{para:integral Galois}, the morphism $\Shh_{K_p K^{\prime p}} \to \Shh_{K_p K^p}$ is a $K^p/ K^{\prime p}$-torsor for any normal open subgroup $K^{\prime p}$ of $K^p$. Hence on $\Shh_K$ there is a canonical $K^p$-level structure of $\cA$, and this has a canonical trivialization over $\Shh_{K_p}$. 
 
 \subsection{}
 \label{para:s_p} Let $\kappa$ be a field extension of $E$, with a fixed algebraic closure $\ol{\kappa}$. For $z \in \Sh_{K}(\kappa)$, we write\index[Oct2026-sorted-notation]{v 0 03 00 v subscript p@$\mathcal V_p$} $$\cV_{p}(z) : = \coh^1_{\et}(\cA_{z, \bar \kappa} , \Z_p) ,$$ which is a finite free $\Z_p$-module equipped with a continuous $\Gal(\bar \kappa/ \kappa)$-action. The canonical $K_p$-level structure explained in \S \ref{para:universal av} gives rise to a distinguished  $\Gal(\bar \kappa /\kappa)$-stable $\cG(\Z_p)$-orbit of $\Z_p$-module isomorphisms $\cV_p(z) \isom V^*_{\Z_p}$. Using any of these isomorphisms, we transport the tensor $s_{\alpha}$ on $V^*_{\Z_p}$  to a tensor on $\cV_{p}(z)$, denoted by\index[Oct2026-sorted-notation]{s 1 00 00 s subscript alpha p cdot@$s_{\alpha,p,\cdot}$} $s_{\alpha, p, z}$,  
 for each $\alpha \in \bm \alpha$. By the aforementioned $\Gal(\bar \kappa /\kappa)$-stability, each $s_{\alpha, p, z}$  is stable under $\Gal(\bar \kappa/\kappa)$. The distinguished $\cG(\Z_p)$-orbit of isomorphisms $\cV_p(z) \isom V^*_{\Z_p}$ is in turn characterized as the set of tensor-preserving isomorphisms. Transporting via such an isomorphism we also obtain a $\Gal(\bar \kappa / \kappa)$-action on $V^*_{\Z_p}$ which factors through a continuous homomorphism $\rho(z): \Gal(\bar \kappa/ \kappa) \to \cG(\Z_p)$. The construction of $\rho(z)$ is canonical up to conjugation by $\cG(\Z_p)$.
 
 Suppose $\tau: \kappa \to \mathbb C$ is an $E$-algebra map. Let $z_\tau \in \Sh_K(\CC)$ be the point induced by $z$. As explained in \cite[\S 1.4.11]{kisin2012modp},
 the $\Z_p$-coefficient Betti cohomology of $\cA_{z_\tau}$ admits an isomorphism to $V^*_{\Z_p}$ which is canonical up to the $K_p$-action; the key point is that the complex uniformization $G(\QQ)\backslash X\times G(\A_f)/K \cong \Sh_K(\CC)$ can be refined via the natural isomorphism $$G(\QQ)\backslash X\times G(\A_f)/K \cong (G(\QQ)\cap K_p) \backslash X \times (K_p \times G(\A_f^p))/K,$$ which holds since $G(\Q)K_p = G(\Q_p)$. Under the Artin isomorphism, this $K_p$-orbit of isomorphisms between Betti cohomology and $V^*_{\Z_p}$ is carried to the $K_p$-orbit of isomorphisms between \'etale cohomology and $V^*_{\Z_p}$ in the last paragraph; this follows from the analogous statement for the Siegel modular variety, which is a consequence of the moduli interpretation. In view of this fact, our current definition of the tensors $s_{\alpha, p, z}$ on $\cV_p(z)$ agrees with the definition in \cite[\S 2.2]{kisin2010integral} or \cite[\S 1.3.6]{kisin2012modp}. 
 
 Similarly, if $\kappa$ is a field extension of the residue field of $E_{\fkp}$, for every $x \in \Shh_K(\kappa)$ we write \index[Oct2026-sorted-notation]{v 0 03 00 v superscript p@$\mathcal V^p$}
 $$\cV^p (x) :  = \coh^1_{\et}(\cA_{x, \bar \kappa}, \A_f^p)$$ and equip it with $\Gal(\bar \kappa/\kappa)$-stable tensors\index[Oct2026-sorted-notation]{s 1 00 00 s subscript alpha a subscript f superscript p cdot@$s_{\alpha,\mathbb A_f^p,\cdot}$} $s_{\alpha, \A_f^p, x}$. A choice of lift of $x$ to $\Shh_{K_p}$ gives rise to a tensor-preserving  isomorphism $\cV^p(x) \isom V^*_{\A_f^p} = V^*_{\QQ} \otimes_{\QQ} \A_f^p$. 
 \subsection{}\label{para:review of s_0}
 Let $x$ be a closed point of the special fiber of $\Shh_{K}$, which gives rise to the abelian variety $\cA_x$ up to prime-to-$p$ isogeny over the residue field $k(x)$. As a key construction in \cite{kisin2010integral} and \cite{kisin2012modp}, for each $\alpha \in \bm \alpha$ we have a tensor\index[Oct2026-sorted-notation]{s 1 00 00 s subscript alpha 0 cdot@$s_{\alpha,0,\cdot}$} $s_{\alpha, 0 , x}$ on the $W(k(x))$-module $\coh^1_{\cris}(\cA_x/W(k(x)))$, which is furthermore $\varphi$-invariant (after inverting $p$; $\varphi$ denoting the crystalline Frobenius). Below, we review the construction of these tensors from \cite[\S 1.3]{kisin2012modp}.

 Let $F$ be a finite extension of $E_{\fkp}$ in $\Qpbar$  whose residue field $k$ contains $k(x)$. Let $\tilde x  \in \Sh_{K}(F)$ be a point that extends to a point in $\Shh_{K} (\oo_F)$ which specializes to $x$. (Clearly such a pair $(F,\tilde x)$ exists for any prescribed $x$ and finite extension $k/k(x)$; one can moreover take $F$ to be $W(k)[1/p]$.) Since $\cA_{x,k}$ and $\cA_{\accentset{\sim}{x}}  $ are the special and generic fibers of an abelian scheme over $\oo_F$,  the $\Gal(\Qpbar/F)$-representation $\cV_p(\tilde x)[1/p] \cong  \coh^1_{\et}(\cA_{\asim{x}, \Qpbar}, \QQ_p)$ is crystalline. Moreover, by the\index{integral comparison isomorphism} \emph{integral comparison isomorphism} we have a canonical isomorphism of integral $F$-isocrystals
 \begin{align}\label{eq:integral comparison0} M_{\cris}(\cV_p(\tilde x)) =
 	M_{\cris} (\coh^1_{\et}(\cA_{ \asim{x}, \Qpbar}, \ZZ_p) \big )  \isom \coh^1_{\cris}(\cA_{x,k}/W(k)),
 \end{align}
 which refines the crystalline comparison isomorphism
 $$ D_{\cris}\big( \coh^1_{\et}(\cA_{\asim{x}, \Qpbar}, \QQ_p) \big) \isom \coh^1_{\cris}(\cA_{x,k}/W(k)) [1/p].$$ Here the functor $M_{\cris}$ is defined in \S \ref{gSmodules}: recall that $M_{\cris}(L)$ is a $W(k)$-lattice in the $W(k)[1/p]$-vector space $D_{\cris}(L\otimes \Q_p)$, for any $L \in \Rep_{\Gamma_K}^{\cris, \circ}$. This result is proved in \cite[Thm.~1.4.2]{kisin2010integral} (cf.~\cite[Thm.~1.1.6]{kisin2012modp} for a correction in the normalization) for $p>2$, and proved in \cite[Prop.~4.2]{Kim12} for $p=2$. An independent proof valid for all $p$ is given by Lau \cite{Lau14,Lau19}.\footnote{In all these references, the integral comparison isomorphism is proved more generally for $p$-divisible groups over $\oo_F$. See the proof of \cite[Thm.~2.12]{KMP16} for a historical account of the different proofs. A vast generalization has been provided by the work of Bhatt--Morrow--Scholze \cite[Thm.~14.6.3 (iii)]{BMS}.}
 
 Now since $s_{\alpha, p ,\asim{x}}$ is $\Gal(\Qpbar/F)$-invariant, we have the $\varphi_{M_{\cris}(\cV_p(\asim{x}))}$-invariant tensor $M_{\cris}(s_{\alpha, p ,\asim{x}})$ over $M_{\cris}(\cV_p(\tilde x))$ by the functoriality of $M_{\cris}$. Under the isomorphism (\ref{eq:integral comparison0}), $M_{\cris} (s_{\alpha, p , \asim{x}})$ corresponds to a $\varphi$-invariant tensor $s_{\alpha, 0 , \asim{x}}$ over $\coh^1_{\cris}(\cA_{x,k}/W(k))$.  It is shown in the proof of \cite[Prop.~2.3.5]{kisin2010integral} (cf.~\cite[Prop.~1.3.9]{kisin2012modp}) that $s_{\alpha, 0 ,\asim{x}}$ in fact comes from a tensor $s_{\alpha, 0 ,x}$  on $\coh^1_{\cris}(\cA_{x}/W(k(x)))$ (via base change along $W(k(x)) \to W(k)$) which depends only on $x$ and not on the choices of $F$ and $\tilde x$.

 \begin{lem}\label{lem:apply KL}
 	There exists a $W(k(x))$-module isomorphism $$
 	V_{\Z_p}^* \otimes_{\Z_p} W(k(x)) \isom  \coh^1_{\cris}(\cA_{x}/W(k(x))) $$ 
 	which takes $s_{\alpha}$ to $s_{\alpha,0,x}$ for each $\alpha \in \bm \alpha$.
 \end{lem}
 \begin{proof} This is proved in \cite[Cor.~1.3.4]{kisin2010integral}; for convenience, we recast 
 	the argument in the framework of \S\ref{subsec:crystalline lattices}. Let $F = W(k(x))[1/p]$, and let $\tilde x \in \Shh_{K}(F)$ be a lift of $x$ as above. As in \S \ref{para:s_p}, we fix an isomorphism $V^*_{\Z_p} \isom \cV_p(\tilde x)$ taking $s_{\alpha}$ to $s_{\alpha,p,\tilde x}$, and obtain a  representation $\rho = \rho (\tilde x):\Gal(\Qpbar/F) \to \cG(\Z_p)$, which is crystalline in this case. It suffices to find a  $W(k(x))$-module isomorphism
 	$\omega_{\rho, 0} (V^*_{\Z_p})  \isom V^*_{\Z_p} \otimes _{\Z_p} W(k(x)) $ which takes each $\omega_{\rho, 0} (s_{\alpha})$ to $s_{\alpha}$ (see \S \ref{para:crystallinerepII} for $\omega_{\rho, 0}$). For this, it is enough to know that $\omega_{\rho, 0}$ is $\otimes$-isomorphic with $\idfunc_{W(k(x))}$, which holds by \S \ref{para:KL} and Lemma \ref{keylemma}.
 \end{proof}
 \subsection{}
 \label{para:isocGstructure}
 For an abelian variety $A$ up to prime-to-$p$ isogeny over $\Fpbar$, we define\index[Oct2026-sorted-notation]{v 0 03 00 v subscript 0 a@$\mathcal V_0(A)$}
 \begin{align*}
 	\cV_{0}(A)  & : =  \coh^1_{\cris} (A_0/W(k)) \otimes _{W(k)} \Z_p^{\ur},
 \end{align*}
 where $A_0$ is a model of $A$ over some finite field $k \subset \Fpbar$. Then $\cV_0(A)$ is an object in $\Isoc_{\Q_p^{\ur}}^\circ$ (see \S\ref{para:K_0 and sigma}), and it is independent of the choices of $k$ and $A_0$ up to canonical isomorphism. For $x  \in \Shh_{K_p}(\Fpbar)$, we write\index[Oct2026-sorted-notation]{v 0 03 00 v subscript 0@$\mathcal V_0$} $\cV_0(x)$ for $\cV_0(\cA_{x})$, and write\index[Oct2026-sorted-notation]{s 1 00 00 s subscript alpha 0 cdot@$s_{\alpha,0,\cdot}$} $s_{\alpha, 0 , x}$ for the tensor on $$\cV_0(x) \cong  \coh^1_{\cris}(\cA_{y}/W(k(y))) \otimes_{W(k(y))} \ZZ_p^{\ur}$$ induced by $s_{\alpha, 0 , y}$, where $y$ is the closed point underlying the image of $x$ in $\Shh_{K}(\Fpbar)$.  By Lemma \ref{lem:apply KL}, the object $\cV_0(x)$ in $\Isoc_{\Q_p^{\ur}}^\circ$, together with the tensors $s_{\alpha, 0 ,x}$ on it, satisfies the hypotheses in \S \ref{para:setting for explicit isoc} with respect to the defining datum $(V^*_{\Z_p}, (s_{\alpha})_{\alpha \in \bm \alpha})$ for $\cG$.  By Lemma \ref{lem:Isocexplicitfiberfunctor}, there exists a pair $(\tomega_x, \iota_x)$, unique up to unique isomorphism, consisting of an integral $F$-isocrystal with $\cG$-structure\index[Oct2026-sorted-notation]{zz20 upsilon 0 00 00 upsilon subscript x@$\Upsilon_x$}
 $$ \tomega_x: \Rep_{\Z_p} \cG \To \Isoc^\circ_{\Q_p^{\ur}} $$
 and an isomorphism $\iota_x: \tomega_x(V_{\Z_p}^*) \isom  \cV_0(x)$ in $\Isoc^{\circ}_{\Q_p^{\ur}}$ taking  $\tomega_x(s_{\alpha})$ to $s_{\alpha,0,x}$ for each $\alpha \in \bm \alpha$. In the sequel, we fix such a pair once and for all. 
 \begin{lem}\label{lem:reduction}
 	Let $x\in \Shh_{K_p}(\Fpbar)$, and let $\tilde x \in \Sh_K(F)$ be a lift of the closed point underlying the image of $x$ in $\Shh_{K}(\Fpbar)$ as in \S \ref{para:review of s_0}. Choose a tensor-preserving isomorphism $t: V^*_{\ZZ_p} \isom  \cV_p(\tilde x)$, leading to $\rho(\tilde x) : \Gal(\Qpbar/F) \to \cG(\ZZ_p)$ as in \S \ref{para:s_p}. This choice determines an isomorphism $\tomega_x\isom \tomega_{[\rho(\asim x)]}$, where the right hand side is defined in \S \ref{fiberfunctoranew}. 
 \end{lem}
 \begin{proof}
 	The choice of $t$ together with the canonical isomorphism (\ref{eq:integral comparison0}) determines an isomorphism $\tomega_{[\rho(\asim x)]} (V^*_{\Z_p})  \isom \cV_0(x)$ taking $\tomega_{[\rho(\asim x)]}(s_\alpha)$ to $s_{\alpha, 0 ,x}$. We then obtain the isomorphism $\tomega_x\isom \tomega_{[\rho(\asim x)]}$ from the uniqueness of $(\tomega_x, \iota_x)$. 
 \end{proof}
 \subsection{}\label{para:defn of isog class}
 Let $x,x' \in \Shh_{K_p} (\Fpbar)$. Let $R$ be a $\Q$-algebra. Any invertible morphism in the $R$-isogeny category $$f \in \Hom ( \cA_x , \cA_{x'})\otimes_{\Z_{(p)}} R$$ induces an $\A_f^p \otimes_{\Q} R$-linear isomorphism $$\cV^p(f): \cV^p(x')\otimes_{\Q} R \isom \cV^p(x) \otimes_{\Q} R ,$$ and a $\Q_p^{\ur} \otimes_{\Q} R $-linear isomorphism
 $$ \cV_0(f) : \cV_0(x') \otimes_{\Z_p^{\ur}} \Q_p^{\ur} \otimes_{\Q} R \isom  \cV_0(x) \otimes_{\Z_p^{\ur}} \Q_p^{\ur} \otimes_{\Q} R .$$ We say that $f$ \emph{preserves $G$-structures}, if $\cV^p(f)$ takes $s_{\alpha, \A_f^p ,x'}$ to $s_{\alpha, \A_f^p,x}$ and  $\cV_0(f)$ takes $s_{\alpha, 0 ,x'}$ to $s_{\alpha, 0 ,x}$ (after inverting $p$) for each $\alpha \in \bm \alpha$. We denote by\index[Oct2026-sorted-notation]{i 0 00 00 i subscript x x prime@$I_{x,x'}$}\index[Oct2026-sorted-notation]{i 0 00 00 i subscript x@$I_x$} $I_{x,x'} (R)$ the set of all such $f$ preserving $G$-structures. 
 
 The functor $R \mapsto I_{x,x'}(R)$ is represented by a $\Q$-scheme $I_{x,x'}$. If $x = x'$, then we write $I_x$ for $I_{x,x'}$, which is a $\Q$-algebraic group. Two points $x,x' \in \Shh_{K_p}(\Fpbar)$ are said to be \emph{isogenous}, if $I_{x,x'} (\Q) \neq \emptyset$. This defines an equivalence relation on $\Shh_{K_p} (\Fpbar)$, and the equivalence classes are called \emph{isogeny classes}.\index{isogeny class} Clearly each isogeny class is stable under the $G(\A_f^p)$-action on $\Shh_{K_p} (\Fpbar)$. 
 
 By \cite[Prop.~1.4.15]{kisin2012modp}\footnote{The requirement in \textit{loc.~cit.}~that $f$ respect weak polarizations is understood up to $\Q^\times$, and is already implied by the preservation of $l$-adic tensors.} , the current definition of isogeny classes is equivalent to that in \cite[\S 1.4.14]{kisin2012modp}.
 
 \section{Strict Kottwitz tuples}
 \label{recall of defns} In this section we define the notion of \emph{strict Kottwitz tuples}, as a variant of the Kottwitz triples considered in \cite{kisin2012modp}. 
 
 \subsection{} \label{subsubsec:associated data}

 For $n \in \NN$ we write\index[Oct2026-sorted-notation]{t 0 04 00 t subscript n@$\mathfrak T_n$} $$\mathfrak T_n = G(\A_f^p) \times  G(\QQ_{ p^n}),$$ and set\index[Oct2026-sorted-notation]{t 0 04 00 t superscript str@$\mathfrak T^{\mathrm{str}}$}
 $$\Tstr :  = \varinjlim_{n \in \NN} \mathfrak T_n,$$ where $\NN$ is   directed  under divisibility, and the transition maps are given by $$\mathfrak T_n \To \mathfrak T_{nt}, \quad ( \gamma, \delta) \longmapsto (\gamma^t, \delta). $$
 
 For any field $F$ of characteristic zero and any $g\in G(F)$, we say that $g$ \emph{has stable centralizer} if $G_{g^n} = G_g$ for all $n \in \NN$. Clearly this is satisfied by all sufficiently divisible powers of $g$. 
 For $\delta \in G(\Q_{p^n})$, we write\index[Oct2026-sorted-notation]{n 0 00 00 n subscript n cdot@$N_n(\cdot)$} $N_n(\delta)$ for $\delta \sigma(\delta)\cdots \sigma^{n-1}(\delta) \in G(\Q_{p^n})$.   Associated with each element $\fkk \in \Tstr$ we have a $\QQ_v$-group $I_v = I_v(\fkk)$ for every finite place $v$, defined as follows.  Pick a representative $(\gamma, \delta) \in \mathfrak T_n$ of $\fkk$. Then define $I_l$ to be $G_{\gamma_l^{t_l}}$ for $l\neq p$, and $I_p$ to be the $\sigma$-centralizer of $\delta$ in $\Res_{\QQ_{p^{nt_p}}/\QQ_p} G$, where $t_v\in \NN$ is chosen sufficiently divisible (for each finite place $v$) such that $\gamma_l^{t_l}$ and $N_{nt_p}(\delta) = N_n(\delta)^{t_p}$ all have stable centralizers in $G$. Note that the canonical projection $$(\Res_{\QQ_{p^{nt_p}}/\QQ_p} G)_{\QQ_{p^{nt_p}}} \To G_{\QQ_{p^{nt_p}}}$$ identifies $I_{p, \QQ_{p^{nt_p}}}$ with the centralizer of $N_n(\delta)^{t_p}$ in $G_{\QQ_{p^{nt_p}}}$, so $I_p$ is also the  $\sigma$-centralizer of $\delta$ in the ind-algebraic $\QQ_p$-group\index[Oct2026-sorted-notation]{res 0 01 00 res subscript q p ur q subscript p g@$\mathrm{Res}_{\Qpur/\Q_p}G$} $\Res_{\Qpur/\Q_p}G = \varinjlim_m \Res_{\QQ_{p^m}/\QQ_p}G$. We record here that $I_{p, \QQ_{p^{nt_p}}} \hookrightarrow G_{\QQ_{p^{nt_p}}}$ is a pre-inner twisting of $\QQ_p$-groups (Definition \ref{defn:inner transfer}) since it differs from the $\sigma$-twist of itself by conjugation by $\delta$. The above definition of $I_v$ is independent of the choice of representative or $t_v$.
 
 \index{Kottwitz tuple!refinement}\index{refinement|see{Kottwitz tuple}}A \emph{refinement} of an element $\fkk \in \Tstr$  is the datum of a (connected) reductive $\Q$-group $I$ and $\QQ_v$-isomorphisms $\iota_v: I_{\QQ_v} \isom I_v(\fkk)$ for all finite $v$ satisfying the following condition: There exists a semi-simple element $\gamma_0 \in G(\QQ)$ with stable centralizer, an element $g \in G(\bar \A_f)$, an inner twisting of $\QQ$-groups $\iota_0 : I_{\Qbar} \isom G_{\gamma_0,\Qbar}$, and a  representative $(\gamma,\delta)\in \mathfrak T_n$ of $\fkk$ such that \begin{itemize}
 	\item $I/\iota_0^{-1}(Z_G) (\RR)$ is compact,\footnote{This implies that $\gamma_0$ is $\RR$-elliptic in $G$.}
 	\item for each finite $v$ we have $\iota_v = \Int(g_v) \circ \iota_0$,
 	\item we have  $$\Int(g) (\gamma_0) = (\gamma, N_n(\delta)) \in G(\A_f^p) \times G(\Q_{p^n}) \subset G(\bar \A_f).$$ 
 \end{itemize} 
 
 Suppose $(I,(\iota_v)_v)$ and $(I',(\iota_v')_v)$ are two refinements of $\fkk \in \Tstr$. Then there exist a $\QQ$-isomorphism $f: I \isom I'$ and an element $\tau \in I^{\ad}(\A_f)$ such that $\iota'_v \circ f = \iota_v \circ \Int(\tau_v)$ for all finite $v$.  Indeed, we may choose a common $\gamma_0$ for both refinements  $(I,(\iota_v)_v)$ and $(I',(\iota_v')_v)$, and then this  statement follows from the Hasse
 principle for $G_{\gamma_0}^{\ad}$. Clearly the pair $(f,\tau)$ is unique up to being replaced by $(f\circ \Int(\tau_0), \tau \tau_0)$ for some $\tau_0 \in I^{\ad}(\QQ)$.  
 
 \index[Oct2026-sorted-notation]{kt 0 04 00 kt superscript str@$\mathfrak{KT}^{\mathrm{str}}$}\index{Kottwitz tuple}\index{Kottwitz tuple!strict}\index{strict Kottwitz tuple|see{Kottwitz tuple}}We call refinable elements of $\Tstr$ \emph{strict Kottwitz tuples}, and denote the set of them by $\KTstr$. We use the symbol\index[Oct2026-sorted-notation]{a 0 06 00 a subscript f superscript star@$\mathbb A_f^*$} $\A^*_f $ to denote $ \A_f^p \times \QQ_p^{\ur}.$ 
 The right action of $G(\A_f^*)$ on $\Tstr$ given by
 \[
 (\gamma,\delta)\cdot(g^p,g_p)
 =((g^p)^{-1}\gamma g^p,g_p^{-1}\delta\sigma(g_p))
 \] (for representatives $(\gamma,\delta)\in \mathfrak T_n$ for sufficiently divisible $n$) 
 stabilizes the subset $\KTstr$. We use the quotients
 \[
 \KT=\KTstr/\cG(\Z_p^{\ur}),\qquad
 \KT/{\sim}=\KTstr/G(\A_f^*).
 \]
 \index[Oct2026-sorted-notation]{kt 0 04 00 kt@$\mathfrak{KT}$}\index{Kottwitz tuple}Elements of $\KT$ are called \emph{Kottwitz tuples}; the relation $\sim$ on either $\KTstr$ or $\KT$ means equality in the final quotient by $G(\A_f^*)$. 
 
 \subsection{}\label{para:tilde I for Kott tuple}
 Let $\fkk \in \KTstr$. For a finite place $v$, let $I_v = I_v(\fkk)$. Then $I_v$ is a $\QQ_v$-subgroup of $G_{\QQ_v}$ for $v\neq p$, and $I_{p,\Qpur}$ is a $\Qpur$-subgroup of $G_{\Qpur}$. The existence of a refinement implies that for all $v$, $I_{v,\ol \QQ_v}$ is a Levi subgroup of $G_{\ol \QQ_v}$. We define a (connected) reductive $\QQ_v$-group $\tilde I_v$ together with a $\QQ_v$-map $\tilde I_v \to I_v$ as follows. For $v\neq p$, we define $\tilde I_v$ to be the inverse image of $I_v$ in $G_{\sconn, \QQ_v}$. For $v=p$, the embedding $I_{p,\Qpbar} \to G_{\Qpbar}$ is a pre-inner twisting as noted in \S \ref{subsubsec:associated data}, and the image is a Levi subgroup. We then apply the construction in \S \ref{para:construction of tilde} to obtain $\tilde I_p$, which is equipped with a $\QQ_p$-map to $I_p$ and a $\Qpbar$-embedding to $G_{\sconn,\Qpbar}$. Note that the embeddings $I_{p,\Qpbar} \to G_{\Qpbar}$ and $\tilde I_{p, \Qpbar} \to G_{\sconn,\Qpbar}$ are both defined over $\Qpur$. In particular there is a natural $\Q_p$-embedding $\tilde I_p \hookrightarrow \Res_{\Qpur/\Q_p} G_{\sconn}$  sitting above $I_p \hookrightarrow \Res_{\Qpur/\Q_p} G$. 
 \subsection{}\label{para:special kott trip}
 Let $\fks = (T,i,h) \in \spd(G,X)$ be a special point datum for $(G,X)$ (see Definition \ref{defn:special point datum}). By Lemma \ref{lem:G=G^c}, $T$ is a cuspidal torus.
 
 Let $\mu_h \in X_*(T)$ be the Hodge cocharacter of $h.$ By Lemma \ref{lem:motelts}, we obtain from the image of $-\mu_h$ in $X_*(T)_{\Gamma_{p,0}}$ a canonical $\circsim$-equivalence class in $T(\Q_p^{\ur})^{\mot}$. Let $\delta_T$ be an element of this equivalence class. Then for $n$ sufficiently divisible, the element $$ \gamma_{0,T,n} : = N_n(\delta_T) = \delta_T \sigma (\delta_T) \cdots \sigma^{n-1}(\delta_T)$$ lies in $T(\Q)$. The image of the element 
 $$
 ((i(\gamma_{0,T,n} ))_{l \neq p}, ~  i(\delta_T))\in \mathfrak T_n  
 $$ in $\Tstr$ depends only on $\fks$ and $\delta_T$, not on $n$. By \cite[Lem.~4.3.11]{kisin2012modp}, this image lies in $\KTstr$. We denote it by\index[Oct2026-sorted-notation]{k 1 04 00 k s delta subscript t@$\mathfrak k(\mathfrak s, \delta_T)$} $$\fkk (\fks, \delta_T) = \fkk (T,i,h, \delta_T) \in \KTstr.$$
 
 Note that the $\circsim$-equivalence class of $\delta_{T}$ is determined by $\fks$. If we $\sigma$-conjugate $\delta_T$ by an element of $\cT^{\circ}(\Z_p^{\ur})$ (or even $T(\Qpur)$), the element $\gamma_{0,T,n}$ remains unchanged as long as $n$ is sufficiently divisible. It follows that the image of $\fkk(\fks, \delta_T)$ in $\KT/{\sim}$ is a well-defined invariant of $\fks \in \spd(G,X)$.\footnote{Note that the same cannot be said for the image of $\fkk(\fks, \delta_T)$ in $\KT = \KTstr  / \cG(\Z_p^{\ur})$. This is because two elements of $G(\Qpur)$ that are $\sigma$-conjugate by an element of $\cT^\circ(\Z_p^{\ur})$ need not be $\sigma$-conjugate by an element of $\cG(\Z_p^{\ur})$.} Elements of $\KT/{\sim}$ arising in this way are called \emph{special}. 
 
 \section{Connected components}\label{sec:conn comp}This section has two goals. Firstly, we discuss the set parametrizing connected components of the Shimura variety and introduce an enlargement of it. Secondly, we define a refinement $\fkH(\phi)$ of the group $\cH(\phi)$ appearing in the Langlands--Rapoport--$\tau$ Conjecture, which is closely related to the set of connected components. 
 \subsection{}
 We use $G(\RR)_+$ to denote the preimage of $G^{\ad}(\RR)^+$ in $G(\RR)$, and write $G(\QQ)_+$ for $G(\QQ) \cap G(\RR)_+$. As in \cite[Lem.~3.6.2]{kisin2012modp}, we set\index[Oct2026-sorted-notation]{zz16 pi 1 00 00 pi g@$\pi(G)$}
 $$\pi (G): = G(\QQ) _+^- \backslash G(\A_f)/
 K_p,$$ where $G(\Q)_+^-$ denotes the closure of $G(\QQ)_+$ in $G(\A_f)$. (Here $\pi(G)$ depends on $K_p$, but we suppress it from the notation.) Since $K_p$ is compact, the projection $G(\A_f) \to G(\A_f)/K_p$ is a closed map. Hence the image of $G(\Q)_+^-$ in $G(\A_f)/K_p$ is closed. It easily follows that we have canonical isomorphisms \begin{align}\label{eq:pi(G)} \pi(G)
 	\cong \varprojlim_{K^p} G(\Q)_+^- \backslash G(\A_f)/ K_p K^p \cong  \varprojlim_{K^p} G(\Q)_+ \backslash G(\A_f)/ K_p K^p ,
 \end{align}
 where $K^p$ runs through compact open subgroups of $G(\A_f^p)$.
 
 \begin{lem}\label{lem:A to pi}
 	The set $\pi(G)$ has the natural structure of an abelian group. The natural map $G(\A_f^p) \to \pi(G)$ is a surjective group homomorphism.
 \end{lem}
 \begin{proof}
 	By strong approximation (see \cite[\S 2.5.1]{deligne1979varietes}), $G(\Q)_+^-$ contains the image of $G_{\sconn} (\A_f) \to G(\A_f)$. Since this image is a normal subgroup of $G(\A_f)$ with abelian quotient (see \cite[\S 2.0.2]{deligne1979varietes}), we see that $\pi(G)$ is naturally an abelian group. The second statement follows from \cite[Lem.~2.2.6]{kisin2010integral} (which uses that $K_p$ is hyperspecial).
 \end{proof}
 \begin{lem}\label{lem:about closure}The subgroup $G(\QQ)_+^-$ of $G(\A_f)$ is generated by $G(\QQ)_+$ and the image of $G_{\sconn}(\A_f) \to G(\A_f).$ \end{lem}
 \begin{proof}
 	Write $Q$ for the image of $G_{\sconn}(\A_f) \to G(\A_f).$ We first show that $G(\QQ) Q$ is closed in $G(\A_f)$ using an argument similar to \cite[\S 2.0.15]{deligne1979varietes}. By Lemma \ref{lem:G=G^c}, $G^{\ab}$ is a cuspidal torus. By Lemma \ref{lem:cusp TFAE}, $G^{\ab} (\QQ) $ is discrete in $G^{\ab}(\A_f)$. It follows that $G_{\der} (\QQ) Q$ is an open subgroup of $G(\QQ) Q$. On the other hand, by \cite[Cor.~2.0.8]{deligne1979varietes}, $G_{\der} (\QQ) Q$ is closed in $G_{\der} (\A_f)$. Hence $G(\QQ ) Q$ is closed in $G(\A_f)$.
 	
 	Consequently, $G(\QQ) _+^- \subset G(\QQ) Q$. Now by strong approximation (see for instance \cite[Thm.~7.12]{plantonov-rapinchuk}) applied to $G_{\sconn}$ and by the connectedness of $G_{\sconn} (\RR)$ (Cartan's theorem), we know that $Q \subset G(\QQ) _+^-$. Thus we have reduced to showing that
 	$$G(\QQ) \cap G(\QQ)_+^- \subset G(\QQ) _+.$$ Let $g \in G(\QQ) \cap G(\QQ) _+^-$. By \cite[Cor.~2.0.7]{deligne1979varietes}, there exists an open subgroup $U$ of $G(\A_f)$ such that $U \cap G(\QQ) \subset G(\RR) ^+$. Since $g \in G(\QQ)_+^-$, there exists
 	$g' \in g U \cap  G(\QQ) _+$. Then $g^{-1} g' \in U \cap G(\QQ) \subset G(\RR)^+$, and so $g \in g' G(\RR) ^+$. Since $g \in G(\Q)$ and $g' \in G(\Q)_+$, we conclude that $g \in G(\QQ) _+$, as desired. 	
 \end{proof}
 \subsection{}
 \label{para:pi(G,X)}
 Now consider the following set (cf.~\cite[\S 2.1.3]{deligne1979varietes}; again it depends on $K_p$, but we suppress it from the notation)\index[Oct2026-sorted-notation]{zz16 pi 1 00 00 pi g x@$\pi(G,X)$} $$\pi (G,X) :=\varprojlim_{K^p} \pi_0(\Sh_{K_pK^p}(\CC)) =  \varprojlim_{K^p} G(\QQ) \backslash \left( \pi_0 (X) \times G(\A_f)/ K^pK_p \right) ,$$ where $K^p$ runs through compact open subgroups of $G(\A_f^p)$.
 There is a natural map $\pi_0(X) \to \pi(G,X)$, induced by $$\pi_0(X) \To  \pi_0 (X) \times G(\A_f), \quad C \longmapsto (C,1). $$ This allows us to speak of the image of an element of $X$, or an element of  $\pi_0(X)$, inside $\pi(G,X)$.
 
 By (\ref{eq:pi(G)}) we see that $\pi(G,X)$ is a $\pi(G)$-torsor. We have a natural $G(\A_f^p)$-equivariant map
 $$\Shh_{ K_p}(\Fpbar) \To \pi (G,X),$$ defined as the inverse limit of the natural maps\footnote{The bijection here follows from the smoothness of $\Shh_{K_pK^p}$ over $\O_{E_\fkp}$. In fact, by the results in \cite{Madapusi19} we have $\pi_0(\Shh_{K_pK^p, \Fpbar}) \cong \pi_0 (\Shh_{K_pK^p, \O_{\ol{E_{\mathfrak p}}}})$, but we shall not need this result.} $$\Shh_{K_pK^p}(\Fpbar) \To \pi_0 (\Shh_{K_pK^p, \O_{\ol{E_{\mathfrak p}}}}) \cong  \pi_0 (\Sh_{K_pK^p} (\CC)) .$$ In particular, for each isogeny class $\mathscr I$ in $\Shh_{K_p} (\Fpbar)$, we have a natural $G(\A_f^p)$-equivariant map\index[Oct2026-sorted-notation]{c 1 00 00 c subscript i@$c_{\mathscr I}$}
 $$ c_{\mathscr I} : \mathscr I \To \pi (G,X),
 $$ which is surjective by Lemma \ref{lem:A to pi}.
 
 \begin{defn}\label{defn:pi^*} We set\index[Oct2026-sorted-notation]{zz16 pi 1 00 00 pi superscript star g@$\pi^*(G)$} $$\pi^*(G)   : = G(\QQ)_+ \backslash G(\A_f^*) / \cG(\Z_p^{\ur}) G_{\sconn}(\A_f^*). $$ This is naturally an abelian group, and is a quotient group of the abelian group $G(\A_f^*)/ G_{\sconn}(\A_f^*)$ (cf.~\cite[\S 2.0.2]{deligne1979varietes}).
 \end{defn}
 \subsection{} We have a natural homomorphism of abelian groups $$G(\Q_p)/\cG(\Z_p) G_{\sconn} (\Q_p) \To G(\Qpur)/\cG(\Z_p^{\ur}) G_{\sconn}(\Qpur).$$ This map is injective for any parahoric model $\cG$ of $G$. Indeed, the inverse image of $\cG(\Zpur)$ in $G_{\sconn}(\Qpur)$ is $\cG_{\sconn}(\Zpur)$, the parahoric subgroup of $G_{\sconn}(\Qpur)$ defined by the facet corresponding to $\cG$. The map $\cG_{\sconn}(\breve \ZZ_p) \to \cG_{\sconn}(\breve \ZZ_p), g \mapsto g \sigma(g)^{-1}$ is surjective (Greenberg's theorem), so any $g\in G_{\sconn}(\Qpur)$ satisfying $g \sigma(g)^{-1} \in \cG_{\sconn}(\Zpur)$ must lie in $\cG_{\sconn}(\Zpur) G_{\sconn}(\QQ_p)$, from which the injectivity follows.
 
 Consequently, the natural inclusion $G(\A_f) \hookrightarrow G(\A_f^*)$ induces an injective group homomorphism $\pi(G) \hookrightarrow \pi^*(G)$. We denote by $\pi^*(G,X)$\index[Oct2026-sorted-notation]{zz16 pi 1 00 00 pi superscript star g x@$\pi^*(G,X)$} the $\pi^*(G)$-torsor obtained by pushing out the $\pi(G)$-torsor $\pi(G,X)$ along $\pi(G) \hookrightarrow \pi^*(G)$, and view $\pi(G,X)$ canonically as a subset of $\pi^*(G,X)$. 
 \subsection{}\label{subsec:frak H}
 We now move on to the second goal of this section. Let $\phi : \Qf \to \G_G$ be an admissible morphism (see Definition \ref{defn of adm morphism}). In \S \ref{para:defn H} and Lemma \ref{lem:cohomological H(phi)} we defined the abelian group $$\mathcal H(\phi) = I_{\phi}(\A_f)\backslash I_{\phi}^{\ad}(\A_f) / I_{\phi}^{\ad} (\Q) \cong \D (Z_{I_{\phi}}, I_{\phi} ;\A_f) / \mathfrak D(Z_{I_{\phi,\sconn}} , I_{\phi,\sconn} ; \QQ) .$$ Our next goal is to construct an abelian group $\mathfrak H(\phi)$ equipped with a surjection to $\mathcal H(\phi)$, using $\pi(G)$ introduced above. In the rest of this section the only special property of $G$ that we will use is the conclusion of Lemma \ref{lem:about closure}. Thus the discussion applies to any Shimura datum $(G,X)$ such that $G^{\ab}$ is cuspidal. 
 
 We first define a homomorphism from the abelian group $I_{\phi,\sconn}(\A_f) \backslash I_{\phi}(\A_f)$ (cf.~\cite[\S 2.0.2]{deligne1979varietes}) to $\pi (G)$.  In \S \ref{para:transfer I_phi} we noted that the canonical embedding $I_{\phi, \Qbar} \to G_{\Qbar}$ is a pre-inner twisting. Thus we have homomorphisms $\coh^i_{\ab}(\A_f, I_\phi) \to \coh^i_{\ab}(\A_f, G)$ for all $i$, as in \S \ref{subsubsec:notation of Sha}.   For $i = 0$, $\coh^0_{\ab}(\A_f, H)$ is identified with the (abelian) quotient group $H_{\sconn}(\A_f) \backslash H(\A_f)$ for any reductive group $H$ over $\QQ$. Indeed, by the same proof as for \cite[Prop.~1.6.8]{Lab99}, $\coh^0_{\ab}(\A_f, H)$ is identified with $\coh^0(\Gamma,\cdot)$ applied to the crossed module $(H_{\sconn}(\bar \A_f) \to H(\bar \A_f))$. The latter is then identified with $H_{\sconn}(\A_f) \backslash H(\A_f)$ by the triviality of $\coh^1(\Gamma, H_{\sconn}(\bar \A_f) ) = \coh^1(\A_f, H_{\sconn})$ (Kneser's theorem). Thus we obtain a  canonical homomorphism \begin{align}\label{eq:map betw H0ab}
 	\coh^0_{\ab}(\A_f, I_\phi) \cong I_{\phi,\sconn}(\A_f) \backslash I_{\phi}(\A_f) ~\To~ \coh^0_{\ab}(\A_f, G) \cong G_{\sconn}(\A_f) \backslash G(\A_f).
 \end{align} As $\pi(G)$ is a quotient group of the right hand side, we obtain a canonical homomorphism \begin{align}\label{eq:map to pi(G)}\coh^0_{\ab}(\A_f, I_\phi) \cong  I_{\phi,\sconn}(\A_f) \backslash I_{\phi}(\A_f)  \To \pi (G).  
 \end{align} We denote its kernel by $$ \coh^0_{\ab}(\A_f, I_\phi)^{\pi}.$$

 We define \index[Oct2026-sorted-notation]{h 0 04 00 h phi@$\mathfrak H(\phi)$}
 $$ \mathfrak H(\phi) : =   \ker (I_{\phi}(\A_f) \to \pi(G)) \backslash I_{\phi}^{\ad}(\A_f) / I_{\phi}^{\ad}(\QQ).  $$ This is equipped with the natural surjection to $\mathcal H(\phi)$.
 
 If we identify $I_{\phi, \sconn} (\A_f) \backslash I_{\phi}^{\ad}(\A_f)$ with $\mathfrak D(Z_{I_{\phi, \sconn}}, I_{\phi,\sconn} ; \A_f) = \coh^1(\A_f, Z_{I_{\phi,\sconn}})$, then as in the proof of Lemma \ref{lem:cohomological H(phi)} the image of $I_{\phi}^{\ad}(\Q)$ in $\coh^1(\A_f, Z_{I_{\phi,\sconn}})$ is identified with the image of the subgroup  $\mathfrak D(Z_{I_{\phi, \sconn}}, I_{\phi,\sconn} ; \QQ) \subset \coh^1(\Q, Z_{I_{\phi, \sconn}})$. It follows that we have a canonical isomorphism of abelian groups 
 \begin{align}\label{eq:pres fkH}
 	\mathfrak H(\phi) \cong   \coh^0_{\ab} (\A_f, I_{\phi})^{\pi} \backslash \coh^1(\A_f, Z_{I_{\phi,\sconn}}) / \mathfrak D(Z_{I_{\phi, \sconn}}, I_{\phi,\sconn} ; \QQ),
 \end{align} where the map $ \coh^0_{\ab} (\A_f, I_{\phi})^{\pi} \to  \coh^1(\A_f, Z_{I_{\phi,\sconn}}) $ is the restriction of the natural map $$ \coh^0_{\ab}(\A_f, I_{\phi}) = \coh^0(\A_f, Z_{I_{\phi,\sconn}} \to Z_{I_{\phi}}) \To  \coh^0(\A_f, Z_{I_{\phi,\sconn}} \to 1)  = \coh^1(\A_f, Z_{I_{\phi,\sconn}}). $$  By this presentation and Lemma \ref{lem:KHC theorem}, $\mathfrak H(\phi)$ depends on $\phi$ only via the complex  $Z_{I_{\phi, \sconn}} \to Z_{I_{\phi}} $  and the map $\pi_1(I_\phi) \to \pi_1(G)$ (used in the definition of $\coh^0_{\ab} (\A_f, I_{\phi})^{\pi}$). By the same argument as in \S \ref{subsubsec:cocycle relation}, for any admissible morphisms $\phi \approx \phi_1$  the abelian groups $\fkH(\phi)$ and $\fkH(\phi_1)$ are canonically isomorphic.  
 
 \subsection{}\label{subsec:refined tori-rational}
 Let $T$ be a maximal torus in $I_{\phi}$. In Definition \ref{defn:rectification} we considered the natural map $\mathcal H(\phi) \to \coh^1(\A_f,T)/\Sha_G^{\infty}(\QQ,T)$. In the following we construct a lift of this map. 
 
 The composite map $T_{\Qbar} \hookrightarrow I_{\phi,\Qbar} \hookrightarrow G_{\Qbar}  \to G^{\ab}_{\Qbar}$ is defined over $\QQ$. We denote its kernel by $T^{\dagger} \subset T$. As usual, denote by $\cG^{\ab}$ the connected N\'eron model of $G^{\ab}_{\QQ_p}$ over $\ZZ_p$. The map $Z_{I_{\phi,\sconn}} \to I_{\phi}$ factors through $T^{\dagger}$, and the map $\coh^1(\A_f, Z_{I_{\phi,\sconn}} ) \to \coh^1(\A_f, T^{\dagger})$ sits in the following 
 commutative diagram:  
 $$ 
 \begin{tikzcd}[
 	row sep={2.2em,between origins},
 	column sep=2.8em,
 	nodes={inner sep=2pt}
 	]
 	\coh^0_{\ab}(\A_f, I_{\phi})^\pi \arrow[rr] \arrow[ddd] & &
 	\coh^1(\A_f, Z_{I_{\phi,\sconn}}) \arrow[ddd] &
 	\mathfrak{D}(Z_{I_{\phi, \sconn}}, I_{\phi,\sconn}; \QQ) \arrow[l] \arrow[ddd, "\gamma"] \\[1.8ex]
 	& \mathcal{G}(\Z_p) \arrow[d] \arrow[ldd] & & \\
 	& \mathcal{G}^{\ab}(\Z_p) \arrow[d] & & \\[1.8ex]
 	\coh^0_{\ab}(\A_f, G) \arrow[r] &
 	G^{\ab}(\A_f) \arrow[r, "\partial"] &
 	\coh^1(\A_f, T^{\dagger}) &
 	\Sha_{G_{\der}}^{\infty}(\QQ, T^{\dagger}) \arrow[l] \\[2.5ex]
 	G(\QQ)_+ \arrow[u] \arrow[r] &
 	G^{\ab}(\QQ) \arrow[u] \arrow[r, "\partial'"] &
 	\coh^1(\QQ, T^\dagger) \arrow[u] &
 \end{tikzcd} $$
 where $\partial$ and $\partial'$ are associated with the short exact sequence $1\to T^{\dagger} \to T \to G^{\ab} \to 1$, and $\gamma$ is induced by the natural map $Z_{I_{\phi,\sconn}} \to T^{\dagger}$; the well-definedness of $\gamma$ follows from Lemma \ref{lem:KHC theorem}. By definition and Lemma \ref{lem:about closure}, the image of $\coh^0_{\ab}(\A_f, I_\phi)^\pi$ in $\coh^0_{\ab} (\A_f, G)$ is contained in the subgroup generated by the images of $G(\Q)_+$ and $\mathcal G (\Z_p)$. By Proposition \ref{prop:G(Q)_+}, the image of $G(\Q)_+$ under the composition of the last row is contained in the image of $\ker(\coh^1(\QQ, Z_{G_{\der}}) \to \coh^1(\QQ, G_{\der}) \times \coh^1(\RR, Z_{G_{\der}}))$ under $\coh^1(\QQ, Z_{G_{\der}}) \to \coh^1(\QQ, T^{\dagger})$. Hence the image of $G(\QQ)_+$ in $\coh^1(\A_f, T^{\dagger})$ is contained in the image of $\Sha_{G_{\der}}^{\infty} (\Q, T^{\dagger})$.   We conclude, in view of (\ref{eq:pres fkH}), that there is a canonical homomorphism 
 $$\mathfrak H(\phi) \To \cok\big( \cG^{\ab}(\Z_p) \xrightarrow{\partial} \coh^1(\A_f, T^{\dagger})/ \Sha^{\infty}_{G_{\der}} (\Q, T^{\dagger})\big) . $$ We denote the right hand side by\index[Oct2026-sorted-notation]{h 0 04 00 h phi t@$\mathfrak H(\phi,T)$} $\mathfrak H(\phi, T)$. Note that  $\mathfrak H(\phi, T)$ naturally maps to $\coh ^1 (\A_f, T) / \Sha _{G} ^{\infty} (\QQ, T)$, and we have a commutative diagram 
 \begin{equation} \label{diag:lifting tori-rational}
 	\xymatrixcolsep{5pc} \xymatrix{ \mathfrak H(\phi) \ar@{->>}[d] \ar[r] & \mathfrak H(\phi, T) \ar[d] \\  	\mathcal H( \phi) \ar[r]^-{(\ref{eq:defn of tori rational})} & \coh ^1 (\A_f, T) / \Sha _{G} ^{\infty} (\QQ, T).  }
 \end{equation} 
 
 \section{Uniformization of an isogeny class}\label{subsec:unif} In this section, we introduce a canonical formal enlargement $\II^*$ of an isogeny class $\II \subset \Shh_{ K_p}(\Fpbar)$, and obtain compatible uniformizations of $\II$ and $\II^*$. Refined strict Kottwitz tuples are attached to points in the uniformizations. We also extend the connected component map $\II \to \pi(G,X)$ to a map $\II^* \to \pi^*(G,X)$.  
 \begin{defn}\label{defn:torsor over groupoid}
 	Let $\II$ be a small connected groupoid category. Let $H$ be a group. By a\index{torsor over a groupoid} \emph{right $H$-torsor over $\II$}, we mean a functor from $\II$ to the category of right $H$-torsors. Let $Y$ be such a functor. Define $$ \bar Y (\II) : = \mathrm{colim}_{\II} Y.$$ This is a right $H$-set and canonically isomorphic to $\Aut(x)\backslash Y(x)$ for any $x\in \II$.\index[Oct2026-sorted-notation]{y 0 00 20 bar y@$\bar Y$}
 \end{defn}
 \subsection{}\label{para:defn of Y(x)}
 For $x \in \Shh_{K_p}(\Fpbar)$, we define\index[Oct2026-sorted-notation]{y 0 00 00 y x@$Y(x)$}
 $$ Y(x) : = Y_p(x) \times Y^p(x),$$
 where $Y_p(x)$ is the right $G(\Q_p^{\ur})$-torsor $Y(\tomega_x)$, and $Y^p(x)$ is the right $G(\A_f^p)$-torsor consisting of $\A_f^p$-module isomorphisms $V^*_{\A_f^p} \isom \cV^p(x)  $ taking $s_{\alpha}$  to $s_{\alpha, \A_f^p, x}$ for each $\alpha \in \bm \alpha$. (By \S \ref{para:s_p}, $Y^p(x)$ has a canonical trivialization.) Thus $Y(x)$ is a right $G(\A_f^*)$-torsor.

 The set $Y_p(x)$ can be interpreted without reference to $\tomega_x$ as follows. On identifying $\tomega_x(V^*_{\Z_p})$ with $\cV_0 (x)$ as in \S \ref{para:isocGstructure}, we obtain from each element of $Y(\tomega_x)$ a $\Qpur$-linear isomorphism $V^*_{\Q_p^{\ur}} \isom \cV_0(x)[1/p]$ taking $s_{\alpha}$ to $s_{\alpha, 0 ,x}$. This defines a bijection from $Y(\tomega_x)$ to the set of such isomorphisms, and the subset $Y(\tomega_x)^{\circ} \subset Y(\tomega_x)$ corresponds to those isomorphisms that map  $V^*_{\Z_p^{\ur}}$ onto $\cV_0(x)$.
 
 Now let $y = (y_p, y^p) \in Y(x)$. Let $n \in \ZZ_{\geq 1}$ be sufficiently divisible such that $\FF_{p^n} \supset \oo_{E, \fkp}/ \fkp$ and the image of $x$ in $\Shh_{K}(\Fpbar)$ comes from an $\FF_{p^n}$-rational point $x_n$. Then the $p^n$-Frobenius acting on $\cV^p(x)$ fixes $s_{\alpha, \A_f^p ,x}$ for all $\alpha \in \bm \alpha$ by \S \ref{para:s_p}, and hence corresponds via $y^p: V^*_{\A_f^p} \isom \cV^p(x)$ to an element $\gamma_n \in G(\A_f^p)$ acting on $V^*_{\A_f^p}$. On the other hand, attached to $y_p \in Y_p(x) = Y(\tomega_x)$ we have the element $\delta_{y_p} \in G(\Q_p^{\ur})$ as in \S \ref{para:delta assoc to isoc}. This element also has the following  concrete interpretation: Under the isomorphism $V^*_{\Q_p^{\ur}} \isom \cV_0(x)[1/p]$ associated with $y_p$  as in the above paragraph, the crystalline Frobenius acting on the right hand side corresponds to $\delta_{y_p} \sigma$ acting on the left hand side.
 
 Let $\fkk(y) \in \Tstr$ be the image of $(\gamma_n,\delta_{y_p}) \in \mathfrak T_n$ (assuming $n$ is sufficiently divisible such that $\delta_{y_p} \in G(\QQ_{p^n})$); this is independent of $n$. By \cite[\S 2.3]{kisin2012modp}, we have $\fkk(y) \in \KTstr$. Thus we have obtained a map
 \begin{align}\label{eq:Kott triple on Y(x)}
 	Y(x) \To \KTstr, \quad y \longmapsto \fkk(y),  
 \end{align} which is clearly $G(\A_f^*)$-equivariant.  
 
 \subsection{}\label{para:inclusionIx}
 Let $y \in Y(x)$. The left $I_{x} (\Q)$-action on the right $G(\A_f^*)$-torsor $Y(x)$ gives rise to a homomorphism
 \begin{align}\label{eq:iota_y}
 	\iota_y: I_x(\Q) \To G(\A_f^*)
 \end{align} defined by $$ j \cdot y = y\cdot \iota_y(j), \quad \forall j \in I_x(\Q).$$ Thus we have a map $\iota_{y,v} : I_x(\Q) \to G(\Q_v)$ for each prime $v \neq p$, and a map $\iota_{y,p} : I_x(\Q) \to G(\Q_p^{\ur})$. By the results in \cite[\S 2.3]{kisin2012modp}, for each finite $v$ the map $\iota_{y,v}$ is induced by an isomorphism of $\Q_v$-groups $I_{x,\Q_v} \isom I_v$, which we still denote by $\iota_{y,v}$.\index[Oct2026-sorted-notation]{zz09 iota 1 00 00 iota subscript y v@$\iota_{y,v}$}   Moreover, $(I_x, (\iota_{y,v})_v)$ is a refinement of $\fkk(y)$.

 \subsection{}\label{para:functoriality of Y}
 Let $\II \subset \Shh_{K_p}(\Fpbar)$ be an isogeny class. We view $\II$ as a small connected groupoid, where the set of morphisms between $x$ and $x'\in \II$ is given by $I_{x,x'}(\Q)$ defined in \S \ref{para:defn of isog class}.  The association $\II \ni x \mapsto Y(x)$ is naturally a right $G(\A_f^*)$-torsor over $\II$, from which we obtain the right $G(\A_f^*)$-set $\bar Y(\II)$, as in Definition \ref{defn:torsor over groupoid}.  For every $x \in \II$ we have a canonical isomorphism
 \begin{align}\label{eq:can iso from bar Y(I)}
 	I_x(\Q) \backslash Y(x) \isom 	\bar Y (\II).  
 \end{align}
 
 For $x, x' \in \II$ and $f\in I_{x,x'}(\Q)$, the bijection $Y(f) : Y(x) \isom Y(x')$ induced by $f$ commutes with the maps $Y(x) \to \KTstr$ and $Y(x') \to \KTstr $ as in (\ref{eq:Kott triple on Y(x)}). Hence the latter maps descend to a canonical map 
 \begin{align}\label{eq:Kott triple on bar Y}
 	\bar Y(\II) \To \KTstr,  \quad \bar y \longmapsto \fkk (\bar y) ,
 \end{align}
 which is independent of the choice of $x\in \II$. 
 
 \subsection{} \label{subsubsec:Y(I)}  Let $\II \subset \Shh_{K_p}(\Fpbar)$ be an isogeny class. Set\index[Oct2026-sorted-notation]{i 0 05 00 i superscript star@$\mathscr I^*$} $$\II^* = \bar Y(\II)/\cG(\Z_p^{\ur}). $$ Then we have a natural right $G(\A_f^p)$-action on $\II^*$, and a map
 \begin{align}
 	\label{eq:Kott triple on I*}
 	\II^* \To \KT = \KTstr/ \cG(\Z_p^{\ur})
 \end{align}
 induced by	(\ref{eq:Kott triple on bar Y}).
 
 For each $x\in \II$, recall that the $G(\A_f^p)$-torsor $Y^p(x)$ has a canonical trivialization. Also inside the $G(\Qpur)$-torsor $Y_p(x) = Y(\tomega_x)$ there is a canonical $\cG(\Z_p^{\ur})$-torsor $Y(\tomega_x)^\circ \subset Y(\tomega_x)$. Hence $Y(x)/\cG(\Z_p^{\ur})$ has a canonical base point, and this determines an element $x^* \in I_x(\Q) \backslash Y(x) / \cG(\Z_p^{\ur})$. We have a canonical isomorphism $   I_x(\Q) \backslash Y(x) / \cG(\Z_p^{\ur}) \cong  \II^*  $ induced by (\ref{eq:can iso from bar Y(I)}), under which we view $x^*$ as an element of  $\II^*$. Sending $x$ to $x^*$, we have obtained a map
 \begin{align}\label{eq:I to I^*}
 	\II \To \II^*. \end{align}
 By \cite[Prop.~2.1.3]{kisin2012modp}, this map is injective and  $G(\A_f^p)$-equivariant. In the sequel we shall often view it as inclusion. We now describe its image. 
 
 As in \S \ref{subsubsec:SD}, we choose $\mu \in \dmu_X^{\cG}$, namely a cocharacter in $\dmu_X(E_{\mathfrak p})$ which extends to a cocharacter of $ \cG_{\cO_{E,\mathfrak p}}$.  Let $x\in \II$. Recall from Definition \ref{defn:intrinsic ADLV} that inside $Y(\tomega_x)$ we have the $\cG(\Z_p^{\ur})$-stable subset $\tilde X_{-\mu} (\tomega_x)$, and the quotient $\tilde X_{-\mu}(\tomega_x)/\cG(\Z_p^{\ur})$ is denoted by $X_{-\mu}(\tomega _x)$. The subset $\cG(\Z_p^{\ur}) p^{-\mu}\cG(\Z_p^{\ur}) \subset G(\Q_p^{\ur})$ depends on $\mu$ only via its $\cG(\Z_p^{\ur})$-conjugacy class, which is canonical. Hence $\tilde X_{-\mu}(\tomega_x)$ and $X_{-\mu}(\tomega_x)$ are independent of the choice of $\mu$ (cf.~the independence of $\mu$ in \S \ref{subsubsec:SD}). Moreover, if $x'$ is another element of $\II$ and if $f \in I_{x, x'}(\Q)$, the bijection $Y(\tomega_x) \isom Y(\tomega_{x'})$  induced by $f$ restricts to a bijection $\tilde X_{-\mu}(\tomega_x) \isom \tilde X_{-\mu}(\tomega_{x'})$. Therefore inside $\bar Y(\II)$ we have a canonical subset of the form\index[Oct2026-sorted-notation]{y 0 00 20 bar y i superscript natural@$\bar Y(\mathscr I)^{\natural}$} $$ \bar Y(\II)^{\natural} \cong I_x(\Q) \backslash \tilde X_{-\mu}(\tomega_x) \times Y^p(x), $$
 which is independent of the choice of $x$.
 
 \begin{prop}\label{prop:isog class reformulated}
 	The image of (\ref{eq:I to I^*}) is equal to the image of $\bar Y(\II)^{\natural}$ under the projection $\bar Y(\II) \to \II^*$. In particular we obtain a uniformization $\II \cong \bar Y(\II)^{\natural}/ \cG(\Zpur)$. 
 \end{prop}
 \begin{proof} To show that the image of (\ref{eq:I to I^*}) is contained in that of $\bar Y(\II)^{\natural}$, we need to show that
 	\begin{align}\label{eq:integral pts in ADLV}
 		Y(\tomega_x)^{\circ} \subset \tilde X_{-\mu} (\tomega_x), \quad \forall x\in \II.
 	\end{align} 
 	Let $y\in 	Y(\tomega_x)^{\circ}$, corresponding to  a tensor-preserving  isomorphism $V^*_{\ZZ_p^{\ur}} \isom \cV_0(x)$ (see \S \ref{para:defn of Y(x)}). Under this isomorphism, the Frobenius on $\cV_0(x)[1/p]$ corresponds to $
 	\delta_y \in G(\Qpur)$. Under the Cartan decomposition, $\delta_y$ must lie in a double coset $\cG(\Zpur) p^{\mu_y} \cG(\Zpur)$ for some minuscule cocharacter $\mu_y$ of $G$, since the analogous statement holds with $(G,\cG)$ replaced by $(\GL(V^*_{\QQ_p}), \GL(V^*_{\ZZ_p}))$, and whether a cocharacter of $G$ is minuscule can be checked in the larger group $\GL(V^*_{\QQ_p})$. As in Lemma \ref{lem:reduction}, pick a lift $\tilde x \in \Sh_K(F)$ of the image of $x$ in $\Shh_K$ and choose a tensor-preserving isomorphism $V^*_{\Z_p} \isom \cV_p (\tilde x)$, leading to a crystalline representation $\rho = \rho(\tilde x) :  \Gal(\Qpbar/F) \to \cG(\Z_p)$ and an isomorphism $\tomega_{[\rho]} \isom \tomega_{x}$. By Corollary \ref{cor:comparing delta and mu}, $w_G(\delta_y) \in \pi_1(G)_{\Gamma_{p,0}} = \pi_1(G)$ is equal to the image of $\mu_\rho$, the cocharacter part of $D^G_{\cris}(\rho) \in \RZ_G/{\sim}$. By the proof of  \cite[Prop.~1.3.7]{kisin2012modp}, we may assume $\mu_{\rho} = -\mu$. On the other hand, $w_G(\delta_y)$ is equal to the image of $\mu_y$, so the minuscule cocharacters $\mu_y$ and $-\mu$ have the same image in $\pi_1(G)$. It follows that they are conjugate, and in particular $\delta_y \in \cG(\Zpur) p^{-\mu} \cG(\Zpur)$. This proves (\ref{eq:integral pts in ADLV}).
 	
 	The containment of $\bar Y(\II)^{\natural}$ in the image of (\ref{eq:I to I^*}) follows from the definition of isogeny classes in \cite[\S 1.4.14]{kisin2012modp}, which we have seen in \S \ref{para:defn of isog class} to be equivalent to ours. The apparent Frobenius twist as in $\upsilon = \sigma(-\mu)$ in \cite[Cor.~1.4.13]{kisin2012modp} is explained in Remark \ref{rem:isog class reformulated} below.
 \end{proof}
 \begin{rem}\label{rem:isog class reformulated} 
 	From Proposition \ref{prop:isog class reformulated}, we see that the choice of an element $x\in \II$ gives rise to a bijection
 	\begin{align*}I_x(\Q) \backslash
 		X_{-\mu}(\tomega_x)  \times Y^p(x) \isom
 		\II.
 	\end{align*}
 	If we further choose  $y\in		\tilde X_{-\mu}(\tomega_x)$, then in view of   Remark \ref{rem:intrinsic ADLV} and \S \ref{para:inclusionIx} we obtain a bijection 
 	\begin{align}\label{eq:presentation of II}
 		\iota_y(I_x(\Q)) \backslash
 		X_{-\mu}(\delta_y)  \times G(\A_f^p) \isom
 		\II. 
 	\end{align} 
 	This bijection is  essentially the same as the map \cite[(2.1.4)]{kisin2012modp}. More precisely, the element $b$ in  \textit{loc.~cit.~}is $\sigma(\delta_y)$. This Frobenius twist is due to the use of Dieudonn\'e crystal as opposed to crystalline cohomology in \textit{loc.~cit.} We have a canonical identification $X_{-\mu}(\delta_y) \isom X_{\sigma(-\mu)} (b)$, and under this  the map (\ref{eq:presentation of II}) is identified with \cite[(2.1.4)]{kisin2012modp}. 
 \end{rem}
 \subsection{}
 
 We summarize the various constructions in the following commutative diagram.
 $$ \xymatrixcolsep{4pc}\xymatrixrowsep{4pc} \xymatrix{  I_x(\Q) \backslash \tilde X_{-\mu}(\tomega_x) \times Y^p(x) \ar@{^{(}->}[r]\ar[d]_{\cong} ^{\forall x \in \II}  & I_x(\Q) \backslash Y(x) \ar[d]_{\cong} ^{ (\ref{eq:can iso from bar Y(I)}), \forall x \in \II} \\
 	\bar Y(\II)^{\natural}  \ar@{^{(}->}[r]  \ar@{->>}[d]_{\text{quot.~by~}\cG(\Z_p^{\ur}) }   &    \bar Y(\II) \ar@{->>}[d]_{\text{quot.~by~}\cG(\Z_p^{\ur}) }  \ar[r]^{(\ref{eq:Kott triple on bar Y})} &  \KTstr \ar@{->>}[d]^{\text{quot.~by~}\cG(\Z_p^{\ur}) }  \\ \II \ar@{^{(}->}[r]^{(\ref{eq:I to I^*})} & \II^* \ar[r]^{(\ref{eq:Kott triple on I*})} &  \KT  }$$

 \begin{lem}\label{lem:extra lem on conn} Let $\II \subset \Shh_{K_p}(\Fpbar)$ be an isogeny class. Let $x, x' \in \II$. Suppose that they have lifts $y, y'\in \bar Y(\II)^{\natural}$ such that $y' = yr$ for $r\in G(\A_f^*)$. Then the image of $x'$ under the composite map
 	$$   \II \xrightarrow{c_{\II}} \pi (G,X)\hookrightarrow \pi^*(G,X)$$ is equal to the image of $x$ right-multiplied by the image of $r$ under the natural homomorphism $G(\A_f^*) \to \pi^*(G)$.
 \end{lem}\begin{proof} By $G(\A_f^p)$-equivariance we reduce to the case where the component of $r$ in $G(\A_f^p)$ is trivial. In the following we view $r$ as in $G(\Qpur)$. Under the identification 
 	$$ \bar Y(\II)^{\natural} \cong I_x(\Q) \backslash \tilde X_{-\mu}(\tomega_x) \times Y^p(x) = I_x(\Q) \backslash \tilde X_{-\mu}(\tomega_x) \times G(\A_f^p) , $$
 	the elements $y$ and $y'$ are represented by $(y_p ,1)$ and $(y_pr,1)$ respectively.
 	
 	Let $\delta = \delta_{y_p} \in G(\Q_p^{\ur})$. As in Remark \ref{rem:intrinsic ADLV}, we have an identification $$X_{-\mu}(\tomega_x) = \tilde X_{-\mu}(\tomega_x)/\cG(\Z_p^{\ur}) \isom X_{-\mu} (\delta), $$ under which the image of $y_p$ (resp.~$y_pr$) in $X_{-\mu} (\tomega_x)$  corresponds to the element $1$ (resp.~$r$) of $X_{-\mu}(\delta) \subset  G(\Q_p^{\ur})/ \cG(\Z_p^{\ur})$. Let $C$ denote the composite map \begin{multline*}
 		X_{-\mu}(\delta) \cong X_{-\mu} (\tomega_x) \xrightarrow{z \mapsto (z,1)}  I_x(\Q) \backslash X_{-\mu} (\tomega_x) \times G(\A_f^p) \cong \II  \to  \pi^*(G,X). 
 	\end{multline*} It remains to show that $C(r)$ and $C(1)$ differ by $r$.

 	We follow \cite[\S\S 1.2--1.4]{kisin2012modp} closely. Let $\mathscr G$ be the $p$-divisible group $\cA_x[p^\infty]$. By \cite[\S 1.2.16, Lem.~1.2.18]{kisin2012modp}, for any sufficiently divisible positive integer $n$, any finite totally ramified extension $F/\QQ_{p^n}$, and any $\cG_{ \ZZ_{p^n}}$-adapted lifting $\tilde {\mathscr G}$ of $\mathscr G$ to $\oo_F$, there is a map 	$$ f_{\tilde {\mathscr G},y_p}: G(\QQ_p)/\cG(\ZZ_p) \To X_{-\mu} (\delta) $$ associated with the element  $y_p$.  We briefly recall its construction\footnote{Again we adjust the normalization so we do not have a Frobenius twist on $\mu$; see Remark \ref{rem:isog class reformulated}.} in the language of \S \ref{subsec:crystalline lattices}.   The $p$-adic Tate-module $T_p\tilde {\mathscr G}^* (-1)$ gives rise to an element $[\rho] \in \Crys_G$, together with a canonical identification $\tomega_{[\rho]} \isom  \tomega_x$ (cf.~\cite[Prop.~1.1.13]{kisin2012modp} and \S \ref{para:isocGstructure}). For each $\lambda \in G(\Q_p)/\cG(\Z_p)$, we have a canonical identification $Y(\tomega_{[\rho]} ^{\lambda})\isom Y(\tomega_{[\rho]})$ as explained in  \S \ref{defnYanew}. For $z_0 \in Y(\tomega_{[\rho]}^{\lambda})^{\circ}$, write its image under the composite isomorphism 
 	$$ Y(\tomega_{[\rho]}^{\lambda}) \isom Y(\tomega_{[\rho]}) \isom Y(\tomega_x) $$
 	as $y_p \lambda_0$ for $\lambda_0 \in G(\Qpur)$. The element $\lambda_0 \cG(\Z_p^{\ur}) \in G(\Qpur)/\cG(\Z_p^{\ur})$ is then independent of the choice of $z_0$. We define this element to be $f_{\tilde {\mathscr G},y_p}(\lambda)$, and it in fact lies in $X_{-\mu}(\delta)$.

 	By \cite[Prop.~1.2.23]{kisin2012modp}, we can choose $\tilde {\mathscr G}$ such that the composite map
 	$$  G(\QQ_p)/\cG(\ZZ_p) \xrightarrow{f_{\tilde {\mathscr G},y_p}}  X_{-\mu} (\delta) \To \pi_0 (X_{-\mu} (\delta) )$$ is surjective. Here $\pi_0 (X_{-\mu} (\delta))$ is the set of \emph{connected components} of the affine Deligne--Lusztig set $X_{-\mu} (\delta)$, defined in \cite{CKV}  (cf.~\cite[\S 1.2]{kisin2012modp}). We fix such a choice of $ \tilde {\mathscr G}$, and choose $g \in G(\QQ_p)$ such that $g_0:=f_{\tilde {\mathscr G},y_p}(g)$ lies in the same connected component of $X_{-\mu} (\delta)$ as $r$.  
 	
 	The argument in \cite[\S 1.4.10]{kisin2012modp} is still valid if the integral model of the Shimura variety is replaced by a connected component of the geometric special fiber, and it then shows that $C$ is constant on each connected component of $X_{-\mu} (\delta)$.  Hence $C(r) = C(g_0)$.
 	By \cite[Cor.~1.4.12]{kisin2012modp}
 	we have $C (g_0) = C(1) \cdot g.$
 	Hence it suffices to show that $g$ and $r$ have the same image in $\pi^*(G)$. For this, it suffices to show that $w_G(g)=w_G(r)$ in view of Lemma \ref{lem:ker of Kott}.
 	This is true since $w_G (r) = w_G (g_0)$ by \cite[Lem.~2.3.6]{CKV}, and $w_G(g_0) = w_G (g)$ by \cite[Lem.~1.2.18]{kisin2012modp} or Proposition \ref{prop:compbuildingpoints} .  \end{proof}
 
 \begin{prop}\label{mark:1.1} Let $\II \subset \Shh_{K_p}(\Fpbar)$ be an isogeny class. There is a unique map\index[Oct2026-sorted-notation]{c 1 00 00 c subscript i superscript star@$c_{\mathscr I}^*$}
 	$$c_{\II}^*: \II^*  \To  \pi^*(G,X)$$
 	extending $c_{\II}: \II \to \pi(G,X)$ such that the composite map 
 	$$ \bar Y(\II) \To \II^* \xrightarrow{c_{\II}^*} \pi ^*(G,X)$$ is $G(\A_f^*)$-equivariant. Here $G(\A_f^*)$ acts on $\pi^*(G,X)$ via the natural homomorphism $G(\A_f^*) \to \pi^*(G)$. 	
 \end{prop}
 \begin{proof} Fix $x \in \II$. Consider $Y(x)^{\natural} = \tilde X_{-\mu}(\tomega_x) \times Y^p(x)$, namely the inverse image of $\II \subset \II^*$ in $Y(x)$. By Lemma \ref{lem:extra lem on conn} and the fact that  $Y(x)^{\natural} \cdot G(\Qpur) = Y(x)$, the composite map $Y(x)^{\natural} \to \II \to \pi^*(G,X)$ extends uniquely 
 	to a $G(\A_f^*)$-equivariant map $Y(x) \to \pi^*(G,X)$. This already shows the uniqueness in the proposition. To show the existence, we need to check that the last map descends to the right quotient by $\cG(\Zpur)$ and the left quotient by $I_x(\QQ)$. The former descent is obvious since the map is $G(\A_f^*)$-equivariant and the subgroup $\cG(\Zpur) \subset G(\A_f^*)$ acts trivially on the target. For the latter descent, again by $G(\A_f^*)$-equivariance it suffices to show that for a fixed $y\in Y(x)^{\natural}$ the composite homomorphism $I_x(\Q) \xrightarrow{\iota_y} G(\A_f^*) \to \pi^*(G)$ is trivial, where $\iota_y$ is as in (\ref{eq:iota_y}). But this holds since we already know that $Y(x)^{\natural}$ is $I_x(\QQ)$-stable and the map $Y(x)^{\natural} \to \pi^*(G,X)$ factors through  $I_x(\QQ) \backslash Y(x)^{\natural} / \cG(\Zpur) = \II$. \end{proof}

 \section{Integral special points}\label{subsec:intsppt} In this section, we construct certain \emph{integral special points} in the uniformization in \S \ref{subsec:unif}. The main result is Proposition \ref{prop:special points on the geom side}, which computes the strict Kottwitz tuple and connected component of such a point. We also review the CM lifting results from \cite{kisin2012modp} in Theorem \ref{thm:geometric speciality}. 
 \subsection{}\label{para:x_fks}
 Let $\fks = (T,i,h) \in \spd(G,X)$ be a special point datum. Let \index[Oct2026-sorted-notation]{x 1 00 23 tilde x subscript s@$\tilde x_{\fks}$}\index[Oct2026-sorted-notation]{x 1 00 22 hat x subscript s@$\hat x_{\fks}$}
 $$ \hat x_{\fks} : = (i\circ h,1) \in \Sh_{K_p} (\CC) ,$$ and let $\tilde x_{\fks}$  be the image of 
 $\hat x_{\fks}$ in $\Sh_K(\CC)$. Then $\hat x_{\fks}$ is actually defined over $\Qbar$, and $\tilde x_{\fks}$ is defined over a finite extension of $\QQ$. Let $F/\Q_p$ be a finite extension containing the field of definition of $\tilde x_{\fks}$. Then $\hat x_{\fks}$ is defined over the maximal unramified extension $F^{\ur}$ of $F$. This follows from the reciprocity law describing the canonical model $\Sh(T,h)$ as in \S \ref{para:Shimura variety}, the discreteness of $T(\Q)$ in $T(\A_f)$ (by Lemmas \ref{lem:cusp TFAE} and \ref{lem:G=G^c}), and the neatness of $i^{-1}(K^p) \subset T(\A_f^p)$. We leave the details to the reader as an exercise. Note that this statement holds for any Shimura variety satisfying the statement about maximal tori in Lemma \ref{lem:G=G^c}, while in the current setting of Hodge type it essentially corresponds to the fact that the $\A_f^p$-Tate module of any CM abelian variety is finitely ramified over $\QQ_p$. 
 
 Applying the extension property\footnote{Here we are using the following version of the extension property which is different from that in \cite{kisin2010integral}: For any DVR $R \supset \mathcal O_{E_{\fkp}}$ of characteristic $(0,p)$, the map $\Shh_{K_p}(R) \to \Shh_{K_p}(R[1/p])$ is a bijection; cf.~\cite{KisinPappas}. This is just a direct consequence of the N\'eron-Ogg-Shafarevich criterion.} of the canonical integral model to the DVR $\oo_{F^{\ur}}$, we see that the point $\hat x_{\fks} \in \Sh_{K_p}(F^{\ur})$ extends to a point in $\Shh_{K_p} (\oo_{F^{\ur}})$. We denote by\index[Oct2026-sorted-notation]{x 1 00 00 x subscript s@$x_{\fks}$} $x_{\fks} \in \Shh_{K_p}(\Fpbar)$ the reduction.

 In the following we write $ \tilde x$ for $\tilde x_{\fks}$, and regard $i:T\to G$ as inclusion. The point $(i\circ h,1) \in X \times G(\A_f)$ gives rise to a canonical trivialization of Betti cohomology $$ V^*_{\ZZ_{(p)}} \isom  \coh^1_B(\cA_{\asim{x} }(\CC), \ZZ_{(p)} ) , $$ under which the Mumford--Tate group (after inverting $p$) is identified with a subgroup of $T$ acting on $V^*_{\QQ}$. Under comparison isomorphisms this gives rise to trivializations \index[Oct2026-sorted-notation]{triv 1 01 00 triv subscript s p@$\mathrm{triv}_{\fks,p}$}\index[Oct2026-sorted-notation]{triv 1 01 00 triv subscript s a subscript f superscript p@$\mathrm{triv}_{\fks,\A_f^p}$}
 \begin{align*}
 	\mathrm{triv}_{\fks, p}	& :  V^*_{\Z_p} \isom  \cV_p(\tilde x) , \\ \mathrm{triv}_{\fks, \A_f^p} & :
 	V^*_{\A_f^p} \isom \cV^p(x_\fks).
 \end{align*}These take $s_{\alpha}$ to $s_{\alpha,p,\asim x}$ and to $s_{\alpha, \A_f^p, x_{\fks}}$ respectively (cf.~\S \ref{para:s_p}), and the latter coincides with the canonical trivialization associated with $x_{\fks}$. 
 
 Up to enlarging $F$, we assume that $\cA_{\asim x}$ has good reduction over $\mathcal O_F$, and that the reduction of $\tilde x$ in $\Shh_K(k_F)$ agrees with the image of $x_{\fks} \in \Shh_{K_p}(\Fpbar)$.   
 Since the Mumford--Tate group of $\cA_{\asim{x}}$ is contained in $T$, there is an embedding 
 of $T$ into the group of self-isogenies of $\cA_{\asim{x}},$ and hence into the group of self-isogenies of 	$\cA_{x_{\fks}},$ the reduction of  $\cA_{\asim {x}}.$
 The isomorphisms $\mathrm{triv}_{\fks,p}[1/p]$ and $\mathrm{triv}_{\fks, \A_f^p}$ are $T(\QQ)$-equivariant, where $T(\QQ)$ acts on the right hand sides via the self-isogenies. 
 
 Write $U$ for $i^{-1}(K) \subset T(\A_f)$. 	Since $\tilde x$ comes from the Shimura variety $\Sh_U(T,h)$, a similar argument as in \S \ref{para:universal av} shows that up to enlarging $F$ (so as to contain the reflex field of $(T,h)$) the abelian variety $\cA_{\asim x}$ over $F$ has a canonical $U$-level structure over $F$, and $\mathrm{triv}_{\fks, p}[1/p]$ is compatible with that level structure. It follows, as in \S \ref{para:s_p}, that  $\mathrm{triv}_{\fks, p}[1/p]$ transports the $\Gal(\Qpbar/F)$-action on $\cV_p(\tilde x)[1/p]$ to a $T(\Q_p)$-valued Galois representation 
 \index[Oct2026-sorted-notation]{zz17 rho 1 00 00 rho subscript t h@$\rho_{(T,h)}$} $$\rho_{(T,h)}: \Gal(\Qpbar/F) \To T(\Q_p)$$ with $T(\QQ_p)$ acting on $V^*_{\QQ_p}$. 
 Hence the $T(\Q)$-action on $\cV_p(\tilde x)[1/p]$ preserves the $\Gal(\Qpbar/F)$-module structure (since $T$ is commutative) as well as the tensors $s_{\alpha, p,\asim x}$ (since $T\subset G$). Applying the functor $D_{\cris}$, we see that the $T(\Q)$-action on $\cV_0(x_{\fks})[1/p]$ via self-isogenies of $\cA_{x_\fks}$ preserves the tensors $s_{\alpha, 0, x_\fks}$. On the other hand, it immediately follows from the $T(\QQ)$-equivariance and tensor-preservation of $\mathrm{triv}_{\fks, \A_f^p}$ that $T(\QQ)$ preserves $s_{\alpha, \A_f^p, x_{\fks}}$. Hence we obtain a canonical embedding
 \begin{align}\label{eq:T to I_x}
 	T \hookrightarrow I_{x_{\fks}}.
 \end{align}
 
 We record that $\rho_{(T,h)}$ has the same germ as $$ r(\mu)_{U,p, \mathrm{loc}} : \Gal(\Qpbar/E_{\mu, U, v}) \To T(\Q_p) $$ in \S \ref{para:r(mu)_U}, with $\mu = \mu_h \in X_*(T)$. This description of $[\rho_{(T,h)}]$ follows from the reciprocity law describing the canonical model $\Sh(T,h)$ as in \S \ref{para:Shimura variety}. See also~\cite[\S (5.5)]{pink1992ladic}. By Proposition \ref{prop:local component}, the element $[\rho_{(T,h)}] \in \Crys_{T_{\Q_p}}$ is equal to $\motinv_T^{-1} (- \mu_h) \in \Mot_T$.
 
 \subsection{}\label{para:tomega_T,h}
 From $[\rho_{(T,h)}] \in \Crys_{T_{\Q_p}}$, we obtain an element $[i\circ \rho_{(T,h)}] \in \Crys_{G_{\Q_p}}$. In the sequel, given any torus $T$ over $\Q$ or $\Q_p$, we write $\cT^{\circ}$ for the connected N\'eron model of $T_{\Q_p}$ over $\Z_p$.  Applying \S \ref{fiberfunctoranew} to $[\rho] = [i \circ \rho_{(T,h)}]$ and $[\rho_T]= [\rho_{(T,h)}]$, we obtain $$\tomega_{[\rho_{(T,h)}]} \in \cT^\circ \text{-} \Isoc_{\Q_p^{\ur}}^\circ \quad \text{and} \quad  \tomega_{[i\circ \rho_{(T,h)}]} \in \cG\text{-} \Isoc_{\Q_p^{\ur}}^\circ.$$ To simplify notation we denote them by\index[Oct2026-sorted-notation]{zz20 upsilon 0 00 00 upsilon subscript t h@$\Upsilon_{(T,h)}$}\index[Oct2026-sorted-notation]{zz20 upsilon 0 00 00 upsilon subscript s@$\Upsilon_{\mathfrak s}$} $\tomega_{(T,h)}$ and $\tomega_{\fks}$ respectively. As in  \S \ref{para:rho_T}, we have a natural injection from the $T(\Qpur)$-torsor $Y(\tomega_{(T,h)})$ to the $G(\Qpur)$-torsor $Y(\tomega_{\fks})$.
 
 By Lemma \ref{lem:reduction}, the canonical $\mathrm{triv}_{\fks,p}$ determines a canonical isomorphism $
 \tomega_{\fks}\cong \tomega_{x_{\fks}}$. Using this we view the canonical injection
 $Y(\tomega_{(T,h)}) \hookrightarrow Y(\tomega_{\fks})$ as an injection
 \begin{align}
 	\label{eq:map from Y_T(x)}
 	Y(\tomega_{(T,h)}) \hookrightarrow Y(\tomega_{x_\fks}) .
 \end{align}
 Recall from \S \ref{para:delta assoc to isoc} that inside the $T(\Qpur)$-torsor $Y(\tomega_{(T,h)})$ we have the subset of integral points $Y(\tomega_{(T,h)})^\circ$, which is a $\cT^\circ(\Z_p^{\ur})$-torsor.
 
 \begin{defn}\label{defn:integral special point}  Let $\fks =(T,i,h) \in \spd(G,X)$, with associated $x_{\fks } \in \Shh_{K_p}(\Fpbar)$. Let $y_p \in Y_p(x_\fks)$ be an element of the image of $Y(\tomega_{(T,h)})^\circ$ under (\ref{eq:map from Y_T(x)}), and let $y = (y_p ,1) \in Y(x_\fks) = Y(\tomega_{x_\fks}) \times G(\A_f^p)$. We call $y$ an\index{integral special point} \emph{integral special point} associated with $\fks$.
 \end{defn}
 \begin{rem}
 	Let $\II$ denote the isogeny class of $x_{\fks}$. The image of an integral special point under the uniformization map $Y(x_{\fks}) \to \II^*$ need not coincide with $x_{\fks}$. Such images are what was referred to as integral special points in the Introduction.
 \end{rem} 
 \begin{prop}\label{prop:special points on the geom side} In the notation of Definition \ref{defn:integral special point}, the following statements hold.
 	\begin{enumerate}
 		\item Let $\delta_{y_p}$ be the image of $y_p$ under the map $Y(\tomega_{(T,h)})  \to T(\Qpur), z \mapsto \delta_z$ as in \S \ref{para:delta assoc to isoc}. Then $\delta_{y_p}$ lies in $T(\Qpur)^{\mot}$, and the $\circsim$-equivalence class of $\delta_{y_p}$ corresponds to the image of $-\mu_h$ in $X_*(T)_{\Gamma_{p,0}}$ under the bijection (\ref{eq:bij for mot}).
 		
 		\item The element $\fkk(y)\in \KTstr$ associated with $y$ as in (\ref{eq:Kott triple on Y(x)}) is equal to $\fkk (\fks, \delta_{y_p}) $ as in \S \ref{para:special kott trip}. 
 		\item Let $\II$ be the isogeny class of $x_{\fks}$. The image of $y$ under the composite map
 		$$ Y(x_{\fks}) \To I_{x_\fks}(\Q) \backslash Y(x_\fks)/ \cG(\Z_p^{\ur}) \cong  \II^* \xrightarrow{c^*_{\II}} \pi^*(G,X)$$ is equal to the image of $i\circ h\in X$ in $\pi(G,X) \subset \pi^*(G,X)$. (See Proposition \ref{mark:1.1} for $c^*_{\II}$.)
 	\end{enumerate}
 \end{prop}
 \begin{proof}
 	We have seen in \S \ref{para:x_fks} that $[\rho_{(T,h)}] = \motinv_T^{-1}(-\mu_h) \in \Mot_T$. Statement (i) follows from this fact and Corollary \ref{cor:inv of T-valued rep}.
 	
 	For (ii), tracing the definitions we see that the $\delta$-component of $\fkk(y)$ is $i(\delta_{y_p})$. 
 	By the Shimura--Taniyama reciprocity law, for sufficiently divisible $n$ the geometric $p^n$-Frobenius in $I_x(\Q)$ is given by an element $\gamma_{0,T,n}'$ of $T(\Q)$, embedded inside $I_x(\QQ)$ as in (\ref{eq:T to I_x}). We claim that $\gamma_{0,T,n}'$ is equal to $\gamma_{0,T,n}$ introduced in \S \ref{para:special kott trip}. (Recall that $\gamma_{0,T,n}$ depends only on $(T,h)$ and $n$.) 
 	
 	We prove the claim. As in \S \ref{para:defn of Y(x)}, the element $y_p \in Y_p(x_{\fks})$ determines an isomorphism $f: V^*_{\Qpur} \isom \mathcal V_0(x_{\fks})[1/p]$. The fact that $y_p$ comes from $Y(\tomega_{(T,h)})$ implies that $f$ is $T(\QQ_p)$-equivariant, where $T(\QQ_p)$ acts on the left hand side by $i: T(\QQ_p) \to G(\QQ_p)$ and acts on the right hand side by  $T(\QQ_p) \hookrightarrow I_{x_\fks}(\QQ_p)$  as in (\ref{eq:T to I_x}). Indeed, let $t \in T(\Q_p)$, and denote by $\alpha(t)$ (resp.~$\beta(t)$) the induced endomorphism of $V^*_{\QQ_p}$ (resp.~$\mathcal V_0(x_\fks)[1/p]$). We view $\alpha(t)$ as a $T$-invariant tensor on $V^*_{\QQ_p}$ (since $T$ is commutative). Under the canonical identification of $\mathcal V_0(x_\fks)[1/p]$ with $D_{\cris}(V^*_{\Q_p})$ (where $V^*_{\Q_p}$ is equipped with the Galois representation $\rho_{(T,h)}$ and canonically identified with $\cV_p(\tilde x_\fks)[1/p]$ via $\mathrm{triv}_{\fks,p}$) we have $\beta(t) = D_{\cris} (\alpha(t))$. Hence $f$ must send the tensor $\alpha(t)$ to the tensor $\beta(t)$, i.e., $f$ is $T(\Q_p)$-equivariant. 
 	Now for sufficiently divisible $n$, the $n$-fold self-composition of the crystalline Frobenius on $\mathcal V_0(x_\fks)[1/p]$ equals $\gamma'_{0,T,n} \sigma^n$, by the usual relation between the crystalline Frobenius and the Frobenius isogeny for an abelian variety over a finite field. Transporting this equality via $f^{-1}$ and using the $T(\Q_p)$-equivariance of $f$, we see that $\gamma_{0,T,n}' = N_n(\delta_{y_p})$. By this equality and part (i), we conclude that $\gamma'_{0,T,n} = \gamma_{0,T,n}$ by the unique characterization in \S \ref{para:special kott trip}.
 	The claim is proved.

 	Now the composition
 	$$ T_ {\A_f^p} \xrightarrow{(\ref{eq:T to I_x})} I_{x_\fks, \A_f^p} \To \GL(\cV^p(x_\fks)) \underset{\sim}{\xrightarrow{\mathrm{triv}_{\fks, \A_f^p} }} \GL(V^*_{\A_f^p}) $$   factors through the base change to $\A_f^p$ of $i: T \to G$. Hence the $\gamma$-component of $\fkk(y)$ is represented by $i(\gamma_{0,T,n}')$ at level $n$. By our previous claim, this is $i(\gamma_{0,T,n})$.   Statement (ii) follows.
 	
 	For (iii), let $y_p'$ be an element of $Y(\tomega_{x_\fks})^\circ \subset Y(\tomega_{x_\fks}) = Y_p(x_\fks)$, and set $y' :  = (y_p' , 1) \in Y_p(x_\fks) \times Y^p(x_\fks) = Y(x_\fks).$ Then the image of $y'$ in $\II^*$ is equal to the image of $x_\fks$ under the canonical injection $\II \hookrightarrow \II^*$; see \S \ref{subsubsec:Y(I)}. By the construction of $x_\fks$, it has the same image in $\pi(G,X)$ as $i\circ h \in X$. Therefore we only need to find $y' = (y'_p, 1)$ as above such that $y$ and $y'$ have the same image in $\pi^*(G,X)$. By the second statement in Proposition \ref{prop:intstrcomp}, we can find $y'_p$ such that it lies in the $\varrho (G_{\sconn}(\Q_p^{\ur}))$-orbit of $y_p$. But then $y$ and $y'$ have the same image in $\pi^*(G,X)$, since the map $Y(x) \to \pi^*(G,X)$ in question is $G(\A_f^*)$-equivariant, and the composition $G_{\sconn}(\A_f^*) \to G(\A_f^*) \to \pi^*(G)$ is trivial.
 \end{proof}
 
 \subsection{}\label{para:x-admissible}
 Let $x  \in \Shh_{K_p} (\Fpbar)$, and let $T$ be a maximal torus in the $\Q$-reductive group $I_x$. For each $\mu \in X_*(T)$, we define $\bar \mu ^{T_{\Q_p}} \in X_*(T)_{\QQ}$ as in \S \ref{para:bar mu} (with respect to the $\Q_p$-torus $T_{\Q_p}$).
 
 Let $y \in Y(x)$. Then we have a $\Q_p$-isomorphism $\iota_{y,p} : I_{x, \Q_p} \isom I_p$, where $I_p = I_p (\fkk(y))$ is the $\sigma$-centralizer in $\Res_{\Qpur/\Q_p} G$ of the component $\delta\in G(\Qpur)$ of $\fkk(y) \in \KTstr$. The canonical embedding $I_p \hookrightarrow J_{\delta}$ identifies $I_p$ with the centralizer in $J_{\delta}$ of  $\delta \sigma(\delta)\cdots \sigma^{t-1} (\delta) \in J_{\delta}(\Q_p)$ for a sufficiently divisible $t$. Hence the Newton cocharacter $\nu_\delta$ of $\delta \in G(\Qpur)$, viewed as a fractional central cocharacter of $J_\delta$, gives rise to a fractional central cocharacter  $\iota_{y,p}^{-1} \circ \nu_{\delta}$ of $I_x$. We say that a cocharacter $\mu \in X_*(T)$ is\index{admissible!x-admissible@$x$-admissible}\index{x-admissible@$x$-admissible|see{admissible}} \emph{$x$-admissible} if
 the composition
 \begin{align}\label{eq:mutoG}
 	\GG_{m, \Qpbar} \overset{\mu}{\To} T_{\Qpbar} \hookrightarrow I_{x, \Qpbar} \xrightarrow{ \iota_{y,p}} I_{p, \Qpbar} \hookrightarrow G_{\Qpbar}
 \end{align}
 lies in $ \dmu_X(\Qpbar)$, and if $\bar \mu^{T_{\Q_p}} = -(\iota_{y,p}^{-1} \circ \nu_{\delta}) \in X_*(T) \otimes \QQ$. This  definition is independent of the choice of $y$, similarly as in \S \ref{para:admissible cochar}.  The following theorem is the geometric analogue of Theorem \ref{thm:precise speciality}. 
 
 \begin{thm}[{\cite{kisin2012modp}}]\label{thm:geometric speciality}
 	Let $x  \in \Shh_{K_p} (\Fpbar)$, and let $T$ be a maximal torus in $I_x$ defined over $\Q$. The following statements hold.
 	\begin{enumerate}
 		\item There exists $\mu \in X_*(T)$ that is $x$-admissible.
 		\item Let $\mu \in X_*(T)$ be an $x$-admissible cocharacter. Then there exists a special point datum of the form $\fks = (T,i,h)$ satisfying the following conditions:
 		\begin{enumerate}
 			\item $\mu = \mu_h$.
 			\item The points $x$ and $x_{\fks}$ lie in the same isogeny class. Moreover, there exists $g \in I_{x,x_{\fks}}(\Q)$ such that the composition of $T  \hookrightarrow I_{x}$ followed by the isomorphism $I_x \isom I_{x_{\fks}}$ induced by $g$ is equal to the canonical embedding  (\ref{eq:T to I_x}).
 		\end{enumerate}
 	\end{enumerate}
 \end{thm}
 \begin{proof}
 	Part (i) is proved in \cite[Lem.~2.2.2]{kisin2012modp}, and part (ii) is proved in \cite[Cor.~2.2.5]{kisin2012modp}.
 \end{proof}
 \section{Analogous structures on the gerb side}\label{subsec:parallel}
 In this section, we develop the structures parallel to those in \S\S \ref{subsec:unif}--\ref{subsec:intsppt} for  admissible morphisms.  
 
 \subsection{}\label{para:defn of Y(phi)}
 Let $\phi: \Qf \rightarrow \G_G$ be an admissible morphism (Definition \ref{defn of adm morphism}). Set\index[Oct2026-sorted-notation]{y 0 00 00 y phi@$Y(\phi)$}
 \begin{align*}
 	Y_p(\phi)  & := \mathcal {UR}( \phi(p)\circ \zeta_p) ,\\ Y^p(\phi) & : = X^p(\phi),\\ Y(\phi) & :=  Y_p(\phi) \times  Y^p(\phi) .
 \end{align*}
 See Definition \ref{defn:unram morph} and \S \ref{subsubsec:defn of S(phi)} for the notations. Thus $Y_p(\phi)$ is a right $G(\Q_p^{\ur})$-torsor, $Y^p(\phi)$ is a right $G(\A_f^p)$-torsor, and $Y(\phi)$ is a right $G(\A_f^*)$-torsor.
 
 The construction in \cite[\S 4.3.9, \S 4.5.1]{kisin2012modp} gives rise to a  map $$ 
 Y(\phi)/\cG(\Z_p^{\ur}) \To \KT = \KTstr/\cG(\Z_p^{\ur}). $$
 We now recall the construction and make a key modification. While our modified map differs from the one in \textit{loc.~cit.},  it has the advantage of admitting a canonical\footnote{That is, the construction depends only on choices independent of the Shimura datum $(G,X)$ and $\phi$.} lift $Y(\phi) \to \KTstr$.

 Fix once and for all an increasing cofinal sequence\index[Oct2026-sorted-notation]{l 0 00 00 l subscript r@$L_r$}
 $L_1\subset L_2\subset\cdots$ of CM Galois fields in $\Qbar$. For each $r \geq 1$, recall the torus $P^{L_r}$ defined in \S \ref{subsec:pseudomotivic}, and the elements 
 $\delta_n \in P^{L_r}(\Q)$ defined in \S \ref{subsec:deltan} for $n$ sufficiently divisible, which we also assume is chosen such that $\nu(p)^{L_r}_n \in X_*(P^{L_r}).$  
 Using \S \ref{subsec:deltan} (iv), one sees that there exists $c_r \in P^{L_r}(\breve \Q_p)$ such that 
 \begin{align}\label{eq:def prop of c}
 	c_r\nu(p)^{L_r}_n(p^{-1})\sigma^n(c_r)^{-1} = \delta_n 
 \end{align} for all sufficiently divisible $n$. Note that the coset $c_r P^{L_r}(\Qpur)$ is uniquely determined by this condition. Hence we may and shall choose $c_r$ in $P^{L_r,\circ}(\breve \Z_p),$ where $P^{L_r,\circ}$ denotes the connected N\'eron model of $P^{L_r}_{\QQ_p}$ over $\Z_p$. We do not impose any compatibility conditions on $c_r$ for varying $r$.

 For any map $a:P_{\Qbar}\to H$ where $H$ is a linear algebraic group over $\Qbar$ occurring below, 
 put
 \[
 r_a:=\min\{r\mid a\text{ factors through }P^{L_r}_{\Qbar}\}. 
 \]
 By Lemma \ref{lem:G=G^c} and
 \cite[Lem.~B3.9]{reimann1997zeta}, each admissible morphism $\phi:\mathfrak Q\to \G_G$ 
 factors through $\mathfrak P$. We regard $\phi^\Delta$ as a morphism $P_{\Qbar} \to G_{\Qbar}$ and write\index[Oct2026-sorted-notation]{r 1 00 00 r subscript phi@$r_\phi$} $r_\phi$ for $r_{\phi^{\Delta}}$. We denote by\index[Oct2026-sorted-notation]{zz21 phi 1 00 00 phi superscript square@$\phi^{\square}$} $\phi^{\square}$ the unique factorization $P^{L_{r_\phi}}_{\Qbar} \to G_{\Qbar}$.

 We now define the map $Y(\phi) \to \KTstr$ by employing the construction in \cite[\S 4.5.1]{kisin2012modp}, but replacing the element $c c'^{-1} \in P^L (\breve \Q_p)$ in \textit{loc.~cit.~}by our fixed $c_{r_\phi} \in P^{L_{r_\phi}}(\breve \Q_p)$. Thus for each $y = (y_ p , y^p) \in Y(\phi) = Y_p(\phi) \times Y^p(\phi)$ and for a sufficiently divisible $n \in \Z_{\geq 1}$, let $\gamma = (\gamma_l)_{l\neq p} \in G(\A_f^p)$ be given by $(y^p)^{-1} \phi^\square(\delta_n) y^p$ and let $$b_{y_p} = b_{\Int (y_p) ^{-1} \circ \phi(p) \circ \zeta_p} \in G(\Qpur)$$ in the notation of Definition \ref{defn:b_theta} and \S \ref{subsubsec:recall of defns in LR}, i.e., $b_{y_p} \rtimes \sigma $ is the image of $d_\sigma$ under the morphism $ \D \to \G_{G_{\QQ_p}}^{\ur}$ underlying the unramified morphism $\Int(y_p)^{-1} \circ \phi(p) \circ \zeta_p : \G_p \to \G_G (p)$.  The morphism $\Int (y_p)^{-1} \circ \phi^{\square} : P^{L_{r_\phi}}_{\Qpbar} \to G_{\Qpbar}$ is actually defined over $\Qpur$ (since $\Int(y_p)^{-1} \circ \phi(p) \circ \zeta_p$ is unramified), and\index[Oct2026-sorted-notation]{c 1 00 00 c subscript y subscript p@$c_{y_p}$} we set\footnote{Our construction of $\delta$ differs from that in \cite[\S 4.5.1]{kisin2012modp} only in the substitution of $c_{r_\phi}$ for $cc'^{-1}$, modulo unraveling the definition of $\gP^{L}(p)^{\ur}$ used in \textit{loc.~cit.}} 
 \begin{align*}
 	c_{y_p} & : = \Int y_p^{-1} \circ \phi^{\square} (c_{r_\phi}) \in G(\breve \QQ_p), \\
 	\delta & :  = c_{y_p} b_{y_p} \sigma (c_{y_p})^{-1} \in G(\breve \QQ_p).
 \end{align*} 
 
 By the same argument as in \cite[\S 4.5.1]{kisin2012modp}, we have $\delta \in G(\Qpur)$, and the image of $(\gamma,\delta) \in \mathfrak T_n$ in $\Tstr$ lies in $\KTstr$, which we denote by $\fkk(y) \in \KTstr$. We thus obtain the map 
 \begin{align}\label{eq:Kott triple on Y(phi)}
 	Y(\phi) \To \KTstr, \quad y \longmapsto \fkk(y).  
 \end{align} This map is clearly $G(\A_f^*)$-equivariant.  
 
 Let $y\in Y(\phi)$. By \cite[\S 4.5.1]{kisin2012modp}, for each finite place $v$ we have an isomorphism $\iota_{y,v}: I_{\phi,\Q_v} \isom I_v(\fkk(y))$ induced by $\Int (y_v)^{-1}$,\footnote{Here $I_p(\fkk(y))$ is viewed as a $\Qpbar$-subgroup of $G$. The fact that $\Int(y_p)^{-1}$ induces an isomorphism $I_{\phi,\QQ_p} \isom I_p (\fkk(y))$ is independent of our choice of $c_{r_\phi}$ in \S \ref{para:defn of Y(phi)}, even though modifying $c_{r_\phi}$ alters the $\delta$-component of $\fkk(y)$. Indeed, $\phi^{\square}(c_{r_\phi})$  commutes with $I_{\phi, \QQ_p}$, so changing $c_{r_\phi}$ (which multiplies its image in $P^{L_{r_\phi}}(\breve \QQ_p)$ by an element of $P^{L_{r_\phi}}(\Qpur)$) replaces $\delta$ by its $\sigma$-conjugate under an element of $G(\Qpur) = (\Res_{\Qpur/\QQ_p} G)(\QQ_p)$ that commutes with $\Int (y_p)^{-1}(I_{\phi, \QQ_p} ) = I_p(\fkk(y)) \subset \Res_{\Qpur/\QQ_p} G$. Thus the $\sigma$-centralizer $I_p(\fkk(y)) $ of $\delta$ remains unchanged.} and  the tuple $(I_\phi, (\iota_{y,v})_v)$ is a refinement of $\fkk(y)$.  	As in \S\ref{para:inclusionIx}, this refinement controls the left $I_\phi(\QQ)$-action on the right $G(\A_f^*)$-torsor $Y(\phi)$ in the following sense: for any $j \in I_\phi(\QQ)$, viewing $\iota_y(j) = (\iota_{y,v}(j))_v$ as an element of $G(\A_f^*)$, we have
 $$ j \cdot y = y \cdot \iota_y (j).$$   
 
 \subsection{}\label{para:functoriality of Y(J)}
 Let $\mathscr J$ be a conjugacy class of admissible morphisms $\Qf \to \G_G$. Then $\mathscr J$ is a small connected groupoid, where morphisms $\phi \to \phi'$ are given by elements $g \in G(\Qbar)$ such that $\phi' = \Int(g) \circ \phi$. The association $\JJ \ni \phi \mapsto Y(\phi)$ is naturally a right $G(\A_f^*)$-torsor over $\mathscr J$, from which we obtain $\bar Y(\JJ)$ that is canonically isomorphic to $ 
 I_\phi(\Q) \backslash Y(\phi)$ for any $\phi \in \JJ$. Since the value $r_\phi$ is constant on $\JJ$, the map (\ref{eq:Kott triple on Y(phi)}) descends to a 
 canonical map 
 \begin{align}\label{eq:Kott triple on bar Y(phi)}
 	\bar Y(\JJ) \To \KTstr,  \quad \bar y \longmapsto \fkk (\bar y) .
 \end{align}
 
 \subsection{}\label{subsec:map to pi}
 Let $\phi: \Qf \to \G_G$ be an admissible morphism.  Recall from \S \ref{subsubsec:recall of defns in LR} and \S \ref{subsubsec:SD} that $$X_p(\phi) =  X_{-\mu} (\phi(p) \circ \zeta_p) = \tilde X_{-\mu} (\phi(p) \circ \zeta_p) / \cG(\Z_p^{\ur}), $$ where $\mu$ is an element of $\dmu_X^{\cG}$ as in \S \ref{subsubsec:SD}, and $\tilde X_{-\mu} (\phi(p) \circ \zeta_p)$ is a subset of $\mathcal {UR} (\phi(p) \circ \zeta_p) = Y_p(\phi)$ which is independent of the choice of $\mu$. We define\index[Oct2026-sorted-notation]{y 0 00 00 y phi superscript natural@$Y(\phi)^{\natural}$}
 $$ Y(\phi)^{\natural} : = \tilde X_{-\mu} (\phi(p) \circ \zeta_p) \times Y^p(\phi) \subset Y(\phi). $$
 In particular, $X(\phi) = Y(\phi)^{\natural}/\cG(\Zpur)$.  
 
 Recall from \cite[\S\S  3.6--3.7]{kisin2012modp} that associated with $\phi$ we have a $\pi(G)$-torsor $\pi(G,\phi)$, together with a map\index[Oct2026-sorted-notation]{c 1 00 23 tilde c subscript phi@$\tilde c_{\phi}$}\index[Oct2026-sorted-notation]{zz16 pi 1 00 00 pi g phi@$\pi (G,\phi)$} $$\tilde c_{\phi}: Y(\phi)^{\natural} \To \pi (G, \phi). $$ For any $\tau \in I_{\phi}^{\ad} (\A_f)$, we write $I_{\phi}(\QQ)_{\tau}$ for $\Int(\tau) (I_{\phi}(\QQ)) \subset I_\phi (\A_f)$ as in \S \ref{subsubsec:defn of S(phi)}. The map $\tilde c_{\phi}$ enjoys the following properties:
 \begin{enumerate}
 	\item It factors through $I_{\phi}(\QQ)_\tau \backslash Y(\phi)^{\natural}$ for every $\tau \in I_{\phi}^{\ad}(\A_f)$ (see \cite[Lem.~3.7.4]{kisin2012modp}).
 	\item For $y,y'\in Y(\phi)^{\natural}$ and $r\in G(\A_f^*)$ with $y' = yr$, we have $\tilde c_\phi (y') = \tilde c_\phi(y) r$, where $r$ acts on the right hand side via its image in $\pi^*(G)$, which necessarily lies in $\pi(G)$. Indeed, this follows from the definition of $\tilde c_\phi$ and the observation in the proof of \cite[Lem.~3.6.2(2)]{kisin2012modp} that the coset $r \varrho(G_{\sconn}(\Qpur)) \cG(\Zpur) \subset G(\Qpur)$ contains an element of $G(\Q_p)$. 
 \end{enumerate} 
 
 Using these properties and arguing as in the proof of Proposition \ref{mark:1.1}, we conclude that $\tilde c_{\phi}$ extends uniquely to a $G(\A_f^*)$-equivariant map 
 $$ \tilde c_{\phi}: Y(\phi) \To \pi (G, \phi)\times^{\pi(G)} \pi^*(G) $$ which, moreover, factors through $I_{\phi}(\QQ)_\tau \backslash Y(\phi)$ for every $\tau \in I_{\phi}^{\ad}(\A_f)$.

 By \cite[Prop.~3.6.10]{kisin2012modp}, there is a canonical isomorphism of $\pi(G)$-torsors\index[Oct2026-sorted-notation]{zz08 theta 1 00 00 vartheta subscript phi@$\vartheta_{\phi}$}
 \begin{align}\label{eq:vartheta}
 	\vartheta_{\phi}:
 	\pi(G,\phi) \cong \pi (G,X).
 \end{align}
 We shall henceforth make this identification tacitly, thereby viewing $\tilde c_\phi$ as a map $Y(\phi) \to \pi^*(G,X)$.

 \subsection{}\label{subsubsec:setting for Y_T(phi)} Let $\fks = (T,i,h) \in \spd(G,X)$. From $\fks$ we obtain a morphism $\Psi_{T, \mu_h} : \Qf \to \G_T$ as in \S \ref{para:Psi_mu}, which we shall simply denote by\index[Oct2026-sorted-notation]{zz23 psi 0 00 00 psi subscript t h@$\Psi_{T,h}$} $\Psi_{T,h}$. 	We make the following definitions which are analogous to \S \ref{para:defn of Y(phi)}.\index[Oct2026-sorted-notation]{y 0 00 00 y psi subscript t h@$Y(\Psi_{T,h})$}
 \begin{align*}
 	Y_p(\Psi_{T,h})  & := \mathcal {UR}( \Psi_{T,h}(p)\circ \zeta_p) ,\\ Y^p(\Psi_{T,h}) & : = X^p(\Psi_{T,h}),\\ Y(\Psi_{T,h}) & := Y_p(\Psi_{T,h})\times Y^p(\Psi_{T,h}).
 \end{align*} By Lemma \ref{lem:unram criterion}, $Y_p(\Psi_{T,h})$ is a $T(\Qpur)$-torsor. The definition of $X^p(\Psi_{T,h})$ is as in \S \ref{subsubsec:defn of S(phi)}, but with $G$ replaced by $T$.  By \cite[Prop.~3.6.7]{kisin2012modp}, $X^p(\Psi_{T,h})$ is non-empty, and hence a $T(\A_f^p)$-torsor.

 Each $y \in Y_p(\Psi_{T,h})$ determines an element $b_y ^T \in T(\Qpur)$ such that the morphism $\D \to \G_T^{\ur}$ underlying the unramified morphism $\Int(y^{-1}) \circ \Psi_{T,h}(p) \circ \zeta_p$ maps $d_{\sigma}$ to $b_y^T \rtimes \sigma$, cf.~Definition \ref{defn:b_theta}. Define\index[Oct2026-sorted-notation]{y 0 00 00 y subscript p psi subscript t h superscript natural@$Y_p(\Psi_{T,h})^\natural$}  
 $$ Y_p(\Psi_{T,h})^\natural : = \set{y \in Y_p(\Psi_{T,h}) \mid w_{T_{\Q_p}} (b_y^T) = [-\mu_h] \in X_*(T)_{\Gamma_{p,0}}}.$$ By  \cite[Prop.~3.6.7]{kisin2012modp}, this is a $T(\Q_p)\cT^\circ(\Z_p^{\ur})$-torsor, where the product is taken inside $T(\Qpur)$. We set\index[Oct2026-sorted-notation]{y 0 00 00 y psi subscript t h superscript natural@$Y(\Psi_{T,h})^\natural$}
 \begin{align*}
 	Y(\Psi_{T,h})^\natural & := Y_p(\Psi_{T,h})^\natural \times X^p(\Psi_{T,h}) \subset Y(\Psi_{T,h}).
 \end{align*}
 Thus $Y(\Psi_{T,h})^{\natural}$ is a torsor under $T(\A_f) \cT^\circ(\Z_p^{\ur})$ (the product inside $T(\A_f^*)$).

 By Lemmas \ref{lem:cusp TFAE} and \ref{lem:G=G^c}, $T(\QQ)$ is discrete in $T(\A_f)$, and hence closed in $T(\A_f)$. 
 Hence by  \cite[Prop.~3.6.7]{kisin2012modp}, there is a canonical $T(\A_f)$-equivariant bijection
 $$ T(\QQ)  \backslash Y(\Psi_{T,h})^{\natural} /\cT^\circ(\Z_p^{\ur}) \isom T(\Q)  \backslash T(\A_f)/\cT^\circ(\Z_p).$$ We denote by\index[Oct2026-sorted-notation]{y 0 00 00 y psi subscript t h subscript ntr@$Y(\Psi_{T,h})_{\ntr}$} $Y(\Psi_{T,h})_{\ntr}$ the $T(\Q) \cT^\circ (\Z_p^{\ur})$-orbit in $Y(\Psi_{T,h})^{\natural}$ corresponding to the image of $1\in T(\A_f)$ in the right hand side.

 As in Definition \ref{defn:special morphisms}, we write $\phi(\fks) = \phi(T,i,h)$ for the morphism $i \circ \Psi_{T,h} : \Qf \to \G_G$. This morphism is admissible, as recalled in Theorem \ref{thm:kisin special}. There is a canonical injection
 \begin{align}\label{eq:inj from Y_T}
 	Y(\Psi_{T,h})  \hookrightarrow Y(\phi(\fks))\end{align}
 induced by $i$. Elements of the image of $Y (\Psi_{T,h})_{\ntr}$ under the above map play a role analogous to the integral special points in Definition \ref{defn:integral special point}.  Their fundamental properties  are given in the following proposition, which is an analogue of Proposition \ref{prop:special points on the geom side}. 
 
 \begin{prop}\label{prop:special pts on the gerb side} Let $y_T\in Y(\Psi_{T,h})_{\ntr}.$ Write  $y \in Y(\phi(\fks))$ for the image of $y_T$ under 
 	$(\ref{eq:inj from Y_T})$. The following statements hold.
 	\begin{enumerate}
 		\item There exists a unique element $\delta_T$ lying in the $\circsim$-equivalence class in $T(\Qpur)^{\mot}$ determined by $-\mu_h$ such that the element $\fkk(y) \in   \KTstr$ associated with $y$ as in (\ref{eq:Kott triple on Y(phi)}) is equal to $ \fkk(\fks, \delta_T)$.     Moreover, $\delta_T$ depends only on $(T,h)$ and $y_T$, not on $i: T \hookrightarrow G$.
 		\item  The image of $y$ under $\tilde c_{\phi(\fks)}: Y(\phi(\fks)) \to \pi^*(G,X)$  is equal to the image of $i\circ h \in X$ in $\pi(G,X) \subset \pi^*(G,X)$.  	 
 	\end{enumerate}
 \end{prop}
 \begin{proof} For part (i), the uniqueness is clear since $i(\delta_T) \in G(\Qpur)$ must be the $\delta$-component of $\fkk(y)$. We show the existence. Note that $r_{\Psi^{\Delta}_{T,h}} = r_{\phi(\fks)}$. Write $r$ for this number, which is independent of $i$. Write $\Psi^{\square}_{T,h}: P^{L_r} \to T$ for the unique factorization of $\Psi_{T,h}^{\Delta}$, which is defined over $\QQ$. Following the construction in \S \ref{para:defn of Y(phi)} (also cf.~\cite[\S 4.3.9]{kisin2012modp}), we see that $\fkk(y)$ is represented, for sufficiently divisible $n$, by  $$((i(\gamma_{0,T,n}))_{l\neq p},i(\delta_T)) \in \mathfrak T_n,$$   with $\gamma_{0,T,n} = \Psi_{T,h}^{\square}(\delta_n) \in T(\QQ)$ and  
 	\begin{align}\label{eq:delta_T}
 		\delta_T = c_T b_{y_{T,p}}^T \sigma(c_T)^{-1} \in T(\QQ_{p^n})
 	\end{align} where $c_T =  \Psi_{T,h}^{\square} (c_r) \in T(\LL)$. By (\ref{eq:def prop of c}) we have $$\gamma_{0,T,n} = \delta_T \sigma(\delta_T) \cdots \sigma^{n-1}(\delta_T) \in T(\Q) \subset T(\Q_p).$$   By the definition of $\delta_n$, we know that  $\gamma_{0,T,n}$ is a $p$-unit in $T(\Q)$. In particular, we have  $\delta_T\in T(\Qpur)^{\mot}$ and  $\fkk(y) = \fkk(\fks, \delta_T)$. It remains to show that $w_{T_{\QQ_p}} (\delta_T) = [-\mu_h] \in X_*(T)_{\Gamma_{p,0}}$. Since $c_r \in P^{L_r,\circ} (\breve \Z_p)$ and $\Psi_{T,h}^{\square}$ maps  $P^{L_r,\circ} (\breve \Z_p)$ into $\cT^{\circ}(\breve \Z_p)$, we have $c_T \in \cT^{\circ}(\breve \Z_p)$. Hence $w_{T_{\QQ_p}} (\delta_T) = w_{T_{\QQ_p}} (b^T_{y_{T,p}})$, and this is $[-\mu_h]$ 	
 	by the definition of $Y _p (\Psi_{T,h})^{\natural}$. The existence in part (i) is proved. 
 	The independence of $i$ follows from (\ref{eq:delta_T}).  
 	
 	Part (ii) follows from  \cite[Prop.~3.6.10 (2)]{kisin2012modp}.
 \end{proof}
 
 \section{Markings and amicable pairs}\label{subsec:markingsamic}
 
 In this section, we introduce \emph{amicable pairs}, the central concept of this chapter. These are pairs consisting of an isogeny class $\II \subset \Shh_{K_p}(\Fpbar)$ and a conjugacy class $\JJ$ of admissible morphisms $\Qf \to \G_G$ satisfying a compatibility condition with respect to strict Kottwitz tuples and connected components. Given an amicable pair $(\II, \JJ)$ and $\phi \in \JJ$, a choice of \emph{marking} yields a bijection $\II \isom S_\tau (\phi)$ for a corresponding element $\tau \in I_{\phi}^{\ad}(\A_f)/ I_{\phi}^{\ad}(\Q)$; see Proposition \ref{prop:marked bij}. While $\tau$ depends on the choice of marking, all such choices yield elements with the same image in the group $\fkH(\phi)$ introduced in \S \ref{sec:conn comp}. Consequently, we obtain a well-defined invariant in $\fkH(\phi)$ of the amicable pair; see Definition \ref{defn:matchable}.
 
 \subsection{}
 Let $\II$ be an isogeny class in $\Shh_{K_p}(\Fpbar)$, and $\mathscr J$ a conjugacy class of admissible morphisms $\gQ \to \gG_G$. For $x\in \II$ and $\phi \in \JJ$, we have maps $Y(x) \to \KTstr$ and $Y(\phi) \to \KTstr$ as in (\ref{eq:Kott triple on Y(x)}) and (\ref{eq:Kott triple on Y(phi)}), and maps $\II^* \to \pi^*(G,X)$ and $Y(\phi) \to \pi^*(G,X)$ as in Proposition \ref{mark:1.1} and \S \ref{subsec:map to pi}. Since $$\II^* = \bar Y(\II)/\cG(\Zpur) \cong I_x(\QQ)\backslash Y(x) /\cG(\Zpur), $$ we have an induced map $Y(x) \to \pi^*(G,X)$. We define the sets of\index{marking}\index{marking!weak} \emph{weak markings} and \emph{markings} by \index[Oct2026-sorted-notation]{mar 0 03 00 mar x phi superscript w@$\mathcal Mar(x,\phi)^w$}\index[Oct2026-sorted-notation]{mar 0 03 00 mar x phi@$\mathcal Mar(x,\phi)$} $$\tmar(x,\phi)^w : = Y(x) \times_{\KTstr} Y(\phi), \quad \tmar(x,\phi)  : = Y(x) \times_{\KTstr \times \pi^*(G,X)} Y(\phi) .$$  
 We say that the pair $(\II, \JJ)$ is \emph{weakly amicable}\index{amicable pair!weakly} (resp.~\emph{amicable}\index{amicable pair}) if $\tmar(x,\phi)^w \neq \emptyset$ (resp.~$\tmar(x,\phi) \neq \emptyset$). These conditions are independent of the choices of $x$ and $\phi$, since the map $Y(x) \to \KTstr \times \pi^*(G,X)$ factors through a map $\bar Y(\II) \to \KTstr \times \pi^*(G,X)$ independent of $x$, and likewise $Y(\phi) \to \KTstr \times \pi^*(G,X)$ factors through a map $\bar Y(\JJ) \to \KTstr \times \pi^*(G,X)$ independent of $\phi$. By abuse of language, we also say that the pair $(x,\phi)$ is weakly amicable (resp.~amicable). Moreover, both $\tmar(x,\phi)^w$ and $\tmar(x,\phi)$ are stable under the diagonal $G(\A_f^*)$-action on $Y(x) \times Y(\phi)$.

 \subsection{}\label{subsubsec:marking} 
 Let $(\II,\JJ)$ be a weakly amicable pair. Fix $x\in \II$ and $\phi \in \JJ$. Let  $\fkm = (y, y') \in \tmar(x,\phi)^w$, and let $\fkk = \fkk (y) = \fkk(y') \in \KTstr$.  Recall from \S \ref{para:inclusionIx} and \S \ref{para:defn of Y(phi)} that we have refinements  $(I_{x}, (\iota_{ y,v})_v)$ and $(I_\phi, (\iota_{y',v})_v)$ of $\fkk$. Hence by \S \ref{subsubsec:associated data}, there exist a $\QQ$-isomorphism $f: I_{x} \isom I_{\phi}$ and an element $\tau= (\tau_{v})_v \in I_{\phi}^{\ad}(\A_f)$ such that for each finite place $v$  we have $$\iota_{y, v} \circ f^{-1} = \iota_{y',v} \circ \Int (\tau_{v}) : I_{\phi, \Q_v} \overset{\sim}{\To} I_v(\fkk). $$ Clearly $f$ is uniquely determined by $\fkm$ up to $I_{\phi}^{\ad}(\QQ)$-conjugation, and the image \index[Oct2026-sorted-notation]{zz19 tau 1 00 00 tau subscript m@$\tau_{\mathfrak m}$}
 $$ \tau_{\fkm} \in I_{\phi}^{\ad}(\A_f)/ I_{\phi}^{\ad}(\Q) $$ of $\tau$  is determined by $\fkm \in \tmar(x,\phi)^w $ alone (independently of $f$). 
 
 Recall from \S \ref{subsec:frak H} that we have a surjective map of abelian groups \begin{multline*}
 	\fkH(\phi)  := \ker (I_{\phi}(\A_f) \to \pi(G)) \backslash I_{\phi}^{\ad}(\A_f) / I_{\phi}^{\ad}(\QQ) \\ \To  \mathcal H(\phi) := I_{\phi}(\A_f)   \backslash I_{\phi}^{\ad}(\A_f) / I_{\phi}^{\ad}(\QQ). 
 \end{multline*} 
 \begin{lem} \label{lem:tau for amicable pair} The $I_{\phi}^{\ad}(\QQ)$-conjugacy class of $f$ and the image of $\tau_{\fkm}$ in $\cH(\phi)$ are both independent of $\fkm \in  \tmar(x,\phi)^w$. If $(\II, \JJ)$ is amicable, then the images of $\tau_{\fkm}$ in $\mathfrak H(\phi)$ for all $\fkm \in \tmar(x,\phi)$ are all equal. Moreover, in this case every element of $ I_{\phi}^{\ad}(\A_f) / I_{\phi}^{\ad}(\Q)$ above this common image is realized as $\tau_{\fkm}$ for some $\fkm \in \tmar(x,\phi)$.  
 \end{lem}
 
 \begin{proof}For the first statement, using  the $G(\A_f^*)$-action we easily reduce to comparing $\fkm_1,\fkm_2\in  \tmar(x,\phi)^w$ with common second coordinate. Write $\fkm_i = (y_i, y')$, and write $\fkk$ for $\fkk(y') = \fkk(y_1) = \fkk(y_2)$.  		Let $h = (h_v)_v \in G(\A_f^*)$ be such that 
 	$ y_1 =  y_2 h$. Then $h_v\in I_v(\fkk)(\Q_v)$ by the $G(\A_f^*)$-equivariance of (\ref{eq:Kott triple on Y(x)}). Let  $s_v  =  (\iota_{y',v})^{-1}(h_v) \in  I_{\phi}(\Q_v)$. The element $s = (s_v)_v$ lies in $I_{\phi}(\A_f)$, since $s_v = \Int(y'_v)(h_v)$ and $y^{\prime p} = (y'_v)_{v\neq p} \in X^p(\phi)$ is an element of $G(\bar \A_f^p)$ such that $\Int (y^{\prime p })^{-1}$ is an $\A_f^p$-map $I_{\phi, \A_f^p} \to G_{\A_f^p}$.  Now suppose we have a $\Q$-isomorphism $f_1: I_x \isom I_{\phi}$  and an element  $\tau_1 = (\tau_{1,v})_v  \in I_{\phi}^{\ad}(\A_f)$  satisfying $$\iota_{y_1,v} \circ f_1^{-1} = \iota_{y',v} \circ \Int (\tau_{1,v}) : I_{\phi, \Q_v} \To I_v(\fkk). $$
 	Then we have
 	\begin{align}\label{eq:f_2=f_1}
 		\iota_{y_2,v} \circ f_1^{-1} = \Int(h_v) \circ \iota_{ y_1,v} \circ f_1^{-1} =  \iota_{y',v} \circ \Int (s_v\tau_{1,v}) : I_{\phi, \Q_v} \To I_v(\fkk).
 	\end{align}It follows that  $ \tau_{\fkm_2} = s \tau_{\fkm_1} \in I_{\phi}^{\ad}(\A_f)/ I_{\phi}^{\ad}(\Q)$, and that $f_1$ is still associated with $\fkm_2$. This proves the first statement.

 	To prove the second statement, it suffices to show that $s$ has trivial image in $\pi(G)$ under the map (\ref{eq:map to pi(G)}) under the assumption  that $\fkm_1, \fkm_2 \in \tmar(x,\phi)$. This assumption implies that the image of $h$ in $\pi^*(G) \supset \pi(G)$ is trivial. Thus it suffices to show that the image $s_G$ of $s$ under the composite map  $$I_{\phi}(\A_f) \xrightarrow{(\ref{eq:map betw H0ab})} G_{\sconn}(\A_f) \backslash G(\A_f) \To G_{\sconn}(\A_f^*) \backslash G(\A_f^*) $$ is equal to the image of $h$. 
 	
 	Let $\eta:\pi_1(I_\phi) 
 	\to \pi_1(G)$ be the map induced by the inclusion $I_{\phi, \Qbar} \hookrightarrow G_{\Qbar}$. 
 	For each finite place $v\neq p$, the $\Q_v$-map $\iota_{y',v} = \Int(y_v^{\prime -1} ) : I_{\phi, \Q_v} \to G_{\Q_v}$  induces $\eta$ at the $\pi_1$-level. Hence  the component of $s_G$ in $ G_{\sconn}(\A_f^p)\backslash G(\A_f^p)$ is equal to the image of $h^p \in G(\A_f^p)$. 
 	Fix $n \geq 1$ such that $I_p(\fkk)$ is a subgroup of $\Res_{\Q_{p^n}/\Q_p} G$, and fix $\tilde \sigma \in \Gal(\Qpbar/\Q_p)$ lifting the Frobenius.   
 	Consider the following two composite maps:
 	\begin{align*}
 		q_1 &: I_{\phi , \Q_p} \xrightarrow{\iota_{y',p}} I_p\hookrightarrow \Res _{\Q_{p^n}/\Q_p} G,\\
 		q_2 & :  I_{\phi, \Qpbar} \hookrightarrow G_{\Qpbar} \xrightarrow{\text{diagonal}}  (\Res_{\Q_{p^n}/\Q_p} G)_{\Qpbar}\end{align*} 
 	They are both $\Q_p$-rational at $\pi_1$-level, and hence induce maps \begin{multline*}
 		q_{i, *} : \coh^0_{\ab}(\Q_p, I_{\phi}) \cong I_{\phi,\sconn}(\Q_p) \backslash I_{\phi}(\Q_p) \\ \To \coh^0_{\ab}(\Q_p, \Res_{\Q_{p^n}/\Q_p} G) \cong  G_{\sconn}(\Q_{p^n}) \backslash G(\Q_{p^n}).
 	\end{multline*} Note that $q_{1,*}(s_p)$ is the image of $h_p \in G(\QQ_{p^n})$, while  $q_{2,*}(s_p)$ maps to the component of $s_G$ in $G_{\sconn}(\Qpur) \backslash G(\Qpur)$. To finish the proof, it suffices to show that $q_1$ and $q_2$ induce the same map at $\pi_1$-level.  		Write $j:I_{\phi,\Qpbar}\hookrightarrow G_{\Qpbar}$ for the
 	natural inclusion, and identify
 	$
 	(\Res_{\Q_{p^n}/\Q_p}G)_{\Qpbar}
 	\cong\prod_{i=0}^{n-1}G_{\Qpbar}$.  The first component of $q_1$ is
 	$\Int(y_p^{\prime -1})\circ j$, and the other components are its
 	Galois conjugates. Since $j$ is a pre-inner twisting, all these
 	components induce the same map
 	$\eta:\pi_1(I_\phi)\to\pi_1(G)$ introduced above. 
 	Each component of $q_2$ also induces $\eta$. This concludes the proof of the second statement.

 	We now prove the last statement in the lemma. Let $\fkm_1 = (y_1, y_1') \in \tmar(x,\phi)$. Then an arbitrary element of $ I_{\phi}^{\ad}(\A_f) / I_{\phi}^{\ad}(\Q)$ above the   image of $\tau_{\fkm_1}$ in $\fkH(\phi)$ is of the form  $s \tau_{\fkm_1}$ with $s \in \ker (I_{\phi}(\A_f) \to \pi(G))$. For each finite place $v$, let $h_v = \iota_{y_1',v} (s_v) \in I_v(\fkk(y_1'))(\QQ_v)$. Then $h = (h_v)_v \in \ker (G(\A_f^*) \to \pi^*(G))$. Let $\fkm' = (y_1 h^{-1}, y_1')$. Then  we have $\fkm' \in \tmar(x,\phi)$ and, by the computation in the proof of the first statement, $\tau_{\fkm'} = s \tau_{\fkm_1}$. 
 \end{proof}
 
 Let $\phi_1$ and $\phi_2$ be conjugate admissible morphisms $\gQ \to \G_G$. Then we have  canonical isomorphisms of abelian groups $\cH(\phi_1) \cong \cH(\phi_2)$ and  $\fkH(\phi_1) \cong \fkH(\phi_2)$   induced by $\Int (g)$ for any $g\in G(\Qbar)$ conjugating $\phi_1$ to $\phi_2$. These are special cases of the isomorphisms discussed in \S \ref{subsubsec:cocycle relation} and \S \ref{subsec:frak H}, as $\phi_1 \approx \phi_2$.
 \begin{prop}\label{prop:well def of tau}	Let $(\II, \JJ)$ be a weakly amicable pair. For $i=1,2$, let $x_i\in \II, \phi_i \in \JJ, \fkm_i\in \tmar(x_i,\phi_i)^w$. \begin{enumerate}
 		\item The canonical isomorphism $\cH(\phi_1) \cong \cH(\phi_2)$ sends the image of  $\tau_{\fkm_1}$ to the image of $\tau_{\fkm_2}$. 
 		\item If $(\II, \JJ)$ is amicable and $\fkm_i\in \tmar(x_i,\phi_i)$ for $i=1,2$, then the canonical isomorphism $\fkH(\phi_1) \cong \fkH(\phi_2)$ sends the image of  $\tau_{\fkm_1}$ to the image of $\tau_{\fkm_2}$.
 	\end{enumerate} 
 \end{prop}  
 \begin{proof}We may assume either $\phi_1 = \phi_2$ or $x_1 = x_2$; the general case follows by combining these two cases.
 	
 	First assume $\phi_1 = \phi_2$, and denote it by $\phi$. 	Choose an isomorphism $f: x_1 \to x_2$ in the category $\II$. This induces an identification $Y(x_1) \isom Y(x_2)$ and hence a bijection $\alpha: \tmar(x_1,\phi)^w \isom \tmar(x_2,\phi)^w$, which maps $\tmar(x_1,\phi)$ onto $\tmar(x_2,\phi)$. By Lemma \ref{lem:tau for amicable pair}, we may assume that $\fkm_2 = \alpha (\fkm_1)$. Then $\fkm_1,\fkm_2$ are of the form $\fkm_i = (y_i , y')$, and for each finite place $v$ the isomorphism $\iota_{y_1, v} : I_{x_1, \Q_v} \isom I_v(\fkk(y'))$ is equal to $\iota_{y_2, v} : I_{x_2, \Q_v} \isom I_v(\fkk(y'))$ pre-composed with the $\QQ$-isomorphism $I_{x_1} \isom I_{x_2} $ induced by $f$. It follows that $\tau_{\fkm_1} = \tau _{\fkm_2} \in I_{\phi}^{\ad}(\A_f)/I_{\phi}^{\ad}(\Q)$. 
 	
 	The case $x_1  = x_2$ is proved in the same way by choosing  $g\in G(\Qbar)$ such that $\phi_2 = \Int(g ) \circ  \phi_1$ and using the induced bijection $\tmar(x, \phi_1)^w \isom \tmar(x,\phi_2)^w$. 
 \end{proof}
 \begin{defn}\label{defn:matchable}Let $(\II, \JJ)$ be an amicable pair. Let $\phi \in \JJ$. We define\index[Oct2026-sorted-notation]{zz19 tau 1 00 00 tau superscript h i j@$\tau^{\mathfrak H}(\mathscr I, \mathscr J)$}   $$\tau^{\fkH} (\II, \JJ) \in \fkH(\phi)$$ to be the image of $\tau_{\fkm}$ for arbitrary $x \in \II$ and $\fkm \in \tmar(x,\phi)$. This is  well defined by Proposition \ref{prop:well def of tau}, and moreover it is independent of  $\phi$  if we canonically identify $\fkH(\phi)$ for all $\phi \in \JJ$. 
 \end{defn}
 
 Let $\phi: \Qf \to \G_G$ be an admissible morphism, and let $\tau \in I_{\phi}^{\ad}(\A_f)$. Recall from \S \ref{subsubsec:defn of S(phi)} that $ S_{\tau}(\phi) : = \varprojlim_{K^p} I_\phi(\QQ)_{\tau}\backslash X(\phi)/ K^p$.   
 The following lemma holds in our current setting of Hodge type. 
 \begin{lem}\label{lem:HT free}
 	For any neat compact open subgroup $K^p \subset G(\A_f^p)$, the action of $I_\phi (\Q)_{\tau}$ on $X(\phi)/K^p$ is free. Moreover, the natural map  
 	$$ I_{\phi}(\Q)_{\tau}\backslash  X(\phi) \To S_\tau(\phi)$$  is a bijection.
 \end{lem}
 \begin{proof}
 	Suppose $\gamma \in I_{\phi} (\Q)_{\tau}$ has a fixed point in $X(\phi)/K^p$. By Lemma \ref{lem:check app licable} (i), we have $\gamma \in Z_G(\Q) \cap K_p K^p$. By Lemmas \ref{lem:Z_s} and \ref{lem:G=G^c}, we have $Z_G(\Q) \cap K_p K^p = \set{1}$. Hence $\gamma =1$. This proves the first part.  
 	
 	It follows that the natural surjection
 	$$  I_{\phi}(\Q)_{\tau} \backslash  \bigg (\varprojlim_{K^p} X(\phi)/K^p \bigg)  \To S_\tau(\phi)$$ is a bijection. Since $X^p(\phi)$ is a torsor under the locally profinite group $G(\A_f^p)$, the natural map $X(\phi)\to \varprojlim_{K^p} X(\phi)/K^p$ is a bijection. This proves the second part. 
 \end{proof}

 \subsection{}\label{subsubsec:markingII} Let $(\II,\JJ)$ be a weakly amicable pair, $x\in \II, \phi \in \JJ$. For each $\fkm = (y,y') \in \tmar(x,\phi)^w$, we shall   construct a bijection 
 $$ \II \isom S_{\tau_{\fkm}}(\phi).$$ Here the right hand side is well defined since $\tau_{\fkm}$ is a well-defined element of $I_{\phi}^{\ad}(\A_f)/I_{\phi}^{\ad}(\Q)$. The construction is essentially the same as the composition of the bijections in \cite[Prop.~4.4.14, Cor.~4.5.4]{kisin2012modp}, except that the element $c c'^{-1}$ used to define $y_{\phi}$ in \cite[Cor.~4.5.4]{kisin2012modp} is replaced by our $c_{r_\phi}$; see \S \ref{para:defn of Y(phi)}.

 Let $\xi:Y(x) \isom G(\A_f^*)$ (resp.~$\xi' : Y(\phi) \isom G(\A_f^*)$) be the unique $G(\A_f^*)$-equivariant bijection  sending $y$ (resp.~$y'$) to $1$. 
 Recall from \S \ref{para:defn of Y(phi)} that the component $\delta$ of $\fkk(y) = \fkk(y')$ is given by 
 $$ \delta = c_{y'_p} b_{y'_p} \sigma (c_{y'_p})^{-1} .$$  Let $\mathbf f: G(\LL)/ \cG(\breve \ZZ_p) \to G(\LL)/ \cG(\breve \ZZ_p)$ be the bijection given by left-multiplication by $c^{-1}_{y'_p}$. Via the natural bijection $G(\Qpur)/\cG(\Z_p^{\ur}) \cong G(\LL)/ \cG(\breve \ZZ_p)$ (by the density of $G(\Qpur)$ in $G(\LL)$), we view $\mathbf f$ as a bijection $G(\Qpur)/\cG(\Z_p^{\ur}) \to  G(\Qpur)/\cG(\Z_p^{\ur})$.

 \begin{prop}\label{prop:marked bij} 
 	The composite bijection 
 	$$ Y(x)/\cG(\Z_p^{\ur}) \underset{\sim}{\overset{\xi}{\To}} G(\A_f^*)/\cG(\Z_p^{\ur})  \underset{\sim} {\xrightarrow{(\id_{G(\A_f^p)},~ \mathbf f)} } G(\A_f^*)/\cG(\Z_p^{\ur})  \underset{\sim} { \xrightarrow{\xi^{\prime -1}}  } Y(\phi)/ \cG(\Z_p^{\ur}) $$ descends to a bijection
 	$$
 	\mathscr B_{\fkm}: 	\II^* \isom I_{\phi}(\Q)_{\tau_\fkm} \backslash  Y(\phi)/\cG(\Z_p^{\ur}).
 	$$  Moreover, $\mathscr B_{\fkm}$ restricts to a bijection $ \II  \isom S_{\tau_\fkm}(\phi)$  which is compatible with the actions of $G(\A_f^p)$ and the $q$-Frobenius $\Phi$ on the two sides. (Here $q$ is the residue cardinality of the fixed prime $\fkp$ of $E$ over $p$.) If $\mathfrak m = (y,y')$ lies in $\tmar(x,\phi)$, then $\mathscr B_{\fkm}$ (resp.~$\mathscr B_{\fkm}|_{\mathscr I}$) commutes with the maps from the two sides to $\pi^*(G,X)$ (resp.~$\pi(G,X)$). 
 \end{prop}
 \begin{proof}Write $\iota$ for the composite map $I_x(\Q) \hookrightarrow I_x(\A_f) \xrightarrow{(\iota_{y,v})_v} G(\A_f^*).$ The map $\xi$ induces a bijection
 	$$ \bar \xi: I_x(\Q) \backslash Y(x)/ \cG(\Z_p^{\ur}) \isom \iota (I_{x}(\Q)) \backslash G(\A_f^*)/ \cG(\Z_p^{\ur}). $$ Similarly, $\xi'$ induces a bijection
 	$$ \bar \xi': I_\phi(\Q)_{\tau_\fkm} \backslash Y(\phi) / \cG(\Z_p^{\ur})  \isom \iota (I_{x}(\Q)) \backslash G(\A_f^*) / \cG(\Z_p^{\ur}) , $$ since the image of $I_{\phi}(\Q)_{\tau_\fkm} \hookrightarrow I_{\phi}(\A_f) \xrightarrow{(\iota_{y',v})_v} G(\A_f^*)$ is equal to $\iota(I_x(\Q))$. To prove the first part of the proposition, it suffices to show that $\mathbf f(hg) = h  \mathbf f(g)$ for all $g\in G(\Qpur)/ \cG(\Z_p^{\ur})$ and all $h\in I_p(\Q_p)$. This holds since $c_{y'_p}$ commutes with   $I_p(\Q_p)$ (because $\phi^{\square}(c_{r_\phi})$ commutes with $I_{\phi, \Q_p}$).
 	
 	We now show that $\mathscr B_{\fkm}(\II) = S_{\tau_\fkm}(\phi)$, where $S_{\tau_\fkm}(\phi)$ is viewed as a subset of $I_{\phi}(\Q)_{\tau_\fkm} \backslash  Y(\phi)/\cG(\Z_p^{\ur})$ via Lemma \ref{lem:HT free}. By Proposition \ref{prop:isog class reformulated},   we have 
 	\begin{align}\label{eq:model 1}
 		\bar \xi (\II) = \iota(I_x(\QQ)) \backslash X_{-\mu} (\delta) \times G(\A_f^p),
 	\end{align} where $\mu \in \dmu_X^{\cG}$. By  \S \ref{subsubsec:recall of defns in LR} and Lemma \ref{lem:HT free}, we have 
 	\begin{align}\label{eq:model 2}
 		\bar \xi' (S_{\tau_\fkm}(\phi))= \iota(I_x(\QQ)) \backslash X_{-\mu} (b_{y'_p}) \times G(\A_f^p).
 	\end{align}  Hence the claim boils down to $ \mathbf f(
 	X_{-\mu}(\delta))  =  X_{-\mu}(b_{y'_p})$, which is clear.
 	
 	The compatibility of $\mathscr B_{\fkm}$ with $G(\A_f^p)$ is clear. The compatibility with $\Phi$ boils down to the following three compatibilities. Firstly, the bijection from $\II$ to (\ref{eq:model 1}) induced by $\bar \xi$ takes the operator $\Phi$ on $\II$ to the operator $(\delta \rtimes \sigma )^{r}$ on  $X_{-\mu}(\delta)$ (with $r = [\FF_q: \FF_p]$). This bijection is just the inverse of (\ref{eq:presentation of II}), which is the map \cite[(2.1.4)]{kisin2012modp} (modulo an identification, see Remark \ref{rem:isog class reformulated}), so the compatibility follows from \cite[Prop.~1.4.4]{kisin2012modp}, cf.~\cite[Cor.~1.4.13, Prop.~2.1.3, Prop.~4.4.14]{kisin2012modp}. Secondly, the bijection $\mathbf f: 	X_{-\mu}(\delta) \isom X_{-\mu}(b_{y'_p}) $ takes	$(\delta \rtimes \sigma )^{r}$ on the left hand side to $(b_{y'_p} \rtimes \sigma)^r$ on the right hand side, which is clear. Thirdly, the bijection from $S_{\tau_\fkm}(\phi)$ to (\ref{eq:model 2}) induced by $\bar \xi' $ takes $\Phi$ on $S_{\tau_\fkm}(\phi)$ to $(b_{y'_p} \rtimes \sigma)^r$ on $X_{-\mu}(b_{y'_p})$, which follows from the discussion in \S \ref{subsubsec:recall of defns in LR}.
 	
 	Finally, suppose $y$ and $y'$ have the same image in $\pi^*(G,X)$. We show that $\mathscr B_{\fkm}$ is over $\pi^*(G,X)$. Suppose $\mathscr B_{\fkm}(t) =s$. Let $t$ be represented by $yg \in Y(x)$ with $g = (g^p, g_p)\in G(\A_f^*)$. Choose $g_{0} \in G(\Qpur)$ sufficiently close to $c^{-1}_{y'_p}$. Then $\mathbf f (g_p) = g_{0}  g_p  \cG(\Zpur)$, and so $s$ is represented by $y' (g^p, g_0 g_p) \in Y(\phi)$. Since $y$ and $y'$ have the same image in $\pi^*(G,X)$, the images of $t$ and $s$ in $ \pi^*(G,X)$ differ by translation by the image of $g_0$ in $\pi^*(G)$.  Since $c_{y'_p}$ is the image under a $\Qpur$-group homomorphism of a $\breve \ZZ_p$-point of the connected N\'eron model of a torus, functoriality gives $w_G(c_{y'_p}) = 0$, and
 	hence $w_G(g_0) = 0$. It then follows from Lemma \ref{lem:ker of Kott} that $g_0$ has trivial image in $\pi^*(G)$. Hence $s$ and $t$ have the same image in $\pi^*(G,X)$.  
 \end{proof}  
 \section{Amicable special forks}\label{subsec:gauges} 
 In this section, we introduce and study \emph{amicable special forks}. The main results are Propositions \ref{prop:LD1} and \ref{prop:LD2},  which are arguably the most technical in this chapter. They should be viewed as intermediate steps that will be reformulated in greater generality in \S\ref{subsec:Galcohprop}.

 \begin{defn}\label{defn:LD}
 	A\index{special fork} \emph{special fork} is an ordered pair $(\fks, \fks')$ consisting of two special point data $\fks , \fks ' \in \spd(G,X)$ of the form $\fks = (T,i,h)$ and $\fks' = (T,i',h)$, satisfying the following conditions:
 	\begin{enumerate}
 		\item The points $i\circ h$ and $i'\circ h$ lie in the same connected component of $X$.
 		\item The maps $i: T \to G$ and $i': T \to G$ are conjugate by $G^{\ad}(\Qbar)$.
 	\end{enumerate}To make the defining data explicit, we also write $(T,h, i , i')$ for a special fork as above.  
 \end{defn}

 \subsection{}\label{subsubsec:T^der for prelinking}
 Given a special fork $(T,h,i,i')$, the two maps $T \rightarrow G^{\ab}$ induced by $i$ and $i'$ are equal, and we denote their kernel by $T^{\dagger}$.  
 The two maps $ \Sha^{\infty} (\Q, T^{\dagger}) \to \Sha^{\infty}(\Q, G_{\der})$ induced by $i$ and $i'$ are equal (since $\Sha^{\infty}(\Q, G_{\der}) \cong \Sha^{\infty}_{\ab}(\Q, G_{\der})$ and since $i,i'$ are conjugate by $G^{\ad}(\Qbar)$). We denote the kernel by $\Sha^{\infty}_{G_{\der}} (\Q, T^{\dagger})$. Similarly, we define $\Sha_{G}^{\infty}(\Q, T)$ to be the kernel of $\Sha^{\infty}(\Q, T) \to \Sha^{\infty}(\Q, G)$ induced by either $i$ or $i'$.
 
 Write $T'$ for $i'(T) \subset G$ and write  $T'^\dagger $ for $T' \cap G_{\der} = i'(T^\dagger)$. Since $T'$ is its own centralizer in $G$, for any $g\in G_{\der}(\Qbar)$ such that $\Int(g) \circ i = i'$ the cocycle $(g \lix^{\rho} g^{-1})_{\rho \in \Gamma}$ defines an element  $\alpha_{i,i'} \in \coh^1(\Q,  T'^\dagger)$ which is independent of the choice of $g$.

 \begin{lem}\label{lem:trivial at infty}The element  $\alpha_{i,i'} $ lies in $\Sha^{\infty} _{G_{\der}} (\Q, T'^{\dagger})$.
 \end{lem}
 
 \begin{proof}The only non-trivial condition to check is that $\alpha_{i,i'}$ has trivial image in $\coh^1(\RR, T'^\dagger)$. The argument is similar to the proof of \cite[Prop.~4.4.13]{kisin2012modp}. Let $K_{\infty}'$ be the $\RR$-algebraic group that is the stabilizer of $i' \circ h$ in $G_{\der , \RR}$. Since $i\circ h$ and $i'\circ h$ lie in the same connected component of $X$, they are conjugate by an element of $G_{\der} (\RR)$ (or even $G_{\sconn} (\RR)$). It follows that $\alpha_{i,i'}$ has trivial image in $\coh ^1 (\RR, K_{\infty}')$. By \cite[Lem.~4.4.5]{kisin2012modp} applied to $H' = T'^{\dagger}$ and $H = K_{\infty}'$, every non-trivial element of $\coh ^1(\RR, T'^{\dagger})$ has a non-trivial image in $\coh ^1 (\RR, K_{\infty}')$. Hence $\alpha_{i,i'}$ must have trivial image in $\coh ^1 (\RR, T'^{\dagger})$, as desired. \end{proof}
 \begin{defn}
 	Let $\fks = (T,i,h) \in \spd(G,X)$. We write $\II_{\fks}$ for the isogeny class of $x_{\fks} \in \Shh_{K_p} (\Fpbar)$, and write $\JJ_{\fks}$ for the conjugacy class of the admissible morphism $\phi(\fks) : \Qf \to \G_G$. We say that a special fork $(\fks, \fks')$ is\index{special fork!amicable}  \emph{amicable} if $(\II_{\fks}, \JJ_{\fks'})$ is an amicable pair. 
 \end{defn}

 \begin{defn}An\index{special fork!extended} \emph{extended special fork} is a tuple $\fkg = (T,h,i,i', y_T,y'_T),$ where
 	\begin{itemize}
 		\item  $ (T,h,i,i')$ is a special fork,
 		\item $y_T$ is an element of the $\cT^{\circ}(\Z_p^{\ur})$-torsor $Y(\tomega_{(T,h)})^\circ$ (see \S \ref{para:tomega_T,h}),
 		\item $y'_T$ is an element of the $T(\Q) \cT^{\circ}(\Zpur)$-torsor $Y(\Psi_{T,h})_{\ntr} $ (see \S \ref{subsubsec:setting for Y_T(phi)})
 	\end{itemize} such that the element of $T(\Qpur)^{\mot}$ associated with $y_T$ as in Proposition \ref{prop:special points on the geom side} (i) is also associated with $y'_T$ as in
 	Proposition \ref{prop:special pts on the gerb side} (i); we denote this element by $\delta_{\fkg}$.  
 \end{defn}
 
 \begin{lem}\label{lem:getting gauge}
 	Let $(T,h)$ be a torus Shimura datum, and let $y'_T \in Y(\Psi_{T,h})_{\ntr}$. Then there exists $y _T\in Y(\tomega_{(T,h)})^{\circ}$ such that for every special fork for $(G,X)$ of the form $(T,h,i,i')$,  the tuple $(T,h,i, i', y_T, y'_T)$ is an extended special fork.
 \end{lem}
 \begin{proof} Let  $\delta \in T(\Qpur)^{\mot}$ be associated with $y_T'$ as in Proposition \ref{prop:special pts on the gerb side} (i). It lies in the $\circsim$-equivalence class in $T(\Qpur)^{\mot}$ determined by $-\mu_h$. By Proposition \ref{prop:special points on the geom side} (i), the image of $Y(\tomega_{(T,h)})^{\circ}$ under the map $Y(\tomega_{(T,h)})  \to T(\Qpur), z \mapsto \delta_z$ is precisely this  $\circsim$-equivalence class. Hence we can find $y_T \in Y(\tomega_{(T,h)})^{\circ}$ mapping to $\delta$.   
 \end{proof}
 
 \begin{defn}\label{defn:gauge rec} Let $\fkg = (T,h,i,i',y_T, y'_T)$ be an extended special fork. Then $\delta_{\fkg}$ lies in the $\circsim$-equivalence class in $T(\Qpur)^{\mot}$ corresponding to $-\mu_h$. Hence as in \S \ref{para:special kott trip} we obtain 
 	\begin{align*}
 		\fkk _{\fkg} :  = \fkk(T,i,h, \delta_{\fkg}) \in \KTstr  , \qquad
 		\fkk'_{\fkg}   : = \fkk(T,i',h, \delta_{\fkg}) \in \KTstr .
 	\end{align*}
 	A\index{rectification} \emph{rectification of $\fkg$} is an element $u \in \ker (G(\A_f^*) \to \pi^*(G))$ such that $ \fkk_{\fkg} \cdot u  =  \fkk_{\fkg}' .$  
 \end{defn}

 \begin{lem}\label{lem:rec gauge}
 	An extended special fork $\fkg= (T,h,i,i', y_T, y'_T)$ admits a rectification if and only if the underlying special fork $(T,h,i,i')$ is amicable. Moreover, each rectification $u$ gives rise to a marking\index[Oct2026-sorted-notation]{m 1 04 00 m g u@$\mathfrak m(\mathfrak g,u)$} $$\fkm(\fkg,u) = ( y u, y') \in \tmar(x_{(T,i,h)}, \phi(T,i',h)), $$ where $y$ is the image of $y_T$ in $Y(x_{(T,i,h)})$ and $y'$ is the image of $y'_T$ in $Y(\phi(T,i',h))$. 
 \end{lem}
 \begin{proof}
 	By Proposition \ref{prop:special points on the geom side} (iii), Proposition \ref{prop:special pts on the gerb side} (ii), and our assumption that $i\circ h$ and $i'\circ h$ lie on the same connected component of $X$, we know that $y$ and $y'$ have the same image in $\pi^*(G,X)$. Therefore $(x_{(T,i,h)},\phi(T,i',h))$ is amicable if and only if there exists $u \in \ker (G(\A_f^*) \to \pi^*(G))$ such that $\fkk(y) u = \fkk(y')$. But we have $\fkk (y) = \fkk_{\fkg}$ by Proposition \ref{prop:special points on the geom side} (ii) and $\fkk(y') = \fkk'_{\fkg}$ by Proposition \ref{prop:special pts on the gerb side} (i), so the condition on $u$ is precisely that it be a rectification. When this is the case, $(yu, y')$ is a marking. 
 \end{proof}

 \begin{cor}[cf.~Proposition \ref{introprop:amicable} in the Introduction]\label{cor:standard amicable}
 	Let $\fks \in \spd(G,X)$. Then $(\II_{\fks}, \JJ_{\fks})$ is an amicable pair.
 \end{cor}
 \begin{proof}
 	By Lemma \ref{lem:getting gauge}, we can extend the special fork $(\fks, \fks)$ to an extended special fork $\fkg$. Clearly $u=1$ is a rectification of $\fkg$. Therefore $(\fks, \fks)$ is amicable by Lemma \ref{lem:rec gauge}.
 \end{proof}

 \begin{defn}\label{defn:adapted iso}
 	Let $(\fks, \fks ' ) = (T,h, i,i')$ be an amicable special fork. We say that a $\Q$-isomorphism $f: I_{x_\fks} \isom I_{\phi(\fks')}$ is\index{i,i'-compatible@$(i,i')$-compatible} \emph{$(i,i')$-compatible} if it belongs to the canonical $I_{\phi(\fks')}^{\ad}(\Q)$-orbit of isomorphisms $I_{x_\fks} \isom I_{\phi(\fks')}$ as in Lemma \ref{lem:tau for amicable pair}, and if it commutes with the following two maps:
 	\begin{itemize}
 		\item the natural $\Q$-embedding $T \hookrightarrow I_{\phi(\fks')}$ (whose composition with $I_{\phi(\fks'), \Qbar} \hookrightarrow G_{\Qbar}$ is $i'$),
 		\item the natural $\Q$-embedding $T \hookrightarrow I_{x_{\fks}}$  as in (\ref{eq:T to I_x}). 
 	\end{itemize} 
 \end{defn}
 \subsection{}\label{para:setting for comm diags}
 Let $ (\fks, \fks')= (T,h, i,i')$ be an amicable special fork. Write $x$ and $\phi$ for $x_\fks$ and $\phi(\fks')$. By Lemma \ref{lem:getting gauge}, we extend $(\fks, \fks')$ to an extended special fork $\fkg$. By Lemma \ref{lem:rec gauge}, we choose a rectification $u$ and obtain $\fkm (\fkg, u) = (yu, y') \in \tmar(x,\phi). $ For $v$ a finite place, write $I_v$ for $I_v(\fkk(y')) = I_v(\fkk(yu))$, and write $I_v^{\natural}$ for $I_v(\fkk(y))$. Write
 \begin{align*}
 	\iota_v  = \iota_{yu, v} & :  I_{x, \QQ_v} \isom I_v \\
 	\iota'_v = \iota_{y',v} &: I_{\phi,\QQ_v} \isom I_v \\
 	\iota_v^{\natural} = \iota_{y,v} &:  I_{x,\Q_v} \isom I_v^{\natural}. 
 \end{align*}
 We set 
 $$ R_v : = \begin{cases}
 	G_{\QQ_v}, & v \neq p,\\
 	\Res_{\Qpur/\QQ_p} G, & v =p. 
 \end{cases}$$ Then $I_v$ and $I_v^{\natural}$ are $\Q_v$-subgroups of $R_v$. 	Denote the natural embeddings $T \hookrightarrow I_{x}$ and $T \hookrightarrow I_{\phi}$ by $j$ and $j'$ respectively. 
 
 \begin{lem}\label{lem:various comm diags}
 	In the setting of \S \ref{para:setting for comm diags}, the following statements hold.
 	\begin{enumerate}
 		\item We have commutative diagrams  
 		\begin{equation*} \xymatrix{ & I_{x,\Q_v} \ar[r]^{\iota_v^{\natural}} & I_v^{\natural}  \ar @{^(->}[d]  \\ T_{\Q_v} \ar[ru]^{j} \ar[r]^{i} & G_{\Q_v} \ar@{^(->}[r]^{\Delta} &  R_v } \quad 
 			\xymatrix{ & I_{\phi,\Q_v} \ar[r]^{\iota_v'} & I_v \ar @{^(->}[d]  \\ T_{\Q_v} \ar[ru]^{j'} \ar[r]^{i'} & G_{\Q_v} \ar@{^(->}[r] ^{\Delta} &  R_v }   \end{equation*} 	Here $\Delta$ is the identity if $v\neq p$ and the diagonal embedding if $v = p$.
 		\item Suppose $f: I_x \isom I_\phi$ is an $(i,i')$-compatible isomorphism. By definition there exists $\tau \in I_{\phi}^{\ad}(\A_f)$ such that
 		\begin{align}\label{eq:tauandf}
 			\iota_{ v} \circ f^{-1} = \iota'_{v} \circ \Int(\tau_v) : I_{\phi, \Q_v} \To  I_v.\end{align} We have a commutative diagram 
 		\begin{align*}  \xymatrixcolsep{4pc}
 			\xymatrix{I_v^{\natural} \ar[r]^{\Int (u_v)^{-1}}_{\sim} & I_v \ar[r]^{ \Int (\iota_v'(\tau_v))^{-1}  }_{\sim} & I_v \\ & T_{\QQ_v} \ar[ul] ^{i} \ar[ur] _{i' }  }
 		\end{align*} Here $\iota_v' (\tau_v)$ is an element of $I_v^{\ad}(\QQ_v)$. We view the component $u_v$ of $u \in G(\A_f^*)$ at $v$ as an element of $R_v(\QQ_v)$, and $\Int (u_v)$ is the associated inner automorphism of $R_v$. 
 	\end{enumerate}
 \end{lem}
 \begin{proof}
 	For part (i), the first diagram obviously commutes for $v\neq p$. For $v=p$, the commutativity boils down to the $T$-equivariance of the isomorphism $V^*_{\Qpur} \isom \cV_0(x)[1/p]$ associated with $y_p$ (where $T$ acts on the left hand side via $i$), which is established in the proof of Proposition \ref{prop:special points on the geom side} (ii). The second diagram commutes because $\iota'_v$ is induced by $\Int(y'_v)^{-1}$ and we have $y'_v \in i'(T(\ol \QQ_v))$. 
 	
 	We now show part (ii). Clearly $\Int(u_v)^{-1}: R_v \isom R_v$ restricts to an isomorphism $I_v^{\natural} \isom I_v$. By part (i), the image of $i_{\ol \QQ_v}$ is contained in $(I_v^\natural)_{\ol \QQ_v}$ viewed as a subgroup of $G_{\ol \QQ_v}$. Moreover, the map $i_{\ol \QQ_v} : T_{\ol \QQ_v} \to (I_v^\natural)_{\ol \QQ_v}$ is defined over $\QQ_v$. Similarly, we have a $\QQ_v$-map $i': T_{\QQ_v} \to I_v$. 
 	
 	We have the following equalities between maps from $T_{\QQ_v}$ to $I_v$: 
 	\begin{align*}
 		\iota'_{v} \circ \Int(\tau_v) \circ  j' & =   \iota_{v} \circ f^{-1} \circ j'  & \text{by (\ref{eq:tauandf})}  \\ \nonumber & = \iota_{v} \circ j  & \text{since $f$ is $(i,i')$-compatible} \\ \nonumber & = \Int(u_v)^{-1} \circ \iota_{v}^{\natural} \circ j .  
 	\end{align*} Hence $$ \iota'_v \circ j'  = \Int (\iota'(\tau_v))^{-1} \circ \Int(u_v)^{-1} \circ \iota_v^{\natural} \circ j. $$ But $\iota'_v \circ j' $ and $\iota_v^{\natural} \circ j$ are the aforementioned maps $i': T_{\QQ_v} \to I_v$ and $i: T_{\QQ_v} \to I_v^{\natural}$ respectively.  This proves the desired commutative diagram. 
 \end{proof}

 \subsection{}\label{para:LD1}
 Let $ (\fks, \fks')= (T,h, i,i')$ be an amicable special fork. Write $x$ and $\phi$ for $x_\fks$ and $\phi(\fks')$. As in  Definition \ref{defn:matchable} we have a well-defined element
 $$ \tau^{\fkH}(\II_{\fks}, \JJ_{\fks'})  \in   \fkH(\phi)  . $$ Since $T$ is naturally a maximal torus in $I_{\phi}$, we have a natural map $\fkH(\phi) \to \fkH(\phi,T) : =\cok\big( \cG^{\ab}(\Z_p) \xrightarrow{\partial} \coh^1(\A_f, T^{\dagger})/ \Sha^{\infty}_{G_{\der}} (\Q, T^{\dagger})\big)$ as in \S \ref{subsec:refined tori-rational}. 
 \begin{prop}\label{prop:LD1}
 	Assume that there exists an $(i,i')$-compatible isomorphism $f: I_{x} \isom I_{\phi}$. Then the image of $\tau^{\fkH}(\II_{\fks}, \JJ_{\fks'})$ in $\fkH(\phi,T)$ is trivial.
 \end{prop}
 \begin{proof}  As in \S \ref{para:setting for comm diags}, we extend $(\fks, \fks')$ to an extended special fork $\fkg$, choose a rectification $u$, and retain all the notations there. Let $\tau$ be as in (\ref{eq:tauandf}). Then $\tau$ is a representative of $\tau^{\fkH}(\II_{\fks}, \JJ_{\fks'})$. We show that its image in $\fkH(\phi,T)$ is trivial.	
 	
 	\textbf{(I) Constructing an element $t_v \in T(\ol \QQ_v)$. }

 	Let $I_\phi^{\dagger} \hookrightarrow I_\phi$ be as  in \S \ref{para:underline isom in pro case} (cf.~\S \ref{para:choice of q_tau}). 
 	Choose $\tilde \tau_v \in I_{\phi}^{\dagger} (\ol{\Q}_v)$ lifting $\tau_v \in I_{\phi}^{\ad}(\Q_v)$, and write $\hat \tau_v $ for $\iota'_v (\tilde \tau_v ) \in I_v(\ol{\Q}_v) \subset R_v(\ol \QQ_v)$. 	Choose $g\in G_{\der}(\Qbar)$ such that
 	$\Int (g) \circ i = i'$, and view $g$ as an element of $R_v(\ol \QQ_v)$ via $\Delta$. Set $$s_v = u_v \hat \tau_v g \in R_v(\ol \QQ_v).$$ We claim that $s_v$ lies in the image of $\Delta \circ i$ in $R_v$. By Lemma \ref{lem:various comm diags} (ii), we know that $s_v$ commutes with this image. Thus for $v \neq p$ the claim is proved. For $v =p$, by Lemma \ref{lem:various comm diags} (i) we know that the image of $\Delta \circ i$ is contained in $I_p^{\natural} \subset R_p$. This image is thus a maximal torus in $I_p^{\natural}$, so it suffices to show that $s_p$ lies in $I_p^{\natural} (\Qpbar)$.  Let $\sigma_R$ denote the automorphism of $R_p$ corresponding to $\sigma \in \Gal(\Qpur/\Q_p)$, i.e., it acts on $R_p(A) = G(\Qpur\otimes_{\Q_p} A
 	)$ via $\sigma$ on the first factor $\Qpur$ for any $\Q_p$-algebra $A$. In the following we view $i(\delta_{\fkg}), i'(\delta_{\fkg}) \in G(\Qpur)$ naturally as elements of $R_p(\QQ_p)$. Using that $i$ and $i'$ are defined over $\QQ_p$, it is straightforward to check that inside $R_p(\Qpbar)$ we have $g i(\delta_\fkg) g^{-1} = i'(\delta_{\fkg})$. (Recall that $g$ is embedded into $R_p(\Qpbar)$ diagonally.)  We then have
 	\begin{align*}
 		s_p  i(\delta_{\fkg}) \sigma_R (s_p)^{-1}  &   = u_p \hat \tau_p g i(\delta_{\fkg}) g^{-1} \sigma_R (\hat \tau_p)^{-1} \sigma_R (u_p)^{-1} & \text{since }\sigma_R(g) = g \\
 		& = u_p \hat \tau_p  i'(\delta_{\fkg}) \sigma_R (\hat \tau_p)^{-1} \sigma_R (u_p)^{-1}  \\
 		& =  u_p i'(\delta_{\fkg}) \sigma_R(u_p)^{-1} & \text{since } \hat \tau_p \in I_p(\Qpbar) \\ & = i(\delta_{\fkg}) & \text{since $u$ is a rectification}. 
 	\end{align*}
 	This precisely means that $s_p \in I_p^{\natural}(\Qpbar)$, as desired. The claim is proved. We write $$ s_v = i (t_v), \quad t_v \in T(\ol \Q_v). $$ 
 	
 	\textbf{(II) Relationship between $t_v$ and $\tau$.}
 	
 	We show that the image of $\tau$ in $\fkH(\phi,T)$ is represented by $(t_v)_v$ in a suitable sense. 
 	
 	We view the embeddings $j' : T \hookrightarrow I_{\phi}$ and $i' : T \hookrightarrow G$ both as inclusions and omit them from the notation.  In the following we let $\Gamma_p$ act on $R_p(\Qpbar) = G(\Qpur \otimes_{\Q_p} \Qpbar)$ only via the second factor $\Qpbar$. For each finite place $v$, the cocycle
 	$(\hat \tau_v^{-1} \lix^\rho \hat \tau_v )_{\rho \in \Gamma_v}$ is valued in $T^{\dagger}$, and the collection of these cocycles for all $v$ represents the image of $\tau$ in $\fkH(\phi,T)$. For $\rho \in \Gamma_v$, we compute
 	\begin{multline}  \hat \tau_v ^{-1} \lix^{\rho} \hat \tau_v  =   g s_v^{-1} u_v \lix ^{\rho} u_v^{-1} \lix ^{\rho}s_v \lix ^{\rho}g ^{-1}  = g s_v^{-1} \lix ^{\rho}s_v \lix ^{\rho}g ^{-1} \\  = g i (t_v ^{-1}\lix^{\rho} t_v) g^{-1} g \lix^\rho g^{-1} = t_v^{-1} \lix^{\rho} t_v g \lix^\rho g^{-1} . 
 	\end{multline}
 	By Lemma \ref{lem:trivial at infty}, the cocycle $(g \lix^\rho g^{-1})_{\rho \in \Gamma}$ represents an element of $\Sha_{G_{\der}} ^\infty (\Q, T^\dagger)$. Therefore the cocycle $(t_v^{-1} \lix^{\rho} t_v)_{\rho \in \Gamma_v}$ represents an element $\beta_v \in \coh^1(\Q_v, T^{\dagger})$, and the collection $(\beta_v)_v$ represents the image of $\tau$ in $\fkH(\phi,T)$.

 	\textbf{(III) Finishing the proof.}

 	To finish the proof it suffices to show that $(\beta_v )_v$ has trivial image in $\fkH(\phi,T)$. 
 	Denote the map $G \to G^{\ab}$ by $r \mapsto r^{\ab}$. 
 	Since $s_p$ is $\sigma_R$-invariant (by the claim in step (I)), $u_p \in R_p(\QQ_p)$ differs from $\sigma_R(u_p)$ by an element of $(\Res_{\Qpur/\Q_p} G_{\der})(\Qpbar)$. It follows that the element  $u_p^{\ab} \in G^{\ab}(\Qpur)$ lies in $G^{\ab}(\Q_p)$. Thus $u^{\ab} \in G^{\ab} (\A_f)$.	Now $\beta_v$ is equal to the image of $t_v^{\ab}  = s_v^{\ab} = u_v^{\ab} \in G^{\ab} (\Q_v)$ under the coboundary map $G^{\ab} (\Q_v) \to \coh ^1(\Q_v, T^{\dagger})$ associated with the short exact sequence $1 \to T^{\dagger} \to T \to G^{\ab} \to 1$.  Since $u \in \ker (G(\A_f^*) \to \pi^*(G))$, we can write $u = u^{(0)} u^{(1)} $, with  $u^{(0)} \in G(\Q)_+$ and  $u^{(1)} \in  G_{\der}(\A^*_f)\cG(\Z_p^{\ur})$. Since $u^{\ab} \in G^{\ab}(\A_f)$, we have $u^{(1) ,\ab}\in \cG^{\ab} (\Z_p)$, with trivial components away from $p$. By Proposition \ref{prop:G(Q)_+}, the image of $u^{(0),\ab}$ in $\coh^1(\A_f, T^{\dagger})$ comes from $\Sha^{\infty} _{G_{\der}} (\Q, T^{\dagger})$. Hence both $u^{(0),\ab}$ and $u^{(1),\ab}$ vanish in $\fkH(\phi,T)$. 
 \end{proof}
 \subsection{}\label{para:setting for LD2} Let  
 $(T,h, i,i')$ be a special fork.
 The two base changes $T\times_{G} G_{\sconn}$ along $i$ and $i' : T \to G$ are canonically identified as a $\QQ$-torus equipped with a map to $T$, which we denote by  $\tilde T \to T$. Let $\tilde i$ and $\tilde i': \tilde T \hookrightarrow G_{\sconn}$ be the base changes of $i$ and $i'$.  Fix a class $\beta_{0} \in \Sha^{\infty}(\QQ, \tilde T)$ and a representing cocycle $\beta_{0} (\cdot)$. We use $\beta_{0} (\cdot)$ to twist the special fork $(T,h,i,i')$ as follows.   
 
 Since $\Sha^\infty (\QQ, G_{\sconn}) = \Sha^\infty_{\ab} (\QQ, G_{\sconn})$ is trivial, we can choose  $\lambda, \lambda' \in G_{\sconn} (\overline \Q)$ such that
 \begin{align}\label{eq:cons of lambda and lambda'}
 	\tilde i(\beta_0(\rho)) = {\lambda} ^{-1 } \lix^{\rho} \lambda , \quad \tilde i'(\beta_0(\rho)) = {\lambda'} ^{-1 } \lix^{\rho} \lambda ',  \quad \forall {\rho} \in \Gamma.
 \end{align} Since $\beta_0$ is trivial at infinity, we have 
 \begin{align}\label{eq:top of lambda and lambda'}
 	\lambda \in G_{\sconn}(\RR) \tilde i(\tilde T(\CC)), \quad 	\lambda' \in G_{\sconn}(\RR) \tilde i'(\tilde T(\CC)).
 \end{align}

 We define
 $$i_1 : =\Int (\lambda)  \circ  i:   T \to G, \quad i'_1 : =\Int (\lambda')  \circ  i':   T \to G .$$ By (\ref{eq:cons of lambda and lambda'}), $i_1$ and $i_1'$ are defined over $\QQ$. By (\ref{eq:top of lambda and lambda'}), the connectedness of $G_{\sconn}(\RR)$ (Cartan's theorem), and the fact that $i\circ h$ and $i' \circ h$ lie in the same connected component of $X$, we know that all four points $i_1 \circ h, ~i'_1 \circ h, ~i\circ h, ~ i' \circ h$ lie in the same connected component of $X$. It is also clear that $i_1$ and $i_1'$ are conjugate by $G^{\ad}(\Qbar)$. Hence we obtain a new special fork $$(T,h, i_1, i_1'), $$ regarded as the twist of $(T,h,i,i')$ by $\beta_0(\cdot)$. It depends only on the class $\beta_0$ and not on the choices of $\beta_0(\cdot), \lambda,\lambda'$ up to conjugating $i_1$ and $i_1'$ separately by $G(\QQ)$.  
 
 We write $x = x_{(T,i,h)}, \phi = \phi(T,i',h), x_1 = x_{(T,i_1,h)}, \phi_1 = \phi(T, i_1', h)$.  Let 
 \begin{align*}
 	j: T \hookrightarrow I_x, && j': T \hookrightarrow I_{\phi}, &&  j_1: T \hookrightarrow I_{x_1}, && j_1': T \hookrightarrow I_{\phi_1}
 \end{align*}
 be the canonical $\Q$-embeddings. We next construct natural inner twistings of $\QQ$-groups $I_{\phi, \Qbar} \isom I_{\phi_1, \Qbar}$ and $I_{x,\Qbar} \isom I_{x_1, \Qbar}$ and compute their cocycles. 
 
 Since $\phi_1^{\Delta} = \Int(\lambda') \circ \phi^{\Delta}$, we have a $\Qbar$-isomorphism $$\psi: I_{\phi, \Qbar} \isom I_{\phi_1, \Qbar}$$ induced by $\Int(\lambda')$. For each $\rho \in \Gamma$ let $q_{\rho} \in \gQ^{\topo}$ be a lift of $\rho$, and write $\phi(q_{\rho}) = g_{\rho} \rtimes \rho$ and $\phi_1(q_{\rho}) = g_{1,\rho} \rtimes \rho$. Then $g_\rho \in i'(T) (\Qbar)$ and $g_{1,\rho} = \Int(\lambda')(g_{\rho})$. We compute 
 \begin{align*}
 	\lix^{\rho}\psi (\cdot) & = \Int (g_{1,\rho}) \bigg( \rho \bigg ( \Int(\lambda' g_{\rho^{-1}}) \big( \rho^{-1}(\cdot)   \big)\bigg)  \bigg ) \\ &   = \Int(\lambda' g_{\rho} \lambda'^{-1} \cdot \lix^{\rho} \lambda' \cdot \lix^{\rho} g_{\rho^{-1}}) (\cdot ) \\  & = \Int(\lambda' g_{\rho} \lambda'^{-1} \lix^{\rho} \lambda'  g_{\rho}^{-1}) (\cdot )  ,
 \end{align*}
 where in the last step we use that $g_{\rho} \lix^{\rho} g_{\rho^{-1}}$ is central in $I_{\phi}(\Qbar)$. By (\ref{eq:cons of lambda and lambda'}) and the fact that $g_{\rho} \in i'(T)$, we know that $g_\rho$ commutes with $\lambda'^{-1} \lix^{\rho} \lambda'$.  Hence 
 $ \lix^{\rho}\psi(\cdot) = \Int(\lix^{\rho} \lambda') (\cdot) .$ Write $\beta(\rho)$ for the image of $\beta_0(\rho) \in \tilde T(\Qbar) $ in $T(\Qbar)$. Thus $\psi$ is an inner twisting satisfying
 \begin{align}\label{eq:cocycle of psi}
 	\psi^{-1} \circ \lix^\rho \psi  = \Int ( j'(\beta(\rho))) \in \Aut(I_{\phi ,\Qbar}), \quad \forall \rho \in \Gamma. 	
 \end{align}
 We also record that $\psi$ and the map  $\Int(\lambda') : G_{\sconn,\Qbar} \to G_{\sconn,\Qbar} $ 	together define an isomorphism 
 $$\tilde \psi: \tilde I_{\phi,  \Qbar}  = I_{\phi ,\Qbar} \times_{G_{\Qbar}} G_{\sconn, \Qbar} \isom \tilde I_{\phi_1, \Qbar} = I_{\phi_1, \Qbar} \times _{G_{\Qbar}} G_{\sconn,\Qbar}.$$ (See \S \ref{para:underline isom in pro case} and \S \ref{para:choice of q_tau} for $\tilde I_\phi$.) The same computation as above shows that $\tilde \psi$ is an inner twisting satisfying 
 \begin{align}\label{eq:cocycle of tilde psi} \tilde \psi^{-1} \circ \lix^\rho  {\tilde \psi}  = \Int ( j'(\beta(\rho)))  = \Int ( \tilde j'(\beta_0(\rho))) \in \Aut(\tilde I_{\phi ,\Qbar}), \quad \forall \rho \in \Gamma.   
 \end{align}
 
 We now construct an inner twisting $\chi: I_{x, \Qbar} \isom I_{x_1,\Qbar}$. As in \S \ref{para:x_fks}, we let $\tilde x$ and $\tilde x_1$ be points of $\Sh_{K}$ associated with $(T,i,h)$ and $(T,i_1,h)$ respectively.  Under the canonical isomorphisms 
 $$ V^*_{\ZZ_{(p)}} \isom  \coh^1_B(\cA_{\asim{x} }(\CC), \ZZ_{(p)} )  \quad \text{and}  \quad  V^*_{\ZZ_{(p)}} \isom  \coh^1_B(\cA_{\asim{x}_1 }(\CC), \ZZ_{(p)} ) ,$$ 
 the element $\lambda \in G(\Qbar)$ induces an isomorphism of $\Qbar$-Hodge structures $$\lambda_B: \coh^1_B(\cA_{\asim{x} }(\CC), \Qbar )  \isom \coh^1_B(\cA_{\asim{x}_1 }(\CC), \Qbar ).  $$ In particular $\lambda_B$ induces a $\Qbar$-isogeny $\theta$ from $\cA_{x_1}$ to $\cA_{x}$. Via the comparison isomorphisms, $\lambda_B$ induces isomorphisms 
 $$ \cV_p(\tilde x)[1/p] \otimes_{\QQ} \Qbar \isom \cV_p(\tilde x_1)[1/p] \otimes_{\QQ} \Qbar , \quad \cV^p(x)\otimes_{\QQ} \Qbar \isom \cV^p(x_1)\otimes_{\QQ} \Qbar. $$
 The first one sends $s_{\alpha,p, \asim {x}}$ to $s_{\alpha,p, \asim{x}_1}$ and preserves the $\Gal(\Qpbar/F)$-action for a sufficiently large finite extension $F/\QQ_p$, so applying $D_{\cris}$ we see that $\theta$ sends $s_{\alpha, 0, x}$ to $s_{\alpha, 0, x_1}$. The second one sends $s_{\alpha, \A_f^p, x}$ to $s_{\alpha, \A_f^p, x_1}$. It follows that $\theta  \in  I_{x_1,x}(\Qbar)$. Clearly   
 \begin{align}\label{eq:ass on chi}
 	\theta \circ \rho (\theta^{-1})  = j(\beta (\rho)) \in j(T(\Qbar)) \subset  I_x (\Qbar), \quad  \forall \rho \in \Gamma. 
 \end{align}
 Define $$\chi : I_{x,\Qbar} \isom I_{x_1, \Qbar}$$ to be the isomorphism induced by $\theta$. From (\ref{eq:ass on chi}) it immediately follows that $\chi$ is an inner twisting satisfying
 \begin{align}\label{eq:cocycle of chi}
 	\chi^{-1} \circ \lix^{\rho} \chi = \Int(j(\beta(\rho))) \in \Aut(I_{x,\Qbar}), \quad \forall \rho \in \Gamma.
 \end{align}
 \begin{prop}\label{prop:LD2} 
 	Keep the setting and notation of \S \ref{para:setting for LD2}. Assume that $(x,\phi)$ is amicable, and that there  exists an $(i,i')$-compatible isomorphism $f: I_{x} \isom I_{\phi}$. Then $(x_1,\phi_1)$ is amicable, and there exists an $(i_1,i_1')$-compatible isomorphism $f_1 : I_{x_1} \isom I_{\phi_1}$. Moreover, the canonical elements $\tau^{\fkH} \in \fkH(\phi)$ and $\tau^{\fkH}_1 \in \fkH(\phi_1)$ associated with the amicable pairs $(x,\phi)$ and $(x_1,\phi_1)$ respectively (Definition \ref{defn:matchable}) correspond to each other under the canonical identification $\fkH(\phi)\cong \fkH(\phi_1)$ as in  \S\ref{subsec:frak H} (since we have $\phi \approx \phi_1$). 
 \end{prop} 
 \begin{proof} 	
 	
 	By Lemma \ref{lem:getting gauge}, we extend the special forks $(T,h,i,i')$ and $(T,h,i_1,i_1')$ to extended special forks of the form $$\fkg = (T,h,i,i', y_T, y'_T), \quad \fkg_1 = (T,h,i_1,i_1', y_T, y'_T) . $$
 	Since $(x,\phi)$ is amicable, by Lemma \ref{lem:rec gauge} we find a  rectification $u$ of $\fkg$ and a resulting marking $\fkm(\fkg, u) = (yu, y') \in \tmar(x,\phi)$. Also let $y_1$ be the image of $y_T$ in $Y(x_1)$, and $y_1'$ be the image of $y'_T$ in $Y(\phi_1)$.  
 	We have $\delta_{\fkg_1} = \delta_{\fkg} \in T(\Qpur)^{\mot}$, and we denote this element by $\delta_T$. For $v$ a finite place, as in \S \ref{para:setting for comm diags} we set  $ R_v$ to be $G_{\QQ_v}$ for $v\neq p$ and $\Res_{\Qpur/\QQ_p} G$ for $v=p$. Similarly we define $R_{\sconn, v}$ using $G_{\sconn}$.    Let $I_v^{\natural}$, $I_v$, $I_{v,1}^{\natural}$, $I_{v,1}$ be the reductive groups over $\Q_v$ associated with $\fkk_{\fkg}, \fkk'_{\fkg}, \fkk_{\fkg_1}, \fkk'_{\fkg_1}$  respectively. These are all $\QQ_v$-subgroups of $R_v$. Each of them has a tilde companion as defined in \S \ref{para:tilde I for Kott tuple}, which is naturally viewed as a $\QQ_v$-subgroup of $R_{\sconn,v}$.  Let 
 	$$ \begin{cases}
 		\iota_v^{\natural}:   I_{x,\Q_v} \isom I_v^{\natural} , \\
 		\iota_v : I_{x,\Q_v} \isom I_v , \\
 		\iota_v':  I_{\phi,\Q_v} \isom I_v , \\
 		\iota_{v,1}^{\natural} :  I_{x_1,\Q_v} \isom I_{v,1}^{\natural},\\
 		\iota_{v,1}' :  I_{\phi_1, \Q_v} \isom I_{v,1}
 	\end{cases} \quad \text{be associated with} \quad \begin{cases}
 		y \in   Y(x) ,\\
 		y   u \in  Y(x) ,\\
 		y'  \in  Y(\phi),\\
 		y_1 \in  Y(x_1), \\
 		y'_1 \in   Y(\phi_1) 
 	\end{cases} $$
 	respectively.  Write $u_v$ for the component of $u$ in $R(\Q_v)$.  
 	
 	\textbf{(I) A commutative diagram of cohomology.}
 	
 	Let  $\tau = (\tau_v)_v \in I_{\phi}^{\ad}(\A_f)$ be associated with $(f, \iota_v, \iota'_v)$ as in (\ref{eq:tauandf}). Fix a lift $(\tilde u_v)_{v\neq p} \in G_{\sconn}(\bar \A_f^p)$ of the image of $(u_v)_{v\neq p}$ in $G^{\ad}(\A_f^p)$. Also fix a lift $\tilde u_p \in R_{\sconn,p}(\Qpbar)$ of the image of $u_p$ in $(\Res_{\Qpur/\Q_p} G^{\ad})(\Q_p)$. For each finite $v$,  since $\Int(\iota'_{v}(\tau_v))^{-1}$ induces the identity on $\coh ^1(\Q_v, \tilde I_v) \cong \coh ^1_{\ab}(\Q_v, \tilde I_v)$, we obtain from Lemma \ref{lem:various comm diags} (ii) the commutative diagram
 	\begin{align*} 
 		\xymatrix{\coh^1(\QQ_v,\tilde I_v^{\natural}) \ar[rr]^{\Int (\tilde u_v)^{-1}} &  & \coh^1(\QQ_v,\tilde I_v) \\  &\coh^1(\QQ_v,\tilde T) \ar[lu]^{\tilde i} \ar[ru]_{\tilde i' } &    }
 	\end{align*} Hence there exists $\Delta_v\in \tilde I_{v}(\overline \Q_v)$ such that \begin{align}\label{eq:introDelta}
 		\Int	(\tilde u_v	)^{-1} [\tilde i (\beta_0 (\rho))] = \Delta_v \cdot \tilde  i' (\beta_0(\rho)) \cdot \lix ^{\rho} \Delta_v^{-1} , \quad \forall \rho \in \Gamma_v.
 	\end{align}
 	In view of (\ref{eq:cons of lambda and lambda'}), we can rewrite (\ref{eq:introDelta}) as
 	\begin{align}\label{eq:Delta_v}
 		\Int	(	\tilde u_v)	 ^{-1} [\lambda ^{-1} \lix^{\rho} \lambda]  = \Delta_v \cdot  {\lambda'}^{-1} \lix ^{\rho} \lambda '  \cdot \lix ^{\rho}  \Delta _v ^{-1}, \quad \forall \rho \in \Gamma_v.
 	\end{align}
 	This is an equality between elements of $R_{\sconn,v}(\ol \QQ_v)$; when $v=p$ we view $\lambda, \lambda' \in G_{\sconn}(\Qpbar)$ as inside $R_{\sconn,p}(\Qpbar)$ via the diagonal embedding $G_{\sconn} \to R_{\sconn,p}$.
 	
 	\textbf{(II) Constructing $u_1$.}
 	
 	We shall construct a rectification $u_1$ of $\fkg_1$. Denote the map $R_{\sconn} \to R$ by $\varrho_R$.  For each finite place $v$, let $$u_{1,v} :  = \varrho_R (\lambda) u_v  \varrho_R(\Delta_v   \lambda'^{-1}) \in R(\ol \QQ_v). $$ Note that we have $$u_{1,v} = u_v \cdot \varrho_R(w_v), $$ with 
 	$$ w_v : = \tilde u_v^{-1} \lambda  \tilde u_v  \Delta_v   \lambda'^{-1} \in R_{\sconn}(\ol \QQ_v).  $$ It immediately follows from (\ref{eq:Delta_v}) that $w_v \in R_{\sconn}(\Q_v)$.

 	We claim that by suitably choosing the elements $\Delta_v$ we can arrange the element $(w_v)_{v\neq p} \in \prod_{v\neq p}G_{\sconn}(\Q_v)$ to lie in $G_{\sconn}(\A_f^p)$. 
 	Indeed, the isomorphism $\iota'_{v} : I_{\phi,\QQ_v} \isom I_v$ canonically lifts to an isomorphism $\tilde \iota'_{v}: \tilde I_{\phi,\QQ_v} \isom \tilde I_v$, where $\tilde I_{\phi}$ is as in \S \ref{para:underline isom in pro case} (cf.~\S \ref{para:choice of q_tau}). If we apply $\tilde \iota^{\prime -1}_v$ to the equality (\ref{eq:introDelta}) and set $\Omega_v : =\tilde \iota^{\prime -1}_{v}(\Delta_v)$, we see that the choice of $\Delta_v$ is equivalent to the choice of $\Omega_v \in \tilde I_{\phi}(\ol \QQ_v)$ such that 
 	$$ \tilde \iota_v^{\prime -1} \Int(\tilde u_v)^{-1} [ \tilde i(\beta_0(\rho))] = \Omega_v \cdot \tilde j' (\beta_0(\rho)) \cdot \lix^{\rho} \Omega_v^{-1}, \quad \forall \rho \in \Gamma_v, $$ where $\tilde j' : \tilde T \to \tilde I_{\phi}$ is the canonical lifting of $j': T \to I_{\phi}$ (which is defined over $\Q$).  By (\ref{eq:tauandf}) and the fact that $\iota_v^{\natural} = \Int(u_v) \circ \iota_v$, we have $$
 	\iota_v^{\prime -1} \circ \Int(u_v)^{-1} \circ i = \Int (\tau_v) \circ f \circ \iota_v^{-1} \circ \Int(u_v)^{-1} \circ i = \Int (\tau_v) \circ f \circ \iota_v^{\natural -1} \circ i  . $$ By Lemma \ref{lem:various comm diags} and the fact that $j' = f\circ j$, the above is equal to $\Int(\tau_v) \circ j'$. Therefore the equation for $\Omega_v$ is  equivalent to 
 	\begin{align}\label{eq:commute cocycle}
 		\Int(\tau_v) [ \tilde j'(\beta_0(\rho))] = \Omega_v \cdot \tilde j' (\beta_0(\rho)) \cdot \lix^{\rho} \Omega_v^{-1}, \quad \forall \rho \in \Gamma_v.
 	\end{align}
 	The collection of the local cocycles appearing on the left hand side, for all $v\neq p$, forms an adelic cocycle in $Z^1(\Gamma, \tilde I_{\phi}(\bar \A_f^p))$. The cocycle $\tilde j' (\beta_0(\rho))$ comes from a global cocycle in $Z^1(\Gamma, \tilde I_{\phi} (\Qbar))$. Hence the elements $\Omega_v$ can be chosen in such a way that   $(\Omega_v)_{v\neq p} \in  \tilde I_{\phi}( \bar  \A_f^p) \subset \prod_{v \neq p} \tilde I_{\phi}(\ol \Q_v)$. It then follows that $(\Delta_v)_{v\neq p}$ lies in $G_{\sconn}(\bar \A_f^p)$. Then by the definition of $w_v$, we have $(w_v)_{v\neq p} \in G_{\sconn}(\bar \A_f^p)$. Since we already know that $w_v \in G_{\sconn}(\Q_v)$ for each $v\neq p$, we have $(w_v)_{v\neq p} \in G_{\sconn}(\A_f^p)$.  This proves the claim. 
 	
 	Set $u_1 = (u_{1,v})_v$. Then $u_1 \in u\cdot \varrho ( G_{\sconn}(\A_f^*) ) \subset G(\A_f^*)$. 
 	
 	\textbf{(III) Proof that $u_1$ is a rectification of $\fkg_1$.}
 	
 	For a sufficiently divisible $n\in \ZZ_{\geq 1}$, let $$\gamma_{0,T,n} = \delta_T \sigma(\delta_T) \cdots \sigma^{n-1}(\delta_T) \in T(\Q)$$ (see \S \ref{para:special kott trip}).  For a finite place $v \neq p$, we have
 	\begin{align}\label{eq:checkrec1}
 		& \Int(u_{1,v})^{-1} (i_1(\gamma_{0,T,n}))\\ \nonumber   &  = \Int \big (\lambda' \Delta_v^{-1} u_v ^{-1} \lambda^{-1} \big ) (i_1(\gamma_{0,T,n})) \\ \nonumber & =  \Int \big (\lambda' \Delta_v^{-1} u_v ^{-1} \big ) (i(\gamma_{0,T,n}))  & \text{since }i_1 = \Int(\lambda) \circ i   \\ \nonumber &   = \Int \big (\lambda' \Delta_v^{-1} \big ) (i'(\gamma_{0,T,n})) & \text{since } \fkk_{\fkg} \cdot u = \fkk'_{\fkg} \\ \nonumber &    =   \Int \big (\lambda'\big )(i'(\gamma_{0,T,n}))& \text{since } \Delta_v \in 
 		\tilde I_v   \text{ centralizes } i'(\gamma_{0,T, n})  \\  \nonumber &  = i_1' (\gamma_{0,T,n}) . & \text{since }i_1' = \Int(\lambda') \circ i' 
 	\end{align}
 	For $a,b \in G(\Qpur \otimes_{\Q_p} \Qpbar)$, we write $\Int_{\sigma}(a) (b)$ for $a b \sigma_R(a)^{-1}$, where $\sigma_R$ acts on $ G(\Qpur \otimes_{\Q_p} \Qpbar)$  via $\sigma$ on the first factor $\Qpur$. By an analogous computation, we obtain 
 	\begin{align}\label{eq:checkrec2}
 		\Int_{\sigma}(u_{1,p})^{-1} (i_1(\delta_T))    & = \Int_{\sigma} \big (\lambda' \Delta_p ^{-1} u_p ^{-1} \lambda^{-1} \big ) (i_1(\delta_T)) \\ \nonumber &   =  \Int_{\sigma} \big (\lambda' \Delta_p^{-1} u_p ^{-1} \big ) (i(\delta_T))  \\ \nonumber &   = \Int _{\sigma} \big (\lambda' \Delta_p^{-1} \big ) (i'(\delta_T)) \\  \nonumber &   =   \Int_{\sigma} \big (\lambda'\big )(i'(\delta_T)) \\ \nonumber &  = i_1' (\delta_T) . 
 	\end{align}	Since $u_1 \in u\cdot \varrho ( G_{\sconn}(\A_f^*) )$, we have $u_1 \in \ker (G(\A_f^*) \to \pi^*(G))$. This together with   (\ref{eq:checkrec1}) and (\ref{eq:checkrec2}) implies that $u_1$ is a rectification of $\fkg_1$.
 	By Lemma \ref{lem:rec gauge}, $(x_1,\phi_1)$ is amicable. 
 	
 	\textbf{(IV) Constructing  $f_1$ and $\tau_1$.}

 	We define $f_1$ by the commutative diagram
 	$$ \xymatrix{I_{x, \Qbar} \ar[r]^f \ar[d]^{\chi} & I_{\phi, \Qbar} \ar[d]^{\psi} \\ I_{x_1,\Qbar} \ar[r]^{f_1}  & I_{\phi_1, \Qbar}} $$
 	By (\ref{eq:cocycle of psi}), (\ref{eq:cocycle of chi}), and the fact that $f\circ j = j'$, we know that $f_1$ is a $\Q$-isomorphism. Moreover, $\chi$ and $\psi$ commute with the canonical embeddings of $T$. Hence $f_1 \circ j_1 = j_1'$.  
 	In order to show that $f_1$ is an $(i_1,i_1')$-compatible isomorphism, it remains to show that $f_1$ belongs to  the canonical $I_{\phi_1}^{\ad}(\Q)$-orbit of isomorphisms $I_{x_1} \isom I_{\phi_1}$. Since we have already shown that $u_1$ is a rectification of $\fkg_1$, we have 
 	$$\fkm (\fkg_1, u_1)  = ( y_1  u_1, y'_1)  \in \tmar(x_1,\phi_1)$$ as in  Lemma \ref{lem:rec gauge}. 
 	Let $
 	\iota_{v,1} = \iota_{y_1u_1, v}= \Int(u_{1,v})^{-1} \circ \iota^{\natural}_{v,1} : I_{x_1, \Q_v} \isom I_{v,1}.$ Then the desired property of $f_1$ is equivalent to the existence of $\tau_1 \in I_{\phi_1}^{\ad}(\A_f)$ such that 
 	\begin{align}\label{eq:tauandfnew1}
 		\iota_{ v, 1} \circ f_1^{-1} = \iota'_{v,1} \circ \Int(\tau_{1,v}) : I_{\phi_1, \Q_v} \To I_{v,1}. \end{align} From the constructions of $\psi$ and $\chi$, we have the commutative diagrams
 	\begin{align}\label{diag:psi}\xymatrixcolsep{4pc}
 		\xymatrix{ I_{\phi,\ol \Q_v} \ar[r]^{\iota'_{v}} \ar[d]^{\psi} & I_{v,\ol \Q_v} \ar[d]^{\Int (\lambda')} \\ I_{\phi_1,\ol \Q_v} \ar[r]^{\iota'_{ v ,1} } & I_{v,1,\ol \Q_v}}
 	\end{align}
 	and
 	\begin{align}\label{diag:chi} \xymatrixcolsep{4pc}
 		\xymatrix{ I_{x,\ol \Q_v} \ar[r]^{\iota_v^{\natural}} \ar[d]^{\chi} & I_{v,\ol \Q_v}^{\natural} \ar[d]^{\Int (\lambda)} \\ I_{x_1,\ol \Q_v} \ar[r]^{\iota_{v,1}^{\natural} } & I_{v,1, \ol \Q_v} ^{\natural}}
 	\end{align}
 	Let $\Omega_v\in \tilde I_{\phi}(\ol \Q_v)$ be as in part (II) of the proof. Let $\tilde \tau_v \in \tilde I_{\phi}(\ol \Q_v)$ be a lift of $\tau_v$, and set $$\tilde \tau_{1,v}  : = \tilde \psi( \Omega_v^{-1}  \tilde \tau_v) \in \tilde I_{\phi_1}(\ol \Q_v).$$ For $\rho \in \Gamma_v$, write $\tilde \beta_{\rho}$ for  $ \tilde j'(\beta_0(\rho)) \in \tilde  I_{\phi} (\ol \QQ_v)$ .  
 	We compute
 	\begin{align} \label{eq:tau_1}
 		\tilde \tau_{1,v}^{-1} \cdot  \lix^{\rho}  \tilde \tau_{1,v}  & = \tilde \psi(  \tilde \tau_v^{-1}  \Omega_v )  \cdot (\lix ^{\rho} \tilde \psi) (\lix^{\rho}    \Omega_v^{-1} \lix^\rho  \tilde \tau_v)\\ \nonumber
 		& = \tilde \psi(  \tilde \tau_v^{-1}
 		\Omega_v  ) \cdot \tilde \psi (\tilde \beta_\rho  \lix^{\rho}   \Omega_v^{-1} \lix^\rho  \tilde \tau_v  \tilde \beta_{\rho} ^{-1} ) & \text{by (\ref{eq:cocycle of tilde psi}) }\\ \nonumber
 		& =  \tilde \psi (    \tilde \tau_v^{-1} \Omega_v  \tilde  \beta_\rho \lix^{\rho}  \Omega_v^{-1} \lix^\rho \tilde \tau_v
 		\tilde \beta_{\rho} ^{-1} ) \\ \nonumber
 		& = \tilde \psi( \tilde \beta_{\rho} \tilde \tau_v^{-1} \lix^{\rho}  \tilde \tau_v  \tilde \beta_{\rho}^{-1})  & \text{by (\ref{eq:commute cocycle})}  \\ \nonumber
 		& =\tilde  \psi(\tilde \tau_v^{-1} \lix^{\rho}  \tilde \tau_v)  & \text{because $\tilde \tau_v^{-1} \lix^{\rho}  \tilde\tau_v \in Z_{\tilde I_\phi}$} .
 	\end{align}
 	It follows that the image of $ \tilde \tau_{1,v}$  in $I_{\phi_1}^{\ad} (\ol \Q_v)$ lies in $I_{\phi_1}^{\ad} ( \Q_v) $. We denote it by $\tau_{1,v} \in I_{\phi_1}^{\ad}(\Q_v)$. Recall from part (II) that $(\Omega_v)_{v\neq p}$ comes from an element of $\tilde I_{\phi}(\bar \A_f^p)$. It easily follows that $\tau_1: = (\tau_{1,v})_v \in \prod_v I_{\phi_1}^{\ad}(\Q_v)$ is in fact an element of $I_{\phi_1}^{\ad}(\A_f)$.
 	
 	We now check that $\tau_1 $ satisfies  (\ref{eq:tauandfnew1}). We compute
 	\begin{align*}
 		\iota'_{v,1}  & = \Int (\lambda') \circ \iota'_v \circ \psi^{-1}  & \text{ by (\ref{diag:psi})} \\  & =  \Int (\lambda' \Delta_v^{-1}) \circ \iota'_v \circ  \Int (\Omega_v) \circ \psi^{-1}  & \text{since }\Omega_v = \tilde \iota_v^{\prime -1}(\Delta_v) \\ &
 		= \Int (\lambda'\Delta_v^{-1}) \circ \iota_v \circ f^{-1} \circ \Int(\tilde \tau_v^{-1}   \Omega_v) \circ  \psi^{-1} & \text{by (\ref{eq:tauandf})} \\  & =  \Int (\lambda'\Delta_v^{-1}) \circ \iota_v \circ f^{-1} \circ \psi^{-1} \circ \Int(\tau_{1,v}^{-1}) & \text{by definition of }\tau_{1,v} \\ &   =  \Int (\lambda'\Delta_v ^{-1} u_v^{-1}) \circ \iota_v^{\natural} \circ f^{-1} \circ \psi^{-1} \circ \Int(\tau_{1,v}^{-1})  &  \\ &
 		= \Int (u_{1,v} ^{-1} \lambda) \circ \iota_v^{\natural} \circ f^{-1} \circ \psi^{-1} \circ \Int(\tau_{1,v}^{-1}) & \text{by definition of }u_{1,v} \\  & =  \Int (u_{1,v}^{-1}) \circ \iota_{v,1}^{\natural} \circ \chi \circ f^{-1} \circ \psi^{-1} \circ \Int(\tau_{1,v}^{-1}) & \text{by (\ref{diag:chi})} \\
 		&  =  \Int (u_{1,v}^{-1}) \circ \iota_{v,1}^{\natural} \circ  f_1^{-1} \circ \Int (\tau_{1,v}^{-1}) & \text{by definition of } f_1\\ &
 		= \iota_{v,1} \circ f_1^{-1} \circ \Int (\tau_{1,v}^{-1}) &  
 	\end{align*}as desired.

 	\textbf{(V) Relation between $\tau$ and $\tau_1$.} 
 	
 	It remains to prove the last statement in the proposition. The elements $ \tau^{\fkH}$ and $ \tau^{\fkH}_1$ are the images of $\tau$ and $\tau_1$ in $\fkH(\phi)$ and  $\fkH(\phi_1)$ respectively. By (\ref{eq:tau_1}), the image of $\tau$ in $\coh^1(\A_f, Z_{\tilde I_\phi})$ corresponds to the image of $\tau_1$ in $\coh^1(\A_f, Z_{\tilde I_{\phi_1}})$ under the canonical identification of $\Q$-groups $Z_{\tilde I_\phi} \cong Z_{\tilde I_{\phi_1}}$ (specifically, the one induced by any inner automorphism of $G_{\Qbar}$ taking $\phi^{\Delta}$ to $\phi_1^{\Delta}$). To complete the proof, it suffices to show that the natural map $I_{\phi}^{\ad}(\A_f) \to \fkH(\phi)$ factors through the coboundary map $I_{\phi}^{\ad}(\A_f) \to \coh^1(\A_f, Z_{\tilde I_\phi})$, and likewise with $\phi_1$ in place of $\phi$. We focus on the case of $\phi$.
 	
 	If two elements of $I_{\phi}^{\ad}(\A_f)$ have the same image in $\coh^1(\A_f, Z_{\tilde I_\phi})$, they differ by left multiplication by an element in the image of $\tilde I_{\phi}(\A_f)$. Consequently, it is enough to show that the composition $\tilde I_{\phi} (\A_f) \to I_{\phi} (\A_f) \to \pi(G)$ is trivial. For this, it suffices to show that the composition $\coh^0_{\ab}(\A_f, \tilde I_{\phi}) \to \coh^0_{\ab}(\A_f, I_{\phi}) \to \coh^0_{\ab}(\A_f, G)$ is trivial. Indeed, this follows from the fact that the composition $\pi_1 (\tilde I_\phi) \to \pi_1(I_\phi) \to \pi_1(G)$ is trivial (see Lemma \ref{lem:tilde conn}).  
 \end{proof}
 \section{Galois cohomological properties of amicable pairs}\label{subsec:Galcohprop}
 In this section, we reformulate and extend the results in \S \ref{subsec:gauges} for general amicable pairs, in Theorems \ref{thm:matchable pair 1} and \ref{thm:matchable pair 2}. We then establish a bijection in Theorem \ref{thm:matchability of bij} as a preliminary to the proof of the Langlands--Rapoport--$\tau$ Conjecture.  
 
 \begin{lem}\label{lem:reduction to special fork}
 	Let $(\II, \JJ)$ be an amicable pair. Let $
 	\phi \in \JJ$ and let $T\subset I_{\phi}$ be a maximal torus. Then there exists an amicable special fork $(\fks, \fks') = (T,h,i,i')$ satisfying the following conditions:
 	\begin{enumerate}
 		\item We have $\II = \II_{\fks}$ and $\JJ=\JJ_{\fks'}$.
 		\item There exists $g\in G(\Qbar)$ such that $\phi(\fks') = \Int(g) \circ \phi$ and the composition $T \hookrightarrow I_\phi \xrightarrow{\Int(g)} I_{\phi(\fks')}$ is the canonical embedding $i'$.
 		\item There exists an $(i,i')$-compatible isomorphism (Definition \ref{defn:adapted iso}) $I_{x_{\fks}} \isom I_{\phi(\fks')}$.  
 	\end{enumerate}
 \end{lem}
 \begin{proof}   Choose $x\in \II$, and choose an isomorphism $f: I_x \isom I_{\phi}$ as in Lemma \ref{lem:tau for amicable pair}. Let $j':T \hookrightarrow I_{\phi}$ be the inclusion map, and let $j = f^{-1} \circ j' : T \hookrightarrow I_x$. Via $j$ we view $T$ as a maximal torus in $I_{x}$.
 	
 	We claim that a cocharacter $\mu \in X_*(T)$ is $\phi$-admissible in the sense of \S \ref{para:admissible cochar} if and only if it is $x$-admissible in the sense of \S \ref{para:x-admissible}. Let $(y,y') \in \tmar(x,\phi)^w$, and let  $\delta\in G(\Qpur)$ be the $p$-adic component of $\fkk(y) = \fkk(y')$. Clearly for a $\mu \in X_*(T)$, the composition (\ref{eq:mutoG})  lies in $\dmu_X(\Qpbar)$ if and only if $j' \circ \mu$ lies in $\dmu_X(\Qbar)$. Hence we only need to check that the fractional central cocharacter $-(\iota_{y,p}^{-1}\circ \nu_{\delta})$ of $I_x$ is taken to $(\phi(p) \circ \zeta_p)^{\Delta}$ of $I_{\phi}$  under $f$. Since $\iota_{y',p}\circ f$ differs from $\iota_{y,p}$ by an inner automorphism, it suffices to check that $-\nu_{\delta} = \iota_{y',p} \circ (\phi(p) \circ \zeta_p)^{\Delta}$. By \S \ref{para:defn of Y(phi)}, we have  $\delta = c' b \sigma(c')^{-1}$ with $b = b_{y'_p}, c' = c_{y'_p}$. By Proposition \ref{prop:twist at p}, we know that $\iota_{y',p} \circ (\phi(p) \circ \zeta_p)^{\Delta} =- \nu_{b}$, and that $c'$ commutes with $\nu_b$. The desired result follows. 
 	
 	By the above claim, Theorem \ref{thm:precise speciality}, and Theorem \ref{thm:geometric speciality}, we can find the following objects: 
 	\begin{itemize}
 		\item Two special point data of the form $\fks = (T,i,h), \fks'= (T,i',h) \in \spd(G,X)$ such that $\II = \II_{\fks}$ and $\JJ = \JJ_{\fks'}$.
 		\item An isomorphism $w : I_x \isom I_{x_{\fks}}$ which is induced by some element of $I_{x, x_{\fks}}(\Q)$, satisfying that the composition $w\circ j: T \to I_{x_{\fks}}$ is the canonical embedding  (\ref{eq:T to I_x}).
 		\item An isomorphism $w' : I_\phi \isom I_{\phi (\fks')}$ of the form $\Int(g)$ for some $g\in G(\Qbar)$ such that $\Int(g) \circ \phi= \phi(\fks')$, satisfying that the composition $w'\circ j': T \to I_{\phi (\fks')}$ is the canonical embedding $i'$.
 	\end{itemize}
 	We claim that we can choose the above objects such that $(T,h,i,i')$ is a special fork. Indeed, the canonical embeddings $T \to I_{x_{\fks}}$ and $T \to I_{\phi(\fks')}$ are related by the composite isomorphism 
 	\begin{align}\label{eq:candidate isom}
 		I_{x_{\fks}}   \underset{w^{-1}} {\isom}  I_x  \underset{f}{\isom} I_{\phi} \underset{w'}{\isom} I_{\phi(\fks')}   ,	\end{align}
 	which lies in the canonical $I_{\phi(\fks')}^{\ad}(\QQ)$-conjugacy class of isomorphisms as in Lemma \ref{lem:tau for amicable pair}. Thus by the discussion in \S \ref{subsubsec:marking}, the canonical $G(\QQ_l)$-conjugacy classes of embeddings $I_{x_{\fks}, \QQ_l} \to G_{\QQ_l}$ and  $I_{\phi (\fks'), \QQ_l} \to G_{\QQ_l}$ give rise to the same $G(\ol \QQ_l)$-conjugacy class of  embeddings $T_{\ol \QQ_l} \to G_{\ol\QQ_l}$. It follows that $i$ and $i'$ are conjugate by $G(\Qlbar)$, and hence by $G(\Qbar)$. By the real approximation theorem, we can choose $a \in G(\Q)$ such that $\Int(a) \circ i' \circ h$ lies in the same connected component of $X$ as $i \circ h$. By replacing $i'$ by $\Int(a)\circ i'$, we can arrange that $(T,h,i,i')$ is a special fork, without affecting the existence of $w$ and $w'$. This proves the claim. Now note that (\ref{eq:candidate isom}) lies in the canonical class and intertwines $i$ and $i'$, so it is an $(i,i')$-compatible isomorphism.  
 \end{proof}

 \begin{thm}\label{thm:matchable pair 1}	Let $(\II,\JJ)$ be an amicable pair. Let $
 	\phi \in \JJ$ and let $T\subset I_{\phi}$ be a maximal torus. The image of $\tau^{\fkH}(\II,\JJ)$ under $\fkH(\phi) \to \fkH(\phi,T)$ is trivial.
 \end{thm}
 \begin{proof}
 	This follows directly from Proposition \ref{prop:LD1} and Lemma \ref{lem:reduction to special fork}. 
 \end{proof}
 \subsection{}\label{para:intro of twisting}
 Throughout the rest of this chapter, we write $\ica$\index[Oct2026-sorted-notation]{i 0 06 00 i@$\mathbb I$} for the set of isogeny classes in $\Shh_{ K_p}(\Fpbar)$ and $\cca = \AM/\mathrm{conj}$\index[Oct2026-sorted-notation]{j 0 06 00 j@$\mathbb J$} for the set of conjugacy classes of admissible morphisms $\Qf \to \G_G$.  
 
 Let $\II \in \ica$ and $\JJ \in \cca$. For any $x,x'\in \II$, the abelian groups $\Sha_G^{\infty}(\QQ, I_{x})$ and $\Sha_G^{\infty} (\QQ, I_{x'})$ are canonically isomorphic, where the isomorphism is induced by any isomorphism $I_x \isom I_{x'}$ arising from an element of $I_{x,x'}(\QQ)$, since they can be defined using abelianized Galois cohomology. We denote the identified group by\index[Oct2026-sorted-notation]{zzz sha 0 00 00 sha subscript g superscript infty i@$\Sha_G^{\infty} (\mathscr I)$} $\Sha_G^\infty(\II)$. Similarly, we define\index[Oct2026-sorted-notation]{zzz sha 0 00 00 sha subscript g superscript infty j@$\Sha_G^{\infty} (\mathscr J)$} $\Sha_G^\infty (\JJ)$, where the allowed isomorphisms $I_\phi \isom I_{\phi'}$ are those induced by $\Int(g)$ for $g\in G(\Qbar)$ such that $\Int(g) \circ \phi = \phi'$. 
 
 For $\beta\in \Sha_G^\infty (\II)$, we have the twisted isogeny class\index[Oct2026-sorted-notation]{i 0 05 00 i superscript beta@$\mathscr I^{\beta}$} $\II^{\beta} \in \ica$  constructed in \cite[\S 4.4.7, Prop.~4.4.8]{kisin2012modp}. The class $\beta$ is uniquely determined by $\II$ and $\II^{\beta}$ as it is the class of the $I_{x_1}$-torsor $I_{x_1,x_2}$ for any $x_1 \in \II, x_2 \in \II^\beta$. Moreover, for such $x_1$ and $x_2$ there is a canonical equivalence class of inner twistings between the $\QQ$-groups $I_{x_1}$ and $I_{x_2}$ induced by elements of $I_{x_1, x_2}(\Qbar)$. Since $\beta$ is trivial at infinity, the induced equivalence class of inner twistings between the $\RR$-groups $I_{x_1, \RR}$ and $I_{x_2,\RR}$ is trivial, i.e., it contains an $\RR$-isomorphism. We thus obtain a canonical isomorphism $\Sha_G^\infty (\II) \isom \Sha_G^\infty (\II^\beta)$, which depends only on $\II$ and $\II^\beta$. 
 
 Analogously,  for $\beta\in \Sha_G^\infty (\JJ)$, we have the twisted conjugacy class\index[Oct2026-sorted-notation]{j 0 05 00 j superscript beta@$\mathscr J^{\beta}$} $\JJ^\beta\in \cca$ as in  Definition \ref{defn:twist morphism} and Proposition \ref{prop:twist adm morph}. Again $\beta$ is uniquely determined by $\JJ$ and $\JJ^{\beta}$, and we have a canonical isomorphism $\Sha_G^\infty (\JJ) \isom \Sha_G^\infty (\JJ^{\beta})$ which depends only on $\JJ$ and $\JJ^{\beta}$.

 \subsection{}   We take this opportunity to correct a mistake in \cite[\S 3.1]{kisin2010integral}, and we freely use the notation 
 introduced there. The twisting of isogeny classes discussed above uses a twisting construction defined in \cite[\S3.1]{kisin2010integral}. The statement of Lemma 3.1.5 of {\em loc.~cit.~}should include the condition that the $Z$-torsor $\mathcal P$ is trivial over $\mathbb R$, i.e., that $\mathcal P_{\RR}$ is a trivial $Z_{\RR}$-torsor. The result is not true without this condition, as the $\Q$-isogeny $\lambda^{\mathcal P}$ 
 need not be a weak polarization. 
 
 This error was imported into \cite[\S 4.1.6]{kisin2012modp} and \cite[Lem.~4.4.8]{KisinPappas} via citation. 
 Fortunately, in all instances where this construction is applied to Shimura varieties in these papers, the condition of triviality at 
 $\infty$ holds. Moreover, the result is always applied to give a moduli theoretic description of a construction defined via complex uniformization. Thus, logically, the fact that $\lambda^{\mathcal P}$ is a weak polarization is never used, but rather follows 
 {\em a posteriori} in all these cases.
 
 Nevertheless, let us explain why $\lambda^{\mathcal P}$ is a weak polarization if $\mathcal P_{\mathbb R}$ is trivial. 
 Recall that for any $\Q$-algebra $R,$ an $R$-isogeny 
 $\phi: \cA \rightarrow \cA^*$ is an element of $\Hom(\cA,\cA^*)\otimes R$ which has an inverse in $\Hom(\cA^*,\cA)\otimes R,$ cf.~\cite[\S 9]{kottwitz1992points}. If $R$ is a subring of $\RR,$ we say that $\phi$ is an $R$-polarization if it is an 
 $R$-linear combination of polarizations, with positive coefficients. We say $\phi$ is a weak $R$-polarization if $\phi$ or $-\phi$ is an $R$-polarization. Thus a weak $\Q$-polarization is the same thing as a weak polarization. 
 Now if $\phi_0: \cA \rightarrow \cA^*$ is a 
 $\Q$-isogeny, then $\phi_0$ is a weak polarization if and only if it is a weak $\RR$-polarization when viewed as an 
 $\RR$-isogeny; this follows from the fact that the set of $\Q$-polarizations is a convex cone in  $\Hom(\cA,\cA^*)\otimes \Q.$ Moreover, the set of weak $\RR$-polarizations is stable under multiplication by $\RR^\times.$ 
 
 We now return to the explanation that $\lambda^{\mathcal P}$ is a weak polarization if $\mathcal P_{\mathbb R}$ is trivial. 
 Thus, suppose $\mathcal P_{\mathbb R}$ is trivial, and let $x \in \mathcal P(\mathbb R).$ Specializing the commutative 
 diagram \cite[(3.1.6)]{kisin2010integral} by $x,$ we see that $f_c(x)^{-1}\lambda^{\mathcal P}$ is a weak $\RR$-polarization. 
 Thus, by what we just saw, $\lambda^{\mathcal P}$ is a weak $\R$-polarization and hence a weak polarization.

 \subsection{}
 \label{subsubsec:setting for two thms} Let $(\II, \JJ)$ be a weakly amicable pair. Let $x\in \II$ and $\phi \in \JJ$. Recall from Lemma \ref{lem:tau for amicable pair} that we have an isomorphism $I_x \isom I_\phi$ which is canonical up to $I_{\phi}^{\ad}(\Q)$-conjugation. This induces a canonical isomorphism of abelian groups $\Sha_G^{\infty} (\II) \isom \Sha_G^\infty (\JJ)$ (which is independent of $x$ and $\phi$). Let $\beta$ be an element of the identified group. By \S \ref{subsec:frak H}, the abelian groups $\fkH(\phi)$ and $\fkH(\phi^{\beta})$ are canonically identified, since $\phi \approx \phi^\beta$. (Here $\phi^\beta$ is defined up to conjugacy, but this ambiguity does not affect the well-definedness of $\fkH(\phi^{\beta})$.)
 
 \begin{thm}[cf.~Proposition \ref{introprop:equivariance} in the Introduction]\label{thm:matchable pair 2} If $(\II, \JJ)$ is an amicable pair, then $(\II^{\beta}, \JJ^{\beta})$ is also amicable. Moreover, the elements
 	$$ \tau^{\fkH}(\II, \JJ)  \in \fkH(\phi)   \quad  \text{and} \quad  \tau^{\fkH}(\II^{\beta}, \JJ^{\beta})  \in \fkH(\phi^{\beta})  $$   correspond to each other under the canonical identification 	
 	$  \fkH(\phi)  \cong   \fkH(\phi^{\beta}) $. \end{thm}
 \begin{proof}
 	Let $T$ be a maximal torus in $I_{\phi }$ such that $\beta$ comes from a class $\beta_T$ in $\coh ^1(\QQ, T)$; such $T$ always exists by \cite[Thm.~5.10]{borovoi}. Let $j':T \hookrightarrow I_{\phi}$ be the inclusion map. By Lemma \ref{lem:reduction to special fork}, we may assume that $\II = \II_{\fks}$ and $\phi = \phi(\fks')$ for an amicable special fork $(\fks, \fks')$, that $j'$ is the canonical embedding (i.e.~$i'$), and that there is an $(i,i')$-compatible isomorphism $f:  I_{x_{\fks}} \isom I_{\phi}$.

 	By \cite[Lem.~4.4.5]{kisin2012modp}, the map of pointed sets $\coh^1(\RR, T) \to \coh^1(\RR, I_{\phi})$ induced by $j'$ has trivial kernel. Hence $\beta _T$ lies in $\Sha^{\infty}_{G} (\QQ, T)$ (defined as in \S \ref{subsubsec:T^der for prelinking}). Let  $\tilde T$  be as in \S \ref{para:setting for LD2}.   By Proposition \ref{prop:lift to sc}, there exists $\beta_{0} \in \Sha^{\infty}(\QQ, \tilde T)$ lifting $\beta_T$. As in \S \ref{para:setting for LD2}, we choose $\beta_0(\cdot), \lambda,\lambda'$ with respect to $\beta_0$ and obtain a twisted special fork $(T,h,i_1,i_1')$. The same argument as in the proof of \cite[Cor.~4.6.5]{kisin2012modp} shows that $(\II^{\beta}, \JJ^{\beta}) = ( \II_{(T,i_1,h)}, \JJ_{(T,i_1',h)})$. Hence the current theorem follows from Proposition \ref{prop:LD2}.
 \end{proof}

 \begin{lem}\label{lem:twist twice phi}
 	Let $\II \in \ica$ and $\JJ\in \cca$. For $\beta, \beta'\in \Sha_G^{\infty} (\II)$, we have $\II^{\beta+\beta'} = (\II^{\beta})^{\beta'}$. Here we use the canonical identification $\Sha_G^\infty(\II) \cong \Sha_G^\infty(\II^{\beta})$ to define the twist by $\beta'$ on the right hand side. Similarly, for $\beta, \beta'\in \Sha_G^{\infty} (\JJ)$, we have $\JJ^{\beta+\beta'} = (\JJ^{\beta})^{\beta'}$. 
 \end{lem}
 \begin{proof} We prove the two statements simultaneously, following almost the same ideas. Let $x\in \II$ and $\phi \in \JJ$. By \cite[Thm.~5.10]{borovoi}, there exists a maximal torus $T\subset I_x$ (resp.~$T \subset I_{\phi}$) such that $\beta$ and $\beta'$ come from elements $\beta_T$ and $\beta'_T$ of $\coh^1(\Q, T)$. By Theorem \ref{thm:geometric speciality} (resp.~Theorem \ref{thm:precise speciality}), we reduce to the case where $x= x_{(T,i,h)}$ (resp.~$\phi = \phi(T,i,h)$) for some $(T,i,h) \in \spd(G,X)$, and  the inclusion $T \hookrightarrow I_x$ (resp.~$T\hookrightarrow I_\phi$) is the canonical embedding as in (\ref{eq:T to I_x}) (resp.~the canonical embedding $i$). By the same argument as in the proof of Theorem \ref{thm:matchable pair 2}, we have $\beta_T,\beta'_T \in \Sha_G^{\infty} (\Q, T)$. Now fix cocycles $\beta_T(\cdot)$ and $\beta_T'(\cdot)$ representing $\beta_T$ and $\beta_T'$, and find $ \lambda_1 ,\lambda_2\in G(\Qbar)$ such that
 	$$ \lambda_1^{-1} \lix^{\rho} \lambda_1 = i(\beta_T(\rho)), \quad \lambda_2^{-1} \lix^{\rho} \lambda_2 = i( \beta_T(\rho) \beta'_T(\rho)) , \quad \forall \rho \in \Gamma. $$
 	Define $i_1 = \Int(\lambda_1) \circ i$ and $i_2 = \Int(\lambda_2) \circ i$. Then $(T,i_1,h) , (T,i_2,h) \in \spd(G,X)$. For the twisting of isogeny classes, the same argument as in the proof of \cite[Cor.~4.6.5]{kisin2012modp} shows that $\II^\beta = \mathscr I_{(T,i_1,h)}$ and $\mathscr I^{\beta +\beta'} = \mathscr I_{(T,i_2,h)}$. For the twisting of admissible morphisms, one directly verifies that   $\JJ^\beta = \JJ_{(T,i_1,h)}$ and $\JJ^{\beta +\beta'} = \JJ_{(T,i_2,h)}$.
 	
 	To finish the proof, we need to show that $$\II_{(T,i_2,h)}=\II_{(T,i_1, h)} ^{\beta'}, \quad \text{resp.}\quad \JJ_{(T,i_2,h)}=\JJ_{(T,i_1, h)} ^{\beta'}. $$ It suffices to find $\lambda \in G(\Qbar)$ satisfying the following conditions:
 	\begin{enumerate}
 		\item We have $\lambda^{-1} \lix^{\rho} \lambda = i_1 (\beta_T'(\rho)) ,~\forall \rho \in \Gamma$.
 		\item We have $i_2 = \Int (\lambda) \circ i_1$.
 	\end{enumerate}
 	Set $\lambda = \lambda_2 \lambda_1^{-1}$. Then it clearly satisfies (ii). To check (i), we compute
 	\begin{align*}
 		\lambda^{-1} \lix^{\rho} \lambda & = \lambda_1 \lambda_2^{-1} \lix^{\rho} \lambda _2 \lix^{\rho} \lambda_1^{-1} = \lambda_1 i \big(\beta_T(\rho) \beta_T'(\rho) \big) \lix^{\rho} \lambda_1^{-1} \\ &  = \lambda_1 \lambda_1^{-1} \lix^{\rho} \lambda_1 i (\beta_T'(\rho)) \lix^{\rho} \lambda_1^{-1} = \lix^{\rho} \lambda_1 i (\beta_T'(\rho)) \lix^{\rho} \lambda_1^{-1},
 	\end{align*}
 	where the second and third equalities follow from the defining properties of $\lambda_2$ and $\lambda_1$ respectively.
 	Since $i$ is defined over $\Q$, the above is equal to
 	$$ \lix^{\rho} \bigg( \lambda_1 i\big({\rho^{-1}} (\beta_T'(\rho)) \big) \lambda_1^{-1}\bigg) = \lix^{\rho} \bigg( i_1\big({\rho^{-1}} (\beta_T'(\rho)) \big)\bigg). $$ But $i_1$ is also defined over $\Q$, so the above is equal to $i_1(\beta_T'(\rho))$, as desired.
 \end{proof}

 \subsection{}\label{para:bij}
 In \S \ref{subsubsec:cocycle relation}, we defined an equivalence relation $\approx$ on the set $\AM$ of admissible morphisms $\Qf \to \G_G$. This descends to an equivalence relation on $\cca$, still denoted by $\approx$. For $\JJ_1, \JJ_2 \in \cca $, we have $\JJ_1 \approx \JJ_2$ if and only if $\JJ_2  = \JJ_1^{\beta}$ for a (unique) $\beta \in \Sha_G^{\infty} (\JJ_1)$. 
 Similarly, we define a binary relation $\approx$ on $\ica$ by declaring $\mathscr I_1 \approx \mathscr I_2$ when there exists $\beta \in \Sha_G^{\infty} (\II_1)$ such that $\mathscr I_2 = \mathscr I_1^{\beta}$. This is an equivalence relation: reflexivity is obvious, while symmetry and transitivity follow easily from Lemma \ref{lem:twist twice phi}.

 \begin{thm}[cf.~Theorem \ref{introprop:Hodge} in the Introduction]\label{thm:matchability of bij}There exists a   bijection
 	$\bij : \cca \isom \ica$
 	satisfying the following conditions:
 	\begin{enumerate}
 		\item  Whenever $\II = \bij (\JJ)$, the pair $(\II, \JJ)$ is amicable. In particular we have a canonical isomorphism $\Sha_G^\infty (\II) \cong \Sha_G^{\infty} (\JJ)$ (see \S \ref{subsubsec:setting for two thms}).
 		\item  Suppose  $\II= \bij(\JJ) $. For any  $\beta \in \Sha_G^\infty (\II) \cong \Sha_G^{\infty} (\JJ)$, we have $\bij (\JJ ^{\beta}) = \II^{\beta}.$
 	\end{enumerate}
 \end{thm}
 
 \begin{proof} 
 	For $\II \in \ica$ and $\JJ \in \cca$, the composite maps $\bar Y(\II) \to \KTstr \to \KT/{\sim}$ and  $\bar Y(\JJ) \to \KTstr \to \KT/{\sim}$ are both constant. Thus we have natural maps $\ica \to \KT/{\sim}$ and $\cca \to \KT/{\sim}$.   By Theorem \ref{thm:precise speciality}, Theorem \ref{thm:geometric speciality}, Proposition \ref{prop:special points on the geom side}, and Proposition \ref{prop:special pts on the gerb side}, the images of the last two maps are both the set of special elements (see \S \ref{para:special kott trip}), which we denote by $\mathfrak {SP}$.
 	By \cite[Lem.~4.4.11, Prop.~4.4.13, Lem.~4.5.6, Prop.~4.5.7]{kisin2012modp},  there is an equivalence relation $\approx$ on $\mathfrak{SP}$ whose pull-back to $\cca$ is $\approx$ and whose pull-back to $\ica$ is $\approx$. Hence we obtain natural bijections $$\cca/{\approx} ~\isom ~ \mathfrak{SP}/{\approx} ~ \isom ~ \ica/{\approx}. $$
 	It remains to construct a bijection from each $\approx$-equivalence class $A \subset \cca$ to the corresponding  $\approx$-equivalence class $B \subset \ica$ satisfying (i) and (ii). 
 	
 	Pick an arbitrary element $\JJ_0 \in A$. By Theorem \ref{thm:precise speciality}, there exists $\fks \in \spd(G,X)$ such that $\JJ_0 = \JJ_{\fks}$. Let $\II_0 = \II_{\fks}$. Then $\II_0$ is an element of $B$. The pair $(\II_0 , \JJ_0)$  is amicable by Corollary \ref{cor:standard amicable}. 
 	In particular, we have a canonical identification $ \Sha_G^{\infty} (\II_0) \cong \Sha_G^{\infty} (\JJ_0)$; we denote the identified group by $\Sha_0$. By the discussion in \S \ref{para:intro of twisting} and \S \ref{para:bij}, we have bijections $\Sha_0 \isom A, \beta \mapsto \JJ_0^\beta$ and $\Sha_0 \isom B, \beta \mapsto \II_0^{\beta}$. Composing them we obtain the bijection $A \isom B$. It satisfies (i) by Theorem \ref{thm:matchable pair 2}, and satisfies (ii) by Lemma \ref{lem:twist twice phi}. 
 \end{proof} 
 \begin{rem}\label{rem:difference in bij} 
 	In the proof of \cite[Cor.~4.6.5]{kisin2012modp}, a bijection $\cca \isom \ica$ is implicitly constructed by a different method than in Theorem~\ref{thm:matchability of bij}. There, one starts by selecting certain pairs $(\JJ_{\fks}, \II_{\fks})$ for $\fks \in \spd(G,X)$ to belong to the graph, and then produces the full bijection by twisting \emph{only} by elements $\beta \in \Sha_G^{\infty}(\QQ, I_\phi)$ coming from $\coh^1(\QQ, Z_{I_{\phi}})$. As a result, every pair in the graph is known to be associated with a common special point datum (a property stronger than being amicable). However, condition~(ii) in Theorem~\ref{thm:matchability of bij} is not easily accessible from that approach.
 \end{rem} 
 
 \section{A strengthened version of the Langlands--Rapoport--\texorpdfstring{$\tau$}{tau} Conjecture} \label{subsec:LR str}
 In this section, we prove a strengthened version (Theorem \ref{thm:LR str}) of the Langlands--Rapoport--$\tau$ Conjecture in the current setting of Hodge type. The strengthening takes into consideration connected components and the equivariance under the action of $\mathscr{A}(G_{\mathbb{Z}_{(p)}})$, which are necessary structures for passing from Hodge type to abelian type, as well as gives more refined information about the twisting elements $\tau$.
 \subsection{}
 The reductive model $\mathcal G$ of $G_{\QQ_p}$ determines a reductive model $\mathcal G^{\ad}$ over $\Z_p$ of $G^{\ad}_{\Q_p}$, and we have an induced  unramified Shimura datum $(G^{\ad}, X^{\ad}, p, \cG^{\ad})$. 
 For any admissible morphism $\phi : \gQ \to \gG_G$ for $(G,X,p,\cG)$, the composition  $$\phi^{\ad} : \gQ \overset{\phi}{\To} \G_G \To \G_{G^{\ad}}$$ is an admissible morphism for $(G^{\ad}, X^{\ad}, p , \cG^{\ad})$. We denote by  $\AM(G)$ and $\AM(G^{\ad})$ respectively  the sets of admissible morphisms $\gQ\to \gG_G$ and $\gQ \to \gG_{G^{\ad}}$. Thus we have a natural map $\AM(G) \to \AM(G^{\ad}), \phi \mapsto \phi^{\ad}$. This further induces a map $\AM(G)/{\approx} \to \AM(G^{\ad})/ {\approx}$.\footnote{Despite the same notation $\approx$, the  equivalence relations on $\AM(G)$ and $\AM(G^{\ad})$ are defined separately.}

 Let  $\phi_0 \in \AM(G^{\ad})$. The inverse image\index[Oct2026-sorted-notation]{i 0 00 00 i subscript phi subscript 0 overline q superscript g@$I_{\phi_0, \overline \QQ}^G$} $I_{\phi_0, \Qbar}^G$ of $I_{\phi_0 ,\Qbar} \subset G^{\ad}_{\Qbar}$ in $G_{\Qbar}$ is a reductive group over $\Qbar$, since $I_{\phi_0 ,\Qbar}$ is the centralizer of a torus in $G^{\ad}_{\Qbar}$ (cf.~\S \ref{para:choice of q_tau}).  If $\phi \in \AM(G)$ is such that $\phi^{\ad}  = \phi_0$, then $I_{\phi}$ is a $\Q$-form of $I_{\phi_0, \Qbar}^G$. In particular, the natural $\Q$-map $I_\phi \to I_{\phi_0}$ is surjective with kernel $Z_G$, and it identifies the adjoint groups of $I_\phi$ and $I_{\phi_0}$.

 \begin{lem}\label{lem:lifting phi_0}
 	Let $\phi_0\in \AM(G^{\ad})$. Let\index[Oct2026-sorted-notation]{u 0 03 00 u phi subscript 0@$\mathcal U(\phi_0)$}  $\mathcal U(\phi_0)$ be the set of conjugacy classes of $\phi \in \AM(G)$ such that $\phi^{\ad}$ is conjugate to $\phi_0$. Then $\mathcal U(\phi_0)$ is non-empty and is acted on transitively by the abelian group $\ker(\Sha ^{\infty} (\QQ, Z_G) \to \coh^1(\Q,G))$, where the action is given by the usual twisting construction (see Proposition \ref{prop:twist adm morph}).
 \end{lem}
 \begin{proof}
 	As in the proof of \cite[Lem.~3.4.8]{kisin2012modp}, the group in the lemma is equal to  $$\Sha^{\infty} (\QQ, Z_G) \cap \im \big(\coh^1(\QQ, Z_{ G_{\sconn}}) \to \coh ^1(\QQ, Z_G)  \big). $$ The lemma follows from this observation and \cite[Lem.~3.4.8, Prop.~3.4.11]{kisin2012modp}.
 \end{proof}

 \begin{lem}\label{lem:bij AM} The natural map $\AM(G)/{\approx} \to \AM(G^{\ad})/{\approx}$ is a bijection.
 \end{lem}
 \begin{proof}
 	The surjectivity follows from Lemma \ref{lem:lifting phi_0}. We show injectivity. Let $\phi, \phi' \in \AM(G)$ be such that $\phi^{\ad} \approx \phi^{\prime ,\ad}$. Then  there exists $\beta _0 \in \Sha_{G^{\ad}}^{\infty} (\QQ, I_{\phi^{\ad}})$ such that $\phi^{\prime, \ad}$ is conjugate to $(\phi^{\ad}) ^{\beta_0}$. The reductive group $\tilde I_\phi$ equipped with a pre-inner twisting into $G_{\sconn}$ associated with the pre-inner twisting $I_{\phi, \Qbar} \to G_{\Qbar}$ (see \S \ref{para:underline isom in pro case}, cf.~\S \ref{para:choice of q_tau}) is canonically identified with the analogous object associated with $I_{\phi^{\ad}, \Qbar} \to G^{\ad}_{\Qbar}$. By Proposition \ref{prop:lift to sc}, the natural maps $\Sha^{\infty}_{\ab} (\QQ, \tilde I_\phi) \to \Sha_{G}^{\infty} (\QQ, I_{\phi}) $ and $\Sha^{\infty}_{\ab} (\QQ, \tilde I_\phi) \to \Sha_{G^{\ad}}^{\infty} (\QQ, I_{\phi^{\ad}}) $ are surjective. Hence  
 	$  \Sha_{G}^{\infty} (\QQ, I_{\phi}) \to    \Sha_{G^{\ad}}^{\infty} (\QQ, I_{\phi^{\ad}}) $ is surjective. Thus we can find $\beta \in \Sha_G^{\infty} (\QQ, I_{\phi})$ lifting $\beta_0$. Then $(\phi^{\beta})^{\ad}$ is conjugate to $\phi^{\prime, \ad}$.\footnote{Here $\phi^{\beta}$ is only well defined up to conjugacy, but the ambiguity does not affect the conjugacy class of $(\phi^{\beta})^{\ad}$.} By the transitivity in Lemma \ref{lem:lifting phi_0}, we have $\phi^{\beta} \approx \phi'$. But $\phi \approx \phi^{\beta}$, so $\phi \approx \phi'$ as desired. 
 \end{proof}

 \subsection{} \label{subsec:A(G)}
 For $\phi_0 \in \AM(G^{\ad})$, let\index[Oct2026-sorted-notation]{i 0 00 00 i subscript phi subscript 0 q bar superscript natural@$I_{\phi_0} (\Qbar)^{\natural}$}  $I_{\phi_0} (\Qbar)^{\natural}$ 
 denote the inverse image of $I_{\phi_0}(\Q)$  in $I_{\phi_0 ,\Qbar}^G(\Qbar)$. Then $I_{\phi_0} (\Qbar)^{\natural}$ acts by left multiplication on the set 
 $$\coprod_{ \phi \in \AM(G), \phi^{\ad} = \phi_0} X(\phi) . $$ Moreover, for any $\tau \in I_{\phi_0}^{\ad} (\A_f)$, we can twist the action by $\Int(\tau)$, and denote the quotient for the twisted action by 
 \index[Oct2026-sorted-notation]{s 0 00 00 s subscript tau g phi subscript 0@$S_{\tau} (G, \phi_0)$}
 $$ S_{\tau} (G, \phi_0) : =  I_{\phi_0} (\Qbar)^{\natural}  \backslash  \coprod_{ \phi \in \AM(G), \phi^{\ad} = \phi_0} X(\phi) .$$
 Here the twisted action is given by 
 $$ X(\Int(a^{-1})\phi) \To X(\phi), \quad x \longmapsto \Int(\tau)(a)  x = \Int(\tau)(a) a^{-1} ax, \quad a \in I_{\phi_0} (\Qbar)^{\natural},  $$  and it is well defined since the element $\Int(\tau)(a) a^{-1} \in I_{\phi_0, \Qbar}^G(\bar \A_f)$ lies in $I_{\phi}(\A_f)$ for any $\phi \in \AM(G)$ such that $\phi^{\ad}=\phi_0$. 
 The action of $\Phi^{\Z} \times G(\A_f^p)$ on $\coprod X(\phi)$ descends to $S_{\tau}(G,\phi_0)$. In addition, the maps $\tilde c_{\phi}: X(\phi) \to \pi(G,X)$ induce a map $S_{\tau}(G,\phi_0) \to \pi(G,X)$; see \cite[Lem.~3.7.4]{kisin2012modp}. 
 
 For each $\phi \in \AM(G)$ such that $\phi^{\ad} = \phi_0$, we can view $\tau \in I_{\phi_0}^{\ad}(\A_f)$ as an element of $I_{\phi}^{\ad}(\A_f)$, and as such define $S_{\tau}(\phi)$ as in \S \ref{subsubsec:defn of S(phi)} (cf.~Lemma \ref{lem:HT free}).  By \cite[Cor.~3.4.16]{kisin2012modp}, we have a natural bijection 
 $$ \coprod_{\phi \in \mathcal U(\phi_0)} S_{\tau}(\phi) \isom S_{\tau} (G, \phi_0).  $$ This bijection is $\Phi^{\Z} \times G(\A_f^p)$-equivariant, and commutes with the maps from the two sides to $\pi(G,X)$. 
 
 As in \cite[\S 3.7.3]{kisin2012modp}, we define \index[Oct2026-sorted-notation]{a 0 05 00 a g subscript z subscript p@$\mathscr A(G_{\ZZ_{(p)}})$} 
 \begin{align*}
 	\mathscr A(G_{\ZZ_{(p)}})     & = G(\A_f^p) * _{G(\Z_{(p)})_+ } G^{\ad}(\Z_{(p)})^+ \\ &  : =G(\A_f^p) \rtimes G^{\ad}(\Z_{(p)})^+  /  \{(\text{image of }\gamma,~ \text{image of }\gamma^{-1}) \mid \gamma \in G(\ZZ_{(p)} )_+ \} ,
 \end{align*}  where \index[Oct2026-sorted-notation]{g 0 00 00 g z subscript p subscript@$G(\mathbb Z_{(p)})_+$} \index[Oct2026-sorted-notation]{g 0 00 00 g superscript ad z subscript p superscript@$G^{\ad}(\mathbb Z_{(p)})^+$} $$  G(\Z_{(p)})_+  :  = \mathcal G(\Z_p) \cap G(\Q)_+,  \quad  G^{\ad}(\Z_{(p)})^+ : =  \mathcal G^{\ad}(\Z_p) \cap G^{\ad}(\QQ)^+,$$ and the action of $G^{\ad}(\Z_{(p)})^+ $ on $G(\A_f^p)$ is the restriction of the adjoint action of $G^{\ad}(\Q)$. By \cite[Lem.~3.7.2]{kisin2012modp}, $\mathscr A(G_{\ZZ_{(p)}}) $ acts on $S_{\tau}(G,\phi_0)$ on the right, where the factor $G(\A_f^p)$ acts in the way already discussed above, and each  $g\in G^{\ad}(\Z_{(p)} )^+$
 acts by right multiplication by $\tilde g$, where $\tilde g$ is any lift of $g$ in $G(\Qbar)$. The action of $\mathscr A(G_{\ZZ_{(p)}}) $ commutes with that of $\Phi^{\Z}$. 
 
 Taking the disjoint union of the maps $S_\tau(\phi) \to \pi(G,X)$ (see \S \ref{subsec:map to pi}) over $\phi \in \mathcal U(\phi_0)$, we obtain a map $S_\tau (G,\phi_0)  \to \pi(G,X) $. 
 By \cite[Lem.~3.7.4]{kisin2012modp}, this map is $\mathscr A(G_{\ZZ_{(p)}})$-equivariant, where $\mathscr A(G_{\ZZ_{(p)}})$ acts on the target by the natural map $\mathscr A(G_{\ZZ_{(p)}}) \to \pi (G)$ sending $(h,g) \in \mathscr A(G_{\ZZ_{(p)}}) = G(\A_f^p) * _{G(\Z_{(p)})_+ } G^{\ad}(\Z_{(p)})^+$ to the class of $h$.
 
 In summary,  $S_{\tau}(G,\phi_0)$ is a set equipped with 
 \begin{itemize}
 	\item an action of $\Phi^{\ZZ} \times \mathscr A(G_{\ZZ_{(p)}})$,
 	\item an $\mathscr A(G_{\ZZ_{(p)}})$-equivariant map $S_{\tau}(G,\phi_0) \to \pi(G,X)$. 
 \end{itemize} 
 
 \subsection{} In \cite[\S 4.2.8]{kisin2012modp}, the right action of $\mathscr A (G_{\ZZ_{(p)}})$ on the scheme $\Shh_{K_p}$ is defined. The factor $G(\A_f^p)$ acts in the usual way, while the $G^{\ad}(\ZZ_{(p)} )^+$-action is induced, via the extension property of the integral model, by the action on the generic fiber defined by the twisting construction on abelian varieties. 
 The $\mathscr A (G_{\ZZ_{(p)}})$-action on $\Shh_{K_p} (\Fpbar)$ permutes the isogeny classes. Specifically, $G(\A_f^p)$ stabilizes each isogeny class, while $\omega \in G^{\ad}(\ZZ_{(p)})^+$ sends each isogeny class $\mathscr I_x$ to $\mathscr I_x^{\beta}$ (with notation as in \S \ref{subsubsec:setting for two thms}) where $\beta$ is the image of $\omega^{-1}$ under the composite map $$ G^{\ad}(\Z_{(p)})^{+} \To \coh^1(\Q, Z_G)  \To \coh^1(\Q, I_x). $$(Note that $\beta \in \Sha_G^\infty(\Q, I_x)$.) Moreover, the natural map $\Shh_{K_p}(\Fpbar) \to \pi(G,X)$ is $\mathscr A (G_{\ZZ_{(p)}})$-equivariant, for the same $\mathscr A (G_{\ZZ_{(p)}})$-action on $\pi(G,X)$ as in \S \ref{subsec:A(G)}; see \cite[\S 3.7.6]{kisin2012modp}.

 For each isogeny class $\mathscr I \subset \Shh_{K_p} (\Fpbar ) $, define \index[Oct2026-sorted-notation]{i 0 05 00 i subscript 0@$\mathscr I_0$}
 $$ \mathscr I_{0} : = \bigcup _{g \in G^{\ad}(\QQ)^+} \mathscr I^{[g]} \subset \Shh_{K_p} (\Fpbar ),$$
 where $[g]$ denotes the image of $g$ in $\Sha_G^\infty(\QQ, I_x)$. Then $\mathscr I_{0}$ is a subset of $\Shh_{K_p} (\Fpbar )$ stable under the $\mathscr A(G_{\ZZ_{(p)}})$-action. Hence it is also a set equipped with  
 \begin{itemize}
 	\item an action of $\Phi^{\ZZ} \times \mathscr A(G_{\ZZ_{(p)}})$,  
 	\item an $\mathscr A(G_{\ZZ_{(p)}})$-equivariant map $\mathscr I_{0}  \to \pi(G,X)$. 
 \end{itemize} 
 
 \begin{prop}\label{prop:equiv bij}
 	Let $(\mathscr I, \mathscr J)$ be an amicable pair, and $\phi \in \mathscr J$. Let $\tau \in I_\phi ^{\ad} (\A_f)$ be a lift of $\tau^{\fkH}(\II, \JJ)  \in \fkH(\phi) $. Then there exists a $\Phi^{\ZZ} \times \mathscr A(G_{\ZZ_{(p)}})$-equivariant bijection $\mathscr I_0 \isom S_{\tau}(G,\phi^{\ad})$ that commutes with the maps from the two sides to $\pi(G,X)$. 
 \end{prop}
 \begin{proof}
 	Let $x\in \mathscr I$. Under the canonical identification $\Sha_G^\infty (\QQ, I_\phi) \cong \Sha_G^\infty (\QQ, I_x)$, the images of the composite maps  
 	$$ G^{\ad}(\QQ)^+ \To \ker (\Sha^\infty(\QQ, Z_G) \to \coh^1(\QQ, G)) \To \Sha_G^\infty (\QQ, I_\phi)   $$ and $$ G^{\ad}(\QQ)^+ \To \ker (\Sha^\infty(\QQ, Z_G) \to \coh^1(\QQ, G)) \To \Sha_G^\infty (\QQ, I_x)   $$ are  identified. We denote this image by $G^{\ad}(\QQ)^{+\prime}$, which is an abelian group. The $\mathscr A (G_{\ZZ_{(p)}})$-action on $S_\tau(G,\phi^{\ad})$ (resp.~on $\mathscr I_0$) permutes the subsets $S_\tau(\phi') \subset S_\tau(G,\phi^{\ad})$ for $\phi' \in \mathcal U(\phi^{\ad})$ (resp.~the isogeny classes $\mathscr I^{[g]} \subset \mathscr I_0$ for $g \in G^{\ad}(\QQ)^+$.) There are bijections 
 	$$ G^{\ad}(\QQ)^{+\prime} \isom \mathcal U(\phi^{\ad}), \quad \beta \mapsto \phi^{\beta}$$ and $$ G^{\ad}(\QQ)^{+\prime} \isom  \{\text{isogeny classes $\mathscr I'  \subset \mathscr I_0$}\}, \quad \beta \mapsto \mathscr I^\beta. $$ For either of the above bijections, two elements of the right hand side are in the same $\mathscr A(G_{\ZZ_{(p)} })$-orbit if and only if the corresponding elements of $G^{\ad}(\QQ)^{+\prime}$ are in the same coset under $G^{\ad}(\ZZ_{(p)})^{+\prime} : = \im (G^{\ad}(\ZZ_{(p)})^{+} \to G^{\ad}(\QQ)^{+\prime})$. 
 	
 	Fix a set of representatives $\{\beta_\lambda\}_{\lambda\in \Lambda}$ of the $G^{\ad}(\ZZ_{(p)})^{+\prime}$-cosets in $G^{\ad}(\QQ)^{+\prime}$. 
 	Write $\phi_\lambda : = \phi^{\beta_\lambda}, \mathscr I_\lambda : = \mathscr I^{\beta_\lambda}$. It suffices to construct, for each $\lambda$, a $\Phi^{\ZZ}\times \mathscr A( G_{\ZZ_{(p)} })$-equivariant bijection 
 	$$\bigcup_{\xi \in \mathscr A( G_{\ZZ_{(p)} })}  \mathscr I_\lambda  \cdot \xi   \isom  \bigcup_{\xi \in \mathscr A( G_{\ZZ_{(p)} })}  S_{\tau}(\phi_\lambda) \cdot \xi   $$ over $\pi(G,X)$. Here the unions are not disjoint. More precisely, write $H$ for the kernel of $G^{\ad}(\Z_{(p)})^{+} \to G^{\ad}(\Z_{(p)})^{+\prime} $. Then  $S_{\tau}(\phi_\lambda)$ and $\mathscr I_{\lambda}$ have the same stabilizer $ G(\A_f^p) * _{G(\Z_{(p)})_+ }  H $ in $\mathscr A( G_{\ZZ_{(p)} })$.  It remains to construct a $\Phi^{\ZZ} \times  (G(\A_f^p) * _{G(\Z_{(p)})_+ }  H)$-equivariant bijection  $  \mathscr I_\lambda \isom S_{\tau}(\phi_\lambda) $ over $\pi(G,X)$, as then induction from $G(\A_f^p) * _{G(\Z_{(p)})_+ }  H$ to $\mathscr A( G_{\ZZ_{(p)} })$ would give the desired bijection.  
 	
 	Write $\mathscr J_\lambda$ for the conjugacy class of $\phi_\lambda$. By Theorem \ref{thm:matchable pair 2}, $(\mathscr I_\lambda, \mathscr J_\lambda)$ is an amicable pair, and the element $\tau \in I_{\phi}^{\ad}(\A_f) \cong I_{\phi_\lambda}^{\ad}(\A_f)$ is a lift of $\tau^{\fkH}(\mathscr I_\lambda, \mathscr J_\lambda)$. The latter property implies, in view of Lemma \ref{lem:tau for amicable pair}, that $\tau = \tau _{\fkm}$ for some $\fkm \in \tmar(x_\lambda,\phi_\lambda), x_\lambda \in \mathscr I_\lambda$. Then by Proposition \ref{prop:marked bij}, associated with $\mathfrak m$ we have a $\Phi^{\ZZ} \times G(\A_f^p)$-equivariant  bijection $$\mathscr B: \mathscr I_\lambda \isom S_\tau (\phi_\lambda)$$ over $\pi(G,X)$. It remains to show that $\mathscr B$ is equivariant under $H$. 
 	
 	Let $x_\lambda' \in \II_\lambda$ and $h \in H$ be arbitrary, and we show that $  \mathscr B(x_\lambda') \cdot h = \mathscr B( x_\lambda' \cdot h) $. On choosing an isomorphism $x_\lambda\isom x'_{\lambda}$ in the category $\II_\lambda$, we have a bijection $\tmar(x_\lambda,\phi_\lambda) \isom \tmar(x_\lambda', \phi_\lambda)$. The bijection $\II_{\lambda} \isom S_\tau (\phi_\lambda)$ associated with the image of $\fkm \in \tmar(x_\lambda,\phi_\lambda)$ in $\tmar(x_\lambda', \phi_\lambda)$ is the same as the original bijection $\mathscr B$. Hence we may assume that $x_\lambda' = x_\lambda$.   We identify $I_{x_\lambda}$ with $I_{\phi_\lambda}$ as in \S \ref{subsubsec:marking}. Let $\tilde h  \in  G(\Qbar)$ be a lift of $h$, and let $i_h \in I_{x_\lambda}(\Qbar)$ be such that $i_h^{-1} \rho(i_h)= \tilde h^{-1} \rho(\tilde h)$ for all $\rho \in \Gal(\Qbar/\QQ)$. Let $\xi, \xi' , \mathbf  f$ be as in \S \ref{subsubsec:markingII} with respect to our current $x_{\lambda}, \phi_\lambda, \fkm = (y,y')$. Then the proof of \cite[Prop.~4.4.14]{kisin2012modp}  shows that $x_\lambda$ and $x_{\lambda} \cdot h$ have lifts in $Y(x_\lambda)/\cG(\Zpur)$ whose images under $\xi : Y(x_\lambda)/\cG(\Zpur) \isom G(\A_f^*) / \cG(\Zpur)$ are of the form $g$ and $\iota(i_h^{-1}) g \tilde h$ respectively, for some $g = (g^p, g_p)\in G(\A_f^*)$. Here $\iota$ is the embedding $I_{x_{\lambda}, \A_f} \to G_{\A_f^p} \times \Res_{\Qpur/ \Q_p} G$ associated with $y$. We view both $\iota(i_h^{-1})$ and $\tilde h$ as elements of $G(\A_f^* \otimes _{\QQ} \Qbar)$, so $\iota(i_h^{-1}) g \tilde h$ is \textit{a priori} an element of $G(\A_f^* \otimes _{\QQ} \Qbar)$, but it turns out to be in $G(\A_f^*)$.      
 	
 	Since $\iota_p(i_h^{-1})$ commutes with $c_{y'_p}$, we  have 
 	$$\mathbf f ( \iota (i_h^{-1}) g_p \tilde h)  =  \iota (i_h^{-1})  \cdot \mathbf f( g_p)  \cdot \tilde h   . $$ Write $u$ for $(g^p, \mathbf f (g_p)) \in G(\A_f^*)$. Then the above discussion shows that $\mathscr B(x_\lambda)$ and $\mathscr B( x_{\lambda} \cdot h)$ have lifts $t_1$ and $t_2$ in $Y(\phi_\lambda)/\cG(\Zpur)$ whose images under $\xi' : Y(\phi_\lambda)/\cG(\Zpur) \isom G(\A_f^*) / \cG(\Zpur)$ are of the form $u$ and $\iota(i_h^{-1}) u \tilde h$ respectively. Thus $t_1 = y'u$ and $t_2 = y' \iota(i_h^{-1}) u \tilde h$. Since $\tau$ is associated with $\fkm = (y,y')$, we have $$\iota(i_h^{-1}) = \Int (y^{\prime -1}) \Int (\tau) (i_h^{-1}),$$  and so  $t_2 = \Int(\tau)(i_h^{-1}) \cdot t_1 \cdot \tilde h. $ As shown in the proof of \cite[Lem.~4.5.9]{kisin2012modp}, this formula for $t_2$ implies that the image of $t_2$ in $S_{\tau}(\phi_\lambda)$ is equal to the  image of $t_1$ acted on by $h$, as desired.  
 \end{proof} 
 \subsection{}
 For each  $\phi \in \AM(G)$, let\index[Oct2026-sorted-notation]{h 0 04 23 widetilde h phi@$\widetilde {\mathfrak H}(\phi)$} $$\widetilde{\fkH}(\phi) : = I_{\phi,\sconn} (\A_f) \backslash I_{\phi}^{\ad}(\A_f)/I_{\phi}^{\ad}(\Q).$$ In view of the triviality of $\coh^1(\A_f, I_{\phi,\sconn})$ and Lemma \ref{lem:KHC theorem}, the coboundary map $I_{\phi}^{\ad}(\A_f) \to \coh^1(\A_f, Z_{I_{\phi, \sconn}})$ induces an isomorphism $$\widetilde{\fkH}(\phi) \isom  \D (Z_{I_{\phi,\sconn}}, I_{\phi,\sconn}; \A_f) /\D (Z_{I_{\phi,\sconn}} , I_{\phi, \sconn} ;\Q) = \coh^1(\A_f, Z_{I_{\phi,\sconn}}) / \Sha^\infty (\Q, Z_{I_{\phi,\sconn}}). $$ Thus $\widetilde{\fkH}(\phi)$ is an abelian group and, similarly as in \S \ref{subsubsec:cocycle relation}, for any  $\phi ,\phi_1 \in\AM(G)$ with $\phi \approx \phi_1$ we have a canonical isomorphism $\widetilde {\fkH}(\phi) \cong \widetilde{\fkH} (\phi_1)$. In \S \ref{subsec:frak H} we have seen the similar property of the abelian group $\fkH(\phi)$. Hence the assignments $\phi \mapsto \widetilde{\fkH}(\phi)$ and $\phi \mapsto \fkH(\phi)$ define sheaves of abelian groups $\widetilde{\fkH}$ and $\fkH$ on $\AM(G)$ which descend to  sheaves $\widetilde{\fkH}_{\approx}$ and $\fkH_{\approx}$ on $\AM(G)/{\approx}$ respectively. Similarly as in Definition \ref{defn:system of obstructions} we write $\Gamma(\fkH)_0$ (resp.~$\Gamma(\widetilde {\fkH})_0$) for the group of those global sections of $\fkH$ (resp.~$\widetilde{\fkH}$) which descend to global sections of $\fkH_{\approx}$ (resp.~$\widetilde{\fkH}_{\approx}$). We have natural surjections of abelian sheaves $\widetilde{\fkH} \to \fkH \to \mathcal H$ which give rise to  surjections of abelian groups $\Gamma(\widetilde {\fkH})_0 \to \Gamma (\fkH)_0 \to \Gamma(\cH)_0.$
 
 \begin{thm}\label{thm:LR str}
 	Let $\{\phi_{0,\varepsilon}\}_{\varepsilon\in \mathcal E}$ be a set of representatives of the conjugacy classes in $\AM(G^{\ad})$. 
 	There exists a family $\{ \tau_{\varepsilon}\}_{\varepsilon\in \mathcal E}$ with $\tau_{\varepsilon} \in I_{\phi_{0,\varepsilon}}^{\ad} (\A_f)$ satisfying the following conditions:
 	\begin{enumerate}
 		\item For each $\phi\in \AM(G)$, let $\varepsilon$ be such that $\phi^{\ad}$ is conjugate  to $\phi_{0,\varepsilon}$, and let $\tilde \tau (\phi)$ be the image of  $\tau_{\varepsilon}$ under the composite map\footnote{Here the isomorphism $I_{\phi_{0,\varepsilon}}^{\ad} \cong I_\phi^{\ad}$ is canonical up to conjugation by $I_{\phi}^{\ad}(\Q)$, so this composite map is canonical.} $I_{\phi_{0,\varepsilon}}^{\ad}(\A_f) \cong I_{\phi}^{\ad}(\A_f) \to \widetilde{\fkH}(\phi)$.  Then the assignment $\phi \mapsto \tilde \tau(\phi)$ is an element of $\Gamma(\widetilde{\fkH})_0$. 
 		\item For every   $\phi\in \AM(G)$ and maximal torus $T \subset I_\phi$, the image of $\tilde \tau(\phi)$ under $\widetilde{\fkH} (\phi) \to \fkH(\phi) \to \fkH(\phi, T)$ is trivial. 
 		\item We have a $\Phi^{\ZZ} \times \mathscr A(G_{\ZZ_{(p)}})$-equivariant bijection 
 		$$  \Shh_{K_p} (\Fpbar ) \isom \coprod_{\varepsilon \in \mathcal E} S_{\tau_{\varepsilon}}(G, \phi_{0,\varepsilon})  $$ which commutes with the maps from the two sides to $\pi(G,X)$.   
 	\end{enumerate}
 \end{thm}
 \begin{proof}
 	Fix a bijection $\bij$ as in Theorem \ref{thm:matchability of bij}. For $\phi \in \AM(G)$, write  $\tau (\phi)$ for the element $\tau^{\fkH}(\bij(\JJ),\JJ) \in \fkH(\phi)$ as in Definition \ref{defn:matchable}, where $\JJ$ is the conjugacy class of $\phi$. By the two conditions satisfied by $\bij$ in Theorem \ref{thm:matchability of bij} and by Theorem \ref{thm:matchable pair 2}, the assignment $\phi \mapsto \tau (\phi)$ is  an element of $\Gamma(\fkH)_0$. Moreover, by Theorem \ref{thm:matchable pair 1}, the image of $\tau(\phi)$ in $\fkH(\phi, T)$ is trivial for every  maximal torus $T \subset I_\phi$. 
 	
 	Fix a lift  $(\phi \mapsto \tilde \tau(\phi)) \in \Gamma(\widetilde{\fkH})_0$ of $(\phi \mapsto \tau(\phi))$  along the surjection $\Gamma(\widetilde {\fkH})_0 \to \Gamma (\fkH)_0$. For each $\varepsilon \in \mathcal E$, we fix once and for all $\phi_{\varepsilon} \in \AM(G)$ such that $\phi_{\varepsilon}^{\ad} = \phi_{0,\varepsilon}$. Now let $\tau_{\varepsilon}$ be a preimage of $\tilde \tau(\phi_\varepsilon)$ along the canonical surjection $I_{\phi_{0,\varepsilon}} ^{\ad} (\A_f) \to \widetilde \fkH (\phi_\varepsilon)$ as in (i). Note that if $\phi$ is any other element of $\AM(G)$ such that $\phi^{\ad}$ is conjugate to $\phi_{0,\varepsilon}$, then $\phi \approx \phi_\varepsilon$ and the elements $\tilde \tau(\phi)$ and $\tilde \tau(\phi_{\varepsilon})$ are identified under the canonical isomorphism $\widetilde{\fkH}(\phi) \cong \widetilde{\fkH}(\phi_\varepsilon)$ (since $(\phi \mapsto \tilde \tau(\phi))$ lies in $\Gamma(\widetilde{\fkH})_0$). It follows that the image of $\tau_{\varepsilon}$ under $I_{\phi_{0,\varepsilon}} ^{\ad} (\A_f) \to \widetilde \fkH (\phi)$ is $\tilde \tau (\phi)$. Hence the family $\{\tau_\varepsilon\}$ satisfies (i) and (ii).  
 	
 	For each $\varepsilon$, let $\JJ_\varepsilon$ be the conjugacy class of $\phi_\varepsilon$, and let $\II_{\varepsilon} = \bij(\JJ_\varepsilon)$. Then $(\II_\varepsilon, \JJ_\varepsilon)$ is an amicable pair, and $\tau_\varepsilon \in I_{\phi_{0,\varepsilon}}^{\ad} (\A_f)$ is a lift of $\tau^{\fkH}( \II_\varepsilon, \JJ_\varepsilon)$. Hence by Proposition \ref{prop:equiv bij}, there is a $\Phi^{\ZZ} \times \mathscr A(G_{\ZZ_{(p)}})$-equivariant bijection 
 	$$ \mathscr B_\varepsilon:    \II_{\varepsilon,0}  \isom S_{\tau_{\varepsilon}}(G, \phi_{0,\varepsilon}) $$ over $\pi(G,X)$.  Now an isogeny class $\II$ in $\Shh_{ K_p}(\Fpbar)$ is contained in $\II_{\varepsilon,0}$ if and only if $\II = \II_{\varepsilon}^{[g]}$ for some $g \in G^{\ad}(\QQ)^+$. In view of condition (ii) in Theorem \ref{thm:matchability of bij}, this happens if and only if $\bij^{-1} (\II) = \JJ_{\varepsilon}^{[g]}$. But the $\JJ_{\varepsilon}^{[g]}$ for $g\in G^{\ad} (\QQ)^+$ are exactly the set of conjugacy classes in $\AM(G)$ which induce the conjugacy class of $\phi_{0,\varepsilon}$ (see \cite[Prop.~3.4.11]{kisin2012modp}, cf.~Lemma \ref{lem:lifting phi_0}). Since $\bij$ is a bijection and the map $\AM(G) \to \AM(G^{\ad})$ is surjective (see Lemma \ref{lem:lifting phi_0}), it follows that the disjoint union of $\II_{\varepsilon,0}$ over all $\varepsilon \in \mathcal E$ is exactly $\Shh_{K_p} (\Fpbar)$. We thus obtain the bijection in (iii) by taking the disjoint union of $\mathscr B_\varepsilon$ over all $\varepsilon \in \mathcal E$. 
 \end{proof}
 \begin{rem}
 	The weaker version of Theorem \ref{thm:LR str} without conditions (i) and (ii) is proved in \cite[Cor.~4.6.5]{kisin2012modp} for $p>2$, assuming that $Z_G$ is connected.  
 \end{rem}
 \begin{rem}\label{rem:from str} In the setting of Theorem \ref{thm:LR str}, as   $\Phi^{\mathbb Z} \times G(\A_f^p)$-sets we have $$\coprod_{\varepsilon \in \mathcal E} S_{\tau_{\varepsilon}}(G, \phi_{0,\varepsilon}) \cong \coprod_{\phi \in \AM(G)/\mathrm{conj}} S_{\dertau(\phi)}(\phi), $$ where $\dertau(\phi)$ is the image of $\tilde \tau(\phi)$ under $\tilde {\fkH}(\phi) \to \fkH(\phi) \to \cH(\phi)$. By conditions (i) and (ii) in the theorem and the commutative diagram (\ref{diag:lifting tori-rational}), $(\phi \mapsto \dertau(\phi))$ is an element of $\Gamma(\cH)_0$ and is tori-rational. Hence Theorem \ref{thm:LR str} (together with Theorem \ref{thm:lan-stroh}) immediately implies the Langlands--Rapoport--$\tau$ Conjecture (Conjecture \ref{hypo about LR}) in the case of Hodge type. 
 \end{rem}

 \cleardoublepage
 
 \chapter{The Langlands--Rapoport--\texorpdfstring{$\tau$}{tau} Conjecture in the case of abelian type} \label{sec:abelian type}
In this chapter we prove the Langlands--Rapoport--$\tau$ Conjecture (Conjecture \ref{hypo about LR}) for all unramified Shimura data of abelian type. The proof is based on Theorem \ref{thm:LR str}.  Throughout this chapter we fix a prime $p$. 

 \section{The proof in the presence of a nice dominating Hodge-type datum}	
 In this section, we prove the Langlands--Rapoport--$\tau$ Conjecture for an unramified Shimura datum of abelian type when there exists a sufficiently nice auxiliary Shimura datum of Hodge type. Recall that a Shimura datum $(G_2,X_2)$ is of \emph{abelian type}\index{Shimura datum!abelian type} if there exists a Shimura datum $(G,X)$ of Hodge type and a central isogeny $G_{\der} \to G_{2,\der}$ inducing an isomorphism of Shimura data $(G^{\ad}, X^{\ad}) \isom (G_2^{\ad}, X_2^{\ad})$.
 \begin{defn}\label{defn:dominance}
 	Let $(G,X,p,\cG)$ be an unramified Shimura	datum of Hodge type. We say that an unramified Shimura datum $(G_2, X_2, p , \cG_2)$ is\index{Shimura datum!domination} \emph{dominated by} $(G,X, p , \cG)$, if there is a central isogeny $G_{\der} \to G_{2,\der}$ inducing an isomorphism of Shimura data $(G^{\ad}, X^{\ad}) \isom (G_2^{\ad}, X_2^{\ad})$ and if the following conditions are satisfied:
 	\begin{enumerate} 
 		\item We have $E(G,X)_{\fkp} = E(G^{\ad}, X^{\ad})_{\fkp'}$, where $\fkp$ and $\fkp'$ are the $p$-adic places determined as in \S \ref{para:unramified SD}.  
 		\item  	The isomorphism	 $ G_{\Q_p}^{\ad} \isom G_{2,\Q_p}^{\ad}$ extends to an isomorphism $\cG^{\ad} \isom \cG_2^{\ad} $.
 	\end{enumerate} In particular, $(G_2,X_2)$ is of abelian type. 
 \end{defn}
 \subsection{} Fix  $(G,X,p,\cG)$ and $(G_2, X_2, p , \cG_2)$ as in the above definition. 	Write $K_p = \cG(\ZZ_p)$ and $K_{2,p}  = \cG_2 (\ZZ_p)$.  Let $\Shh_{K_p} = \Shh_{K_p}(G,X)$ and $\Shh_{K_{2,p}} = \Shh_{K_{2,p}}(G_2,X_2)$ be the canonical integral models recalled in Theorem \ref{thm:integral model} (also cf.~Theorem \ref{thm:lan-stroh}). Write $\Phi$ for the $p^{[E(G,X)_{\fkp}:\QQ_p]}$-Frobenius, and similarly define $\Phi_2$. 
 
 Let $\phi_{0,\varepsilon}$ and $\tau_{\varepsilon}$ be as in Theorem \ref{thm:LR str} applied to   $(G,X,p,\cG)$. Then we have  a $\Phi^{\ZZ} \times \mathscr A(G_{\ZZ_{(p)}})$-equivariant bijection 
 $$  \Shh_{K_p} (\Fpbar ) \isom  \coprod_{\varepsilon \in \mathcal E} S_{\tau_{\varepsilon}}(G, \phi_{0,\varepsilon})  $$ over  $\pi(G,X)$. As in the proof of \cite[Thm.~4.6.7]{kisin2012modp}, by taking the fibers over a fixed element of $\pi(G,X)$ and applying  $[\mathscr A(G_{2, \ZZ_{(p)}}) \times \cdot ] / (\mathscr A(G_{\ZZ_{(p)}}))^\circ$ to both sides, we obtain a $\Phi_2^{\ZZ} \times \mathscr A(G_{2, \ZZ_{(p)}})$-equivariant bijection 
 $$  \Shh_{K_{2,p}} (\Fpbar ) \isom \coprod_{\varepsilon \in \mathcal E} S_{\tau_{\varepsilon}}(G_2, \phi_{0,\varepsilon}) .   $$
 Now as $\Phi_2^{\ZZ} \times G_2(\A_f^p)$-sets we have (cf.~Remark \ref{rem:from str})
 $$\coprod_{\varepsilon \in \mathcal E} S_{\tau_{\varepsilon}}(G_2, \phi_{0,\varepsilon}) \cong \coprod_{\phi_2 \in \AM(G_2)/\mathrm{conj}} S_{\dertau(\phi_2)}(\phi_2), $$ where $\dertau(\phi_2) \in \cH(\phi_2)$ for each $\phi_2$ is defined as follows: Let $\varepsilon$ be such that $\phi_2^{\ad} \in \AM(G_2^{\ad}) \cong \AM(G^{\ad})$ is conjugate to $\phi_{0,\varepsilon}$. We set $\dertau(\phi_2)$ to be the image of $\tau_{\varepsilon}$ under the (canonical) composite map $I_{\phi_{0,\varepsilon}}^{\ad}(\A_f) \cong I_{\phi_2}^{\ad}(\A_f) \to \cH(\phi_2)$.

 Denote by $\widetilde {\fkH}_2, \cH_2$ the sheaves on $\AM(G_2)$ that are the counterparts of $\widetilde{\fkH}, \cH$ on $\AM(G)$ respectively. In order to prove the Langlands--Rapoport--$\tau$ conjecture for $\Shh_{K_{2,p}}$, it remains to show that the family $(\dertau(\phi_2))$ is an element of $\Gamma(\cH_2)_0$ and is tori-rational. 
 
 By Lemma \ref{lem:bij AM} (this result and Lemma \ref{lem:lifting phi_0} used in its proof are valid for arbitrary unramified Shimura data) applied to $(G,X,p,\cG)$ and $(G_2,X_2,p,\cG_2)$, the following spaces are identified:
 $$\AM(G)/{\approx}, \quad  \AM(G^{\ad})/{\approx}, \quad  \AM(G_2)/{\approx}, \quad  \AM(G_2^{\ad})/{\approx}. $$
 On the identified space, the sheaves $\widetilde{\fkH}_{\approx}$ and $\widetilde{\fkH}_{2,\approx}$ are identified, since $\widetilde \fkH(\phi) = I_{\phi,\sconn}(\A_f) \backslash I_{\phi}^{\ad} (\A_f) / I_{\phi}^{\ad} (\Q)$ depends only on   $I_{\phi}^{\ad}$ and not on $G$ and similarly for $\widetilde\fkH(\phi_2)$. Via this identification, the family 
 $$(\tilde \tau (\phi))_{\phi \in \AM(G)} \in \Gamma (\widetilde{\fkH})_0 $$ in Theorem \ref{thm:LR str} (i) determines a family  $$(\tilde \tau (\phi_2))_{\phi_2 \in \AM(G_2)} \in \Gamma (\widetilde{\fkH}_2)_0. $$ Clearly the element $\dertau(\phi_2)$ for each $\phi_2 \in \AM(G_2)$ is the image of $\tilde \tau(\phi_2)$ under $\widetilde{\fkH}(\phi_2) \to \cH(\phi_2)$. Hence we have 
 $$(\dertau (\phi_2))_{\phi_2 \in \AM(G_2)} \in \Gamma (\cH_2)_0. $$

 It remains only to show that the above family is tori-rational. 
 For the moment, we show this under an additional assumption on $(G,X,p,\cG)$.  
 Let $\mu^{\ab} \in X_*(G^{\ab})$ denote the composite cocharacter $$ \GG_m \overset{\mu}{\To} G_{\Qbar} \To G^{\ab}_{\Qbar} , $$ where $\mu$ is any element of $\dmu_X(\Qbar)$ (with notation in \S \ref{subsubsec:SD}). Clearly $\mu^{\ab}$ is independent of the choice of $\mu$. 
 
 \begin{prop}\label{prop:transporting tori-rationality}
 	Assume  that $X_*(G^{\ab})$ is generated by $\mu^{\ab}$ as a $\Gal(\Qbar/\Q)$-module. Then $(\dertau (\phi_2)) \in \Gamma (\cH_2)_0$ is tori-rational. 
 \end{prop}
 \begin{proof} 
 	Let $\phi_2 \in \AM(G_2) $ and let $T_2$ be a maximal torus in $I_{\phi_2}$. We need to show that the image of $\dertau(\phi_2) \in \cH(\phi_2)$ in $\coh^1(\A_f,T_2)/\Sha_{G_2}^{\infty} (\Q, T_2)$ under the map (\ref{eq:defn of tori rational}) is trivial. Let $\phi\in \AM(G)$ be such that $\phi^{\ad} \in \AM(G^{\ad}) \cong \AM(G_2^{\ad})$ is conjugate to $\phi_2^{\ad}$. We choose an identification of $\QQ$-groups $I_{\phi^{\ad}} \cong I_{\phi_2^{\ad}}$ (which is canonical up to inner automorphisms) and denote the identified group by $I_0$. 
 	
 	Let $T_0$ be the image of $T_2$ under $I_{\phi_2} \to I_0$, and let $T$ be the preimage of $T_0$ along $I_\phi \to I_0$.  Then $T_0$ (resp.~$T$) is a maximal torus in $I_0$ (resp.~$I_{\phi}$).
 	As usual, set
 	\begin{align*}
 		T^{\dagger} & : = \ker(T \to G^{\ab}) = T \cap I_{\phi}^{\dagger} \\ T_2^{\dagger} & : = \ker (T_2 \to G_2^{\ab}) = T_2 \cap I_{\phi_2}^{\dagger}. 
 	\end{align*}
 	Then there is a natural $\QQ$-map $u:T^{\dagger}\to T_2^{\dagger}$ induced by the central isogeny $G_{\der} \to G_{2,{\der}}$.

 	Define the torus $U: = (T \times _{T_0} T_2)^0$. Denote the natural maps $U \to T$ and $U \to T_2$ by $p_1$ and $p_2$ respectively. The inclusion map $T^{\dagger} \hookrightarrow T$ and the   map $u: T^{\dagger} \to T_2^{\dagger} \subset T_2$ together define a map $T^{\dagger} \to U$, which we denote by $\Delta$. Let $V$ be the quotient torus $U/\Delta(T^{\dagger})$. Then we have a commutative diagram with exact rows in the category of tori over $\QQ$: 
 	\begin{align}\label{diag:three rows}
 		\xymatrix { 1 \ar[r] & T^{\dagger} \ar[r] \ar@{=}[d] & T \ar[r] & G^{\ab} \ar[r] & 1 \\
 			1 \ar[r] & T^{\dagger} \ar[r]^{\Delta}  \ar[d]^{u} & U \ar[d]^{p_2} \ar[u]_{p_1} \ar[r] & V \ar[u]_{\bar p_1} \ar[d]^{\bar p_2} \ar[r] & 1 \\
 			1 \ar[r] & T_2^{\dagger} \ar[r] & T_2 \ar[r] & G_2^{\ab} \ar[r] & 1
 		}
 	\end{align}
 	Here $\bar p_1$ and $\bar p_2$ are induced by $p_1$ and $p_2$ respectively.
 	
 	We claim that the image of the map $X_*(V) \to X_*(G^{\ab})$ induced by $\bar p_1$ contains $\mu^{\ab}$. To show the claim, let $\mu_T$ be an arbitrary element of $\dmu_X(\Qbar)$ that factors through $T_{\Qbar}$ (which is embedded into $G_{\Qbar}$ via $I_{\phi, \Qbar} \hookrightarrow G_{\Qbar}$). Similarly, pick $\mu_{T_2}\in \dmu_{X_2} (\Qbar)$ that factors through $T_{2,\Qbar
 	}$. The two cocharacters of $T_{0 , \Qbar}$ induced by $\mu_T$ and $\mu_{T_2}$ respectively are conjugate by $G^{\ad}(\Qbar)$, and are hence in the same orbit under the Weyl group of $(G^{\ad}_{\Qbar}, T_{0,\Qbar})$ (cf.~\cite[Lem.~(1.1.3)]{kottwitztwisted}). Since the Weyl group of $(G^{\ad}_{\Qbar}, T_{0,\Qbar})$ and that of $(G_{2,\Qbar} , T_{2,\Qbar})$ are canonically isomorphic, we can replace $\mu_{T_2}$ by a Weyl-conjugate and assume that $\mu_T$ and $\mu_{T_2}$ induce the same cocharacter of $T_{0, \Qbar}$. We then obtain from $
 	\mu_T $ and $\mu_{T_2}$ an element of $X_*(U)$. The image of this element under the composite map $X_*(U) \to X_*(V) \to X_*(G^{\ab})$ is $\mu^{\ab}$, by the upper right commutative square in (\ref{diag:three rows}). The claim is proved.

 	By the claim and by our assumption on $X_*(G^{\ab})$, the map $X_*(V) \to X_*(G^{\ab})$ is surjective.
 	It follows that the kernel of $\bar p_1: V \to G^{\ab}$ is a torus, which we denote by $V^{\dagger}$. Clearly the map $(\bar p_1, \bar p_2): V \to G^{\ab} \times G_2^{\ab}$ is an isogeny between tori over $\Q$, and so $V_{\QQ_p}$ and $V^{\dagger}_{\QQ_p}$ are both unramified tori. 
 	Let $\cV$ denote the connected N\'eron model of $V_{\Q_p}$ over $\ZZ_p$. The kernel of $\cV \to \cG^{\ab}$ is the connected N\'eron model of $V^{\dagger}_{\Q_p}$. Hence by Lang's theorem,  the map $\cV (\Z_p ) \to \cG^{\ab} (\Z_p)$ is surjective.

 	Using this surjectivity result and the commutative diagram (\ref{diag:three rows}), we see that the natural map
 	$$ \coh^1(\A_f, T^{\dagger})/ \Sha^{\infty}_{G_{\der}} (\Q, T^{\dagger}) \To \coh^1(\A_f, T_2^{\dagger})/ \Sha^{\infty}_{G_{2,\der}} (\Q, T_2^{\dagger})$$
 	induced by $u: T^{\dagger} \to T_2^{\dagger}$ descends to a map
 	\begin{multline*}
 		u_*:    \fkH(\phi, T) = 		\cok \bigg( \cG^{\ab}(\Z_p) \rightarrow \coh^1(\A_f, T^{\dagger})/\Sha_{G_{\der}}^{\infty}(\Q, T^{\dagger}) \bigg) \To  \\ \fkH(\phi_2, T_2) =  \cok \bigg( \cG_2^{\ab}(\Z_p) \rightarrow \coh^1(\A_f, T_2^{\dagger})/\Sha_{G_{2,\der}}^{\infty}(\Q, T_2^{\dagger}) \bigg).
 	\end{multline*}
 	(Indeed, the first cokernel does not change if $\cG^{\ab}(\Z_p)$ is replaced by $\cV(\Z_p)$.) Moreover, there is a natural commutative diagram 
 	$$ \xymatrix{ \widetilde{\fkH}(\phi) \ar[d]^{\cong}\ar[r] & \fkH(\phi) \ar[r] & \fkH(\phi,T) \ar[d] ^{u_*} \\ \widetilde{\fkH}(\phi_2) \ar[rr] \ar[d] &  & \fkH(\phi_2,T_2) \ar[d] \\ \cH(\phi_2) \ar[rr]^-{(\ref{eq:defn of tori rational})} && \coh^1(\A_f,T_2)/\Sha_{G_2}^{\infty} (\Q, T_2)  .}$$
 	From this and condition (ii) in Theorem \ref{thm:LR str}, we see that the element  $\dertau(\phi_2)$, being the image of  $\tilde \tau (\phi)$ under $\widetilde{\fkH}(\phi) \isom  \widetilde{\fkH}(\phi_2) \to \cH(\phi_2)$,  has trivial image in $\coh^1(\A_f,T_2)/\Sha_{G_2}^{\infty} (\Q, T_2)$ as desired.  \end{proof}
 In conclusion, we have proved the following result. 
 \begin{thm}\label{thm:partial} Let $(G,X,p,\cG)$ be an unramified Shimura datum of Hodge type such that $X_*(G^{\ab})$ is generated by $\mu^{\ab}$ as a $\Gal(\Qbar/\Q)$-module. Then Conjecture \ref{hypo about LR} holds for any unramified Shimura datum $(G_2, X_2, p , \cG_2)$  dominated by $(G,X, p , \cG)$. \qed
 \end{thm}
 
 \section{Construction of a nice  dominating Hodge-type datum}\label{subsec:Proof of LRtau}
 In view of Theorem \ref{thm:partial}, in order to prove Conjecture \ref{hypo about LR} for \emph{all} unramified Shimura data of abelian type, it remains to prove the following result.

 \begin{prop}\label{prop:aux datum} 
 	Let $(G_2,X_2,p ,\cG_2)$ be an unramified Shimura datum of abelian type. Then there exists an unramified Shimura datum $(G,X, p ,\cG)$ of Hodge type dominating $(G_2,X_2,p ,\cG_2)$ such that $X_*(G^{\ab})$ is generated by $\mu^{\ab}$ as a $\Gal(\Qbar/\Q)$-module. 
 \end{prop}

 \subsection{}\label{para:Deligne's category} As preparation we recall some constructions from \cite[\S A]{deligne1982motifs}.  
 We denote by $\iota$ the complex conjugation in $\Gal(\Qbar/\Q)$, defined with respect to our fixed embedding $\Qbar \hookrightarrow \CC$. For any $\Q$-torus $T$ and $\mu \in X_*(T)$, we write $E_{\mu}$ for the field of definition of $\mu$ inside $\Qbar$. 
 
 Following \cite[\S A, (a)]{deligne1982motifs}, we consider the category\index[Oct2026-sorted-notation]{i 0 03 00 i deligne prime s category@$\mathcal I$ (Deligne's category)} $\cI$ as follows. The objects are the pairs $(T,\mu)$, where $T$ is a $\Q$-torus and $\mu \in X_*(T)$ such that $w : = \mu + \iota(\mu)$ is defined over $\Q$ and $(T/w(\GG_m))_{\RR}$ is anisotropic. Morphisms $(T,\mu) \to (T',\mu')$ are $\Q$-homomorphisms $T \to T'$ taking $\mu$ to $\mu'$. For each $(T,\mu)$ in $\cI$, we know that $T$ is a cuspidal torus since it satisfies condition (ii) in Lemma \ref{lem:cusp TFAE}. By condition (vii) in that lemma, $T$ splits over a CM or totally real field, so $E_{\mu}$ is either CM or totally real. In the latter case we have $E_{\mu} = \Q$, since $\mu + \iota(\mu)$ is defined over $\Q$ and is equal to $2\mu$. 
 
 As in \textit{loc.~cit.}, for every $(T,\mu)$ in $\cI$, there exists an object in $\cI$ of the form $(S^L, \mu^L)$ which maps to $(T, \mu)$. Here $L$ is a CM field with maximal totally real subfield $L_0$, and $S^L$ is the $\Q$-torus $$ \Res_{L/\QQ} \GG_m / \ker (\N_{L_0/\Q} : \Res_{L_0/\QQ}\GG_m \to \GG_m ). $$ The cocharacter $ \mu^L \in X_*(S^L)$ is  induced by the cocharacter of $\Res_{L/\QQ} \GG_m$ corresponding to the canonical embedding $L \hookrightarrow \ol \QQ$.  In fact, as discussed above, $E_\mu$ is either a CM field or $\Q$. We can take $L$ to be any CM field containing $E_{\mu}$, and take the homomorphism $S^L \to T$ to be the one induced by the composite homomorphism 
 $$ \Res_{L/\QQ} \GG_m \xrightarrow{\Res_{L/\QQ} \mu }  \Res_{L/\QQ} T   \xrightarrow{\N_{L/\QQ}}  T   .$$ 
 Using the defining conditions of $\cI$ satisfied by $(T,\mu)$, one easily checks that the above homomorphism factors through $S^L$.

 \begin{proof}[Proof of Proposition \ref{prop:aux datum}]  Let $(G_1, X_1, p, \cG_1)$  be an unramified Shimura datum of Hodge type dominating  $(G_2, X_2, p , \cG_2)$; this exists by \cite[Lem.~4.6.6]{kisin2012modp} and the proof of \cite[Cor.~3.4.14]{kisin2010integral}.\footnote{In the proof of  \cite[Lem.~4.6.6]{kisin2012modp} the following fact is used: Every totally real field $F$ admits a totally imaginary quadratic extension $K/F$ such that every prime of $F$ over $p$ splits in $K$. This fact is an immediate consequence of Thm.~5 or Thm.~6 in \cite[\S X.2]{AT}, regardless of the parity of $p$. In the proof of \cite[Cor.~3.4.14]{kisin2010integral}, the assumption that $p$ is odd is not used for the construction of $\cG$.}  Let $\mu_1^{\ab} \in X_*(G_1^{\ab})$ be induced by  $\mu_1 \in \dmu_{X_1}(\Qbar)$.
 	Since $(G_1,X_1)$ is of Hodge type,  its weight cocharacter $w_1$ is defined over $\QQ$, and $(Z_{G_1}^0/ w_1(\GG_m))_{\RR}$ is anisotropic. Since $G_1^{\ab}$ is isogenous to $Z_{G_1}^0$ and since $\mu_1^{\ab} + \iota(\mu_1^{\ab}) \in X_*(G_1^{\ab})$ is induced by $w_1 \in X_*(Z_{G_1}^0)$, the pair $(G_1^{\ab}, \mu_1^{\ab})$ is in the category $\cI$. 
 	
 	We define a field $L$ as follows. 
 	As discussed in \S \ref{para:Deligne's category}, the field $E_{\mu_1^{\ab}}$ is either CM or $\QQ$. In the former case, we let $L$ be $E_{\mu_1^{\ab}}$. In the latter case, we let $L$ be an arbitrary imaginary quadratic field in which $p$ splits. We construct the morphism $(S^L ,\mu^L) \to (G_1^{\ab}, \mu_1^{\ab})$ in $\cI$ as in \S \ref{para:Deligne's category}. Since we have an unramified Shimura datum $(G_1, X_1, p ,\cG_1)$, the torus $G_1^{\ab}$ is unramified over $\Q_p$. Hence every conjugate of the subgroup $\Gamma_{p,0} \subset \Gal(\Qbar/\Q)$ fixes $\mu_1^{\ab}$. It follows that the Galois closure of $E_{\mu_1^{\ab}}/\QQ$ is unramified over $p$, and the same holds for the Galois closure of $L/\QQ$. Hence the torus $S^L$ is unramified over $\Q_p$ since it splits over the Galois closure of $L/\QQ$.

 	Recall that   specifying a Shimura datum for a given torus over $\Q$ amounts to specifying the Hodge cocharacter, which can be an arbitrary cocharacter over $\Qbar$.
 	We thus obtain from $(S^L,  \mu^L)$ a Shimura datum $(S^L, \set{h^L})$. By the paragraph preceding \cite[\S A, (a)]{deligne1982motifs} and by \cite[Lem.~A.2]{deligne1982motifs}, we know that $S^L$ admits a faithful representation over $\Q$ such that the Hodge structure on that representation determined by $h^L$ is of type $\set{(-1,0), (0,-1)}$. By \cite[Prop.~2.3.2]{deligne1979varietes}, the last fact implies that $(S^L, \set{h^L})$ is a Shimura datum of Hodge type.
 	
 	Let $G = G_1 \times_{G_1^{\ab}}  S^L $. Fix an element $h_1 \in X_1$. Note that $h_1$ and $h^L$ induce the same map $\mathbb S \to G_{1,\RR}^{\ab}$. Hence we obtain a homomorphism $$h= (h_1, h^L) : \mathbb S \To G_{\RR}. $$  Let $X$ be the $G(\RR)$-conjugacy class of $h$, which is clearly a Shimura datum for $G$.\footnote{Here $X$ may depend on the choice of $h_1$, but this does not matter for our proof.} Since $(G_1,X_1)$ and $(S^L, \set{h^L})$ are both of Hodge type, so is $(G,X)$. Indeed, choose faithful symplectic similitude representations of $(G_1,X_1)$ and
 	$(S^L,\{h^L\})$, and pull them back to $G$.
 	The ratio of their similitude characters is a character of $G$
 	which factors through $G^{\ab}=S^L$.
 	This ratio is trivial on $\mu^L$, since the image of $\mu^L$ under $X_*(S^L) \to X_*(G_1^{\ab})$ is $\mu_1^{\ab}$. 
 	It is therefore trivial on every Galois conjugate of $\mu^L$.
 	These conjugates generate $X_*(S^L)$, so the two similitude
 	characters agree on $G$.
 	Hence	by taking the direct sum of the two representations we obtain a faithful symplectic similitude
 	representation of $(G,X)$.

 	Since there is a central isogeny $G_{1,\der} \to G_{2,\der}$ inducing an isomorphism $(G_1^{\ad}, X_1^{\ad}) \isom (G_2^{\ad}, X_2^{\ad})$, there is also a central isogeny $G_{\der} \to G_{2,\der}$ inducing an isomorphism $(G^{\ad}, X^{\ad}) \isom (G_2^{\ad}, X_2^{\ad})$. For any number field $F$ in $\Qbar$ we write $F_p$ for the completion of $F$ with respect to the fixed embedding $\Qbar \hookrightarrow \Qpbar$. Note that $E_{\mu^L}$ is contained in $L$, and $L_p$ is equal to $(E_{\mu_1^{\ab}})_p$, which is a subfield of $E(G_1,X_1)_p$. It follows that $E(G,X)_p = E(G_1, X_1)_p$. Since $E(G_1, X_1)_p = E(G_1^{\ad}, X_1^{\ad})_p$ by assumption, we have $E(G,X)_p = E(G^{\ad}, X^{\ad})_p$. Since $G_1$ and $S^L$ are both unramified over $\Q_p$, so is $G$. Thus as in the proof of \cite[Cor.~3.4.14]{kisin2010integral} we can extend $(G,X)$ to an unramified Shimura datum $(G,X,p, \cG)$ such that the  isomorphism	 $ G_{\Q_p}^{\ad} \isom G_{2,\Q_p}^{\ad}$ extends to an isomorphism $\cG^{\ad} \isom \cG_2^{\ad} $. In conclusion, the unramified Shimura datum $(G,X,p,\cG)$ of Hodge type dominates $(G_2,X_2,p,\cG_2)$.  
 	
 	It remains to show that $X_*(G^{\ab})$ is generated by $\mu^{\ab}$ as a $\Gal(\Qbar/\Q)$-module. This is clear since $G^{\ab}= S^L$ and $\mu^{\ab}=\mu^L$. 
 \end{proof}

 \begin{thm}[cf.~Theorem \ref{introthm:LRtau} in the Introduction]\label{thm:LRtau for abelian type} Conjecture \ref{hypo about LR} holds for all unramified Shimura data of abelian type. 
 \end{thm}
 \begin{proof}
 	This follows from Theorem \ref{thm:partial} and Proposition \ref{prop:aux datum}. 
 \end{proof}

 \begin{thm}[cf.~Theorem \ref{thm:intro point counting} in the Introduction] \label{thm:announcement 1}
 	Conjecture \ref{conj:point counting formula} holds for all Shimura varieties of abelian type.
 \end{thm}
 \begin{proof}
 	This follows from Theorems \ref{thm:main thm in point counting} and \ref{thm:LRtau for abelian type}. 
 \end{proof} 
 
 \cleardoublepage
 
 \part{~Stabilization}
 \cleardoublepage
 
 \chapter{Preliminaries for stabilization}\label{sec:prelim-for-stab}
 	This chapter recalls the basics of the trace formula and endoscopy that are needed to stabilize the point counting formula for Shimura varieties. This includes local base-change endoscopy, which is a special case of twisted endoscopy pertaining to twisted orbital integrals at $p$. Along the way, we discuss the trace formula with fixed central character and miscellaneous results in endoscopy since references are often hard to locate in the desired generality.

 \section{Central character data and the trace formula}\label{sec:TF-fixed-central}
 
 To stabilize the point counting formula for Shimura varieties (\ref{pcf}) in general, it is necessary to work with fixed central characters. 
 To this end, we are going to introduce the formalism of central character data following \cite[\S 3.1]{Arthur}.
 The point counting formula can be understood through the lens of a particular central character datum, but it is useful to allow flexible central character data to accommodate $z$-extensions during the stabilization process.
 
 We devote \S\ref{sec:TF-fixed-central} to discussing the invariant trace formula with fixed central character. Though this is not logically needed for the stabilization in \S\ref{sec:stabilization}, it puts central character data into context, motivates the definition of stable distributions with fixed central characters, and also has an application in \S\ref{sec:spectral interpretation} below.
 
 \subsection{}
 
 Throughout \S\ref{sec:TF-fixed-central}, let $G$ be a reductive  group over $\Q$ with center $Z$. Write $A_Z$ for the maximal $\Q$-split torus in $Z$ and set\index[Oct2026-sorted-notation]{a 0 00 00 a subscript z infty@$A_{Z,\infty}$} $A_{Z,\infty}:=A_Z(\R)^0$.

 \begin{defn}
 	A\index{central character datum} \emph{central character datum} for $G$ is a pair $(\fkX,\chi)$, where $\fkX$ is a closed subgroup of $Z(\A)$ containing $A_{Z,\infty}$ such that $Z(\Q)\fkX$ is closed in $Z(\A)$, with a choice of Haar measure on $\fkX$ (often implicit),
 	and  $\chi:\fkX\ra\CC^\times$ a continuous character which is trivial on\index[Oct2026-sorted-notation]{x 0 04 00 x subscript q@$\mathfrak X_{\mathbb Q}$} $\fkX_{\Q}:=\fkX\cap Z(\Q)$. 
 \end{defn}
 
 \begin{rem}
 	The two extreme cases where $\fkX=A_{Z,\infty}$ or $\fkX=Z(\A)$ are already interesting, but it is important to allow intermediate groups for our purposes.
 \end{rem}
 
 \subsection{}\label{central-character-data}
 
 Let $\Gamma_{\el}(G)$\index[Oct2026-sorted-notation]{zz03 gamma 0 00 00 gamma subscript ell g@$\Gamma_{\el}(G)$} denote the set of elliptic conjugacy classes in $G(\Q)$. Given a central character datum $(\fkX,\chi)$, denote by $\Gamma_{\el,\fkX}(G)$\index[Oct2026-sorted-notation]{zz03 gamma 0 00 00 gamma subscript ell x g@$\Gamma_{\el,\fkX}(G)$} the set of $\fkX_\Q$-orbits in $\Gamma_{\el}(G)$ with respect to the multiplication action.
 Write\index[Oct2026-sorted-notation]{stab 0 01 00 stab subscript x gamma@$\mathrm{Stab}_{\fkX}(\gamma)$} $\mathrm{Stab}_{\fkX}(\gamma)$ for the stabilizer subgroup of $\fkX_\Q$ fixing $\gamma\in \Gamma_{\el}(G)$. It is not hard to see that $\mathrm{Stab}_{\fkX}(\gamma)$ is finite, for instance by reducing to the case of a product of general linear groups via a (possibly reducible) faithful representation of $G$.
 
 Fix a maximal compact subgroup\index[Oct2026-sorted-notation]{k 0 00 00 k subscript infty@$K_\infty$} $K_\infty\subset G(\R)$.
 Let $v$ be a place of $\Q$, and\index[Oct2026-sorted-notation]{z 0 04 00 z@$\mathfrak Z$} $\fkZ$ a closed subgroup of $Z(\Q_v)$.
 For a continuous character $\omega:\fkZ\ra \CC^\times$ define\index[Oct2026-sorted-notation]{h 0 03 00 h g q subscript v omega superscript minus 1@$\mathcal H(G(\Q_v),\omega^{-1})$} $\cH(G(\Q_v),\omega^{-1})$ to be the Hecke algebra of smooth functions on $G(\Q_v)$ which transform under $\fkZ$ by $\omega^{-1}$ and have compact support modulo $\fkZ$; we also require $K_\infty$-finiteness if $v=\infty$. Let $\pi$ be an admissible representation of $G(\Q_v)$ which has central character on $\fkZ$ equal to $\omega$. For $f\in \cH(G(\Q_v),\omega^{-1})$ define
 \begin{eqnarray}
 	\pi(f)(u) &:=& \int_{G(\Q_v)/\fkZ} f(g)\pi(g)u \cdot dg,\quad u\in \pi.\nonumber
 \end{eqnarray}
 The trace of the trace-class operator $\pi(f)$ is denoted by $\tr(f\mid \pi)$ or $\tr \pi(f)$. The orbital integrals for $f\in \cH(G(\Q_v),\omega^{-1})$ are defined by the same formula as for $\cH(G(\Q_v))$, cf.~\S\ref{subsubsec:setting for point count}.
 These definitions obviously extend to the adelic setting.

 \subsection{}\label{sub:inv-dist}
 
 We recall the invariant distributions\index[Oct2026-sorted-notation]{i 0 00 00 i subscript star chi subscript 0@$I_{\star,\chi_0}$} \index[Oct2026-sorted-notation]{t 0 00 00 t subscript star chi subscript 0@$T_{\star,\chi_0}$}
 $$I_{\geom,\chi_0},~I_{\spec,\chi_0},~I_{\el,\chi_0},~I_{\disc,\chi_0},~T_{\el,\chi_0},~T_{\disc,\chi_0}$$
 in the classical setup where the central character datum consists of $\fkX=A_{Z,\infty}$ and $\chi_0:A_{Z,\infty}\ra\CC^\times$.
 First off, $I_{\geom,\chi_0}$ and $I_{\spec,\chi_0}$ are Arthur's invariant distributions given in sections 3 and 4 of \cite{Art88b}, respectively. Define $I_{\el,\chi_0}$ to be the $M=G$ part of formula (3.3) and $I_{\disc,\chi_0}$ to be formula (4.4), both referenced to \textit{loc.~cit}. All four distributions are distributions on $\cH(G(\A),\chi_0^{-1})$.\footnote{Arthur defines them as distributions on a certain space of functions on $G(\A)^1$ named $\cH(G(\A)^1)$, but this space is isomorphic to $\cH(G(\A),\chi_0^{-1})$ via the product decomposition $G(\A)= G(\A)^1 \times A_{Z,\infty}$. We do not mention $\cH(G(\A)^1)$ again.}
 
 For $\gamma\in \Gamma_{\el}(G)$, write $I_\gamma$ for the connected centralizer of $\gamma$ in $G$.
 We put
 \begin{eqnarray*}
 	T_{\el,\chi_0}(f)&:=&\sum_{\gamma\in \Gamma_{\el}(G)} \iota(\gamma)^{-1}\vol(I_\gamma(\Q)\bs I_\gamma(\A)/A_{Z,\infty}) O_{\gamma}(f), \\
 	T_{\disc,\chi_0}(f)&:=&\tr(f \mid  L^2_{\disc,\chi_0}(G(\Q)\bs G(\A))),\quad\quad f\in \cH(G(\A),\chi_0^{-1}). 
 \end{eqnarray*}
 In general $I_{\star}(f)$ is more complicated than $T_{\star}(f)$ for $\star\in \{\el,\disc\}$. Arthur's invariant trace formula is the equality 
 \begin{equation}\label{e:ArthurITF}
 	I_{\geom,\chi_0}=I_{\spec,\chi_0}.
 \end{equation}
 When $G/Z$ is anisotropic over $\Q$,
 \begin{equation}\label{e:ArthurITFanis}
 	T_{\el,\chi_0}=I_{\el,\chi_0}=I_{\geom,\chi_0}=I_{\spec,\chi_0}=I_{\disc,\chi_0}=T_{\disc,\chi_0}.
 \end{equation}

 \subsection{}\label{sub:trace-formula-central-character}
 
 Next we introduce the trace formula with respect to a fixed character on a closed central subgroup. This must be well known to experts, but for lack of a convenient reference, we state the formula here.\footnote{Section 3.1 of \cite{Arthur} discusses such a variant in the discrete part of the trace formula. Sections 2 and 3 of \cite{ArtSTF1} present both the spectral and geometric expansions of the trace formula with fixed central character on an induced torus. We treat a more general case than \textit{loc.~cit.} but proceed in a similar spirit.}

 Let $(\fkX,\chi)$ be a central character datum for $G$. Suppose that $\chi_0:A_{Z,\infty}\ra\CC^\times$ is the restriction of $\chi$. We will obtain the $\chi$-versions of the above invariant distributions by averaging.

 \begin{lem}\label{lem:torus-cocompact}
 	Let $D$ be a multiplicative group over $\Q$, $A_D$ its maximal $\Q$-split subtorus, and $A_{D,\infty}:=A_D(\R)^0$. Then $D(\Q)\bs D(\A)/A_{D,\infty}$ is compact.
 \end{lem}
 
 \begin{proof}
 	Replacing $D$ by $D^0$, we may assume that $D$ is a torus. Via a closed embedding, we reduce to the case where $D$ is a finite product of tori of the form $\Res_{F/\Q}\GG_m$ for a finite extension $F$ over $\Q$. When $D= \Res_{F/\Q}\GG_m$, the lemma is clear since $F^\times \bs \A_F^1$ is compact, where $\A_F^1$ denotes the group of id\`eles of norm $1$.
 \end{proof}
 
 \begin{cor}
 	The quotient $\fkX_\Q\bs \fkX/A_{Z,\infty}$ is compact.
 \end{cor}
 
 \begin{proof}
 	This is clear since the inclusion $\fkX\subset Z(\A)$ induces a closed embedding
 	from $\fkX_\Q\bs \fkX/A_{Z,\infty}$ into $Z(\Q)\bs Z(\A)/A_{Z,\infty}$.	(The image is closed since $ Z(\Q)\fkX$ is closed in $Z(\A)$.)
 \end{proof}
 \subsection{}
 There is a surjection $\cH(G(\A),\chi_0^{-1})\ra \cH(G(\A),\chi^{-1})$ sending $f$ to the function
 $$g\mapsto \ol{f}_\chi(g):=\int_{\fkX/A_{Z,\infty}} f(gz)\chi(z) dz,$$
 where the integral converges because $z\mapsto f(gz)$ has compact support in $\fkX/A_{Z,\infty}$.
 Given a function $f$ on $G(\A)$ and $z\in Z(\A)$, write $f_z$ for the translated function $g\mapsto f(gz)$. For each $\star\in \{\geom,\spec,\el,\disc\}$ define\index[Oct2026-sorted-notation]{i 0 00 00 i subscript star chi@$I_{\star,\chi}$} \index[Oct2026-sorted-notation]{t 0 00 00 t subscript star chi@$T_{\star,\chi}$}
 \begin{equation}\label{eq:def-of-T_chi}
 	I_{\star,\chi}(f):=\frac{1}{\vol(\fkX_\Q\bs \fkX/A_{Z,\infty})}\int_{\fkX_\Q\bs \fkX/A_{Z,\infty}} \chi(z) I_{\star,\chi_0}(f_z)dz,\quad f\in \cH(G(\A),\chi_0^{-1}),
 \end{equation}
 as well as $T_{\el,\chi}$ and $T_{\disc,\chi}$ in a similar manner.
 In the special case where $G/Z$ is anisotropic over $\Q$, it is clear from \eqref{e:ArthurITFanis} that
 $$T_{\el,\chi}(f)=I_{\el,\chi}(f)=I_{\geom,\chi}(f)=I_{\spec,\chi}(f)=I_{\disc,\chi}(f)=T_{\disc,\chi}(f).$$
 These averaged distributions factor through
 $f\mapsto\bar f_\chi$. We use the same notation for the
 resulting distributions on $\cH(G(\A),\chi^{-1})$.
 
 Write 
 $L^2_{\disc,\chi}(G(\Q)\bs G(\A))$
 for the discrete spectrum in the $L^2$-space of complex-valued functions $\phi$ on $G(\Q)\bs G(\A)$ such that
 \begin{itemize}
 	\item  $\phi(gz)=\chi(z)\phi(g)$ for every $g\in G(\A)$ and every $z\in \fkX$,
 	\item $\int_{G(\Q)\bs G(\A)^1/\fkX\cap G(\A)^1}|\phi(g)|^2 dg <\infty$ (for any Haar measure), where $G(\A)^1$ denotes the ``norm one'' subgroup of $G(\A)$ as
 	defined in \cite[p.~917]{Art78}.
 \end{itemize}

 \begin{prop}\label{prop:Tell,Tdisc} For $f\in \cH(G(\A),\chi_0^{-1})$ the following equalities hold.
 	\begin{eqnarray*}
 		T_{\el,\chi}(f)&=&\sum_{\gamma\in \Gamma_{\el,\fkX}(G)} |\mathrm{Stab}_{\fkX}(\gamma)|^{-1}\iota(\gamma)^{-1} \vol(I_\gamma(\Q)\bs I_\gamma(\A)/\fkX) O_{\gamma}(\ol{f}_\chi), 
 		\\		T_{\disc,\chi}(f)&=&\tr (\ol{f}_\chi \,|\, L^2_{\disc,\chi}(G(\Q)\bs G(\A))). 
 	\end{eqnarray*}
 \end{prop}
 
 \begin{proof}
 	We compute $T_{\el,\chi}(f)$ as follows:
 	
 	\begin{eqnarray*}
 		& &\int_{\fkX_\Q\bs \fkX/A_{Z,\infty}}\sum_{\gamma\in \Gamma_{\el}(G)}  \chi(z)\iota(\gamma)^{-1} \frac{\vol(I_\gamma(\Q)\bs I_\gamma(\A)/A_{Z,\infty})}{\vol} O_{\gamma z}(f)dz\\
 		&=&\sum_{\gamma\in \Gamma_{\el}(G)} \int_{\fkX_\Q\bs \fkX/A_{Z,\infty}} \chi(z)\iota(\gamma)^{-1} \vol(I_\gamma(\Q)\bs I_\gamma(\A)/\fkX) O_{\gamma z}(f)dz\\
 		&=& \sum_{\gamma\in \Gamma_{\el,\fkX}(G)} \frac{\vol(I_\gamma(\Q)\bs I_\gamma(\A)/\fkX)}{|\mathrm{Stab}_{\fkX}(\gamma)| \iota(\gamma)} \int_{\fkX/A_{Z,\infty}} \chi(z) O_{\gamma z}(f)dz\\
 		&=&\sum_{\gamma\in \Gamma_{\el,\fkX}(G)}  \frac{\vol(I_\gamma(\Q)\bs I_\gamma(\A)/\fkX)}{|\mathrm{Stab}_{\fkX}(\gamma)| \iota(\gamma)}  O_{\gamma}(\ol{f}_\chi).
 	\end{eqnarray*}
 	
 	The equality for $T_{\disc,\chi}(f)$ follows from the following two observations. First, given an admissible representation $\pi$ of $G(\A)$ with central character $\chi_\pi$ on $\fkX$,
 	$$\int_{\fkX_\Q\bs \fkX/A_{Z,\infty}} \chi(z) \pi(f_z)dz= \pi(f) \int_{\fkX_\Q\bs \fkX/A_{Z,\infty}} \chi(z)\chi_{\pi}^{-1}(z)dz,$$
 	which equals 0 if $\chi\neq \chi_\pi$ and $\vol(\fkX_\Q\bs \fkX/A_{Z,\infty})\pi(f)$ if $\chi=\chi_\pi$. Second, if $\chi=\chi_\pi$,
 	\begin{eqnarray}
 		\pi(\ol{f}_\chi)&=&\int_{G(\A)/\fkX}\pi(g)\int_{\fkX/A_{Z,\infty}} f(gz)\chi(z)dzdg \nonumber\\
 		&=&\int_{G(\A)/\fkX}\int_{\fkX/A_{Z,\infty}} \pi(gz)f(gz)dz dg = \pi(f).\nonumber
 	\end{eqnarray}

 \end{proof}
 
 \subsection{}
 
 We record a simplification of the trace formula with fixed central character when the test function is a stable cuspidal function at the real place.
 
 Given a central character datum $(\fkX,\chi)$, suppose that $\fkX=\fkX^\infty\times \fkX_\infty$ with $\fkX^\infty\subset Z(\A_f)$ and $\fkX_\infty\subset Z(\R)$. (In particular $\fkX_\infty$ contains $A_{Z,\infty}$.) Accordingly we decompose $\chi=\chi^\infty\chi_\infty$. We also assume that $G(\R)$ contains a maximal torus which is compact modulo $\fkX_\infty$. 
 
 Let $\xi$ be an irreducible algebraic representation of $G_\CC$. The inverse of the central character of $\xi$ on $Z(\R)$ is denoted by $\chi_\xi$. Let\index[Oct2026-sorted-notation]{zz16 pi 0 00 00 pi subscript 2 xi@$\Pi_2(\xi)$} $\Pi_2(\xi)$ denote the set of isomorphism classes of (irreducible) discrete series representations of $G(\R)$ whose central and infinitesimal characters are the same as those of the contragredient of $\xi$. Define\index[Oct2026-sorted-notation]{f 1 00 00 f subscript xi@$f_\xi$} $f_\xi\in \cH(G(\R),\chi^{-1}_\xi)$ to be the sum of pseudo-coefficients of $\pi_\infty$ as $\pi_\infty$ runs over $\Pi_2(\xi)$, cf.~\cite[Lem.~3.1]{Art89}, \cite{CD90}.
 
 A\index{stable cuspidal function} \emph{stable cuspidal function} on $G(\R)$ (relative to $\fkX_\infty$) is defined to be $f_\infty\in \cH(G(\R),\chi_\infty^{-1})$ such that for every irreducible tempered representation $\pi_\infty$ of $G(\R)$ whose central character restricts to $\chi_\infty$ on $\fkX_\infty$,
 we have (i) $\tr \pi_\infty(f_\infty)=0$ unless $\pi_\infty$ is a discrete series representation, and (ii) $\tr \pi_\infty(f_\infty)$ has a constant value as $\pi_\infty$ runs over each discrete series $L$-packet, cf. \cite[\S4, p.~270]{Art89}. 
 An example is $f_\xi$ in the last paragraph. In general, when $A_{G_{\R}}(\R)^0\subset \fkX_\infty$, a stable cuspidal function is a finite linear combination of character twists of functions of the form $f_\xi$ (for  different $\xi$'s), up to a function whose orbital integrals are identically zero.

 \begin{prop}\label{prop:simple-TF}
 	Let $f=f^\infty f_\infty\in \cH(G(\A),\chi^{-1})$ with $f_\infty$ stable cuspidal. Then
 	$$I_{\spec,\chi}(f)=I_{\disc,\chi}(f)=T_{\disc,\chi}(f).$$
 \end{prop}
 
 \begin{proof}
 	This is proved in \cite[\S3]{Art89} when $\fkX_\infty=A_{Z,\infty}$. The same proof extends.
 \end{proof}
 
 \begin{rem}
 	When $f_\infty$ is stable cuspidal,
 	a simple expansion for $I_{\geom,\chi_0}(f)$ is obtained in  \cite[Thm.~6.1]{Art89}.
 	A similar expansion for $I_{\geom,\chi}(f)$ is given in \cite[\S 6.4, \S 6.5]{Dalal} for more general central character data.
 \end{rem}

 \section{Endoscopic data and \texorpdfstring{$z$}{z}-extensions}\label{sub:endoscopic-data}
 
 Langlands--Shelstad \cite[\S1.2]{LS87} and Kottwitz--Shelstad \cite[\S2.1]{KS99} have defined endoscopic data and related notions in the untwisted and twisted settings. Here we recall the untwisted case as well as a specific kind of local twisted endoscopy (generalizing the unramified base change) as studied in \cite[App.~ A]{Mor10}. We also prove some lemmas about $z$-extensions to facilitate reduction to the case of groups whose derived subgroups are simply connected.
 
 \subsection{} Throughout the rest of \S \ref{sec:prelim-for-stab}, let $G$ be a reductive  group over a local or global field $F$ of characteristic zero.
 
 \begin{defn}\label{def:end-datum}
 	Let $F,G$ be as above. An\index{endoscopic datum} \emph{endoscopic datum} for $G$ is a quadruple\index[Oct2026-sorted-notation]{e 1 04 00 e h h s eta@$\mathfrak e=(H,\cH,s,\eta)$} $\fke=(H,\cH,s,\eta)$, where
 	\begin{itemize}
 		\item $H$ is a quasi-split reductive group over $F$,
 		\item $\cH$ is a split extension of $W_F$ by $\hat H$ such that the $L$-action of $W_F$ on $\hat H$ determined by $\cH$ coincides with the $L$-action of the $L$-group $^L H$,
 		\item $s$ is a semi-simple element of $\hat G$,
 		\item $\eta: \cH \ra {}^L G$ is an $L$-morphism inducing an isomorphism $\hat H\cong \Cent(s,\hat G)^0$ (via $\eta$ we view $s$ also as an element of $\hat H$),
 	\end{itemize}
 	such that $\Int(s)\circ \eta=a\cdot \eta$ for a 1-cocycle $a:W_F\ra Z(\hat G)$ which is trivial (resp.~locally trivial) if $F$ is local (resp.~global). In this case $H$ is said to be an\index{endoscopic group} \emph{endoscopic group} for $G$.
 	
 	The datum $\fke$ is said to be\index{endoscopic datum!elliptic} \emph{elliptic} if $\eta(Z(\hat H)^\Gamma)^0\subset Z(\hat G)$. When $F$ is non-archimedean, we say that $\fke$ is\index{unramified!endoscopic datum} \emph{unramified} if $H$ and $G$ are unramified groups over $F$ and if $\eta$ is inflated from an $L$-morphism with respect to the Weil group of $F^{\ur}$ over $F$. An \emph{isomorphism} from $\fke=(H,\cH,s,\eta)$ to $\fke'=(H',\cH',s',\eta')$ is an element $g\in \hat G$ such that $g \eta(\cH) g^{-1}=\eta'(\cH')$ in $^L G$ and $g s g^{-1}=s'$ in $\hat G/Z(\hat G)$.
 \end{defn}
 
 \subsection{}\label{Z(G)-Z(H)}

 Automorphisms of an endoscopic datum $\fke$ induce outer automorphisms of $H$ as in \cite[(2.1.8)]{KS99}. By\index[Oct2026-sorted-notation]{out 0 01 00 out subscript f e@$\mathrm{Out}_F(\fke)$} $\Out_F(\fke)$ we mean the image subgroup of the outer automorphism group $\Out_F(H):=\Aut_F(H)/H_{\textrm{ad}}(F)$. Set\index[Oct2026-sorted-notation]{zz11 lambda 1 00 00 lambda e@$\lambda(\mathfrak e)$}
 $$\lambda(\fke):=|\Out_F(\fke)|\in \Z_{>0}.$$
 The set of endoscopic data is denoted by\index[Oct2026-sorted-notation]{e 0 00 00 e g@$E(G)$}\index[Oct2026-sorted-notation]{e 0 00 00 e subscript ell g@$E_{\el}(G)$} $E(G)$. Write $E_{\el}(G)$ for the subset of elliptic endoscopic data. Write\index[Oct2026-sorted-notation]{e 0 03 00 e g@$\mathcal E(G)$}\index[Oct2026-sorted-notation]{e 0 03 00 e subscript ell g@$\mathcal E_{\el}(G)$} $\cE(G)$ and $\cE_{\el}(G)$ for the corresponding sets of isomorphism classes. 
 
 For $\fke=(H,\cH,s,\eta)\in E(G)$, there exists a canonical injection over $F$
 \begin{equation}\label{eq:ZGZH}
 	Z_G \hra Z_H
 \end{equation}
 as we now explain. Since an inner twisting induces a canonical isomorphism of centers, we may assume that $G$ is quasi-split over $F$. Choose a maximal torus $T_H$ of $H$ over $F$. Then there exist a maximal torus $T$ of $G$ and an isomorphism $T_H\cong T$, both defined over $F$ \cite[p.~226]{LS87} such that the composite embedding $T\cong T_H \subset H$ is canonical up to $H(\ol F)$-conjugacy. Restricting from $T_H\cong T$, we get the desired map $Z_G\hra Z_H$, which is independent of the choices involved.

 \subsection{}\label{para:z-ext}
 If $G_{\der}=G_{\textup{sc}}$ then by \cite[Prop.~1]{Lan79}, every $\fke\in \cE(G)$ is represented by $(H,{}^L H,s,\eta)$, that is, we can take $\cH={}^L H$.  In general there is no guarantee that this is possible, so we use $z$-extensions to reduce to this case.
 
 Recall that a\index{z-extension@$z$-extension} $z$-\emph{extension} of $G$ over $F$ is defined as a reductive  group $G_1$ equipped with a short exact sequence
 \begin{equation}\label{eq:z-ext-defined}
 	1 \To Z_1 \To G_1 \To G \To 1
 \end{equation}
 such that $G_{1,\textup{der}}=G_{1,\textup{sc}}$, $Z_1\subset Z_{G_1}$, and $Z_1\cong \prod_{i\in I} \Res_{F_i/F}\GG_m$ for  finite extensions $F_i$ of $F$ over a finite index set $I$. We also call such a short exact sequence itself a $z$-extension of $G$.

 \begin{lem}\label{lem:z-ext-choice}
 	If $F$ is non-archimedean and if $G$ is unramified, then there exists a $z$-extension $G_1$ of $G$ that is unramified. 
 	Similarly, if $F$ is a number field and if $G$ is unramified at a finite set of finite places $S$, then there exists a $z$-extension $G_1$ of $G$ that is  unramified at $S$. 
 \end{lem}
 
 \begin{proof}
 	The lemma follows from \cite[Prop.~3.1]{MS82}, possibly except for the point, pertaining to the global case, that there exists a maximal torus $T$ of $G$ such that $T$ splits over an extension of $F$ unramified at $S$. Let us verify this.
 	
 	Since $G$ is unramified over $F_v$ for each $v\in S$, so is $G_{\sconn}$. Thus there exists an unramified maximal torus $T_{\sconn,v}$ in $G_{\sconn,F_v}$. Write $T_{\sconn,v}(F_v)_{\textup{rs}}\subset T_{\sconn,v}(F_v)$ for the open subset consisting of regular semi-simple elements.
 	Then the non-empty subset
 	$$Y_v:=\bigcup_{g_v\in G_{\sconn}(F_v)} g_v \cdot T_{\sconn,v}(F_v)_{\textup{rs}} \cdot g_v^{-1} \subset  G_{\sconn}(F_v)$$  
 	is open in $G_{\sconn}(F_v)$ by Harish-Chandra's submersion principle \cite{HCsubmersion}.
 	By weak approximation for $G_{\sconn}$, there exists an element $\gamma_0\in G_{\sconn}(F)\cap (\prod_{v\in S} Y_v)$. 	Let $T'$ denote the connected centralizer in $G_{\sconn}$ of $\gamma_0$. Then $T'$ is unramified at $S$ as it is conjugate to $T_{\sconn,v}$ at each $v\in S$. The obvious image of $T'\times Z_G^0$ in $G$ is then a  maximal torus of $G$ which is unramified at $S$.
 \end{proof}
 
 \begin{lem}\label{lem:z-ext-H}
 	Fix a $z$-extension of $G$ as in \eqref{eq:z-ext-defined}.
 	For each $\fke=(H,\cH,s,\eta)\in E(G)$, there exists a central extension $1\ra Z_{1}\ra H_1\ra H\ra 1$ over $F$, with $H_1$ reductive , such that the induced short exact sequence $1\ra Z_1 \ra Z_{H_1} \ra Z_H\ra 1$ fits in the following row-exact commutative diagram:
 	$$ \xymatrix{ 1 \ar[r] & Z_{1} \ar[r] \ar@{=}[d] & Z_{G_1} \ar[r] \ar[d] & Z_G \ar[d]^-{\eqref{eq:ZGZH}} \ar[r] & 1
 		\\ 1 \ar[r] & Z_{1} \ar[r] & Z_{H_1} \ar[r] & Z_H \ar[r] & 1. }$$
 	Moreover, when $F$ is non-archimedean, we can choose $H_1$ to be unramified if $H$ and $G_1$ are unramified. When $F$ is a number field, $H_1$ can be chosen to be unramified at a finite set of places $S$ if $H$ and $G_1$ are unramified at $S$.
 \end{lem}
 
 \begin{rem}
 	We do not claim that $H_1$ is a $z$-extension of $H$, that is, the derived subgroup of $H_1$ need not be simply connected.
 \end{rem}
 
 \begin{proof}
 	We may assume that $G_1$ and $G$ are quasi-split. Indeed, we can replace $G_1$ by a quasi-split inner form $G_1^*$ and identify $Z_1\subset G_1$ as a central $F$-subgroup of $G_1^*$ to take the quotient $G^*:=G_1^*/Z_1$. Then the top row in the above diagram is identified with a short exact sequence $1\ra Z_1\ra Z_{G_1^*}\ra Z_{G^*}\ra 1$.
 	
 	Choose maximal tori $T_H\subset H$ and $T\subset G$ over $F$ together with an isomorphism $T_H\cong T$ over $F$ as in \S\ref{Z(G)-Z(H)}. We pull back the resulting embedding $Z_H\subset T_H \cong T \subset G$ via the surjection $G_1\rightarrow G$ to obtain the preimage $\tilde Z_1 \subset G_1$ fitting in the exact sequence
 	\begin{equation}\label{e:ZZZ}
 		1\To Z_1\To  \tilde Z_1\To Z_H\To 1.
 	\end{equation}

 	Let $T_{H_{\sconn}},T_{\sconn}$ denote the preimages of $T_H,T$ in $H_{\sconn},G_{\sconn}$, respectively. The cocharacter lattices of $T_{H_{\sconn}},T_{\sconn}$ are identified with the coroot lattices of $T_H,T$ in $H,G$, respectively. Since the coroots of $T_H$ in $H$ form a subset in the set of coroots of $T$ in $G$ via $T_H\cong T$, we have a map $T_{H_{\sconn}}\ra T_{\sconn}$. The composite map $Z(H_{\sconn})\hookrightarrow T_{H_{\sconn}}\ra T_{\sconn}\hookrightarrow G_{\sconn}=G_{1,\der}$ has image in $\tilde Z_1$ because its image in $G$ lies in $Z_H$.
 	
 	Recall from \cite[\S 2.0.1]{deligne1979varietes} and the notation therein that $H=H_{\sconn} *_{Z(H_{\sconn})} Z_H$. Setting
 	$H_1:=H_{\sconn} *_{Z(H_{\sconn})} \tilde Z_1$ via the map $Z(H_{\sconn})\ra \tilde Z_1$ in the preceding paragraph,
 	we have a surjection $H_1\ra H$ induced by $\tilde Z_1\ra Z_H$, whose kernel is identified with $Z_1$. The construction yields $Z_{H_1}=\tilde Z_1$ and the commutative diagram of the lemma. To verify the last two assertions of the lemma, it suffices to choose the torus $T_H$ at the start of the proof to be unramified if $F$ is local and unramified at $S$ if $F$ is global, as in the proof of Lemma \ref{lem:z-ext-choice}. (E.g., if $F$ is local, since $Z_1$ and $T$ are unramified, the preimage of $T$ in $G_1$ is unramified. Hence $\tilde Z_1$ is also unramified.)
 \end{proof}
 
 \subsection{}\label{G1H1}
 Let $F$ be a local or global field of characteristic $0$.
 Fix a $z$-extension of $G$ as in \eqref{eq:z-ext-defined}.
 Let us explain how each $\fke=(H,\cH,s,\eta)\in E(G)$ gives rise to an endoscopic datum for $G_1$.
 
 Fix inner twistings to quasi-split inner forms $H_{\ol F}\cong H^*_{\ol F}$ and $G_{\ol F}\cong G^*_{\ol F}$ together with $F$-pinnings for $H^*$ and $G^*$. Along the central extension $H_1\ra H$ provided by the preceding lemma, we can lift inner twistings to obtain $H_{1,\ol F}\cong H_{1,\ol F}^*$ and $G_{1,\ol F}\cong G_{1,\ol F}^*$ as well as $F$-pinnings for $H_1^*$ and $G_1^*$. As explained in \cite[1.8]{Kot84a}, we obtain
 $\Gamma$-equivariant maps $\hat H\ra \hat H_1$, $\hat G\ra \hat G_1$, $\hat H_1\ra \hat Z_1$, and $\hat G_1\ra \hat Z_1$. We will consider the natural extension of the last three maps to $L$-morphisms $\zeta_{G_1}:{}^L G\ra {}^L G_1$, $^L H_1\ra {}^L Z_1$, and $^L G_1 \ra {}^L Z_1$, respectively. (As for $\hat H\ra \hat H_1$, we have $\zeta_{H_1}$ in the lemma below.) The composition $\hat H\ra \hat G\ra \hat G_1$ factors through the embedding $\hat H \ra \hat H_1$ to yield the following commutative diagram.
 \begin{equation}\label{e:H-H_1-G-G_1}
 	\xymatrix{ \hat H \ar[r] \ar[d]  & {}\hat H_1 \ar[d] \\ \hat G \ar[r] & \hat G_1}
 \end{equation}

 \begin{lem}\label{lem:endoscopic-datum-z-ext}
 	Retain the notation of \S\ref{G1H1}.
 	\begin{enumerate}
 		\item The embedding $\hat H \hra \hat H_1$ can be extended to an $L$-morphism $\zeta_{H_1}:\cH\ra {}^L H_1$ such that $\zeta_{H_1}$ induces a homeomorphism from $\cH$ onto its image.
 		\item There exists an $L$-morphism $\eta_1:{}^L H_1\ra {}^L G_1$ such that the following is a commutative diagram extending \eqref{e:H-H_1-G-G_1}.
 		$$ \xymatrix{ \cH \ar[r]^-{\zeta_{H_1}} \ar[d]^-{\eta} & {}^L H_1 \ar[d]^-{\eta_1} \\ ^L G \ar[r]^-{\zeta_{G_1}} & ^L G_1}$$
 		\item The quadruple $\fke_1:=(H_1,{}^L H_1,s_1,\eta_1)$, with $s_1:=\zeta_{G_1}(s)$, is an endoscopic datum for $G_1$. The isomorphism class of $\fke_1$ is independent of the choices in (i) and (ii).
 	\end{enumerate}
 \end{lem}
 
 \begin{proof}
 	To verify (i), consider the split extension of $W_F$ by $\hat H_1$ given by $\cH_1:=(Z(\hat H_1)\rtimes \cH)/Z(\hat H)$, where $Z(\hat H)$ embeds in the semi-direct product diagonally. The assignment $h\mapsto (1\rtimes h)$ induces an embedding $\cH\hra \cH_1$ extending the map $\hat H\hra \hat H_1$. As remarked in \S\ref{para:z-ext}, $\cH_1\cong {}^L H_1$ since $G_{1,\textup{sc}}=G_{1,\der}$, so we obtain the desired map $\zeta_{H_1}$ by composition. Let us verify (ii). Since $G_{1,\textup{sc}}=G_{1,\der}$, we can extend $\hat H_1\ra \hat G_1$ to $\eta_0:{}^L H_1\ra {}^L G_1$ by \cite[Prop.~1]{Lan79}. The two $L$-morphisms $\zeta_{G_1}\eta$ and $\eta_0\zeta_{H_1}$ differ by a 1-cocycle $a:W_F\ra C$, where $C:=\Cent(\hat H,\hat G_1)$. Clearly $Z(\hat H_1)\subset C$. We also note that $Z(\hat H_1)\cap (C\cap \hat G)=Z(\hat H_1)\cap \hat G=Z(\hat H)$. Thus
 	$$\hat Z_1=Z(\hat H_1)/Z(\hat H)\subset C/C\cap \hat G \subset \hat G_1/\hat G=\hat Z_1,$$
 	implying that $C=Z(\hat H_1)$.
 	As $a$ is valued in $Z(\hat H_1)$, one can twist $\eta_0$ by $a$ to obtain $\eta_1$, which then makes the diagram commute.	Lastly (iii) is a routine check.
 \end{proof}
 
 \subsection{}\label{para:lambda1}

 Given a central extension $H_1$ of $H$ as above, choose a splitting $W_F\ra \cH$ to consider the composition
 $$W_F \To \cH \xrightarrow{\zeta_{H_1}} {}^L H_1 \To  {}^L Z_1.$$
 Write $\lambda_1:Z_1(F)\ra \CC^\times$ if $F$ is local, or $\lambda_{1}:Z_1(F)\bs Z_1(\A_F)\ra \CC^\times$ if $F$ is global, for the corresponding continuous character via the (geometrically normalized) Langlands correspondence for tori. Then $\lambda_1$ is independent of the choice of splitting. This character naturally shows up in endoscopic transfer. (Our $\lambda_1$ is exactly the character $\lambda_{H_1}$ described in \cite[p.23]{KS99}.)
 
 We show that the assignment $\fke\ra \fke_1$ admits an inverse map.

 \begin{lem}\label{lem:fke_1=fke}
 	The map $\fke\mapsto \fke_1$ defined by Lemma \ref{lem:endoscopic-datum-z-ext} induces a bijection from $\mathcal E(G)$ onto $\mathcal E(G_1)$.
 \end{lem}

 \begin{proof}
 	The injectivity is easy to see since $\hat G_1$ is generated by $\hat G$ and $Z(\hat G_1)$.
 	
 	To prove the surjectivity, let $\fke_1=(H_1,{}^L H_1,s_1,\eta_1)\in E(G_1)$.
 	Replacing $\fke_1$ by an isomorphic datum, we may assume that $\eta_1(W_F)$ lies in the subgroup $^L G$ of $^L G_1$. Indeed, consider the exact sequence of continuous cohomology
 	$$\coh^1(W_F,\hat G) \ra \coh^1(W_F,\hat G_1)\ra \coh^1(W_F,\hat Z_1).$$
 	The image of $\eta_1$ under the second map lifts to a 1-cocycle $c$ valued in $Z(\hat G_1)$, up to a 1-coboundary, via the map $\coh^1(W_F,Z(\hat G_1))\ra \coh^1(W_F,\hat Z_1)$, which is surjective by \cite[Lem.~4]{Lan79}. Then $c^{-1}\cdot \eta_1$ comes from a 1-cocycle valued in $\hat G$ up to a 1-coboundary. 
 	
 	We define $H$ to be the cokernel of the composite map $Z_1\hra Z_{G_1}\hra Z_{H_1}\hra H_1$.
 	We can write $s_1=sz$ for some $s\in \hat G$ and $z\in Z(\hat G_1)$. By pulling back $\eta_1:{}^L H_1\ra {}^L G_1$ via $^L G\hookrightarrow {}^L G_1$, we obtain an injection $\eta:\cH\ra {}^L G$ and see that $\cH$ is a split extension of $W_F$ by $\hat H$ (as $\cH$ is generated by $\eta_1(W_F)$ and $\hat H$).
 	
 	It is enough to verify that $\fke:=(H,\cH,s,\eta)$ is an endoscopic datum for $G$, since it would then be obvious by construction that $\fke\mapsto\fke_1$, and we will be done. The main point is to show that $\Int(s)\circ \eta=a\eta$ with trivial (resp.~locally trivial) 1-cocycle $a:W_F\ra Z(\hat G)$ if $F$ is local (resp.~global). Let us check this when $F$ is global, as the local case is simpler. Since $\fke_1\in \cE(G_1)$ we know that $\Int(s)\circ \eta_1=\Int(s_1)\circ \eta_1=a_1\eta_1$ with a locally trivial 1-cocycle $a_1:W_F\ra Z(\hat G_1)$. Since $\eta_1(W_F)\subset {}^L G$ we have that $a_1(W_F)\subset Z(\hat G_1)\cap \hat G=Z(\hat G)$ and that $\Int(s)\circ \eta=a_1\eta$. As $G_1$ is a $z$-extension of $G$, the map $\coh^1(W_{F_v},Z(\hat G))\ra \coh^1(W_{F_v},Z(\hat G_1))$ is injective at each place $v$ (e.g., by \cite[Cor.~2.3]{Kot84a}). Hence $a_1$ is locally trivial as a cocycle valued in $Z(\hat G)$.
 \end{proof}

 \subsection{}\label{subsub:twisted-end-at-p}
 
 From here until the end of \S\ref{sub:endoscopic-data}, $F$ is assumed to be non-archimedean. Write $F_m$ for the unramified extension of $F$ of degree $m\in \Z_{\ge 1}$ in a fixed algebraic closure $\ol F$. Denote by $\sigma\in \Gal(F_m/F)$ the arithmetic Frobenius generator. Fix a $z$-extension $1\ra Z_1\ra G_1\ra G\ra 1$. Set $R:=\Res_{F_m/F} G$ and $R_1:=\Res_{F_m/F} G_1$. Write $\theta$ (resp.~$\theta_1$) for the automorphism of $R$ (resp.~$R_1$) induced by $\sigma$. Identify
 $$\hat R=\hat G^{\Hom(F_m,\ol{F})}= \prod_{j=0}^{m-1} \hat G$$
 such that the $j$-th component corresponds to the inclusion $F_m\subset \ol{F}$ pre-composed by $\sigma^j$, and similarly for $\hat R_1$. There are unique embeddings $i:{}^L G\hra {}^L R$ and $i_1:{}^L G_1\hra {}^L R_1$ such that the maps are diagonal embeddings on the dual groups and the identity map on the Weil groups.

 The following is a variant of Lemma \ref{lem:z-ext-H}. In practice $G_1$ and $H_1$ over $F$ will come from central extensions over $\Q$ independently of $m$. By contrast, the extensions $G'_1$ and $H'_1$ below depend on $m$ and will be considered only in a local setting.

 \begin{lem}\label{lem:z-ext2}
 	Let $\fke=(H,\cH,s,\eta)\in E(G)$. Consider central extensions $1\ra Z_{1}\ra G_1\ra G\ra 1$ and $1\ra Z_{1}\ra H_1\ra H\ra 1$ as constructed in Lemma \ref{lem:z-ext-H} (disregarding the last assertion). Let $F'/F$ be a finite extension. Then there exist 
 	\begin{enumerate}
 		\item a $z$-extension $1\ra Z'_1\ra G'_1\ra G\ra 1$ such that $Z'_1\cong \prod_{j\in J}\Res_{F'_j/F}\GG_m$ with $J$ a finite index set and $F'_j \supset F'$, and
 		\item a central extension $1\ra Z'_1\ra H'_1\ra H\ra 1$ arising from (i) as in Lemma \ref{lem:z-ext-H},
 	\end{enumerate}	
 	and such that
 	there is a commutative diagram with an injective middle vertical arrow
 	\begin{equation}\label{eq:z-ext2}
 		\xymatrix{ 1 \ar[r] & Z_1 \ar[r] \ar[d] & G_1 \ar[r] \ar@{^(->}[d] & G \ar@{=}[d] \ar[r] & 1
 			\\ 1 \ar[r] & Z'_1 \ar[r] & G'_1 \ar[r] &G \ar[r] & 1 },
 	\end{equation} 
 	as well as the analogous diagram with $H,H_1,H'_1$ in place of $G,G_1,G'_1$. Moreover, if $F'$, $G_1$ and $\fke$ (thus also $H_1$) are unramified over $F$, then $G'_1$ and $H'_1$ are unramified over $F$.
 \end{lem}
 
 \begin{proof}
 	We have $Z_1\cong \prod_{j\in J} \Res_{E_j/F}\GG_m$ for finite extensions $E_j/F$. Take $Z'_1:=\prod_{j\in J} \Res_{E_jF'/F}\GG_m$. Along the canonical inclusion $Z_1\hra Z'_1$, we make a pushout diagram from the top row of Lemma \ref{lem:z-ext-H} as follows:
 	$$ \xymatrix{ 1 \ar[r] & Z_1 \ar[r] \ar[d] & Z_{G_1} \ar[r] \ar[d] & Z_G \ar@{=}[d] \ar[r] & 1
 		\\ 1 \ar[r] & Z'_1 \ar[r] & Z_{G_1}' \ar[r] & Z_G \ar[r] & 1 }$$
 	With the bottom row in place of \eqref{e:ZZZ} we construct a $z$-extension $1\ra Z'_1 \ra G'_1\ra G \ra 1$ such that $Z_{G'_1}=Z_{G_1}'$ as in the proof of Lemma \ref{lem:z-ext-H}. By construction, we have (i) and \eqref{eq:z-ext2}.
 	Applying Lemma \ref{lem:z-ext-H} to $G'_1$, we obtain (ii). 
 	Since $Z_{H_1}$ and $Z_{H'_1}$ are preimages of $Z_H\subset G$ under $G_1\rightarrow G$ and $G'_1\rightarrow G$, we have $Z_{H_1}\subset Z_{H'_1}$. This in turn induces $H_1=H_{\textup{sc}} *_{Z(H_{\textup{sc}})} Z_{H_1} \subset H'_1=H_{\textup{sc}} *_{Z(H_{\textup{sc}})} Z_{H'_1}$. With this inclusion as the middle vertical arrow, we see that there is a commutative diagram analogous to \eqref{eq:z-ext2} for $H,H_1,H'_1$. Finally, if $F'$, $G_1$, and $\fke$ are unramified over $F$, then $E_j/F$ are unramified in particular, and the torus $Z'_1$ is unramified by construction. Since $H,G$ are unramified, it follows that $H'_1,G'_1$ are unramified. 
 \end{proof}
 
 \subsection{}\label{para:twisted-datum}
 
 Under the temporary assumption on $\fke=(H,\cH,s,\eta)\in E(G)$ that
 $$s\in Z(\hat H)^{\Gamma_F}$$
 (in general $s\in Z(\hat H)^{\Gamma_F} Z(\hat G)$), we construct some twisted endoscopic data to be used in the stabilization (\S\ref{sec:stabilization}). Put $\tilde s:=(s,1,\ldots,1)\in \hat R$,  which lies in $\cZ:=\Cent(i\eta(\hat{H}),\hat R)$. Equip $\cZ$ with a $W_F$-action by $w\cdot z= \Int(i\eta(h_w))z$ for any choice of $h_w\in \cH$ lying over $w\in W_F$; the action is independent of the choice. Define the $L$-morphism $\tilde \eta: \cH\ra {}^L R$ to be the twist of $i\eta$ by the unramified 1-cocycle $a:W_F\ra \cZ$ mapping $\sigma$ to $\tilde s$. Then\index[Oct2026-sorted-notation]{e 1 04 23 tilde e h h tilde s tilde eta@$\tilde{\mathfrak e}=(H,\cH,\tilde s, \tilde \eta)$} $\tilde\fke=(H,\cH,\tilde s, \tilde \eta)$ is checked to be a twisted endoscopic datum  for $(R,\theta)$, cf.~ \cite[A.1.3, A.2.6]{Mor10} or \cite[\S7]{Kot90}. Replacing $\fke$ by $\fke_1$ (noting that the temporary assumption is still satisfied for $\fke_1$, i.e., $s_1\in Z(\hat H_1)^{\Gamma_F}$) we construct a twisted endoscopic datum $\tilde\fke_1=(H_1,{}^L H_1,\tilde s_1,\tilde \eta_1)$ for $(R_1,\theta_1)$. With $H'_1$ playing the role of $H_1$, we also construct $\fke'_1$ and $\tilde\fke'_1$.

 \section{Cohomological lemmas}\label{sub:cohomological}
 
 This section studies certain Galois cohomological invariants with respect to $z$-extensions.
 
 \subsection{}
 Let $F$, $G$, and $\fke$ be as in \S\ref{sub:endoscopic-data}. (The field $F$ is either local or global.) Take a $z$-extension $G_1$ as in Lemma \ref{lem:z-ext-H}. By Hilbert 90 the map $G_1(F)\ra G(F)$ is onto. Let $\gamma\in G(F)_{\semi}$ and choose a lift $\gamma_1\in G_1(F)$. We have a commutative diagram of reductive groups
 \begin{equation}\label{e:Z-I-G}
 	\xymatrix{ 1 \ar[r] & Z_{1} \ar[r] \ar@{=}[d] & I_{\gamma_1} \ar[r] \ar[d] & I_{\gamma} \ar[r] \ar[d] & 1\\ 1 \ar[r] & Z_{1} \ar[r] & G_1 \ar[r] & G \ar[r] & 1,}
 \end{equation}
 which gives rise to a $\Gamma$-equivariant commutative diagram by \cite[1.8]{Kot84a}:
 \begin{equation}\label{e:xymatrix-Z-I-G}
 	\xymatrix{ 1 \ar[r] & Z(\hat I_{\gamma}) \ar[r]  & Z(\hat I_{\gamma_1}) \ar[r]  &  \hat Z_{1}\ar[r]  & 1\\ 1 \ar[r] & Z( \hat G ) \ar[r] \ar[u] & Z( \hat G_1 ) \ar[u] \ar[r] &  \hat Z_{1} \ar[r] \ar@{=}[u] & 1.}
 \end{equation}
 In particular we get a $\Gamma$-equivariant isomorphism $$Z(\hat I_\gamma)/Z(\hat G)\cong Z(\hat I_{\gamma_1})/Z(\hat G_1). $$
 
 Recall the definition of $\mathfrak K(I_\gamma/F)$ from \cite[\S4]{Kot86} (which depends on the ambient group $G$). The exact sequence
 $
 1\rightarrow Z(\widehat G)
 \rightarrow Z(\widehat I_\gamma)
 \rightarrow Z(\widehat I_\gamma)/Z(\widehat G)
 \rightarrow1
 $
 gives a connecting homomorphism
 \[
 \partial_F:
 \pi_0\!\left(
 [Z(\widehat I_\gamma)/Z(\widehat G)]^{\Gamma_F}
 \right)
 \longrightarrow H^1(F,Z(\widehat G)).
 \]
 If $F$ is local, set
 $
 \mathfrak K(I_\gamma/F):=\ker(\partial_F).
 $
 If $F$ is global, 
 $
 \mathfrak K(I_\gamma/F)
 :=
 \partial_F^{-1}(\Sha(F,Z(\widehat G))).
 $
 For $F=\mathbf Q$, this agrees with the definition recalled in Corollary \ref{cor:compare with Kottwitz}.
 
 \begin{lem}\label{lem:z-ext-wrt-K-and-D}
 	We have 
 	\begin{enumerate}
 		\item If $F$ is global or local, there is a canonical isomorphism 
 		$$\fkK(I_\gamma/F)\cong \fkK(I_{\gamma_1}/F).$$
 		\item If $F$ is local, there is a canonical bijection 
 		$$\fkD(I_{\gamma_1},G_1;F)\cong \fkD(I_{\gamma},G;F).$$
 	\end{enumerate}
 \end{lem}
 
 \begin{proof}
 	Let us check (i) in the global case. We obtain the following commutative diagram from \cite[Cor.~2.3]{Kot84a} and diagram \eqref{e:xymatrix-Z-I-G}:
 	$$\xymatrix{ \pi_0((Z(\hat I_\gamma)/Z(\hat G))^{\Gamma_F}) \ar[r]  \ar[d]^-{\sim} & \coh^1(F,Z(\hat G)) \ar[d] \\ \pi_0((Z(\hat I_{\gamma_1})/Z(\hat G_1))^{\Gamma_F}) \ar[r] & \coh^1(F,Z(\hat G_1)). }$$
 	The left vertical map is an isomorphism. The right vertical map is injective by \cite[Cor.~2.3]{Kot84a} since $\hat {Z}_1^{\Gamma_F}$ is connected (a product of copies of $\CC^\times$).
 	Further observe that the injective right vertical map induces an isomorphism $\Sha(F,Z(\hat G))\cong \Sha(F,Z(\hat G_1))$. To see this, notice that the cokernel is mapped injectively into $\Sha(F,\hat Z_1)$, which is trivial since $Z_1$ is a product of induced tori.
 	Since $\fkK(I_\gamma/F)$, $\fkK(I_{\gamma_1}/F)$ are the preimages of $\Sha(F,Z(\hat G))$, $\Sha(F,Z(\hat G_1))$ under the horizontal arrows, the left vertical map induces the desired isomorphism of (i). In the case of local fields, the same argument works if we replace the $\Sha$-groups with the trivial group.

 	The two pointed sets in (ii) fit in the following commutative diagram, where the row exact sequences come from \cite[p.42]{Lab99}, and the vertical maps are induced by surjections $G_1(F)\twoheadrightarrow G(F)$ and $I_{\gamma_1}(F)\twoheadrightarrow I_\gamma(F)$ as well as an $F$-isomorphism $I_{\gamma_1}\backslash G_1 \isom I_\gamma\backslash G$ coming from \eqref{e:Z-I-G} by Hilbert 90. 
 	$$\xymatrix{ 1\ar[r] & I_{\gamma_1}(F)\backslash G_1(F) \ar[d]^-{\sim}\ar[r] & \coh^0(F,I_{\gamma_1}\backslash G_1) \ar[r] \ar[d]^-{\sim} & \fkD(I_{\gamma_1},G_1;F) \ar[d]^-{\mathrm{pr}}\ar[r] & 1\\
 		1\ar[r] & I_{\gamma}(F)\backslash G(F) \ar[r] & \coh^0(F,I_{\gamma}\backslash G) \ar[r]  & \fkD(I_{\gamma},G;F) \ar[r] & 1.}
 	$$  Hence $\mathrm{pr}$ is a surjection whose kernel is a singleton. A twisting argument shows that the fiber of $\mathrm{pr}$ is a singleton over every element  $a\in \fkD(I_{\gamma},G;F)$. That is, suppose $\mathrm{pr}(a_1)=a$. We can represent $a_1$ by $\sigma\mapsto g_1^{-1}{}^\sigma g_1$ for some $g_1\in G_1(\ol{F})$; then $a$ is represented by $g^{-1}{}^\sigma g$, where $g\in G(\ol F)$ is the image of $g_1$. Then $\gamma'_1:=g_1 \gamma_1 g_1^{-1}$ and $\gamma':=g \gamma g^{-1}$ are $F$-rational. Applying the bijection $(b_\sigma)_{\sigma}\mapsto (g b_\sigma {}^\sigma g^{-1})_{\sigma}$ on $\fkD(I_{\gamma},G;F)$ and similarly on $\fkD(I_{\gamma_1},G_1;F)$, we carry $\mathrm{pr}^{-1}(a)$ bijectively onto the kernel of $\fkD(I_{\gamma'_1},G_1;F)\ra \fkD(I_{\gamma'},G;F)$, and the latter is a singleton by the preceding argument. 
 \end{proof}

 \subsection{}\label{Kottwitz parameters z-ext}
 
 Let $(G,X)$ be a Shimura datum.
 We study Kottwitz parameters and their invariants with respect to a $z$-extension
 $1\ra Z_1 \ra G_1 \ra G\ra 1$ over $\Q$.
 Let $T\subset G_{\R}$ be an elliptic maximal torus.
 There exists $h\in X$ such that $\mu_h:\GG_m\ra G_\CC$ factors through $T_{\CC}$.
 In the notation of \S\ref{subsubsec:SD}, $\mu_h\in \dmu_X(\CC)$.
 Write $T_1\subset G_{1,\R}$ for the preimage of $T$. Then $\mu_h:\GG_m\ra T_\CC$ lifts to a cocharacter $\mu_1:\GG_m\ra T_{1,\CC}$, which we fix henceforth and view also as a cocharacter of $G_1$ over $\CC$. As noted in \cite[3.4]{MS82}, the conjugacy class of $\mu_1$ comes from a Shimura datum $(G_1,X_1)$ for a suitable $X_1$, that is, $\mu_1\in \dmu_{X_1}(\CC)$.
 In particular the discussion of cohomological invariants (\S\S\ref{defn of Kottwitz invariant}--\ref{subsub:comparison_Kottwitz_triples}) applies to $(G_1,X_1)$ and $\mu_1$.
 
 Let $\gamma_{0,1}\in G_1(\A_f^p)$ and $\gamma_{0}\in G(\A_f^p)$ such that $\gamma_{0,1}$ maps to $\gamma_0$. Write $I_{0,1}$ and $I_0$ for the connected centralizers of $\gamma_{0,1}$ and $\gamma_0$ in $G_1$ and $G$ over $\A_f^p$, respectively.

 \begin{lem}
 	Suppose that $\gamma_0\in G(\A_f^p)$ is the image of an element $\gamma_{0,1}\in G_1(\A_f^p)_{\textup{ss}}$. Then the natural map $\fkD(I_{0,1},G_1;\A_f^p)\ra \fkD(I_0,G;\A_f^p)$ is a bijection.
 	
 \end{lem}
 
 \begin{proof}
 	This follows from \eqref{eq:adelic coh} and part (ii) of Lemma \ref{lem:z-ext-wrt-K-and-D}.
 \end{proof}
 
 \subsection{}\label{para:Kottwitz-triples-z-ext-1}
 
 Now consider a $z$-extension $1\ra Z_1 \ra G_1 \ra G\ra 1$ over $\Q_p$ (which need not come from $\Q$-groups via base change). Let $\mu:\GG_m\ra G_{\Qpbar}$ be a cocharacter over $\Qpbar$. Let $\gamma_{0,1}\in G_1(\Q_p)_{\textup{ss}}$, and $\gamma_0\in G(\Q_p)$ the image of $\gamma_{0,1}$. As usual, $I_0$ and $I_{0,1}$ are the connected centralizers of $\gamma_0$ and $\gamma_{0,1}$ in $G$ and $G_1$ over $\Q_p$, respectively. Fix a level $n\in \Z_{>0}$. Set $R:=\Res_{\Q_{p^n}/\Q_p}G$ and $R_1:=\Res_{\Q_{p^n}/\Q_p}G_1$.
  
 Denote by\index[Oct2026-sorted-notation]{d 0 04 00 d subscript n gamma subscript 0 g q subscript p@$\mathfrak D_n(\gamma_0,G;\Q_p)$} $\fkD_n(\gamma_0,G;\Q_p)$ the set of all $[b]\in \B(I_{0})$ satisfying condition \textbf{KP1} of Definition \ref{pn-adm} with the given $\gamma_0$ and $n$. (Here we do not require $\textbf{KP0}$, which  will come into play in Corollary \ref{cor:transfer-at-p} below.) This means that for some (thus every) representative $b\in I_0(\LL)$ of $[b]$, there exists $c\in G(\LL)$ such that $$c^{-1} \gamma_0 c = c^{-1} b \sigma( b) \cdots {\sigma^{n-1}} (b) {\sigma^n}(c). $$
 Given $[b]\in \fkD_n(\gamma_0,G;\Q_p)$, we will often write $\delta_{[b]}\in  G(\Q_{p^n})$ for the element arising from $[b]$ (well defined up to $\sigma$-conjugacy in $G(\Q_{p^n})$) as in Lemma \ref{b decides delta}. Then $\gamma_0$ is a degree $n$ norm of $\delta_{[b]}$.
 
 Likewise we define $\fkD_n(\gamma_{0,1},G_1;\Q_p)$. The natural map $I_{0,1}\ra I_0$ induces a map
 \begin{equation}\label{eq:D01-to-D0}
 	\fkD_n(\gamma_{0,1},G_1;\Q_p)
 	\To  \fkD_n(\gamma_{0},G;\Q_p).
 \end{equation}

 \begin{lem}\label{lem:Kottwitz-triples-z-ext-1}
 	Let $\gamma_0,\gamma_{0,1}$ be as above. Suppose that $Z_1\cong \prod_{j\in J} \Res_{F_j/\Q_p}\GG_m$ with all $F_j$ containing $\Q_{p^n}$ and $J$ a finite index set.
 	Then the map \eqref{eq:D01-to-D0} is a bijection.
 \end{lem}
 
 \begin{proof}
 	The map \eqref{eq:D01-to-D0} fits in the following commutative diagram
 	$$\xymatrix{ \fkD_n(\gamma_{0,1},G_1;\Q_p) ~\ar@{^(->}[r] \ar[d]_{\eqref{eq:D01-to-D0}} & \B(I_{0,1}) \ar[r] \ar[d]& \B(R_1) \ar[d] \\
 		\fkD_n(\gamma_{0},G;\Q_p) ~\ar@{^(->}[r] & \B(I_{0}) \ar[r]  & \B(R),}$$
 	where the vertical maps are induced by the natural maps $I_{0,1}\ra I_0$ and $R_1\ra R$.
 	
 	We may assume that $\fkD_n(\gamma_{0},G;\Q_p)$ is non-empty since the lemma is vacuously true otherwise.
 	We claim that $\fkD_n(\gamma_{0,1},G_1;\Q_p)$ is also non-empty. To see this, fix $[b]\in \fkD_n(\gamma_{0},G;\Q_p)$ and pick a lift $\delta'_1\in G_1(\Q_{p^n})$ of $\delta_{[b]}$. Let $\gamma'_{0,1}\in G_1(\Q_p)$ be a degree $n$ norm of $\delta'_1$. Then the norm of $z\delta'_1$ is $z{}^\sigma z \cdots {}^{\sigma^{n-1}} z \gamma'_{0,1}$ with $z\in Z_1(\Q_{p^n})$. The norm map $Z_1(\Q_{p^n})\ra Z_1(\Q_p)$ is onto by the hypothesis on $Z_1$, so we may choose $z$ such that $\gamma_{0,1}$ is a norm of $\delta_1:=z\delta'_1$. This implies that there is $c_1\in G_1(\LL)$ such that $c_1^{-1} \gamma_{0,1} c_1 = \delta_1 {}^\sigma \delta_1 \cdots {}^{\sigma^{n-1}} \delta_1$.  Setting $b_1:=c_1\delta_1 {}^\sigma c_1^{-1}$, we see that $b_1\in I_{0,1}(\breve \Q_p)$ (since $G_1$ has simply connected derived subgroup) and that $[b_1]\in \fkD_n(\gamma_{0,1},G_1;\Q_p)$, proving the claim. 
 	
 	We fix $[b]$, $\delta_1$, and $[b_1]$ as in the last paragraph. In particular $\delta_1$ maps to $\delta:=\delta_{[b]}$, and $[b_1]$ to $[b]$.
 	By Corollary \ref{b is basic}, $[b_1]\in \B(I_{0,1})$ and $[b]\in \B(I_0)$ are basic.  Define $\fkD([b_1],R_1)$ to be the set of basic elements $[b'_1]\in \B(I_{0,1})$ such that $\kappa_{I_{0,1}}([b'_1])-\kappa_{I_{0,1}}([b_1])$ lies in $\ker(\pi_1(I_{0,1})_{\Gamma_p,\tors}\to \pi_1(R_1)_{\Gamma_p,\tors})$. Define\index[Oct2026-sorted-notation]{d 0 04 00 d b r leftbracket rightbracket@$\mathfrak D([b],R)$}\index[Oct2026-sorted-notation]{d 0 04 00 d b subscript 1 r subscript 1 leftbracket rightbracket@$\mathfrak D([b_1],R_1)$} $\fkD([b],R)$ exactly in the same way with $I_0$ and $R$ in place of $I_{0,1}$ and $R_1$.
 	
 	We claim that $\fkD_n(\gamma_{0},G;\Q_p)$ is a subset of $\fkD([b],R)$. Indeed, consider an element $[b']\in \fkD_n(\gamma_{0},G;\Q_p)$, which gives rise to $\delta'\in G(\Q_{p^n})$ as in Lemma \ref{b decides delta}. Corollary \ref{b is basic} implies that $\kappa_{I_{0}}([b'])-\kappa_{I_{0}}([b])$ is a torsion element in $\pi_1(I_{0})_{\Gamma_p}$ via the right half of the commutative diagram in \cite[p.~162]{rapoportrichartz}. We know from Lemma \ref{b decides delta} that $\delta$ and $\delta'$ have the same norm, implying that $\delta$ and $\delta'$ are stably $\sigma$-conjugate by \cite[Prop.~5.7]{Kot82}. Therefore $[\delta]$ and $[\delta']$ are equal in $\B(R)$. Since $b$ and $b'$ are $\sigma$-conjugate to $\delta$ and $\delta'$ in $G(\LL)$, respectively, it follows that $[b]=[b']$ in $\B(R)$. Hence $\kappa_{I_{0}}([b'])$ and $\kappa_{I_{0}}([b])$ have the same image in $\pi_1(R)_{\Gamma_p}$, completing the proof of the claim.
 	
 	What we have shown is summarized in the following commutative diagram.
 	\begin{equation}\label{eq:D0,1-D_0}
 		\xymatrix{ \fkD_n(\gamma_{0,1},G_1;\Q_p) ~ \ar@{^(->}[r] \ar[d] & \fkD([b_1],R_1) \ar[d]^-{\mathrm{bij}.}  \\
 			\fkD_n(\gamma_{0},G;\Q_p) ~ \ar@{^(->}[r] & \fkD([b],R)}
 	\end{equation}
 	The right vertical map is a bijection by the proof of \cite[Lem.~5.6.(2)]{Kot82}. (This relies on the assumption of the lemma on  $Z_1$.) Indeed, the proof there shows a canonical bijection from $$\ker(\coh^1(\Q_p,I_{\delta_1})\ra \coh^1(\Q_p,R_1))$$ to $$\ker(\coh^1(\Q_p,I_{\delta})\ra \coh^1(\Q_p,R)), $$ but this is exactly the right vertical map above via the functorial bijection between $\coh^1(\Q_p,H)$ and $\pi_1(H)_{\Gamma_p,\tors}$ for an arbitrary reductive  group $H$ over $\Q_p$, cf.~\cite[Prop.~1.6.7]{Lab99} and \cite[Thm.~1.15.(i)]{rapoportrichartz}.
 	
 	Now we verify that the top horizontal map in \eqref{eq:D0,1-D_0} is a bijection. As we have seen in the last paragraph, $\fkD([b_1],R_1)$ is in a canonical bijection with $\ker(\coh^1(\Q_p,I_{\delta_1})\ra \coh^1(\Q_p,R_1))$, which in turn is canonically bijective onto the set of $\sigma$-conjugacy classes in the stable $\sigma$-conjugacy class of $\delta_1$ in $G_1(\Q_{p^n})$. Recall that $\gamma_{0,1}$ is a norm of $\delta_1$. Therefore each $\delta'_1\in G_1(\Q_{p^n})$ stably $\sigma$-conjugate to $\delta_1$ gives rise to an element of $\fkD_n(\gamma_{0,1},G_1;\Q_p)$ as described in the second paragraph of the proof of the current lemma (more detailed on p.~167 of \cite{Kot90}). It is routine to check that this map gives the inverse of the top horizontal map in \eqref{eq:D0,1-D_0}.
 	
 	Going back to diagram \eqref{eq:D0,1-D_0}, it is now clear that the left vertical map (as well as the lower horizontal map) is a bijection. \end{proof}

 \begin{cor}\label{cor:bijection}
 	Assume that $\fkD_n(\gamma_{0},G;\Q_p)$ is non-empty. Then the map
 	$[b]\mapsto \delta_{[b]}$ gives a surjection from $\fkD_n(\gamma_{0},G;\Q_p)$ onto the set of $\sigma$-conjugacy classes in
 	the stable $\sigma$-conjugacy class in $G(\Q_{p^n})$ whose norm is $\gamma_0$.
 	If $G$ has simply connected derived subgroup then this map is a bijection.
 \end{cor}
 
 \begin{proof}
 	Thanks to Lemma \ref{lem:z-ext2}, we can choose a $z$-extension $G_1$ such that the condition on $Z_1$ in Lemma \ref{lem:Kottwitz-triples-z-ext-1} is satisfied.
 	Then the corollary follows from  the bijectivity of the lower horizontal map of \eqref{eq:D0,1-D_0}, together with the interpretation of $\fkD([b],R)$ in terms of $\sigma$-conjugacy classes below \eqref{eq:D0,1-D_0}.
 \end{proof}
 
 \subsection{}\label{real place z-ext}
 
 Now we turn to the real place.
 Let $(G,X)$ and $(G_1,X_1)$ be as in \S\ref{Kottwitz parameters z-ext}. 
 Following \S\ref{defn of Kottwitz invariant} and \S\ref{subsub:local_Kottwitz_inv}, for each elliptic element $\gamma_0\in G(\R)$, we have $\tilde\beta_\infty(\gamma_0)\in \pi_1(I_0)=X^*(Z(\hat I_0))$ mapping to $[\mu]\in \pi_1(G)$ for any $\mu\in \dmu_X(\Qbar)$. The definition involves an extra choice, but the restriction of $\tilde\beta_\infty(\gamma_0)$ to $Z(\hat I_0)^{\Gamma_\infty}Z(\hat G)$ is independent of the choice. Analogously, for each elliptic element $\gamma_{0,1}\in G_1(\R)$, we have $\tilde\beta_\infty(\gamma_{0,1})\in  \pi_1(I_{0,1})=X^*(Z(\hat I_{0,1}))$ mapping to $[\mu_1]\in \pi_1(G_1)$ with $\mu_1\in \dmu_{X_1}(\Qbar)$.

 \begin{lem}
 	Assume that $\gamma_{0,1}\in G_1(\R)$ and $\gamma_0\in G(\R)$ are elliptic elements such that $\gamma_{0,1}$ maps to $\gamma_0$. Then the image of $\tilde\beta_\infty(\gamma_{0,1})$ under the natural map $\pi_1(I_{0,1}) \ra \pi_1(I_{0})$ coincides with $\tilde\beta_\infty(\gamma_0)$ on $Z(\hat I_0)^{\Gamma_\infty}Z(\hat G)$.
 \end{lem}
 
 \begin{proof}
 	In \S\ref{defn of Kottwitz invariant} and \S\ref{subsub:local_Kottwitz_inv} (adapted to $(G_1,X_1)$), choose any $\R$-elliptic maximal torus $T_1\subset G_1$ containing $\gamma_{0,1}$ and $h_1\in X_1$ factoring through $T_1$ to define $\tilde\beta_\infty(\gamma_{0,1})$. Take $T\subset G$ and $h\in X$ to be the images of $T_1$ and $h_1$ in the definition of $\tilde\beta_\infty(\gamma_0)$.
 	Then the lemma is true on the nose.
 \end{proof}

 \section{Langlands--Shelstad--Kottwitz transfer}\label{sub:LSK-transfer}
 
 Our stabilization relies on the Langlands--Shelstad transfer in untwisted endoscopy and a special kind of the Kottwitz--Shelstad transfer in twisted endoscopy for local unramified base change. We recall the statements, including the fundamental lemma in each context, and slightly generalize them via $z$-extensions as necessary.
 
 \subsection{}
 Let $F$, $G$, and $\fke=(H,\cH,s,\eta)$ be as in Definition \ref{def:end-datum}.
 Throughout \S\ref{sub:LSK-transfer}, assume that $F$ is a local field (of characteristic $0$).
 Let\index[Oct2026-sorted-notation]{zz23 psi 1 00 00 psi g superscript star subscript bar f ra g subscript bar f@$\psi:G^*_{\ol F}\ra G_{\ol F}$} $\psi:G^*_{\ol F}\ra G_{\ol F}$ be an inner twisting of $F$-groups with $G^*$ quasi-split over $F$. We will use the following notation for $G$ (and likewise for other reductive groups).
 \begin{itemize}
 	\item $\Gamma(G)=\Gamma(G(F))$ is\index[Oct2026-sorted-notation]{zz03 gamma 0 00 00 gamma g@$\Gamma(G)$} the set of semi-simple $G(F)$-conjugacy classes in $G(F)$,
 	\item $\Sigma(G)=\Sigma(G(F))$ is\index[Oct2026-sorted-notation]{zz18 sigma 0 00 00 sigma g@$\Sigma(G)$} the set of stable semi-simple conjugacy classes in $G(F)$,
 	\item $\cH(G)=\cH(G(F))$ and $\cH(G,\omega^{-1})=\cH(G(F),\omega^{-1})$ as\index[Oct2026-sorted-notation]{h 0 03 00 h g@$\mathcal H(G)$}\index[Oct2026-sorted-notation]{h 0 03 00 h g omega superscript minus 1@$\mathcal H(G,\omega^{-1})$} in \S\ref{central-character-data},
 	\item When $F$ is non-archimedean and $G$ is unramified, we fix a hyperspecial subgroup $K\subset G(F)$ and define\index[Oct2026-sorted-notation]{h 0 03 00 h superscript ur g@$\mathcal H^{\ur}(G)$}\index[Oct2026-sorted-notation]{h 0 03 00 h superscript ur g omega superscript minus 1@$\mathcal H^{\ur}(G,\omega^{-1})$} $\cH^{\ur}(G)\subset \cH(G)$ and $\cH^{\ur}(G,\omega^{-1}) \subset \cH(G,\omega^{-1})$ to be the subalgebra consisting of $K$-bi-invariant functions.
 \end{itemize}

 \subsection{}\label{subs:trans-conjugacy}
 
 We explain the transfer of conjugacy classes in endoscopy following \cite[\S 3.1]{Kot86}. See also \cite[\S 1.3]{LS87}.
 
 Let $\gamma_H\in H(F)_{\semi}$. Choose a maximal torus $T_H$ of $H$ over $F$ containing $\gamma_H$. There exists a canonical $G^*(\ol{F})$-conjugacy class of embeddings $j:T_H\ra G^*$ over $\ol F$. Fix a choice of $j$ and put $T^*:=j(T_H)$. Denoting the set of absolute roots of $T_H$ in $H$ (resp.~$T^*$ in $G^*$) by\index[Oct2026-sorted-notation]{r 0 00 00 r t subscript h h@$R(T_H,H)$}\index[Oct2026-sorted-notation]{r 0 00 00 r t superscript star g superscript star@$R(T^*,G^*)$} $R(T_H,H)$ (resp.~$R(T^*,G^*)$), we have $R(T_H,H)\subset R(T^*,G^*)$. The element $\gamma_H$ is said to be\index{G,H-regular@$(G,H)$-regular} \emph{$(G,H)$-regular} if $\alpha(\gamma_H)\neq 1$ for all $\alpha\in R(T^*,G^*)\backslash R(T_H,H)$. The definition depends only on the $H(\ol{F})$-conjugacy class of $\gamma_H$ and not on the extra choices. The $(G,H)$-regular subsets of $H(F)_{\semi}$ and $\Sigma(H)$ will be denoted by\index[Oct2026-sorted-notation]{g 0 00 00 g h minus reg@$(G,H)\text{-reg}$} $H(F)_{(G,H)\text{-} \reg}$ and $\Sigma(H)_{(G,H)\text{-} \reg}$, respectively.
 
 We make the temporary assumption that $G_{\der}=G_{\sconn}$.
 Let $\gamma_H\in H(F)_{(G,H)\text{-} \reg}$. The $\ol{F}$-embedding $\psi\circ j:T_H\ra G$ is canonical up to $G(\ol{F})$-conjugacy, so $\gamma_H$ determines a semi-simple $G(\ol{F})$-conjugacy class in $G(\ol{F})$ defined over $F$.  If this conjugacy class contains an element $\gamma$ of $G(F)$, then we take the stable conjugacy class of $\gamma$ to be the image of $\gamma_H$. Otherwise the image of $\gamma_H$ is formally denoted by the empty set symbol $\emptyset$. 
 (Such a $\gamma$ always exists by \cite[Thm.~4.4]{Kot82} if $G$ is quasi-split and has simply connected derived subgroup, but not in general.)
 To summarize, we obtain a map
 \begin{equation}\label{e:Sigma(H)-Sigma(G)}
 	\Sigma(H)_{(G,H)\text{-} \reg} \To \Sigma(G)\cup \{\emptyset\}.
 \end{equation}
 We say that $\gamma_H$ and $\gamma\in G(F)_{\semi}$ have\index{matching conjugacy classes}\index{matching conjugacy classes!image} \emph{matching conjugacy classes}, or simply that $\gamma$ is an \emph{image} of $\gamma_H$. Since the centralizer of $\gamma$ in $G$ is connected (because $G_{\der}=G_{\sconn}$), the centralizer of $\gamma_H$ in $H$ is also connected by \cite[Lem.~3.2]{Kot86}.
 In fact, the notion of matching conjugacy classes as above works as long as  the centralizer of $\gamma$ in $G$ is connected, e.g., if $\gamma$ is strongly regular (without requiring $G_{\der}=G_{\sconn}$).
 
 Now we drop the assumption that $G_{\der}=G_{\sconn}$. Choose a $z$-extension
 $G_1$ and its associated central extension $H_1$ as in \eqref{eq:z-ext-defined} and Lemma \ref{lem:z-ext-H}, and let $Z_1$ denote their common kernel. Then the above construction can be performed for $G_1$ and $H_1$ in place of $G$ and $H$. Then the analogue of \eqref{e:Sigma(H)-Sigma(G)} is equivariant for the translation by elements of $Z_1(F)$; hence it descends along the surjections $G_1\rightarrow G$ and $H_1\rightarrow H$ to yield the following commutative diagram (where the symbol $\emptyset$ maps to itself under the right vertical map). 
 $$
 \xymatrix{
 	\Sigma(H_1)_{(G_1,H_1)\text{-} \reg} \ar[r]\ar[d] & \Sigma(G_1)\cup \{\emptyset\} \ar[d] \\
 	\Sigma(H)_{(G,H)\text{-} \reg} \ar[r] & \Sigma(G)\cup \{\emptyset\}
 }
 $$
 The bottom horizontal map is independent of the choice of $z$-extensions, as can be shown by using \cite[Lem.~2.4.4]{Kot84a} (which produces a $z$-extension mediating two arbitrary $z$-extensions of a group).
 
 By slight abuse of language (as $H_1$ is not an endoscopic group of $G$), $\gamma\in G(F)_{\semi}$ is said to be an image of $\gamma_{H_1}\in H_1(F)_{(G_1,H_1)\text{-} \reg}$ if the stable conjugacy class of $\gamma_{H_1}$ maps to that of $\gamma$ in the above diagram.

 \subsection{}
 
 We introduce $\kappa$-orbital integrals, of which stable orbital integrals are the special case.
 Let us assume that $F$ is local for convenience. The main definitions here extend to the adelic setting in the obvious manner.
 
 Let $\gamma\in G(F)_{\semi}$ and $x\in G(\ol{F})$. Suppose that $\gamma_x:=x^{-1}\gamma x\in G(F)$ and $x \rho(x)^{-1}\in I_\gamma(\ol{F})$ for every $\rho\in \Gamma_F$. Then $x$ and the 1-cocycle $\rho\mapsto x \rho(x)^{-1}$ define an element of $\coh^0(F,I_\gamma\bs G)$, to be denoted by $\dot x$. The map $x\mapsto \gamma_x$ factors through $\coh^0(F,I_\gamma\bs G)$, namely there is an induced map\index[Oct2026-sorted-notation]{zz03 gamma 1 00 00 gamma subscript dot x@$\gamma_{\dot x}$}
 $$\coh^0(F,I_\gamma\bs G)\To G(F),\quad \dot x\longmapsto \gamma_{\dot x}.$$
 Recall that there is a short exact sequence
 $$1\To I_\gamma(F)\bs G(F)\To \coh^0(F,I_\gamma\bs G) \To \fkD(I_\gamma,G;F) \To 1,$$
 coming from a long exact sequence. Given Haar measures on $I_\gamma(F)$ and $G(F)$ and the counting measure on $\fkD(I_\gamma,G;F)$, there is a unique way to equip $\coh^0(F,I_\gamma\bs G)$ with a compatible measure. The map $\dot x\mapsto \gamma_{\dot x}$ induces a map\index[Oct2026-sorted-notation]{zz03 gamma 1 00 00 gamma subscript x leftbracket rightbracket@$\gamma_{[x]}$}
 $$\fkD(I_\gamma,G;F)\To \Gamma(G),\qquad [x]\longmapsto \gamma_{[x]},$$
 whose image consists of conjugacy classes in the stable conjugacy class of $\gamma$. If $G_\gamma=\Cent(\gamma,G)$ is connected (so that it equals $I_\gamma$) then $[x]\mapsto \gamma_{[x]}$ is a bijection onto the image. In general the fiber over the conjugacy class of $\gamma'\in G(F)$ (stably conjugate to $\gamma$) is in bijection with $\ker(\coh^1(F,I_{\gamma'})\ra \coh^1(F,G_{\gamma'}))$.
 Recall the map $\fkD(I_\gamma,G;F)\ra \fkE(I_\gamma,G;F)=\fkK(I_\gamma/F)^D$ from \S\ref{para:two groups}. Thus we have a pairing
 $$\lg \cdot,\cdot \rg: ~\fkD(I_\gamma,G;F)\times \fkK(I_\gamma/F)\To \CC^\times.$$
 Given $\kappa\in \fkK(I_\gamma/F)$ and $f\in \cH(G)$ we define the\index{orbital integral!kappa@$\kappa$-}\index{kappa-orbital integral@$\kappa$-orbital integral|see{orbital integral}} $\kappa$-orbital integral of $\gamma$ by\index[Oct2026-sorted-notation]{o 0 00 00 o superscript g f kappa subscript gamma f@$O^{G(F),\kappa}_\gamma(f)$}
 \begin{eqnarray}
 	O^{G(F),\kappa}_\gamma(f)&:=&
 	\int_{\coh^0(F,I_\gamma\backslash G)} e(I_{\gamma_{\dot x}})\lg \dot x,\kappa\rg f(\gamma_{\dot x}) d\dot x\nonumber\\
 	&=&	\sum_{[x]\in \fkD(I_\gamma,G;F)} e(I_{\gamma_{[x]}})\lg [x],\kappa\rg O^{G(F)}_{\gamma_{[x]}}(f).
 	\nonumber
 \end{eqnarray} Here\index[Oct2026-sorted-notation]{e 1 00 00 e cdot kottwitz sign@$e(\cdot)$ (Kottwitz sign)} $e(\cdot)$ is the\index{Kottwitz sign} Kottwitz sign \cite{Kot83}. 
 When $\kappa$ is trivial, one has the\index{orbital integral!stable}\index{stable orbital integral|see{orbital integral}} stable orbital integral\index[Oct2026-sorted-notation]{so 0 00 00 so superscript g f subscript gamma f@$SO^{G(F)}_\gamma(f)$}
 $$SO^{G(F)}_\gamma(f):=O^{G(F),1}_\gamma(f).$$
 The superscript $G(F)$ will be omitted if there is no danger of confusion. The above definition of $\kappa$-orbital integrals works verbatim for $f\in \cH(G(F),\omega^{-1})$.

 \subsection{}
 \label{subsub:untwisted-transfer}
 We recall the\index{transfer!Langlands--Shelstad} Langlands--Shelstad transfer in the local untwisted case. See \S\ref{subsub:local-twisted-transfer} below for the twisted case.
 
 When $\cH={}^L H$ in the endoscopic datum (this assumption will be removed via $z$-extensions), Langlands--Shelstad \cite{LS87,LS90} define the transfer factor\index[Oct2026-sorted-notation]{zz04 delta 0 00 00 delta cdot cdot@$\Delta(\cdot,\cdot)$}\index{transfer!factor}
 $$\Delta(\cdot,\cdot): H(F)_{(G,H)\text{-} \reg}\times G(F)_{\semi} \To \CC$$
 which vanishes on $(\gamma_{H},\gamma)$ unless $\gamma$ is an image of $\gamma_{H}$. 
 
 When $G$ is quasi-split, the canonical  transfer factor\index[Oct2026-sorted-notation]{zz04 delta 0 00 00 delta subscript 0@$\Delta_0$} $\Delta_0$ was given by Langlands--Shelstad depending only on a choice of $F$-pinning for $G$. Another natural normalization in the quasi-split case is the Whittaker normalization in \cite[\S5.3]{KS99}. While there is no direct analogue of either when $G$ is not quasi-split, the Whittaker normalization can be extended to $G$ given a suitable rigidification of an inner twisting of $G$ against its quasi-split inner form. See \cite{KalLoc,KalGlo} for Kaletha's notion of rigid inner forms and a discussion of other rigidifications. In this work, we do not attempt to choose a rigidification or a canonical normalization of transfer factor at every place. However when $G$ is defined over $\Q$, we may and will always choose a global normalization as in \cite[\S 6.4]{LS87} so that the product formula (Corollary 6.4.B therein) holds true.
 
 In fact there is more than one sign convention for (untwisted and twisted) transfer factors as explained in \cite{KS12}. We have chosen to work with the factor $\Delta'$ in \textit{loc.~cit.}, which we believe is the easiest to compare with the transfer factors in the literature:  $\Delta'$ coincides with the convention in \cite{Kot90} (see p.~178 therein) but differs from that of \cite{LS87,LS90} by the map $(H,\cH,s,\eta)\mapsto(H,\cH,s^{-1},\eta)$.
 A comparison with \cite[App.~A]{Mor10} is detailed in \cite[\S5.6]{KS12}. A downside of $\Delta'$ is that its definition is based on the arithmetic normalization of the Langlands correspondence for tori; see \cite[\S\S4--5]{KS12} and \cite[App.~F]{AGIKMS} for details on the normalization issues. Since our convention for the Langlands correspondence is geometric throughout this work, this results in an extra inverse in the formula for transfer factors under translation by central elements.

 \subsection{}\label{ss:lambdaH}
 
 There exists a smooth character\index[Oct2026-sorted-notation]{zz11 lambda 1 00 00 lambda subscript h@$\lambda_H$} $\lambda_H:Z_G(F)\ra \CC^\times$ such that
 \begin{equation}\label{eq:trans-Z-equiv}
 	\Delta(z\gamma_H,z\gamma) = \lambda_H(z) \Delta(\gamma_H,\gamma),\qquad z\in Z_G(F).
 \end{equation}
 This is \cite[Lem.~3.5.A]{LS90}. We can describe $\lambda_H$ explicitly on $Z_G^0(F)$ as follows. The lemma remains equivalent if we replace $\lambda_H$ with $\lambda^{-1}_H$ and adopt the arithmetic normalization. 
 
 \begin{lem}\label{lem:lambdaH}
 	When $\cH={}^L H$, the restriction of $\lambda^{-1}_H$ to $Z_G^0(F)$ corresponds  (via the geometrically normalized Langlands correspondence for tori) to the composite Langlands parameter 
 	$$W_F \To {}^L H \overset{\eta}{\To} {}^L G \overset{\zeta}{\To} {}^L Z_G^0,$$
 	where the first map is the distinguished splitting, and $\zeta$ is dual to the inclusion $Z_G^0\hra G$.
 \end{lem}
 
 \begin{proof}
 	Consider $\gamma_H\in H(F)$, $\gamma\in G(F)$, and maximal tori $T\subset G$ and $T_H\subset H$ as in \cite[\S3]{LS87}. In particular we are given an isomorphism $i:T\cong T_H$, inducing an isomorphism $^L i: {}^L T_H\cong {}^L T$.
 	They construct $L$-morphisms $\xi_{T_H}: {}^L T_H \ra {}^L H$ and $\xi_T: {}^L T \ra {}^L G$ (depending on some additional choices) as well as a 1-cocycle $\mathbf{a}:W_F\ra \hat T_H$ such that $\eta \xi_{T_H}=\mathbf{a}\cdot \xi_T {}^L i$ as $L$-morphisms from $^L T_H$ to $^L G$. Restricting the equality via the splitting $W_F\ra {}^L T_H$, we obtain 
 	\begin{equation}\label{eq:trans-factor-equiv}
 		\eta(\xi_{T_H}(w)) = \mathbf{a}(w) \xi_T(w) \in {}^L G,\qquad w\in W_F.
 	\end{equation}
 	The first paragraph of \cite[p.~253]{LS87} tells us that $\lambda_H|_{Z_G^0(F)}$ corresponds (under the arithmetic normalization) to the composition of $\mathbf{a}$ with $\hat T_H \cong\hat T \rightarrow \hat Z_G^0$, which is dual to $Z_G^0\subset T \cong T_H$. To prove the lemma, it is thus enough to verify that the composite map in the lemma is $w\mapsto \zeta(\mathbf{a}(w))\rtimes w$.
 	
 	To this end, write $\xi_{T_H}(w)=b(w)\rtimes w$ and $\xi_T(w)=c(w)\rtimes w$. From the construction of $\xi_T$ in \cite[\S 2.6]{LS87}, it follows that $\zeta(c(w))=1$.  Indeed, the two main points are that the image of every morphism $\SL_2\ra \hat G$ maps trivially to $\hat Z_G^0$ (thus also $n(\omega_T(\sigma))$ therein) and that the coroots of $\hat G$ map trivially to $\hat Z_G^0$ (thus also $r_p(w)$ therein). Similarly, $b(w)$ maps to $1\in \hat Z_H^0$, so  $\zeta(b(w))=1$. Now we apply $\zeta$ to \eqref{eq:trans-factor-equiv} to see that
 	$\zeta(\eta(w))=\zeta(\mathbf{a}(w))\rtimes w$, as desired.
 \end{proof}
 
 \subsection{}\label{para:trans-factor-general}
 We introduce transfer factors in general by reducing to the case $G_{\der}=G_{\textup{sc}}$, in which case we can always assume that $\cH={}^L H$, cf.~\S\ref{para:z-ext} and \S\ref{subsub:untwisted-transfer}.
 Let $\fke=(H,\cH,s,\eta)\in E(G)$. Take a $z$-extension $1\ra Z_1\ra G_1\ra G\ra 1$ and define $\fke_1=(H_1,{}^L H_1,s_1,\eta_1)\in E(G_1)$ as in \S\ref{G1H1}.
 Then we have $\Delta(\gamma_{H_1},\gamma_1)\in \CC$ defined on $\gamma_{H_1}\in H_1(F)_{(G_1,H_1)\text{-} \reg}$ and $\gamma_1\in G_1(F)_{\semi}$. Langlands--Shelstad define\index[Oct2026-sorted-notation]{zz04 delta 0 00 00 delta cdot cdot@$\Delta(\cdot,\cdot)$}
 $$\Delta(\cdot,\cdot): H_1(F)_{(G_1,H_1)\text{-} \reg}\times G(F)_{\semi} \To \CC$$
 as follows. Set $\Delta(\gamma_{H_1},\gamma)=0$ unless $\gamma$ is an image of $\gamma_{H_1}$ (\S\ref{subs:trans-conjugacy}), that is, unless $\gamma$ lifts to $\gamma_1\in G_1(F)$ which is an image of $\gamma_{H_1}$, in which case $\Delta(\gamma_{H_1},\gamma):=\Delta(\gamma_{H_1},\gamma_1)$. By \eqref{eq:trans-Z-equiv}, 
 \begin{equation}\label{eq:central-equiv-H1} \Delta(z\gamma_{H_1},\ol{z}\gamma) = \lambda_{H_1}(z) \Delta(\gamma_{H_1},\gamma),\qquad z\in Z_{G_1}(F),\end{equation}
 where $\ol{z}$ denotes the image of $z$ in $Z_G(F)$.
 Notice that $Z_1\subset Z_{G_1}^0$. Lemma \ref{lem:lambdaH} implies that $\lambda_{H_1}|_{Z_1(F)}=\lambda_1$, where $\lambda_1$ was given in \S\ref{para:lambda1} as explained in \cite[p.23]{KS99} and \cite[Rem.~5.3.1]{KS12}. (The references tell us that the two characters are mutual inverses, but note that our sign convention for $\lambda_{H_1}$ is flipped, cf.~Lemma \ref{lem:lambdaH}.)

 \subsection{}\label{subsub:untwisted-unram-Hecke-transfer}
 
 Continue to assume that $F$ is a local field. We state the\index{transfer!Langlands--Shelstad} Langlands--Shelstad transfer and the\index{fundamental lemma} fundamental lemma for reductive  groups $G_1$ with $G_{1,\der}=G_{1,\textup{sc}}$. Let $\fke_1=(H_1,{}^L H_1,s_1,\eta_1)\in E(G_1)$.
 If $G_1$ and $\fke_1$ are unramified, then $\eta_1$ induces a $\CC$-algebra map via the Satake transform:
 \begin{equation}\label{eq:eta1}
 	\eta_1^*: \cH^{\ur}(G_1)\To \cH^{\ur}(H_1).
 \end{equation}

 \begin{prop}\label{prop:LS-transfer}
 	For each $f_1\in \cH(G_1)$, there exists $f_1^{H_1}\in \cH(H_1)$ enjoying the following property: If $\gamma_{H_1}\in H_1(F)_{(G_1,H_1)\text{-} \reg}$ has no image in $G_1(F)_{\semi}$ then $SO_{\gamma_{H_1}}(f_1^{H_1})=0$. If $\gamma_{H_1}$ has an image $\gamma_1\in G_1(F)_{\semi}$ then
 	\begin{equation}\label{eq:FL}
 		SO_{\gamma_{H_1}}(f^{H_1}_1) =
 		\sum_{[x]\in \fkD(I_{\gamma_1},G_1;F)} e(I_{\gamma_{1,[x]}}) \Delta(\gamma_{H_1},\gamma_{1,[x]}) O_{\gamma_{1,[x]}}(f_1) .
 	\end{equation}
 	Moreover, the fundamental lemma (FL) holds true, i.e., if $G_1$ and $\fke_1$ are unramified and if $f_1\in \cH^{\ur}(G_1)$ then the above holds with $f^{H_1}_1=\eta_1^* f_1$.
 	
 \end{prop}
 
 \begin{proof}
 	The second assertion (FL) follows from work of Ng\^o as well as Cluckers-Loeser, Hales, and Waldspurger  \cite{CL10,Hal95,Ngo10,Wal06}. The first assertion (the transfer conjecture) is implied by FL \cite{Wal97} if $F$ is non-archimedean; see \cite{Shelstad2012} in the archimedean case. 
 \end{proof}
 
 \subsection{}\label{para:transfer-with-central-character}
 Let us adapt Proposition \ref{prop:LS-transfer} to the setting with fixed central characters.
 Let $\fkX$ be a closed subgroup of $Z(F)$ and $\chi:\fkX\ra \CC^\times$ a continuous character. We view $\fkX$ also as a closed subgroup of $Z_H(F)$ via $Z\hra Z_H$. Denote by $\fkX_{H_1}$ and $\fkX_1$ the preimages of $\fkX$ under $H_1(F)\rightarrow H(F)$ and $G_1(F)\rightarrow G(F)$, respectively, so that $\fkX_{H_1}\cong \fkX_1$ canonically. Choose compatible Haar measures on $\fkX_{H_1}$ and $\fkX_1$. Let $\chi_1$ denote the character of $\fkX_{H_1}$ or $\fkX_{1}$ pulled back from $\chi:\fkX\ra \CC^\times$.  Restricting $\lambda_{H_1}$ as in \S\ref{para:trans-factor-general}, we obtain another character
 $$\lambda_{H_1}|_{\fkX_{H_1}}:\fkX_{H_1}\To\CC^\times.$$
 By slight abuse of notation, we will often denote the above still by $\lambda_{H_1}$ (as we will not consider $\lambda_{H_1}$ on a larger domain). In the special case when $\fkX=\{1\}$, notice that $\fkX_1=\fkX_{H_1}=Z_1(F)$, $\chi_1=1$, and $\lambda_{H_1}=\lambda_1$.
 
 If $F$ is non-archimedean, $\chi:\fkX\ra \CC^\times$ is said to be \emph{unramified} if the character is trivial on the maximal compact subgroup of $\fkX$. The same definition works with $\fkX_{H_1}$ in place of $\fkX$. When $G$, $\fke$, $\chi$, and the $z$-extension $G_1$ are unramified, we have $G_1$ and $\fke_1$ also unramified. In this case, the map $\eta_1^*:\cH^{\ur}(G_1)\ra \cH^{\ur}(H_1)$ from \eqref{eq:eta1} is averaged to give a map
 \begin{equation}\label{eq:eta1-variant}
 	\cH^{\ur}(G,\chi^{-1})\To \cH^{\ur}(H_1,\chi_1^{-1}\lambda_{H_1})
 \end{equation}
 as follows. To get started, we claim that $\lambda_{H_1}|_{\fkX_{H_1}}$ is unramified, i.e., $\lambda_{H_1}(z)=1$ for every $z$ in the maximal compact subgroup of $\fkX_{H_1}$. For every $G_1$-regular element $\gamma_{H_1}$, the central equivariance of transfer factors implies that \[\lambda_{H_1}(z)SO_{\gamma_{H_1}}(\eta_1^* f_1)=SO_{z\gamma_{H_1}}(\eta_1^* f_1)=SO_{\gamma_{H_1}}(\eta_1^* f_1),\qquad f_1\in \cH^{\ur}(G_1).\]
 Since $SO_{\gamma_{H_1}}(\eta_1^* f_1)\neq 0$ if $f_1$ is the unit element and $\gamma_{H_1}$ lies in the fixed hyperspecial subgroup of $H_1(F)$, we see that  $\lambda_{H_1}(z)=1$. The claim is proved.
 
 Now let $f_1\in \cH^{\ur}(G_1)$ be any lift of $f$ along the natural surjective map, where the first arrow is averaging against $\chi_1$:
 $$\cH^{\ur}(G_1)\To\cH^{\ur}(G_1,\chi_1^{-1}) \cong \cH^{\ur}(G,\chi^{-1}).$$
 For $z\in \fkX_{H_1}$, write $f_{1,z}(h):=f_1(zh)$. Define a function $f^{H_1}$ on $H_1(F)$ by $$f^{H_1}(h):=\int_{\fkX_{H_1}} \chi_{1}(z)\eta_1^*f _{1,z}(h) dz=\int_{\fkX_{H_1}} \chi_1(z)\lambda^{-1}_{H_1}(z) \eta_1^*f _{1}(zh) dz.$$
 Then $f^{H_1}$ belongs to  $\cH(H_1,\chi_1^{-1}\lambda_{H_1})$ and is independent of the choice of $f_1$. The resulting map \eqref{eq:eta1-variant} is again denoted by $\eta_1^*$ as there is little danger of confusion.

 \begin{prop}\label{prop:untwisted-transfer}
 	For each $\chi:\fkX\ra \CC^\times$ and $f\in \cH(G,\chi^{-1})$, there exists
 	$$f^{H_1}\in \cH(H_1,\chi_1^{-1}\lambda_{H_1})$$
 	such that for every $\gamma_{H_1}\in H_1(F)_{(G_1,H_1)\text{-} \reg}$, if $\gamma\in G(F)_{\semi}$ is an image of $\gamma_{H_1}$ then
 	$$ SO^{H_1(F)}_{\gamma_{H_1}}(f^{H_1})=\sum_{[x]\in \fkD(I_\gamma,G;F)} e(I_{\gamma_{[x]}}) \Delta(\gamma_{H_1},\gamma_{[x]}) O^{G(F)}_{\gamma_{[x]}}(f).$$
 	If $\gamma_{H_1}$ admits no image in $G(F)$ then $ SO^{H_1(F)}_{\gamma_{H_1}}(f^{H_1})=0$.

 	Moreover when $F$ is non-archimedean, if $G_1$, $\fke$, and $\chi$ are unramified, and if $f\in \cH^{\ur}(G,\chi^{-1})$, then $\chi_1,\lambda_{H_1}|_{\fkX_{H_1}}$ are unramified and the above holds true with 
 	$f^{H_1}=\eta_1^*(f)$.
 \end{prop}
  
\begin{proof}
Given $f\in \cH(G,\chi^{-1})$, choose a lift $f_1\in\cH(G_1)$ under the surjective composite map
$\cH(G_1) \rightarrow \cH(G_1,\chi_1^{-1})=\cH(G,\chi^{-1})$.
Let $ f_1^{H_1}\in \cH(H_1)$ be a transfer of $f_1$ as in Proposition \ref{prop:LS-transfer}. 	 Define $f^{H_1}\in \cH(H_1,\chi_1^{-1}\lambda_{H_1})$ by
$$f^{H_1}(\gamma_{H_1}):= \int_{\fkX_{H_1}} f^{H_1}_1 (z\gamma_{H_1}) \chi_1(z)\lambda_{H_1}^{-1}(z)dz,\quad \gamma_{H_1}\in H_1(F),$$
so that for $\gamma_{H_1}\in H_1(F)_{\semi}$,
\begin{equation}\label{eq:SO-average}
	SO_{\gamma_{H_1}}(f^{H_1})= \int_{\fkX_{H_1}} \chi_1(z)\lambda_{H_1}^{-1}(z) SO_{z \gamma_{H_1}}( f^{H_1}_1 ) dz.
\end{equation}
If $\gamma_{H_1}$ has no image in $G(F)$ then it has no image in $G_1(F)$ either, so $SO_{\gamma_{H_1}}(f^{H_1})$ vanishes.
Otherwise, let $\gamma_1\in G_1(F)$ be an image of  $\gamma_{H_1}$. Then $z\gamma_1$ is an image of $z\gamma_{H_1}$ for each $z\in \fkX_{H_1}$. Applying Proposition \ref{prop:LS-transfer} to \eqref{eq:SO-average}, we see that $SO_{\gamma_{H_1}}(f^{H_1})$ equals
$$\int_{\fkX_{H_1}} \chi_1(z) \sum_{[x]\in \fkD(I_{\gamma_1},G_1;F)} e(I_{\gamma_{1,[x]}}) \frac{\Delta(z\gamma_{H_1},z\gamma_{1,[x]})}{\lambda_{H_1}(z)} O_{z\gamma_{1,[x]}}(f_1) dz.$$
Lemma \ref{lem:z-ext-wrt-K-and-D} gives a bijection $\fkD(I_{z\gamma_{1}},G_1;F)=\fkD(I_{\gamma_{1}},G_1;F)\cong \fkD(I_{\gamma},G;F)$.
(We view $[x]$ also as an element of $\fkD(I_{z\gamma_{1}},G_1;F)$ or $\fkD(I_{\gamma},G;F)$.)
By definition and \eqref{eq:trans-Z-equiv}, we have $$\Delta(z\gamma_{H_1},z\gamma_{1,[x]})=\Delta(z\gamma_{H_1},\ol{z}\gamma_{[x]})= \lambda_{H_1}(z) \Delta(\gamma_{H_1},\gamma_{[x]}),\qquad z\in \fkX_{H_1},$$
where $\ol{z}\in \fkX$ is the image of $z$.
By \cite[Cor.~(2)]{Kot83} $e(I_{\gamma_{1,[x]}})=e(I_{\gamma_{[x]}})$. Altogether, with both sums running over $[x]\in \fkD(I_\gamma,G;F)$ below,
\begin{eqnarray}
	SO_{\gamma_{H_1}}(f^{H_1})	&=& \sum_{[x]} e(I_{\gamma_{[x]}}) \Delta(\gamma_{H_1},\gamma_{[x]})
	\int_{\fkX_1} \chi_1(z)O_{z\gamma_{1,[x]}}(f_1)dz\nonumber
	\\ &=& \sum_{[x]} e(I_{\gamma_{[x]}}) \Delta(\gamma_{H_1},\gamma_{[x]}) O_{\gamma_{[x]}}(f). \nonumber
\end{eqnarray}

It remains to prove the last assertion when $G$, $\fke_1$, and $\chi$ are unramified. The unramifiedness of $\lambda_{H_1}$ is proved as in \S\ref{para:transfer-with-central-character}. The fact that we can take $f^{H_1}=\eta_1^*(f)$ follows from the earlier part of the current proof, where we can pick $f_1\in \cH^{\ur}(G_1)$ and choose $f_1^{H_1}$ to be the image of $f_1$ under \eqref{eq:eta1}.
\end{proof}

\subsection{}
\label{subsub:local-twisted-transfer} 
Here we work out a small generalization of some results in local twisted endoscopy by Morel and Kottwitz \cite[\S9, App.~ A]{Mor10} to the setting where $\cH$ in the endoscopic datum cannot be taken to be an $L$-group.
We put ourselves in the setting of \S\ref{subsub:twisted-end-at-p} with $F=\Q_p$.
Let $\fke=(H,\cH,s,\eta)\in E(G)$, which gives rise to $\fke_1=(H_1,{}^L H_1,s_1,\eta_1)\in E(G_1)$ as in Lemma \ref{lem:endoscopic-datum-z-ext}. We make an additional hypothesis that 
$$s\in Z(\hat H)^{\Gamma_p}.$$
The group $G$ is assumed quasi-split over $\Q_p$ so that \cite[App.~ A]{Mor10} applies.
If $G$ and $H$ are unramified (which we are not assuming) then we may and will choose $G_1$ and $H_1$ to be also unramified.

In twisted endoscopy, the norm map is defined by\index{transfer!Kottwitz--Shelstad} Kottwitz and Shelstad \cite{KS99} for strongly regular elements and by Labesse \cite{Lab04} for elliptic elements (which may not be strongly regular). The norm map in untwisted endoscopy is simply the transfer of conjugacy classes as in \S\ref{subs:trans-conjugacy}. In the special case of base change, Kottwitz \cite{Kot82} defines the norm map for general elements. (These norm maps  coincide when more than one definition is available in a given setting.) For our purpose, we define the norm from $R_1(\Q_p)$ to $H_1(\Q_p)$ as follows: $\gamma_{H_1}$ is said to be a norm of $\delta_1\in R_1(\Q_p)$ if an image of $\gamma_{H_1}$ in $G_1(\Q_p)$ is a degree $n$ base-change norm of $\delta_1$. It is an exercise to check that this is consistent with the norm map by Kottwitz--Shelstad and Labesse.

Let $\delta_1\in  G_1(\Q_{p^n})=R_1(\Q_p)$. It has a degree $n$ norm $\gamma_{0,1}\in G_1(\Q_p)$. We assume that $\gamma_{0,1}$ is semi-simple. Then $\delta_1$ corresponds to $[b_1]\in \fkD_n(\gamma_{0,1},G_1;\Q_p)$ according to Corollary \ref{cor:bijection}. So we write $\delta_{[b_1]}$ for $\delta_1=\delta_{[b_1]}$. Recall $\beta_p(\gamma_{0,1},[b_1])=\kappa_{I_{0,1}}([b_1])$. 
Suppose that $\gamma_{H_1}\in H_1(\Q_p)_{\semi}$ is a norm of $\delta_1$. Then $\gamma_{0,1}$ is an image of $\gamma_{H_1}$. 

When $\gamma_{H_1}$ is strongly $G_1$-regular (equivalently when $\gamma_{0,1}$ is strongly regular), Kottwitz \cite[Cor.~A.2.10]{Mor10} proved that\index[Oct2026-sorted-notation]{zz04 delta 0 00 00 delta subscript 0@$\Delta_0$}
\begin{equation}\label{e:tw-trans-fac-at-p}
\Delta_0(\gamma_{H_1},\delta_{[b_1]})=\Delta_0(\gamma_{H_1},\gamma_{0,1}) \lg \beta_p(\gamma_{0,1},[b_1]),s_1\rg,
\end{equation}
where we have taken the sign correction of \cite[\S5.6]{KS12} into account.
To make sense of the pairing, we view $s_1\in Z(\hat H_1)^{\Gamma_p}$ as an element of $Z(\hat I_{0,1})^{\Gamma_p}$ via $Z(\hat H_1)^{\Gamma_p}\subset Z(\hat I_{\gamma_{H_1}})^{\Gamma_p} = Z(\hat I_{0,1})^{\Gamma_p}$, cf. \cite[A.3.11]{Mor10}.

When $\gamma_{H_1}$ is $(G_1,H_1)$-regular but not strongly $G_1$-regular, we take \eqref{e:tw-trans-fac-at-p} as the definition of $\Delta_0(\gamma_{H_1},\delta_{[b_1]})$, cf.~\cite[(A.3.11.1)]{Mor10}.

Let $\gamma_0\in G(\Q_p)$ and $[b]\in \fkD_n(\gamma_{0},G;\Q_p)$ be the images of $\gamma_{0,1}$ and $[b_1]$, respectively.
Then $[b]$ gives rise to $\delta_{[b]}\in G(\Q_{p^n})$ via Corollary \ref{cor:bijection}, and $\gamma_0$ is a degree $n$ norm of $\delta_{[b]}$. If the hypothesis on $Z_1$ in Lemma \ref{lem:Kottwitz-triples-z-ext-1} is satisfied, and if $[b_1]\in\fkD_n(\gamma_{0,1},G_1;\Q_p)$ and $[b]\in \fkD_n(\gamma_{0},G;\Q_p)$ correspond under the bijection of that lemma, then the $\sigma$-conjugacy class of $\delta_{[b_1]}$ maps to that of $\delta_{[b]}$.
We defined $\beta_p(\gamma_{0},[b])$ to be $\kappa_{I_{0}}([b])$, which equals the image of $\beta_p(\gamma_{0,1},[b_1])$. Thus we have
$$\lg \beta_p(\gamma_{0,1},[b_1]),s_1\rg=\lg \beta_p(\gamma_{0},[b]),s\rg.$$

In case $G_1$ and $\fke$ are unramified, the twisted datum $\tilde\fke_1=(H_1,\cH_1,\tilde s_1, \tilde \eta_1)$ of \S\ref{para:twisted-datum} is also unramified, and $\tilde\eta_1$
induces $\CC$-algebra morphisms
$$\cH^{\ur}(R_1)\ra \cH^{\ur}(H_1)\quad\mbox{and}\quad\cH^{\ur}(R)\ra \cH^{\ur}(H_1,\lambda_{H_1})$$
as in the untwisted case, cf.~\S\ref{subsub:untwisted-unram-Hecke-transfer} and \S\ref{para:transfer-with-central-character}. By slight abuse of notation we call both maps $\tilde\eta_1^*$, so that if the averaging map 
$\cH^{\ur}(R_1)\ra \cH^{\ur}(R)$ sends $f_1\mapsto f$, then $\tilde\eta_1^*(f)$ is the average of $\tilde\eta_1^*(f_1)$ against $\lambda_{H_1}^{-1}$:
\begin{equation}\label{eq:two-tilde-eta1}
\tilde\eta_1^*(f)(\gamma_{H_1})=\int_{Z_1(\Q_p)} \lambda^{-1}_{H_1}(z) \tilde\eta_1^*(f_1)(z\gamma_{H_1}) dz,\quad \gamma_{H_1}\in H_1(\Q_p).
\end{equation}

The following is a twisted analogue of Proposition \ref{prop:untwisted-transfer}.

\begin{prop}\label{prop:transfer-at-p}
Let $f\in \cH(R(\Q_p))$. Then there exists $f^{H_1}\in \cH(H_1,\lambda_{H_1})$ such that the following holds:
Let $\gamma_{H_1}\in H_1(\Q_p)_{(G_1,H_1)\textup{-reg}}$. We have $$SO_{\gamma_{H_1}}(f^{H_1})=0$$ if $z\gamma_{H_1}$ is not a norm from $R_1(\Q_p)$ for any $z\in Z_1(\Q_p)$. If $z\gamma_{H_1}$ is a norm of $\delta_1\in R_1(\Q_p)=G_1(\Q_{p^n})$ for some $z\in Z_1(\Q_p)$ then whenever $\gamma_{0}\in G(\Q_p)$ is an image of $\gamma_{H_1}$, we have
\begin{equation}\label{e:tw-orb-id}
	SO_{\gamma_{H_1}}(f^{H_1})=\Delta_0(\gamma_{H_1},\gamma_0) \sum_{[b]\in \fkD_n(\gamma_0,G;\Q_p)} e(I_{\delta_{[b]}})  \lg \beta_p(\gamma_{0},[b]),s\rg TO_{\delta_{[b]}}(f).
\end{equation}

Moreover if $G_1$ and $\fke$ are unramified, and if $f\in \cH^{\ur}(R)$ then the above is true for $f^{H_1}=\tilde\eta_1^* (f)\in \cH^{\ur}(H_1,\lambda_{H_1})$.		

\end{prop}

\begin{rem}\label{rem:TFL}
In the essential case when $G_1$ and $\fke$ are unramified with $f$ in the unramified Hecke algebra, if $\gamma_{H_1}$ is assumed to be strongly $G_1$-regular, the proposition is a special case of the\index{fundamental lemma!twisted} twisted fundamental lemma (TFL) for the full unramified Hecke algebra. 
When the residue characteristic $p$ is large, TFL for the unit element is true thanks to Ng\^o, Waldspurger, Cluckers-Loeser, and others (\cite{Ngo10,Wal08} with \cite{Wal06} or \cite{CL10}). 
Together with \cite{Wal97,Wal08}, these results imply the existence of twisted endoscopic transfer for arbitrary test functions, with orbital integral identities on the strongly regular locus, for every residue characteristic and without an unramifiedness assumption.

Returning to the unramified setting and TFL, for general elements of unramified Hecke algebras, including the case of small $p$, TFL was recently established by Lemaire, Moeglin, and Waldspurger \cite{LMW-TFL} and \cite{LW-TFL}. Note that TFL states an orbital integral identity only for strongly regular elements; in fact the twisted transfer factors are not defined for other elements in general. However we want an orbital integral identity for $(G_1,H_1)$-regular semi-simple elements which may not be strongly regular. The transfer factors for such elements are available in our setting by \cite[App.~ A]{Mor10}, which also explains how to obtain Proposition \ref{prop:transfer-at-p} once we know it for strongly regular elements, under the hypotheses that $G_{\der}=G_{\sconn}$ and that $\cH={}^L H$ in the endoscopic datum. 
So the point of our proof below is to remove the hypotheses of \textit{loc.~cit}.
The basic idea is to employ $z$-extensions, but there is a technical problem: no single $z$-extension $G_1$ satisfies the assumption of Lemma \ref{lem:Kottwitz-triples-z-ext-1} for all (sufficiently large) $n$. We choose an auxiliary $z$-extension as in Lemma \ref{lem:z-ext2} to get around the issue.
\end{rem}
 
\begin{proof}[Proof of Proposition \ref{prop:transfer-at-p}]
We present a proof when $G_1$ and $\fke$ are unramified assuming that $f\in \cH^{\ur}(R)$. 
The general case will be taken care of in the last paragraph of this proof.

Put $\tilde Z_1:=\Res_{\Q_{p^n}/\Q_p} (Z_1)_{\Q_{p^n}}$. Choose $f_1$ to be any preimage of $f$ under the averaging map $\cH^{\ur}(R_1)\rightarrow \cH^{\ur}(R)$ over $\ker(R_1(\Q_p)\ra R(\Q_p))=\tilde Z_1(\Q_p)$. Suppose that $z\gamma_{H_1}$ is not a norm from $R_1(\Q_p)$ for any $z\in Z_1(\Q_p)$. Then TFL for general elements of unramified Hecke algebras (together with its extension to $(G_1,H_1)$-regular elements by \cite[Prop.~A.3.14]{Mor10}; see Remark \ref{rem:TFL}) tells us that $SO_{z\gamma_{H_1}}(\tilde\eta_1^*(f_1))=0$. Hence
$$SO_{\gamma_{H_1}}(\tilde\eta_1^*(f))=\int_{Z_1(\Q_p)} \lambda^{-1}_{H_1}(z) SO_{z\gamma_{H_1}}(\tilde\eta_1^*(f_1))dz=0.$$

Now assume the existence of $z\in Z_1(\Q_p)$ such that $\gamma'_{H_1}:=z\gamma_{H_1}\in H_1(\Q_p)$ is a norm of some element $\delta'_1\in R_1(\Q_p)$.
Let $\gamma_{0,1}\in G_1(\Q_p)$ be an image of $\gamma'_{H_1}$. (Such a $\gamma_{0,1}$ exists since $G_1$ is quasi-split and has simply connected derived subgroup, cf.~\cite[Thm.~4.4]{Kot82}.) Write $\gamma_H\in H(\Q_p)$ and $\gamma_0\in G(\Q_p)$ for the respective images of $\gamma_{H_1}$ and $\gamma_{0,1}$. Take central extensions $G'_1$ and $H'_1$ as in Lemma \ref{lem:z-ext2} and write $Z'_1$ for the common kernel of the surjections $G'_1\rightarrow G$ and $H'_1\rightarrow H$. Thus $Z'_1\cong \prod_{i=1}^r \Res_{F_i/\Q_p}\GG_m$ for some $r\in \Z_{\ge 1}$ and $F_i\supset \Q_{p^n}$.
Set $$\tilde Z'_1:=\Res_{\Q_{p^n}/\Q_p}(Z'_1)_{\Q_{p^n}}\cong \prod_{i=1}^r \prod_{\Q_{p^n}\hra F_i} \Res_{F_i/\Q_p}\GG_m,$$
where the second product runs over the set of $\Q_p$-embeddings. Fix an embedding $\Q_{p^n}\hra F_i$ for each $i$ and define $Y_{1}$ to be the subtorus of $\tilde Z'_1$ whose components outside the set of fixed embeddings are trivial. The Frobenius automorphism $\sigma\in \Gal(\Q_{p^n}/\Q_p)$ acts on $\tilde Z'_1$ by permuting the $\Q_p$-embeddings. 

The norm map $N:\tilde Z'_1(\Q_p)\ra Z'_1(\Q_p)$ is obviously onto and restricts to an isomorphism
$$N:Y_1(\Q_p)\isom Z'_1(\Q_p).$$
From $\fke$, we have $\fke'_1=(H'_1,{}^L H'_1,s'_1,\eta'_1)\in E(G'_1)$ as in Lemma \ref{lem:endoscopic-datum-z-ext}. Since $s\in Z(\hat H)^\Gamma$ we have that $s'_1\in Z(\hat H'_1)^\Gamma$. Consider $R'_1:=\Res_{\Q_{p^n}/\Q_p} (G'_1)_{\Q_{p^n}}$, equipped with a $\Q_p$-automorphism $\theta'_1$ induced by $\sigma$. One builds a twisted endoscopic datum $\tilde\fke'_1$ for $(R'_1,\theta'_1)$ as in \S\ref{subsub:twisted-end-at-p} and \S\ref{para:twisted-datum}.

Choose $f'_1\in \cH^{\ur}(R'_1)$ to be a preimage of $f$ under the averaging surjection $ \cH^{\ur}(R'_1)\rightarrow \cH^{\ur}(R)$.
We have two functions:
$$(\tilde\eta'_1)^*(f'_1)\in \cH^{\ur}(H'_1)\quad\mbox{and}\quad f^{H'_1}:=(\tilde\eta'_1)^*(f)\in \cH^{\ur}(H'_1,\lambda_{H'_1}),$$
following the paragraph above \eqref{eq:two-tilde-eta1}.
The inclusion $H_1\hra H'_1$ induces an isomorphism $\cH^{\ur}(H_1,\lambda_{H_1})\cong \cH^{\ur}(H'_1,\lambda_{H'_1})$
since the character $\lambda_{H'_1}$ restricts to $\lambda_{H_1}$.
The functions $f^{H'_1}$ and $f^{H_1}$ correspond under the isomorphism. Clearly
$$SO_{\gamma'_{H_1}}(f^{H'_1})= SO_{\gamma'_{H_1}}(f^{H_1})=\lambda_{H_1}(z)SO_{\gamma_{H_1}}(f^{H_1}),$$
where $\gamma'_{H_1}$ is considered as an element of either $H_1(\Q_p)$ or $H'_1(\Q_p)$ via $H_1\subset H'_1$.
By TFL for unramified Hecke algebras and its extension (Remark \ref{rem:TFL} again),
\begin{equation}\label{eq:Mor10-951-A314}
SO_{\gamma'_{H_1}}((\tilde\eta'_1)^*(f'_1))=\sum_{\delta'_1} e(I_{\delta'_1}) \Delta_0(\gamma'_{H_1},\delta'_1) TO_{\delta'_1}(f'_1),
\end{equation}
where the sum runs over a set of representatives for the $\sigma$-conjugacy classes in $G'_1(\Q_{p^n})=R'_1(\Q_p)$ whose norm is $\gamma_{0,1}$. (Again $\gamma'_{H_1}\in H_1(\Q_p)$ is viewed as an element of $H'_1(\Q_p)$ via $H_1\subset H'_1$. Similarly $\gamma_{0,1}\in G_1(\Q_p)$ is viewed as an element of $G'_1(\Q_p)$.)
Let $\delta\in R(\Q_p)$ denote the image of $\delta'_1$, and $y\in \tilde Z'_1(\Q_p)$ an arbitrary element. We collect the following facts.
\begin{itemize}
\item the set of representatives $\delta'_1$ in the sum is in bijection with $\fkD_n(\gamma_{0,1},G'_1;\Q_p)$, and also with $\fkD_n(\gamma_0,G;\Q_p)$ (Lemma \ref{lem:Kottwitz-triples-z-ext-1} and Corollary \ref{cor:bijection}),
\item $e(I_{y\delta'_1})=e(I_\delta)$ by \cite[Cor.~(2)]{Kot83},
\item $N(y)\gamma_{0,1}$ (resp.~$N(y)\gamma'_{H_1}$) is a norm of $y\delta'_1$ in $G'_1$ (resp. $H'_1$),
\item \eqref{eq:Mor10-951-A314} holds with $N(y)\gamma_{0,1}$ and $N(y)\gamma'_{H_1}$ in place of $\gamma_{0,1}$ and $\gamma'_{H_1}$,
\item $\Delta_0(N(y)\gamma'_{H_1},y\delta'_1)=\lambda_{H'_1}(N(y))\Delta_0(\gamma'_{H_1},\delta'_1)$ by \cite[Lem.~5.1.C and p.~53]{KS99},\footnote{Unlike the formula in \emph{loc.~cit.}, we do not put an inverse over $\lambda_{H'_1}$. Their character $\lambda^{-1}_{H_1}$ equals our $\lambda_1$ (since their $\lambda_{H_1}$ and our $\lambda_1$ correspond to the same $L$-parameter but via opposite normalizations), and our $\lambda_{H'_1}$ equals $\lambda_1$ on $Z_1(F)$ by \S\ref{para:trans-factor-general}.}
\item  $\Delta_0(\gamma'_{H_1},\gamma_{0,1})$ remains the same whether it is viewed with respect to the transfer between $H'_1$ and $G'_1$ or between $H_1$ and $G_1$.\footnote{To see this, one reduces to the strongly regular case by \cite[\S 4.3]{LS87} and \cite[\S 2.4]{LS90}. The same $a$-data and $\chi$-data may be chosen in the two cases to compare transfer factors, as the centralizer of $\gamma_{0,1}$ determines the same root system in $G_1$ and $G'_1$ respectively. Then one sees that each of $\Delta_{I}$, $\Delta_{II}$, $\Delta_{III_1}$, $\Delta_{III_2}$, and $\Delta_{IV}$ is the same by inspecting the definition in \cite[\S 3]{LS87}.}
\end{itemize}
We put all this together below, where the first equality uses \eqref{eq:two-tilde-eta1} with $Z'_1,H'_1$ in place of $Z_1,H_1$; the second last equality holds because unfolding the central fiber
$\tilde Z'_1(\Q_p)/Z'_1(\Q_p)$ via
$t\mapsto t^{-1}\theta'_1(t)$ supplies integration over $\ker N$, which combines with $\int_{Y_1(\Q_p)}$ to give the full central average defining $f$ by Fubini.
\begin{eqnarray*}
SO_{\gamma'_{H_1}}(f^{H'_1})&=&\int_{Z'_1(\Q_p)} \lambda^{-1}_{H'_1}(z') SO_{z'\gamma'_{H_1}}((\tilde\eta'_1)^*(f'_1))dz'\nonumber\\
&=&\int_{Y_1(\Q_p)}\lambda^{-1}_{H'_1}(N(y)) \sum_{\delta'_1} e(I_{y\delta'_1}) \Delta_0(N(y)\gamma'_{H_1},y\delta'_1) TO_{y\delta'_1}(f'_1)dy.\nonumber\\
&=&\sum_{\delta'_1} e(I_{\delta'_1}) \Delta_0(\gamma'_{H_1},\delta'_1) \int_{Y_1(\Q_p)}TO_{y\delta'_1}(f'_1)dy.\nonumber\\
&=&\sum_{\substack{[b'_1]\in \\ \fkD_n(\gamma_{0,1},G'_1;\Q_p)}} e(I_{\delta_{[b'_1]}}) \Delta_0(\gamma'_{H_1},\gamma_{0,1}) \lg \beta_p(\gamma_{0,1},[b'_1]),s'_1\rg TO_{\delta_{[b]}}(f).\nonumber\\
&=&\sum_{[b]\in \fkD_n(\gamma_0,G;\Q_p)} e(I_{\delta_{[b]}})\lambda_{H_1}(z)\Delta_0(\gamma_{H_1},\gamma_{0}) \lg \beta_p(\gamma_{0},[b]),s\rg  TO_{\delta_{[b]}}(f).\nonumber
\end{eqnarray*}
With $\lambda_{H_1}(z)$ canceled out, this finishes the proof in the unramified case.

The general case works in the same way. If $z\gamma_{H_1}$ is not a norm from $R_1(\Q_p)$ for any $z\in Z_1(\Q_p)$ then $SO_{\gamma_{H_1}}(f^{H_1})=0$ as before. Otherwise we assume that $\gamma'_{H_1}:=z\gamma_{H_1}$ is a norm for some $z$. Then we repeat the preceding argument via a choice of $f'_1\in \cH(R'_1)$ as before, with the difference occurring in the choice of $f^{H'_1}$. Namely by Remark \ref{rem:TFL}, there exists $f_0^{H'_1}\in \cH(H'_1(\Q_p))$ such that \eqref{eq:Mor10-951-A314} holds true with $f_0^{H'_1}$ in place of $(\tilde\eta'_1)^*(f'_1)$, and such that $SO_{\gamma'_{H_1}}(f^{H'_1}_0)=0$ if $\gamma'_{H_1}$ is not a norm. Setting $f^{H'_1}(h'):=\int_{Z'_{1}(\Q_p)} \lambda^{-1}_{H'_1}(z') f_0^{H'_1}(z'h') dz'$, we see from the above computation that Proposition \ref{prop:transfer-at-p} holds. The proof is complete.
\end{proof}

\subsection{}\label{para:H-unramified}
The last part of Proposition \ref{prop:transfer-at-p} can be slightly generalized, following \cite[p.~181]{Kot90}. Assume that $H$ and $G$ are unramified, and choose $H_1$ and $G_1$ to be unramified, but allow $\fke_1$ to be ramified. In that case, one can write $\eta_1=c \cdot \eta^\circ_1$ with a continuous 1-cocycle $c:W_{\Q_p}\ra Z(\hat H_1)$ such that $(H_1,\cH_1,s_1,\eta^\circ_1)$ is unramified.  (To see this, apply \cite[Prop.~1]{Lan79} for an unramified extension to find $\eta^\circ_1$, and observe that $\eta_1$ and $\eta^\circ_1$ must differ by a continuous 1-cocycle valued in $Z(\hat H_1)$.) Via local class field theory with geometric normalization, $c$ determines a smooth character $\chi_c:H_1(\Q_p)\ra \CC^\times$. Writing $f^{H_1,\circ}:=\tilde\eta^{\circ,*}_1(f)$, we see that $f^{H_1}:=\chi^{-1}_c \cdot f^{H_1,\circ}$ satisfies the desired condition of Proposition \ref{prop:transfer-at-p}. This follows from the fact that $\Delta_0(\gamma_{H_1},\gamma_{0,1})$ gets multiplied by $\chi^{-1}_c(\gamma_{H_1})$ when changing from $\eta^\circ_1$ to $\eta_1$. (This factor comes from $\Delta_2$ in \cite[\S 3.5]{LS87} as the 1-cocycle $a$ there is replaced with $ac$.)

\begin{prop}\label{prop:ramified-transfer-at-p=0}
In the setting of Proposition \ref{prop:transfer-at-p}, assume that $G$ is unramified over $\Q_p$ and $f\in \cH^{\ur}(R(\Q_p))$. If $H$ is ramified over $\Q_p$ then we can take $f^{H_1}=0$ (i.e., the right hand side of \eqref{e:tw-orb-id} always vanishes).
\end{prop}

\begin{proof}
The untwisted analogue of this proposition is proven in \cite[Prop.~7.5]{Kot86}. We adapt this proof, referred to as ``\emph{loc.~cit.}'' below, to our twisted setting. (A proof in the twisted case is alluded to on p.~189 of \cite{Kot90}. We are elaborating on the details.) We write $G'$ for the group $G_1$ in \textit{loc.~cit.}, as we reserve the symbol $G_1$ to stand for a $z$-extension.
Via $z$-extensions we reduce to the case where $G_{\der}=G_{\sconn}$ and $\cH={}^L H$. Recall we are also assuming $s\in Z(\hat H)^{\Gamma_p}$.

Fix $[b]\in \fkD_n(\gamma_0,G;\Q_p)$ and a representative $b\in G(\LL)$ of $[b]$. Fix $c\in G(\LL)$ as in  condition \textbf{KP1} of Definition \ref{pn-adm} so that
$c^{-1}\gamma_0 c = \delta \sigma(\delta) \cdots {\sigma^{n-1}} ( \delta) $
with $\delta=c^{-1} b\sigma (c) \in G(\Q_{p^n})$ (which is also denoted $\delta_{[b]}$). As in the proof of Lemma \ref{lem:Kottwitz-triples-z-ext-1}, the following sets are naturally in bijection with each other:
\begin{enumerate}
\item the set of $\sigma$-conjugacy classes in the stable $\sigma$-conjugacy class of $\delta$,
\item $\ker(H^1(\Q_p,I_\delta)\ra H^1(\Q_p,R))$, and
\item $\fkD_n(\gamma_0,G;\Q_p)$.
\end{enumerate}
To go from (i) to (ii), let $\delta' = x \delta \theta(x)^{-1}$ with $x\in R(\Qpbar)$, where $\theta$ denotes the automorphism of $R$ induced by $\sigma\in \Gal(\Q_{p^n}/\Q_p)$. Sending $ \delta'$ to the cocycle $z_{\delta,\delta'}:\tau\mapsto x^{-1} \lix{^\tau} x$ induces the bijection from (i) onto (ii). To go from (i) to (iii), take 
$c'\in G(\LL)$ such that
$(c')^{-1}\gamma_0 c' = \delta' \sigma(\delta') \cdots {\sigma^{n-1}} ( \delta') $.
Then we send $\delta'$ to $[c' \delta' \sigma(c')^{-1}]\in B(I_0)$, which lies in the set (iii).
Moreover we have the following compatibility: it follows from the bottom commutative diagram in \cite[p.~273]{kottwitzisocrystal2} (with $I_\delta$ and $I_0$ in place of $J^h$ and $H$ there, and the cocycle $h$ determined by $\sigma\mapsto c\delta \sigma(c)^{-1}$) that the image of $z_{\delta,\delta'}$ under the composite map
$$H^1(\Q_p,I_\delta) \cong \pi_1(I_\delta)_{\Gamma_p,\tors}\cong \pi_1(I_0)_{\Gamma_p,\tors}\hra  \pi_1(I_0)_{\Gamma_p}$$
coincides with $\kappa_{I_0}([b'])-\kappa_{I_0}([b])$ under the composite map
$$\fkD_n(\gamma_0,G;\Q_p)\hra \B(I_0) \xrightarrow{\kappa_{I_0}} \pi_1(I_0)_{\Gamma_p}.$$
Write $\textup{inv}(\delta,\delta')$ for the image in $\pi_1(I_0)_{\Gamma_p}$.
With the preparation so far, we follow \emph{loc.~cit.} to construct an exact sequence of unramified reductive groups over $\Q_p$
$$1 \ra G \ra G' \ra C \ra 1,$$
where $C$ is a non-trivial unramified torus. Define $R':=\Res_{\Q_{p^n}/\Q_p} G'$. Via the dual map $\hat G'\ra \hat G$, we can pull-back $\eta: \hat H \ra \hat G$ to define $\hat H'$ equipped with an embedding $\hat H' \hookrightarrow \hat G'$. We equip $\hat H'$ with a $\Gamma_p$-action as in \emph{loc.~cit.}

Write $I'_{\delta}$ for the connected $\sigma$-centralizer of $\delta$ in $R'$. As in \emph{loc.~cit.~}we have an exact sequence
$1\ra I_\delta \ra I'_\delta \ra C \ra 1$, whose dual exact sequence fits in the following commutative diagram, where rows are $\Gamma_p$-equivariant and exact:
$$\xymatrix{
1 \ar[r] & \hat C \ar@{=}[d]  \ar[r] & Z(\hat H') \ar[d] \ar[r] & Z(\hat H) \ar[r] \ar[d] & 1\\
1 \ar[r] & \hat C  \ar[r] & Z(\hat I'_\delta) \ar[r] & Z(\hat I_\delta) \ar[r] & 1
}$$
We have $s\in Z(\hat H)^{\Gamma_p}\subset Z(\hat I_{\delta})^{\Gamma_p}$. Write $\chi_s$ for the image of $s$ in $H^1(\Q_p,\hat C)$ under the connecting homomorphism arising from the first row. Then $\chi_s$ determines a smooth character $C(\Q_p)\ra \CC^\times$, still denoted by $\chi_s$. As in the proof of \cite[Prop.~7.5]{Kot86}, the character $\chi_s$ is \emph{ramified}, i.e., non-trivial on $C(\Z_p)$. (As $C$ is unramified, it extends uniquely to a $\Z_p$-torus.) This is crucial for our proof.

Now consider $g_1\in G'(\Q_p)$ and put $ \delta' := g_1 \delta  \theta( g_1)^{-1}=g_1 \delta g_1^{-1}\in G(\Q_{p^n})$. We claim that $\delta'$ is stably $\sigma$-conjugate to $\delta$. Indeed, writing $c_1\in C(\Q_p)$ for the image of $g_1$, take a lift $g_2\in I'_{\delta}(\LL)$ of $c_1^{-1}$ via the surjection $I'_{\delta}\ra C$. Set $g:=g_1 g_2 \in R(\LL)$. Then we have $\delta'= g \delta  \theta( g)^{-1}$, which proves the claim. 

Let $[b']$ denote the image of $\delta'$ under the bijection from (i) to (iii) above. 
Then we assert that
$$\langle \beta_p(\gamma_{0},[b']),s\rangle \langle \beta_p(\gamma_{0},[b]),s\rangle ^{-1}
= \langle \textup{inv}(\delta,\delta'),s\rangle= \chi_s(c_1^{-1}).$$
The first equality follows from the aforementioned compatibility. The second equality 
comes from \cite[Lem.~1.6]{Kot86} (applied to $I=C$ and $G=I_\delta$), 
where the image of $c_1^{-1}$ in $H^1(\Q_p,I_\delta)$ is represented by the cocycle $z_{\delta,\delta'}$. Indeed, after applying the injection from $H^1(\Q_p,I_\delta)$ into $B(I_\delta)=H^1(W_{\Q_p},I_\delta(\ol{\breve{\Q_p}}))$, they are represented by the same cocycle since $ g_2^{-1} {}^\tau g_{2}=g^{-1} {}^\tau g$ for $\tau\in W_{\Q_p}$.

Choose $T$ to be the generic fiber of the torus in a Borel pair of $\cG$ and set $T'=\Cent(T,G')$. With the compatible integral models, $T'(\Z_p)\to C(\Z_p)$ is surjective. Choose $c_1\in C(\Z_p)$ with $\chi_s(c_1)\ne 1$ and lift it to $g_1\in T'(\Z_p)$. For $f\in \cH^{\ur}(R(\Q_p))$, define $f_0\in \cH^{\ur}(R(\Q_p))$ as
$x\mapsto f(g_1^{-1} x g_1)$ with $g_1$ as above. (The analogue of $f_0$ is denoted by $f_1$ in \emph{loc.~cit.})
Write $f^H$ and $f_0^H$ for their twisted transfers to $H$. On the one hand, we have $f=f_0$ by the argument of \emph{loc.~cit.}, so we can take $f^H=f_0^H$. On the other hand, comparing the right hand sides of \eqref{e:tw-orb-id}, we have
$$
SO_{\gamma_H}(f_0^H) = \langle \beta_p(\gamma_{0},[b']),s\rangle
\langle \beta_p(\gamma_{0},[b]),s\rangle ^{-1} SO_{\gamma_H}(f^H)= \chi_s(c_1^{-1})SO_{\gamma_H}(f^H)
$$
if $\gamma_H$ is a norm of some $\delta$ as above, hence $SO_{\gamma_H}(f^H)=0$. If $\gamma_H$ is not a norm of any $\delta$, $SO_{\gamma_H}(f^H)=0$ by Proposition~\ref{prop:transfer-at-p}. We conclude that the stable orbital integral of $f^H$ vanishes everywhere.
\end{proof}

\subsection{}\label{para:fn}

Finally we drop the assumption that $s\in Z(\hat H)^{\Gamma_p}$ and consider the general case where $s\in Z(\hat H)^{\Gamma_p}Z(\hat G)$.
Write $s=s's''$ with $s'\in Z(\hat H)^{\Gamma_p}$ and $s''\in Z(\hat G)$. As in \S\ref{para:twisted-datum}, $(H,\cH,s',\eta)\in E(G)$ yields a twisted endoscopic datum $(H_1,{}^L H_1,\tilde s'_1,\tilde \eta'_1)$ for $(R_1,\theta_1)$, with $s'$ playing the role of $s$.

Consider the setting of  \S\ref{subsubsec:setting for point count}, where $\phi_n\in \cH^{\ur}(R(\Q_{p}))$ was introduced. If $\fke_1$ is unramified, then take\index[Oct2026-sorted-notation]{f 1 00 00 f subscript n superscript h subscript 1@$f_n^{H_1}$} $f_n^{H_1}:=\mu^{-1}(s'')\tilde\eta_1^{\prime,*} (\phi_n)$. If $H$ is unramified but $\fke_1$ is ramified, write
$\eta_1=c\eta_1^\circ$ as in \S\ref{para:H-unramified}
and take
$
f_n^{H_1}:=\mu(s'')^{-1}\chi_c^{-1}
\tilde\eta_1^{\prime,\circ,*}(\phi_n),
$
where $\tilde\eta_1^{\prime,\circ}$ is constructed from
$s'$ and $\eta_1^\circ$. If $H$ (thus also $H_1$) is ramified over $\Q_p$ then take $f_n^{H_1}:=0$. We check that Proposition \ref{prop:transfer-at-p} extends to this case.

\begin{cor}\label{cor:transfer-at-p}
Let $\gamma_{H_1}\in H_1(\Q_p)_{(G_1,H_1)\textup{-reg}}$. We have $SO_{\gamma_{H_1}}(f_n^{H_1})=0$ if $z\gamma_{H_1}$ is not a norm from $R_1(\Q_p)$ for any $z\in Z_1(\Q_p)$. If $z\gamma_{H_1}$ is a norm of $\delta_1\in R_1(\Q_p)$ for some $z\in Z_1(\Q_p)$ then denoting by $\gamma_{0}\in G(\Q_p)$ an image of $\gamma_{H_1}$, we have
$$
SO_{\gamma_{H_1}}(f_n^{H_1})=\sum_{[b]\in \fkD_n(\gamma_0,G;\Q_p)} e(I_{\delta_{[b]}}) \Delta_0(\gamma_{H_1},\gamma_{0}) \lg \tilde \beta_p(\gamma_{0},[b]),s\rg TO_{\delta_{[b]}}(\phi_n),
$$
where the summand is understood to be zero if $TO_{\delta_{[b]}}(\phi_n)=0$. 
If $TO_{\delta_{[b]}}(\phi_n)\neq~0$ then $[b]$ satisfies \textbf{KP0} and \textbf{KP1} in Definition \ref{pn-adm} so $\tilde \beta_p(\gamma_{0},[b])$ is well defined as explained in \S\ref{defn of Kottwitz invariant}.
\end{cor}

\begin{proof}

By definition $[b]\in \fkD_n(\gamma_0,G;\Q_p)$ always satisfies \textbf{KP1}. Let $\mu \in \dmu_X^{\cG}$ as in \S\ref{subsubsec:SD}.
If $TO_{\delta_{[b]}}(\phi_n)\neq 0$ for some $[b]\in \fkD_n(\gamma_0,G;\Q_p)$ then $X_{-\mu}(\delta_{[b]})$ is non-empty, so $[\delta_{[b]}]\in \B(G,-\mu)$ by \cite{rapoportrichartz}. Since $[b]=[\delta_{[b]}]$ in $\B(G)$, we see that \textbf{KP0} holds true.
Thus 
$$ \lg \tilde \beta_p(\gamma_{0},[b]),s\rg =\mu^{-1}(s'')  \lg \tilde \beta_p(\gamma_{0},[b]),s'\rg = \mu^{-1}(s'')  \lg  \beta_p(\gamma_{0},[b]),s'\rg.$$
(The second equality holds because $s'\in Z(\hat H)^{\Gamma_p}\subset Z(\hat I_0)^{\Gamma_p}$ and $\beta_p(\gamma_{0},[b])$ is an element of $\pi_1(I_0)_{\Gamma_p}$.)
Now the proof follows from Proposition \ref{prop:transfer-at-p} with $\phi_n$ and $s'$ in place of $f$ and $s$, respectively, as $\mu^{-1}(s'')\in \CC^\times$ cancels out.
\end{proof}

\cleardoublepage
\chapter{Stabilization}\label{sec:stabilization} 

We return to the point counting formula for Shimura varieties. Taking Conjecture \ref{conj:point counting formula} for granted, we carry out the stabilization of the formula \eqref{pcf}, with a view towards a representation-theoretic description of the cohomology.

\section{Initial steps}\label{sub:initial}

We start by fixing a central character datum, rewrite the coefficients in \eqref{pcf}, and apply a Fourier transform on the finite abelian group $\fkK(I_0/\Q)$.

\subsection{}

We freely use the setting of \S \ref{subsec:stating pcf}. Throughout stabilization, we fix an unramified Shimura datum $(G,X,p,\cG)$, which determines $K_p=\cG(\Z_p)$, and an open compact subgroup $K^p\subset G(\A_f^p)$ such that $K=K_pK^p$ is a neat subgroup.

\subsection{}\label{para:initial-steps-measures}

Recall that $Z$ is the center of $G$. We endow $A_{Z,\infty}=A_{G,\infty}$, which is isomorphic to a finite product of copies of $\R^\times_{>0}$, with the standard multiplicative Haar measure. 
Fix Haar measures on $Z(\A_f)$ and $Z(\R)$, thereby also on $\fkX:=(Z(\A_f)\cap K)\cdot Z(\R)$ and $\fkX_{\Q}\bs \fkX/ A_{Z,\infty}$, relative to the counting measure on $\fkX_{\Q}$ (which is discrete and cocompact in $ \fkX/ A_{Z,\infty}$). With respect to the set of places of $\Q$, we have the decomposition $\fkX=\fkX^{p,\infty} \fkX_p \fkX_\infty$.
We put a Haar measure on $Z(\Q)\bs Z(\A)/\fkX$ via the following exact sequence of topological groups
$$1\To \fkX_{\Q}\bs \fkX/ A_{Z,\infty} \To Z(\Q)\bs Z(\A)/A_{Z,\infty} \To Z(\Q)\bs Z(\A)/\fkX \To 1,$$
where the space in the middle is given the\index{Tamagawa measure} Tamagawa measure. 

We promote $\fkX$ to a central character datum $(\fkX,\chi)$ by defining the finite part $\chi^\infty$ to be trivial on $Z(\A_f)\cap K$ and the infinite part $\chi_\infty:=\omega_{\xi}^{-1}$, where $\omega_{\xi}$ is the central character of $\xi$ on $Z(\R)$.

For a reductive  subgroup $G_0$ of $G$ over $\Q$ containing $Z$ such that $A_{G_0}=A_Z$, we use the quotient measure to define\index[Oct2026-sorted-notation]{zz19 tau 1 00 00 tau subscript x@$\tau_{\mathfrak X}$}
$$\tau_\fkX(G_0):=\vol(G_0(\Q)\bs G_0(\A)/\fkX)$$
by viewing the double coset space as the quotient of $G_0(\Q)\bs G_0(\A)/A_{G_0,\infty}$ by $\fkX_\Q\bs \fkX/A_{Z,\infty}$, where the former is equipped with the Tamagawa measure and the latter with the measure explained above. For any inner form $G_0'$ of $G_0$ over $\Q$ we make sense of $\tau_\fkX(G_0')$ by transporting the measure. Since the Tamagawa volumes are equal for $G_0$ and $G_0'$ by \cite{Kot88}, we have $\tau_\fkX(G_0)=\tau_\fkX(G_0')$.

Once and for all, fix a $z$-extension $G_1$ of $G$ over $\Q$ which is unramified over $\Q_p$. This is possible by Lemma \ref{lem:z-ext-choice} (recalling that $G$ is unramified over $\Q_p$). Write $\fkX_1\subset Z_{G_1}(\A)$ for the preimage of $\fkX$ under $G_1\rightarrow G$.
Write $\chi_1$ for the pull-back of $\chi$ to $\fkX_1$ (thus the finite part of $\chi_1$ is trivial). Then $(\fkX_1,\chi_1)$ is a central character datum for $G_1$, and $\fkX_1=\fkX_1^{p,\infty} \fkX_{1,p}\fkX_{1,\infty}$ analogously as the decomposition of $\fkX$. We put the unique Haar measure on $ \fkX_1$, and thus also on $\fkX_{1,\Q} \backslash \fkX_1/A_{G_1,\infty}$ similarly as above, such that
\begin{equation}\label{eq:Haar-measure-on-X1}
\vol(\fkX_{1,\Q} \backslash \fkX_1/A_{G_1,\infty})/\vol(\fkX_{\Q} \backslash \fkX/A_{G,\infty})
=\tau(G_1)/\tau(G),
\end{equation}
so that
\begin{equation}\label{eq:tauG1=tauG}
\tau_{\fkX_1}(G_1)=\tau_{\fkX}(G).
\end{equation}

Let $\fkc=(\gamma_0,a,[b])\in \KP_{\adm}(p^n)$ for some $n\in\Z_{\ge1}$. Recall from \S\ref{local inner forms} that $\fkc$ determines an inner form $I_v$ of $I_0$ over $\Q_v$. Note that $I_p$ depends only on $\gamma_0$ and $[b]$ (not on $a$). Since $I_\infty$ is always a compact-mod-center inner form of $I_{0,\R}$ we often write $I_0^{\cpt}$ for $I_\infty$. If the Kottwitz invariant $\alpha(\fkc)$ (defined in \S\ref{defn of Kottwitz invariant}) vanishes, then we have an inner form $I$ of $I_0$ over $\Q$ which localizes to $I_v$ at each $v$ (Proposition \ref{existence of I}). In this case, choose the Haar measure $di_\infty$ on
$I_\infty(\R)$ so that $di_p\,di^p\,di_\infty$ induces the Tamagawa
measure on $I(\Q)\bs I(\A)/A_{Z,\infty}$, using counting measure on
$I(\Q)$ and the fixed measure on $A_{Z,\infty}$. Equip
$Z(\R)\bs I_\infty(\R)$ with the quotient measure from $di_\infty$
and the fixed measure on $Z(\R)$. We defined the constants $c_1(\fkc,K^p,di_p di^p)$ and $c_2(\gamma_0)$ in \S\ref{para:finite c_1}.

\begin{lem}\label{lem:c_1c_2}
If $\alpha(\fkc)$ vanishes, then for the measures chosen above,
$$c_1(\fkc,K^p,di_p di^p)c_2(\gamma_0)=\tau_{\fkX}(G)\cdot  |\fkK(I_0/\Q)|\cdot \vol(Z(\R)\bs I_{\infty}(\R))^{-1}.$$
\end{lem}

\begin{proof}
With the choice of measures as above,
\begin{eqnarray}
c_1(\fkc,K^p,di_p di^p)&=&\vol(I(\Q)\bs I(\A)/Z_K I_\infty(\R))\nonumber \\
&=&\frac{\vol(I(\Q)\bs I(\A)/Z_K Z(\R))}{\vol(Z(\R)\bs I_\infty(\R))}=\frac{\tau_{\fkX}(I_0)}{\vol(Z(\R)\bs I_\infty(\R))}.
\nonumber
\end{eqnarray}
On the other hand one deduces as in \cite[p.~395]{Kot86} that
$$c_2(\gamma_0)=\tau(I_0)^{-1}\tau(G)|\fkK(I_0/\Q)| = \tau_{\fkX}(I_0)^{-1}\tau_{\fkX}(G)|\fkK(I_0/\Q)|.$$
We conclude by taking product of the two equations.
\end{proof}

\subsection{}\label{para:initial}

Let $\fkc=(\gamma_0,a,[b])\in \KP_{\adm}(p^n)$. Define $e(\fkc):=\prod_v e(I_v) \in \{\pm1\}$, the product of  Kottwitz signs over all places. If $\alpha(\fkc)$ vanishes, then $e(\fkc)=1$ since $I_v$'s come from a $\Q$-group $I$. Hence

\begin{equation}\label{e:finite-Fourier}
\sum_{\kappa\in \fkK(I_0/\Q)} e(\fkc)\lg \alpha(\fkc),\kappa\rg=\left\{ \begin{array}{cl}
|\fkK(I_0/\Q)|, & \mbox{if}~\alpha(\fkc)=0,\\
0, & \mbox{otherwise}.\end{array}\right.
\end{equation}

Applying Lemma \ref{lem:c_1c_2} and \eqref{e:finite-Fourier} to \eqref{pcf}, we have
\begin{equation}\label{e:initial-steps} T(\Phi_{\mathfrak p}^m, f^p dg^p)  
= \tau_\fkX(G)\sum_{\gamma_0\in \Sigma_{\RR \text{-} \el,\fkX}(G)}\ol\iota_G(\gamma_0)^{-1}\sum_{\kappa\in \fkK(I_0/\Q)} \sum_{(a,[b])} N(\gamma_0,\kappa,a,[b]),
\end{equation}
where the third sum runs over $\fkD(I_0,G;\A_f^p)\times \fkD_n(\gamma_0,G;\Q_p)$, and  we set $N(\gamma_0,\kappa,a,[b])=0$ if $TO_{\delta_{[b]}}(\phi_n)=0$; if $TO_{\delta_{[b]}}(\phi_n)\neq 0$ then \textbf{KP0} and \textbf{KP1} hold by Corollary \ref{cor:transfer-at-p}, and we set
\begin{equation}
N(\gamma_0,\kappa,a,[b])=\lg \alpha(\gamma_0,a,[b]),\kappa\rg e(\gamma_0,a,[b]) \frac{O_{\gamma_{a}}(f^p) TO_{\delta_{[b]}}(\phi_n) \tr \xi(\gamma_0)}{\vol(Z(\R)\bs I_\infty(\R))}.\nonumber
\end{equation}

\section{Local transfer of orbital integrals}\label{sub:local-transfer}

Here we prove adelic orbital integral identities for each elliptic endoscopic datum for $G$ by dividing into three cases: away from $\{p,\infty\}$, at $p$, and at $\infty$.

\subsection{}
Let $\fke=(H,\cH,s,\eta)\in E_{\el}(G)$. The fixed $z$-extension $G_1$ (which is unramified over $\Q_p$) gives rise to a central extension $H_1$ of $H$ and an endoscopic datum $\fke_1$ for $G_1$ via Lemmas \ref{lem:z-ext-H} and \ref{lem:endoscopic-datum-z-ext}. In particular, if $H$ is unramified over $\Q_p$ then so is $H_1$.
As before, $\fkX_{H_1}$ is the preimage of $\fkX\subset Z_G(\A)\subset Z_H(\A)$ under $Z_{H_1}\rightarrow Z_H$.
For every place $v$, let $\Lambda_{H_1,v}: Z_{G_1}(\Q_v) \ra \CC^\times$ be the character in \S\ref{ss:lambdaH} applied to $\fke_1$, which is unramified at almost all $v$. We obtain a continuous character \[\Lambda_{H_1}:=\textstyle\prod_v \Lambda_{H_1,v}:Z_{G_1}(\A)\ra \CC^\times.\] By Lemma \ref{lem:lambdaH}, the restriction of $\Lambda_{H_1}^{-1}$ to $Z^0_{G_1}(\A)$ corresponds to the parameter \[W_{\Q}\ra {}^L H_1 \stackrel{\eta_1}{\ra} {}^L G_1 \ra {}^L Z_{G_1}^0.\] Moreover, $\Lambda_{H_1}=1$ on $Z_{G_1}(\Q)$. To see this, we reduce to a quasi-split inner form
of $G_1$ by inner twisting. Then we apply \eqref{eq:trans-Z-equiv} to $z\in Z_{G_1}(\Q)$ and strongly regular $\gamma\in G_1(\Q)$ and $\gamma_H\in H_1(\Q)$ with matching conjugacy classes at every $v$. Taking the product over all $v$, we conclude $\Lambda_{H_1}(z)=1$ by the product formula for transfer factors. Now we restrict $\Lambda_{H_1}$ to $\fkX_1$ to define
\begin{equation}\label{eq:lambdaH1}
\lambda_{H_1}:\fkX_{1,\Q}\backslash \fkX_{1}\To  \CC^\times,
\end{equation}
which can be viewed as a character of $\fkX_{H_1}$ via $\fkX_1\cong \fkX_{H_1}$. According as $\fkX_1=\fkX_1^{p,\infty}\fkX_{1,p}\fkX_{1,\infty}$ we decompose $\lambda_{H_1}=\lambda_{H_1}^{p,\infty}\lambda_{H_1,p}\lambda_{H_1,\infty}$. With the Haar measure on $\fkX_1$ transferred to $\fkX_{H_1}$ via the isomorphism, the analogue of \eqref{eq:tauG1=tauG} holds.

\begin{lem}\label{lem:tauH1=tauH}
We have $\tau_{\fkX_{H_1}}(H_1) = \tau_{\fkX}(H)$.
\end{lem}

\begin{proof}
Once we prove that $\tau(G)/\tau(G_1)=\tau(H)/\tau(H_1)$, the lemma follows in the same way as \eqref{eq:tauG1=tauG} is implied by \eqref{eq:Haar-measure-on-X1}. 
By \cite[(5.2.3), \S5.3]{Kot84a} ($\tau_1(\cdot)$ therein is $\tau(\cdot)$ by \cite{Kot88}), we have 
\begin{align*}
\tau(G)/\tau(G_1) &=\left|\cok(X_*( Z(\hat G_1))^\Gamma \ra X_*( \hat Z_1)^\Gamma)\right|,\\
\tau(H)/\tau(H_1)&=\left|\cok(X_*( Z(\hat H_1))^\Gamma \ra X_*( \hat Z_1)^\Gamma)\right|.
\end{align*}
So it is enough to show that the two cokernels are isomorphic. 
Consider the commutative diagram below, where the rows are coming from the exact sequence of \cite[Cor.~2.3]{Kot84a}. 

$$ \xymatrixcolsep{1.6pc} \xymatrix{
1 \ar[r] &  X_*(Z(\hat G))^\Gamma \ar[r]^-{i_G} \ar@{^(->}[d] &  X_*(Z(\hat G_1))^\Gamma \ar[r] \ar@{^(->}[d] & X_*(\hat Z_1)^\Gamma  \ar@{=}[d]  \\
1 \ar[r] &  X_*(Z(\hat H))^\Gamma \ar@{->>}[d]  \ar[r]^-{i_H} &  X_*(Z(\hat H_1))^\Gamma \ar[r] \ar@{->>}[d]  &  X_*(\hat Z_1)^\Gamma \\
& X_*({Z(\hat H)}/{Z(\hat G)})^\Gamma \ar@{=}[r] & X_*({Z(\hat H_1)}/{Z(\hat G_1)})^\Gamma 
}$$
We obtain $\cok i_G\cong \cok i_H$ by diagram chase, and thereby the two cokernels above are isomorphic.
\end{proof}

\subsection{}\label{trans-away-from-p}

For the moment we assume that $s\in Z(\hat H)^{\Gamma_p}$, and we are going to drop this assumption in \S\ref{para:fH1} below.
By Proposition \ref{prop:untwisted-transfer} there exists a transfer\index[Oct2026-sorted-notation]{f 1 00 00 f superscript h subscript 1 p infty@$f^{H_1,p,\infty}$} $$f^{H_1,p,\infty}\in \cH(H_1(\A_f^{p,\infty}), \chi_1^{-1}\lambda_{H_1}^{p,\infty})$$ of $f^{p,\infty}\in \cH(G(\A_f^p){\sslash} K^p)$ with the following property. For each $$\gamma_{H_1}\in H_1(\A_f^p)_{(G_1,H_1)\textup{-reg}},$$ if it has no image in $G_1(\A_f^p)_{\semi}$ then $ SO_{\gamma_{H_1}}(f^{H_1,p,\infty})=0$.
If there exists an image $\gamma_{0,1}\in G_1(\A_f^p)_{\semi}$ of $\gamma_{H_1}$, then writing $\gamma_0\in G(\A_f^p)$ for the projection of $\gamma_{0,1}$, we have
\begin{eqnarray}
SO_{\gamma_{H_1}}(f^{H_1,p,\infty})&=&\sum_{a_1} e(a_1) \Delta(\gamma_{H_1},\gamma_{0,1,a_1}) O_{\gamma_{0,1, a_1}}(f^{p,\infty})
\label{e:SO-away-p-infty}\nonumber\\
&=& \sum_{a_1}  e(a_1) \Delta(\gamma_{H_1},\gamma_{0,1}) \lg  \beta^{p,\infty}(\gamma_{0,1},a_1) ,s_1 \rg O_{\gamma_{0,1, a_1}}(f^{p,\infty})\nonumber \\
&=& \sum_{a} \left(\prod_{v\neq p,\infty} e(I_v)\right) \Delta(\gamma_{H_1},\gamma_{0,1}) \lg  \beta^{p,\infty}(\gamma_{0},a) ,s \rg O_{\gamma_{0, a}}(f^{p,\infty}),
\end{eqnarray}
where the sums for $a_1$ and $a$ run over $\fkD(I_{0,1},G_1;\A_f^p)$ and $\fkD(I_{0},G;\A_f^p)$  respectively. Recall that $\beta^{p,\infty}(\cdot,\cdot)$ was introduced in \S\ref{subsub:local_Kottwitz_inv}.
The second equality above follows from a basic property of transfer factors regarding the change of $\gamma_{0,1}$ within its stable conjugacy class as stated in \cite[Conj.~5.5]{Kot86}, which can be proved by arguing as in the proof of \cite[Lem.~4.1.C ]{LS87} in the $G_1$-regular case and extended to the $(G_1,H_1)$-regular case by \cite{LS90}.

\subsection{}\label{trans-at-p}

From \S\ref{para:fn} and Corollary \ref{cor:transfer-at-p}, we obtain\index[Oct2026-sorted-notation]{f 1 00 00 f superscript h subscript 1 subscript p@$f^{H_1}_{p}$} $f^{H_1}_{p}\in \cH(H_1(\Q_p),\lambda_{H_1,p})$ (renaming $f^{H_1}_n$) with the following property. Let $\gamma_{H_1}\in H_1(\Q_p)_{(G_1,H_1)\textup{-reg}}$. If $z\gamma_{H_1}$ is not a norm from an element of $R_1(\Q_p)$ for any $z\in Z_1(\Q_p)$ then $SO_{\gamma_{H_1}}(f^{H_1}_p)=0$. If $z\gamma_{H_1}$ is a norm for some $z\in Z_1(\Q_p)$, there exists $\gamma_{0,1}\in G_1(\Q_p)$ whose conjugacy class matches $z\gamma_{H_1}$. Writing $\gamma_0\in G(\Q_p)$ for the projection of $\gamma_{0,1}$, we have
\begin{equation}\label{e:SO-p}
SO_{\gamma_{H_1}}(f^{H_1}_p)=\sum_{[b]\in \fkD_n(\gamma_0,G;\Q_p)} e(I_{\delta_{[b]}}) \Delta_0(\gamma_{H_1},\gamma_{0}) \lg \tilde \beta_p(\gamma_0,[b]), s \rg TO_{\delta_{[b]}}(\phi_n).
\end{equation}
By definition $I_{\delta_{[b]}}=I_p$ if $[b]$ comes from $\fkc=(\gamma_0,a,[b])\in \KP_{\adm}(p^n)$.
Thanks to Proposition \ref{prop:ramified-transfer-at-p=0} we may and will take $f^{H_1}_p=0$ if $H$ is ramified over $\Q_p$.

\subsection{}\label{h-infty}

Our starting point for real orbital integrals is the argument of \cite[\S7]{Kot90} based on Shelstad's real endoscopy \cite{She82} and the pseudo-coefficients of  Clozel--Delorme \cite{CD90}. We incorporate central characters and $z$-extensions.

We have the characters $\chi_{1,\infty}$ and $\lambda_{H_1,\infty}$ on $\fkX_{H_1,\infty}$. We will drop $\infty$ from the subscript and write $\lambda_{H_1}$ and $\chi_1$ when the context is clearly local archimedean. Let $\xi_1$ be the irreducible representation of $G_1$ obtained from $\xi$ via the surjection $G_1\rightarrow G$.
As in \S\ref{real place z-ext}, fix a Shimura datum $(G_1,X_1)$ and $\mu_1\in \dmu_{X_1}(\Qbar)$. For each elliptic $\gamma_{0,1}\in G_1(\R)$ we have an element $\tilde\beta_\infty(\gamma_{0,1})\in \pi_1(I_{0,1})$ which maps to $[\mu_1]\in \pi_1(G_1)$, cf.~\S\ref{subsub:local_Kottwitz_inv}.

We simply set $f^{H_1}_\infty:=0$ if elliptic maximal tori of $G_{\R}$ do not come from those of $H_{\R}$ via transfer (equivalently, if elliptic maximal tori of $G_{1,\R}$ do not come from those of $H_{1,\R}$). In particular $f^{H_1}_\infty=0$ if $H_{\R}$ contains no elliptic maximal tori. From now on, we assume that $H_{\R}$ contains elliptic maximal tori and they transfer to elliptic maximal tori of $G_{\R}$. In particular $Z_{G_{\R}}^0$ and $Z_{H_{\R}}^0$ have the same $\R$-rank. 

Following Kottwitz, we construct a smooth function\index[Oct2026-sorted-notation]{f 1 00 00 f superscript h subscript 1 subscript infty@$f^{H_1}_\infty$} $f^{H_1}_\infty$ on $H_1(\R)$ that is compactly supported modulo center, by taking a suitable finite linear combination of pseudo-coefficients of discrete series representations (explicitly as on p.~186 of \cite{Kot90} with our $H_1$ and $G_1$ playing the roles of his $H$ and $G$) such that the following hold:
if $\gamma_{H_1}\in H_1(\R)_{\semi}$ is elliptic and $(G_1,H_1)$-regular, which by the above assumption has an elliptic element $\gamma_{0,1}\in G_1(\R)_{\semi}$ as image, then
\begin{equation}\label{e:terms-infty-1}
SO_{\gamma_{H_1}}(f^{H_1}_\infty)=\frac{ e(I_{0,1}^{\cpt})\lg \tilde\beta_\infty(\gamma_{0,1}),s_1\rg \Delta_\infty(\gamma_{H_1},\gamma_{0,1})
\tr \xi_1(\gamma_{0,1})}{\vol(Z_{G_1}(\R)\bs I_{0,1}^{\cpt}(\R))},
\end{equation}
whereas if $\gamma_{H_1}$ is non-elliptic and $(G_1,H_1)$-regular then $SO_{\gamma_{H_1}}(f^{H_1}_\infty)=0$. Here we choose the pseudo-coefficients in
$\cH(H_1(\R),\chi_{1,\infty}^{-1}\lambda_{H_1,\infty})$,
with the quotient measures used in \eqref{e:terms-infty-1}. If $\gamma_{H_1}\in H_1(\R)_{\semi}$ is not $(G_1,H_1)$-regular, then by \cite[Prop.~3.3.4, Rem 3.3.5]{Mor10} adapted to our setting we have 
$SO_{\gamma_{H_1}}(f^{H_1}_\infty)=0$.

We claim that if $\gamma_{H_1}$ is elliptic and $(G_1,H_1)$-regular then \eqref{e:terms-infty-1} implies that
\begin{equation}\label{e:terms-infty-2}
SO_{\gamma_{H_1}}(f^{H_1}_\infty)= \frac{e(I_{0}^{\cpt}) \lg \tilde\beta_\infty(\gamma_{0}),s\rg \Delta_\infty(\gamma_{H_1},\gamma_{0}) \tr \xi(\gamma_0)}{ \vol(Z_G(\R)\bs I_0^{\cpt}(\R))},
\end{equation}
where $\gamma_0$ is the image of $\gamma_{0,1}$ in $G(\R)$. It is routine to check term-by-term equalities between the right hand sides of \eqref{e:terms-infty-1} and \eqref{e:terms-infty-2}. We illustrate the idea by showing that $\lg \tilde\beta_\infty(\gamma_{0,1}),s_1\rg=\lg \tilde\beta_\infty(\gamma_{0}),s\rg$, leaving the rest to the reader. Let us write $s$ for $s_1$ since the latter is the image of $s$ under the inclusion $\hat G\ra \hat G_1$. We have a decomposition $s=s's''$ with $s'\in Z(\hat H)^{\Gamma_\infty}$ and $s''\in Z(\hat G)$. By definition (\cite[\S2]{Kot90}) the character $\tilde\beta_\infty(\gamma_0)$ (resp.~$\tilde\beta_\infty(\gamma_{0,1})$) restricts to a character on $Z(\hat I_0)^{\Gamma_\infty}$ (resp.~$Z(\hat I_{0,1})^{\Gamma_\infty}$) and a character on $Z(\hat G)$ (resp.~$Z(\hat G_1)$), each of which is determined by $\mu$ (resp.~$\mu_1$). They are related via the following commutative diagrams, from which it is obvious that $\lg \tilde\beta_\infty(\gamma_{0,1}),s's''\rg = \lg \tilde\beta_\infty(\gamma_{0}),s's''\rg$.
$$\xymatrix{ s'\in Z(\hat H)^{\Gamma_\infty} \ar[r] \ar[d] & Z(\hat I_0)^{\Gamma_\infty} \ar[r]^-{\tilde\beta_\infty(\gamma_0)} \ar[d] & \CC^\times & s''\in Z(\hat G) \ar[r]^-{\tilde\beta_\infty(\gamma_0)} \ar[d] & \CC^\times \\  Z(\hat H_1)^{\Gamma_\infty} \ar[r] & Z(\hat I_{0,1})^{\Gamma_\infty} \ar[ur]_-{\tilde\beta_\infty(\gamma_{0,1})} & & Z(\hat G_1) \ar[ur]_-{\tilde\beta_\infty(\gamma_{0,1})}}$$

\subsection{}\label{para:fH1}

So far we have constructed the function\index[Oct2026-sorted-notation]{f 1 00 00 f superscript h subscript 1@$f^{H_1}$} $f^{H_1}:=f^{H_1,p,\infty}f^{H_1}_p f^{H_1}_\infty$ on $H_1(\adele)$.
If $H$ is ramified over $\Q_p$ or if elliptic tori of $G_{\R}$ do not come from those of $H_{\R}$, then we have $f^{H_1}=0$ since $f^{H_1}_p=0$ or $f^{H_1}_\infty=0$ in each case. When neither is the case,
$f^{H_1}$ depends only on the image of $s$ modulo $Z(\hat G)$. To see this, suppose that $s$ is replaced with $sz$ for $z\in Z(\hat G)$. Then $f^{H_1,p,\infty}$ as well as the identity \eqref{e:SO-away-p-infty} remains unchanged as $\prod_{v\neq p,\infty}\tilde\beta_v(\gamma_0,a)$ is trivial on $Z(\hat G)$. The function $f^{H_1}_p$ is multiplied by $\mu(z)^{-1}$ according to \S\ref{para:fn}. Since $\lg \tilde\beta_\infty(\gamma_0),s\rg$ is the only term in  \eqref{e:terms-infty-2}
to change and it is multiplied by $\lg \tilde\beta_\infty(\gamma_0),z\rg=\mu(z)$, the function $f^{H_1}_\infty$ is multiplied by  $\mu(z)$ to keep \eqref{e:terms-infty-2} valid. All in all, $f^{H_1}$ indeed remains invariant.

\subsection{}

Let us summarize the above results in terms of adelic orbital integral identities. To this end, we slightly extend the definition of $N(\gamma_0,\kappa,a,[b])$ in \S\ref{para:initial} from rational to adelic elements. Let $\fke=(H,\cH,s,\eta)\in \cE_{\el}(G)$.
Let $\gamma_0\in G(\A)$ be $\R$-elliptic and suppose that $\gamma_0$ is an image of $\gamma_{H_1}\in H_1(\A)_{(G_1,H_1)\textup{-reg}}$ at every place $v$. Define a similar quantity $N'(\gamma_{0},\gamma_{H_1},s, a ,[b])$ for the fixed $s$ of the endoscopic datum, $a\in \fkD(I_0,G;\A_f^p)$, and $[b]\in \fkD_n(\gamma_0,G;\Q_p)$ by $N'(\gamma_{0},\gamma_{H_1},s, a ,[b])=0$ if $TO_{\delta_{[b]}}(\phi_n)$  vanishes; if it does not vanish, in light of the last statement of Corollary \ref{cor:transfer-at-p}, set
\begin{multline*}
N'(\gamma_{0},\gamma_{H_1},s, a ,[b]) =  \lg \tilde\beta^{p,\infty}(\gamma_0,a),s\rg \lg \tilde\beta_p(\gamma_0,[b]),s \rg \lg \tilde\beta_\infty(\gamma_0),s\rg \Delta_\A(\gamma_{H_1},\gamma_{0})  \\ \times 
(\textstyle\prod_v e(I_v)) O_{\gamma_{a}}(f^p) TO_{\delta_{[b]}}(\phi_n) \tr \xi(\gamma_0) \vol(Z(\R)\bs I_\infty(\R))^{-1}. 
\end{multline*}
To compare with \S\ref{para:initial}, if $\gamma_0\in G(\Q)_{\R\textup{-ell}}$ is an image of $\gamma_{H_1}\in H_1(\Q)_{(G_1,H_1)\textup{-reg}}$ then
\begin{equation}\label{eq:N'=N}
N'(\gamma_{0},\gamma_{H_1},s,a,[b])=N(\gamma_0,s,a,[b])
\end{equation}
by definition of
$\alpha(\gamma_0,a,[b])$ and the product formula that $\Delta_\A(\gamma_{H_1},\gamma_{0})=1$.

\begin{lem}\label{lem:adelic-orb-int-identity}  Let $\gamma_{H_1}\in H_1(\A)_{(G_1,H_1)\textup{-reg}}$. If $\gamma_{0}\in G(\A)_{\semi}$ is $\R$-elliptic and an image of $\gamma_{H_1}$ at every place $v$ then
$$SO_{\gamma_{H_1}}(f^{H_1})=\sum_{a\in \fkD(I_0,G;\A_f^p) \atop [b]\in \fkD_n(\gamma_0,G;\Q_p)} N'(\gamma_{0},\gamma_{H_1},s,a,[b]).$$
If no $\R$-elliptic $\gamma_{0}\in G(\A)_{\semi}$ is an image of $\gamma_{H_1}$ then
$SO_{\gamma_{H_1}}(f^{H_1})=0$.
If $\gamma_{H_1}\in H_1(\Q)_{\semi}$ is not $(G_1,H_1)$-regular then again $SO_{\gamma_{H_1}}(f^{H_1})=0$.
\end{lem}

\begin{proof}
If elliptic maximal tori of $G_{\R}$ do not come from those of $H_{\R}$ then no $\R$-elliptic element of $G(\R)$ is an image of an element of $H(\R)$. Then $f^{H_1}_\infty=0$ by construction, so the lemma holds. If $\gamma_{H_1,\infty}$ is not $(G_1,H_1)$-regular then we saw $SO_{\gamma_{H_1,\infty}}(f^{H_1}_\infty)=0$, so in particular $SO_{\gamma_{H_1}}(f^{H_1})=0$ for non-$(G_1,H_1)$-regular elements $\gamma_{H_1}\in H_1(\Q)_{\semi}$.

From now on, let $\gamma_{H_1}\in H_1(\A)_{(G_1,H_1)\textup{-reg}}$, and assume that elliptic maximal tori of $G_{\R}$ come from those of $H_{\R}$. If $\gamma_{H_1}$ is not $\R$-elliptic (or equivalently if $\gamma_{H_1}$ has no image in $G(\R)_{\el}$) then 
$SO_{\gamma_{H_1,\infty}}(f^{H_1}_\infty)=0$ by \cite[Lem.~3.1]{Kot92b} (asserting that the orbital integrals of pseudo-coefficients of discrete series representations are supported on elliptic elements). Thus the lemma is verified in this case.
If $\gamma_{H_1,v}$ does not have an image in $G_1(\Q_v)_{\semi}$ at some finite place $v$ then $SO_{\gamma_{H_1}}(f^{H_1}_v)=0$ by \S\ref{trans-away-from-p} and \S\ref{trans-at-p} so the lemma is again true.
In the remaining case, there exists $\gamma_{0}$ as in the lemma.
Then the desired equality follows from \eqref{e:SO-away-p-infty}, \eqref{e:SO-p}, and \eqref{e:terms-infty-2}.
\end{proof}

\section{Final steps}\label{sub:final}

Applying the adelic transfer of orbital integrals to the outcome of initial steps, we finish stabilization.

\subsection{}
Let us resume from the formula \eqref{e:initial-steps}.
To reparametrize the sum over $(\gamma_0,\kappa)$ in \eqref{e:initial-steps},
consider the set of equivalence classes\index[Oct2026-sorted-notation]{zz18 sigma 0 00 00 sigma k subscript ell g@$\Sigma{\mathfrak K}_{\el}(G)$}
$$\Sigma\fkK_{\el}(G):=\{(\gamma_0,\kappa) 
\mid \gamma_0\in G(\Q)_{\el},~\kappa\in \fkK(I_0/\Q)\}/{\sim},$$
where $(\gamma_0,\kappa)\sim(\gamma'_0,\kappa')$ if there exists $g\in G(\ol\Q)$ such that (i) $g \gamma_0 g^{-1}=\gamma'_0$, (ii) $g^{-1}{}^{\tau}g\in I_0(\ol{\Q})$ for every $\tau\in \Gamma$, and (iii) the inner twisting $I_{0,\ol \Q}\cong I_{\gamma'_0,\ol \Q}$ induced by $\mathrm{Int}(g)$ carries $\kappa$ to $\kappa'$. Define another set of equivalence classes\index[Oct2026-sorted-notation]{e 0 03 00 e sigma subscript ell g@$\mathcal{E}\Sigma_{\el}(G)$}
$$\cE\Sigma_{\el}(G):=\{(\fke,\gamma_H) \mid \fke=(H,\cH,s,\eta)\in E_{\el}(G),~\gamma_H\in \Sigma_{\el}(H)_{(G,H)\text{-} \reg}\}/{\sim},$$
where $\Sigma_{\el}(H)_{(G,H)\text{-} \reg}$ is the set of stable conjugacy classes of $(G,H)$-regular semi-simple elements of $H(\Q)$, and $(\fke,\gamma_H)\sim(\fke',\gamma'_H)$ if there is an isomorphism between endoscopic data $\fke=(H,\cH,s,\eta)$ and $ \fke'=(H',\cH',s',\eta')$ such that $\gamma_H$ is carried to $\gamma'_H$. (In particular $(\fke,\gamma_H)\sim (\fke,\gamma'_H)$ if there is an outer automorphism of $\fke$ mapping $\gamma_H$ to $\gamma'_H$.) Define an analogous set with the obvious surjection
\begin{equation}\label{e:cEtilde-to-cE}\index[Oct2026-sorted-notation]{e 0 03 00 e sigma superscript sim subscript ell g@$\mathcal{E}\Sigma^{\sim}_{\el}(G)$}
\cE\Sigma^\sim_{\el}(G):=\left\{(\fke,\gamma_H) \mid \begin{array}{c}\fke=(H,\cH,s,\eta)\in \cE_{\el}(G),\\\gamma_H\in \Sigma_{\el}(H)_{(G,H)\text{-} \reg}\end{array}  \right\}
\To  \cE\Sigma_{\el}(G).
\end{equation}
The outer automorphism group $\Out_{\Q}(\fke)$ acts transitively on each fiber of the map.
Let us define a map\index[Oct2026-sorted-notation]{e 0 04 23 tilde e@$\tilde{\mathfrak E}$}
$$\tilde\fkE:\cE\Sigma^\sim_{\el}(G)\To \Sigma\fkK_{\el}(G)\cup \{\emptyset\}.$$
We explained in \S\ref{subs:trans-conjugacy} how
$(\fke,\gamma_H)$ determines either a stable conjugacy class of $\gamma_0\in G(\Q)$ or $\emptyset$ (when there is no matching conjugacy class in $G(\Q)$).
In the latter case, set $\tilde\fkE:(\fke,\gamma_H)\mapsto \emptyset$. In the former, we map $(\fke,\gamma_H)$ to $(\gamma_0,\kappa)$, where $\kappa\in \fkK(I_0/\Q)$ is determined by the image of $s$ under the composition $Z(\hat H)\hra Z(\hat I_{\gamma_H}) \cong Z(\hat I_0)$. Here the canonical isomorphism comes from the fact that $I_0$ is an inner form of $I_{\gamma_H}$.  The ellipticity of $\gamma_0$ follows from that of $\gamma_H$. Letting $Z(\Q)$ act on each of $\Sigma_{\el}(G)$ and $\Sigma_{\el}(H)$ by multiplication, we see that $\tilde\fkE$ is $Z(\Q)$-equivariant.

By \cite[Lem.~9.7]{Kot86}, when $G_{\der}=G_{\sconn}$, the map $\tilde\fkE$ factors through a unique map $\fkE: \cE\Sigma_{\el}(G)\ra \Sigma\fkK_{\el}(G)\cup \{\emptyset\}$, which is a bijection over $\Sigma\fkK_{\el}(G)$. In particular, the image of $\tilde\fkE$ contains $\Sigma\fkK_{\el}(G)$. When $G_{\der}$ is not simply connected, write $\tilde\fkE_1$ and $\fkE_1$ for the analogous maps for the $z$-extension $G_1$. We have a commutative diagram (\textit{a priori} without $\fkE$)
\begin{equation}\label{eq:ESigma-SigmaK}
\xymatrix{
\cE\Sigma^\sim_{\el}(G_1) \ar@{->>}[r] \ar[d] & \cE\Sigma_{\el}(G_1) \ar[r]^-{\fkE_1} \ar[d] & \Sigma\fkK_{\el}(G_1)\cup \{\emptyset\} \ar[d]\\
\cE\Sigma^\sim_{\el}(G) \ar@/_1.0pc/[rr]_-{\tilde \fkE}  \ar@{->>}[r] & \cE\Sigma_{\el}(G) \ar@{-->}[r]^-{\exists!\, \fkE}   & \Sigma\fkK_{\el}(G)\cup \{\emptyset\},
}
\end{equation}
where the vertical maps are induced by $G_1\rightarrow G$, using $\cE_{\el}(G_1)\cong \cE_{\el}(G)$ (Lemma \ref{lem:fke_1=fke}) and part (i) of Lemma \ref{lem:z-ext-wrt-K-and-D}. We add that the right vertical is required to send $\emptyset$ to itself. By the $Z_{G_1}(\Q)$-equivariance of $\fkE_1$, there exists a unique $Z_G(\Q)$-equivariant map\index[Oct2026-sorted-notation]{e 0 04 00 e@$\mathfrak E$} $\fkE$ making the entire diagram commute. (The image of $\fkE$ does not contain $\emptyset$ if $G$ is quasi-split over $\Q$, since every $\gamma_H$ then has an image in $G(\Q)$.)

Write $n(\gamma_0,\kappa)$ for the size of the
$(G_{\gamma_0}/I_0)(\Q)$-orbit of $\kappa$
in $\fkK(I_0/\Q)$ under conjugation.

\begin{lem}\label{lem:endo-correspondence}
The map $\fkE$ is $Z(\Q)$-equivariant and bijective when restricted to the preimage of $\Sigma\fkK_{\el}(G)$. Moreover for $(\gamma_0,\kappa)\in \Sigma\fkK_{\el}(G)$ and $(\fke,\gamma_{H,0})\in \cE\Sigma^\sim_{\el}(G)$ in the preimage,
$$n(\gamma_0,\kappa)\ol\iota_G(\gamma_0)^{-1}=\lambda(\fke)^{-1}\sum_{\gamma_H \in \Sigma(H)~\textup{s.t.}\atop \fkE:(\fke,\gamma_H)\mapsto (\gamma_0,\kappa) } \ol\iota_H(\gamma_H)^{-1}.$$
\end{lem}

\begin{rem} We remark that, in the strongly regular case the lemma follows from \cite[Lem.~7.2.A]{KS99}, which covers twisted endoscopy.
\end{rem}

\begin{proof}
The $Z(\Q)$-equivariance was observed above. 
Since $\fkE_1$ is bijective over $\Sigma\fkK_{\el}(G_1)$ and $Z_{G_1}(\Q)$-equivariant, this property descends along the vertical surjections to tell us that $\fkE$ is bijective over $\Sigma\fkK_{\el}(G)$.
Finally the equality asserted in the lemma follows from \cite[Cor.~IV.3.6]{Lab04}.
\end{proof}

\subsection{}\label{para:fH1-final}

Recall that $K^p$ is a neat subgroup and that $\fkX=(Z(\A_f)\cap K)\cdot Z(\R)$. In particular
\begin{equation}\label{e:X-cap-Gder-trivial}
\fkX\cap G_{\der}(\ol{\Q})=\{1\}.
\end{equation}
Let $\cE\Sigma^\sim_{\el,\fkX}(G)$ denote\index[Oct2026-sorted-notation]{e 0 03 00 e sigma superscript sim subscript ell x g@$\mathcal{E}\Sigma^{\sim}_{\el,\fkX}(G)$} the quotient set of $\cE\Sigma^\sim_{\el}(G)$ by the obvious multiplication of $\fkX_\Q=\fkX\cap Z(\Q)$. Likewise, define\index[Oct2026-sorted-notation]{zz18 sigma 0 00 00 sigma subscript x h@$\Sigma_{\fkX}(H)$} $\Sigma_{\fkX}(H)$ and\index[Oct2026-sorted-notation]{zz18 sigma 0 00 00 sigma k subscript ell x g@$\Sigma{\mathfrak K}_{\el,\fkX}(G)$} $\Sigma\fkK_{\el,\fkX}(G)$ from $\Sigma(H)$ and $\Sigma\fkK_{\el}(G)$.
Let\index[Oct2026-sorted-notation]{e 0 04 23 tilde e subscript x@$\tilde{\mathfrak E}_{\fkX}$} $$\tilde\fkE_{\fkX}:\cE\Sigma^\sim_{\el,\fkX}(G)\ra \Sigma\fkK_{\el,\fkX}(G)\cup \{\emptyset\}$$ denote the map $\fkE$ pre-composed with \eqref{e:cEtilde-to-cE} modulo the $\fkX_\Q$-action.

\begin{cor}\label{cor:endo-correspondence}
The image of $\tilde{\fkE}_\fkX$ is equal to $\Sigma\fkK_{\el,\fkX}(G)$ if $G$ is quasi-split and contains $\Sigma\fkK_{\el,\fkX}(G)$ in general. For $(\gamma_0,\kappa)\in \Sigma\fkK_{\el,\fkX}(G)$ and for $(\fke,\gamma_{H,0})$ in the preimage of $(\gamma_0,\kappa)$ under $\tilde{\fkE}_\fkX$,
$$n(\gamma_0,\kappa)\ol\iota_G(\gamma_0)^{-1}=\lambda(\fke)^{-1}\sum_{\gamma_H \in \Sigma_\fkX(H)~\textup{s.t.}\atop \fkE_\fkX:(\fke,\gamma_H)\mapsto (\gamma_0,\kappa) } \ol\iota_H(\gamma_H)^{-1}.$$
\end{cor}

\begin{proof}
The only non-trivial point is to deduce the equality from Lemma \ref{lem:endo-correspondence}. Let $\fke=(H,\cH,s,\eta)\in E_{\el}(G)$, $\gamma_H,\gamma'_H\in H(\Q)_{(G,H)\text{-} \reg}$. Suppose that $\fkE$ maps both $(\fke,\gamma_H)$ and $(\fke,\gamma'_H)$ to $(\gamma_0,\kappa)$. What we need to show is that, if $\gamma'_H$ is stably conjugate to $z\gamma_H$ for some $z\in \fkX_{\Q}$, then $\gamma'_H$ is stably conjugate to $\gamma_H$.

By assumption, $\fkE$ maps $(\fke,z\gamma_H)$ to $(z\gamma_0,\kappa)$. Hence $z\gamma_0$ is stably conjugate to $\gamma_0$. In particular $gz\gamma_0g^{-1}=\gamma_0$ for some $g\in G(\ol{\Q})$. Thanks to \eqref{e:X-cap-Gder-trivial} the element $z$ is trivial, so $\gamma'_H$ is stably conjugate to $\gamma_H$ as desired.
\end{proof}

\begin{lem}\label{lem:Q-conj-in-A-conj}
Let $f^{H_1}$ be as in \S\ref{para:fH1}.
Assume that $\tilde\fkE_\fkX$ maps $(\fke,\gamma_H)$ to $\emptyset$. 
Then $SO_{\gamma_{H_1}}(f^{H_1})=0$ for every lift $\gamma_{H_1}\in H_1(\Q)$ of $\gamma_H$.
\end{lem}

\begin{proof}
Suppose that $SO_{\gamma_{H_1}}(f^{H_1})\neq 0$. By Lemma \ref{lem:adelic-orb-int-identity} there exists $\gamma_{0,1}\in G_1(\A)$ such that $\gamma_{0,1,v}$ is an image of $\gamma_{H_1,v}$ at every place $v$. To show that $(\fke,\gamma_H)$ is not mapped to $\emptyset$ under $\fkE$, it suffices to show the existence of an element of $G_1(\Q)$ which is stably conjugate to $\gamma_{0,1}$ in $G_1(\Q_v)$ for every place $v$.
Indeed, $(\fke_1,\gamma_{H_1})$ then does not map to $\emptyset$, thus $(\fke,\gamma_H)$ does not either, cf.~\eqref{eq:ESigma-SigmaK}. Thus the proof boils down to the case where $G=G_1$ and $H=H_1$ (with $G_{\der}=G_{\sconn}$). Henceforth we drop the subscript $1$.

Let $G^*$ denote a quasi-split inner form of $G$. Write $\gamma^*_0\in G^*(\Q)_{\semi}$ for the image of $\gamma_{H}$ under \eqref{e:Sigma(H)-Sigma(G)}.
Now recall that Labesse (see \cite[\S 2.6]{Lab99}, with $L= G, H= G^*$) constructs a non-empty subset $$ \mathrm{obs}_{\gamma_{0}^*} (\gamma_0)\subset \E (I_{\gamma_0^*}, G^*; \adele/\QQ)  \xlongequal{\mathrm{Cor.~\ref{cor:compare with Kottwitz}}} \fkK(I_{\gamma^*_0}/\Q)^D , $$ generalizing the construction of Kottwitz in \cite{Kot86}.

As in the second paragraph of \cite[p.~188]{Kot90}, the Chebotarev density theorem implies that the natural map $\fkK(I_{\gamma^*_0}/\Q_w)^D\ra \fkK(I_{\gamma^*_0}/\Q)^D$ is surjective for some finite place $w$.
In Labesse's construction, if we twist $\gamma_{0 ,w}$ within its stable conjugacy class by a class $c \in \fkK(I_{\gamma^*_0}/\Q_w)^D $, then $\mathrm{obs}_{\gamma_0^*} (\gamma_0)$ gets shifted by the image of $c$ in the abelian group $\fkK(I_{\gamma^*_0}/\Q)^D $. Hence, if we replace $\gamma_{0,w}$ by a stably conjugate element and keep the components of $\gamma_0$ outside $w$ unchanged, then we may arrange that
$ 0 \in\mathrm{obs}_{\gamma_0^*} (\gamma_0) \subset \fkK(I_{\gamma^*_0}/\Q)^D$. By \cite[Thm.~2.6.3]{Lab99}, the $G(\A)$-conjugacy class of $\gamma_0$ contains an element of $G(\Q)$, as desired. 
\end{proof}

\subsection{}\label{para:STH1}

For each $\fke=(H,\cH,s,\eta)\in \cE_{\el}(G)$, set\index[Oct2026-sorted-notation]{zz09 iota 1 00 00 iota g h@$\iota(G,H)$}
\begin{equation}\label{eq:iota}
\iota(G,H):=\tau(G)\tau(H)^{-1}\lambda(\fke)^{-1}=\tau_\fkX(G)\tau_\fkX(H)^{-1}\lambda(\fke)^{-1}.
\end{equation}
For $\gamma_H\in \Sigma_{\el,\fkX}(H)$, choose a representative $\tilde\gamma_{H}\in\Sigma_{\el}(H)$ and define $\mathrm{Stab}_{\fkX}(\gamma_H)$ to be the group of $z\in \fkX_\Q$ such that $z\tilde\gamma_H=\tilde\gamma_H$ in $ \Sigma_{\el}(H)$; this is independent of the choice of $\tilde\gamma_{H}$. The same argument as in \S\ref{central-character-data} shows finiteness of $\mathrm{Stab}_{\fkX}(\gamma_H)$.

Let us introduce the stable analogue of $T_{\el}$ (see \S\ref{sub:trace-formula-central-character}) for $H_1$ with respect to the central character datum $(\fkX_{H_1},\chi_{H_1})$, where $\chi_{H_1}:=\chi_1 \lambda_{H_1}^{-1}$. To check that it is indeed a central character datum, note that both $\chi_1$ and $\lambda_{H_1}$ are trivial on $\fkX_{H_1,\Q}$ by construction, cf.~\S\ref{para:initial-steps-measures} and \eqref{eq:lambdaH1}.
For $ h\in \cH(H_1(\A),\chi_{H_1}^{-1})$, set\index[Oct2026-sorted-notation]{st 0 00 00 st superscript h subscript 1 subscript ell chi subscript h subscript 1 h@$ST^{H_1}_{\el, \chi_{H_1}}(h)$}

\begin{equation*}
ST^{H_1}_{\el, \chi_{H_1}}(h) :=   \tau_{\fkX_{H_1}}(H_1) \sum_{\gamma_{H_1}\in \Sigma_{\el,\fkX_{H_1}}(H_1)} \ol\iota(\gamma_{H_1})^{-1} |\mathrm{Stab}_{\fkX_{H_1}}(\gamma_{H_1})|^{-1} SO_{\gamma_{H_1}}(h).
\end{equation*}
In practice, we have $\ol\iota(\gamma_{H_1})=1$ whenever $\gamma_{H_1}$ is $(G_1, H_1)$-regular
and has an image in $G_1$ with $G_{1,\textup{sc}}=G_{1,\textup{der}}$ (even though $H_{1,\textup{der}}$ need not be simply connected), cf.~\cite[Lem.~3.2]{Kot86}.

\begin{lem}\label{lem:compare-stabilizers} 
For $\gamma_{H_1}\in H_1(\Q)_{(G_1,H_1)\textup{-reg}}$,
$$|\mathrm{Stab}_{\fkX_{H_1}}(\gamma_{H_1})|=|\mathrm{Stab}_{\fkX}(\gamma_{H})|\ol\iota_H(\gamma_H).$$
\end{lem}

\begin{proof}
It suffices to construct a short exact sequence of groups
$$1\ra (H_{\gamma_H}/H^0_{\gamma_H})(\Q) \ra \mathrm{Stab}_{\fkX_{H_1}}(\gamma_{H_1}) \ra \mathrm{Stab}_{\fkX}(\gamma_{H}) \ra 1.$$
The third arrow is induced by the projection $\fkX_{H_1}\rightarrow \fkX$ and is clearly surjective. To construct the second arrow, given $\ol h\in (H_{\gamma_H}/H^0_{\gamma_H})(\Q)$, choose a lift $h\in H_{\gamma_H}(\ol\Q)$ and a further lift $h_1\in H_1(\ol\Q)$. Then $x_1:=h_1 \gamma_{H_1} h_1^{-1} \gamma_{H_1}^{-1} $ belongs to $Z_1(\Q)$, and moreover $x_1\in \mathrm{Stab}_{\fkX_{H_1}}(\gamma_{H_1})$ since $h_1 \gamma_{H_1} h_1^{-1}  = x_1 \gamma_{H_1}$. The assignment $\ol h \mapsto x_1$ is a well-defined homomorphism. To check this map is injective, suppose $h_1 \gamma_{H_1} h_1^{-1}=\gamma_{H_1}$. Then $h_1$ lies in the centralizer of $\gamma_{H_1}$ in $H_1$, which is connected. (To see this, choose a place $v$ such that $G_{1,\Q_v}$ is quasi-split, so that $\gamma_{H_1,v}$ has an image $\gamma_{1,v}\in G_{1,\Q_v}$. Since the centralizer of $\gamma_{1,v}$ is connected, and since $\gamma_{H_1}$ is $(G_1,H_1)$-regular, the same is true for $\gamma_{H_1,v}$ by \cite[Lem.~3.2]{Kot86}.) Thus the image of $h_1$ in $H_{\gamma_H}$ lies in $H^0_{\gamma_H}$, implying that $\ol h$ is trivial.

The composition of the two maps above is clearly trivial. 
Finally suppose that $x_1\in \mathrm{Stab}_{\fkX_{H_1}}(\gamma_{H_1})$ maps trivially into $ \mathrm{Stab}_{\fkX}(\gamma_{H})$. Then $x_1\in Z_1(\Q)$. To check that $x_1$ comes from $ (H_{\gamma_H}/H^0_{\gamma_H})(\Q) $, choose $h_1\in H_1(\ol\Q)$ such that $h_1 \gamma_{H_1} h_1^{-1}  = x_1 \gamma_{H_1}$. Its image $h$ in $H(\ol\Q)$ clearly centralizes $\gamma_H$. Writing $\ol h$ for the image of $h$ in $(H_{\gamma_H}/H^0_{\gamma_H})(\ol\Q)$, we see that $\ol h$ is $\Q$-rational and maps to $x_1$ by construction. 
\end{proof}

\begin{rem}
In our situation $|\mathrm{Stab}_{\fkX}(\gamma_H)|=1$. Indeed the stabilizer group is a finite subgroup of $Z(\Q)_K$, which is torsion-free as $K^p$ is neat.
\end{rem}

\begin{thm}[cf.~Theorem \ref{thm:intro stabilization} in the Introduction]\label{thm:end-of-stabilization}
Assume that Conjecture \ref{conj:point counting formula} is true (cf.~Remark \ref{rem:Deligne's conj}). For all sufficiently large $m$, with $f^{H_1}$ constructed as in \S\ref{para:fH1} for each $\fke\in \cE_{\el}(G)$, we have
$$
\sum_{i} (-1)^i \tr \bigg({\Phi}_{\mathfrak p}^m \times (f^p dg^p)  \mid \coh_c^i(\Sh_{\ol{E}},\xi)^{K_p} \bigg)  = \sum_{\fke\in \cE_{\el}(G)} \iota(G,H) ST^{H_1}_{\el,\chi_{H_1}}(f^{H_1}). $$   
\end{thm}

\begin{proof}
We compute the right hand side as follows.
\begin{eqnarray}
& &  \sum_{\fke\in \cE_{\el}(G)} \iota(G,H) ST^{H_1}_{\el,\chi_{H_1}}(f^{H_1})\nonumber\\
& = &\tau_{\fkX}(G)  \sum_{\fke\in \cE_{\el}(G)}  \lambda(\fke)^{-1} \sum_{\gamma_{H_1}\in \Sigma_{\el,\fkX_{H_1}}(H_1)} |\mathrm{Stab}_{\fkX_{H_1}}(\gamma_{H_1})|^{-1} SO_{\gamma_{H_1}}(f^{H_1})
\nonumber\\
& = & \tau_{\fkX}(G)  \sum_{\fke\in \cE_{\el}(G)}  \lambda(\fke)^{-1} \sum_{\gamma_{H}\in \Sigma_{\el,\fkX}(H)} \ol\iota(\gamma_H)^{-1} SO_{\gamma_{H_1}}(f^{H_1})
\nonumber\\
& = &
\tau_\fkX(G)\sum_{(\fke,\gamma_H)\in \cE\Sigma^\sim_{\el,\fkX}(G)} \lambda(\fke)^{-1}\ol\iota(\gamma_H)^{-1}SO_{\gamma_{H_1}}(f^{H_1})
\nonumber\\
& =&
\tau_\fkX(G) \sum_{(\gamma_0,\kappa)\in \Sigma\fkK_{\el,\fkX}(G)} \sum_{(\fke,\gamma_H)\in \cE\Sigma^\sim_{\el,\fkX}(G)
	\atop \fkE_\fkX:(\fke,\gamma_H)\mapsto (\gamma_0,\kappa) } \lambda(\fke)^{-1}\ol\iota(\gamma_H)^{-1}SO_{\gamma_{H_1}}(f^{H_1}) \nonumber\\
&  = &
\tau_\fkX(G) \sum_{(\gamma_0,\kappa)\in \Sigma\fkK_{\el,\fkX}(G)\atop \gamma_0:\,\R\textrm{-elliptic}} \sum_{a \in \fkD(I_0,G;\A_f^p) \atop [b]\in \fkD_n(\gamma_0,G;\Q_p)} n(\gamma_0,\kappa) \ol\iota(\gamma_0)^{-1} N(\gamma_0,\kappa,a,[b])
\nonumber\\
& = &
\tau_\fkX(G)
\sum_{\gamma_0\in\Sigma_{\RR\text{-}\el,\fkX}(G)}
\ol\iota(\gamma_0)^{-1}
\sum_{\kappa\in\fkK(I_0/\Q)}
\sum_{a\in\fkD(I_0,G;\A_f^p)
	\atop [b]\in\fkD_n(\gamma_0,G;\Q_p)}
N(\gamma_0,\kappa,a,[b]).
\nonumber
\end{eqnarray}
In the third, fourth, and fifth lines, $\gamma_{H_1}\in H_1(\Q)_{\semi}$ is an arbitrary lift of $\gamma_H\in H(\Q)_{\semi}$. Each summand is independent of the choice since $f^{H_1}$ transforms under the character $\chi_{H_1}^{-1}$, which is trivial on $Z_1(\Q)$. The last expression above is exactly the left hand side of the theorem by \S\ref{para:initial}.  It remains to  justify these equalities. 

The first equality uses Lemma \ref{lem:tauH1=tauH} and \eqref{eq:iota}.
The next one is based on Lemma \ref{lem:compare-stabilizers} and the bijection $\Sigma_{\el,\fkX_{H_1}}(H_1)\ra \Sigma_{\el,\fkX}(H)$ induced by the surjection $H_1(\Q)\rightarrow H(\Q)$. We also used $|\mathrm{Stab}_{\fkX}(\gamma_H)|=1$, and the vanishing of the summand if $\gamma_{H_1}$ is not $(G_1,H_1)$-regular by Lemma \ref{lem:adelic-orb-int-identity}.
To continue, the third equality is justified by Lemma \ref{lem:adelic-orb-int-identity} telling us that only $(G,H)$-regular $\gamma_H$ contributes to the sum.
The fourth equality follows from Lemma~\ref{lem:Q-conj-in-A-conj}. The second last equality is deduced from Lemma~\ref{lem:adelic-orb-int-identity},~Corollary~\ref{cor:endo-correspondence}, and~\eqref{eq:N'=N}, noting that $\gamma_{0,1}$ can be taken from $G_1(\Q)$ whenever $SO_{\gamma_{H_1}}(f^{H_1})\neq 0$ as shown in the proof of Lemma \ref{lem:Q-conj-in-A-conj}. The last equality follows from Proposition \ref{funct of alpha conj}, by transporting $(a,[b])$ under conjugation and expanding each $(G_{\gamma_0}/I_0)(\Q)$-orbit of $\kappa$ into its $n(\gamma_0,\kappa)$ elements. 
\end{proof}
\begin{thm}[cf.~Theorem \ref{thm:intro main} in the Introduction] \label{thm:announcement 2} With  notation as in Theorem \ref{thm:end-of-stabilization},  for $(G,X)$ of abelian type and all sufficiently large $m$, we have
$$
\sum_{i} (-1)^i \tr \bigg({\Phi}_{\mathfrak p}^m \times (f^p dg^p)  \mid \coh_c^i(\Sh_{\ol{E}},\xi)^{K_p} \bigg)  = \sum_{\fke\in \cE_{\el}(G)} \iota(G,H) ST^{H_1}_{\el,\chi_{H_1}}(f^{H_1}).$$ 
\end{thm}
\begin{proof}
This follows from Theorems \ref{thm:announcement 1} and  \ref{thm:end-of-stabilization}.
\end{proof}  

\cleardoublepage

\chapter{Spectral interpretation}\label{sec:spectral interpretation}

In order to read off spectral information from Theorem \ref{thm:end-of-stabilization}, we need to turn the geometric stable distribution into a spectral expansion.
After discussing the stable trace formula in \S\ref{sub:STF} when the test function is stable cuspidal at $\infty$, we will remark on the prospect for unconditional spectral interpretation in \S\ref{sub:speculation}.

\section{The stable trace formula}\label{sub:STF}

\subsection{}

Let $G$ be a quasi-split reductive  group over $\Q$ with a fixed $z$-extension $G_1$. (We do not assume that $G$ is part of a Shimura datum in \S\ref{sub:STF}.)
Let $(\fkX,\chi)$ be a central character datum such that $\fkX\supset A_{Z,\infty}$. For each elliptic endoscopic datum $\fke=(H,\cH,s,\eta)\in \cE_{\el}(G)$, choose a central extension $H_1$ and define $\fke_1\in \cE_{\el}(G_1)$ as well as characters $\chi_1$ and $\lambda_{H_1}$ as at the start of \S\ref{sub:local-transfer}. We have a central character datum $(\fkX_{H_1},\chi_{H_1})$ with $\chi_{H_1}:=\chi_1\lambda_{H_1}^{-1}$ as in \S\ref{para:STH1}.
Let us recall relationships between certain stable distributions on $G(\A)$.

\subsection{}
For every real split torus $S$, give $S(\R)^0$ the Lebesgue measure, namely the Haar
measure which, under logarithms, gives $X_*(S)$ covolume $1$.
Use the resulting quotient measures throughout.
Set $\fkX_0:=A_{Z,\infty}$. Define $\chi_0:A_{Z,\infty}\ra\CC^\times$ to be the restriction of $\chi$ to $\fkX_0$.
Then $(\fkX_0,\chi_0)$ is a central character datum.
Let $S_{\chi_0}=S^G_{\chi_0}$ denote\index{Arthur's stable distribution} Arthur's stable distribution on $\cH(G(\A),\chi_0^{-1})$ inductively defined in \cite[\S9]{ArtSTF1}. Write $S_{\disc,\chi_0}$ for the discrete part of $S_{\chi_0}$ (see (7.11) in \textit{loc.~cit.}). Define $S_{\chi}$ and $S_{\disc,\chi}$ in terms of $S_{\chi_0}$ and $S_{\disc,\chi_0}$ exactly as in \eqref{eq:def-of-T_chi}. The equality defining $S^G_{\chi_0}$ inductively leads to the analogous equality (compare with \cite[(3.2.3)]{Arthur})
\begin{equation}\label{eq:S=I-S}
S^G_{\chi}(f)=I^G_\chi(f)-\sum_{\fke=(H,\cH,s,\eta)\in \cE_{\el}(G)\atop H\neq G} \iota(G,H) S^{H_1}_{\chi_{H_1}}(f^{H_1}),
\end{equation}
where $I^G_\chi$ means either $I^G_{\spec,\chi}$ or $I^G_{\geom,\chi}$, which are equal, and the function $f^{H_1}\in \cH(H_1(\A),\chi^{-1}_{H_1})$ denotes a Langlands--Shelstad transfer of $f$ (Proposition \ref{prop:untwisted-transfer}). If $f_\infty$ is stable cuspidal then $f^{H_1}_\infty$ can be chosen stable cuspidal (possibly trivial). Indeed, this is reduced via $z$-extensions to the case that $G_{\der}=G_{\sconn}$, where this fact follows from work of Shelstad and Clozel--Delorme by the argument in \cite[pp.~182--186]{Kot90}. (This argument is also at the basis of constructing $f^{H_1}_\infty$ in \S\ref{h-infty}.) 

Likewise the analogue of \eqref{eq:S=I-S} holds true with $S_{\disc}$ and $I_{\disc}$ in place of $S$ and $I$.
What follows is the stable version of Proposition \ref{prop:simple-TF}.

\begin{lem}\label{lem:S=Sdisc} Let $f=f^\infty f_\infty\in \cH(G(\A),\chi^{-1})$ with $f_\infty$ stable cuspidal. Then
$$S^G_{\chi}(f)=S^G_{\disc,\chi}(f).$$
\end{lem}

\begin{proof}
This follows from Proposition \ref{prop:simple-TF}, which implies that  $I_{\spec,\chi}(f)=I_{\disc,\chi}(f)$
via the inductive definition above.
\end{proof}

\subsection{}

We are going to state a stabilization of the geometric side. Let $(\fkX,\chi)$ be a central character datum for $G$. Write $A_{G_{\R}}$ for the maximal $\R$-split torus in $Z_{G_{\R}}$. (In general $A_{G_{\R}}\neq (A_G)_{\R}$.)
Consider the following hypotheses:
\begin{itemize}
\item[(H1)] $G_{\R}$ contains an elliptic maximal torus,
\item[(H2)] $\fkX=\fkX^\infty\times \fkX_\infty$ with $\fkX^\infty \subset Z(\A_f)$ and $A_{G_\R}(\R)^0\subset \fkX_\infty\subset Z(\R)$.
\end{itemize}
The two conditions are satisfied by the groups contributing to the right hand side of Theorem \ref{thm:end-of-stabilization}, so (H1) and (H2) are harmless to assume for our purpose.

We adapt the definition of the stable distribution $ST^G_M$ in \cite[\S5.4]{Mor10} to the case of fixed central character.
Let $T_{\infty}$ be an elliptic maximal torus of $G_{\R}$ and write $T_{\textup{sc},\infty}\subset G_{\textup{sc},\R}$ for the preimage. Write $G^{\textup{cpt}}$ for an inner form of $G_{\R}$ which is anisotropic modulo $A_{G_{\R}}$; such a  $G^{\textup{cpt}}$ exists by (H1).
The Haar measure on $G^{\textup{cpt}}(\R)$ is always chosen to be compatible with that of $G(\R)$. Equip $\fkX_\infty/A_{G_{\R}}(\R)^0$ with the quotient measure.
Define\index[Oct2026-sorted-notation]{k 1 00 00 k g subscript r@$k(G_{\RR})$} \index[Oct2026-sorted-notation]{v 1 00 20 overline v g subscript r@$\overline v (G_{\RR})$}
\begin{eqnarray}
k(G_{\R})&:=&|\textup{im}(\coh^1(\R,T_{\textup{sc},\infty}) \ra H^1(\R,T_\infty))|, \nonumber\\
\ol{v}(G_{\R})&:=&e(G^{\textup{cpt}}) \vol(G^{\textup{cpt}}(\R)/A_{G_{\R}}(\R)^0). \nonumber
\end{eqnarray}
The two numbers depend only on $G_{\R}$ and a Haar measure on $G(\RR)$. 

Let $M_{\R}\subset G_{\R}$ be an $\R$-rational Levi subgroup containing an elliptic maximal torus. (So the torus is anisotropic modulo $A_{M_{\R}}$.) Let $\Pi$ be a discrete series $L$-packet of $G(\R)$ with fixed central character $\chi_\infty$ on $\fkX_\infty$, and write $\Theta_\Pi$ for the associated stable character (either as a function on regular elements or as a distribution on the space of test functions). Let\index[Oct2026-sorted-notation]{d 0 00 00 d superscript g subscript r subscript m subscript r@$D^{G_{\R}}_{M_{\R}}$} $D^{G_{\R}}_{M_{\R}}$ denote the Weyl discriminant. Write\index[Oct2026-sorted-notation]{zz21 phi 0 00 00 phi superscript g subscript r subscript m subscript r@$\Phi^{G_{\R}}_{M_{\R}}$} $\Phi^{G_{\R}}_{M_{\R}} (\cdot,\Theta_\Pi)$ for the unique continuous function on the set of elliptic elements in $M(\R)$ which extends the function $\gamma \mapsto |D^{G_{\R}}_{M_{\R}}(\gamma)|^{1/2} \Theta_\Pi(\gamma)$ on $M(\R)\cap G(\R)_{\textup{reg}}$; such an extension exists by \cite[Lem.~4.2]{Art89}. 
Let $f_\infty\in \cH(G(\R),\chi_\infty^{-1})$.
For elliptic elements $\gamma$ in $M(\R)$, define\index[Oct2026-sorted-notation]{s 0 00 00 s phi superscript g subscript r subscript m subscript r@$S\Phi^{G_{\R}}_{M_{\R}}$}
$$ S\Phi^{G_{\R}}_{M_{\R}} (\gamma,f_\infty)
:= (-1)^{\dim A_{M_{\R}}/A_{G_{\R}}} \cdot \ol v(M^0_{\R,\gamma})^{-1} \frac{k(M_\R)}{k(G_\R)} \vol(\fkX_\infty/A_{G_{\R}}(\R)^0)
$$
$$\times\sum_{\Pi} \Phi^{G_{\R}}_{M_{\R}}(\gamma^{-1},\Theta_{\Pi}) \Theta_{\Pi}(f_\infty),
$$
where the sum runs over discrete series $L$-packets with central character $\chi_{\infty}$. Set $S\Phi^{G_{\R}}_{M_{\R}} (\gamma,f_\infty)=0$ if $\gamma\in M(\R)$ is not elliptic.

Turning back to the global setting, assuming (H1) and (H2) for $G$, let $M$ be a $\Q$-rational Levi subgroup of $G$, which is said to be\index{G-cuspidal@$G$-cuspidal} \emph{$G$-cuspidal} if $M_{\R}$ contains an elliptic maximal torus and if $\dim A_M/A_G = \dim A_{M_{\R}}/A_{G_{\R}}$.
This relativizes the notion of cuspidal reductive groups over $\Q$. Note that $M=G$ is always $G$-cuspidal even if $G$ is not cuspidal over $\Q$. 
Let $f=f^\infty f_\infty$ with $f^\infty\in \cH(G(\A_f),(\chi^\infty)^{-1})$ and $f_\infty\in \cH(G(\R),\chi_\infty^{-1})$. Denote by $f^\infty_M\in \cH(M(\A_f),(\chi^\infty)^{-1})$ the constant term of $f^\infty$ defined by \cite[(7.13.2)]{GKM97}. (The same definition works regardless of fixed central character. As explained therein, it is not $f^\infty_M$ itself but its orbital integral that has well-defined values.) 

If $M$ is $G$-cuspidal, set $\fkX_M:=\fkX A_{M,\infty}$, $\tau_{\fkX_M}(M):=\vol(M(\Q)\backslash M(\A)/\fkX_M)$, normalized by the Tamagawa measure for $M$, the Haar measure on $\fkX$, and the Lebesgue measure on $A_{M,\infty}$. On $\fkX_M/A_{M_{\R}}(\R)^0\simeq
\fkX/A_{G_{\R}}(\R)^0$ retain the given quotient measure. Put\index[Oct2026-sorted-notation]{st 0 00 00 st superscript g subscript m chi@$ST^G_{M,\chi}$}
$$ST^G_{M,\chi}(f):=\tau_{\fkX_M}(M) \sum_\gamma |\operatorname{Stab}_{\mathfrak X_{\mathbf Q}}(\gamma)|^{-1}
\ol\iota_M(\gamma)^{-1} SO_\gamma(f^\infty_M) S\Phi^{G_{\R}}_{M_{\R}} (\gamma,f_\infty),
$$
where the sum runs over $\fkX_{\Q}$-orbits on the set of stable semi-simple conjugacy classes in $M(\Q)$. (The summand is zero unless $\gamma$ is elliptic in $M(\R)$. The sum can be taken over $(\fkX_M)_{\Q}$-orbits as one can show that $\fkX_{\Q}=(\fkX_M)_{\Q}$.) Define $ST^G_{M,\chi}$ to be identically zero if $M$ is not $G$-cuspidal. Finally, define\index[Oct2026-sorted-notation]{st 0 00 00 st superscript g subscript chi@$ST^G_{\chi}$}
$$ ST^G_{\chi}(f):=\sum_M |(N_G(M)/M)(\Q)|^{-1}
ST^G_{M,\chi}(f),$$
where $M$ runs over  $G(\Q)$-conjugacy classes of $\Q$-rational Levi subgroups.

\begin{lem}\label{lem:S=ST}  Let $G$ be a quasi-split reductive group over $\Q$ with central character datum $(\fkX,\chi)$.
Assume (H1) and (H2) above, and let $f=f^\infty f_\infty$ as above with $f_\infty$ stable cuspidal. 
Then
$$S^G_{\chi}(f)=ST^G_{\chi}(f).$$
\end{lem}

\begin{proof}
Since $Z_G(\R)/\fkX_\infty$ is compact and $f_\infty$ is $K_\infty$-finite, we may assume that $f_\infty\in \cH(G(\R),\omega^{-1})$ for a continuous character $\omega:Z_G(\R)\to\CC^\times$ such that $\omega|_{\fkX_\infty}=\chi_\infty$.
Choose $\tilde f=\tilde f^\infty \tilde f_\infty\in
\cH(G(\A),\chi_0^{-1})$ whose central
averages over $\fkX^\infty$ and $\fkX_\infty/A_{G,\infty}$ are $f^\infty$ and $f_\infty$, respectively. For $z\in \fkX$, write $\tilde f_z$ for the $z$-translate of $\tilde f$.
By \eqref{eq:def-of-T_chi},
\[
S^G_\chi(f)=
\frac{1}{\vol(\fkX_\Q\backslash\fkX/A_{G,\infty})}
\int_{\fkX_\Q\backslash\fkX/A_{G,\infty}}
\chi(z)S^G_{\chi_0}(\tilde f_z)\,dz.
\]
We apply Arthur's stable geometric expansion \cite[Lem.~7.2(b), Thm.~10.1(b), (10.5)]{ArtSTF1} to $S^G_{\chi_0}(\tilde f_z)$
with semisimple descent \cite[(1.7), (2.9)]{ArtSTF2} as well as the stable splitting and descent formulas
\cite[Thms.~6.1(b), 7.1(b)]{ArtTransferWeighted}.
At the real place, take Peng's induced central torus
$Z=A_{G_\R}$ with
$\zeta=\omega|_{A_{G_\R}(\R)}$.
Since $\mathfrak a_Z=\mathfrak a_{G_\R}$,
\cite[Thm.~7.4, (7.6), and the calculation leading to
Thm.~8.1]{Peng}, applied to the stable cuspidal function
$f_\infty$, leave only semisimple terms from
$G$-cuspidal rational Levi subgroups. Grouping the stable conjugacy classes of $M(\Q)$ into $\fkX_\Q$-orbits and unfolding gives
\[
\begin{aligned}
S^G_\chi(f)
={}&
\sum_M
\frac{
	e^G_{M_\R}(M,G_\R)\,
	\vol(M(\Q)\backslash M(\A)/A_{M,\infty})
}{
	|(N_G(M)/M)(\Q)|\,
	\vol(\fkX_\Q\backslash\fkX/A_{G,\infty})
}
\\ \times
&
\sum_\gamma
\frac{1}{
	|\operatorname{Stab}_{\fkX_\Q}(\gamma)|
	\bar\iota_M(\gamma)}
\left(
\int_{\fkX^\infty}\chi^\infty(z)
SO_\gamma((\tilde f^\infty_z)_M)\,dz
\right)
\frac{S^{G_\R}_{M_\R}(\gamma,f_\infty)}{|D^M(\gamma)|_\infty^{1/2}}
,
\end{aligned}
\]
where $M$ runs over $G(\Q)$-conjugacy classes of $G$-cuspidal Levi subgroups, $\gamma$ over $\fkX_\Q$-orbits of stable semisimple classes in $M(\Q)$, $D^M$ is the Weyl discriminant for $M$ (with respect to $M_{\gamma}$), and $S^{G_\R}_{M_\R}(\gamma,\cdot)$ denotes Arthur's normalized stable weighted distribution.
The inverse stabilizer order comes from integration over
$\operatorname{Stab}_{\fkX_\Q}(\gamma)\backslash
\fkX/A_{G,\infty}$.
The discriminant factor comes from the product formula when Arthur's normalized finite orbital integrals are
replaced by our unnormalized stable orbital integrals. The coefficient $e^G_{M_\R}(M,G_\R)$ is that of
\cite[(7.1), (7.3)]{ArtTransferWeighted}, expressed using our quotient measures.

Since $A_{G_\R}(\R)\cap A_{M_\R}(\R)^0
=A_{G_\R}(\R)^0$, passing from quotients by
$A_{G_\R}(\R)$ to quotients by $A_{G_\R}(\R)^0$
multiplies both the packet trace and its signed
compact-volume denominator by
$|A_{G_\R}(\R)/A_{G_\R}(\R)^0|$, leaving their ratio
unchanged. The character formula of Peng, with our quotient measures,
therefore gives
\[
|D^M(\gamma)|_\infty^{-1/2}
S^{G_\R}_{M_\R}(\gamma,f_\infty)
=
S\Phi^{G_\R}_{M_\R}(\gamma,f_\infty).
\]

Finally, Arthur's coefficient in our quotient measures is
\[
e^G_{M_\R}(M,G_\R)=
\frac{\vol(\fkX_\Q\backslash\fkX/A_{G,\infty})}
{\vol((\fkX_M)_\Q\backslash\fkX_M/A_{M,\infty})}.
\]
Thus the coefficient preceding the inner sum becomes
$\tau_{\fkX_M}(M)/|(N_G(M)/M)(\Q)|$.
Since constant terms commute with central averaging
and $\tilde f^\infty$ averages to $f^\infty$, the displayed integral over $\fkX^\infty$ equals $SO_\gamma(f_M^\infty)$.
The resulting expression precisely tells us that $S^G_\chi(f)=ST^G_\chi(f)$.
\end{proof}

\section{Speculations}\label{sub:speculation} 

\subsection{}
We return to the setting of Theorem \ref{thm:end-of-stabilization} for the compactly supported cohomology of Shimura varieties, where (H1) and (H2) hold true for $G$ and $\fkX$ as well as for $H_1$ and $\fkX_{H_1}$ contributing non-trivially to the right hand side.

The analogue of Theorem \ref{thm:end-of-stabilization} is expected to be true for the intersection cohomology of the Baily--Borel compactification. Writing $T^{\textup{IH}}(\Phi^m_{\fkp},f^p dg^p)$ for the intersection cohomology analogue of $T(\Phi^m_{\fkp},f^p dg^p)$, the conjectural stabilization should have the form (cf.~\cite[(10.1)]{Kot90})
\begin{equation}\label{eq:trace-on-IH}
T^{\textup{IH}}(\Phi^m_{\fkp},f^p dg^p)=\sum_{\fke\in \cE_{\el}(G)} \iota(G,H) ST^{H_1}_{\chi_{H_1}}(f^{H_1}).
\end{equation}
The point is that the non-elliptic terms in $ST^{H_1}_{\chi_{H_1}}$ (coming from proper Levi subgroups) should be accounted for exactly by the boundary strata of the Baily--Borel compactification. For non-proper Shimura varieties, \eqref{eq:trace-on-IH} is known for certain special orthogonal groups and unitary similitude groups in addition to general symplectic groups in \cite{Montreal92,Mor08,Mor10,Mor11,Zhu-orthogonal}.
On the other hand, Lemmas \ref{lem:S=Sdisc} and \ref{lem:S=ST} imply that
\begin{equation}\label{eq:ST=S}
ST^{H_1}_{\chi_{H_1}}(f^{H_1}) = S^{H_1}_{\disc,\chi_{H_1}}(f^{H_1}).
\end{equation}
Combined with \eqref{eq:trace-on-IH}, this yields a trace formula for the intersection cohomology in terms of the stable distributions $S^{H_1}_{\disc,\chi_{H_1}}$, which are of a spectral nature. Then one can follow Kottwitz \cite[\S\S9--10]{Kot90} to unravel $S^{H_1}_{\disc,\chi_{H_1}}(f^{H_1})$ to obtain a conjectural description of the intersection cohomology (in each degree, by purity) as a $G(\A_f)\times \Gal(\ol E/E)$-module, in terms of automorphic representations of $G(\A)$ and their endoscopic classification; see p.~201 therein.\footnote{We do not reproduce Kottwitz's argument or his conjectural description here. We content ourselves with remarking that the destabilization process in \cite{Kot90} may also be carried out by applying the conjectural stable multiplicity formula, cf.~\cite[Thm.~4.1.2, (4.8.5)]{Arthur}. Still the key computation at $p$ and $\infty$ of \cite[\S9]{Kot90} is irreplaceable as it reflects the features of test functions at $p$ and $\infty$ specific to the context of Shimura varieties.}
The endoscopic classification for classical groups is worked out in \cite{Arthur,Mok,KMSW,Taibi}, in the quasi-split case and some more. However, little is known for groups of higher rank beyond classical groups, except for partial results on general symplectic and orthogonal groups in \cite{Xu-GSpGSO,Xu-GSpGSO2}.

\subsection{}

We return to compactly supported cohomology. 
In the special case that $G/Z$ is anisotropic over $\Q$, the Shimura variety $\Sh_K$ is proper over $E$ for each $K$. Thus the intersection cohomology coincides with the compactly supported cohomology. 
In particular $T^{\textup{IH}}(\Phi^m_{\fkp},f^p dg^p)=T(\Phi^m_{\fkp},f^p dg^p)$, and the above considerations suggest that
\begin{equation}\label{eq:ST=STell}
ST^{H_1}_{\chi_{H_1}}(f^{H_1}) \stackrel{?}{=} ST^{H_1}_{\el,\chi_{H_1}}(f^{H_1}).
\end{equation}
We stress that this equality is not intrinsic to $H_1$. Indeed, a quasi-split inner form $G^*$ of $G$ over $\Q$ shares the same elliptic endoscopic data as $G$. When $G^*$ can be promoted to a Shimura datum,  \eqref{eq:ST=STell} need not hold for $f^{H_1}$ constructed in the context of the Shimura variety for $G^*$ (since the latter does have a non-empty boundary).
Once \eqref{eq:ST=STell} is verified, we obtain \eqref{eq:trace-on-IH}, and the preceding paragraph explains how to extract the spectral information for the compactly supported cohomology in this case.

If $G/Z$ is isotropic over $\Q$, then the description of the $G(\A_f)\times \Gal(\ol E/E)$-module structure on the compactly supported cohomology is expected to be very complicated. Indeed, this is confirmed by Franke's formula \cite{Fra98} even if the Galois action is forgotten. See also \cite{Lau97} for the case of $\textup{GSp}_4$. Moreover, there may be cancellations between different degrees since the compactly supported cohomology need not be pure. We think that it is better to study the intersection cohomology by proving \eqref{eq:trace-on-IH} in this case.

\subsection{}

We end by summarizing the prospect of unconditional results on the cohomology of Shimura varieties associated with $(G,X)$. In the case of abelian type, our main result is that the identity in Theorem \ref{thm:end-of-stabilization} holds unconditionally whenever $K_p$ is a hyperspecial subgroup of $G(\Q_p)$ and $m$ is sufficiently large. When $G/Z$ is anisotropic over $\Q$, in order to make the conjectural description in the style of \cite[p.~201]{Kot90} unconditional, the two main missing ingredients are 
the endoscopic classification of automorphic representations (for $G$ and the groups $H_1$ contributing to the stabilization) and the equality \eqref{eq:ST=STell}. 
When $G/Z$ is isotropic, instead of \eqref{eq:ST=STell}, one should attempt to prove \eqref{eq:trace-on-IH} by extending the methods of  \cite{Mor10} and  \cite{Zhu-orthogonal}. On top of that, the same endoscopic classification is needed to arrive at the final description of cohomology.

\cleardoublepage

\bibliographystyle{alpha}

\bibliography{myref} 

@article {BKH,
	AUTHOR = {Borovoi, Mikhail and Kaletha, Tasho},
	TITLE = {Galois cohomology of reductive groups over global fields},
	NOTE = {With an appendix by Vladimir Hinich},
	JOURNAL = {Mem. Amer. Math. Soc.},
	FJOURNAL = {Memoirs of the American Mathematical Society},
	VOLUME = {319},
	YEAR = {2026},
	NUMBER = {1621},
	PAGES = {vii+110},
	ISSN = {0065-9266,1947-6221},
	ISBN = {978-1-4704-8046-2; 978-1-4704-8668-6},
	MRCLASS = {11E72 (20G10 20G25 20G30)},
	MRNUMBER = {5072372},
	DOI = {10.1090/memo/1621},
	URL = {https://doi.org/10.1090/memo/1621},
}

@book {Youcis,
	AUTHOR = {Youcis, Alexander Frank},
	TITLE = {The {L}anglands-{K}ottwitz {M}ethod and {D}eformation {S}paces
	of p-{D}ivisible {G}roups of {A}belian {T}ype},
	NOTE = {Thesis (Ph.D.)--University of California, Berkeley},
	PUBLISHER = {ProQuest LLC, Ann Arbor, MI},
	YEAR = {2019},
	PAGES = {192},
	ISBN = {978-1085-79410-7},
	MRCLASS = {99-05},
	MRNUMBER = {4035310},
	URL =
	{http://gateway.proquest.com/openurl?url_ver=Z39.88-2004&rft_val_fmt=info:ofi/fmt:kev:mtx:dissertation&res_dat=xri:pqm&rft_dat=xri:pqdiss:13885192},
}

@article {varshavsky07,
	AUTHOR = {Varshavsky, Yakov},
	TITLE = {Lefschetz-{V}erdier trace formula and a generalization of a
	theorem of {F}ujiwara},
	JOURNAL = {Geom. Funct. Anal.},
	FJOURNAL = {Geometric and Functional Analysis},
	VOLUME = {17},
	YEAR = {2007},
	NUMBER = {1},
	PAGES = {271--319},
	ISSN = {1016-443X,1420-8970},
	MRCLASS = {14F20 (11G25)},
	MRNUMBER = {2306659},
	MRREVIEWER = {Tam\'as\ Szamuely},
	DOI = {10.1007/s00039-007-0596-9},
	URL = {https://doi.org/10.1007/s00039-007-0596-9},
}

@article {Xu,
	AUTHOR = {Xu, Yujie},
	TITLE = {Normalization in integral models of {S}himura varieties of
	abelian type},
	JOURNAL = {Math. Ann.},
	FJOURNAL = {Mathematische Annalen},
	VOLUME = {394},
	YEAR = {2026},
	NUMBER = {4},
	PAGES = {Paper No. 95, 15},
	ISSN = {0025-5831,1432-1807},
	MRCLASS = {14G35},
	MRNUMBER = {5042416},
	DOI = {10.1007/s00208-026-03338-3},
	URL = {https://doi.org/10.1007/s00208-026-03338-3},
}

@article {vanHoftenLR,
	AUTHOR = {van Hoften, Pol},
	TITLE = {Mod {$p$} points on {S}himura varieties of parahoric level},
	JOURNAL = {Forum Math. Pi},
	FJOURNAL = {Forum of Mathematics. Pi},
	VOLUME = {12},
	YEAR = {2024},
	PAGES = {Paper No. e20, 67},
	ISSN = {2050-5086},
	MRCLASS = {11G18 (14G35)},
	MRNUMBER = {4818941},
	MRREVIEWER = {C\'icero\ Carvalho},
	DOI = {10.1017/fmp.2024.22},
	URL = {https://doi.org/10.1017/fmp.2024.22},
}

@article {Lee,
	AUTHOR = {Lee, Dong Uk},
	TITLE = {Lefschetz number formula for {S}himura varieties of {H}odge
	type},
	JOURNAL = {Asian J. Math.},
	FJOURNAL = {Asian Journal of Mathematics},
	VOLUME = {28},
	YEAR = {2024},
	NUMBER = {1},
	PAGES = {103--196},
	ISSN = {1093-6106,1945-0036},
	MRCLASS = {11G18 (11G25 11S40 14G35)},
	MRNUMBER = {4787153},
	DOI = {10.4310/ajm.2024.v28.n1.a5},
	URL = {https://doi.org/10.4310/ajm.2024.v28.n1.a5},
}

@incollection {FontaineGroth,
	AUTHOR = {Fontaine, Jean-Marc},
	TITLE = {Repr\'{e}sentations {$p$}-adiques des corps locaux. {I}},
	BOOKTITLE = {The {G}rothendieck {F}estschrift, {V}ol. {II}},
	SERIES = {Progr. Math.},
	VOLUME = {87},
	PAGES = {249--309},
	PUBLISHER = {Birkh\"{a}user Boston, Boston, MA},
	YEAR = {1990},
	ISBN = {0-8176-3428-2},
	MRCLASS = {11S23 (14F30 14L05)},
	MRNUMBER = {1106901},
	MRREVIEWER = {Rutger\ Noot},
}

@article {LiLiu,
	AUTHOR = {Li, Chao and Liu, Yifeng},
	TITLE = {Chow groups and {$L$}-derivatives of automorphic motives for
	unitary groups},
	JOURNAL = {Ann. of Math. (2)},
	FJOURNAL = {Annals of Mathematics. Second Series},
	VOLUME = {194},
	YEAR = {2021},
	NUMBER = {3},
	PAGES = {817--901},
	ISSN = {0003-486X,1939-8980},
	MRCLASS = {11G18 (11G40 11G50 14C15)},
	MRNUMBER = {4334978},
	MRREVIEWER = {Stefano\ Vigni},
	DOI = {10.4007/annals.2021.194.3.6},
	URL = {https://doi.org/10.4007/annals.2021.194.3.6},
}

@incollection {langlandsmarchen,
	AUTHOR = {Langlands, R. P.},
	TITLE = {Automorphic representations, {S}himura varieties, and motives.
	{E}in {M}\"archen},
	BOOKTITLE = {Automorphic forms, representations and {$L$}-functions
	({P}roc. {S}ympos. {P}ure {M}ath., {O}regon {S}tate {U}niv.,
	{C}orvallis, {O}re., 1977), {P}art 2},
	SERIES = {Proc. Sympos. Pure Math., XXXIII},
	PAGES = {205--246},
	PUBLISHER = {Amer. Math. Soc., Providence, R.I.},
	YEAR = {1979},
	MRCLASS = {12A67 (10D40)},
	MRNUMBER = {546619},
	MRREVIEWER = {R. Howe},
}

@article {KneserII,
	AUTHOR = {Kneser, Martin},
	TITLE = {Galois-{K}ohomologie halbeinfacher algebraischer {G}ruppen
	\"{u}ber {${\mathfrak p}$}-adischen {K}\"{o}rpern. {II}},
	JOURNAL = {Math. Z.},
	FJOURNAL = {Mathematische Zeitschrift},
	VOLUME = {89},
	YEAR = {1965},
	PAGES = {250--272},
	ISSN = {0025-5874},
	MRCLASS = {14.50},
	MRNUMBER = {188219},
	MRREVIEWER = {R. Steinberg},
	DOI = {10.1007/BF02116869},
	URL = {https://doi-org.proxy-um.researchport.umd.edu/10.1007/BF02116869},
}

@article {borovoi83,
	AUTHOR = {Borovo\u{\i}, M. V.},
	TITLE = {Langlands' conjecture concerning conjugation of connected
	{S}himura varieties},
	NOTE = {Selected translations},
	JOURNAL = {Selecta Math. Soviet.},
	FJOURNAL = {Selecta Mathematica Sovietica},
	VOLUME = {3},
	YEAR = {1983/84},
	NUMBER = {1},
	PAGES = {3--39},
	ISSN = {0272-9903},
	MRCLASS = {11G18 (11F55 14K10 14K22 22E40)},
	MRNUMBER = {732450},
	MRREVIEWER = {Thomas Zink},
}

@incollection {keller,
	AUTHOR = {Keller, Bernhard},
	TITLE = {Derived categories and their uses},
	BOOKTITLE = {Handbook of algebra, {V}ol. 1},
	SERIES = {Handb. Algebr.},
	VOLUME = {1},
	PAGES = {671--701},
	PUBLISHER = {Elsevier/North-Holland, Amsterdam},
	YEAR = {1996},
	MRCLASS = {18E30},
	MRNUMBER = {1421815},
	MRREVIEWER = {Jeremy Rickard},
	DOI = {10.1016/S1570-7954(96)80023-4},
	URL = {https://doi-org.proxy-um.researchport.umd.edu/10.1016/S1570-7954(96)80023-4},
}

@inproceedings {Tits,
	AUTHOR = {Tits, J.},
	TITLE = {Reductive groups over local fields},
	BOOKTITLE = {Automorphic forms, representations and {$L$}-functions
	({P}roc. {S}ympos. {P}ure {M}ath., {O}regon {S}tate {U}niv.,
	{C}orvallis, {O}re., 1977), {P}art 1},
	SERIES = {Proc. Sympos. Pure Math., XXXIII},
	PAGES = {29--69},
	PUBLISHER = {Amer. Math. Soc., Providence, R.I.},
	YEAR = {1979},
	MRCLASS = {20G25 (20G10)},
	MRNUMBER = {546588},
}

@article {Shimura63,
	AUTHOR = {Shimura, Goro},
	TITLE = {On analytic families of polarized abelian varieties and
	automorphic functions},
	JOURNAL = {Ann. of Math. (2)},
	FJOURNAL = {Annals of Mathematics. Second Series},
	VOLUME = {78},
	YEAR = {1963},
	PAGES = {149--192},
	ISSN = {0003-486X},
	MRCLASS = {14.51 (14.40)},
	MRNUMBER = {156001},
	MRREVIEWER = {A. Borel},
	DOI = {10.2307/1970507},
	URL = {https://doi-org.proxy-um.researchport.umd.edu/10.2307/1970507},
}

@article {Shimura64,
	AUTHOR = {Shimura, Goro},
	TITLE = {On the field of definition for a field of automorphic
	functions},
	JOURNAL = {Ann. of Math. (2)},
	FJOURNAL = {Annals of Mathematics. Second Series},
	VOLUME = {80},
	YEAR = {1964},
	PAGES = {160--189},
	ISSN = {0003-486X},
	MRCLASS = {10.23 (10.68)},
	MRNUMBER = {167466},
	MRREVIEWER = {Andr\'{e} Weil},
	DOI = {10.2307/1970497},
	URL = {https://doi-org.proxy-um.researchport.umd.edu/10.2307/1970497},
}

@article {Shimura65,
	AUTHOR = {Shimura, Goro},
	TITLE = {On the field of definition for a field of automorphic
	functions. {II}},
	JOURNAL = {Ann. of Math. (2)},
	FJOURNAL = {Annals of Mathematics. Second Series},
	VOLUME = {81},
	YEAR = {1965},
	PAGES = {124--165},
	ISSN = {0003-486X},
	MRCLASS = {10.23 (14.51)},
	MRNUMBER = {171783},
	MRREVIEWER = {Andr\'{e} Weil},
	DOI = {10.2307/1970388},
	URL = {https://doi-org.proxy-um.researchport.umd.edu/10.2307/1970388},
}

@article {Shimura66,
	AUTHOR = {Shimura, Goro},
	TITLE = {On the field of definition for a field of automorphic
	functions. {III}},
	JOURNAL = {Ann. of Math. (2)},
	FJOURNAL = {Annals of Mathematics. Second Series},
	VOLUME = {83},
	YEAR = {1966},
	PAGES = {377--385},
	ISSN = {0003-486X},
	MRCLASS = {10.68 (10.20)},
	MRNUMBER = {190100},
	MRREVIEWER = {M. Eichler},
	DOI = {10.2307/1970437},
	URL = {https://doi-org.proxy-um.researchport.umd.edu/10.2307/1970437},
}

@article {Shimura67,
	AUTHOR = {Shimura, Goro},
	TITLE = {Construction of class fields and zeta functions of algebraic
	curves},
	JOURNAL = {Ann. of Math. (2)},
	FJOURNAL = {Annals of Mathematics. Second Series},
	VOLUME = {85},
	YEAR = {1967},
	PAGES = {58--159},
	ISSN = {0003-486X},
	MRCLASS = {14.49 (10.68)},
	MRNUMBER = {204426},
	MRREVIEWER = {M. Eichler},
	DOI = {10.2307/1970526},
	URL = {https://doi-org.proxy-um.researchport.umd.edu/10.2307/1970526},
}

@article {Shimura67b,
	AUTHOR = {Shimura, Goro},
	TITLE = {Algebraic number fields and symplectic discontinuous groups},
	JOURNAL = {Ann. of Math. (2)},
	FJOURNAL = {Annals of Mathematics. Second Series},
	VOLUME = {86},
	YEAR = {1967},
	PAGES = {503--592},
	ISSN = {0003-486X},
	MRCLASS = {10.65 (14.00)},
	MRNUMBER = {222048},
	MRREVIEWER = {M. Eichler},
	DOI = {10.2307/1970613},
	URL = {https://doi-org.proxy-um.researchport.umd.edu/10.2307/1970613},
}

@article {Shimura70,
	AUTHOR = {Shimura, Goro},
	TITLE = {On canonical models of arithmetic quotients of bounded
	symmetric domains},
	JOURNAL = {Ann. of Math. (2)},
	FJOURNAL = {Annals of Mathematics. Second Series},
	VOLUME = {91},
	YEAR = {1970},
	PAGES = {144--222},
	ISSN = {0003-486X},
	MRCLASS = {10.68},
	MRNUMBER = {257031},
	MRREVIEWER = {K. Iyanaga},
	DOI = {10.2307/1970604},
	URL = {https://doi-org.proxy-um.researchport.umd.edu/10.2307/1970604},
}

@article {Shimura70b,
	AUTHOR = {Shimura, Goro},
	TITLE = {On canonical models of arithmetic quotients of bounded
	symmetric domains. {II}},
	JOURNAL = {Ann. of Math. (2)},
	FJOURNAL = {Annals of Mathematics. Second Series},
	VOLUME = {92},
	YEAR = {1970},
	PAGES = {528--549},
	ISSN = {0003-486X},
	MRCLASS = {10D05},
	MRNUMBER = {292758},
	MRREVIEWER = {I. Piatetski-Shapiro},
	DOI = {10.2307/1970631},
	URL = {https://doi-org.proxy-um.researchport.umd.edu/10.2307/1970631},
}

@book {Springerbook,
	AUTHOR = {Springer, T. A.},
	TITLE = {Linear algebraic groups},
	SERIES = {Progress in Mathematics},
	VOLUME = {9},
	EDITION = {Second},
	PUBLISHER = {Birkh\"{a}user Boston, Inc., Boston, MA},
	YEAR = {1998},
	PAGES = {xiv+334},
	ISBN = {0-8176-4021-5},
	MRCLASS = {20G15 (14L10)},
	MRNUMBER = {1642713},
	MRREVIEWER = {Andy R. Magid},
	DOI = {10.1007/978-0-8176-4840-4},
	URL = {https://doi-org.proxy-um.researchport.umd.edu/10.1007/978-0-8176-4840-4},
}

@book {AT,
	AUTHOR = {Artin, Emil and Tate, John},
	TITLE = {Class field theory},
	NOTE = {Reprinted with corrections from the 1967 original},
	PUBLISHER = {AMS Chelsea Publishing, Providence, RI},
	YEAR = {2009},
	PAGES = {viii+194},
	ISBN = {978-0-8218-4426-7},
	MRCLASS = {11-01 (11R37)},
	MRNUMBER = {2467155},
	DOI = {10.1090/chel/366},
	URL = {https://doi-org.proxy-um.researchport.umd.edu/10.1090/chel/366},
}

@incollection{deligne1982motifs,
	title={Motifs et groupes de {T}aniyama},
	author={Deligne, Pierre},
	booktitle={Hodge cycles, motives, and Shimura varieties},
	pages={261--279},
	year={1982},
	publisher={Springer}
}

@inproceedings {LanglandsAntwerp,
	AUTHOR = {Langlands, R. P.},
	TITLE = {Modular forms and {$\ell $}-adic representations},
	BOOKTITLE = {Modular functions of one variable, {II} ({P}roc. {I}nternat.
	{S}ummer {S}chool, {U}niv. {A}ntwerp, {A}ntwerp, 1972)},
	PAGES = {361--500. Lecture Notes in Math., Vol. 349},
	YEAR = {1973},
	MRCLASS = {12A70 (10D15 14G10 22E50)},
	MRNUMBER = {0354617},
	MRREVIEWER = {I. I. Pjateckii-Sapiro},
}

@article {Liucompatibility,
	AUTHOR = {Liu, Tong},
	TITLE = {Compatibility of {K}isin modules for different uniformizers},
	JOURNAL = {J. Reine Angew. Math.},
	FJOURNAL = {Journal f\"{u}r die Reine und Angewandte Mathematik. [Crelle's
	Journal]},
	VOLUME = {740},
	YEAR = {2018},
	PAGES = {1--24},
	ISSN = {0075-4102},
	MRCLASS = {11F80 (14F30)},
	MRNUMBER = {3824780},
	MRREVIEWER = {Takuya Yamauchi},
	DOI = {10.1515/crelle-2015-0074},
	URL = {https://doi-org.proxy-um.researchport.umd.edu/10.1515/crelle-2015-0074},
}

@article {BMS,
	AUTHOR = {Bhatt, Bhargav and Morrow, Matthew and Scholze, Peter},
	TITLE = {Integral {$p$}-adic {H}odge theory},
	JOURNAL = {Publ. Math. Inst. Hautes \'{E}tudes Sci.},
	FJOURNAL = {Publications Math\'{e}matiques. Institut de Hautes \'{E}tudes
	Scientifiques},
	VOLUME = {128},
	YEAR = {2018},
	PAGES = {219--397},
	ISSN = {0073-8301},
	MRCLASS = {14F30},
	MRNUMBER = {3905467},
	DOI = {10.1007/s10240-019-00102-z},
	URL = {https://doi-org.proxy-um.researchport.umd.edu/10.1007/s10240-019-00102-z},
}

@article {steinbergtorsion,
	AUTHOR = {Steinberg, Robert},
	TITLE = {Torsion in reductive groups},
	JOURNAL = {Advances in Math.},
	FJOURNAL = {Advances in Mathematics},
	VOLUME = {15},
	YEAR = {1975},
	PAGES = {63--92},
	ISSN = {0001-8708},
	MRCLASS = {20G15},
	MRNUMBER = {354892},
	MRREVIEWER = {S. I. Gel\cprime fand},
	DOI = {10.1016/0001-8708(75)90125-5},
	URL = {https://doi-org.proxy-um.researchport.umd.edu/10.1016/0001-8708(75)90125-5},
}

@article {ZZ,
	AUTHOR = {Zhou, Rong and Zhu, Yihang},
	TITLE = {Twisted orbital integrals and irreducible components of affine
	{D}eligne-{L}usztig varieties},
	JOURNAL = {Camb. J. Math.},
	FJOURNAL = {Cambridge Journal of Mathematics},
	VOLUME = {8},
	YEAR = {2020},
	NUMBER = {1},
	PAGES = {149--241},
	ISSN = {2168-0930},
	MRCLASS = {11G18 (22E35)},
	MRNUMBER = {4085434},
	DOI = {10.4310/CJM.2020.v8.n1.a3},
	URL = {https://doi-org.proxy-um.researchport.umd.edu/10.4310/CJM.2020.v8.n1.a3},
}

@article{conradlifting,
	title={Lifting global representations with local properties},
	author={Conrad, Brian},
	journal  = {preprint},
	note = {Available at \url{http://math.stanford.edu/~conrad/papers/locchar.pdf}},
	year = {2011}
}

@book {harristaylor,
	AUTHOR = {Harris, Michael and Taylor, Richard},
	TITLE = {The geometry and cohomology of some simple {S}himura
	varieties},
	SERIES = {Annals of Mathematics Studies},
	VOLUME = {151},
	NOTE = {With an appendix by Vladimir G. Berkovich},
	PUBLISHER = {Princeton University Press, Princeton, NJ},
	YEAR = {2001},
	PAGES = {viii+276},
	ISBN = {0-691-09090-4},
	MRCLASS = {11G18 (11F70 11S37 14G35 22E45)},
	MRNUMBER = {1876802},
	MRREVIEWER = {James Milne},
}

@article {Wintenberger97,
	AUTHOR = {Wintenberger, Jean-Pierre},
	TITLE = {Propri\'{e}t\'{e}s du groupe tannakien des structures de {H}odge
	{$p$}-adiques et torseur entre cohomologies cristalline et
	\'{e}tale},
	JOURNAL = {Ann. Inst. Fourier (Grenoble)},
	FJOURNAL = {Universit\'{e} de Grenoble. Annales de l'Institut Fourier},
	VOLUME = {47},
	YEAR = {1997},
	NUMBER = {5},
	PAGES = {1289--1334},
	ISSN = {0373-0956},
	MRCLASS = {14F30 (14D07 14F20)},
	MRNUMBER = {1600395},
	MRREVIEWER = {Claudio Pedrini},
	URL = {http://www.numdam.org/item?id=AIF_1997__47_5_1289_0},
}

@book {Saavedra,
	AUTHOR = {Saavedra Rivano, Neantro},
	TITLE = {Cat\'{e}gories {T}annakiennes},
	SERIES = {Lecture Notes in Mathematics, Vol. 265},
	PUBLISHER = {Springer-Verlag, Berlin-New York},
	YEAR = {1972},
	PAGES = {ii+418},
	MRCLASS = {14L15 (18D10 20G05)},
	MRNUMBER = {0338002},
	MRREVIEWER = {P. Abellanas},
}

@incollection {kisin2006crystalline,
	AUTHOR = {Kisin, Mark},
	TITLE = {Crystalline representations and {$F$}-crystals},
	BOOKTITLE = {Algebraic geometry and number theory},
	SERIES = {Progr. Math.},
	VOLUME = {253},
	PAGES = {459--496},
	PUBLISHER = {Birkh\"{a}user Boston, Boston, MA},
	YEAR = {2006},
	MRCLASS = {11S20 (14F30)},
	MRNUMBER = {2263197},
	MRREVIEWER = {Martin C. Olsson},
	DOI = {10.1007/978-0-8176-4532-8_7},
	URL = {https://doi-org.proxy-um.researchport.umd.edu/10.1007/978-0-8176-4532-8_7},
}

@incollection{deligne1982tannakian,
	title={Tannakian categories},
	author={Deligne, Pierre and Milne, James S},
	booktitle={Hodge cycles, motives, and Shimura varieties},
	pages={101--228},
	year={1982},
	publisher={Springer}
}

@article {wedhorn2004tannakian,
	AUTHOR = {Wedhorn, Torsten},
	TITLE = {On {T}annakian duality over valuation rings},
	JOURNAL = {J. Algebra},
	FJOURNAL = {Journal of Algebra},
	VOLUME = {282},
	YEAR = {2004},
	NUMBER = {2},
	PAGES = {575--609},
	ISSN = {0021-8693},
	MRCLASS = {18D10 (16W30 18E15 20G35)},
	MRNUMBER = {2101076},
	MRREVIEWER = {Ulrich G\"{o}rtz},
	DOI = {10.1016/j.jalgebra.2004.07.024},
	URL = {https://doi-org.proxy-um.researchport.umd.edu/10.1016/j.jalgebra.2004.07.024},
}

@misc{milneAGS,
	author={Milne, J. S.},
	title={Basic Theory of Affine Group Schemes},
	year={2012},
	note={Available at \url{https://www.jmilne.org/math/CourseNotes/AGS.pdf}}
}

@misc{milnelettertodeligne,
	author={Milne, J. S.},
	title={Letter to {D}eligne},
	year={1990},
	note={Available at \url{https://www.jmilne.org/math/articles/1990b.pdf}}
}

@article {Serre68,
	AUTHOR = {Serre, Jean-Pierre},
	TITLE = {Groupes de {G}rothendieck des sch\'{e}mas en groupes r\'{e}ductifs
	d\'{e}ploy\'{e}s},
	JOURNAL = {Inst. Hautes \'{E}tudes Sci. Publ. Math.},
	FJOURNAL = {Institut des Hautes \'{E}tudes Scientifiques. Publications
	Math\'{e}matiques},
	volume = {34},
	YEAR = {1968},
	PAGES = {37--52},
	ISSN = {0073-8301},
	MRCLASS = {14.55},
	MRNUMBER = {231831},
	MRREVIEWER = {T. Oda},
	URL = {http://www.numdam.org/item?id=PMIHES_1968__34__37_0},
}

@article {Bor63,
	AUTHOR = {Borel, Armand},
	TITLE = {Some finiteness properties of adele groups over number fields},
	JOURNAL = {Inst. Hautes \'{E}tudes Sci. Publ. Math.},
	FJOURNAL = {Institut des Hautes \'{E}tudes Scientifiques. Publications
	Math\'{e}matiques},
	volume = {16},
	YEAR = {1963},
	PAGES = {5--30},
	ISSN = {0073-8301},
	MRCLASS = {20.75 (10.25)},
	MRNUMBER = {202718},
	MRREVIEWER = {T. Ono},
	URL = {http://www.numdam.org/item?id=PMIHES_1963__16__5_0},
}

@article{HainesRapoport,
	title={On parahoric subgroups},
	author={Haines, Thomas and Rapoport, Michael},
	journal={Advances in Mathematics},
	volume={219},
	number={1},
	pages={188--198},
	year={2008},
	publisher={Elsevier}
}

@book {HartsRD,
	AUTHOR = {Hartshorne, Robin},
	TITLE = {Residues and duality},
	SERIES = {Lecture notes of a seminar on the work of A. Grothendieck,
	given at Harvard 1963/64. With an appendix by P. Deligne.
	Lecture Notes in Mathematics, No. 20},
	PUBLISHER = {Springer-Verlag, Berlin-New York},
	YEAR = {1966},
	PAGES = {vii+423},
	MRCLASS = {14.55},
	MRNUMBER = {0222093},
	MRREVIEWER = {R. L. Knighten},
}

@incollection {Mil83,
	AUTHOR = {Milne, J. S.},
	TITLE = {The action of an automorphism of {${\mathbf C}$} on a {S}himura
	variety and its special points},
	BOOKTITLE = {Arithmetic and geometry, {V}ol. {I}},
	SERIES = {Progr. Math.},
	VOLUME = {35},
	PAGES = {239--265},
	PUBLISHER = {Birkh\"{a}user Boston, Boston, MA},
	YEAR = {1983},
	MRCLASS = {11G18 (14K15 14K22)},
	MRNUMBER = {717596},
	MRREVIEWER = {Paul B. Garrett},
}

@incollection {Mil88,
	AUTHOR = {Milne, J. S.},
	TITLE = {Canonical models of (mixed) {S}himura varieties and
	automorphic vector bundles},
	BOOKTITLE = {Automorphic forms, {S}himura varieties, and {$L$}-functions,
	{V}ol. {I} ({A}nn {A}rbor, {MI}, 1988)},
	SERIES = {Perspect. Math.},
	VOLUME = {10},
	PAGES = {283--414},
	PUBLISHER = {Academic Press, Boston, MA},
	YEAR = {1990},
	MRCLASS = {11G18 (11F03 11F12 14K22)},
	MRNUMBER = {1044823},
	MRREVIEWER = {Min Ho Lee},
}

@incollection {lanstrohII,
	AUTHOR = {Lan, Kai-Wen and Stroh, Beno\^{\i}t},
	TITLE = {Nearby cycles of automorphic \'{e}tale sheaves, {II}},
	BOOKTITLE = {Cohomology of arithmetic groups},
	SERIES = {Springer Proc. Math. Stat.},
	VOLUME = {245},
	PAGES = {83--106},
	PUBLISHER = {Springer, Cham},
	YEAR = {2018},
	MRCLASS = {11G18 (11F75 11G15 14F20 14G35)},
	MRNUMBER = {3848816},
	MRREVIEWER = {Giovanni Rosso},
}

@article {RV14,
	AUTHOR = {Rapoport, Michael and Viehmann, Eva},
	TITLE = {Towards a theory of local {S}himura varieties},
	JOURNAL = {M\"unster J. Math.},
	FJOURNAL = {M\"unster Journal of Mathematics},
	VOLUME = {7},
	YEAR = {2014},
	PAGES = {273--326},
	ISSN = {1867-5778},
	MRCLASS = {11G18 (14G35)},
	MRNUMBER = {3271247},
}

@article {debacker2006,
	AUTHOR = {DeBacker, Stephen},
	TITLE = {Parameterizing conjugacy classes of maximal unramified tori
	via {B}ruhat-{T}its theory},
	JOURNAL = {Michigan Math. J.},
	FJOURNAL = {Michigan Mathematical Journal},
	VOLUME = {54},
	YEAR = {2006},
	NUMBER = {1},
	PAGES = {157--178},
	ISSN = {0026-2285},
	MRCLASS = {22E35},
	MRNUMBER = {2214792},
	MRREVIEWER = {David A. Renard},
	DOI = {10.1307/mmj/1144437442},
	URL = {https://doi-org.proxy-um.researchport.umd.edu/10.1307/mmj/1144437442},
}

@article {CIMdecomp,
	AUTHOR = {Alper, Jarod and Heinloth, Jochen and Halpern-Leistner,
	Daniel},
	TITLE = {Cartan-{I}wahori-{M}atsumoto decompositions for reductive
	groups},
	JOURNAL = {Pure Appl. Math. Q.},
	FJOURNAL = {Pure and Applied Mathematics Quarterly},
	VOLUME = {17},
	YEAR = {2021},
	NUMBER = {2},
	PAGES = {593--604},
	ISSN = {1558-8599,1558-8602},
	MRCLASS = {14L24 (13A50 14L35)},
	MRNUMBER = {4257594},
	MRREVIEWER = {P.\ E.\ Newstead},
	DOI = {10.4310/PAMQ.2021.v17.n2.a1},
	URL = {https://doi.org/10.4310/PAMQ.2021.v17.n2.a1},
}

@book {SGA1,
	author = {Grothendieck, Alexander},
	TITLE = {Rev\^{e}tements \'{e}tales et groupe fondamental ({SGA} 1)},
	SERIES = {Documents Math\'{e}matiques (Paris)},
	VOLUME = {3},
	NOTE = {S\'{e}minaire de g\'{e}om\'{e}trie alg\'{e}brique du Bois Marie 1960--61. 
	Directed by A. Grothendieck,
	With two papers by M. Raynaud,
	Updated and annotated reprint of the 1971 original [Lecture
	Notes in Math., 224, Springer]},
	PUBLISHER = {Soci\'{e}t\'{e} Math\'{e}matique de France, Paris},
	YEAR = {2003},
	PAGES = {xviii+327},
	ISBN = {2-85629-141-4},
	MRCLASS = {14E20 (14-06 14F35)},
	MRNUMBER = {2017446},
}

@incollection {deligneCT,
	AUTHOR = {Deligne, P.},
	TITLE = {Cat\'{e}gories tannakiennes},
	BOOKTITLE = {The {G}rothendieck {F}estschrift, {V}ol. {II}},
	SERIES = {Progr. Math.},
	VOLUME = {87},
	PAGES = {111--195},
	PUBLISHER = {Birkh\"{a}user Boston, Boston, MA},
	YEAR = {1990},
	MRCLASS = {14A99 (12H05 18A99)},
	MRNUMBER = {1106898},
	MRREVIEWER = {James Milne},
}

@article {KisinPappas,
	AUTHOR = {Kisin, M. and Pappas, G.},
	TITLE = {Integral models of {S}himura varieties with parahoric level
	structure},
	JOURNAL = {Publ. Math. Inst. Hautes \'{E}tudes Sci.},
	FJOURNAL = {Publications Math\'{e}matiques. Institut de Hautes \'{E}tudes
	Scientifiques},
	VOLUME = {128},
	YEAR = {2018},
	PAGES = {121--218},
	ISSN = {0073-8301},
	MRCLASS = {11G18},
	MRNUMBER = {3905466},
	MRREVIEWER = {Martin Orr},
	DOI = {10.1007/s10240-018-0100-0},
	URL = {https://doi-org.proxy-um.researchport.umd.edu/10.1007/s10240-018-0100-0},
}

@article {broshi2013,
	AUTHOR = {Broshi, Michael},
	TITLE = {{$G$}-torsors over a {D}edekind scheme},
	JOURNAL = {J. Pure Appl. Algebra},
	FJOURNAL = {Journal of Pure and Applied Algebra},
	VOLUME = {217},
	YEAR = {2013},
	NUMBER = {1},
	PAGES = {11--19},
	ISSN = {0022-4049},
	MRCLASS = {14L15},
	MRNUMBER = {2965898},
	MRREVIEWER = {Andy R. Magid},
	DOI = {10.1016/j.jpaa.2012.01.011},
	URL = {https://doi-org.proxy-um.researchport.umd.edu/10.1016/j.jpaa.2012.01.011},
}

@article {Anschutz,
	AUTHOR = {Ansch\"utz, Johannes},
	TITLE = {Extending torsors on the punctured {${\mathrm{ Spec}}(A_{\inf})$}},
	JOURNAL = {J. Reine Angew. Math.},
	FJOURNAL = {Journal f\"ur die Reine und Angewandte Mathematik. [Crelle's
	Journal]},
	VOLUME = {783},
	YEAR = {2022},
	PAGES = {227--268},
	ISSN = {0075-4102,1435-5345},
	MRCLASS = {20G05 (14M15)},
	MRNUMBER = {4373246},
	MRREVIEWER = {Nguy\cftil en Qu\^oc Th\'ang},
	DOI = {10.1515/crelle-2021-0077},
	URL = {https://doi.org/10.1515/crelle-2021-0077},
}

@article {rapoportguide,
	AUTHOR = {Rapoport, Michael},
	TITLE = {A guide to the reduction modulo {$p$} of {S}himura varieties},
	NOTE = {Automorphic forms. I},
	JOURNAL = {Ast\'{e}risque},
	FJOURNAL = {Ast\'{e}risque},
	volume = {298},
	YEAR = {2005},
	PAGES = {271--318},
	ISSN = {0303-1179},
	MRCLASS = {11G18 (11G40 14G35)},
	MRNUMBER = {2141705},
	MRREVIEWER = {Ulrich G\"{o}rtz},
}

@article {Greenberg,
	AUTHOR = {Greenberg, Marvin J.},
	TITLE = {Schemata over local rings. {II}},
	JOURNAL = {Ann. of Math.},
	FJOURNAL = {Annals of Mathematics. Second Series},
	VOLUME = {78},
	NUMBER = {2},  
	YEAR = {1963},
	PAGES = {256--266},
	ISSN = {0003-486X},
	MRCLASS = {14.55 (14.08)},
	MRNUMBER = {0156855},
	MRREVIEWER = {M. Artin},
	DOI = {10.2307/1970342},
	URL = {https://doi.org/10.2307/1970342},
}

@incollection {Fontaine79,
	AUTHOR = {Fontaine, Jean-Marc},
	TITLE = {Modules galoisiens, modules filtr\'{e}s et anneaux de
	{B}arsotti-{T}ate},
	BOOKTITLE = {Journ\'{e}es de {G}\'{e}om\'{e}trie {A}lg\'{e}brique de {R}ennes. ({R}ennes,
	1978), {V}ol. {III}},
	SERIES = {Ast\'{e}risque},
	VOLUME = {65},
	PAGES = {3--80},
	PUBLISHER = {Soc. Math. France, Paris},
	YEAR = {1979},
	MRCLASS = {14L05 (12B20 14G20)},
	MRNUMBER = {563472},
	MRREVIEWER = {Lucile B\'{e}gueri},
}

@article {CF,
	AUTHOR = {Colmez, Pierre and Fontaine, Jean-Marc},
	TITLE = {Construction des repr\'{e}sentations {$p$}-adiques semi-stables},
	JOURNAL = {Invent. Math.},
	FJOURNAL = {Inventiones Mathematicae},
	VOLUME = {140},
	YEAR = {2000},
	NUMBER = {1},
	PAGES = {1--43},
	ISSN = {0020-9910},
	MRCLASS = {11S23 (14F30 14F40)},
	MRNUMBER = {1779803},
	MRREVIEWER = {Mark Kisin},
	DOI = {10.1007/s002220000042},
	URL = {https://doi-org.ezproxy.cul.columbia.edu/10.1007/s002220000042},
}

@incollection {Serre79,
	AUTHOR = {Serre, Jean-Pierre},
	TITLE = {Groupes alg\'{e}briques associ\'{e}s aux modules de {H}odge-{T}ate},
	BOOKTITLE = {Journ\'{e}es de {G}\'{e}om\'{e}trie {A}lg\'{e}brique de {R}ennes. ({R}ennes,
	1978), {V}ol. {III}},
	SERIES = {Ast\'{e}risque},
	VOLUME = {65},
	PAGES = {155--188},
	PUBLISHER = {Soc. Math. France, Paris},
	YEAR = {1979},
	MRCLASS = {14L05 (12B25 14G20)},
	MRNUMBER = {563476},
	MRREVIEWER = {Kenneth A. Ribet},
}

@article {CKV,
	AUTHOR = {Chen, Miaofen and Kisin, Mark and Viehmann, Eva},
	TITLE = {Connected components of affine {D}eligne-{L}usztig varieties
	in mixed characteristic},
	JOURNAL = {Compos. Math.},
	FJOURNAL = {Compositio Mathematica},
	VOLUME = {151},
	YEAR = {2015},
	NUMBER = {9},
	PAGES = {1697--1762},
	ISSN = {0010-437X},
	MRCLASS = {14G35 (20G25)},
	MRNUMBER = {3406443},
	MRREVIEWER = {Zongbin Chen},
	DOI = {10.1112/S0010437X15007253},
	URL = {https://doi.org/10.1112/S0010437X15007253},
}

@article {Madapusi19,
	AUTHOR = {Madapusi Pera, Keerthi},
	TITLE = {Toroidal compactifications of integral models of {S}himura
	varieties of {H}odge type},
	JOURNAL = {Ann. Sci. \'{E}c. Norm. Sup\'{e}r. (4)},
	FJOURNAL = {Annales Scientifiques de l'\'{E}cole Normale Sup\'{e}rieure.
	Quatri\`eme S\'{e}rie},
	VOLUME = {52},
	YEAR = {2019},
	NUMBER = {2},
	PAGES = {393--514},
	ISSN = {0012-9593,1873-2151},
	MRCLASS = {11G18 (14G35)},
	MRNUMBER = {3948111},
	MRREVIEWER = {Brandon\ Levin},
	DOI = {10.24033/asens.2391},
	URL = {https://doi.org/10.24033/asens.2391},
}

@incollection {milne92,
	AUTHOR = {Milne, J. S.},
	TITLE = {The points on a {S}himura variety modulo a prime of good
	reduction},
	BOOKTITLE = {The zeta functions of {P}icard modular surfaces},
	PAGES = {151--253},
	PUBLISHER = {Univ. Montr\'eal, Montreal, QC},
	YEAR = {1992},
	MRCLASS = {11G18 (14G35)},
	MRNUMBER = {1155229},
	MRREVIEWER = {Takeshi Ooe},
}

@article {pinkcalc,
	AUTHOR = {Pink, Richard},
	TITLE = {On the calculation of local terms in the {L}efschetz-{V}erdier
	trace formula and its application to a conjecture of
	{D}eligne},
	JOURNAL = {Ann. of Math. (2)},
	FJOURNAL = {Annals of Mathematics. Second Series},
	VOLUME = {135},
	YEAR = {1992},
	NUMBER = {3},
	PAGES = {483--525},
	ISSN = {0003-486X},
	MRCLASS = {14F20 (11G25 14G15)},
	MRNUMBER = {1166642},
	MRREVIEWER = {Takeshi Ooe},
	DOI = {10.2307/2946574},
	URL = {http://dx.doi.org/10.2307/2946574},
}

@article {varshavsky,
	AUTHOR = {Varshavsky, Yakov},
	TITLE = {A proof of a generalization of {D}eligne's conjecture},
	JOURNAL = {Electron. Res. Announc. Amer. Math. Soc.},
	FJOURNAL = {Electronic Research Announcements of the American Mathematical
	Society},
	VOLUME = {11},
	YEAR = {2005},
	PAGES = {78--88},
	ISSN = {1079-6762},
	MRCLASS = {14F20 (11G25)},
	MRNUMBER = {2176068},
	MRREVIEWER = {Kanetomo Sato},
	DOI = {10.1090/S1079-6762-05-00150-2},
	URL = {http://dx.doi.org/10.1090/S1079-6762-05-00150-2},
}

@article {fujiwara,
	AUTHOR = {Fujiwara, Kazuhiro},
	TITLE = {Rigid geometry, {L}efschetz-{V}erdier trace formula and
	{D}eligne's conjecture},
	JOURNAL = {Invent. Math.},
	FJOURNAL = {Inventiones Mathematicae},
	VOLUME = {127},
	YEAR = {1997},
	NUMBER = {3},
	PAGES = {489--533},
	ISSN = {0020-9910},
	MRCLASS = {14F20 (14G20)},
	MRNUMBER = {1431137},
	DOI = {10.1007/s002220050129},
	URL = {http://dx.doi.org/10.1007/s002220050129},
}

@BOOK{RZ96,
	title = {Period spaces for {$p$}-divisible groups},
	publisher = {Princeton University Press, Princeton, NJ},
	year = {1996},
	author = {Rapoport, M. and Zink, Th.},
	volume = {141},
	pages = {xxii+324},
	series = {Annals of Mathematics Studies},
	doi = {10.1515/9781400882601},
	isbn = {0-691-02782-X; 0-691-02781-1},
	mrclass = {14G20 (11G18 14F30 14L05 14M15 20G05 20G25)},
	mrnumber = {1393439},
	mrreviewer = {Robert E. Kottwitz},
	url = {http://dx.doi.org/10.1515/9781400882601}
}

@book {borel1991,
    AUTHOR = {Borel, Armand},
     TITLE = {Linear algebraic groups},
    SERIES = {Graduate Texts in Mathematics},
    VOLUME = {126},
   EDITION = {Second},
 PUBLISHER = {Springer-Verlag},
   ADDRESS = {New York},
      YEAR = {1991},
     PAGES = {xii+288},
      ISBN = {0-387-97370-2},
   MRCLASS = {20-01 (20Gxx)},
  MRNUMBER = {1102012 (92d:20001)},
MRREVIEWER = {F. D. Veldkamp},
       DOI = {10.1007/978-1-4612-0941-6},
       URL = {http://dx.doi.org/10.1007/978-1-4612-0941-6},
}

@article {borovoi,
    AUTHOR = {Borovoi, Mikhail},
     TITLE = {Abelian {G}alois cohomology of reductive groups},
   JOURNAL = {Mem. Amer. Math. Soc.},
  FJOURNAL = {Memoirs of the American Mathematical Society},
    VOLUME = {132},
      YEAR = {1998},
    NUMBER = {626},
     PAGES = {viii+50},
      ISSN = {0065-9266},
     CODEN = {MAMCAU},
   MRCLASS = {20G10 (11E72 14E20 18G50)},
  MRNUMBER = {1401491 (98j:20061)},
MRREVIEWER = {James E. Humphreys},
       DOI = {10.1090/memo/0626},
       URL = {http://dx.doi.org/10.1090/memo/0626},
}

@article {Lab99,
    AUTHOR = {Labesse, Jean-Pierre},
     TITLE = {Cohomologie, stabilisation et changement de base},
      NOTE = {Appendix A by Laurent Clozel and Labesse, and Appendix B by
              Lawrence Breen},
   JOURNAL = {Ast\'erisque},
  FJOURNAL = {Ast\'erisque},
    volume = {257},
      YEAR = {1999},
     PAGES = {vi+161},
      ISSN = {0303-1179},
   MRCLASS = {22E55 (11F70 18G30)},
  MRNUMBER = {1695940 (2001k:22040)},
MRREVIEWER = {Nguy{\cftil{e}}n Qu{\^o}c Th{\'a}ng},
}

@article {kottwitztwisted,
    AUTHOR = {Kottwitz, Robert E.},
     TITLE = {Shimura varieties and twisted orbital integrals},
   JOURNAL = {Math. Ann.},
  FJOURNAL = {Mathematische Annalen},
    VOLUME = {269},
      YEAR = {1984},
    NUMBER = {3},
     PAGES = {287--300},
      ISSN = {0025-5831},
     CODEN = {MAANA},
   MRCLASS = {11G18 (11F70 11G40 22E50)},
  MRNUMBER = {761308 (87b:11047)},
MRREVIEWER = {Ernst-Wilhelm Zink},
       DOI = {10.1007/BF01450697},
       URL = {http://dx.doi.org/10.1007/BF01450697},
}

@article {kisin2012modp,
	AUTHOR = {Kisin, Mark},
	TITLE = {Mod {$p$} points on {S}himura varieties of abelian type},
	JOURNAL = {J. Amer. Math. Soc.},
	FJOURNAL = {Journal of the American Mathematical Society},
	VOLUME = {30},
	YEAR = {2017},
	NUMBER = {3},
	PAGES = {819--914},
	ISSN = {0894-0347},
	MRCLASS = {11G18 (11G10 14G35)},
	MRNUMBER = {3630089},
	DOI = {10.1090/jams/867},
	URL = {http://dx.doi.org/10.1090/jams/867},
}

@book{lan2013arithmetic,
 author = {Lan, Kai-Wen},
 title = {Arithmetic compactifications of {PEL}-type
 {Shimura} varieties},
 series = {London Mathematical Society Monographs},
 volume = {36},
 year = {2013},
 publisher = {Princeton University Press}
}

@article {rapoportrichartz,
	AUTHOR = {Rapoport, M. and Richartz, M.},
	TITLE = {On the classification and specialization of {$F$}-isocrystals
	with additional structure},
	JOURNAL = {Compositio Math.},
	FJOURNAL = {Compositio Mathematica},
	VOLUME = {103},
	YEAR = {1996},
	NUMBER = {2},
	PAGES = {153--181},
	ISSN = {0010-437X},
	MRCLASS = {14F30 (22E50)},
	MRNUMBER = {1411570},
	MRREVIEWER = {Abdellah Mokrane},
	URL = {http://www.numdam.org/item?id=CM_1996__103_2_153_0},
}

@book {plantonov-rapinchuk,
	AUTHOR = {Platonov, Vladimir and Rapinchuk, Andrei},
	TITLE = {Algebraic groups and number theory},
	SERIES = {Pure and Applied Mathematics},
	VOLUME = {139},
	NOTE = {Translated from the 1991 Russian original by Rachel Rowen},
	PUBLISHER = {Academic Press, Inc., Boston, MA},
	YEAR = {1994},
	PAGES = {xii+614},
	ISBN = {0-12-558180-7},
	MRCLASS = {11E57 (11-02 20Gxx)},
	MRNUMBER = {1278263},
}

@article {kottwitz1992points,
    AUTHOR = {Kottwitz, Robert E.},
     TITLE = {Points on some {S}himura varieties over finite fields},
   JOURNAL = {J. Amer. Math. Soc.},
  FJOURNAL = {Journal of the American Mathematical Society},
    VOLUME = {5},
      YEAR = {1992},
    NUMBER = {2},
     PAGES = {373--444},
      ISSN = {0894-0347},
   MRCLASS = {11G18 (11G25 14G10 14G35)},
  MRNUMBER = {1124982 (93a:11053)},
MRREVIEWER = {Min Ho Lee},
       DOI = {10.2307/2152772},
       URL = {http://dx.doi.org/10.2307/2152772},
}

@article {kottwitzisocrystal,
    AUTHOR = {Kottwitz, Robert E.},
     TITLE = {Isocrystals with additional structure},
   JOURNAL = {Compositio Math.},
  FJOURNAL = {Compositio Mathematica},
    VOLUME = {56},
      YEAR = {1985},
    NUMBER = {2},
     PAGES = {201--220},
      ISSN = {0010-437X},
     CODEN = {CMPMAF},
   MRCLASS = {14L25 (14F30 14L15 20G25)},
  MRNUMBER = {809866 (87i:14040)},
MRREVIEWER = {K. F. Lai},
       URL = {http://www.numdam.org/item?id=CM_1985__56_2_201_0},
}

@article {kottwitzisocrystal2,
    AUTHOR = {Kottwitz, Robert E.},
     TITLE = {Isocrystals with additional structure. {II}},
   JOURNAL = {Compositio Math.},
  FJOURNAL = {Compositio Mathematica},
    VOLUME = {109},
      YEAR = {1997},
    NUMBER = {3},
     PAGES = {255--339},
      ISSN = {0010-437X},
     CODEN = {CMPMAF},
   MRCLASS = {20G25 (11S25 14F30 14L05)},
  MRNUMBER = {1485921 (99e:20061)},
MRREVIEWER = {Guy Rousseau},
       DOI = {10.1023/A:1000102604688},
       URL = {http://dx.doi.org/10.1023/A:1000102604688},
}

@article {langlands1977shimura,
    AUTHOR = {Langlands, R. P.},
     TITLE = {Shimura varieties and the {S}elberg trace formula},
   JOURNAL = {Canad. J. Math.},
  FJOURNAL = {Canadian Journal of Mathematics. Journal Canadien de
              Math\'ematiques},
    VOLUME = {29},
      YEAR = {1977},
    NUMBER = {6},
     PAGES = {1292--1299},
      ISSN = {0008-414X},
   MRCLASS = {10D20 (14G10)},
  MRNUMBER = {0498400 (58 \#16528)},
MRREVIEWER = {Rita Hall},
}

@article {langlands1979zeta,
    AUTHOR = {Langlands, R. P.},
     TITLE = {On the zeta functions of some simple {S}himura varieties},
   JOURNAL = {Canad. J. Math.},
  FJOURNAL = {Canadian Journal of Mathematics. Journal Canadien de
              Math\'ematiques},
    VOLUME = {31},
      YEAR = {1979},
    NUMBER = {6},
     PAGES = {1121--1216},
      ISSN = {0008-414X},
     CODEN = {CJMAAB},
   MRCLASS = {10D40 (12A67 14G25 14K15)},
  MRNUMBER = {553157 (82h:10042)},
MRREVIEWER = {O. M. Fomenko},
       DOI = {10.4153/CJM-1979-102-1},
       URL = {http://dx.doi.org/10.4153/CJM-1979-102-1},
}

@article {langlands1987gerben,
    AUTHOR = {Langlands, R. P. and Rapoport, M.},
     TITLE = {Shimuravariet\"aten und {G}erben},
   JOURNAL = {J. Reine Angew. Math.},
  FJOURNAL = {Journal f\"ur die Reine und Angewandte Mathematik},
    VOLUME = {378},
      YEAR = {1987},
     PAGES = {113--220},
      ISSN = {0075-4102},
     CODEN = {JRMAA8},
   MRCLASS = {11G18 (14G10 14G15 14K22)},
  MRNUMBER = {895287 (88i:11036)},
MRREVIEWER = {James Milne},
}

@incollection {deligne1971traveaux,
    AUTHOR = {Deligne, Pierre},
     TITLE = {Travaux de {S}himura},
 BOOKTITLE = {S\'eminaire {B}ourbaki, 23\`eme ann\'ee (1970/71), {E}xp.
              {N}o. 389},
     PAGES = {123--165. Lecture Notes in Math., Vol. 244},
 PUBLISHER = {Springer},
   ADDRESS = {Berlin},
      YEAR = {1971},
   MRCLASS = {14G25 (10D25 14K15)},
  MRNUMBER = {0498581 (58 \#16675)},
MRREVIEWER = {Author's review},
}

@incollection {deligne1979varietes,
    AUTHOR = {Deligne, Pierre},
     TITLE = {Vari\'et\'es de {S}himura: interpr\'etation modulaire, et
              techniques de construction de mod\`eles canoniques},
 BOOKTITLE = {Automorphic forms, representations and {$L$}-functions
              ({P}roc. {S}ympos. {P}ure {M}ath., {O}regon {S}tate {U}niv.,
              {C}orvallis, {O}re., 1977), {P}art 2},
    SERIES = {Proc. Sympos. Pure Math., XXXIII},
     PAGES = {247--289},
 PUBLISHER = {Amer. Math. Soc.},
   ADDRESS = {Providence, R.I.},
      YEAR = {1979},
   MRCLASS = {10D20 (14D20 14G25 14K15)},
  MRNUMBER = {546620 (81i:10032)},
MRREVIEWER = {J. S. Milne},
}

@book {faltingschai,
    AUTHOR = {Faltings, Gerd and Chai, Ching-Li},
     TITLE = {Degeneration of abelian varieties},
    SERIES = {Ergebnisse der Mathematik und ihrer Grenzgebiete (3) [Results
              in Mathematics and Related Areas (3)]},
    VOLUME = {22},
      NOTE = {With an appendix by David Mumford},
 PUBLISHER = {Springer-Verlag},
   ADDRESS = {Berlin},
      YEAR = {1990},
     PAGES = {xii+316},
      ISBN = {3-540-52015-5},
   MRCLASS = {14K10 (11G10 14D20 14K25)},
  MRNUMBER = {1083353 (92d:14036)},
MRREVIEWER = {Min Ho Lee},
}

@book {pink1989compactification,
    AUTHOR = {Pink, Richard},
     TITLE = {Arithmetical compactification of mixed {S}himura varieties},
    SERIES = {Bonner Mathematische Schriften [Bonn Mathematical
              Publications], 209},
      NOTE = {Dissertation, Rheinische Friedrich-Wilhelms-Universit{\"a}t
              Bonn, Bonn, 1989},
 PUBLISHER = {Universit\"at Bonn Mathematisches Institut},
   ADDRESS = {Bonn},
      YEAR = {1990},
     PAGES = {xviii+340},
   MRCLASS = {11G18 (14G35)},
  MRNUMBER = {1128753 (92h:11054)},
MRREVIEWER = {Min Ho Lee},
}

@article {pink1992ladic,
    AUTHOR = {Pink, Richard},
     TITLE = {On {$l$}-adic sheaves on {S}himura varieties and their higher
              direct images in the {B}aily-{B}orel compactification},
   JOURNAL = {Math. Ann.},
  FJOURNAL = {Mathematische Annalen},
    VOLUME = {292},
      YEAR = {1992},
    NUMBER = {2},
     PAGES = {197--240},
      ISSN = {0025-5831},
     CODEN = {MAANA},
   MRCLASS = {11G18 (11F75 14F20)},
  MRNUMBER = {1149032 (93c:11042)},
MRREVIEWER = {Min Ho Lee},
       DOI = {10.1007/BF01444618},
       URL = {http://dx.doi.org/10.1007/BF01444618},
}

@article {kisin2010integral,
    AUTHOR = {Kisin, Mark},
     TITLE = {Integral models for {S}himura varieties of abelian type},
   JOURNAL = {J. Amer. Math. Soc.},
  FJOURNAL = {Journal of the American Mathematical Society},
    VOLUME = {23},
      YEAR = {2010},
    NUMBER = {4},
     PAGES = {967--1012},
      ISSN = {0894-0347},
   MRCLASS = {11G18 (14G35)},
  MRNUMBER = {2669706 (2011j:11109)},
MRREVIEWER = {Jeffrey D. Achter},
       DOI = {10.1090/S0894-0347-10-00667-3},
       URL = {http://dx.doi.org/10.1090/S0894-0347-10-00667-3},
}

@book {reimann1997zeta,
    AUTHOR = {Reimann, Harry},
     TITLE = {The semi-simple zeta function of quaternionic {S}himura
              varieties},
    SERIES = {Lecture Notes in Mathematics},
    VOLUME = {1657},
 PUBLISHER = {Springer-Verlag},
   ADDRESS = {Berlin},
      YEAR = {1997},
     PAGES = {viii+143},
      ISBN = {3-540-62645-X},
   MRCLASS = {11G18 (11G40 14G10 14G35)},
  MRNUMBER = {1470457 (98k:11076)},
MRREVIEWER = {Rutger Noot},
}

@article {hainesbasechange,
    AUTHOR = {Haines, Thomas J.},
     TITLE = {The base change fundamental lemma for central elements in
              parahoric {H}ecke algebras},
   JOURNAL = {Duke Math. J.},
  FJOURNAL = {Duke Mathematical Journal},
    VOLUME = {149},
      YEAR = {2009},
    NUMBER = {3},
     PAGES = {569--643},
      ISSN = {0012-7094},
     CODEN = {DUMJAO},
   MRCLASS = {22E50 (20G25)},
  MRNUMBER = {2553880},
MRREVIEWER = {Anne-Marie H. Aubert},
       DOI = {10.1215/00127094-2009-045},
       URL = {http://dx.doi.org/10.1215/00127094-2009-045},
}

@article {Art78,
    AUTHOR = {Arthur, James },
     TITLE = {A trace formula for reductive groups. {I}. {T}erms associated
              to classes in {$G({\mathbf Q})$}},
   JOURNAL = {Duke Math. J.},
  FJOURNAL = {Duke Mathematical Journal},
    VOLUME = {45},
      YEAR = {1978},
    NUMBER = {4},
     PAGES = {911--952},
      ISSN = {0012-7094},
   MRCLASS = {10D40 (22E55)},
  MRNUMBER = {518111},
MRREVIEWER = {Stephen Gelbart},
       URL = {http://projecteuclid.org.libproxy.berkeley.edu/euclid.dmj/1077313104},
}

@article {Art88b,
 	AUTHOR = {Arthur, James},
 	TITLE = {The invariant trace formula. {II}. {G}lobal theory},
 	JOURNAL = {J. Amer. Math. Soc.},
 	FJOURNAL = {Journal of the American Mathematical Society},
 	VOLUME = {1},
 	YEAR = {1988},
 	NUMBER = {3},
 	PAGES = {501--554},
 	ISSN = {0894-0347},
 	MRCLASS = {22E55 (11F72)},
 	MRNUMBER = {939691 (89j:22039)},
 	MRREVIEWER = {A. B. Venkov},
 	DOI = {10.2307/1990948},
 	URL = {http://dx.doi.org/10.2307/1990948},
 }

@article {Art89,
 	AUTHOR = {Arthur, James},
 	TITLE = {The {$L^2$}-{L}efschetz numbers of {H}ecke operators},
 	JOURNAL = {Invent.~Math.},
 	FJOURNAL = {Inventiones Mathematicae},
 	VOLUME = {97},
 	YEAR = {1989},
 	NUMBER = {2},
 	PAGES = {257-290},
 }

@article {ArtSTF1,
 	AUTHOR = {Arthur, James},
 	TITLE = {A stable trace formula. {I}. {G}eneral expansions},
 	JOURNAL = {J. Inst. Math. Jussieu},
 	FJOURNAL = {Journal of the Institute of Mathematics of Jussieu. JIMJ.
 	Journal de l'Institute de Math\'ematiques de Jussieu},
 	VOLUME = {1},
 	YEAR = {2002},
 	NUMBER = {2},
 	PAGES = {175--277},
 	ISSN = {1474-7480},
 	MRCLASS = {11F72 (11R39 22E55)},
 	MRNUMBER = {1954821 (2004b:11066)},
 	MRREVIEWER = {Volker J. Heiermann},
 	DOI = {10.1017/S1474-748002000051},
 	URL = {http://dx.doi.org/10.1017/S1474-748002000051},
 }

@article {ArtSTF2,
 	AUTHOR = {Arthur, James},
 	TITLE = {A stable trace formula. {II}. {G}lobal descent},
 	JOURNAL = {Invent. Math.},
 	FJOURNAL = {Inventiones Mathematicae},
 	VOLUME = {143},
 	YEAR = {2001},
 	NUMBER = {1},
 	PAGES = {157--220},
 	ISSN = {0020-9910},
 	CODEN = {INVMBH},
 	MRCLASS = {11F72 (22E45)},
 	MRNUMBER = {1802795 (2002m:11043)},
 	MRREVIEWER = {Volker J. Heiermann},
 	DOI = {10.1007/s002220000107},
 	URL = {http://dx.doi.org/10.1007/s002220000107},
 }

@book {Arthur,
 	AUTHOR = {Arthur, James},
 	TITLE = {The endoscopic classification of representations},
 	SERIES = {American Mathematical Society Colloquium Publications},
 	VOLUME = {61},
 	NOTE = {Orthogonal and symplectic groups},
 	PUBLISHER = {American Mathematical Society, Providence, RI},
 	YEAR = {2013},
 	PAGES = {xviii+590},
 	ISBN = {978-0-8218-4990-3},
 	MRCLASS = {22E55 (11F66 11R37 20G25 22E50)},
 	MRNUMBER = {3135650},
 }

@article {Dalal,
	AUTHOR = {Dalal, Rahul},
	TITLE = {Sato-{T}ate equidistribution for families of automorphic
	representations through the stable trace formula},
	JOURNAL = {Algebra Number Theory},
	FJOURNAL = {Algebra \& Number Theory},
	VOLUME = {16},
	YEAR = {2022},
	NUMBER = {1},
	PAGES = {59--137},
	ISSN = {1937-0652,1944-7833},
	MRCLASS = {11F55 (11F70 11F72 11F75 22E50 22E55)},
	MRNUMBER = {4384564},
	MRREVIEWER = {Henry\ H.\ Kim},
	DOI = {10.2140/ant.2022.16.59},
	URL = {https://doi.org/10.2140/ant.2022.16.59},
}

@ARTICLE{BM,
 	author={Brumley, Farrell and Milicevic, Djordje},
 	title={Counting cusp forms by analytic conductor},
 	journal={preprint, https://www.math.univ-paris13.fr/~brumley/},
 	
 	year={2017},
 	
 }

@article {CD90,
    AUTHOR = {Clozel, Laurent and Delorme, Patrick},
     TITLE = {Le th\'eor\`eme de {P}aley-{W}iener invariant pour les groupes de
              {L}ie r\'eductifs. {II}},
   JOURNAL = {Ann. Sci. \'Ecole Norm. Sup. (4)},
  FJOURNAL = {Annales Scientifiques de l'\'Ecole Normale Sup\'erieure. Quatri\`eme
              S\'erie},
    VOLUME = {23},
      YEAR = {1990},
    NUMBER = {2},
     PAGES = {193--228},
      ISSN = {0012-9593},
   MRCLASS = {22E30 (22E45 22E47)},
  MRNUMBER = {1046496},
MRREVIEWER = {Michael Cowling},
       URL = {http://www.numdam.org/item?id=ASENS_1990_4_23_2_193_0},
}

@article {CL10,
 	AUTHOR = {Cluckers, Raf and Loeser, Fran{\c{c}}ois},
 	TITLE = {Constructible exponential functions, motivic {F}ourier
 	transform and transfer principle},
 	JOURNAL = {Ann. of Math. (2)},
 	FJOURNAL = {Annals of Mathematics. Second Series},
 	VOLUME = {171},
 	YEAR = {2010},
 	NUMBER = {2},
 	PAGES = {1011--1065},
 	ISSN = {0003-486X},
 	CODEN = {ANMAAH},
 	MRCLASS = {14E18},
 	MRNUMBER = {2630060 (2011g:14036)},
 	MRREVIEWER = {Tommaso De Fernex},
 	DOI = {10.4007/annals.2010.171.1011},
 	URL = {http://dx.doi.org/10.4007/annals.2010.171.1011},
 }

@article {EGA4-4,
	AUTHOR = {Grothendieck, A.},
	TITLE = {\'{E}l\'{e}ments de g\'{e}om\'{e}trie alg\'{e}brique. {IV}. \'{E}tude locale des
	sch\'{e}mas et des morphismes de sch\'{e}mas {IV}},
	JOURNAL = {Inst. Hautes \'{E}tudes Sci. Publ. Math.},
	FJOURNAL = {Institut des Hautes \'{E}tudes Scientifiques. Publications
	Math\'{e}matiques},
	volume = {32},
	YEAR = {1967},
	PAGES = {5--361},
	ISSN = {0073-8301},
	MRCLASS = {14.55},
	MRNUMBER = {238860},
	MRREVIEWER = {J. P. Murre},
	URL = {http://www.numdam.org/item?id=PMIHES_1967__32__361_0},
}

@article {Fra98,
 	AUTHOR = {Franke, Jens},
 	TITLE = {Harmonic analysis in weighted {$L_2$}-spaces},
 	JOURNAL = {Ann. Sci. \'Ecole Norm. Sup. (4)},
 	FJOURNAL = {Annales Scientifiques de l'\'Ecole Normale Sup\'erieure.
 	Quatri\`eme S\'erie},
 	VOLUME = {31},
 	YEAR = {1998},
 	NUMBER = {2},
 	PAGES = {181--279},
 	ISSN = {0012-9593},
 	CODEN = {ASENAH},
 	MRCLASS = {11F75 (11F70 22E41)},
 	MRNUMBER = {1603257 (2000f:11065)},
 	DOI = {10.1016/S0012-9593(98)80015-3},
 	URL = {http://dx.doi.org/10.1016/S0012-9593(98)80015-3},
 }

@ARTICLE{GKM97,
 	AUTHOR =       {Goresky, M. and Kottwitz, R. and MacPherson, R.},
 	TITLE =        {{D}iscrete series characters and the {L}efschetz formula for {H}ecke operators},
 	JOURNAL =      {Duke Math. J.},
 	YEAR =         {1997},
 	volume =       {89},
 	number =       {3},
 	pages =        {477-554},
 	
 }

@article {Gro99,
 	AUTHOR = {Gross, Benedict H.},
 	TITLE = {Algebraic modular forms},
 	JOURNAL = {Israel J. Math.},
 	FJOURNAL = {Israel Journal of Mathematics},
 	VOLUME = {113},
 	YEAR = {1999},
 	PAGES = {61--93},
 	ISSN = {0021-2172},
 	CODEN = {ISJMAP},
 	MRCLASS = {11F55 (11F80 20G30 22E55)},
 	MRNUMBER = {1729443 (2001b:11037)},
 	MRREVIEWER = {Stefan K{\"u}hnlein},
 	DOI = {10.1007/BF02780173},
 	URL = {http://dx.doi.org/10.1007/BF02780173},
 }

@article {Hal95,
	AUTHOR = {Hales, Thomas C.},
	TITLE = {On the fundamental lemma for standard endoscopy: reduction to
	unit elements},
	JOURNAL = {Canad. J. Math.},
	FJOURNAL = {Canadian Journal of Mathematics. Journal Canadien de
	Math\'{e}matiques},
	VOLUME = {47},
	YEAR = {1995},
	NUMBER = {5},
	PAGES = {974--994},
	ISSN = {0008-414X},
	MRCLASS = {22E50 (20C05 22E35)},
	MRNUMBER = {1350645},
	MRREVIEWER = {Mark Reeder},
	DOI = {10.4153/CJM-1995-051-5},
	URL = {https://doi-org.proxy-um.researchport.umd.edu/10.4153/CJM-1995-051-5},
}

@article {Lau14,
	AUTHOR = {Lau, Eike},
	TITLE = {Relations between {D}ieudonn\'{e} displays and crystalline
	{D}ieudonn\'{e} theory},
	JOURNAL = {Algebra Number Theory},
	FJOURNAL = {Algebra \& Number Theory},
	VOLUME = {8},
	YEAR = {2014},
	NUMBER = {9},
	PAGES = {2201--2262},
	ISSN = {1937-0652},
	MRCLASS = {14L05 (14F30)},
	MRNUMBER = {3294388},
	MRREVIEWER = {Alexander Boris Ivanov},
	DOI = {10.2140/ant.2014.8.2201},
	URL = {https://doi-org.proxy-um.researchport.umd.edu/10.2140/ant.2014.8.2201},
}

@article {Lau19,
	AUTHOR = {Lau, Eike},
	TITLE = {Displayed equations for {G}alois representations},
	JOURNAL = {Nagoya Math. J.},
	FJOURNAL = {Nagoya Mathematical Journal},
	VOLUME = {235},
	YEAR = {2019},
	PAGES = {86--114},
	ISSN = {0027-7630},
	MRCLASS = {14L05 (11F80 14F30)},
	MRNUMBER = {3986712},
	MRREVIEWER = {Davide Lombardo},
	DOI = {10.1017/nmj.2018.3},
	URL = {https://doi-org.proxy-um.researchport.umd.edu/10.1017/nmj.2018.3},
}

@article {Kim12,
	AUTHOR = {Kim, Wansu},
	TITLE = {The classification of {$p$}-divisible groups over 2-adic
	discrete valuation rings},
	JOURNAL = {Math. Res. Lett.},
	FJOURNAL = {Mathematical Research Letters},
	VOLUME = {19},
	YEAR = {2012},
	NUMBER = {1},
	PAGES = {121--141},
	ISSN = {1073-2780},
	MRCLASS = {14L15 (11S20 14F30)},
	MRNUMBER = {2923180},
	MRREVIEWER = {Marc-Hubert Nicole},
	DOI = {10.4310/MRL.2012.v19.n1.a10},
	URL = {https://doi-org.proxy-um.researchport.umd.edu/10.4310/MRL.2012.v19.n1.a10},
}

@article {KMP16,
 	AUTHOR = {Kim, Wansu and Madapusi Pera, Keerthi},
 	TITLE = {2-adic integral canonical models},
 	JOURNAL = {Forum Math. Sigma},
 	FJOURNAL = {Forum of Mathematics. Sigma},
 	VOLUME = {4},
 	YEAR = {2016},
 	PAGES = {e28, 34},
 	ISSN = {2050-5094},
 	MRCLASS = {14G35 (11G18)},
 	MRNUMBER = {3569319},
 	MRREVIEWER = {Su-ion Ih},
 	DOI = {10.1017/fms.2016.23},
 	URL = {http://dx.doi.org/10.1017/fms.2016.23},
 }

@article {KMSW,
 	AUTHOR = {Kaletha, Tasho and Mínguez, Alberto and Shin, Sug Woo and White, Paul-James},
 	TITLE = {Endoscopic Classification of Representations: Inner Forms of Unitary Groups},
 	EPRINT = {arXiv:1409.3731},
 
        YEAR = {2014},
        note = {arXiv:1409.3731,
        \url{https://arxiv.org/abs/1409.3731}}
	}

@article {Kot82,
 	AUTHOR = {Kottwitz, Robert E.},
 	TITLE = {Rational conjugacy classes in reductive groups},
 	JOURNAL = {Duke Math. J.},
 	FJOURNAL = {Duke Mathematical Journal},
 	VOLUME = {49},
 	YEAR = {1982},
 	NUMBER = {4},
 	PAGES = {785--806},
 	ISSN = {0012-7094},
 	CODEN = {DUMJAO},
 	MRCLASS = {20G15 (22E55)},
 	MRNUMBER = {683003 (84k:20020)},
 	MRREVIEWER = {Andy R. Magid},
 	URL = {http://projecteuclid.org/euclid.dmj/1077315531},
 }

@article {Kot83,
 	AUTHOR = {Kottwitz, Robert E.},
 	TITLE = {Sign changes in harmonic analysis on reductive groups},
 	JOURNAL = {Trans. Amer. Math. Soc.},
 	FJOURNAL = {Transactions of the American Mathematical Society},
 	VOLUME = {278},
 	YEAR = {1983},
 	NUMBER = {1},
 	PAGES = {289--297},
 	ISSN = {0002-9947},
 	CODEN = {TAMTAM},
 	MRCLASS = {22E35 (22E30)},
 	MRNUMBER = {697075 (84i:22012)},
 	MRREVIEWER = {Martin L. Karel},
 	DOI = {10.2307/1999316},
 	URL = {http://dx.doi.org/10.2307/1999316},
 }

@article {Kot84a,
 	AUTHOR = {Kottwitz, Robert E.},
 	TITLE = {Stable trace formula: cuspidal tempered terms},
 	JOURNAL = {Duke Math. J.},
 	FJOURNAL = {Duke Mathematical Journal},
 	VOLUME = {51},
 	YEAR = {1984},
 	NUMBER = {3},
 	PAGES = {611--650},
 	ISSN = {0012-7094},
 	CODEN = {DUMJAO},
 	MRCLASS = {11R39 (11F70 11F72 22E55)},
 	MRNUMBER = {757954 (85m:11080)},
 	MRREVIEWER = {Jean-Pierre Labesse},
 	DOI = {10.1215/S0012-7094-84-05129-9},
 	URL = {http://dx.doi.org/10.1215/S0012-7094-84-05129-9},
 }

@article {Kot86,
 	AUTHOR = {Kottwitz, Robert E.},
 	TITLE = {Stable trace formula: elliptic singular terms},
 	JOURNAL = {Math. Ann.},
 	FJOURNAL = {Mathematische Annalen},
 	VOLUME = {275},
 	YEAR = {1986},
 	NUMBER = {3},
 	PAGES = {365--399},
 	ISSN = {0025-5831},
 	CODEN = {MAANA},
 	MRCLASS = {22E55 (11F70 11F72)},
 	MRNUMBER = {858284 (88d:22027)},
 	MRREVIEWER = {J. Szmidt},
 	DOI = {10.1007/BF01458611},
 	URL = {http://dx.doi.org/10.1007/BF01458611},
 }

@article {Kot88,
 	AUTHOR = {Kottwitz, Robert E.},
 	TITLE = {Tamagawa numbers},
 	JOURNAL = {Ann. of Math. (2)},
 	FJOURNAL = {Annals of Mathematics. Second Series},
 	VOLUME = {127},
 	YEAR = {1988},
 	NUMBER = {3},
 	PAGES = {629--646},
 	ISSN = {0003-486X},
 	CODEN = {ANMAAH},
 	MRCLASS = {11F70 (11E72 22E55)},
 	MRNUMBER = {942522 (90e:11075)},
 	MRREVIEWER = {Stephen Gelbart},
 	DOI = {10.2307/2007007},
 	URL = {http://dx.doi.org/10.2307/2007007},
 }

@article {Kot92b,
 	AUTHOR = {Kottwitz, Robert E.},
 	TITLE = {On the {$\lambda$}-adic representations associated to some
 	simple {S}himura varieties},
 	JOURNAL = {Invent. Math.},
 	FJOURNAL = {Inventiones Mathematicae},
 	VOLUME = {108},
 	YEAR = {1992},
 	NUMBER = {3},
 	PAGES = {653--665},
 	ISSN = {0020-9910},
 	MRCLASS = {11G18 (11F70 11F85 14G35 22E55)},
 	MRNUMBER = {1163241},
 	MRREVIEWER = {Min Ho Lee},
 	DOI = {10.1007/BF02100620},
 	URL = {http://dx.doi.org/10.1007/BF02100620},
 }

@incollection {Kot05,
 	AUTHOR = {Kottwitz, Robert E.},
 	TITLE = {Harmonic analysis on reductive {$p$}-adic groups and {L}ie
 	algebras},
 	BOOKTITLE = {Harmonic analysis, the trace formula, and {S}himura varieties},
 	SERIES = {Clay Math. Proc.},
 	VOLUME = {4},
 	PAGES = {393--522},
 	PUBLISHER = {Amer. Math. Soc.},
 	ADDRESS = {Providence, RI},
 	YEAR = {2005},
 	MRCLASS = {22E35 (17B99)},
 	MRNUMBER = {2192014 (2006m:22016)},
 	MRREVIEWER = {David A. Renard},
 }

@article {KS99,
 	AUTHOR = {Kottwitz, Robert E. and Shelstad, Diana},
 	TITLE = {Foundations of twisted endoscopy},
 	JOURNAL = {Ast\'erisque},
 	FJOURNAL = {Ast\'erisque},
 	volume = {255},
 	YEAR = {1999},
 	PAGES = {vi+190},
 	ISSN = {0303-1179},
 	MRCLASS = {22E55 (11F70 11R34 22-02 22E50)},
 	MRNUMBER = {1687096 (2000k:22024)},
 	MRREVIEWER = {Volker J. Heiermann},
 }

@article {KS-GSp,
	AUTHOR = {Kret, Arno and Shin, Sug Woo},
	TITLE = {Galois representations for general symplectic groups},
	JOURNAL = {J. Eur. Math. Soc. (JEMS)},
	FJOURNAL = {Journal of the European Mathematical Society (JEMS)},
	VOLUME = {25},
	YEAR = {2023},
	NUMBER = {1},
	PAGES = {75--152},
	ISSN = {1435-9855,1435-9863},
	MRCLASS = {11R39 (11F70 11F80 11G18)},
	MRNUMBER = {4556781},
	MRREVIEWER = {Jack\ Shotton},
	DOI = {10.4171/jems/1179},
	URL = {https://doi.org/10.4171/jems/1179},
}

@article {Lab04,
 	AUTHOR = {Labesse, J.-P.},
 	TITLE = {Stable twisted trace formula: elliptic terms},
 	JOURNAL = {J. Inst. Math. Jussieu},
 	FJOURNAL = {Journal of the Institute of Mathematics of Jussieu. JIMJ.
 	Journal de l'Institute de Math\'ematiques de Jussieu},
 	VOLUME = {3},
 	YEAR = {2004},
 	NUMBER = {4},
 	PAGES = {473--530},
 	ISSN = {1474-7480},
 	MRCLASS = {11F72 (11R34 11R39)},
 	MRNUMBER = {2094449 (2005g:11084)},
 	MRREVIEWER = {Pierre-Henri Chaudouard},
 	DOI = {10.1017/S1474748004000143},
 	URL = {http://dx.doi.org/10.1017/S1474748004000143},
 }

@article {Lan79,
 	AUTHOR = {Langlands, R. P.},
 	TITLE = {Stable conjugacy: definitions and lemmas},
 	JOURNAL = {Canad. J. Math.},
 	FJOURNAL = {Canadian Journal of Mathematics. Journal Canadien de
 	Math\'ematiques},
 	VOLUME = {31},
 	YEAR = {1979},
 	NUMBER = {4},
 	PAGES = {700--725},
 	ISSN = {0008-414X},
 	CODEN = {CJMAAB},
 	MRCLASS = {10D40 (12A67 22E50)},
 	MRNUMBER = {540901 (82j:10054)},
 	MRREVIEWER = {Joe Repka},
 	DOI = {10.4153/CJM-1979-069-2},
 	URL = {http://dx.doi.org/10.4153/CJM-1979-069-2},
 }

@article {Lau97,
 	AUTHOR = {Laumon, G{\'e}rard},
 	TITLE = {Sur la cohomologie \`a supports compacts des vari\'et\'es de
 	{S}himura pour {${\mathrm {GSp}}(4)_{\mathbf Q}$}},
 	JOURNAL = {Compositio Math.},
 	FJOURNAL = {Compositio Mathematica},
 	VOLUME = {105},
 	YEAR = {1997},
 	NUMBER = {3},
 	PAGES = {267--359},
 	ISSN = {0010-437X},
 	CODEN = {CMPMAF},
 	MRCLASS = {11G18 (11F70 11F72 14G35 22E55)},
 	MRNUMBER = {1440724 (2000a:11097)},
 	MRREVIEWER = {{\'E}ric Urban},
 	DOI = {10.1023/A:1000102414055},
 	URL = {http://dx.doi.org/10.1023/A:1000102414055},
 }

@article {LS87,
	AUTHOR = {Langlands, R. P. and Shelstad, D.},
	TITLE = {On the definition of transfer factors},
	JOURNAL = {Math. Ann.},
	FJOURNAL = {Mathematische Annalen},
	VOLUME = {278},
	YEAR = {1987},
	NUMBER = {1-4},
	PAGES = {219--271},
	ISSN = {0025-5831},
	MRCLASS = {11S37 (11F70 11F72 22E35 22E45 22E50)},
	MRNUMBER = {909227},
	MRREVIEWER = {Ernst-Wilhelm Zink},
	DOI = {10.1007/BF01458070},
	URL = {https://doi-org.proxy-um.researchport.umd.edu/10.1007/BF01458070},
}

@incollection {LS90,
	AUTHOR = {Langlands, R. and Shelstad, D.},
	TITLE = {Descent for transfer factors},
	BOOKTITLE = {The {G}rothendieck {F}estschrift, {V}ol. {II}},
	SERIES = {Progr. Math.},
	VOLUME = {87},
	PAGES = {485--563},
	PUBLISHER = {Birkh\"{a}user Boston, Boston, MA},
	YEAR = {1990},
	MRCLASS = {22E50 (11S37)},
	MRNUMBER = {1106907},
	MRREVIEWER = {Ernst-Wilhelm Zink},
}

@article {Mok,
 	AUTHOR = {Mok, Chung Pang},
 	TITLE = {Endoscopic classification of representations of quasi-split
 	unitary groups},
 	JOURNAL = {Mem. Amer. Math. Soc.},
 	FJOURNAL = {Memoirs of the American Mathematical Society},
 	VOLUME = {235},
 	YEAR = {2015},
 	NUMBER = {1108},
 	PAGES = {vi+248},
 	ISSN = {0065-9266},
 	ISBN = {978-1-4704-1041-4; 978-1-4704-2226-4},
 	MRCLASS = {22E55 (11R42 22E50)},
 	MRNUMBER = {3338302},
 	MRREVIEWER = {Neven Grbac},
 	DOI = {10.1090/memo/1108},
 	URL = {http://dx.doi.org/10.1090/memo/1108},
 }

@ARTICLE{Mor08,
 	author={Morel, S.},
 	title={Complexes pond\'{e}r\'{e}s sur les compactifications de {B}aily-{B}orel:
 	{L}e cas des vari\'{e}t\'{e}s de {S}iegel},
 	journal={J. Amer. Math. Soc.},
 	year={2008},
 	volume={21},
 	pages={23-61},
 }

@article {MT,
 	AUTHOR = {Matz, Jasmin and Templier, Nicolas},
 	TITLE = {Sato-{T}ate equidistribution for families of {Hecke-Maass} forms on {$SL(n,\mathbb R)/SO(n)$}},
 	JOURNAL = {preprint, arXiv:1505.07285},
 	year={},
 }

@article {Ngo10,
 	AUTHOR = {Ng{\^o}, Bao Ch{\^a}u},
 	TITLE = {Le lemme fondamental pour les alg\`ebres de {L}ie},
 	JOURNAL = {Publ. Math. Inst. Hautes \'Etudes Sci.},
 	FJOURNAL = {Publications Math\'ematiques. Institut de Hautes \'Etudes
 	Scientifiques},
 	volume = {111},
 	YEAR = {2010},
 	PAGES = {1--169},
 	ISSN = {0073-8301},
 	MRCLASS = {22E35 (11S37 14D23 14G35 22E50)},
 	MRNUMBER = {2653248},
 	DOI = {10.1007/s10240-010-0026-7},
 	URL = {http://dx.doi.org/10.1007/s10240-010-0026-7},
 }

@article {CS79,
 	AUTHOR = {Colliot-Th\'{e}l\`{e}ne, J.-L. and Sansuc, J.-J.},
 	TITLE = {Fibr\'{e}s quadratiques et composantes connexes r\'{e}elles},
 	JOURNAL = {Math. Ann.},
 	FJOURNAL = {Mathematische Annalen},
 	VOLUME = {244},
 	YEAR = {1979},
 	NUMBER = {2},
 	PAGES = {105--134},
 	ISSN = {0025-5831},
 	MRCLASS = {14G05 (10C04)},
 	MRNUMBER = {550842},
 	MRREVIEWER = {Daniel Coray},
 	DOI = {10.1007/BF01420486},
 	URL = {https://doi-org.proxy-um.researchport.umd.edu/10.1007/BF01420486},
 	}

@ARTICLE{She82,
 	author={Shelstad, D.},
 	title={{$L$}-indistinguishability for real groups},
 	journal={Math. Ann.},
 	year={1982},
 	volume={259},
 	pages={385-430},
 }

@article {Shi09,
 	AUTHOR = {Shin, Sug Woo},
 	TITLE = {Counting points on {I}gusa varieties},
 	JOURNAL = {Duke Math. J.},
 	FJOURNAL = {Duke Mathematical Journal},
 	VOLUME = {146},
 	YEAR = {2009},
 	NUMBER = {3},
 	PAGES = {509--568},
 	ISSN = {0012-7094},
 	MRCLASS = {11G18 (11F72 11G15 11R39 14G35)},
 	MRNUMBER = {2484281},
 	MRREVIEWER = {Nguy\'a$\gg$꿳 Qu\~A'c Th\~A!`ng},
 	DOI = {10.1215/00127094-2009-004},
 	URL = {http://dx.doi.org/10.1215/00127094-2009-004},
 }

@article {Shin10,
 	AUTHOR = {Shin, Sug Woo},
 	TITLE = {A stable trace formula for {I}gusa varieties},
 	JOURNAL = {J. Inst. Math. Jussieu},
 	FJOURNAL = {Journal of the Institute of Mathematics of Jussieu. JIMJ.
 	Journal de l'Institut de Math\'ematiques de Jussieu},
 	VOLUME = {9},
 	YEAR = {2010},
 	NUMBER = {4},
 	PAGES = {847--895},
 	ISSN = {1474-7480},
 	MRCLASS = {22Exx (11G18 14G35)},
 	MRNUMBER = {2684263},
 	DOI = {10.1017/S1474748010000046},
 	URL = {http://dx.doi.org/10.1017/S1474748010000046},
 }

@article {Wal97,
 	AUTHOR = {Waldspurger, J.-L.},
 	TITLE = {Le lemme fondamental implique le transfert},
 	JOURNAL = {Compositio Math.},
 	FJOURNAL = {Compositio Mathematica},
 	VOLUME = {105},
 	YEAR = {1997},
 	NUMBER = {2},
 	PAGES = {153--236},
 	ISSN = {0010-437X},
 	CODEN = {CMPMAF},
 	MRCLASS = {22E50 (11F72 11F85 22E35)},
 	MRNUMBER = {1440722 (98h:22023)},
 	MRREVIEWER = {David Manderscheid},
 	DOI = {10.1023/A:1000103112268},
 	URL = {http://dx.doi.org/10.1023/A:1000103112268},
 }

@article {Wal06,
 	AUTHOR = {Waldspurger, J.-L.},
 	TITLE = {Endoscopie et changement de caract\'eristique},
 	JOURNAL = {J. Inst. Math. Jussieu},
 	FJOURNAL = {Journal of the Institute of Mathematics of Jussieu. JIMJ.
 	Journal de l'Institut de Math\'ematiques de Jussieu},
 	VOLUME = {5},
 	YEAR = {2006},
 	NUMBER = {3},
 	PAGES = {423--525},
 	ISSN = {1474-7480},
 	MRCLASS = {22E35 (20G25 22E60)},
 	MRNUMBER = {2241929 (2007h:22007)},
 	MRREVIEWER = {Bertrand Lemaire},
 	DOI = {10.1017/S1474748006000041},
 	URL = {http://dx.doi.org/10.1017/S1474748006000041},
 }

@article {Wal08,
 	AUTHOR = {Waldspurger, J.-L.},
 	TITLE = {L'endoscopie tordue n'est pas si tordue},
 	JOURNAL = {Mem. Amer. Math. Soc.},
 	FJOURNAL = {Memoirs of the American Mathematical Society},
 	VOLUME = {194},
 	YEAR = {2008},
 	NUMBER = {908},
 	PAGES = {x+261},
 	ISSN = {0065-9266},
 	CODEN = {MAMCAU},
 	MRCLASS = {22E50 (11R34 22E35)},
 	MRNUMBER = {2418405 (2011d:22020)},
 	MRREVIEWER = {Anne-Marie H. Aubert},
 }

@article {Xu-GSpGSO,
    AUTHOR = {Xu, Bin},
     TITLE = {L-packets of quasisplit {$GSp(2n)$} and {$GO(2n)$}},
   JOURNAL = {Math. Ann.},
  FJOURNAL = {Mathematische Annalen},
    VOLUME = {370},
      YEAR = {2018},
    NUMBER = {1-2},
     PAGES = {71--189},
      ISSN = {0025-5831},
   MRCLASS = {22E50 (11F70)},
  MRNUMBER = {3747484},
MRREVIEWER = {Ivan Mati\'{c}},
       DOI = {10.1007/s00208-016-1515-x},
       URL = {https://doi.org/10.1007/s00208-016-1515-x},
}

@article {Xu-GSpGSO2,
	AUTHOR = {Xu, Bin},
	TITLE = {Global {$L$}-packets of quasisplit {${\mathrm{ GSp}}(2n)$} and
	{${\mathrm {GO}}(2n)$}},
	JOURNAL = {Amer. J. Math.},
	FJOURNAL = {American Journal of Mathematics},
	VOLUME = {147},
	YEAR = {2025},
	NUMBER = {2},
	PAGES = {401--464},
	ISSN = {0002-9327,1080-6377},
	MRCLASS = {22E50 (11F70)},
	MRNUMBER = {4887967},
	MRREVIEWER = {Neven\ Grbac},
	DOI = {10.1353/ajm.2025.a954647},
	URL = {https://doi.org/10.1353/ajm.2025.a954647},
}

@book {Mor10,
 	AUTHOR = {Morel, Sophie},
 	TITLE = {On the cohomology of certain noncompact {S}himura varieties},
 	SERIES = {Annals of Mathematics Studies},
 	VOLUME = {173},
 	NOTE = {With an appendix by Robert Kottwitz},
 	PUBLISHER = {Princeton University Press, Princeton, NJ},
 	YEAR = {2010},
 	PAGES = {xii+217},
 	ISBN = {978-0-691-14293-7},
 	MRCLASS = {11F75 (11F72 11G18 14G35)},
 	MRNUMBER = {2567740 (2011b:11073)},
 	MRREVIEWER = {Eva Viehmann},
 	DOI = {10.1515/9781400835393},
 	URL = {http://dx.doi.org/10.1515/9781400835393},
 }

@article {Mor11,
 	AUTHOR = {Morel, Sophie},
 	TITLE = {Cohomologie d'intersection des vari\'et\'es modulaires de
 	{S}iegel, suite},
 	JOURNAL = {Compos. Math.},
 	FJOURNAL = {Compositio Mathematica},
 	VOLUME = {147},
 	YEAR = {2011},
 	NUMBER = {6},
 	PAGES = {1671--1740},
 	ISSN = {0010-437X},
 	MRCLASS = {11G18 (14F43 14G32)},
 	MRNUMBER = {2862060},
 	MRREVIEWER = {Keerthi Madapusi Pera},
 	DOI = {10.1112/S0010437X11005409},
 	URL = {http://dx.doi.org/10.1112/S0010437X11005409},
 }

@incollection {Lan76,
 	AUTHOR = {Langlands, R. P.},
 	TITLE = {Some contemporary problems with origins in the {J}ugendtraum},
 	BOOKTITLE = {Mathematical developments arising from {H}ilbert problems
 	({P}roc. {S}ympos. {P}ure {M}ath., {V}ol. {XXVIII}, {N}orthern
 	{I}llinois {U}niv., {D}e {K}alb, {I}ll., 1974)},
 	PAGES = {401--418},
 	PUBLISHER = {Amer. Math. Soc.},
 	ADDRESS = {Providence, R. I.},
 	YEAR = {1976},
 	MRCLASS = {12A65 (10D15)},
 	MRNUMBER = {0437500 (55 \#10426)},
 	MRREVIEWER = {Stephen Gelbart},
 }

@incollection {Kot90,
 	AUTHOR = {Kottwitz, Robert E.},
 	TITLE = {Shimura varieties and {$\lambda$}-adic representations},
 	BOOKTITLE = {Automorphic forms, {S}himura varieties, and {$L$}-functions,
 	{V}ol.\ {I} ({A}nn {A}rbor, {MI}, 1988)},
 	SERIES = {Perspect. Math.},
 	VOLUME = {10},
 	PAGES = {161--209},
 	PUBLISHER = {Academic Press, Boston, MA},
 	YEAR = {1990},
 	MRCLASS = {11G18 (11F70 11R39)},
 	MRNUMBER = {1044820},
 	}

@book {Montreal92,
 	TITLE = {The zeta functions of {P}icard modular surfaces},
 	EDITOR = {Langlands, Robert P. and Ramakrishnan, Dinakar},
 	PUBLISHER = {Universit\'e de Montr\'eal, Centre de Recherches Math\'ematiques,
 	Montreal, QC},
 	YEAR = {1992},
 	PAGES = {xiv+492},
 	ISBN = {2-921120-08-9},
 	MRCLASS = {14-06 (00B25 11-06 14G10 14G35)},
 	MRNUMBER = {1155223},
 	}

@incollection {MS82,
 	AUTHOR =       {Milne, James and Shih, Kuang-Yen},
 	TITLE =        {Conjugates of {S}himura varieties},
 	BOOKTITLE = {Hodge cycles, motives, and {S}himura varieties},
 	SERIES = {Lecture Notes in Mathematics},
 	VOLUME = {900},
 	PAGES = {280--356},
 	PUBLISHER = {Springer-Verlag},
 	YEAR = {1982},
 	
 	}

@article {KalGlo,
 		AUTHOR = {Kaletha, Tasho},
 		TITLE = {Global rigid inner forms and multiplicities of discrete
 		automorphic representations},
 		JOURNAL = {Invent. Math.},
 		FJOURNAL = {Inventiones Mathematicae},
 		VOLUME = {213},
 		YEAR = {2018},
 		NUMBER = {1},
 		PAGES = {271--369},
 		ISSN = {0020-9910},
 		MRCLASS = {11R34},
 		MRNUMBER = {3815567},
 		MRREVIEWER = {Meng Fai Lim},
 		DOI = {10.1007/s00222-018-0791-3},
 		URL = {https://doi-org.proxy-um.researchport.umd.edu/10.1007/s00222-018-0791-3},
 	}

@article {KalLoc,
 	AUTHOR = {Kaletha, Tasho},
 	TITLE = {Rigid inner forms of real and {$p$}-adic groups},
 	JOURNAL = {Ann. of Math. (2)},
 	FJOURNAL = {Annals of Mathematics. Second Series},
 	VOLUME = {184},
 	YEAR = {2016},
 	NUMBER = {2},
 	PAGES = {559--632},
 	ISSN = {0003-486X},
 	MRCLASS = {20G25 (14L17 20G10)},
 	MRNUMBER = {3548533},
 	MRREVIEWER = {Wen-Wei Li},
 	URL = {https://doi.org/10.4007/annals.2016.184.2.6},
 	}

@ARTICLE{KS12,
 		author = {{Kottwitz}, R. and {Shelstad}, D.},
 		title = "{On Splitting Invariants and Sign Conventions in Endoscopic Transfer}",
 		journal = {arXiv e-prints},
 		year = 2012,
 		month = jan,
 		eid = {arXiv:1201.5658},
 		pages = {arXiv:1201.5658},
 		archivePrefix = {arXiv},
 		eprint = {1201.5658},
 		primaryClass = {math.RT},
 		adsurl = {https://ui.adsabs.harvard.edu/abs/2012arXiv1201.5658K}
 	}

@article {LMW-TFL,
	AUTHOR = {Lemaire, Bertrand and Moeglin, Colette and Waldspurger,
	Jean-Loup},
	TITLE = {Le lemme fondamental pour l'endoscopie tordue: r\'{e}duction aux
	\'{e}l\'{e}ments unit\'{e}s},
	JOURNAL = {Ann. Sci. \'{E}c. Norm. Sup\'{e}r. (4)},
	FJOURNAL = {Annales Scientifiques de l'\'{E}cole Normale Sup\'{e}rieure. Quatri\`eme
	S\'{e}rie},
	VOLUME = {51},
	YEAR = {2018},
	NUMBER = {2},
	PAGES = {281--369},
	ISSN = {0012-9593},
	MRCLASS = {11F72},
	MRNUMBER = {3798304},
	MRREVIEWER = {Micha\l  Zydor},
	DOI = {10.24033/asens.2357},
	URL = {https://doi-org.proxy-um.researchport.umd.edu/10.24033/asens.2357},
}

@incollection {LW-TFL,
	AUTHOR = {Lemaire, Bertrand and Waldspurger, Jean-Loup},
	TITLE = {Le lemme fondamental pour l'endoscopie tordue: le cas o\`u le
	groupe endoscopique elliptique non ramifi\'{e} est un tore},
	BOOKTITLE = {Representation theory, number theory, and invariant theory},
	SERIES = {Progr. Math.},
	VOLUME = {323},
	PAGES = {399--468},
	PUBLISHER = {Birkh\"{a}user},
	YEAR = {2017},
	MRCLASS = {22E50},
	MRNUMBER = {3753919},
}

@book{SGA4-3,
	AUTHOR = {Artin, Michael and Grothendieck, Alexander}, 
	TITLE = {Th\'{e}orie des topos et cohomologie \'{e}tale des sch\'{e}mas. {T}ome 3},
	SERIES = {Lecture Notes in Mathematics, Vol. 305},
	NOTE = {S\'{e}minaire de G\'{e}om\'{e}trie Alg\'{e}brique du Bois-Marie 1963--1964 (SGA 4),
	Dirig\'{e} par M. Artin, A. Grothendieck et J. L. Verdier. Avec la
	collaboration de P. Deligne et B. Saint-Donat},
	PUBLISHER = {Springer-Verlag, Berlin-New York},
	YEAR = {1973},
	PAGES = {vi+640},
	MRCLASS = {14-06},
	MRNUMBER = {0354654},
	}

@article {Zhu-orthogonal,
		AUTHOR = {Zhu, Yihang},
		TITLE = {The stabilization of the {F}robenius-{H}ecke traces on the
		intersection cohomology of orthogonal {S}himura varieties},
		JOURNAL = {Ast\'erisque},
		FJOURNAL = {Ast\'erisque},
		volume = {453},
		YEAR = {2024},
		PAGES = {ix+292},
		ISSN = {0303-1179,2492-5926},
		ISBN = {978-2-37905-203-3},
		MRCLASS = {11G18 (11F80 11G40 14G35)},
		MRNUMBER = {4838820},
		MRREVIEWER = {Ivan\ Mati\'c},
		DOI = {10.24033/ast.1225},
		URL = {https://doi.org/10.24033/ast.1225},
	}

@article {Peng,
    AUTHOR = {Peng, Zhifeng},
     TITLE = {Multiplicity formula and stable trace formula},
   JOURNAL = {Amer. J. Math.},
  FJOURNAL = {American Journal of Mathematics},
    VOLUME = {141},
      YEAR = {2019},
    NUMBER = {4},
     PAGES = {1037--1085},
      ISSN = {0002-9327},
   MRCLASS = {11R39 (11F72 22E50)},
  MRNUMBER = {3992572},
MRREVIEWER = {Feng-Wen An},
       DOI = {10.1353/ajm.2019.0027},
       URL = {https://doi-org.libproxy.berkeley.edu/10.1353/ajm.2019.0027},
}

@article {Scholze:LK,
    AUTHOR = {Scholze, Peter},
     TITLE = {The {L}anglands-{K}ottwitz method and deformation spaces of
              {$p$}-divisible groups},
   JOURNAL = {J. Amer. Math. Soc.},
  FJOURNAL = {Journal of the American Mathematical Society},
    VOLUME = {26},
      YEAR = {2013},
    NUMBER = {1},
     PAGES = {227--259},
}

@article {KS-GSO,
	AUTHOR = {Kret, Arno and Shin, Sug Woo},
	TITLE = {Galois representations for even general special orthogonal
	groups},
	JOURNAL = {J. Inst. Math. Jussieu},
	FJOURNAL = {Journal of the Institute of Mathematics of Jussieu. JIMJ.
	Journal de l'Institut de Math\'ematiques de Jussieu},
	VOLUME = {23},
	YEAR = {2024},
	NUMBER = {5},
	PAGES = {1959--2050},
	ISSN = {1474-7480,1475-3030},
	MRCLASS = {11F80 (11F70 11G18 11R39)},
	MRNUMBER = {4821551},
	MRREVIEWER = {Neven\ Grbac},
	DOI = {10.1017/S1474748023000427},
	URL = {https://doi.org/10.1017/S1474748023000427},
}

@article {Shin-MC,
 	AUTHOR = {Mack-Crane, Sander and Shin, Sug Woo},
 	TITLE = {Counting points on {I}gusa varieties of {H}odge type},
 	JOURNAL = {Math. Res. Lett.},
 	FJOURNAL = {Mathematical Research Letters},
 	VOLUME = {32},
 	YEAR = {2025},
 	NUMBER = {5},
 	PAGES = {1597--1665},
 	ISSN = {1073-2780,1945-001X},
 	MRCLASS = {11G25 (11G18)},
 	MRNUMBER = {5009221},
 	DOI = {10.4310/mrl.251215195034},
 	URL = {https://doi.org/10.4310/mrl.251215195034},
 }

@article {HainesRicharz1,
    AUTHOR = {Haines, Thomas J. and Richarz, Timo},
     TITLE = {The test function conjecture for parahoric local models},
   JOURNAL = {J. Amer. Math. Soc.},
  FJOURNAL = {Journal of the American Mathematical Society},
    VOLUME = {34},
      YEAR = {2021},
    NUMBER = {1},
     PAGES = {135--218},
}

@article {HainesRicharz2,
    AUTHOR = {Haines, Thomas J. and Richarz, Timo},
     TITLE = {The test function conjecture for local models of
              {W}eil-restricted groups},
   JOURNAL = {Compos. Math.},
  FJOURNAL = {Compositio Mathematica},
    VOLUME = {156},
      YEAR = {2020},
    NUMBER = {7},
     PAGES = {1348--1404},
}

@article {LTXZZ,
	AUTHOR = {Liu, Yifeng and Tian, Yichao and Xiao, Liang and Zhang, Wei
	and Zhu, Xinwen},
	TITLE = {On the {B}eilinson-{B}loch-{K}ato conjecture for
	{R}ankin-{S}elberg motives},
	JOURNAL = {Invent. Math.},
	FJOURNAL = {Inventiones Mathematicae},
	VOLUME = {228},
	YEAR = {2022},
	NUMBER = {1},
	PAGES = {107--375},
	ISSN = {0020-9910,1432-1297},
	MRCLASS = {11G40 (11G05 11G18 11R34)},
	MRNUMBER = {4392458},
	MRREVIEWER = {A.\ Raghuram},
	DOI = {10.1007/s00222-021-01088-4},
	URL = {https://doi.org/10.1007/s00222-021-01088-4},
}

@article{Cornut,
  AUTHOR = {Cornut, Christophe},
  TITLE = {An {E}uler System of {H}eegner Type},
 journal={preprint},
 note={Available at \url{https://webusers.imj-prg.fr/~christophe.cornut/papers/ESHT.pdf}},
 year = {2018}
}

@article {ZhouParahoric,
    AUTHOR = {Zhou, Rong},
     TITLE = {Mod {$p$} isogeny classes on {S}himura varieties with
              parahoric level structure},
   JOURNAL = {Duke Math. J.},
  FJOURNAL = {Duke Mathematical Journal},
    VOLUME = {169},
      YEAR = {2020},
    NUMBER = {15},
     PAGES = {2937--3031},
}

@article {Matsushima,
    AUTHOR = {Matsushima, Yoz\^{o}},
     TITLE = {A formula for the {B}etti numbers of compact locally symmetric
              {R}iemannian manifolds},
   JOURNAL = {J. Differential Geometry},
  FJOURNAL = {Journal of Differential Geometry},
    VOLUME = {1},
      YEAR = {1967},
     PAGES = {99--109},
      ISSN = {0022-040X},
   MRCLASS = {57.31},
  MRNUMBER = {222908},
MRREVIEWER = {S. Murakami},
       URL = {http://projecteuclid.org/euclid.jdg/1214427883},
}

@incollection {HCsubmersion,
    AUTHOR = {Harish-Chandra},
     TITLE = {A submersion principle and its applications},
 BOOKTITLE = {Geometry and analysis},
     PAGES = {95--102},
 PUBLISHER = {Indian Acad. Sci., Bangalore},
      YEAR = {1980},
   MRCLASS = {22E45 (22E50 43A65)},
  MRNUMBER = {592255},
MRREVIEWER = {P. C. Trombi},
}

@article {Taibi,
    AUTHOR = {Ta{\"i}bi, Olivier},
     TITLE = {Arthur's multiplicity formula for certain inner forms of
              special orthogonal and symplectic groups},
   JOURNAL = {J. Eur. Math. Soc. (JEMS)},
  FJOURNAL = {Journal of the European Mathematical Society (JEMS)},
    VOLUME = {21},
      YEAR = {2019},
    NUMBER = {3},
     PAGES = {839--871},
      ISSN = {1435-9855},
   MRCLASS = {11R39 (22E55)},
  MRNUMBER = {3908767},
MRREVIEWER = {\Dbar \cftil{o} Ng\d{o}c Di\cfudot{e}p},
       DOI = {10.4171/JEMS/852},
       URL = {https://doi.org/10.4171/JEMS/852},
}

@misc{AGIKMS,
	title={Local Intertwining Relations and Co-tempered {A}-packets of Classical Groups}, 
	author={Hiraku Atobe and Wee Teck Gan and Atsushi Ichino and Tasho Kaletha and Alberto Mínguez and Sug Woo Shin},
	year={2026},
	note = {version 3,  arXiv:2410.13504,  \url{https://arxiv.org/abs/2410.13504}}
}

@article {Shelstad2012,
    AUTHOR = {Shelstad, D.},
     TITLE = {On geometric transfer in real twisted endoscopy},
   JOURNAL = {Ann. of Math. (2)},
  FJOURNAL = {Annals of Mathematics. Second Series},
    VOLUME = {176},
      YEAR = {2012},
    NUMBER = {3},
     PAGES = {1919--1985},
      ISSN = {0003-486X,1939-8980},
   MRCLASS = {22E45 (22E50)},
  MRNUMBER = {2979862},
MRREVIEWER = {Anne-Marie\ H.\ Aubert},
       DOI = {10.4007/annals.2012.176.3.9},
       URL = {https://doi.org/10.4007/annals.2012.176.3.9},
}

@article {ArtTransferWeighted,
    AUTHOR = {Arthur, James},
     TITLE = {On the transfer of distributions: weighted orbital integrals},
   JOURNAL = {Duke Math. J.},
  FJOURNAL = {Duke Mathematical Journal},
    VOLUME = {99},
      YEAR = {1999},
    NUMBER = {2},
     PAGES = {209--283},
      ISSN = {0012-7094,1547-7398},
   MRCLASS = {22E55 (11F70 11R39 22E50)},
  MRNUMBER = {1708030},
MRREVIEWER = {Volker\ J.\ Heiermann},
       DOI = {10.1215/S0012-7094-99-09909-X},
       URL = {https://doi-org.libproxy.berkeley.edu/10.1215/S0012-7094-99-09909-X},
}

\cleardoublepage

\printindex
\printindex[Oct2026-sorted-notation]

\end{document}